\documentclass[english]{smfbook}
\usepackage{smfenum,sabbah-twistor_041220}
\xyoption{dvips}
\makeindex

\begin{document}
\selectlanguage{french}
\selectlanguage{english}

\frontmatter
\title{Polarizable twistor $\mathcal D$-modules}
\alttitle{$\mathcal D$-modules avec structure de twisteur polarisable}

\author[C.~Sabbah]{Claude Sabbah}
\address{UMR 7640 du CNRS\\
Centre de Math\'ematiques Laurent Schwartz\\
\'Ecole polytechnique\\
F--91128 Palaiseau cedex\\
France}
\email{sabbah@math.polytechnique.fr}
\urladdr{http://www.math.polytechnique.fr/cmat/sabbah/sabbah.html}

\thanks{The first part of this work was achieved under the Indo-French program IFCPAR, and the other part under the INTAS program 97-1664}

\begin{abstract}
We prove a Decomposition Theorem for the direct image of an irreducible local system on a smooth complex projective variety under a morphism with values in another smooth complex projective variety. For this purpose, we construct a category of polarized twistor $\mathcal{D}$-modules and show a Decomposition Theorem in this category.
\end{abstract}

\begin{altabstract}
Nous montrons un th\'eor\`eme de d\'ecomposition pour l'image directe d'un syst\`eme local irr\'eductible sur une vari\'et\'e projective complexe lisse par un morphisme \`a valeurs dans une autre vari\'et\'e projective complexe lisse. \`A cet effet, nous construisons une cat\'egorie de $\mathcal{D}$-modules avec structure de twisteur polaris\'ee et nous montrons un th\'eor\`eme de d\'ecomposition dans cette cat\'egorie.
\end{altabstract}


\altkeywords{Structure de twisteur, $\mathcal{D}$-module, sp\'ecialisation, polarisation, poids}

\maketitle

\chapterspace{2}
\tableofcontents

\mainmatter

\chapterspace{-3}
\chapter*{Introduction}

Let $X$ be a smooth complex projective manifold and let $\cF$ be a locally constant sheaf of $\CC$-vector spaces of finite dimension on $X$. Assume that $\cF$ is \emph{semisimple}, \ie a direct sum of irreducible locally constant sheaves on $X$. Then it is known that, given any ample line bundle on $X$, the corresponding Hard Lefschetz Theorem holds for the cohomology of $X$ with values in $\cF$: if $\cF$ is constant, this follows from Hodge theory; for general semisimple local systems, this was proved by C\ptbl Simpson \cite{Simpson92} using the existence of a harmonic metric \cite{Corlette88}. The existence of such a metric also allows him to show easily that the restriction of $\cF$ to any smooth subvariety of $X$ remains semisimple.

In this article, we extend to such semisimple local systems other properties known to be true for the constant sheaf, properties usually deduced from Hodge theory. These properties will concern the behaviour with respect to morphisms. They were first proved for the constant sheaf (\cf \cite{Deligne68,DeligneHII,Schmid73,Steenbrink76,B-B-D81,G-N90}) and then, more generally, for local systems underlying a polarizable Hodge Module, as a consequence of the work of M\ptbl Saito \cite{MSaito86}.

Given a local system $\cF$ of finite dimensional $\CC$-vector spaces on a complex manifold $X$, it will be convenient to denote by $\pcF$ the associated perverse complex $\cF[\dim X]$, \ie the complex having $\cF$ as its only nonzero term, this term being in degree $-\dim X$.

The proof of the following results will be given in \T\ref{sec:main}.

\begin{thm}[Decomposition Theorem]\label{thm:decomp}
Let $X$ be a smooth complex projective variety and let $\cF$ be a semisimple local system of finite dimensional $\CC$-vector spaces on~$X$. Let $U$ be an open set of $X$ and let $f:U\to Y$ be a proper holomorphic mapping in a complex manifold $Y$. Fix an ample line bundle on $X$. Then
\begin{enumerate}
\item\label{thm:decomp1}
the relative Hard Lefschetz Theorem holds for the perverse cohomology sheaves ${}^p\!\cH^i(\bR f_*\pcF_{|U})$ of the direct image;
\item\label{thm:decomp2}
the direct image complex $\bR f_*\pcF_{|U}$ decomposes (maybe non canonically) as the direct sum of its perverse cohomology sheaves:
\[
\bR f_*\pcF_{|U}\simeq\ooplus_i {}^p\!\cH^i(\bR f_*\pcF_{|U})[-i];
\]
\item\label{thm:decomp3}
each perverse cohomology sheaf ${}^p\!\cH^i(\bR f_*\pcF_{|U})$ decomposes as the direct sum of intersection complexes supported on closed irreducible analytic subsets $Z$ of $Y$, \ie of the form $\IC^{\cbbullet}(\pcL)$, where $\cL$ is a local system on a smooth open dense set $Z\moins Z'$, with $Z'$ closed analytic in $Z$;
\item\label{thm:decomp4}
if moreover $U=X$ and $Y$ is projective, then each perverse cohomology sheaf ${}^p\!\cH^i(\bR f_*\pcF)$ is \emph{semisimple}, \ie the local systems $\cL$ are semisimple.
\end{enumerate}
\end{thm}

\begin{thm}[Vanishing cycles] \label{thm:vanishcycles}
Let $X$ be a smooth complex projective variety and let $\cF$ be a semisimple local system on $X$. Let $U$ be an open set of $X$ and let $f:U\to\CC$ be a holomorphic function on $U$ which is proper. Then, for any $\ell\in\ZZ$, the perverse complexes $\gr_\ell^{\rM}\psip_f\pcF$ and $\gr_\ell^{\rM}\phip_f\pcF$, obtained by grading with respect to the monodromy filtration the perverse complexes of nearby or vanishing cycles, are semisimple perverse sheaves on $f^{-1}(0)$.
\end{thm}

\begin{Remarques*}
\begin{enumerate}
\item
Notice that $\eqref{thm:decomp1}\implique\eqref{thm:decomp2}$ in Main Theorem \ref{thm:decomp} follows from an argument of Deligne \cite{Deligne68}.

\item
The nearby and vanishing cycles functors $\psi_f$ and $\phi_f$ defined by Deligne \cite{Deligne73} are shifted by $-1$, so that they send perverse sheaves to perverse sheaves. They are denoted by $\psip_f$ and $\phip_f$, following M\ptbl Saito \cite{MSaito86}.

\item
It is known that the Main Theorem \ref{thm:decomp} implies the \emph{local invariant cycle theorem} for the cohomology with coefficients in $\cF$ (\cf \cite[Cor\ptbl 6.2.8 and 6.2.9]{B-B-D81}, see also \cite[Cor\ptbl 3.6 and 3.7]{MSaito90}). If for instance $Y=\CC$ then, for any $k\geq0$ and for $t\neq0$ small enough, there is an exact sequence
\[
H^k(f^{-1}(0),\cF)\to H^k(f^{-1}(t),\cF)\To{T-\id}H^k(f^{-1}(t),\cF),
\]
where $T$ denotes the monodromy. It also implies the exactness of the \emph{Clemens-Schmid sequence}.

\item
Owing to the fact that, if $\cF_\QQ$ is a perverse complex of $\QQ$-vector spaces on a complex analytic manifold, then $\cF_\QQ$ is semisimple if and only if $\cF_\CC=\CC\otimes_\QQ\cF_\QQ$ is so, the previous results apply as well to $\QQ$-local systems, giving semisimple $\QQ$-perverse complexes as a result.

\item
It would be possible to define a category of perverse complexes ``of smooth origin'', obtained after iterating various operations starting from a semisimple local system on a smooth complex projective variety, \eg taking perverse cohomology of a projective direct image, taking monodromy-graded nearby or vanishing cycles relative to a projective holomorphic function, taking sub-quotients of such objects. The perverse complexes in this category are semisimple.

\item
A conjecture of M\ptbl Kashiwara \cite{Kashiwara98} ---which was the main motivation for this work--- asserts in particular that these results should hold when $\cF$ is any semisimple \emph{perverse sheaf} (with coefficients in $\CC$) on $X$. In the complex situation that we consider, they are proved when $\cF$ underlies a \emph{polarizable Hodge Module}, \ie if on a smooth dense open set of its support, the perverse sheaf $\cF$ is (up to a shift) a local system defined over $\QQ$ or $\RR$ underlying a variation of polarized Hodge structures defined over $\QQ$ or $\RR$: this is a consequence of the work of M\ptbl Saito \cite{MSaito86,MSaito87} and \cite{C-K-S87,K-K87}, and of the known fact (see \cite{Deligne84}) that, on a smooth Zariski open set of a projective variety, the local system underlying a variation of complex Hodge structures is semisimple.

Let us indicate that the conjecture of Kashiwara is even more general, as it asserts that analogues of such results should be true for semisimple holonomic $\cD$-modules on smooth complex projective varieties. However, we will not seriously consider nonregular $\cD$-modules in this article.

\item
First were proved the arithmetic analogues of these theorems, \ie for ``pure sheaves'' instead of semisimple sheaves (\cf \cite{B-B-D81}) and they were used to give the first proof of the Decomposition Theorem for the constant sheaf in the complex case. An arithmetic approach to the conjecture of Kashiwara (at least for $\CC$-perverse sheaves) has recently been proposed by V\ptbl Drinfeld \cite{Drinfeld01}.

\item
It should be emphasized that we work with \emph{global} properties on a projective variety, namely, semisimplicity. Nevertheless, the main idea in the proof is to show that these global properties may be expressed by \emph{local} ones, \ie by showing that each irreducible local system on $X$ underlies a variation of some structure, analogous to a polarized Hodge structure, called a \emph{polarized twistor structure}. Extending this to irreducible perverse sheaves is the contents of Conjecture \ref{conj:ssimple}.

\item
It will be more convenient to work with the category of regular holonomic $\cD_X$-modules instead of that of $\CC$-perverse sheaves on $X$. It is known that both categories are equivalent \emph{via} the de~Rham functor, and that this equivalence is compatible with the corresponding direct image functors or with the nearby and vanishing cycles functors. We will freely use this compatibility.
\end{enumerate}
\end{Remarques*}

Let us now give some explanation on the main steps of the proof. We will use three sources of ideas:
\begin{enumerate}
\item
the theory of \emph{twistor structures} developed by C\ptbl Simpson (after ideas of P\ptbl De\-li\-gne),
\item
the techniques developed by M\ptbl Saito in the theory of polarizable Hodge Modules,
\item
the use of distributions and Mellin transform, as inspired by the work of M\ptbl Kashiwara and D\ptbl Barlet.
\end{enumerate}

One of the main objectives, when trying to prove a decomposition theorem, is to develop a notion of \emph{weight} satisfying good properties with respect to standard functors. In other words, the category of semisimple local systems (or, better, semisimple perverse sheaves) should satisfy the properties that one expects for pure sheaves. If the Hodge structure contains in its very definition such a notion, it is not clear \emph{a~priori} how to associate a weight to an irreducible perverse sheaf: one could give it weight~$0$, but one should then explain why ${}^p\!\cH^i(\bR f_*\cF)$ has weight~$i$ for instance. On the other hand, it is natural to expect that, if a notion of ``pure sheaf'' exists in the complex setting, it should be more general than that of polarized Hodge Modules, and even of that of ``pure perverse sheaf''. Indeed, in the arithmetic situation, one is able to treat sheaves with wild ramification (\eg Fourier transform of pure sheaves with moderate ramification).

The very nice idea of a twistor structure allows one to work with the notion of weight. Let us quickly explain it, referring to \cite{Simpson97} (see also \S\T\ref{sec:twst0} and \ref{sec:smtwqc} below) for a more detailed presentation. Let $(V,\nabla)$ be a \emph{flat} holomorphic vector bundle on a smooth manifold $X$, that we view as a $C^\infty$ vector bundle $H$ on $X$ equipped with a flat $C^\infty$ connection $D$, the holomorphic structure being given by the $(0,1)$ part $D''$ of the connection. A twistor structure of weight~$w\in\ZZ$ on $(V,\nabla)$ (or a variation of twistor structures) consists of the datum of a $C^\infty$ vector bundle $\wt\cH$ on $X\times \PP^1$, holomorphic with respect to the variable of $\PP^1$, equipped with relative connections $\wt D', \wt D''$ (\ie there is no derivation with respect to the $\PP^1$ variable), with poles along $X\times\{0\}$ and $X\times\{\infty\}$ respectively, and such that the restriction of $\wt\cH$ to any $\{x_o\}\times\PP^1$ is isomorphic to $\cO_{\PP^1}(w)^{\rg V}$ (see \T\ref{sec:smtwqc} for a more precise definition, in particular for $\wt D', \wt D''$).

Therefore, a variation of twistor structures on $X$ lives on $X\times\PP^1$. One of the main properties required for weights, namely that there is no nonzero morphism from an object of weight~$w$ to an object of weight~$w'<w$, follows from the analogous property for line bundles on $\PP^1$. One may also define the notion of polarization (see \loccit).

The main device to produce a variation of polarized twistor structures on a holomorphic flat bundle $(V,\nabla)$ is given by the construction of a \emph{harmonic metric}. It follows from a theorem of K\ptbl Corlette \cite{Corlette88} and C\ptbl Simpson \cite{Simpson92} that a local system $\cF$ of $\CC$-vector spaces on a compact K\"ahler manifold $X$ ``underlies'' a variation of polarized twistor structures if and only if it is semisimple, because semisimplicity is a necessary and sufficient condition to build on the flat bundle $(V,\nabla)$ associated with~$\cF$ a harmonic metric.

The next step closely follows ideas of M\ptbl Saito \cite{MSaito86}, namely it consists in defining in its own right a category of ``singular variations of polarized twistor structures''. This is done \emph{via} the theory of $\cD$-modules, and more precisely \emph{via} the theory of $\cR$-modules, which is a natural extension to $X\times (\PP^1\moins\{\infty\})$ of the theory of $\cD_X$-modules. In order to keep some control on the coherence properties, we are not allowed to use $C^\infty$ coefficients. Therefore, we modify a little bit the presentation of the object $\wt\cH$ introduced above, as associated to the left $\cD_X$-module $(V,\nabla)$.

Put $\Omega_0=\PP^1\moins\{\infty\}$ with coordinate $\hb$ and $\Omega_\infty=\PP^1\moins\{0\}$. We may view $\wt\cH$ as the result of a $C^\infty$ gluing between
$\wt\cH_{|X\times\Omega_0}$ and $\wt\cH_{|X\times\Omega_\infty}$ on some neighbourhood of $X\times\bS$, where $\bS$ denotes the circle $\module{\hb}=1$. Equivalently, denoting by $\wt\cH^\vee$ the dual bundle, the gluing may be viewed as a nondegenerate \emph{pairing} on $\wt\cH_{|X\times\bS}^\vee\otimes\wt\cH_{|X\times\bS}$ with values in the sheaf of $C^\infty$ functions on some neighbourhood of $X\times\bS$ which are holomorphic with respect to $\hb$. We may restrict this pairing to the holomorphic/antiholomorphic object $\cHS'\otimes\ov{\cHS''}$, where we put
$$
\cHS'=\ker\big[\wt D'':\wt\cH_{|X\times\bS}^\vee\to\wt\cH_{|X\times\bS}^\vee\big]\qqbox{and} \ov{\cHS''}=\ker \big[\wt D':\wt\cH_{|X\times\bS}\to\wt\cH_{|X\times\bS}\big].
$$
Its restriction to any $\{x_o\}\times\PP^1$ should define a bundle $\cO_{\PP^1}(w)^{\rg V}$.

We extend this construction to $\cD_X$-modules as follows: the basic objects are triples $(\cM',\cM'',C)$, where $\cM',\cM''$ are coherent $\cR_{X\times\Omega_0}$-modules (see \S\T\ref{subsec:debutRX}--\ref{num:RX} for the definition of the sheaf $\cR$). To any $\cR_{X\times\Omega_0}$-module $\cM$ is associated a ``conjugate'' object $\ov\cM$, which is a coherent $\cR_{\ov X\times\Omega_\infty}$-module (here, $\ov X$ is the complex conjugate manifold); now, $C$ is a \emph{pairing} on $\cM'_{|X\times\bS}\otimes_{\cO_\bS}\ov{\cM''_{|X\times\bS}}$ which takes values in \emph{distributions} on $X\times\bS$ which are continuous with respect to $\hb$. A \emph{polarization} will then appear as an isomorphism $\cM''\isom\cM'$ of $\cR_{X\times\Omega_0}$-modules.

\begin{exemple*}
Given any $\CC$-vector space $H$, denote by $\ov H$ its complex conjugate and by $H^\vee$ its dual. Define $H'=H^\vee$ and $H''=\ov H$. There is a natural pairing (\ie $\CC$-linear map) $H'\otimes_\CC\ov{H''}\to\CC$, induced by the natural duality pairing $H^\vee\otimes_\CC H \to\CC$.

On the other hand, consider the category of triples $(H',H'',C)$, where $H',H''$ are $\CC$-vector spaces and $C$ is a nondegenerate pairing $H'\otimes_\CC\ov{H''}\to\CC$; morphisms $\varphi:(H_1',H_1'',C_1)\to(H'_2,H_2'',C_2)$ are pairs $\varphi=(\varphi',\varphi'')$ with $\varphi':H'_2\to H'_1$, $\varphi'':H''_1\to H''_2$ such that $C(\varphi'(m'_2),\ov{m''_1})=C(m'_2,\ov{\varphi''(m''_1)})$.

We have constructed above a functor from the category of $\CC$-vector spaces to this category of triples. It is easily seen to be an equivalence.

Under this equivalence, the Hermitian dual $H^*=\ov H^\vee$ of $H$ corresponds to $(H',H'',C)^*\defin(H'',H',C^*)$ with $C^*(m'',\ov{m'})\defin \ov{C(m',\ov{m''})}$, and a sesquilinear form on $H$, which is nothing but a morphism $\cS:H^*\to H$, corresponds to a morphism $\cS:(H',H'',C)^*\to (H',H'',C)$, \ie a pair $(S',S'')$ with $S',S'':H''\to H'$ such that $C(S'm'',\ov{\mu''})=\ov{C(S''\mu'',\ov{m''})}$.
\end{exemple*}

In order to say that the pairing $C$ on $\cM'_{|X\times\bS}\otimes\ov{\cM''_{|X\times\bS}}$ is holomorphic and nondegenerate, and therefore defines a ``gluing'', we should be able to restrict it to $\{x_o\}\times\bS$ for any $x_o\in X$. ``Restriction'' is understood here under the broader sense of ``taking nearby or vanishing cycles''. Hence, in order to ``restrict'' $\cM'$ or $\cM''$, we impose that they have a Malgrange-Kashiwara filtration, \ie admit Bernstein polynomials. In order to restrict the pairing $C$, we use a device developed by D\ptbl Barlet in a nearby context, namely by taking residues of Mellin transforms of distributions.

\medskip
The main technical result is then the construction of the category of regular polarized twistor $\cD$-modules, mimicking that of polarized Hodge Modules \cite{MSaito86}, and the proof of a decomposition theorem in this category (Theorem \ref{th:imdirtwistor}).

To conclude with a proof of M\ptbl Kashiwara's conjecture for semisimple perverse sheaves, one should prove that the functor which associates to each regular polarized twistor $\cD_X$-module $(\cM,\cM,C)$ of weight~$0$ (the polarization is $\id:\cM\to\cM$) the $\cD_X$-module $\cM_{|X\times\{1\}}$ is an equivalence with the subcategory of semisimple regular holonomic $\cD_X$-modules, when $X$ is a complex projective manifold.

We are not able to prove this equivalence in such a generality. However, we prove the equivalence for smooth objects, and also when $X$ is a curve. According to a Zariski-Lefschetz Theorem due to H\ptbl Hamm and L\^e~D.T. \cite{H-L85}, and using the Riemann-Hilbert correspondence, this implies at least that the functor above takes values in the category of semisimple regular holonomic $\cD_X$-modules. This is enough to get the Main Theorems.

\enlargethispage{3mm}
\smallskip
What is the overlap with M\ptbl Saito's theory of polarizable Hodge Modules? The main difference with M\ptbl Saito's theory consists in the way of introducing the polarization.

The method of M\ptbl Saito is ``\`a la Deligne'', using a perverse complex defined over $\QQ$ or~$\RR$ (and the de~Rham functor from holonomic $\cD$-modules to $\CC$-perverse sheaves) to get the rational or real structure. The polarization is then introduced at the topological level (perverse complexes) as a bilinear form, namely the Poincar\'e-Verdier duality. We do not know whether such an approach would be possible for polarizable twistor $\cD$-modules.

Here, we use a purely analytical approach ``\`a la Griffiths'', without paying attention to the possible existence of a $\QQ$- or $\RR$-structure. The polarization is directly introduced as a Hermitian form. In particular, we do not use the duality functor and we do not need to show various compatibilities with the de~Rham functor. This approach uses therefore less derived category techniques than the previous one. Moreover, it is possible (\cf \T\ref{subsec:HodgeDM}) to introduce a category of polarizable Hodge $\cD$-modules as a subcategory of twistor $\cD$-modules, by considering those twistor $\cD$-modules which are invariant under the natural $\CC^*$ action on the category (similarly to what C\ptbl Simpson does for ``systems of Hodge bundles'' \cite{Simpson88,Simpson92}). This gives a generalization of \emph{complex} variations of Hodge structures (without real structure). We do not know if this category is equivalent to the category one gets by M\ptbl Saito's method, but this can be expected. A similar category, that of \emph{integrable} twistor $\cD$-modules, is considered in Chapter~\ref{chap:int}.

\smallskip
What is the overlap with C\ptbl Simpson's study of Higgs bundles \cite{Simpson92}? First, notice that we consider objects which may have complicated singularities, so we do not consider any question concerning moduli. We are mainly interested in the functor sending a twistor $\cD$-module to its associated $\cD$-module by restricting to $\hb=1$ ($\hb$ is the standard name we use for the variable on $\PP^1$). We could also consider its associated Higgs module by restricting to $\hb=0$ (see \T\ref{subsec:Bernstein}). In the first case, we at least know the image category, namely that of semisimple regular holonomic $\cD$-modules. In the second case, we have no idea of how to characterize the image of the functor and if an equivalence could be true, similarly to what is done in the smooth case by C\ptbl Simpson \cite{Simpson92} or in a slightly more general case by O\ptbl Biquard \cite{Biquard97}.

\medskip
Let us now describe with more details the contents of this article.

In Chapter \ref{chap:RX}, we give the main properties of $\cR_\cX$-modules. They are very similar to that of $\cD_X$-modules. The new objects are the sesquilinear pairing $C$, and the category $\RTriples(X)$ (the objects are triples formed with two $\cR_\cX$-modules and a sesquilinear pairing between them), to which we extend various functors. We have tried to be precise concerning signs.

Chapter \ref{chap:smtw} introduces the notion of a (polarized) twistor structure, following C\ptbl Simpson \cite{Simpson97}. We first consider the case where the base $X$ is a point, to get the analogue of a (polarized) Hodge structure. We develop the notion of a Lefschetz twistor structure and adapt to this situation previous results of M\ptbl Saito and P\ptbl Deligne. Last, we develop the notion of a smooth twistor structure on a smooth complex manifold $X$. The main point of this chapter is to express the notion of a twistor structure in the frame of the category $\RTriples$, in order to extend this notion to arbitrary holonomic $\cR_\cX$-modules.

Chapter \ref{chap:spe} extends to $\cR_\cX$-modules the notion of specializability along a hyper\-sur\-face---a notion introduced by B\ptbl Malgrange and M\ptbl Kashiwara for $\cD_X$-module, together with the now called Malgrange-Kashiwara filtration---and analyzes various properties of the nearby and vanishing cycles functors. The specialization of a sesquilinear pairing is then defined by means of the residue of a Mellin transform, in analogy with some works of D\ptbl Barlet. All together, this defines the notion of a specializable object in the category $\RTriples(X)$. The category of S-decomposable objects, introduced in \T\ref{sec:Sdecomp}, is inspired from \cite{MSaito86}.

In Chapter \ref{chap:twistor}, we introduce the category of twistor $\cD$-modules on $X$ as a subcategory of $\RTriples(X)$. We prove various property of the category of (polarized) twistor $\cD$-modules, analogous to that of (polarized) Hodge Modules \cite{MSaito86}. We show that regular twistor $\cD$-modules induce semisimple regular $\cD$-modules by the de~Rham functor $\Xi_{\DR}$.

Chapter \ref{chap:curves} establishes the equivalence between regular twistor $\cD$-modules and semisimple perverse sheaves (or semisimple regular holonomic $\cD$-modules) on compact Riemann surfaces, by expressing the results of C\ptbl Simpson \cite{Simpson90} and O\ptbl Biquard \cite{Biquard97} in the frame of polarized regular twistor $\cD$-modules. In order to establish the equivalence, we also adapt results of D\ptbl Barlet and H.-M\ptbl Maire \cite{B-M89} concerning Mellin transform.

The main theorems are proved in Chapter \ref{chap:decomp}, following the strategy of M\ptbl Saito~\cite{MSaito86}. We reduce the proof to the case where $X$ is a compact Riemann surface and $f$ is the constant map to a point. In this case, we generalize the results of S\ptbl Zucker \cite{Zucker79} to polarizable regular twistor $\cD$-modules.

In Chapter \ref{chap:int}, we consider the category of \emph{integrable} twistor $\cD$-modules. This chapter, written somewhat after the previous ones, is an adaptation to the present theory of the notion of CV-structure considered in \cite{Hertling01}. We mainly prove a ``local unitarity'' statement (the local exponents are real and, in the regular case, the eigenvalues of local monodromies have absolute value equal to one). The interest of such a subcategory should be for the nonregular case, where it should play the role of singular variations of polarized Hodge structures.

In Chapter \ref{chap:8}, we analyze the behaviour of polarized regular twistor $\cD$-modules under a partial (one-dimensional) Fourier-Laplace transform. We generalize to such objects the main result of \cite{Bibi96a}, comparing, for a given function $f$, the nearby cycles at $f=\infty$ and the nearby or vanishing cycles for the partial Fourier-Laplace transform in the $f$-direction (Theorem \ref{th:Fstrict}).

\medskip
Since the first version of this article was written, there has been progress in various directions.
\begin{enumerate}
\item
In the first version of this article, the category of polarized twistor $\cD$-modules was restricted to the \emph{local unitary case}, mainly because of a lack of proof of Theorem \ref{th:L2} in general. This restriction is now unnecessary, due a new proof of this theorem.
\item
The main progress comes from recent work of T\ptbl Mochizuki \cite{Mochizuki03,Mochizuki04}. Continuing \cite{Mochizuki02}, T\ptbl Mochizuki generalizes the contents of Chapter~\ref{chap:curves} in two directions:

-- he considers an arbitrary parabolic structure along the divisor, whereas only a natural parabolic structure is considered here, that we call ``Deligne type''; depending on the point of view, one could call the objects defined by T\ptbl Mochizuki as ``twistor $\cD$-modules with parabolic structure'', or the objects of the present article as ``twistor $\cD$-modules of Deligne type'' (or ``pure imaginary'' after \cite{Mochizuki04}); the category of polarized regular twistor $\cD$-modules that we define here should be (and is, after the work of Mochizuki) equivalent to the category of semisimple perverse sheaves, on a smooth projective variety, whereas twistor $\cD$-modules with parabolic structure give rise to semisimple ``perverse-sheaves-with-parabolic-structure'';

-- he is able to treat the case of the complement of a normal crossing divisor on a smooth complex manifold of arbitrary dimension.

All together, it seems that, in view of the work of O\ptbl Biquard \cite{Biquard97} and J\ptbl Jost and K\ptbl Zuo \cite{J-Z96,J-Z97,Zuo99} (revisited in \cite{Mochizuki04}), the proof of of Conjecture \ref{conj:ssimple}, hence a proof of the conjecture of Kashiwara for perverse sheaves (and even for perverse sheaves ``with parabolic structure'') with analytical methods, is now complete.

\item
On the other hand, according to recent results of G\ptbl Boeckle and C\ptbl Khare \cite{B-K03} or of D\ptbl Gaitsgory \cite{Gaitsgory04}, a proof of the conjecture of de~Jong used by V\ptbl Drinfeld is available; therefore, the arithmetic approach of Drinfeld \cite{Drinfeld01} to the conjecture of Kashiwara (for perverse sheaves) is also complete.

\item
Let us also mention a new proof of the decomposition theorem for the constant sheaf, obtained by M.A\ptbl de\kern2pt Cataldo and L\ptbl Migliorini \cite{C-M03}, with methods completely different from those developed by M\ptbl Saito. We do not know if such methods can be adapted to more general local systems.
\item
The nonregular case of Kashiwara's conjecture is still open. Extending the work of C\ptbl Simpson \cite{Simpson90} and O\ptbl Biquard \cite{Biquard97} to holomorphic bundles on compact Riemann surfaces with meromorphic connections having \emph{irregular} singularities would be a first step. Some results in this direction are obtained in \cite{Bibi98} and \cite{B-B03}.
\end{enumerate}

\subsubsection*{Acknowledgements}
The first part of this work arose from numerous discussions with Nitin Nitsure. It grew up during my visit to the Tata Institute (Bombay), and during a visit of Nitin Nitsure at \'Ecole polytechnique (Palaiseau), both under the Indo-French program IFCPAR. I~also had many discussions with Olivier Biquard, who explained me his work \cite{Biquard97}. M\ptbl Kashiwara and C\ptbl Simpson kindly answered my questions. For Chapter~\ref{chap:int}, I benefitted of many discussions with Claus Hertling. Anonymous referees carefully read the successive versions of this work, suggesting improvements and correcting some mistakes. I~thank all of them.

\setcounter{chapter}{-1}
\begin{chapter}{Preliminaries}
\numberwithin{equation}{subsection}
\def\thesubsection{\thechapter.\arabic{subsection}}

\Subsection{Some signs}\label{subsec:signes}

{\def\theenumi{\alph{enumi}}
\begin{enumerate}
\item\label{subsec:signesb}
We will use the function
\begin{align*}
\ZZ&\To{\varepsilon}\{\pm1\}\\
a&\mto\varepsilon(a)=(-1)^{a(a-1)/2}
\end{align*}
which satisfies in particular
\[
\index{$eps$@$\varepsilon$}\varepsilon(a+1)=\varepsilon(-a)=(-1)^a\varepsilon(a),\quad \varepsilon(a+b)=(-1)^{ab}\varepsilon(a)\varepsilon(b).
\]
Recall that, on $\CC^n$ with coordinates $z_k=x_k+iy_k$ ($k=1,\dots,n$), we have
\[
dz_1\wedge\cdots\wedge dz_n\wedge d\ov z_1\wedge\cdots\wedge d\ov z_n= \varepsilon(n)(dz_1\wedge d\ov z_1)\wedge\cdots\wedge(dz_n\wedge d\ov z_n)
\]
and that $dz_k\wedge d\ov z_k=-2i(dx_k\wedge dy_k)$.

\item\label{subsec:signesc}
We follow the sign convention given in \cite[\S\T0 and 1]{DeligneSGA4}. When we write a multi-complex, we understand implicitly that we take the associated simple complex, ordered as written, with differential equal to the sum of the partial differentials.

Given any sheaf $\cL$, denote by $(\index{$god$@$\God$}\God^{\cbbullet}\cL,\delta)$ the standard semisimplicial resolution of $\cL$ by flabby sheaves, as defined in \cite[Appendice]{Godement64}. For a complex $(\cL^{\cbbullet},d)$, we view $\God^{\cbbullet}\cL^{\cbbullet}$ as a double complex ordered as written, \ie with differential $(\delta_i,(-1)^id_j)$ on $\God^i\cL^j$, and therefore also as the associated simple complex.
\end{enumerate}}

\subsection{}\label{subsec:debutRX}
In this article, $X$ denotes a complex analytic manifold of dimension $\dim X=n$, $\index{$ox$@$\cO_X$}\cO_X$ denotes the sheaf of holomorphic functions on $X$ and $\index{$dx$@$\cD_X$}\cD_X$ the sheaf of linear differential operators with coefficients in $\cO_X$. The sheaf $\cO_X$ is equipped with its natural structure of left $\cD_X$-module and the sheaf $\index{$omegax$@$\omega_X$}\omega_X$ of holomorphic differential $n$-forms with its structure of right $\cD_X$-module.

We denote by $\index{$xbar$@$\ov X$}\ov X$ the complex conjugate manifold, equipped with the structure sheaf $\cO_{\ov X}\defin\ov{\cO_X}$ and by $\index{$xR$@$X_\RR$}X_\RR$ the underlying $C^\infty$ manifold.

The increasing filtration of $\cD_X$ by the order is denoted by $F_\bbullet\cD_X$. Given a filtered object $(M,F_\bbullet M)$ (filtrations are increasing and indexed by $\ZZ$), the associated Rees object $R_FM\defin\oplus_{k\in\ZZ}F_kM\hb^k$ is the graded object constructed with the new variable~$\hb$. In particular, we will consider the Rees ring $R_F\cD_X$: this is a sheaf of rings on $X$. The filtration induced by $F_\bbullet\cD_X$ on $\cO_X$ satisfies $F_k\cO_X=\cO_X$ for $k\geq 0$ and $F_k\cO_X=0$ for $k<0$, so the associated Rees ring $R_F\cO_X$ is equal to $\cO_X[\hb]$.

We will denote by $R_X$ the sheaf $R_F\cD_X$ when we forget its grading, and call it the \emph{differential deformation sheaf}. This is a sheaf of rings on $X$. In local coordinates on~$X$, we denote now by $\index{$dartiall$@$\partiall$}\partiall_{x_i}$ in $R_X$ the element $\hb\partial_{x_i}$ in $R_F\cD_X$. With such a notation, we have
\[
R_X=\cO_X[\hb]\langle\partiall_{x_1},\ldots,\partiall_{x_n}\rangle,
\]
where $\partiall_{x_i}$ satisfy the relations
\[
[\partiall_{x_i},\partiall_{x_j}]=0\qqbox{and}
[\partiall_{x_i},f(x,\hb)]=\hb\frac{\partial f}{\partial x_i}.
\]
One has (forgetting grading)
\[
R_X/\hb R_X=R_F\cD_X/\hb R_F\cD_X=\gr^F\!\cD_X=\cO_X[TX],
\]
where $\cO_X[TX]$ denotes the sheaf of holomorphic functions on the cotangent bundle $T^*X$ which are polynomials with respect to the fibres of $T^*X\to X$; in other words, $\cO_X[TX]=p_*\cO_{\PP(T^*X\oplus\mathbf{1})}(*\infty)$, where $p:\PP(T^*X\oplus \mathbf{1})\to X$ denotes the projection and $\infty$ denotes the section at infinity of the projective bundle.

The left action of $R_X$ on $R_F\cO_X=\cO_X[\hb]$ is defined by
\[
\partiall_{x_i}(f(x,\hb))=\hb\frac{\partial f}{\partial x_i}.
\]

\subsection{}\label{num:RX}
Let $\index{$omega0$@$\Omega_0$}\Omega_0=\CC$ be the complex line equipped with a fixed coordinate $\hb$. We will denote by $\imhb$ the imaginary part of $\hb$. Put $\CC^*=\{\hb\neq0\}$. Denote by $\index{$omega0infty$@$\Omega_\infty$}\Omega_\infty$ the other chart of $\PP^1$, centered at $\hb=\infty$.

Denote by $\index{$xcurl$@$\cX$}\cX$ the product $X\times\Omega_0$ and by $\index{$oxcurl$@$\cO_\cX$}\cO_\cX$ the sheaf of holomorphic functions on it, by $\index{$xcurlcirc$@$\cX^\Cir$}\cX^\Cir$ the product $X\times \CC^*$. Let $\pi:\cX\to X$ denote be the natural projection. If $\cM$ is a sheaf on $\cX$, we denote by $\cM^\Cir$ its restriction to $\cX^\Cir$. We will consider the sheaves
\[
\index{$rxcurl$@$\cR_\cX$}\cR_\cX=\cO_\cX\ootimes_{\cO_X[\hb]}R_F\cD_X\qqbox{and}\cR_{\cX^\Cir}=\cO_{\cX^\Cir}\ootimes_{\cO_X[\hb]}R_F\cD_X.
\]
In local coordinates on $X$, we have $\cR_\cX=\cO_\cX\langle\partiall_{x_1},\dots,\partiall_{x_n}\rangle$.

We will denote by $\index{$thXcurl$@$\Theta_\cX$}\Theta_\cX$ the sheaf of holomorphic vector fields relative to the projection $\cX\to\Omega_0$, which \emph{vanish} at $\hb=0$. This is the $\cO_\cX$-locally free sheaf generated by $\partiall_{x_1},\dots,\partiall_{x_n}$. It is contained in $\cR_\cX$.

Dually, we denote by $\index{$omegaxcurl1$@$\Omega^1_\cX$}\Omega^1_\cX=\hbm\Omega^1_{X\times\Omega_0/\Omega_0}\subset\Omega^1_{X\times\Omega_0/\Omega_0}[\hbm]$ the sheaf of holomorphic $1$-forms on $\cX$ relative to the projection $\cX\to\Omega_0$, which have \emph{a pole of order one at most} along $\hb=0$. We will put $\index{$omegaxcurlk$@$\Omega^k_\cX$}\Omega^k_\cX=\wedge^k\Omega^1_\cX$. The differential $d:\Omega^k_\cX\to\Omega^{k+1}_\cX$ is induced by the relative differential $d=d_{X\times\Omega_0/\Omega_0}$. The natural left multiplication of $\Theta_\cX$ on $\cR_\cX$ can be written as a connection
\[
\nabla:\cR_\cX\to\Omega^1_\cX\ootimes_{\cO_\cX}\cR_\cX
\]
satisfying the Leibniz rule $\nabla(fP)=df\otimes P+f\nabla P$. More generally, a left $\cR_\cX$-module $\cM$ is nothing but a $\cO_\cX$-module with a \emph{flat} connection $\nabla:\cM\to\Omega^1_\cX\otimes_{\cO_\cX}\cM$. Put $\index{$omegaxcurl$@$\omega_\cX$}\omega_\cX\defin\Omega^n_\cX=\hb^{-n}\omega_{X\times\Omega_0/\Omega_0}$. This is naturally a right $\cR_\cX$-module: the action is given by $\omega\cdot\xi=-\cL_\xi\omega$, where $\cL_\xi$ denotes the Lie derivative, here equal to the composition of the interior product $\iota_\xi$ by $\xi$ with the relative differential $d$.

\subsection{}\label{num:bS}
We denote by $\index{$sbold$@$\bS$}\bS$ the circle $\module{\hb}=1$ in $\Omega_0\cap\Omega_\infty$. For a $\cO_\cX$ or $\cR_\cX$-module $\cM$, we denote by~$\cMS$ its sheaf-theoretic restriction to $X\times\bS$. In particular, we will consider the sheaves $\cO_{\cX|\bS}$ and $\cR_{\cX|\bS}$. We will simply denote $\cO_{\Omega_0|\bS}$ by $\index{$osbold$@$\cO_\bS$}\cO_\bS$. We will also use the sheaves\index{$oxxbarsbold$@$\cO_{(X,\ov X)\times\bS}$}\index{$rxxbarsbold$@$\cR_{(X,\ov X)\times\bS}$}
\begin{align*}
\cO_{(X,\ov X)\times\bS}&\defin\cO_{(X\times\Omega_0)|\bS}\otimes_{\cO_\bS}\cO_{(\ov X\times\Omega_0)|\bS}\\
\cR_{(X,\ov X)\times\bS}&\defin\cR_{(X\times\Omega_0)|\bS}\otimes_{\cO_\bS}\cR_{(\ov X\times\Omega_0)|\bS}.
\end{align*}

\subsection{Distributions and currents}\label{num:Dbh}
We will need to consider distributions on $X_\RR\times\bS$ which are regular with respect to the variable on $\bS$, in order to be able to specialize them with respect to the variable of $\bS$. Let us introduce some notation.

Let $T$ be a $C^\infty$ manifold with a fixed volume form $\vol_T$ (we will mainly use $(T,\vol_T)=(\bS,d\arg\hb)$). For $k=0,\dots,\infty$, we denote by $\cC^k_{X_\RR\times T}$ the sheaf of $C^k$ functions on $X_\RR\times T$ and by $\cE^{(n,n)}_{X_\RR\times T/T}$ the sheaf of $C^\infty$ relative (with respect to the projection $X\times T\to T$) $(n,n)$ forms of maximal degree, and we put an index $c$ for those objects with compact support. We denote by $\DbhT{X}^k$ the sheaf on $X_\RR\times T$ of distributions which are $C^k$ with respect to~$T$: by definition, given any open set $W$ of $X_\RR\times T$, an element of $\DbhT{X}^k(W)$ is a $C^\infty(T)$-linear map $\cE^{(n,n)}_{X_\RR\times T/T,c}(W)\to C^k_c(T)$ which is continuous with respect to the usual norm on $C^k_c(T)$ (sup of the modules of partial derivatives up to order $k$) and the family of semi-norms on $\cE^{(n,n)}_{X_\RR\times T/T,c}(W)$ obtained by taking the sup on some compact set of $W$ of the module of partial derivatives up to some order with respect to $X$ and up to order $k$ with respect to $T$. Given a compact set in $W$, the smallest order in $\partial_x$ which is needed is called the \emph{order} of~$u$. Such an element $u$ defines a usual distribution on $X_\RR\times T$ by integration along $T$ with the fixed volume form $\vol_T$.

It is sometimes more convenient to work with currents of maximal degree, which are $C^k$ with respect to $T$. We denote by $\gChT{X}^k$ the corresponding sheaf: a section on $W$ is a continuous $C^\infty(T)$-linear map $\cC^\infty_{X_\RR\times T,c}(W)\to C^k_c(T)$. In particular, when $T=\bS$, $\index{$dbhX$@$\Dbh{X}$}\Dbh{X}^k$ (\resp $\index{$chX$@$\gCh{X}$}\gCh{X}^k$) is a left (\resp right) module over $\cR_{(X,\ov X),\bS}$ defined above.

The regularity of distributions with respect to $T$ is useful in order to get the following:
\begin{itemize}
\item
The restriction to any subvariety $T'$ of an object of $\DbhT{X}^k$ is well defined and is an object of $\Db_{X_\RR\times T'/T'}^k$. In particular, if $T'$ is reduced to a point, we get a ordinary distribution on $X_\RR$.
\item
When $T=\bS$, if $p(\hb)$ is any nonzero polynomial, then
\begin{equation}\label{eq:sanstorsion}
u\in\Dbh{X}^k(W)\text{ and }p(\hb)\cdot u=0\implique u=0.
\end{equation}
\end{itemize}

\begin{Exemples}\label{ex:distran}
\begin{enumerate}
\item\label{ex:distran1}
If $T$ is the Euclidean space $\RR^p$ and $\Delta_T$ is the Laplacian on it, the subsheaf $\ker\Delta_T$ of $T$-harmonic distributions in $\Db_{X_\RR\times T}$ is contained in $\Db_{X_\RR\times T/T}^k$ for any~$k$.

\item\label{ex:distran2}
If $\Omega$ is an open set in $\CC$ with coordinate $\hb$, the sheaf $\index{$dbhXan$@$\Db^\an_{X_\RR\times\Omega}$}\Db^\an_{X_\RR\times\Omega}\defin\ker\ov\partial_\hb$ of distributions on $X_\RR\times\Omega$ which are holomorphic with respect to $\hb$ is contained in $\Db_{X_\RR\times\Omega/\Omega}^k$ for any $k$. We will also denote by $\index{$cxinftyan$@$\cC^{\infty,\an}_{X_\RR\times\Omega}$}\cC^{\infty,\an}_{X_\RR\times\Omega}$ the sheaf of $C^\infty$ functions on $X_\RR\times\Omega$ which are holomorphic with respect to $\hb$ and by $\cCh{\cX}$ the restriction to $\bS$ of $\cC^{\infty,\an}_{\cX}$.

\item\label{ex:distran3}
Let $X=D$ be the open disc of radius $1$ and coordinate $t$. Set $\Omega=\CC\moins(-\NN^*)$. For each $\ell\in\NN$,
\[
U_\ell\defin\frac{\mt^{2s}\module{\log t\ov t}^\ell}{\ell!}
\]
defines a global section of $\Db_{X\times\Omega}^{\an}$ (variable $s$ on $\Omega$). The order of $U_\ell$ is finite on any domain $\Omega\cap \{\reel s\geq -N\}$, $N>0$. If we set $U_\ell=0$ for $\ell<0$, we have
\[
t\partial_t U_\ell=\ov t\partial_{\ov t} U_\ell=sU_\ell+U_{\ell-1}.
\]
\item\label{ex:distran4}
Let $X$ be as above. For $\hb\in\bS$, $t\in X$ and $\ell\geq1$, $\hb/\ov t^\ell-1/\hb t^\ell$ is purely imaginary, so that the function $e^{\hb/\ov t^\ell-1/\hb t^\ell}$ defines a distribution on $X_\RR\times\bS$. This distribution is a section of $\Dbh{X}^0$ (it is even a section of $\Dbh{X}^1$ if $\ell=1$).
\end{enumerate}
\end{Exemples}

Denote by $\cC^{\an,k}_{X\times\bS}$ the subsheaf of $\cC^k_{X\times\bS}$ of functions which are holomorphic with respect to $X$. We then have:

\begin{lemme}[Dolbeault-Grothendieck]\label{lem:DGhb}
$
\ker[\ov\partial_X :\Dbh{X}^k\to\Dbh{X}^{(0,1),k}]=\cC^{\an,k}_{X\times\bS}
$.\qed
\end{lemme}

In the following, we will only use the continuity property with respect to $\bS$, and we will denote by $\Dbh{X}$ the sheaf $\Dbh{X}^0$.

\subsection{Spencer and de~Rham}\label{subsec:spederh}
The de~Rham complex (of $\cO_\cX$) will be usually shifted by $n=\dim X$, with differential $(-1)^nd$. We will denote it by $(\Omega_\cX^{n+\bbullet},(-1)^nd)$.

Given any $k\geq0$, the \emph{contraction} is the morphism
\begin{equation}\label{eq:contraction}
\begin{array}{rcl}
\omega_\cX\otimes_{\cO_\cX}\wedge^k\Theta_\cX&\To{}&\Omega^{n-k}_\cX\\
\omega\otimes\xi&\mto&\varepsilon(n-k)\omega(\xi\wedge\cbbullet).
\end{array}
\end{equation}
The \emph{Spencer complex} $(\index{$spX$@$\Sp_\cX$}\Sp_\cX^{\cbbullet}(\cO_\cX),\delta)$ is the complex $\cR_\cX\otimes_{\cO_\cX}\wedge^{-\bbullet}\Theta_\cX$ (with $\cbbullet\leq0$) of locally free left $\cR_\cX$-modules of finite rank, with differential $\delta$ given by
\begin{align*}
P\otimes\xi_1\wedge\cdots\wedge\xi_k\Mto{\delta}& \sum_{i=1}^{k}(-1)^{i-1}P\xi_i\otimes\xi_1\wedge\cdots\wedge\wh{\xi_i}\wedge\cdots\wedge\xi_k\\ +&
\sum_{i<j}(-1)^{i+j}P\otimes[\xi_i,\xi_j]\wedge\xi_1\wedge\cdots\wedge\wh{\xi_i}\wedge\cdots\wedge\wh{\xi_j}\wedge\cdots\wedge\xi_k.
\end{align*}
It is locally isomorphic to the Koszul complex $K(\cR_\cX,{}\cdot\partiall_{x_1},\dots,{}\cdot\partiall_{x_n})$. It is a resolution of $\cO_\cX$ as a left $\cR_\cX$-module. Under the contraction~\eqref{eq:contraction}
\[
\omega_\cX\ootimes_{\cO_\cX}\wedge^{k}\Theta_\cX\isom\Omega_\cX^{n-k}
\]
the complex $\big(\omega_{\cX}\otimes_{\cR_\cX}\Sp_\cX^{\cbbullet}(\cO_\cX),\delta\big)=\big(\omega_{\cX}\otimes_{\cO_\cX}\wedge^{-\cbbullet}\Theta_{\cX},\delta\big)$ is identified with the complex $(\Omega_\cX^{n+\bbullet},(-1)^nd)$.

Similarly, putting as above $n=\dim X$, the complex $(\Omega_\cX^{n+\bbullet}\otimes_{\cO_\cX}\cR_\cX,\nabla)$, with differential $\nabla$ given by
\[
\omega_{n+\ell}\otimes P\To{\nabla} (-1)^n d\omega_{n+\ell}\otimes P+(-1)^{\ell}\omega_{n+\ell}\wedge\nabla P,
\]
is a resolution of $\omega_\cX$ as a right $\cR_\cX$-module. We will use the notation $\Omega_\cX^{-k}$ for $\wedge^k\Theta_\cX$.

Denote by $(\cE_\cX^{(n+\bbullet,0)}, (-1)^nd')$ the complex $\cC^{\infty,\an}_\cX\otimes_{\cO_\cX}\Omega^{n+\cbbullet}_\cX$ with the differential induced by $(-1)^nd$ (here, we assume $n+\cbbullet\geq0$; recall that the exponent $\an$ means ``holomorphic with respect to $\hb$''). More generally, let $\cE_\cX^{(n+p,q)}=\cE_\cX^{(n+p,0)}\wedge\pi^*\cE_X^{(0,q)}$ (the antiholomorphic part does not produce new poles or zeros along $\hb=0$) and let $d''$ be the usual antiholomorphic differential. For any $p$, the complex $(\cE_\cX^{(n+p,\bbullet)},d'')$ is a resolution of $\Omega_\cX^{n+p}$. We therefore have a complex $(\cE_\cX^{n+\cbbullet}, (-1)^nd)$, which is the single complex associated to the double complex $(\cE_\cX^{(n+\bbullet,\bbullet)}, (-1)^nd',d'')$.

In particular, we have a natural quasi-isomorphism of complexes of right $\cR_\cX$-modules:
\[
(\Omega_\cX^{n+\bbullet}\otimes_{\cO_\cX}\cR_\cX,\nabla)\isom (\cE_\cX^{n+\bbullet}\otimes_{\cO_\cX}\cR_\cX,\nabla)
\]
by sending holomorphic $k$-forms to $(k,0)$-forms. Remark that the terms of these complexes are flat over $\cO_\cX$.

\enlargethispage{3mm}
\subsection{Left and right}\label{subsec:leftright}
At some places in this paper, it is simpler to work with right $\cR_\cX$-modules. The correspondence between both points of view is analogous to that for $\cD_X$-modules. Any left $\cR_\cX$-module $\index{$mcurllr$@$\cM^l$, $\cM^r$}\cM^l$ gives rise to a right one $\cM^r$ by putting (\cf \cite{Castro90} for instance) $\cM^r=\omega_\cX\otimes_{\cO_\cX}\cM^l$ and, for any vector field $\xi$,
\[
(\omega\otimes m)\cdot \xi=\omega\xi\otimes m-\omega\otimes\xi m.
\]
Conversely, put $\cM^l=\cHom_{\cO_\cX}(\omega_\cX,\cM^r)$, which has in a natural way the structure of a left $\cR_\cX$-module. The natural morphisms
\[
\cM^l\to\cHom_{\cO_\cX}(\omega_\cX,\omega_\cX\otimes_{\cO_\cX}\cM^l),\qquad \omega_\cX\otimes_{\cO_\cX}\cHom_{\cO_\cX}(\omega_\cX,\cM^r)\to\cM^r
\]
are isomorphisms of $\cR_\cX$-modules.

If $\cM,\cN$ are two left $\cR_\cX$-modules, we have a natural isomorphism of sheaves of $\CC$-vector spaces
\begin{equation}\label{eq:leflefrig}
\begin{array}{rcl}
\cN^r\otimes_{\cR_\cX}\cM&\isom&\cM^r\otimes_{\cR_\cX}\cN\\
(\omega\otimes n)\otimes m&\mto&(\omega\otimes m)\otimes n,
\end{array}
\end{equation}
which is functorial in $\cM$ and in $\cN$.

\medskip
Notice that $\omega_\cX\otimes_{\cO_\cX}\cR_\cX$ has therefore two structures of right $\cR_\cX$-module, denoted by $\cdotr$ and $\cdott$:
\[
(\omega\otimes P)\cdotr Q=\omega\otimes(PQ)\qqbox{and} (\omega\otimes P)\cdott \xi= \omega\cdot\xi\otimes P-\omega\otimes\xi P,
\]
for $P,Q$ local sections of $\cR_\cX$ and $\xi$ a local section of $\Theta_\cX$.
Recall (\cf\cite[Lemme 2.4.2]{MSaito86}) that there is a unique involution $\iota:\omega_\cX\otimes_{\cO_\cX}\cR_\cX\to\omega_\cX\otimes_{\cO_\cX}\cR_\cX$ which is the identity on $\omega_\cX\otimes1$ and exchanges both structures: it is given by $\omega\otimes P\mto (\omega\otimes1)\cdott P$.

In particular, the isomorphism of right $\cR_\cX$-modules
\begin{align*}
\omega_\cX\otimes_{\cO_\cX}\big(\cR_\cX\otimes_{\cO_\cX}\wedge^{k}\Theta_\cX\big)&\xrightarrow[\ts\sim]{~\ts\iota~}\Omega_\cX^{n-k}\otimes_{\cO_\cX}\cR_\cX\\
\big[\omega\otimes(1\otimes\xi)\big]\cdott P&\mto
\big(\varepsilon(n-k)\omega(\xi\wedge\cbbullet)\big)\otimes P
\end{align*}
where the right structure of the right-hand term is the trivial one and that of the left-hand term is nothing but that induced by the left structure after going from left to right, induces an isomorphism of complexes of right $\cR_\cX$-modules
\begin{equation}\label{eq:lefrigDR}
\iota:\omega_\cX\ootimes_{\cO_\cX}\big(\Sp^\cbbullet_\cX(\cO_\cX),\delta\big)\isom\big(\Omega_\cX^{n+\bbullet}\ootimes_{\cO_\cX}\cR_\cX,\nabla\big).
\end{equation}

\medskip
Similarly, if $\cM$ is any left $\cR_\cX$-module and $\cM^r=\omega_\cX\otimes_{\cO_\cX}\cM$ is the associated right $\cR_\cX$-module, there is an isomorphism
\begin{multline}\label{eq:isoDRrg}
\cM^r\otimes_{\cR_\cX}\big(\Sp^\cbbullet_\cX(\cO_\cX),\delta\big)\simeq \big(\omega_\cX\otimes_{\cO_\cX}\cM\otimes_{\cO_\cX}\wedge^{-\bbullet}\Theta_\cX,\delta\big)\\ \isom
\big(\Omega_\cX^{n+\bbullet}\otimes_{\cO_\cX}\cM,\nabla\big)\simeq\big(\Omega_\cX^{n+\bbullet}\ootimes_{\cO_\cX}\cR_\cX,\nabla\big)\otimes_{\cR_\cX}\cM
\end{multline}
given on $\omega_\cX\otimes_{\cO_\cX}\cM\otimes_{\cO_\cX}\wedge^{k}\Theta_\cX$ by
\[
\omega\otimes m\otimes \xi\mto\varepsilon(n-k)\omega(\xi\wedge\cbbullet)\otimes m.
\]

\medskip
In the same vein, let $\cM$ be a left $\cR_\cX$-module. Then $\cM\otimes_{\cO_\cX}\cR_\cX$ has the structure of a left and of a right $\cR_\cX$-module: $\xi\cdot(m\otimes P)=(\xi m)\otimes P+m\otimes\xi P$, and $(m\otimes P)\cdot \xi=m\otimes(P\xi)$ for any local vector field $\xi$. Similarly, $\cR_\cX\otimes_{\cO_\cX}\cM$ also has such a structure: $\xi\cdot(P\otimes m)=(\xi P)\otimes m$ and $(P\otimes m)\cdot \xi=P\xi\otimes m-P\otimes\xi m$.

Then, there exists a unique isomorphism
\begin{equation}\label{eq:lefrig}
\cM\otimes_{\cO_\cX}\cR_\cX\isom\cR_\cX\otimes_{\cO_\cX}\cM
\end{equation}
of left and right $\cR_\cX$-modules, which induces the identity on $\cM\otimes 1=1\otimes\cM$. If $\xi$ is any local vector field, this isomorphism is given by
\[
m\otimes\xi=(m\otimes 1)\xi\mto(1\otimes m)\xi=\xi\otimes m-1\otimes\xi m.
\]

\subsection{}\label{subsec:defLt}
Let $D$ be an open disc centered at the origin in $\CC$ with complex coordinate $t$. Assume that it has radius $\leq1$ (this will be always the case later on, as we may reduce the size of the disc). The logarithm $\index{$lt$@$\Lt$}\Lt$ is defined as
\[
\Lt=\module{\log\mt^2}=-\log(t\ov t).
\]
It satisfies, for any $k\in\RR$,
\begin{equation}\label{eq:diffLt}
\Lt^kt\;\frac{\partial\lefpar \Lt^{-k}\rigpar}{\partial t} =\Lt^k\ov t\;\frac{\partial\lefpar \Lt^{-k}\rigpar}{\partial \ov t} =k\Lt^{-1}.
\end{equation}

\subsection{}\label{subsec:defstar}
Let $\alpha=\alpha'+i\alpha''$ be a complex number with $\alpha'=\reel(\alpha)$, $\alpha''=\im(\alpha)$. For $\hb\in\CC$, put
\begin{equation}\label{eq:defstar}
\index{$alphastarhb$@$\alpha\star\hb$}\alpha\star\hb=\alpha\hb+i\alpha''(\hb-1)^2/2=\alpha'\hb+i\alpha''(\hb^2+1)/2.
\end{equation}
The following properties are easily verified:
\begin{itemize}
\item
$\alpha\mto \alpha\star\hb$ is $\RR$-linear;
\item
for $\hb\neq0$, the expression
\[
\frac{\alpha\star\hb}{\hb}=\alpha'+i\alpha''(\hb+1/\hb)/2
\]
is ``real'' in the sense of the conjugation defined in \T\ref{subsubsec:conj}, \ie is invariant when we replace $i$ with $-i$ and $\hb$ with $-1/\hb$: indeed, using this notion of conjugation, it is the ``real'' part of $\alpha'+i\hb\alpha''$, namely $\frac12[(\alpha'+i\hb\alpha'')+(\ov{\alpha'+i\hb\alpha''})]$;
\item
if $\alpha''=0$ or if $\hb_o=\pm i$, we have $(\alpha\star\hb_o)/\hb_o=\alpha'$;
\item
We have $\alpha\star\hb_o=0$ if and only if one of the following properties is satisfied:
\begin{enumerate}
\item
$\alpha=0$,
\item
$\alpha\neq0$ is real (\ie $\alpha''=0$ and $\alpha'\neq0$) and $\hb_o=0$,
\item
$\alpha\neq0$ is not real (\ie $\alpha''\neq0$) and $\hb_o=i\big((\alpha'/\alpha'')\pm\sqrt{1+(\alpha'/\alpha'')^2}\big)$ (in particular $\hb_o\in i\RR^*$);
\end{enumerate}
in particular,
\[
\alpha\star\hb_o=0 \text{ and }\alpha\neq0\implique
\begin{cases}
\hb_o\in i\RR,\\
\hb_o=\pm i \iff \alpha\text{ is purely imaginary}.
\end{cases}
\]
\end{itemize}

Let $A\subset\CC$ be a finite subset and put $\Lambda=A+\ZZ$. A complex number $\hb_o\in\Omega_0$ is \emph{singular with respect to $\Lambda$} if there exist $\alpha_1,\alpha_2\in\Lambda\cup\ZZ$ such that $\alpha_1\neq\alpha_2$ and $(\alpha_1-\alpha_2)\star\hb_o=0$. Such a $\hb_o$ is purely imaginary. The set of nonzero $\Lambda$-singular complex numbers is discrete in $i\RR^*$, and $0$ is its only possible limit point in $i\RR$. It is reduced to $\{0\}$ if $\Lambda\subset\RR$. We denote it by $\index{$lambda$@$\Lambda$, $\Sing(\Lambda)$}\Sing(\Lambda)$.

\medskip
For $\hb_o\in\Omega_0$, denote $\imhb_o=\im\hb_o$ and set $\ell_{\hb_o}(\alpha)=\alpha'-\imhb_o\alpha''=\reel(\alpha'+i\hb_o\alpha'')$ (where $\reel$ is taken in the usual sense). The following lemma will be useful:

\begin{lemme}\label{lem:imhb}
Fix $\hb_o\in\Omega_0$. If $\ell_{\hb_o}(\alpha)=0$, then $\alpha\star\hb_o=0\implique\alpha=0$.
\end{lemme}

\begin{proof}
Assume first that $\hb_o=0$. The hypothesis is that $\alpha'=0$. We then have $\alpha\star\hb_o=i\alpha''/2$.

Assume now that $\hb_o\neq0$. The hypothesis is that $\alpha'=\imhb_o\alpha''$ and we have
\[
\frac{\alpha\star\hb_o}{\hb_o}=\alpha''(\imhb_o+i(\hb_o+1/\hb_o)/2).
\]
If $\alpha''\neq0$, this could vanish only if $\hb_o+1/\hb_o$ is purely imaginary, hence only if $\hb_o$ is so, \ie $0\neq\hb_o=i\imhb_o$; but we would have $0=\imhb_o+i(\hb_o+1/\hb_o)/2=(\imhb_o+1/\imhb_o)/2$ with $\imhb_o\in\RR^*$, impossible.
\end{proof}
\end{chapter}

\makeatletter
\def\cl@subsection{}
\makeatother
\numberwithin{equation}{section}
\def\thesubsection{\thesection.\alph{subsection}}
\chapter{Coherent and holonomic $\cR_\cX$-modules}\label{chap:RX}

\section{Coherent and good $\cR_\cX$-modules}\label{subsec:RX}
\subsection{}
The ring $R_X$ is equipped with a natural increasing filtration (locally given by the total degree in $\partiall_{x_i}$) and the associated graded object is naturally identified with the sheaf
\[
\cO_X[\hb][TX]\defin p_*\cO_{\PP(T^*X\oplus \mathbf{1})}(*\infty)[\hb],
\]
where $p$ is as in \T\ref{subsec:debutRX}. By usual arguments, it follows that $R_X$ is a coherent sheaf of rings on $X$.

Analogous results hold for $\cR_\cX$, which is a coherent sheaf of rings on $\cX$, by replacing $\cO_X[\hb]$ with $\cO_\cX$.

We may identify the restriction $\cR_{\cX^\Cir}$ with the sheaf of relative differential operators $\cD_{\cX^\Cir/\CC^*}$ by
\[
\partiall_{x_i}\in\cR_{\cX^\Cir}\mto\partial_{x_i}=\hbm\partiall_{x_i}\in\cD_{\cX^\Cir/\CC^*}.
\]
It follows in particular that, for any $\hb_o\neq0$ and any coherent $\cR_\cX$-module $\cM$, the cohomology modules of the complex of $\cR_\cX$-modules $\bL i^*_{\hb_o}\cM\defin\{\cM\stackrel{\hb-\hb_o}{\longrightarrow}\cM\}$ are coherent $\cD_X$-modules. Following \cite{Simpson90}, we put
$$\index{$xdoldr$@$\Xi_{\Dol}$, $\Xi_{\DR}$}\Xi_{\Dol}(\cM)\defin\cM/\hb\cM \quad \text{and} \quad\Xi_{\DR}(\cM)\defin\cM/(\hb-1)\cM.
$$
Then $\Xi_{\Dol}(\cM)$ is a coherent $\cO_X[TX]$-module and $\Xi_{\DR}(\cM)$ is a coherent $\cD_X$-module.

Notice that the datum of a left $\cD_X$-module $M$ is equivalent to the datum of a $\cO_X$-module $M$ equipped with a \emph{flat} connection $\nabla:M\to\Omega^1_X\otimes_{\cO_X}M$. Similarly, the datum of a (left) $\cO_X[TX]$-module $M$ is equivalent to that of a $\cO_X$-module $M$ equipped with a \emph{$\cO_X$-linear morphism} $\theta:M\to\Omega^1_X\otimes_{\cO_X}M$ satisfying the \emph{Higgs condition} $\theta\wedge\theta=0$.

\subsection{De~Rham and Dolbeault complexes}\label{subsec:derhamdolbeault}
Recall that the de~Rham complex $\index{$dr$@$\DR(M)$}\DR(M)$ of a left $\cD_X$-module $M$ is the complex $\Omega_X^{n+\bbullet}\otimes_{\cO_X}M$ with differential
\[
\nabla(\omega_{n+k}\otimes m)=(-1)^nd\omega_{n+k}\otimes m+(-1)^k\omega_{n+k}\wedge\nabla m.
\]
Similarly, the Dolbeault complex $\Dol(M)$ of a $\cO_X[TX]$-module $M$ is the complex $\Omega_X^{n+\bbullet}\otimes_{\cO_X}M$ with differential $(-1)^\cbbullet{}\wedge\theta$ (\ie $(-1)^n\theta\wedge{}$ when putting the forms on the right).

Let now $\cM$ be a left $\cR_\cX$-module. Its de~Rham complex $\index{$drcurl$@$\cDR(\cM)$}\cDR(\cM)$ is $\Omega_\cX^{n+\bbullet}\otimes_{\cO_\cX}\cM$ with differential
\[
\nabla(\omega_{n+k}\otimes m)=(-1)^nd\omega_{n+k}\otimes m+(-1)^k\omega_{n+k}\wedge\nabla m.
\]
As $\DR\cM$ is a complex of $\cO_{\Omega_0}$-modules, we have $\bL i^*_{\hb_o}(\cDR\cM)=\DR(\bL i^*_{\hb_o}\cM)$ if $\hb_o\neq0$ and $\bL i^*_0(\cDR\cM)=\Dol(\bL i^*_0\cM)$.

\medskip
The de~Rham complex $\cDR(\cM)$ of a left $\cR_\cX$-module is also equal to the complex $(\Omega_\cX^{n+\bbullet}\otimes_{\cO_\cX}\cR_\cX,\nabla)\otimes_{\cR_\cX}\cM$.

We define the de~Rham complex of a right $\cR_\cX$-module $\cN$ as $\cN\otimes_{\cR_\cX}\Sp_{\cX}^\cbbullet(\cO_\cX)$.

Using \eqref{eq:isoDRrg}, we have a functorial isomorphism $\cDR(\cM^r)\isom\cDR(\cM)$ for any left $\cR_\cX$-module $\cM$.

\subsection{}
The sheaf $\cR_\cX$ comes equipped with an increasing filtration by locally free $\cO_\cX$-submodules, indexed by the order in $\partiall_{x_1},\dots,\partiall_{x_n}$. We may therefore define, as usual, the notion of a \emph{good filtration} on a $\cR_\cX$-module.
Following \cite{Sch-Sch94}, we say that a $\cR_\cX$-module $\cM$ is \emph{good} if, for any compact subset $\cK\subset \cX$, there exists in a neighbourhood of $\cK$, a finite filtration of $\cM$ by $\cR_\cX$-modules such that all successive quotients have a \emph{good} $\cR_\cX$-filtration. This implies that $\cM$ is coherent.
\section{The involutivity theorem}\label{subsec:Bernstein}
Let $\cM$ be a coherent $\cR_\cX$-module. The support in $T^*X\times\Omega_0$ of the graded module associated to any local good filtration of $\cM$ does not depend on the choice of such a good filtration and is defined globally (see \eg \cite[Prop\ptbl A:III.3.21]{Bjork93}): this is the characteristic variety $\index{$charMcurl$@$\Char\cM$}\Char\cM$ of $\cM$.

For any $\hb_o\neq0$, denote by $\Char_{\hb_o}(\cM)$ the union of the characteristic varieties of the cohomology $\cD_X$-modules of $\bL i^*_{\hb_o}\cM$. There is a natural inclusion $\Char_{\hb_o}(\cM)\subset\Char\cM\cap(T^*X\times \{\hb_o\})$.

Let $\index{$sigmaMcurl$@$\Sigma(\cM)$}\Sigma(\cM)\subset T^*X$ be the support of $\Xi_{\Dol}(\cM)$. This is a kind of ``characteristic variety'', but may not be homogeneous with respect to the usual $\CC^*$-action on $T^*X$.

It is possible to associate a multiplicity to each irreducible component of $\Char\cM$, $\Char_{\hb_o}(\cM)$ or $\Sigma(\cM)$ to get a characteristic cycle.

The definition of $\Char\cM$, $\Char_{\hb_o}(\cM)$ or $\Sigma(\cM)$ extends to complexes: just take the union of characteristic varieties of the cohomology sheaves.

The \emph{support} $\index{$suppMcurl$@$\supp\cM$}\supp\cM\subset X$ is by definition the closure of the projection of $\Char\cM$ in $X$. This is the smallest closed subset $Z$ of $X$ such that $\cM$ vanishes identically on $(X\moins Z)\times\Omega_0$.

\begin{definition}
A $\cR_\cX$-module $\cM$ is said to be \emph{strict} if it has no $\cO_{\Omega_0}$-torsion. A complex $\cM^{\cbbullet}$ of $\cR_\cX$-modules is said to be \emph{strict} if each of its cohomology modules is so. A morphism $\varphi:\cM\to\cN$ is \emph{strict} if the corresponding complex is so, \ie if $\ker\varphi$ and $\coker\varphi$ are strict.
\end{definition}

Notice that, for a strict $\cR_\cX$-module $\cM$, the restriction $\bL^*_{\hb_o}\cM$ reduces to the degree~$0$ term $\cM/(\hb-\hb_o)\cM$.

\begin{Lemme}\label{lem:Wstrict}
\begin{enumerate}
\item\label{lem:Wstrict1}
Let $\cM$ be a $\cR_\cX$-module equipped with a finite increasing filtration $W_\bbullet\cM$ by $\cR_\cX$-submodules. If each $\gr_k^W\cM$ is strict, then $\cM$ is strict.
\item\label{lem:Wstrict2}
Let $\varphi:\cM\to\cN$ be a morphism of $\cR_\cX$-modules. Assume that $\cM,\cN$ have a finite filtration $W$ by $\cR_\cX$-submodules and that $\varphi$ is strictly compatible with $W$, \ie satisfies $\varphi(W_k\cM)=W_k\cN\cap\varphi(\cM)$. If $\gr^W_k\varphi$ is strict for all $k$, then $\varphi$ is strict.
\end{enumerate}
\end{Lemme}

\begin{proof}
The first point is clear. Let us prove \eqref{lem:Wstrict2}. By strict compatibility, the sequence
\[
0\to\gr_k^W\ker\varphi\to\gr_k^W\cM\To{\gr_k^W\varphi}\gr_k^W\cN\to\gr_k^W\coker\varphi\to0
\]
is exact, putting on $\ker\varphi$ and $\coker\varphi$ the induced filtration. By strictness of $\gr_k^W\varphi$, and applying \eqref{lem:Wstrict1} to $\ker\varphi$ and $\coker\varphi$, one gets \eqref{lem:Wstrict2}.
\end{proof}

\begin{theoreme}\label{th:Gabber}
Let $\cM$ be a strict coherent $\cR_\cX$-module. Then $\Sigma(\cM)$ and $\Char_{\hb_o}(\cM)$ ($\hb_o\in\CC^*$) are involutive in $T^*X$, and $\Char\cM$ is involutive in $T^*X\times\Omega_0$ (with respect to the Poisson bracket $\hb\{\,,\,\}$).
\end{theoreme}

\begin{proof}
This is well-known for the characteristic varieties $\Char_{\hb_o}(\cM)$ and $\Char\cM$ (\cite{Gabber81}, see also \cite[A:III.3.25]{Bjork93}). The proof of Gabber's involutivity theorem also applies to $\Sigma(\cM)$ because $\cM$ is strict (indeed, $\cR_\cX/\hb^2\cR_\cX$ is a Gabber ring, in the sense of \cite[A:III.3]{Bjork93}).
\end{proof}

The restriction of $\Char\cM$ over $\{\hb=0\}$ is not controlled by the involutivity theorem. However, the restriction to $\{\hb=0\}$ of components of $\Char\cM$ for which the fibre at some $\hb_o\neq0$ is Lagrangian is a union of irreducible conical Lagrangian closed analytic subsets of $T^*X$.

\begin{definition}\label{def:holo}
A $\cR_\cX$-module $\cM$ is said to be \emph{holonomic} if it is \emph{good} and there exists a conical Lagrangian variety $\Lambda\subset T^*X$ such that the characteristic variety $\Char\cM$ is contained in $\Lambda\times\Omega_0$. A complex $\cM^{\cbbullet}$ of $\cR_\cX$-modules is said to be \emph{holonomic} if each of its cohomology module is so.
\end{definition}

If $\cM$ is holonomic, then any irreducible component of $\Char\cM$ is equal to $T^*_ZX\times\Omega_0$ or to $T^*_ZX\times\{\hb_{o}\}$ for some closed irreducible analytic subset $Z$ of $X$ and some $\hb_{o}\in \Omega_0$. In particular, the support of $\cM$ is \emph{equal} to the projection of $\Char\cM$ in $X$.

In an exact sequence $0\to\cM'\to\cM\to\cM''\to0$ of good $\cR_\cX$-modules, $\cM$ is holonomic if and only if $\cM'$ and $\cM''$ are so: indeed, we have $\Char(\cM)=\Char(\cM')\cup\Char(\cM'')$.

\begin{proposition}[Restriction to $\hb=\hb_o$]\label{prop:lagrange}
If $\cM$ is holonomic, then for any \hbox{$\hb_o\!\neq\!0$}, the cohomology modules of $\bL^*_{\hb_o}\cM$ are holonomic $\cD_X$-modules. Moreover, if $\cM$ is strict, $\Sigma(\cM)$ is Lagrangian.
\end{proposition}

\begin{proof}
Let $F_\bbullet\cM$ be a good filtration of $\cM$ locally near a point of $X\times\{\hb_o\}$. It induces a good filtration $(F_\bbullet\cM)\cap(\hb-\hb_o)\cM$ on the coherent $\cR_\cX$-submodule $(\hb-\hb_o)\cM$ as well as on the coherent quotient $\cM/(\hb-\hb_o)\cM$, the graded module of which is a quotient of $\gr^F\!\cM/(\hb-\hb_o)\gr^F\!\cM$. Similarly, $\gr^F\big[\ker(\hb-\hb_o)\big]$ is contained in the kernel of $\hb-\hb_o$ acting on $\gr^F\!\cM$. This implies the first point. For $\hb_o=0$, we have
\begin{align*}
\dim X&\leq \dim\Sigma(\cM)\quad\text{($\Sigma(\cM)$ is involutive)}\\
&=\dim\supp\gr^F(\cM/\hb\cM)\quad\text{(conservation of the dimension by grading)}\\
&\leq \dim\supp\gr^F\!\cM/\hb\,\gr^F\!\cM\\
&\leq \dim X \quad\text{($\cM$ is holonomic)}.\qedhere
\end{align*}
\end{proof}

The variety $\Sigma(\cM)$ is well-behaved in exact sequences only for strict objects in general. One has for instance

\begin{proposition}
Let $0\to\cM'\to\cM\to\cM''\to 0$ be an exact sequence of strict holonomic $\cR_\cX$-modules. Then $\Sigma(\cM)=\Sigma(\cM')\cup\Sigma(\cM'')$ and the corresponding Lagrangian cycles behave in an additive way.\qed
\end{proposition}

\begin{remarque}
Analogous results hold for $R_X$-modules. We leave them to the reader.
\end{remarque}

\begin{proposition}\label{prop:CK}
Let $\cM$ be a strict coherent $\cR_\cX$-module, the characteristic variety of which is contained in the zero section $T^*_XX\times\Omega_0$. Then,
\begin{enumerate}
\item\label{prop:CK1}
$\cM$ is $\cO_\cX$-coherent,
\item\label{prop:CK2}
$\cM^\Cir$ is locally (on $\cX^\Cir$) isomorphic to $\cO_{\cX^\Cir}^{d}$ equipped with its natural structure of left $\cR_{\cX^\Cir}$-module, for some integer $d$,
\item\label{prop:CK3}
there exists a nowhere dense closed analytic subset $Z\subset X$ such that $\cM$ is $\cO_\cX$-locally free on $\cX\moins\cZ$.
\end{enumerate}
\end{proposition}

\begin{proof}
The first point is clear. As $\cM^\Cir$ is $\cO_{\cX^\Cir}$-coherent (hence good as a $\cR_{\cX^\Cir}$-module), the second point follows from \cite[Theorem 2.23 (iii)]{Deligne70} and strictness. The third point also follows from strictness.
\end{proof}

\begin{remarque}\label{rem:Stein}
Under the assumption of Proposition \ref{prop:CK}, there exists locally on~$X$ a vector bundle $E$ such that $\cM^{\Cir}=\pi^{\Cir *}E$ as a $\cO_{\cX^\Cir}$-module (indeed, if $U$ is any contractible Stein open set of $X$, any vector bundle on $U\times \CC^*$ is topologically trivial, hence analytically trivial, by Grauert's theorem).
\end{remarque}

\section{Examples}

\subsection{The twistor deformation of an irreducible flat connection to a Higgs bundle \cite{Simpson92,Simpson97}}\label{subsub:a}
Let $X$ be a projective manifold and let $(V,\nabla)$ be a flat holomorphic vector bundle on $X$. The construction explained in \loccit gives an example of (and in fact is a model for) a strict holonomic $\cR_\cX$-module. Let us recall the main definitions.

Let $\index{$dV$@$D_V$}D_V=D'_V+D''_V$ be the flat connection on $H\defin\cC^\infty_X\otimes_{\cO_X} V$, so that $(V,\nabla)=(\ker D''_V,D'_V)$, and let $h$ be a metric on $\index{$hdv$@$(H,D_V)$}(H,D_V)$. There exist connections denoted $\index{$dE$@$D_E$}D'_E$ (of type $(1,0)$) and $D''_E$ (of type $(0,1)$), and a $(1,0)$-form $\theta'_E$ with values in $\Endom(H)$ such that, denoting by $\theta''_E$ the adjoint of $\index{$thE$@$\theta_E$}\theta'_E$ with respect to $h$, we have, for any local sections $u,v$ of $H$,
\begin{align*}
d'h(u,v)&=h(D'_E u,v)+h(u,D''_Ev),\\
d'' h(u,v)&=h(D''_E u,v)+h(u,D'_Ev),\\
h(\theta'_Eu,v)&=h(u,\theta''_Ev),\\
D'_V=D'_E+\theta'_E,&\quad D''_V=D''_E+\theta''_E.
\end{align*}
These objects are uniquely defined by the previous requirements. Notice that, by applying $d'$ or $d''$ to each of the first three lines above, we see that $D_E^{\prime\prime2}$ is adjoint to $D_E^{\prime2}$, $D''_E(\theta'_E)$ is adjoint to $D'_E(\theta''_E)$ and $D'_ED''_E+D''_ED'_E$ is selfadjoint with respect to~$h$.

The triple $(H,D_V,h)$ (or $(V,\nabla,h)$, or simply $h$, if $(V,\nabla)$ is fixed) is said to be \emph{harmonic} if the operator $D''_E+\theta'_E$ has square $0$. By looking at types, this is equivalent to
\[
D_E^{\prime\prime2}=0,\quad D''_E(\theta'_E)=0,\quad \theta'_E\wedge\theta'_E=0.
\]
By adjunction, this implies
\[
D_E^{\prime2}=0,\quad D_E'(\theta''_E)=0,\quad \theta''_E\wedge\theta''_E=0.
\]
Moreover, the flatness of $D_V$ implies then
\[
D'_E(\theta'_E)=0,\quad D''_E(\theta''_E)=0,\quad D'_ED''_E+D''_ED'_E=-(\theta'_E\theta''_E+\theta''_E\theta'_E).
\]

Let $E=\ker D''_E:H\to H$. This is a holomorphic vector bundle equipped with a holomorphic $\endom(E)$-valued $1$-form $\theta'_E$ satisfying $\theta'_E\wedge\theta'_E=0$. It is called a \emph{Higgs bundle} and $\theta'_E$ is its \emph{associated Higgs field}.

Remark that $\theta'_E:E\to E\otimes_\cO\Omega^1_X$ can be viewed as a holomorphic map $\theta'_E:\Theta_X\to\endom(E)$ satisfying $[\theta'_E(\xi),\theta'_E(\eta)]=0$ for any vector fields $\xi,\eta$, defining thus the structure of a $\cO_X[TX]$-module on $E$. Its support in $T^*X$ is therefore a finite ramified covering of $X$.

The previous relations also imply that, if $\hb_o$ is any complex number, the operator $D''_E+\hb_o\theta''_E$ is a complex structure on $H$. Moreover, if $\hb_o\neq0$, the holomorphic bundle $V_{\hb_o}=\ker (D''_E+\hb_o\theta''_E)$ is equipped with a flat holomorphic connection $\nabla_{\hb_o}=D'_E+\hbm_o\theta'_E$. For $\hb_o=1$ we recover $(V,\nabla)$.

Consider the $\cC^{\infty,\an}_\cX$-module $\index{$hcurl$@$\cH$}\cH=\cC^{\infty,\an}_\cX\otimes_{\pi^{-1}\cC^\infty_X}\pi^{-1}H$, equipped with a $d''$ operator
\begin{equation}\label{eq:dbar}
\index{$dHcurl$@$D_\cH$}D''_{\cH}=D''_E+\hb\theta''_E.
\end{equation}
This defines a holomorphic subbundle $\cH'$ (that is, a locally free $\cO_\cX$-submodule such that $\cC^{\infty,\an}_\cX\otimes_{\cO_\cX}\cH'=\cH$), which is thus strict. Moreover, it has the natural structure of a good $\cR_\cX$-module, using the flat connection
\begin{equation}\label{eq:dH}
D'_\cH=D'_E+\hbm\theta'_E.
\end{equation}
One has $\Xi_{\Dol}(\cH')=(E,\theta'_E)$ and $\Xi_{\DR}(\cH')=(V,D'_V)$. Clearly, $\Char\cH'$ is equal to $T^*_XX\times\Omega_0$ (take the trivial filtration). Then $\cH'$ is a strict holonomic $\cR_\cX$-module.

\begin{remarque*}
The support $\Sigma\subset T^*X$ of a holomorphic Higgs bundle $(E,\theta'_E)$ (viewed as a $\cO_X[TX]$-module) coming from a harmonic flat bundle $(H,D_V,h)$ is Lagrangian in $T^*X$, after Proposition \ref{prop:lagrange}. More generally, any Higgs bundle $(E,\theta'_E)$ on a projective manifold $X$ satisfies this property, without referring to the existence of a Hermite-Einstein metric (\ie an associated flat harmonic bundle): indeed, restrict the standard holomorphic Liouville $1$-form on $T^*X$ to $\Sigma$ and then lift it to a resolution $\wt\Sigma$ of the singularities of $\Sigma$, which is a projective manifold, as it is a finite ramified covering of $X$; by standard Hodge theory, the lifted form is closed, hence so is its restriction to the regular part $\Sigma^o$ of $\Sigma$; the restriction to $\Sigma^o$ of the canonical $2$-form on $T^*X$ is thus identically $0$ on $\Sigma^o$.
\end{remarque*}

\subsection{Filtered $\cD_X$-modules}\label{subsub:b}
Let $(M,F)$ be a filtered holonomic $\cD_X$-module and $R_FM$ the associated graded Rees module. Put $\cM=\cO_\cX\otimes_{\cO_X[\hb]}R_FM$. By construction, $\cM$ has no $\hb$-torsion and thus is strict holonomic, because $\Char\cM=\Char(M)\times\Omega_0$.

\subsection{Variations of complex Hodge structures \cite{Simpson92}}\label{subsub:c}
Let $H=\oplus_{p\in\ZZ}H^{p,w-p}$ be a $C^\infty$ vector bundle on $X$, where $w\in\ZZ$ is fixed, equipped with a flat connection $D_V=D'_V+D''_V$ and a flat nondegenerate Hermitian bilinear form $k$ such that the direct sum decomposition of $H$ is $k$-orthogonal, $(-1)^pi^{-w}k$ is a metric on $H^{p,w-p}$, \ie $(-1)^pi^{-w}k$ is positive definite on the fibres of $H^{p,w-p}$ for each $p$, and
\begin{align*}
D'_V(H^{p,w-p})&\subset\lefpar H^{p,w-p}\oplus H^{p-1,w-p+1}\rigpar \otimes_{\cO_X} \Omega_{X}^{1}\\
D''_V(H^{p,w-p})&\subset\lefpar H^{p,w-p}\oplus H^{p+1,w-p-1}\rigpar \otimes_{\cO_{\ov X}}\Omega_{\ov X}^{1}
\end{align*}
where $\ov X$ denotes the complex conjugate manifold.

Denote by $D'_V= D'_E+\theta'_E$ and $D''_V=D''_E+\theta''_E$ the corresponding decomposition. Then the metric $h$ defined as $(-1)^pi^{-w}k$ on $H^{p,w-p}$ and such that the direct sum decomposition of $H$ is $h$-orthogonal is a harmonic metric and the objects $ D'_E$, $D''_E$, $\theta'_E$ and $\theta''_E$ are the one associated with $(h,D_V)$ as in example \ref{subsub:a}.

Put $F^p=\oplus_{q\geq p}H^{q,w-q}$. This bundle is stable under $D''_V$. Let $F^pV=F^p\cap V$ be the corresponding holomorphic bundle. Consider on the Rees module $\oplus F^p\hb^{-p}\subset H[\hb,\hbm]=H\otimes_\CC\CC[\hb,\hbm]$ the holomorphic structure induced by $D''_V$. The holomorphic bundle corresponding to it is the Rees module $\oplus_pF^{p}V\hb^{-p}$ attached as in example \ref{subsub:b} to the filtered $\cD_X$-module $(V,F^{\cbbullet}V)$ (put $F_\bbullet=F^{-\bbullet}$ to get an increasing filtration).

On the other hand, consider on $\CC[\hb]\otimes_\CC H$ the holomorphic structure given by $D''_E+\hb\theta''_E$, as defined in example \ref{subsub:a}.

The natural $\CC[\hb]$-linear map
\begin{align*}
\CC[\hb]\otimes_\CC H&\To{\varphi}\CC[\hb,\hbm]\otimes_\CC H\\
1\otimes (\oplus u_p)&\mto\sum u_p\hb^{-p}
\end{align*}
is an isomorphism onto $\oplus F^p\hb^{-p}$ and the following diagram commutes
$$
\xymatrix{
\CC[\hb]\otimes_\CC H\ar[r]^-{\varphi}\ar[d]_-{D''_E+\hb\theta''_E}&\oplus F^p\hb^{-p}\ar[d]^{D''_V}\\
\CC[\hb]\otimes_\CC H\ar[r]^-{\varphi}&\oplus F^p\hb^{-p}
}
$$
showing that, in case of complex variation of Hodge structures, the construction of examples \ref{subsub:a} and \ref{subsub:b} are isomorphic.

\section{Direct and inverse images of $\cR_\cX$-modules}\label{sec:iminvdir}

\subsection{Direct images of $\cR_\cX$-modules}\label{subsec:imdir}
Let $f:X\to Y$ be a holomorphic map between analytic manifolds and denote also by $f:\cX\to\cY$ the map trivially induced. As in the theory of $\cD_X$-modules, one defines the sheaves $\cR_{\cX\to\cY}$ and $\cR_{\cY\leftarrow \cX}$ with their bimodule structure: the sheaf $\index{$rxy$@$\cR_{\cX\to\cY}$}\cR_{\cX\to\cY}=\cO_\cX\otimes_{f^{-1}\cO_\cY} f^{-1}\cR_\cY$ is a left-right $(\cR_\cX,f^{-1}\cR_\cY)$-bimodule when using the natural right $f^{-1}\cR_\cY$-module structure and the usual twisted left $\cR_\cX$-module structure: for any section $\xi$ of the sheaf $\Theta_\cX$ of vector fields on $\cX$ tangent to the fibres of $\pi$ and vanishing at $\hb=0$ (\cf \T\ref{num:RX}), $Tf(\xi)$ is a local section of $\cO_\cX\otimes_{f^{-1}\cO_\cY}\Theta_\cY$, hence acts by left multiplication on $\cR_{\cX\to\cY}$; put $\xi\cdot(\varphi\otimes P)=\xi(\varphi)\otimes P+Tf(\xi)(\varphi\otimes P)$.

The sheaf $\cR_{\cY\leftarrow \cX}$ is obtained by using the usual left-right transformation (see \eg \cite{Castro90} for details). Recall that, if $f$ is an embedding, the sheaves $\cR_{\cX\to\cY}$ and $\cR_{\cY\leftarrow \cX}$ are locally free over $\cR_\cX$.

Denote by $\index{$spXY$@$\Sp_{\cX\to\cY}$}\Sp_{\cX\to\cY}^{\cbbullet}(\cO_\cX)$ the complex $\Sp_\cX^{\cbbullet}(\cO_\cX)\otimes_{f^{-1}\cO_\cY}f^{-1}\cR_\cY$, where the left $\cR_\cX$ structure for each term is twisted as above (recall that the Spencer complex $\Sp_\cX^{\cbbullet}(\cO_\cX)$ was defined in \T\ref{subsec:spederh}). Then $\Sp_{\cX\to\cY}^{\cbbullet}(\cO_\cX)$ is a resolution of $\cR_{\cX\to\cY}$ as a bimodule, by locally free left $\cR_\cX$-modules.

\begin{exemples}\label{ex:Sp}
For $f=\id:X\to X$, the relative Spencer complex $\Sp_{\cX\to\cX}^{\cbbullet}(\cO_\cX)$ which is nothing but $\Sp_\cX^{\cbbullet}(\cO_\cX) \otimes_{\cO_\cX}\cR_\cX$ is a resolution of $\cR_{\cX\to\cX}=\cR_\cX$ as a left and right $\cR_\cX$-module. For $f:X\to\mathrm{pt}$, the complex $\Sp_{\cX\to\mathrm{pt}}^{\cbbullet}(\cO_\cX)=\Sp_{\cX}^{\cbbullet}(\cO_\cX)$ is a resolution of $\cR_{\cX\to\mathrm{pt}}=\cO_X$. If $X=Y\times Z$ and $f$ is the projection, the complex $\cR_\cX\otimes_{\cO_\cX}\wedge^{-\bbullet}\Theta_{\cX/\cY}$ is also a resolution of $\cR_{\cX\to\cY}$ as a bimodule. We moreover have a canonical quasi-isomorphism as bimodules
\begin{align*}
\Sp_{\cX\to\cY}^{\cbbullet}(\cO_\cX)&=\big(\cR_\cX\otimes_{\cO_\cX}\wedge^{-\bbullet}\Theta_{\cX/\cY}\big)\ootimes_{f^{-1}\cO_\cY}f^{-1}\big(\wedge^{-\bbullet}\Theta_\cY\otimes_{\cO_\cY}\cR_\cY\big)\\
&=\big(\cR_\cX\otimes_{\cO_\cX}\wedge^{-\bbullet}\Theta_{\cX/\cY}\big)\ootimes_{f^{-1}\cR_\cY}f^{-1}\big(\Sp_\cY^{\cbbullet}(\cO_\cY)\otimes_{\cO_\cY}\cR_\cY\big)\\
&\isom\big(\cR_\cX\otimes_{\cO_\cX}\wedge^{-\bbullet}\Theta_{\cX/\cY}\big)\ootimes_{f^{-1}\cR_\cY}f^{-1}\cR_{\cY\to\cY}\\
&=\cR_\cX\otimes_{\cO_\cX}\wedge^{-\bbullet}\Theta_{\cX/\cY}.
\end{align*}
\end{exemples}

Recall that $\God^{\cbbullet}$ denotes the canonical Godement resolution (\cf \T\ref{subsec:signes}\eqref{subsec:signesc}). Remark that, if $\cL$ and $\cF$ are $\cO_\cX$-modules and if $\cF$ is locally free, then the natural inclusion of complexes $\God^{\cbbullet}(\cL)\otimes_{\cO_\cX}\cF\hto\God^{\cbbullet}(\cL\otimes_{\cO_\cX}\cF)$ is a quasi-isomorphism.

\begin{definition}\label{def:fdag}
The direct image with proper support $f_\dag$ is the functor from $\Mod^r(\cR_\cX)$ to $D^+(\Mod^r(\cR_\cY))$ defined by (we take the single complex associated to the double complex)
\[
\index{$fdag$@$f_\dag\cM$}f_\dag\cM= f_!\God^{\cbbullet}\Big(\cM\ootimes_{\cR_\cX} \Sp_{\cX\to\cY}^{\cbbullet}(\cO_\cX)\Big).
\]
It is a realization of $\bR f_!\big(\cM\otimes_{\cR_\cX}^{\bL}\cR_{\cX\to\cY}\big)$.
\end{definition}

\begin{Remarques}\label{rem:imdir}
\begin{enumerate}
\item\label{rem:imdir1}
$f_\dag$ can be extended as a functor from $D^+(\Mod^r(\cR_\cX))$ to $D^+(\Mod^r(\cR_\cY))$.
\item\label{rem:imdir2}
Let $f:X\to Y$ be as above and let $Z$ be another manifold. Put $F=f\times\id_Z:X\times Z\to Y\times Z$. Denote by $f_\dag$ the direct image defined on $\Mod^r(\cR_{\cX\times_{\Omega_0}\cZ})$ using the relative Spencer complex (\ie defined with $\Theta_{\cX\times_{\Omega_0} \cZ/\cZ}$) and viewing the $\cR_\cZ$ action as an extra structure on the terms of $f_\dag\cM$, commuting with the differentials of this complex. Then there is a canonical and functorial isomorphism of functors $f_\dag\isom F_\dag$. We will not distinguish between both functors.
\end{enumerate}
\end{Remarques}

\begin{Proposition}\label{prop:imdir}
\begin{enumerate}
\item Let $f:X\to Y$ and $g:Y\to Z$ be two maps. There is a functorial canonical isomorphism of functors $(g\circ f)_\dag=g_\dag f_\dag$.
\item If $f$ is an embedding, then $f_\dag\cM=f_*(\cM\otimes_{\cR_\cX}\cR_{\cX\to\cY})$.
\item
If $f:X=Y\times Z\to Y$ is the projection, we have
\[
f_\dag\cM=f_!\God^{\cbbullet}\big(\cM\otimes_{\cO_\cX}\wedge^{-\bbullet}\Theta_{\cX/\cY}\big),
\]
and this complex is canonically and functorially isomorphic to the relative Dolbeault complex $f_!\big(\cM\otimes_{\cO_\cX}\cE_{\cX/\cY}^{-\bbullet,\bbullet}\big)$.
\end{enumerate}
\end{Proposition}

\begin{proof}
We have a natural morphism
\begin{multline*}
\Sp_{\cX\to\cY}^{\cbbullet}(\cO_\cX)\ootimes_{f^{-1}\cR_\cY}\Sp_{\cY\to\cZ}^{\cbbullet}(\cO_\cY) \to\Sp_{\cX\to\cY}^{\cbbullet}(\cO_\cX)\ootimes_{f^{-1}\cR_\cY}f^{-1}\cR_{\cY\to\cZ}\\
\isom \Sp_{\cX\to\cZ}^{\cbbullet}(\cO_\cX).
\end{multline*}
Both complexes are a resolution of $\cR_{\cX\to\cZ}$ by locally free $\cR_\cX$-modules: this is clear for the right-hand term; for the left-hand term, remark that it is naturally quasi-isomorphic to
\begin{align*}
\cR_{\cX\to\cY}\ootimes_{f^{-1}\cR_\cY}\Sp_{\cY\to\cZ}^{\cbbullet}&(\cO_\cY) =\cO_\cX\ootimes_{f^{-1}\cO_\cY}\Sp_{\cY\to\cZ}^{\cbbullet}(\cO_\cY) \\
& = \cO_\cX\ootimes^{\bL}_{f^{-1}\cO_\cY}f^{-1}\cR_{\cY\to\cZ}\\
& = \cO_\cX\ootimes_{g^{-1}f^{-1}\cO_\cZ}g^{-1}f^{-1}\cR_\cZ\quad(\cR_{\cY\to\cZ}\text{ is } \cO_\cY \text{ locally free})\\
&= \cR_{\cX\to\cZ}.
\end{align*}
Use now the fact that the natural morphism $g_!f_!\God^{\cbbullet}\to g_!\God^{\cbbullet} f_!\God^{\cbbullet}$ is an isomorphism, as $f_!\God^{\cbbullet}$ is $c$-soft, to get the first point. The second point is easy, as $\cR_{\cX\to\cY}$ is then $\cR_\cX$ locally free. For the third point, use Example \ref{ex:Sp}. The canonical isomorphism is obtained by applying $f_!$ to the diagram
\[
\xymatrix{
\God^{\cbbullet}\big(\cM\otimes_{\cO_\cX}\cE_{\cX/\cY}^{-\bbullet,\bbullet}\big)&\ar[l]^-{\sim}\cM\otimes_{\cO_\cX}\cE_{\cX/\cY}^{-\bbullet,\bbullet}\\
\ar[u]^-{\wr}\God^{\cbbullet}\big(\cM\otimes_{\cO_\cX}\wedge^{-\bbullet}\Theta_{\cX/\cY}\big)
}
\]
and by using the fact that, each term $\cM\otimes_{\cO_\cX}\cE_{\cX/\cY}^{-\bbullet,\bbullet}$ being $c$-soft on each fibre of $f$, its direct images $R^jf_!\big(\cM\otimes_{\cO_\cX}\cE_{\cX/\cY}^{-\bbullet,\bbullet}\big)$ vanish for any $j\neq0$.
\end{proof}

If $f$ is proper, or proper on the support of $\cM$, we have an isomorphism in the category $D^+(\Mod^r(\cR_\cY))$:
\[
\bR f_!\lefpar \cM\otimes_{\cR_\cX}^{\bL}\cR_{\cX\to\cY}\rigpar \isom\bR f_*\lefpar \cM\otimes_{\cR_\cX}^{\bL}\cR_{\cX\to\cY}\rigpar\defin \index{$fplus$@$f_+\cM$}f_+\cM.
\]
If moreover $\cM$ is good, then, for any compact set $\cK$ in $\cY$, it can be expressed in a neighbourhood of $f^{-1}(\cK)$ as a successive extension of modules which admit a resolution by \emph{coherent} induced $\cR_\cX$-modules in the neighbourhood of $f^{-1}(\cK)$. Arguing for instance as in \cite{Malgrange85}, one gets:

\begin{theoreme}\label{th:imdirgood}
Let $\cM$ be a good right $\cR_\cX$-module (or bounded complex) and let $f:X\to Y$ be a holomorphic map which is proper (or proper on the support of $\cM$). Then the object $f_+\cM$ is good and $\Char f_+\cM\subset f\big[(T^*f)^{-1}\Char\cM\big]$.\qed
\end{theoreme}

\begin{corollaire}\label{cor:holo}
If $\cM$ is holonomic and $f$ is proper on the support of $\cM$, then $f_+\cM$ is holonomic.\qed
\end{corollaire}

\begin{remarque}
Let $i:X\hto X'$ be a closed inclusion. Then a $\cR_\cX$-module $\cM$ is coherent (\resp holonomic, \resp strict) if and only if the $\cR_{\cX'}$-module $i_+\cM$ is so. Indeed, this is a local property and it is enough to verify it for the inclusion $X=X\times \{0\}\hto X\times \CC$. Denote by $t$ the coordinate on $\CC$. Then $i_+\cM=\cM\otimes_\CC\CC[\partiall_t]$, the $\cO_{\cX'}$-action on the degree $0$ terms being defined as the action of $\cO_\cX$. The assertion is then clear.
\end{remarque}

\begin{remarque}[Direct image for a left $\cR_\cX$-module] \label{rem:imdirleftmod}
The direct image for left $\cR_\cX$-modules is defined as usual by using the standard left-right transformation. It can be obtained by a formula analogous to that of Definition \ref{def:fdag}, using the bi-module $\cR_{\cY\leftarrow \cX}$.

Assume that $f:X=Z\times Y\to Y$ is the projection and put $n=\dim Z=\dim X/Y$. If $\cM$ is a left $\cR_\cX$-module, the direct image $f_\dag\cM$ can be computed directly with relative de~Rham complex:
\[
f_\dag\cM=f_!\God^{\cbbullet}\big(\Omega_{\cX/\cY}^{n+\bbullet}\otimes_{\cO_\cX}\cM\big),
\]
and, using the Dolbeault resolution,
\[
f_\dag\cM=f_!\big(\cE_{\cX/\cY}^{n+\bbullet}\otimes_{\cO_{\cX}} \cM\big).
\]
\end{remarque}

\begin{remarque}[Restriction to $\hb=\hb_o$]
If $\hb_o\neq0$, one has $\bL^*_{\hb_o}(f_\dag\cM)=f_\dag(\bL^*_{\hb_o}\cM)$, where the right-hand $f_\dag$ denotes the direct image of $\cD$-modules.
\end{remarque}

\subsection{Inverse images of $\cR_\cX$-modules}\label{subsec:iminv}

Let $f:X\to Y$ be a holomorphic map and let $\cM$ be a left $\cR_\cY$-module. The inverse image $\index{$fplush$@$f^+\cM$}f^+\cM$ is the object $$\cR_{\cX\to\cY}\otimes_{f^{-1}\cR_\cY}^{\bL}f^{-1}\cM.$$

In general, we only consider the case where $f$ is smooth (or, more generally, noncharacteristic, \cf\T\ref{sec:nonchar}). Then, $\cO_\cX$ being $f^{-1}\cO_\cY$-flat, we have $f^+\cM=\cO_\cX\otimes_{f^{-1}\cO_\cY}f^{-1}\cM$ with the structure of a left $\cR_\cX$-module defined as for $\cR_{\cX\to\cY}$. It is $\cR_\cX$-good if $\cM$ is $\cR_\cY$-good.

Assume on the other hand that $\cM$ is $\cO_\cY$-locally free of finite rank, but make no assumption on $f$. Then $f^+\cM=f^*\cM$ as an $\cO_\cX$-module, with the structure of a left $\cR_\cX$-module defined as above. It is also $\cO_\cX$-locally free of finite rank.

\section{Sesquilinear pairings on $\cR_{\cX}$-modules} \label{subsec:sesqui}
\subsection{Conjugation}\label{subsubsec:conj}
Denote by $X_\RR$ the $C^\infty$-manifold underlying $X$ and by $\ov X$ the complex analytic manifold conjugate to $X$, \ie\ $X_\RR$ equipped with the structural sheaf $\ov{\cO_X}$ of antiholomorphic functions. Notice that the conjugate of an open set of $X$ is the \emph{same} open set with a different sheaf of holomorphic functions. Recall that the conjugation $\index{$conj$@$\ovvv$ (conjugation)}\ovvv:\cO_X\to\cO_{\ov X}$ makes $\cO_{\ov X}$ a $\cO_X$-module. Given any $\cO_X$-module $\cF$, we denote by $\ov\cF$ its conjugate $\cO_{\ov X}$-module defined by
\[
\ov\cF=\cO_{\ov X}\ootimes_{\cO_X}\cF.
\]
One may extend $\ovvv$ as a ring morphism $\cD_X\to\cD_{\ov X}$ (in local coordinates, $\ov{\partial_{x_i}}=\partial_{\ov x_i}$) and define similarly $\ov\cF$ for $\cD_X$-modules.

On the $\PP^1$ factor, we will define a geometric conjugation that we also denote by $\ovvv$. It is induced by the involution $\hb\mto-1/\hb$. For notational convenience, we denote for a while by $c$ the usual conjugation functor on $\PP^1$. Given any open set $\Omega$ of $\PP^1$, denote by $\ov\Omega$ its image by the previous involution. Then, if $g(\hb)$ is a holomorphic function on $\Omega$, its \emph{conjugate} $\ov g(\hb)$ is by definition the holomorphic function $c\big(g(-1/c(\hb))\big)$ on $\ov\Omega$. Notice that we have
\[
\ov{\Omega_0}=\Omega_\infty,\quad\ov{\Omega_\infty}=\Omega_0\quad\text{and}\quad \ov\bS=\bS.
\]

We will now mix these two notions to get a conjugation functor on $X\times\PP^1$. We continue to denote by $c$ the usual conjugation functor on $X\times\PP^1$, but we keep the notation $\ovvv$ on $X$. Let $\sigma:c\PP^1\to\PP^1$ or $\PP^1\to c\PP^1$ denote the antilinear involution of $\PP^1$ defined by
\[
\sigma(c(\hb))=-1/\hb\qqbox{or} \sigma(\hb)=-1/c(\hb).
\]
Then, for any open set $\Omega\subset\PP^1$, $\sigma$ induces isomorphisms
\[
c(\Omega)\isom \ov\Omega.
\]
We also denote by $\sigma$ the inverse isomorphisms. Define $\sigma$ on $X\times\PP^1$ so that it is the identity on the $X$-factor.

We now have a conjugation functor $\ovvv\defin\sigma^*c$: given any holomorphic function $f(x,\hb)$ on an open set $U\times\Omega$ of $X\times\PP^1$, we put $\ov f(x,\hb)=c\big(f(x,-1/c(\hb))\big)$; therefore, $\ovvv$ defines a functor, also denoted by $\ovvv$ as above, which sends $\cR_{X\times\Omega}$-modules to $\ov{\cR_{X\times\Omega}}=\cR_{\ov X\times \ov\Omega}$-modules and conversely (in particular, $\ov{\partiall_{x_i}}=-\hbm\partial_{\ov x_i}$). Notice also that, if $m$ is a section of $\cM$ on $U\times\Omega$, it is also a section of $c\cM$ on $U\times\Omega$, that we denote by $c(m)$, and it defines a section $\ov m$ of $\ov\cM=\sigma^*c\cM$ on $U\times\sigma^{-1}(\Omega)=U\times\sigma(\Omega)$.

In particular, we have an identification
\[
\ov{\cO_{X\times\Omega}}=\cO_{\ov X\times\ov\Omega},\quad \ov{\cO_{X\times\Omega,(x,\hb)}}=\cO_{\ov X\times\ov\Omega,(x,-1/\hb)}
\]
by putting as above $\ov f(x,\hb)=c\big(f(x,-1/c(\hb))\big)$. Similarly, we have $\ov{\cR_{X\times\Omega}}=\cR_{\ov X\times\ov\Omega}$ when $\Omega\subset\CC^*$, $\ov{\cC^{\infty}_{X_\RR\times\bS}}=\cC^{\infty}_{X_\RR\times\bS}$, $\ov{\Dbh{X}}=\Dbh{X}$, \etc. Be careful, however, that these last identifications are not $\cO_\bS$-linear, but are linear \emph{over} the ``conjugation''
\begin{equation}\label{eq:conjug}
\begin{split}
\cO_\bS&\to\ov{\cO_\bS}=c\sigma^*\cO_\bS\\
\lambda(\hb)&\mto\ov\lambda(\hb)=c\big(\lambda(-1/c(\hb))\big).
\end{split}
\end{equation}

\subsection{Sesquilinear pairings}
We use notation and results of \S\T\ref{num:RX}--\ref{num:Dbh}.

Given two left $\cR_\cX$-modules $\cM'$ and $\cM''$, we define a \emph{sesquilinear pairing} between $\cMS'$ and $\cMS''$ as a $\cR_{(X,\ov X),\bS}$-linear pairing
\[
\index{$c$@$C$ (sesquilinear pairing)}C:\cMS'\ootimes_{\cO_\bS}\ov{\cMS''}\to\Dbh{X},
\]
where the left-hand term is equipped with its natural $\cR_{(X,\ov X),\bS}$-structure. Similarly, a sesquilinear pairing for right $\cR_\cX$-modules takes values in $\gCh{X}$.

\begin{remarque}\label{rem:lefrigsesqui}
It is easy to verify that the various functors ``going from left to right'', for $\cR_{\cX|\bS}$- and $\cR_{(X,\ov X),\bS}$-modules are compatible, and that they are compatible with sesquilinear pairings.
\end{remarque}

\begin{lemme}\label{lem:smoothsesqui}
If $\cM',\cM''$ are strict holonomic $\cR_\cX$-modules having their characteristic varieties $\Char\cM^{\prime\Cir}$, $\Char\cM^{\prime\prime\Cir}$ contained in the zero section, then any sesquilinear pairing $C$ between $\cM'$ and $\cM''$ takes values in $\cCh{\cX}$.
\end{lemme}

\begin{proof}
The assertion is local on $\cX$. According to Proposition \ref{prop:CK}, when restricted to $\cX^{\!\circ}$, $\cM',\cM''$ are $\cO_{\cX}$-locally free of finite rank and, given any $(x_{o},\hb_{o})\in X\times\bS$, we may find bases $\bme',\bme''$ of $\cM'_{(x_o,\hb_o)},\cM''_{(x_o,-\hb_o)}$ satisfying $\partiall\bme'=0$, $\partiall\bme''=0$. Lemma \ref{lem:DGhb} then shows in particular that $C$ takes values in $\cCh{\cX}$ (and more precisely in the subsheaf of functions which are real analytic with respect to $X$).
\end{proof}

\begin{exemple}[Basic holomorphic distributions]\label{ex:distr}
Let $\beta=\beta'+i\beta''$ be a complex number such that $\beta'\not\in-\NN^*$. Notice that there exists an open neighbourhood $\nb_\beta(\bS)$ such that the map $\nb_\beta(\bS)\to\CC$ defined by $\hb\mto (\beta\star\hb)/\hb$ (recall that the operation $\star$ is defined by \eqref{eq:defstar}) takes values in $\CC\moins(-\NN^*)$:
\begin{itemize}
\item
if $\beta''=0$, this function is constant and equal to $\beta'$;
\item
otherwise, $(\beta\star\hb)/\hb=-k\in-\NN^*$ is equivalent to $\hb= i\big((\beta'+k)/\beta''\pm\sqrt{1+[(\beta'+k)/\beta'']^2}\big)$, $k\in\NN^*$; the solutions belong to $i\RR\moins\{\pm i\}$ and do not accumulate at $\pm i$.
\end{itemize}

Putting $s=\beta\star\hb/\hb$ in Example \ref{ex:distran}\eqref{ex:distran3},
we therefore get a section $u_{\beta,\ell}=U_\ell$ of $\Dbh{X}$; we have $u_{\beta,\ell}=\mt^{2(\beta\star\hb)/\hb}\Lt^\ell/\ell!$. Then, the $u_{\beta,\ell}$ satisfy (if we set $u_{\beta,-1}=0$)
\begin{equation}\label{eq:tdtu}
(t\partiall_{t}-\beta\star\hb)u_{\beta,\ell}=u_{\beta,\ell-1}.
\end{equation}
Recall that $\ov{(\beta\star\hb)/\hb}=(\beta\star\hb)/\hb$, \ie $(\beta\star\hb)/\hb$ is ``real'', where $\ovvv$ is the conjugation defined in \T\ref{subsubsec:conj}, or equivalently, $\ov{\beta\star\hb}=\ov \beta\star\ov\hb$. Hence $\ov{u_{\beta,\ell}}=u_{\beta,\ell}$. Consequently, the $u_{\beta,\ell}$ also satisfy
\begin{equation}\label{eq:tdtbu}
(\ov t\partiall_{\ov t}-\ov\beta\star\ov\hb)u_{\beta,\ell}=u_{\beta,\ell-1}.
\end{equation}

Fix $\hb_o\in\bS$. Let $\beta\in\CC$ be such that $\reel(\beta\star\hb_o)/\hb_o=\beta'-\beta''(\hb_o+1/\hb_o)/2>-1$. Then, there exists a neighbourhood $\Deltag$ of $\hb_o$ in $\bS$ on which $u_{\beta,\ell}$ defines an element of $L^1_{\loc}(D\times\Deltag)$. If $B$ is a finite set of complex numbers $\beta$ such that
\begin{align}
\beta\in B&\implique\reel(\beta\star\hb_o)/\hb_o>-1,\label{eq:propB1}\\
\beta_1,\beta_2\in B\text{ and } \beta_1-\beta_2\in\ZZ&\implique\beta_1=\beta_2,\label{eq:propB2}
\end{align}
then the family $(u_{\beta,\ell})_{\beta\in B,\ell\in\NN}$ of elements of $L^1_{\loc}(D\times\Deltag)$ (with $\Deltag$ small enough, depending on $B$), is free over $C^\infty(D\times\Deltag)$: this is seen by considering the order of growth for any $\hb\in\Deltag$.
\end{exemple}

\section{The category $\RTriples(X)$}\label{subsec:Rtriples}

\Subsection{The category $\RTriples(X)$ and Hermitian adjunction} \label{subsec:somedef}
\begin{Definition}[of $\RTriples(X)$]\label{def:RtriplesX}
\begin{itemize}
\item
An object of $\index{$rtriples$@$\RTriples$}\RTriples(X)$ is a triple $\cT=(\cM',\cM'',C)$, where $\cM',\cM''$ are left $\cR_\cX$-modules and
\[
C:\cMS'\ootimes_{\cO_\bS}\ov{\cMS''}\to\Dbh{X}
\]
is a sesquilinear pairing.

\item
A morphism $\varphi:(\cM'_1,\cM''_1,C_1)\to (\cM'_2,\cM''_2,C_2)$ is a pair $(\varphi',\varphi'')$, where $\varphi':\cM'_2\to\cM'_1$ and $\varphi'':\cM''_1\to\cM''_2$ are $\cR_\cX$-linear and compatible with $C_1,C_2$, \ie satisfy
\[
C_1(\varphi'\cbbullet,\ov\star)=C_2(\cbbullet,\ov{\varphi''\star}).
\]

\item
The \emph{Tate twist} of an object of $\RTriples(X)$ is defined, for $k\in\hZZ$, by $$(\cM',\cM'',C)(k)=(\cM',\cM'',\thb^{-2k}C).$$
Any morphism between triples is also a morphism between the twisted triples (with the same twist), and we denote it in the same way.
\end{itemize}
\end{Definition}

The category $\RTriples(X)$ is abelian. If $\index{$tcurl$@$\cT$, $\cT^*$, $\cT(k)$}\cT$ is an object of $\RTriples(X)$ and $\lambda(\hb)\in\cO(\bS)$, the object $\lambda(\hb)\cdot\cT$ is by definition the object $(\cM',\cM'',\lambda(\hb)C)$. If $\varphi:\cT_1\to\cT_2$ is a morphism, then it is also a morphism between $\lambda(\hb)\cdot\cT_1$ and $\lambda(\hb)\cdot\cT_2$.

There are two functors $\cT\mto\cM'$ and $\cT\mto\cM''$, the first one to the category $\Mod(\cR_\cX)^{\rm op}$ (opposite category), and the second one to $\Mod(\cR_\cX)$. The identity morphism $\id_\cT$ is defined as $(\id_{\cM'},\id_{\cM''})$.

\begin{definition}[Adjunction]\label{def:hermdual}
Let $\cT=(\cM',\cM'',C)$ be an object of $\RTriples(X)$. Its \emph{Hermitian adjoint} $\cT^*$ is by definition
\[
\cT^*=(\cM',\cM'',C)^*\defin(\cM'',\cM',C^*),\quad\text{with}\quad \index{$ca$@$C^*$}C^*(\mu,\ov m)\defin\ov{C(m,\ov\mu)}.
\]
If $\varphi=(\varphi',\varphi''):\cT_1\to\cT_2$ is a morphism in $\RTriples(X)$, its \emph{adjoint} $\varphi^*\defin(\varphi'',\varphi')$ is a morphism $\cT_2^*\to\cT_1^*$ in $\RTriples(X)$.
\end{definition}

For $k\in\hZZ$, we choose a canonical isomorphism:
\begin{equation}\label{eq:adjcanonical}
((-1)^{2k}\id_{\cM'},\id_{\cM''}):\cT(k)\isom\cT^*(-k)^*.
\end{equation}
This isomorphism defines by adjunction an isomorphism $\cT^*(-k)\isom\cT(k)^*$ which is compatible with the composition of twists. These isomorphisms are equal to $\id$ if $k\in\ZZ$.

\begin{remarque}\label{rem:leftrightTriple}
We may define similarly the category $\RTriples(X)^r$: the objects are triples $(\cM',\cM'',C)$, where $\cM',\cM''$ are right $\cR_\cX$-modules and $C$ takes values in $\gCh{X}$. The Hermitian adjoint is defined similarly.

Given a left triple $(\cM',\cM'',C)$, the associated right triple is $(\cM^{\prime r},\cM^{\prime\prime r},C^r)$ with $\cM^{\prime r}=\omega_\cX\otimes_{\cO_\cX}\cM'$, and similarly for $\cM''$; moreover, $C^r$ is defined as
\[
C^r(\omega'\otimes m,\ov{\omega''\otimes \mu})=\varepsilon(n)\big(\itwopi \big)^n\cdot C(m,\ov{\mu})\omega'\wedge\ov{\omega''} \qquad(n=\dim X).
\]
Going from left to right is compatible with adjunction: $(C^r)^*=(C^*)^r$.
\end{remarque}

\begin{definition}\label{def:sesqdual}
A \emph{sesquilinear duality} of weight~$w\in\ZZ$ on $\cT$ is a morphism $\cS:\cT\to\cT^*(-w)$.
\end{definition}

Write $\cS=(S',S'')$ with $S',S'':\cM''\to\cM'$. Then $\cS$ is a morphism if and only if $C,S',S''$ satisfy, for local sections $\mu_1,\mu_2$ of $\cM''$,
\[
C(S'\mu_1,\ov\mu_2)=\thb^{2w}C^*(\mu_1,\ov{S''\mu_2}).
\]
Let $k\in\hZZ$. Put $\index{$scurl$@$\cS$, $\cS(k)$}\cS(k)=\big((-1)^{2k}S',S''\big)$. Then $\cS$ is a sesquilinear duality of weight~$w$ on $\cT$ if and only if $\cS(k):\cT(k)\to(\cT(k))^*(-w+2k)$ is a sesquilinear duality of weight~$w-2k$ on $\cT(k)$. In particular, $\cS(w/2)$ is a sesquilinear duality of weight~$0$ on $\cT(w/2)$.

Notice that $\cS(k)$ is obtained by composing $\cS:\cT(k)\to\cT^*(-w+k)$ with the canonical isomorphism chosen above $\cT^*\isom\cT(k)^*(k)$, applied to the $(-w+k)$ twisted objects.

If $\cS:\cT\to\cT^*(-w)$ has weight~$w$, associate to it the sesquilinear pairing on $\bS$:
\begin{equation}\label{eq:sesqupair}
h_{\bS,\cS}\defin \thb^{-w} C\circ(S''\otimes\id):\cMS''\ootimes_{\cO_\bS}\ov{\cMS''} \to\Dbh{X}.
\end{equation}
Notice that $h_{\bS,\cS}=h_{\bS,\cS(k)}$ for any $k\in\hZZ$. We denote this pairing by $h_\bS$ when $\cS$ is fixed.

\begin{definition}\label{def:hermw}
A sesquilinear duality of weight~$w$ on $\cT$ is said to be \emph{Hermitian} if it satisfies
\[
\cS^*=(-1)^w\cS,\quad\text{\ie\ }S'=(-1)^wS''.
\]
\end{definition}

The exponent $w$ is useful to get that, if $\cS$ is Hermitian, then $\cS(k)$ is Hermitian for any $k\in\hZZ$. If $\cS:\cT\to\cT^*(-w)$ is Hermitian, then its associated sesquilinear pairing $h_{\bS,\cS}$ on $\bS$ is Hermitian, \ie $h_\bS^*=h_\bS$.

\begin{remarque}\label{rem:hermdual}
Let $\cS:\cT\to\cT^*$ be a Hermitian duality of weight~$0$. Assume that $\cS$ is an isomorphism. Put $\cM=\cM''$. Then $\cT$ is isomorphic to the triple $(\cM,\cM,h_{\bS})$ which is self-adjoint and, under this isomorphism, $\cS$ corresponds to $(\id_\cM,\id_\cM)$. Indeed, the isomorphism is nothing but
\[
(S'',\id_{\cM''}):(\cM',\cM'',C)\to(\cM'',\cM'',h_{\bS}),
\]
as $S'=S''$ by assumption. This trick, combined with a Tate twist by $(w/2)$, reduces the study of polarized twistor $\cD$-modules (\cf Definition \ref{def:polarization}) to that of objects of the form $[(\cM,\cM,C),(\id_{\cM},\id_{\cM})]$.
\end{remarque}

\subsection{Smooth triples}\label{subsec:smoothtriple}
We say that an object $\cT=(\cM',\cM'',C)$ of $\RTriples(X)$ is \emph{smooth} if $\cM'$ and $\cM''$ are $\cO_{\cX}$ locally free of finite rank. It follows from Lemma \ref{lem:smoothsesqui} that, for a smooth triple $\cT$, the sesquilinear pairing $C$ takes values in $\cC^\infty_{X_\RR\times\bS}$.

\begin{definition}[Inverse image]\label{def:restrsmooth}
Let $f:Y\to X$ be a holomorphic map between complex analytic manifolds $Y$ and $X$. The inverse image by $f$ of the left smooth triple $\cT=(\cM',\cM'',C)$ is the smooth triple $f^+\cT=(f^+\cM',f^+\cM'',f^+C)$, where $f^+\cM=f^*\cM$ is taken in the sense of $\cO$-modules with connections (\cf \T\ref{subsec:iminv}) and $\index{$fplushc$@$f^+C$, $f^+\cT$, $f^+\cS$}f^+C(1\otimes m',\ov{1\otimes m''})=C(m',\ov{m''})\circ f$.
\end{definition}

\begin{Remarques}
\begin{itemize}
\item
The inverse image by $f$ of $C$ by is well defined because $C$ takes values in $C^\infty$ functions on $X$, and not only in distributions on $X$.

\item
The inverse image by $f$ of a morphism is the usual inverse image of each component of the morphism.

\item
The inverse image functor commutes with Tate twist and Hermitian adjunction.
\end{itemize}
\end{Remarques}

The last remark allows one to introduce:

\begin{definition}[Inverse image of a sesquilinear duality]\label{def:restrsesqui}
Assume that $\cT$ is smooth. The inverse image $$f^+\cS:f^+\cT\to f^+\cT^*(-w)$$ of a sesquilinear pairing $\cS=(S',S''):\cT\to\cT^*(-w)$ of weight~$w$ is the morphism $\big(f^*S',f^*S''\big)$.
\end{definition}

A Hermitian sesquilinear duality of weight~$w$ remains Hermitian of weight~$w$ after inverse image.

\subsection{Differential graded $\cR$-triples}\label{subsec:dgtriples}
Consider the category of graded $\cR$-triples $\cT=\oplus_{j\in\ZZ}\cT^j$. Morphisms are graded. We will follow the usual convention when writing indices: $\cT_j=\cT^{-j}$. For $k\in\ZZ$, put $(\cT[k])_j=\cT_{j-k}$ or $(\cT[k])^j=\cT^{j+k}$. The shift $[\cbbullet]$ and the twist $(\cbbullet)$ commute. A differential $d$ is a morphism $\cT\to\cT[1](\varepsilon)$ such that $d\circ d=0$, for some $\varepsilon\in\ZZ$.

The Hermitian adjunction is defined by $\cT^*=\oplus_j(\cT^*)^j$ with the grading $(\cT^*)^j=(\cT^{-j})^*$. We have $(\cT[k])^*=\cT^*[-k]$.

A sesquilinear duality of weight~$w$ on $\cT$ is a (graded) morphism $\cS:\cT\to\cT^*(-w)$, \ie a family of morphisms $\cS^j:\cT^j\to\cT^{-j*}(-w)$.

A morphism $\varphi:\cT\to\cT[k](\ell)$ is \emph{selfadjoint} (\resp \emph{skewadjoint}) with respect to $\cS$ if the following diagram commutes (\resp anticommutes):
\[
\xymatrix{
\cT\ar[r]^-{\cS}\ar[d]_-\varphi&\cT^*(-w)\ar[d]^-{\varphi^*}\\
\cT[k](\ell)\ar[r]^-{\cS}&\cT^*[k](\ell-w)
}
\]
A differential $d$ is selfadjoint with respect to $\cS$ if and only if $\cS$ is a morphism of complexes $(\cT,d)\to(\cT^*(-w),d^*)$.

Filtered objects are defined similarly: a decreasing filtration $F^{\cbbullet}\cT$ of $\cT$ consists of the datum of decreasing filtrations $F^{\cbbullet}\cM',F^{\cbbullet}\cM''$ such that, for any $k\in\ZZ$, we have $C(F^{-k+1}\cM',\ov{F^k\cM''})=0$; then $F^k\cT=(\cM'/F^{-k+1}\cM',F^k\cM'',C)$ is well defined and we have $\gr^k_F\cT=(\gr^{-k}_F\cM',\gr^k_F\cM'',C)$ (where we still denote by $C$ the pairing naturally induced by $C$).

Define the decreasing filtration $F^{\cbbullet}\cT^*$ by $F^k(\cT^*)= (\cM''/F^{-k+1}\cM',F^k\cM',C^*)$. Then $\gr^k_F(\cT^*)=(\gr^{-k}_F\cT)^*$ and, considering the total graded object $\gr_F\cT$, this is compatible with the definition above of adjunction for graded objects.

\begin{lemme}\label{lem:d1}
Let $(\cT,F^{\cbbullet})$ be a filtered $\cR$-triple, equipped with a filtered differential $d$ and a filtered sesquilinear duality $\cS$ of weight~$w$. Assume that $d$ is selfadjoint with respect to $\cS$. Then $\cS$ induces a natural sesquilinear duality of weight~$w$ on $E_1^{j}=\oplus_pH^{j}(\gr_F^p\cT)$ with respect to which the differential $d_1:E_1^j\to E_1^{j+1}$ is selfadjoint.\qed
\end{lemme}

\subsection{Direct images in $\RTriples$}\label{subsec:imdirsesqui}
The purpose of this paragraph is to define, for any holomorphic map $f:X\to Y$ and any object $\cT=(\cM',\cM'',C)$ of $\RTriples(X)$, an object $f_\dag\cT$ in the derived category $D^+(\RTriples(Y))$. Such an object is a complex $((\cN^{\prime\bbullet})^{\rm op},\cN^{\prime\prime\bbullet},C^{\cbbullet})$, where the first term is a complex in the opposite category $\Mod(\cR_\cY)^{\rm op}$ (given a complex $\cN^{\cbbullet}$ in $\Mod(\cR_\cY)$, we put $\cN^{{\rm op}\,k}=\cN^{-k}$) and the second term a complex in $\Mod(\cR_\cY)$. Therefore $C^k$ is a morphism $\cNS^{\prime -k}\otimes\cNS^{\prime\prime k}\to\gCh{Y}$ which is compatible with the differentials, \ie the following diagram commutes:
\[
\xymatrix@C=1.5cm{
\cNS^{\prime -k}\otimes\cNS^{\prime\prime k}\ar@<3ex>[d]^-{d''}\ar[r]^-{C^k}&\gCh{Y}\ar@<-1ex>@{=}[d]\\
\cNS^{\prime -k-1}\otimes\cNS^{\prime\prime k+1} \ar@<5ex>[u]^-{d'}\ar[r]^-{C^{k+1}}&\gCh{Y}}
\]
The complex $f_\dag\cT$ will take the form $((f_\dag\cM')^{\rm op},f_\dag\cM'',f_\dag C)$, where $f_\dag\cM',f_\dag\cM''$ are defined in \T\ref{subsec:imdir} and $(f_\dag\cM')^{\rm op}$ is the corresponding complex in the opposite category. We will therefore obtain a family of sesquilinear pairings
\begin{equation}\label{eq:fdagC}
\index{$fdagC$@$f_\dag C$, $f_\dag\cT$}f^i_\dag C:\cH^{-i}(f_{\dag}\cMS') \otimes_{\cO_\bS}\ov{\cH^i(f_{\dag}\cMS'')}\to\gCh{Y}.
\end{equation}

We will define $f_\dag C$ when $f$ is a projection and when $f$ is an inclusion. For a general $f$, we write it as the composition of its graph inclusion $i_f$ and of the canonical projection $p_f$, and put $f_\dag=p_{f\dag}i_{f\dag}$.

We will prove
\begin{equation}\label{eq:fdag*}
f_\dag^j(\cT^*)=(f_\dag^{-j}\cT)^*
\end{equation}
and, whenever $f$ and $g$ are composable,
\begin{equation}\label{eq:fgdag}
(g\circ f)_\dag\cT=g_\dag(f_\dag\cT)
\end{equation}

\subsubsection*{(a) Case of a projection $f:X=Z\times Y\to Y$ and left triples}
Recall (\cf Remark \ref{rem:imdirleftmod}) that we have $(f_{\dag}\cM)_\bS= f_!\big(\cE_{X_\RR\times\bS/Y_\RR\times\bS}^{n+\bbullet}\otimes_{\cO_{\cX|\bS}}\cMS\big)$. Consider the family of morphisms
\[
(f_\dag C)^j: f_!\big(\cE_{X_\RR\times\bS/Y_\RR\times\bS}^{n-j}\otimes_{\cO_{\cX|\bS}}\cMS'\big)\ootimes_{\cO_\bS} \ov{f_!\big(\cE_{X_\RR\times\bS/Y_\RR\times\bS}^{n+j}\otimes_{\cO_{\cX|\bS}}\cMS''\big)}\to\Dbh{Y}
\]
defined by
\begin{equation}\label{eq:lefrigimdir}
(\eta^{n-j}\otimes m')\otimes \ov{(\eta^{n+j}\otimes m'')}\mto\dfrac{\varepsilon(n+j)}{(2i\pi)^n}\int_fC(m',\ov{m''})\,\eta^{n-j}\wedge \ov{\eta^{n+j}}.
\end{equation}

\begin{Lemme}\label{lem:imdiradjointcomp}
\begin{enumerate}
\item\label{lem:imdiradjointcomp1}
$f_\dag\cT\defin\big((f_\dag\cM')^{\rm op},f_\dag\cM'',f_\dag C\big)$ is an object of $D^+(\RTriples(Y))$.
\item\label{lem:imdiradjointcomp2}
We have $(f_\dag C^*)^{-j}=((f_\dag C)^j)^*$.
\item\label{lem:imdiradjointcomp3}
If $f,g$ are two composable projections, we have $(g\circ f)_\dag C=g_\dag(f_\dag C)$.
\end{enumerate}
\end{Lemme}

\begin{proof}
We have, by $\cR_{(X,\ov X),\bS}$-linearity of $C$, and up to multiplication by $(-1)^n$,
\begin{equation*}
\begin{split}
(f_\dag C)^j\big(\nabla(&\eta^{n-j-1}\otimes m')\otimes \ov{(\eta^{n+j}\otimes m'')}\big)\\ &=\dfrac{\varepsilon(n+j)}{(2i\pi)^n}\int_f \big(C(m',\ov{m''})d\eta^{n-j-1}+C(\nabla m',\ov{m''})\wedge\eta^{n-j-1}\big)\wedge\ov{\eta^{n+j}}\\
&=\dfrac{\varepsilon(n+j)}{(2i\pi)^n}\int_f \big(C(m',\ov{m''})d\eta^{n-j-1}+d'C(m',\ov{m''})\wedge\eta^{n-j-1}\big)\wedge\ov{\eta^{n+j}}
\end{split}
\end{equation*}
and similarly
\begin{equation*}
\begin{split}
(f_\dag C)^{j+1}\big((&\eta^{n-j-1}\otimes m')\otimes \ov{\nabla(\eta^{n+j}\otimes m'')}\big)\\ &=\dfrac{\varepsilon(n+j+1)}{(2i\pi)^n}\int_f \eta^{n-j-1}\wedge\big(C(m',\ov{\nabla m''})\wedge\ov{\eta^{n+j}}+ C(m',\ov{m''})d\ov{\eta^{n+j}}\big)\\
&=\dfrac{\varepsilon(n+j+1)}{(2i\pi)^n}\int_f \eta^{n-j-1}\wedge\big(d''C(m',\ov{m''})\wedge\ov{\eta^{n+j}}+C(m',\ov{m''})d\ov{\eta^{n+j}}\big).
\end{split}
\end{equation*}
Using that $\varepsilon(n+j+1)=(-1)^{n+j}\varepsilon(n+j)=(-1)^{n-j}\varepsilon(n+j)$, Stokes Formula implies that both terms are equal, hence \ref{lem:imdiradjointcomp}\eqref{lem:imdiradjointcomp1}.

\medskip
To prove \ref{lem:imdiradjointcomp}\eqref{lem:imdiradjointcomp2}, remark that we have
\begin{align*}
(f_\dag C^*)^{-j}\big((\eta^{n+j}\otimes m'')\otimes \ov{(\eta^{n-j}\otimes m')}\big)&=
\dfrac{\varepsilon(n-j)}{(2i\pi)^n}\int_fC^*(m'',\ov{m'})\,\eta^{n+j}\wedge \ov{\eta^{n-j}}\\
&= \ov{\dfrac{(-1)^j\varepsilon(n-j)}{(2i\pi)^n}\int_fC(m',\ov{m''})\,\eta^{n-j}\wedge \ov{\eta^{n+j}}}\\
&=((f_\dag C)^j)^*\big((\eta^{n+j}\otimes m'')\otimes \ov{(\eta^{n-j}\otimes m')}\big),
\end{align*}
the last equality following from $(-1)^j\varepsilon(n-j)=\varepsilon(n+j)$.

\medskip
Last, \ref{lem:imdiradjointcomp}\eqref{lem:imdiradjointcomp3} follows from
\[
(-1)^{(m-k)(n+j)}\varepsilon(m+n+j+k)=(-1)^{(m+k)(n+j)}\varepsilon(m+n+j+k)=\varepsilon(n+j)\varepsilon(m+k).\qedhere
\]
\end{proof}

\subsubsection*{(b) Case of a closed inclusion $i:X\hto Y$}
Consider first the case of right triples. In this case, $i_\dag\cM$ is generated by $i_*\cM$ as a $\cR_\cY$-module. The pairing $i_\dag C$ is extended by $\cR_\cY$-linearity from its restriction to $i_*\cMS'\otimes\ov{i_*\cMS''}$, where it is defined as the composition of $C$ with the direct image of currents
\begin{align*}
\gCh{X}&\To{i_{++}}\gCh{Y}\\
u&\mto i_{++}u:\psi\mto\langle u,\psi\circ i\rangle\quad(\psi\in\cC^{\infty}_c(Y\times\bS))
\end{align*}
For left triples, define $i_\dag C$ in such a way that $(i_\dag C)^r=i_\dag(C^r)$, where $C^r$ is defined in Remark \ref{rem:leftrightTriple}. If $X$ is a submanifold of $Y$ defined by $x_j=0$ ($j\in J\subset\{1,\dots,n\}$), we identify $i_\dag\cM$ with $\cM\otimes_\CC\CC[{\partiall_{x_j}}_{j\in J}]$ and we have, for any $\varphi\in\cE^{n,n}_c(Y\times\bS/\bS)$,
\[
\langle(i_\dag C)(m'\otimes 1,\ov{m''\otimes1}),\varphi\rangle=\Big\langle C(m',\ov{m''}),\frac{\varphi}{\prod_{j\in J} \itwopi dx_j\wedge d\ov {x_j}}\Big\rangle.
\]
The conclusions of Lemma \ref{lem:imdiradjointcomp} clearly hold for the case of closed inclusions.

\begin{remarque*}
It may be more convenient to write the form $\itwopi dx_j\wedge d\ov {x_j}$ as
\[
\frac{1}{2i\pi}\,\frac{dx_j}{\hb}\wedge \frac{d\ov {x_j}}{\ov\hb}
\]
when checking the compatibilities below.
\end{remarque*}

\subsubsection*{(c) General case}
For a general $f$, define $f_\dag$ as the composition $p_{f\dag}i_{f\dag}$ as indicated above. One has to check first that this is compatible with the previous definitions when $f$ is a projection or an inclusion. This is mainly reduced to checking this for a closed inclusion $f:X\hto Y$.

Then, to conclude that $(g\circ f)_\dag=g_\dag f_\dag$, and to end the proof of \eqref{eq:fdag*} and \eqref{eq:fgdag}, it is enough to show that, in the following cartesian diagram,
\[
\xymatrix@=1.5cm{
X\times Y\ar@{^{ (}->}[r]^-{\id\times i}\ar[d]^-p&X\times Y'\ar[d]^-{p'}\\
Y\ar@{^{ (}->}[r]^-{i}&Y'
}
\]
we have $p'_\dag\circ(\id\times i')_\dag=i_\dag p_\dag$.

We leave both computations to the reader.\qed

\subsection{The Lefschetz morphism}\label{subsec:lef}
Let $c\in H^2(X,\CC)$ be a real $(1,1)$-class (the Chern class of a holomorphic line bundle for instance). We will show that it induces, for any $\cR_\cX$-module $\cM$, a morphism
\[
L_{c}:\cH^if_\dag\cM\to\cH^{i+2}f_\dag\cM
\]
in such a way that, if $\cT=(\cM',\cM'',C)$ is an object of the category $\RTriples(X)$ or the derived category $D(\RTriples(X))$, then $L_{c}$ induces a functorial morphism in $\RTriples(Y)$:
\begin{equation}\label{eq:lefschetzRtriple}
\index{$lc$@$L_c$, $\cL_{c}$, $\cL_\omega$}\cL_{c}=(L'_c,L''_c)\defin(-L_c,L_c):f_\dag^j\cT\to f_\dag^{j+2}\cT(1).
\end{equation}
It will be enough to apply the following computation to a closed real $(1,1)$-form $\omega\in\Gamma(X,\cE_X^{1,1})$ representing~$c$.

Let $i_f:X\hto X\times Y$ be the graph inclusion of $f$ and let $p:X\times Y\to Y$ be the projection. We have $f_\dag\cT=p_\dag i_{f,\dag}\cT$. So we may replace $\cT$ with $i_{f,\dag}\cT$ and assume that $f$ is the projection $X=Z\times Y\to Y$.

Let $\omega$ be a real closed $(1,1)$-form on $X$. We define $L_\omega:f_\dag\cM\to f_\dag\cM[2]$ as the morphism induced by $\hbm\omega\wedge$ on $f_!\big(\cE^{n+\bbullet}_{\cX/\cY}\otimes_{\cO_\cX}\cM\big)$, where $n=\dim X/Y$ (\cf Remark \ref{rem:imdirleftmod}). It is a morphism because $\omega$ is closed and the morphism induced on the cohomology depends on $c$ only.

We will now show that, as $\omega$ is real $2$-form, we have a morphism
\[
\cL_\omega=(L'_\omega,L''_\omega)\defin(-L_\omega,L_\omega):f_\dag^j\cT\to f_\dag^{j+2}\cT(1).
\]
As $\omega$ is a real $2$-form, we have, by Formula \eqref{eq:lefrigimdir} and using that $\varepsilon(n+j+2)=-\varepsilon(n+j)$,
\begin{equation*}
\begin{split}
f_\dag^jC\big(L'_\omega(\eta^{n-j-2}\otimes m),&\ov{\eta^{n+j}\otimes \mu}\big)=-\dfrac{\varepsilon(n+j)}{(2i\pi)^n}\int_f C(m,\ov{\mu})\hbm\omega\wedge\eta^{n-j-2}\wedge\ov{\eta^{n+j}}\\
&=-\thb^{-2}\dfrac{\varepsilon(n+j)}{(2i\pi)^n}\int_f C(m,\ov{\mu})\eta^{n-j-2}\wedge\ov{\hbm\omega\wedge\eta^{n+j}}\\
&=\thb^{-2}\dfrac{\varepsilon(n+j+2)}{(2i\pi)^n}\int_f C(m,\ov{\mu})\eta^{n-j-2}\wedge\ov{\hbm\omega\wedge\eta^{n+j}}\\
&=\thb^{-2} f_\dag^{j+2}C\big(\eta^{n-j-2}\otimes m,\ov{L''_\omega(\eta^{n+j}\otimes \mu)}\big).
\end{split}
\end{equation*}

\chapter{Smooth twistor structures}\label{chap:smtw}

The notion of a twistor structure has been introduced by C\ptbl Simpson in \cite{Simpson97} (see also \cite{Simpson97b}) in order to extend the formalism of variations of Hodge structures to more general local systems. The purpose of this chapter is to review the basic definitions in the language of $\RTriples$, in order to extend them to $\cR_\cX$-modules. Following Simpson, we express the Hodge theory developed in \cite{Simpson92} in terms of twistors and recall the proof of the Hodge-Simpson theorem \ref{th:HodgeSimpson}. Nevertheless, this theorem will not be used in its full strength for twistor $\cD$-modules, according to the method of M\ptbl Saito. It is useful, however, to understand the main result (Theorem \ref{th:imdirtwistor}).

\section{Twistor structures in dimension $0$}\label{sec:twst0}
In order to explain such a definition, we will first give details on the simplest example, \ie when $X$ is reduced to a point. We will first give definitions for weight~$0$, then give the way to obtain a twistor of weight~$0$ from a twistor of weight~$w$: this is the analogue of the Weil operator $C$ in Hodge theory. The convention taken here will look convenient later on.

\subsection{Twistor structures in dimension $0$ after C\ptbl Simpson \cite{Simpson97}}
A \emph{pure twistor} of rank $d$ and weight~$w$ is a vector bundle on $\PP^1$ isomorphic to $(\cO_{\PP^1}(w))^d$. A~\emph{mixed twistor} is a vector bundle $\wt\cH$ on $\PP^1$ equipped with an increasing filtration $W_\bbullet$ indexed by $\ZZ$ such that, for each $\ell\in\ZZ$, $\gr_\ell^W\wt\cH$ is pure of weight~$\ell$. A pure twistor structure can also be viewed as a mixed twistor structure in a natural way. A morphism of mixed twistor structures is a morphism of vector bundles which respects the filtrations $W$. There are no nonzero morphisms between pure twistor structures when the weight of the source is strictly bigger than the weight of the target. The category of pure twistor structures of weight~$w$ is equivalent to the category of $\CC$-vector spaces, hence it is abelian. The category of mixed twistor structures is therefore abelian and any morphism is strict with respect to $W$ (see \loccit, and also \cite[Th\ptbl 1.2.10]{DeligneHII}, \cite[Lemme 5.1.15]{MSaito86}).

\subsection{Twistor structures in dimension $0$ as objects of $\RTriples$} \label{subsubsec:smtw0}
We will give an equivalent definition of twistor structures which will be extended to arbitrary dimensions in \T\ref{subsec:defDtwist}. A \emph{twistor structure of rank $d$ and of weight~$w$} consists of the data of two free $\cO_{\Omega_0}$-modules $\index{$hcurlprime$@$\cH'$}\cH'$ and $\cH''$ and of a $\cO_\bS$-linear pairing
\begin{equation}\label{eq:twC0}
\index{$c$@$C$ (sesquilinear pairing)}C:\cHS'\ootimes_{\cO_\bS}\ov{\cHS''}\to\cO_\bS
\end{equation}
(in other words, $(\cH',\cH'',C)$ is an object of $\RTriples(\text{pt})$ where $C$ takes values in $\cO_\bS$ instead of only $\cC^0_\bS$), such that \eqref{cond:twistora} and \eqref{cond:twistorb} below are satisfied:

{\def\theenumi{\alph{enumi}}
\begin{enumerate}
\item\label{cond:twistora}
The sesquilinear pairing \eqref{eq:twC0} is \emph{nondegenerate}, \ie its matrix in any local basis of $\cHS'$ and $\cHS''$ is invertible, so that the associated $\cO_\bS$-linear morphism induces an \emph{isomorphism}
\begin{equation}\label{eq:tw0iso}
\ov{\cHS''}\isom \cHS^{\prime\vee}\defin\cHom_{\cO_{\Omega_0}}(\cH',\cO_{\Omega_0})_\bS.
\end{equation}

\item\label{cond:twistorb}
The locally free $\cO_{\PP^1}$-module $\index{$hcurltilde$@$\wt\cH$}\wt\cH$ obtained by gluing $\cH^{\prime\vee}$ (dual of $\cH'$, chart $\Omega_0$) and $\ov{\cH''}$ (conjugate of $\cH''$, chart $\Omega_\infty$; \cf \T\ref{subsubsec:conj}) using \eqref{eq:tw0iso} is \emph{isomorphic to} $\cO_{\PP^1}(w)^d$.
\end{enumerate}}

Having chosen $0$ and $\infty$ on $\PP^1$, the category of twistor structures of weight~$w$ is a full subcategory of $\RTriples(\text{pt})$: a morphism $\varphi:(\cH'_1,\cH''_1,C_1)\to (\cH'_2,\cH''_2,C_2)$ consists of the data of morphisms $\varphi':\cH'_2\to\cH'_1$, $\varphi'':\cH''_1\to\cH''_2$ of $\cO_{\Omega_0}$-modules such that $C_1(\varphi'(m'_2)\otimes\ov{m''_1})=C_2(m'_2\otimes \ov{\varphi''(m''_1)})$ for any $\hb_o\in\bS$ and all $m'_2\in\cH'_{2,\hb_o}$, $m''_1\in\cH''_{1,-\hb_o}$ (recall that $\sigma(\hb_o)=-\hb_o$ when $\hb_o\in\bS$). This category is clearly equivalent to the category of semistable vector bundles of slope $w$ on $\PP^1$, or to the category of finite dimensional $\CC$-vector spaces.

\smallskip
The notion of a \emph{Tate twist} is also well defined: recall that, for $k\in\hZZ$, we put
\begin{equation}\label{eq:Tate}
(\cH',\cH'',C)(k)\defin(\cH',\cH'',\thb^{-2k}C).
\end{equation}
The weight of $(\cH',\cH'',C)(k)$ is $w-2k$, if $(\cH',\cH'',C)$ has weight~$w$.

\subsubsection*{Weil reduction to weight~$0$}
Given a twistor $\cT=(\cH',\cH'',C)$ of weight~$w$, its associated twistor of weight~$0$ is by definition $\wt\cT=\cT(w/2)=(\cH',\cH'',\thb^{-w}C)$. Notice that, if $\varphi:\cT_1\to\cT_2$ is a morphism of twistors of weight~$w$ then $\varphi$ also defines a morphism $\wt\cT_1\to\wt\cT_2$. Notice also that, for any $k\in\hZZ$, the twistors $\cT$ and $\cT(k)$ have the same associated twistor of weight~$0$.

\begin{remarque}
A triple $(\cH',\cH'',C)$ is a pure twistor of weight~$w$ if and only of one can find $\CC$-vector spaces $H'\subset \Gamma(\Omega_0,\cH')$, $H''\subset\Gamma(\Omega_0,\cH'')$ of finite dimension (equal to $\rg\cH'=\rg\cH''$) such that $\cH'=\cO_{\Omega_0}\otimes_\CC H'$, $\cH''=\cO_{\Omega_0}\otimes_\CC H''$, the restriction of $\Gamma(\bS,C)$ to $H'\otimes \ov{H''}$ takes values in $\hb^w\CC\subset\Gamma(\bS,\cO_\bS)$, and induces an isomorphism $\ov{H''}\isom H^{\prime\vee}$.
\end{remarque}

We may also define the category of mixed twistor structures as the category of triples with a finite filtration, such that $\gr_\ell^W$ is a pure twistor structure of weight~$\ell$. If $(\cH',\cH'',C)$ is a mixed twistor structure, then $\cH'$ and $\cH''$ are locally free and $C$ is nondegenerate. This category is equivalent to the category of mixed twistor structures in the sense of Simpson.

\subsection{Hermitian adjunction and polarization}\label{subsec:hermdualpol}
Recall that the \emph{Hermitian adjoint} $\cT^*$ of the triple $\cT=(\cH',\cH'',C)$ is the triple $\cT^*=(\cH'',\cH',C^*)$ (\cf Definition \ref{def:hermdual}). If $\cT$ is a twistor of weight~$w$, then $\cT^*$ has weight~$-w$. Recall also that $(\cT(k))^*=\cT^*(-k)$ for $k\in\ZZ$ and that we have chosen an isomorphism $(\cT(k))^*\isom\cT^*(-k)$ if $k\in\hZZ$.

\begin{definition}
Let $\cT=(\cH',\cH'',C)$ be a twistor structure of pure weight~$w$. A \emph{Hermitian duality} of $\cT$ is an isomorphism $\cS:\cT\to\cT^*(-w)$ which is a Hermitian sesquilinear pairing of weight~$w$ in the sense of Definition \ref{def:hermw}, \ie which satisfies $\cS^*=(-1)^w\cS$.
\end{definition}

The associated sesquilinear pairing (\cf \eqref{eq:sesqupair})
\begin{equation}\label{eq:hermit}
\index{$hsbold$@$h_\bS$}h_\bS\defin \thb^{-w}C\circ (S''\otimes \id): \cHS''\ootimes_{\cO_\bS} \ov{\cHS''}\to\cO_\bS
\end{equation}
is nondegenerate and \emph{Hermitian}, \ie satisfies $h_\bS(m,\ov \mu)=\ov{h_\bS(\mu,\ov m)}$ for local sections $m,\mu$ of $\cHS''$.

\medskip
Assume now that $w=0$. We therefore have $S'=S''\defin S$. Let $\wt \cH$ be the trivial $\cO_{\PP^1}$-bundle defined in (b) above and consider its conjugate $\ov{\wt\cH}$ in the sense of \T\ref{subsubsec:conj}. This is the trivial vector bundle obtained by gluing $\cH''$ (chart $\Omega_0$) and $\ov{\cH^{\prime\vee}}$ (chart $\Omega_\infty$), using the conjugate of the map \eqref{eq:tw0iso} (also denoted by $C$) induced by $C$, that we denote by $\ov C$. We may define a $\cO_{\PP^1}$-linear pairing $\wt h:\ov{\wt\cH}\otimes_{\cO_{\PP^1}}\wt\cH\to\cO_{\PP^1}$ in the following way:
\begin{itemize}
\item
in the chart $\Omega_0$, $\wt h$ is the pairing $\langle S\cbbullet,\cbbullet\rangle:\cH''\otimes_{\cO_{\Omega_0}}\cH^{\prime\vee}\to\cO_{\Omega_0}$, where $\Crochet$ is the standard duality pairing;
\item
in the chart $\Omega_\infty$, $\wt h$ is the pairing $\langle\cbbullet, \ov{S\cbbullet}\rangle:\ov{\cH^{\prime\vee}}\otimes_{\cO_{\Omega_\infty}}\ov{\cH''}\to\cO_{\Omega_\infty}$;
\item
that both definitions agree near $\bS$ follows from the fact that $h_\bS$ is Hermitian; using $C$, one may identify $\wt h_\bS$ with $h_\bS$; for the same reason, $\wt h$ is Hermitian.
\end{itemize}

Denote by $\ov H$ the rank $d$ vector space $H^0(\PP^1,\wt\cH)$. Its conjugate $H$ is canonically identified with $H^0(\PP^1,\ov{\wt\cH})$. Therefore, $\wt h$ induces a Hermitian pairing
\begin{equation}\label{eq:ccP}
h=\pi_*\wt h:H\ootimes_\CC \ov H\to\CC.
\end{equation}

\begin{remarque}\label{rem:h}
We have a canonical inclusion $\ov H\subset \Gamma(\bS,\ov{\cHS''})$ (restriction of sections). We also have a canonical conjugate inclusion $H\subset\Gamma(\bS,\cHS'')$. Then $h$ may be identified with the restriction of $\Gamma(\bS,h_\bS)$ to $H\otimes_\CC \ov H$ by these inclusions.
\end{remarque}

If $\cT$ has weight~$0$, we say that the Hermitian duality $\cS:\cT\to\cT^*$ is a \emph{polarization} of $\cT$ if (c) below is satisfied:

{\def\theenumi{\alph{enumi}}
\begin{enumerate}\setcounter{enumi}{2}
\item\label{cond:twistorc}
The Hermitian pairing $h$ defines on $H$ a \emph{positive definite Hermitian form}.
\end{enumerate}}

\subsubsection*{Reduction to weight~$0$}
If $\cT$ has weight~$w$, we say that a Hermitian duality $\cS:\cT\to\cT^*(-w)$ is a \emph{polarization} of the twistor structure $\cT$ of pure weight~$w$ if $\cS(w/2)$ (defined after Definition \ref{def:sesqdual}) is a polarization of $\cT(w/2)$, \ie the positivity condition \eqref{cond:twistorc} above is satisfied for the Hermitian duality $\cS(w/2)$ of the twistor structure $\cT(w/2)$ of weight~$0$. This is equivalent to saying that, if $H$ is defined as above with $\cT(w/2)$ and $H\subset\cHS$ is the corresponding inclusion, then the restriction of $h_\bS$ defined by \eqref{eq:hermit} to $H\otimes_\CC\ov H$ is positive definite.

Notice that, if $(\cT,\cS)$ is a polarized twistor structure of weight~$w$, then, for any $k\in\hZZ$, $(\cT(k),\cS(k))$ is a polarized twistor structure of weight~$w-2k$.

\medskip
Under the equivalence above, the category of polarized twistor structures of weight~$w$ (the morphisms being the morphisms of twistor structures) is equivalent to the category of $\CC$-vector spaces with a positive definite Hermitian form (the morphisms being all linear maps). In particular we have:

\begin{fact}\label{fact:pol}
Let $\cT_1$ be a subtwistor structure of the polarized twistor structure $(\cT,\cS)$ (the weight is fixed). Then $\cS$ induces a polarization on the subtwistor, which is a direct summand of $(\cT,\cS)$.\qed
\end{fact}

\Subsection{Complex Hodge structures and twistor structures (after \cite{Simpson97})}
\label{subsec:chsmtw0}
Consider the example of \T\ref{subsub:c} with $X$ a point. So, let $H$ be a $\CC$-vector space equipped with a decomposition $H=\oplus_pH^{p,w-p}$, that we call a \emph{complex Hodge structure}. Consider the two decreasing filtrations
\[
F^{\prime p}=\ooplus_{p'\geq p}H^{p',w-p'}\qqbox{and} F^{\prime\prime q} =\ooplus_{q'\geq q}H^{w-q',q'}
\]
and the Rees modules associated to these filtrations
\[
\cH^{\prime \vee}\defin \ooplus_p\ov F^{\prime p}\hb^{-p}\subset\CC[\hb,\hbm]\ootimes_\CC \ov H\qquad
\cH''\defin \ooplus_q F^{\prime\prime q}\hb^{-q} \subset\CC[\hb,\hbm]\ootimes_\CC H.
\]
We will now work with the algebraic variant of twistor structures, where we replace $\cO(\bS)$ with $\CC[\hb,\hbm]$, $\cH',\cH''$ are viewed as free $\CC[\hb]$-modules, and where we replace $\cHS'$ with $\ccH^{\prime\Cir}=\CC[\hb,\hbm]\otimes_{\CC[\hb]}\cH'$, \etc.
We have
\begin{align*}
\cH'&=\ooplus_p\big((\ov H^{p,w-p})^\vee\hb^{p}\CC[\hb]\big)&
\cH^{\prime \vee}&=\ooplus_p\big(\ov H^{p,w-p}\hb^{-p}\CC[\hb]\big)\\
\cH''&=\ooplus_p\big(H^{p,w-p}\hb^{p-w}\CC[\hb]\big)&
\ov{\cH''}&=\ooplus_p\big(\ov H^{p,w-p}\hb^{w-p}\CC[\hbm]\big).
\end{align*}
The inclusions
\[
\cH^{\prime\vee}=\ooplus_p\ov F^{\prime p}\hb^{-p}\subset\CC[\hb,\hbm]\ootimes_\CC \ov H\supset \ooplus_q\ov F^{\prime\prime q}\hb^q=\ov{\cH''}
\]
define a semistable vector bundle of weight~$w$ on $\PP^1$. The pairing $C$ is induced by the natural $\CC$-duality pairing $\Crochet$:
\begin{align*}
(\ov H^{p,w-p})^\vee\hb^{p}\otimes \ov H^{p,w-p}\hb^{w-p}&\to\hb^w\CC[\hb]\\
x^\vee\hb^p\otimes y\hb^{w-p}&\mto\langle x^\vee,y\rangle\hb^w.
\end{align*}

Let us now compare the notion of polarization with the usual one. The natural inclusion $j:H\hto\ccH^{\prime\prime\Cir}$ is induced by $H^{p,w-p}\mto H^{p,w-p}\hb^{p-w}$. The conjugate inclusion $\ov\jmath:\ov H\hto\ov{\ccH^{\prime\prime\Cir}}$ is given by $\ov H^{p,w-p}\mto \ov H^{p,w-p}(-1)^{p-w}\hb^{w-p}$.

Let $k$ be a polarization of the Hodge structure as in \T\ref{subsub:c} (with $\dim X=0$). We view $k$ as an isomorphism $H^{p,w-p}\isom(\ov H^{p,w-p})^\vee$. Define $S':\cH''\to\cH'$ by $S'=\hb^w k:H^{p,w-p}\hb^{p-w}\to (\ov H^{p,w-p})^\vee\hb^{p}$ and $S''=(-1)^wS'=(-\hb)^wk$.

For $x,y\in H^{p,w-p}$, we have
\begin{align*}
h_\bS(j(x),\ov \jmath(\ov y))&=\thb^{-w}C\big(S''({j(x)}),\ov \jmath(\ov y)\big)\\
&=\thb^{-w}C\big(S''(x\hb^{p-w}),(-1)^{p-w}\ov y\hb^{w-p}\big)\\
&=(-1)^p\thb^{-w}C\big(k(x)\hb^p,\ov y\hb^{w-p}\big)\\
&=(-1)^pi^{-w}\langle k(x),\ov y\rangle.
\end{align*}
We therefore recover the usual notion of positivity.

\subsection{Graded Lefschetz twistor structures}\label{subsec:grlefschetz}
Fix $\varepsilon =\pm1$. A \emph{graded Lefschetz twistor structure} $[\cT=(\cH',\cH'',C),\cL]$ of weight~$w\in\ZZ$ and type $\varepsilon$ consists of
\begin{itemize}
\item
a finite family, indexed by $j\in\ZZ$, of twistor structures $\cT_j=(\cH'_{-j},\cH''_j,C_j)$ of weight~$w-\varepsilon j$, so that $\cT=\oplus_j\cT_j$,
\item
a graded morphism $\cL:\cT\to\cT(\varepsilon)$ of degree $-2$, such that, for any $j$, the component $\cL:\cT_j\to\cT_{j-2}(\varepsilon)$ is a morphism of twistor structures of weight~$w-\varepsilon j$, and that, for any $j\geq0$, $\cL^j:\cT_j\to\cT_{-j}(\varepsilon j)$ is an isomorphism.
\end{itemize}
Morphisms $\varphi$ are families $(\varphi_j)_j$ of morphisms of twistor structures, which are compatible with $\cL$. Remark that $\cL=(L',L'')$, with $L':\cH'_j\to\cH'_{j-2}$ and $L'':\cH''_j\to\cH''_{j-2}$ (we forget the index $j$ for $\cL$).

\begin{remarque}\label{rem:dessusdessous}
We will follow the usual convention when writing upper indices: put $\cT^j=\cT_{-j}$; then $\cT^j$ has weight~$w+\varepsilon j$ and $\cL$ is a morphism $\cT^j\to\cT^{j+2}(\varepsilon)$.
\end{remarque}

\medskip
The \emph{primitive part} of index $j\geq0$ of $(\cT,\cL)$ is the twistor structure $P\cT_j$ of weight~$w-\varepsilon j$ defined as $\ker \cL^{j+1}:\cT_j\to\cT_{-j-2}(\varepsilon(j+1))$. We therefore have
\[
P\cT_j=(L^{\prime j}P\cH'_j,P\cH''_j,C_j),
\]
where $P\cH''_j=\ker [L^{\prime\prime j+1}:\cH''_j\to\cH''_{-j-2}]$, and similarly for $P\cH'_j$. The primitive parts allow one to reconstruct $(\cT,\cL)$ up to isomorphism using the \emph{Lefschetz decomposition}: for any $j\geq0$,
\begin{equation}\label{eq:Lefdec}
\cT_j=\ooplus_{k\geq0}\cL^kP\cT_{j+2k}(-\varepsilon k),\quad \cT_{-j}=\ooplus_{k\geq0}\cL^{k+j}P\cT_{j+2k}(-\varepsilon (k+j)).
\end{equation}

Put $\cT'=\cT^*(-w)$, with the grading $\cT'_j=\cT_{-j}^*(-w)$ and define $\cL'=\cL^*$. Then $(\cT',\cL')$ is an object of the same kind as $(\cT,\cL)$. Moreover we have
\[
P\cT'_j=(\cL^jP\cT_j)^*(-w)\xrightarrow[\ts\sim]{~\ts\cL^{*j}~} (P\cT_j)^*(-w+\varepsilon j).
\]

\subsubsection*{Weil reduction to weight~$0$}
Any graded Lefschetz twistor structure $(\cT,\cL)$ of weight~$w$ and type $\varepsilon=\pm1$ gives rise to a graded twistor structure $\wt\cT$ of weight~$0$ and a graded morphism $\wt\cL:\wt\cT_j\to\wt\cT_{j-2}$ by putting $\wt\cT_j=\cT_j((w-\varepsilon j)/2)$ and $\wt \cL=\cL$. We hence have
\[
\wt\cH'_{-j}=\cH'_{-j},\quad\wt\cH''_j=\cH''_j,\quad\wt C_j=(iz)^{-w+\varepsilon j}C_j,\quad \wt L'=L',\quad \wt L''=L''.
\]
The category of graded Lefschetz twistor structures of weight~$0$ is equivalent to the category of graded $\CC$-vector spaces $H=\oplus_j H_j$ of finite dimension, equipped with a graded nilpotent endomorphism $L:H\to H$ of degree $-2$, such that, for any $j\geq0$, $L^j:H_j\to H_{-j}$ is an isomorphism. More precisely, the twistor condition on $\cT_j$ gives vector spaces $H'_j\subset \cH'_j$, $H''_j\subset\cH''_j$ such that $C_j:H'_{-j}\otimes\ov{H''_j}\to z^{w-\varepsilon j}\CC$ is nondegenerate. Put $H_j=H''_j$ and use $\wt C_j$ to identify $H'_{-j}$ with $\ov{H_j^*}$. Define $L:H_j\to H_{j-2}$
as the restriction of $L''$ to $H''_j$.

We call this situation the case of weight~$0$ and type $0$ (this is not exactly obtained by putting $w=0$ and $\varepsilon=0$ in the previous case).

From this equivalence, it is clear that the category of graded Lefschetz twistor structures of weight~$w$ and type $\varepsilon$ is abelian and that any morphism is graded with respect to the Lefschetz decomposition \eqref{eq:Lefdec}.

\subsubsection*{Polarization}
A \emph{Hermitian duality} of $\cT$ is a graded isomorphism $\cS:\cT\to\cT^*(-w)$ which is Hermitian in the sense of Definition \ref{def:hermw}, in other words $(-1)^w$-selfadjoint, \ie a family $(\cS_j)_{j\in\ZZ}$ of isomorphisms
\[
\cS_j=(S'_{-j},S''_j):\cT_j\to\cT_{-j}^*(-w),
\]
satisfying $\cS_j^*=(-1)^w\cS_{-j}$, in other words a family indexed by $j\in\ZZ$ of isomorphisms
\[
S'_j:\cH''_{j}\isom\cH'_{j},\qquad
S''_j:\cH''_j\isom\cH'_{j},
\]
with $S'_j=(-1)^wS''_{j}$, such that, for $x\in\cH''_{-j}$ and $y\in\cH''_j$, we have
\[
C_j(S'_{-j}x,\ov y)=\thb^{2w}C_{-j}^*(x,\ov{S''_jy})\defin \thb^{2w}\ov{C_{-j}(S''_jy,\ov x)}.
\]
We say that this Hermitian duality is compatible with $\cL$ if $\cL$ is skewadjoint with respect to $\cS$, \ie
\begin{equation}\label{eq:skew}
\cL^*\circ\cS_j+\cS_{j-2}\circ \cL=0.
\end{equation}
This implies that $\cL^{*j}\circ\cS_j=(-1)^j\cS_{-j}\circ \cL^j$, which can be written, using the symmetry of $\cS$, as
\[
(\cS_{-j}\circ\cL^j)^*=(-1)^{w-\epsilon j}(\cS_{-j}\circ\cL^j)\quad (\epsilon=\pm1).
\]
This also implies that $\cS$ is completely determined by its restriction to the primitive parts $P\cT_j$ and can be rebuild using the Lefschetz decomposition.

Given a Hermitian duality $\cS$ of $(\cT,\cL)$ (that is, a Hermitian duality of $\cT$ compatible with $\cL$), the composed morphism ($j\geq0$)
\begin{equation}\label{eq:grlefhermdual}
(P\cS)_j\defin\cS_{-j}\circ \cL^j: P\cT_j\To{\cL^j}\cL^jP\cT_j\To{\cS_{-j}}(P\cT_j)^*(-w+\varepsilon j)
\end{equation}
is thus a Hermitian duality of $P\cT_j$. We then say that $\cS$ is a \emph{polarization} of $(\cT,\cL)$ if, for any $j\geq0$, this Hermitian duality $(P\cS)_j$ is a polarization of the twistor structure $P\cT_j$ (of weight $(w-\epsilon j)$) as defined in \T\ref{subsec:hermdualpol}.

\begin{exemple}\label{ex:lefpol}
Assume that $w=0$ and that $\cH'_j=\cH''_j=\cH_j$ and $C^*_{-j}=C_j$ for all $j$, so that $\cT^*_{-j}=\cT_j$. Assume also that $\cS_j=(\id,\id)$ for all $j$. The fact that $\cL$ is skewadjoint with respect to $\cS$ means that $\cL^*=-\cL$, \ie $L'=-L''$. The polarization $(P\cS)_j$ on $P\cT_j$ is the morphism $P\cT_j\to(P\cT_j)^*(\varepsilon j)$ given by
\[
(L^{\prime j},L^{\prime\prime j}):\big(L^{\prime j}P\cH_j,P\cH_j,C_j\big)\to\big(P\cH_j,L^{\prime\prime j}P\cH_j,\thb^{-2\varepsilon j}C_{-j}\big).
\]
The positivity condition is that $\Gamma(\bS,\cbbullet)$ of the Hermitian form \eqref{eq:hermit}
\[
\thb^{\varepsilon j}C_j(L^{\prime\prime j}\cbbullet,\ov\cbbullet): P\cH_{j|\bS}\otimes_{\cO_\bS} \ov{P\cH_{j|\bS}}\to\cO_\bS
\]
takes values in $\CC$ and is positive definite, when restricted to $PH_j\subset \Gamma(\bS,P\cH_{j|\bS})$.

In other words, $P\cT_j(-\varepsilon j/2)$ with its polarization is isomorphic to the twistor structure $\big(P\cH_j,P\cH_j,\thb^{\varepsilon j}C_j(L^{\prime\prime j}\cbbullet,\ov\cbbullet)\big)$ with polarization $(\id,\id)$.
\end{exemple}

The datum of a polarized graded Lefschetz twistor structure $(\cT,\cL,\cS)$ of weight~$w$ and type $\varepsilon$ is therefore equivalent to the datum of one of weight~$0$ and type $0$ plus that of $w,\varepsilon$, \ie to the datum of $(H,L,w,\varepsilon,h)$, where $(H,L,w,\varepsilon)$ is as above and $h$ is a homogeneous sesquilinear form of degree $0$ on $H\otimes\ov H$, \ie corresponds to a family of pairings
\[
h:H_{-j}\ootimes_\CC\ov H_j\to\CC
\]
such that $L$ is skew-adjoint with respect to $h$ and $h(L^j\cbbullet,\ov\cbbullet)$ is a positive definite Hermitian form on $PH_j$ for any $j\geq0$.

\begin{remarque}\label{rem:W}
The datum of $(H,h,L)$ as above is equivalent to the datum of a finite dimensional Hermitian $\CC$-vector space $(H,h)$ with a $\SL_2(\RR)$-action: the torus $\CC^*$-action gives the grading, and $L$ comes from the corresponding $\sld(\RR)$-action. The positivity condition for $h(L^j\cbbullet,\ov\cbbullet)$ on $PH_j$ is then equivalent to the positivity of the form $h(W\cbbullet,\ov \cbbullet)$ on $H$, where $W$ corresponds to $\big(\begin{smallmatrix}0&1\\-1&0\end{smallmatrix}\big)$.
\end{remarque}

\begin{remarque}\label{rem:factpol}
The Fact \ref{fact:pol} also applies to graded Lefschetz twistor structures of weight~$w$ and type $\varepsilon$. Indeed, reduce first to weight~$0$ and type $0$. According to Remark \ref{rem:W}, it is then a matter of proving that, given $(H,h)$ with an action of $\SL_2(\RR)$ such that $k(\cbbullet,\ov \cbbullet)\defin h(W\cbbullet,\ov \cbbullet)$ is positive definite, any $\SL_2(\RR)$-stable subspace $H'$ has a $\SL_2(\RR)$-stable $k$-orthocomplement: this is clear.
\end{remarque}

\begin{remarque}\label{rem:grlefrtriples}
If one forgets the notion of weight and the notion of positivity, one may define the category of \emph{graded Lefschetz $\RTriples$} and the notion of Hermitian duality on objects of this category by changing above the words ``twistor structure'' with the words ``object of $\RTriples$''.
\end{remarque}

\begin{remarque}[Stability by extension]\label{rem:extleftriples}
Assume that we have an exact sequence $0\to\cT'\to\cT\to\cT''\to0$ in the category of graded $\RTriples$ (morphisms are graded of degree $0$); assume that each of these objects is equipped with a graded morphism $\cL$ of degree $-2$ and type $\epsilon$, in a compatible way with the exact sequence; lastly, assume that $(\cT',\cL)$ and $(\cT'',\cL)$ are graded Lefschetz twistor structures of the same weight~$w$. Then so is $(\cT,\cL)$. Indeed, applying the Weil reduction procedure above to all objects, one can assume from the beginning that $\epsilon=0$ and $w=0$. We have exact sequences $0\to\wt\cT'_j\to\wt\cT_j\to\wt\cT''_j\to0$ for any $j$. The locally free $\cO_{\PP^1}$-module corresponding to $\cT_j$ is an extension of two free $\cO_{\PP^1}$-modules, hence is free. We are now reduce to consider exact sequences of graded vector spaces $0\to H'\to H\to H''\to0$ with compatible endomorphisms $L$, such that, for any $j\geq1$, $L_j:H'_j\to H'_{-j}$ and $L_j:H''_j\to H''_{-j}$ are isomorphisms. Then, $L_j:H_j\to H_{-j}$ is an isomorphism.
\end{remarque}

\subsection{Two results on polarized graded Lefschetz twistor structures}
Let $(\cT,\cL,\cS)$ and $(\cT',\cL',\cS')$ be polarized graded Lefschetz twistor structures of type $\varepsilon=\pm1$ and weight~$w$ and $w-\varepsilon$ respectively. Let $c,v$ be graded morphisms of twistors of degree $-1$
\[
c:\cT_{j+1}\to\cT'_j,\qquad v:\cT'_j\to\cT_{j-1}(\varepsilon)
\]
such that $\cL=v\circ c$ and $\cL'=c\circ v$. Assume that $c,v$ are adjoint with respect to $\cS$ and $\cS'$, \ie for any $j$, the following diagram commutes:
\[
\xymatrix{
\cT_j\ar[r]^-{\cS_j}\ar[d]_-{c}&\cT^*_{-j}(-w)\ar[d]^-{v^*}\\
\cT'_{j-1}\ar[r]_-{\cS'_{j-1}}&\cT^{\prime*}_{-j+1}(-w+\varepsilon)
}
\]
and the adjoint diagram anticommutes.

\begin{proposition}[\cf {\cite[lemme 5.2.15]{MSaito86}}]\label{prop:declef}
Under these assumptions, we have a decomposition $\cT'=\im c\oplus\ker v$ as a graded Lefschetz twistor structure.
\end{proposition}

\begin{proof}
Remark first that each term of the decomposition is stable under $\cL'$. Apply the Weil reduction to weight~$0$. Now, $\wt\cT,\wt\cT'$ are graded vector spaces with a nilpotent endomorphism $\wt\cL,\wt\cL'$ of degree $-2$, and with sesquilinear forms $h,h'$ of degree $0$, such that $\wt\cL,\wt\cL'$ are skewadjoint. There are morphisms $\wt c,\wt v$ of degree $-1$ such that $\wt v\circ \wt c=\wt\cL$, $\wt c\circ \wt v=\wt\cL'$, and which are adjoint or skewadjoint to each other, as above. The proof of \cite[lemme 5.2.15]{MSaito86} applies to this case.
\end{proof}

Fix $\varepsilon_1,\varepsilon_2=\pm1$. The notion of a (polarized) bigraded Lefschetz twistor structure $(\cT,\cL_1,\cL_2)$ of weight~$w$ and bi-type $(\varepsilon_1,\varepsilon_2)$ is defined in a natural way, similarly to the single graded case: $\cL_1$ and $\cL_2$ should commute and the \emph{primitive part} in $\cT_{j_1,j_2}$ is by definition the intersection of $\ker\cL_1^{j_1+1}$ and $\ker\cL_2^{j_2+1}$. The following lemma will be useful in \T\ref{sec:n1n0}.

\begin{lemme}\label{lem:n1n0}
Let $(\cT,\cL_1,\cL_2,\cS)$ be a polarized bigraded Lefschetz twistor structure of weight~$w$ and bi-type $(\varepsilon,\varepsilon)$. Put on $\cT$ the grading $\cT_\ell=\oplus_{j+k=\ell}\cT_{j,k}$ and set $\cL=\cL_1+\cL_2$. Then $(\cT,\cL_1+\cL_2,\cS)$ is a polarized graded Lefschetz twistor structure of weight~$w$ and type $\varepsilon$.
\end{lemme}

\begin{proof}
Reduce to weight~$0$ and bi-type $(0,0)$. We therefore have a Hermitian $\CC$-vector space $(H,h)$ equipped with a $\SL_2(\RR)\times\SL_2(\RR)$ action and a positivity condition (\cf Remark \ref{rem:W}). Consider the diagonal $\SL_2(\RR)$-action. Then $W$ acts by $W_2=(W,W)$. Now, the positivity of $h(W\cbbullet,\ov\cbbullet)$ follows from \cite[\T4.3]{G-N90}.
\end{proof}

A \emph{differential} $d$ on $(\cT,\cL_1,\cL_2,\cS)$ is a morphism of bidegree $(-1,-1)$
\[
d:\cT_{j_1,j_2}\to\cT_{j_1-1,j_2-1}(\varepsilon_1+\varepsilon_2)
\]
such that $d\circ d=0$, which commutes with $\cL_1,\cL_2$ and is selfadjoint with respect to~$\cS$.

\begin{proposition}[M\ptbl Saito, P\ptbl Deligne, \cf {\cite[proposition 4.2.2]{MSaito86}, \cite[th\'eor\`eme 4.5]{G-N90}}]\label{prop:lefcoh}
In such a situation, the cohomology $\ker d/\im d$, with the induced $\cL_1,\cL_2,\cS$, is a polarized bigraded Lefschetz twistor structure $(\cT,\cL_1,\cL_2,\cS)$ of weight~$w$ and type $\varepsilon_1,\varepsilon_2$.
\end{proposition}

\begin{proof}
Notice first that $\cS$ induces a Hermitian duality on $\ker d/\im d$. Notice also that the Weil reduction to weight~$0$ commutes with taking cohomology. It is therefore enough to prove the proposition in the case of a bigraded Lefschetz twistor structure of weight~$0$ and bitype $0,0$ defined as above for the single graded case. This is done in \loccit
\end{proof}

\section{Smooth twistor structures in arbitrary dimension} \label{sec:smtwqc}
\subsection{}\label{subsec:smtwqc}
A \emph{smooth twistor structure} $\cT=(\index{$hcurlprime$@$\cH'$}\cH',\cH'',C)$ (or a \emph{variation of twistor structure}) weight~$w$ on $\cX$ is a smooth object of $\RTriples(X)$ in the sense of \T\ref{subsec:smoothtriple} such that its restriction to each $x_o\in X$ (in the sense of Definition \ref{def:restrsmooth}) is a twistor structure of pure weight~$w$. The \emph{rank} of $\cT$ is the rank of the bundles $\cH',\cH''$.

A \emph{polarization} is a Hermitian pairing $\cS:\cT\to\cT^*(-w)$ of weight~$w$ (in the sense of Definition \ref{def:hermw}) which induces a polarization by restriction to any $x_o\in X$ (in the sense of Definition \ref{def:restrsesqui}).

\enlargethispage{1.5\baselineskip}%
\begin{remarque}\label{rem:sans pole}
We have seen that the sesquilinear pairing $C$ takes values in $\cCh{\cX}$, according to Lemma \ref{lem:smoothsesqui}. So the restriction to $x_o$ of each component of the smooth twistor structure is well defined. It is also nondegenerate and gives a gluing of $\cH^{\prime*}$ with $\ov{\cH''}$, defining thus a $\cC^{\infty,\an}_{X\times\PP^1}$-bundle $\wt\cH$ on $X\times \PP^1$.
\end{remarque}

\begin{lemme}[\cite{Simpson97}]\label{lem:twharm}
The datum of a smooth polarized twistor structure of weight~$0$ on~$X$ is equivalent to the datum of a flat holomorphic bundle $(V,\nabla)$ on $X$ with a harmonic metric $h$, or the datum of a holomorphic Higgs bundle $(E,\theta'_E)$ with a harmonic metric~$h$.
\end{lemme}

\begin{proof}
Given $\index{$hdv$@$(H,D_V)$}(H,D'_V,D''_V,h)$ as in \T\ref{subsub:a}, consider the associated operators $\index{$dE$@$D_E$}D'_E$, $D''_E$, $\theta'_E$, $\theta''_E$. Recall also that $\cH=\cC^{\infty,\an}_\cX\otimes_{\cC_X^\infty}H$ is equipped with connections $\index{$dHcurl$@$D_\cH$}D'_\cH$ and $D''_\cH$ as in \eqref{eq:dbar} and \eqref{eq:dH}. Let $\cH'\subset\cH$ be defined as the kernel of $D''_\cH$. As we canonically have $\cC_X^\infty=\ov{\cC_X^\infty}$, we may identify the locally free $\cC_X^\infty$-module $H$ with its conjugate and view $h$ as a $\cC_X^\infty$-linear morphism $H\otimes_{\cC_X^\infty}\ov H\to\cC_X^\infty$. Consider on $\ov H$ the operators
\[
D'_{\ov E}\defin \ov{D''_E},\quad D''_{\ov E}\defin \ov{D'_E},\quad \theta'_{\ov E}\defin \ov{\theta''_E},\quad \theta''_{\ov E}\defin \ov{\theta'_E}.
\]
For local sections $u,v$ of $H$, we therefore have
\begin{align*}
d'h(u,\ov v)&=h(D'_Eu,\ov v)+h(u,D'_{\ov E}\ov v)\\
h(\theta'_Eu,\ov v)&=h(u,\theta'_{\ov E}\ov v)
\end{align*}
and the $(0,1)$ analogues. Extend $h$ as
\[
\index{$hsbold$@$h_\bS$}h_\bS:\cHS\ootimes_{\cCh{\cX}}\ov{\cHS}\to\cCh{\cX}
\]
by $\cCh{\cX}$-linearity. If we define as above $D'_{\ov\cH}=\ov{D''_\cH}$ and $D''_{\ov\cH}=\ov{D'_\cH}$, the previous relations may be written in a more convenient way:
\begin{align*}
d'h_\bS(u,\ov v)&=h_\bS(D'_\cH u,\ov v)+h_\bS(u,D'_{\ov\cH}\ov v)\\
d''h_\bS(u,\ov v)&=h_\bS(D''_\cH u,\ov v)+h_\bS(u,D''_{\ov\cH}\ov v).
\end{align*}
Define $\cH''=\cH'$ and $S''=\id$, $S'=\id$. Let $C$ be the restriction of $h_\bS$ to $\cHS'\otimes_{\cO_\bS}\ov{\cHS'}$. The $\cO_{(X,\ov X),\bS}$-linearity of $C$ is clear from the analogous property of $h_\bS$. Let us verify the $\cR_{(X,\ov X),\bS}$-linearity: for a local section $v$ of $\cHS'$, we have $D''_\cH v=0$, hence \hbox{$D'_{\ov \cH}\ov v=0$}; therefore, given $\hb_o\in\bS$, for local sections $u$ of $\cH'_{\hb_o}$ and $v$ of $\cH'_{-\hb_o}$, we have (using the standard $d'$ operator on functions),
\[
d' C(u,\ov v)=h_\bS(D'_\cH u,\ov v)+ h_\bS(u,D'_{\ov \cH}\ov v)=h_\bS(D'_\cH u,\ov v)=C(D'_\cH u,\ov v).
\]
The $d''$-linearity is obtained similarly, exchanging the roles of $u$ and $v$.

Now we will show that $(\cH',\cH',C)$ is a smooth polarized twistor structure of weight~$0$ on $X$. Denote by $\index{$hcurltilde$@$\wt\cH$}\wt\cH$ the bundle on $X\times\PP^1$ obtained by gluing the dual $\cH^\vee$ (chart $\Omega_0$, coordinate $\hb$) with the conjugate $\ov\cH$ (chart $\Omega_\infty$, coordinate $\hb'$) using the isomorphism $h_\bS:\ov{\cHS}\isom\cHS^\vee$. As $h_\bS$ induces the isomorphism $h:\ov H\isom H^\vee$ (using the natural inclusions $H^\vee=1\otimes H^\vee\subset \pi_*\cHS^\vee$ and $\ov H=1\otimes \ov H\subset \pi_*\ov{\cHS}$), the natural map $\cC^{\infty,\an}_{X\times\PP^1}\otimes_{\cC^\infty_X}\ov H\to\wt\cH$ is an isomorphism and $\pi_*\wt h$ (\cf \eqref{eq:ccP}) is identified with $h$. The restriction of these objects to each $x_o\in X$ gives therefore a polarized twistor structure of weight~$0$. Now, as we have $\cH=\cC^{\infty,\an}_\cX\otimes_{\cO_\cX}\cH'$, the restriction of $\cH'$ to $x_o\in X$ is equal to that of $\cH$, and this shows that $(\cH',\cH', C,\cS)$ with $\cS=(\id,\id)$, is a smooth polarized twistor structure of weight~$0$ on $X$.

Last, let us show that, by restriction to $\hb=1$, one recovers $(V,\nabla)$. We know that $D^{\prime\prime2}_\cH=0$, so the Dolbeault complex $\big(\cH\otimes_{\cC_\cX^\infty}\cE_\cX^{(0,\bbullet)}\big)$ is a resolution of $\cH'$. As $\cH'$ and the terms of this complex are $\cO_{\Omega_0}$-locally free, the restriction to $\hb=1$ of this complex is a resolution of $\cH'/(\hb-1)\cH'$. But we clearly have $\cH/(\hb-1)\cH=H$ and ${D''_\cH}_{|\hb=1}=D''_V$, so $\cH'/(\hb-1)\cH'=\ker D''_V=V$. Conclude by noticing that the restriction of $D'_\cH$ to $\hb=1$ is $D'_V$.

\medskip
Conversely, let $(\cH',\cH'',C,\cS)$ be a polarized twistor structure of weight~$0$. We will assume that $\cH''=\cH'$ and $\cS=(\id,\id)$ (it is not difficult to reduce to this case, \cf Remark \ref{rem:hermdual}). Put $\cH=\cC^{\infty,\an}_\cX\otimes_{\cO_\cX}\cH'$ and denote by $h_\bS:\cHS\otimes_{\cCh{\cX}}\ov{\cHS}\to\cCh{\cX}$ the $\cCh{\cX}$-linear morphism induced by $C$. As it is nondegenerate, we may use it to glue $\cH^\vee$ (on $X\times\Omega_0$) with $\ov\cH$ (on $X\times\Omega_\infty$), and obtain a $\cC^{\infty,\an}_{X\times\PP^1}$-linear bundle $\wt\cH$.

As $(\cH',\cH',C)$ has weight~$0$, the restriction of $\wt\cH$ to any $x_o\in X$ is the trivial bundle on $\PP^1$, thus the natural morphism $\pi^*\pi_*\wt\cH\to\wt\cH$ is an isomorphism and $\ov H\defin\pi_*\wt\cH$ is a $\cC^\infty_X$ locally free sheaf such that $\wt\cH=\cC^{\infty,\an}_{X\times\PP^1}\otimes_{\cC^\infty_X}\ov H$. We define the metric $h$ on $H$ as $\pi_*\wt h$, where $\wt h$ is constructed as in \T\ref{subsec:hermdualpol}. As there is a natural inclusion $\ov H\subset \pi_*\ov{\cHS}$ and a conjugate inclusion $H\subset\pi_*\cHS$ (induced by the natural restriction morphism $\pi_*\to\pi_{|\bS*}$), the metric $h$ may also be defined as the restriction of $\Gamma(\bS,h_{\bS})$ to $H\otimes_{\cC^{\infty}_{X}}\ov H$ (\cf Remark \ref{rem:h}).

The bundle $\wt\cH$ comes equipped with a $(1,0)$-connection relative to $\pi:X\times\PP^1\to\PP^1$, with poles of order at most one along $X\times\{0\}$:
\[
D'_{\wt\cH}=\wt D': \wt\cH\to\Omega^1_{X\times\PP^1/\PP^1}(X\times\{0\}) \ootimes_{\cO_{X\times\PP^1}}\wt\cH.
\]
Indeed, in the chart $\Omega_0$, the $\cR_\cX$-structure on $\cH'$ defines such a $(1,0)$-connection on $\cH$ (\cf \T\ref{num:RX}), hence in a natural way on $\cH^\vee$. Put the trivial $\sigma^*\cR_\cX$-structure on $\ov{\cH'}=\sigma^*c\cH'$ (a bundle which is purely antiholomorphic with respect to $X$ because of $c$). Therefore, $\wt\cH_{|X\times\Omega_\infty}\defin\cC^{\infty,\an}_{X\times\Omega_\infty}\otimes_{\cO_{\ov X\times\Omega_\infty}}\ov{\cH'}$ has a natural connection of type $(1,0)$ induced by $d'$ on $\cC^{\infty,\an}_{X\times\Omega_\infty}$.

Let us verify that both $\hb$-connections on $\cHS^\vee$ and $\ov{\cHS}$ correspond each other under the gluing $h_\bS$. For $m,\mu$ local sections of $\cHS$ and $\varphi,\psi$ local sections of $\cCh{\cX}$, we have
\[
d'h_\bS(\varphi m,\psi\ov \mu)=d'\big(\varphi\psi C(m,\ov \mu)\big)=h_\bS\big(\wt D'(\varphi m),\psi\ov \mu)+h_\bS\big(\varphi m,\wt D'(\psi\ov \mu)\big),
\]
as $C$ is $\cR_{(X,\ov X),\bS}$-linear and $\wt D'(\psi\ov \mu)=d'\psi\otimes\ov\mu$.

A similar definition and construction can be done for the conjugate notion, namely a relative connection of type $(0,1)$, with poles along $X\times\{\infty\}$. We denote it by $\wt D''$.

Composing $\wt D'$ with the residue along $X\times\{0\}$, we get an endomorphism of $\ov H$, that we denote by $-\theta'_{\ov E}$. According to the relative triviality of $\wt\cH$, we may write $\wt D'$ as
\[
\wt D'=D'_{\ov E}-\hbm\theta'_{\ov E},
\]
where $D'_{\ov E}$ is a $(1,0)$-connection on $\ov H=\pi_*\wt\cH$. Similarly, write $\wt D''=D''_{\ov E}+\ov\hbm\theta''_{\ov E}=D''_{\ov E}-\hb\theta''_{\ov E}$. Define also on $\ov{\wt\cH}$
\begin{align*}
D'_{\ov{\wt\cH}}&=\sigma^*c\big(\wt D''\big)=D'_E+\hbm\theta'_E=D'_E-\ov\hb\theta'_E,\\
D''_{\ov{\wt\cH}}&=\sigma^*c\big(\wt D'\big)=D''_E-\ov\hbm\theta''_E=D''_E+\hb\theta''_E.
\end{align*}
The $\cR_{(X,\ov X),\bS}$-linearity of $C$ implies that $h_\bS$ is compatible with $D'_{\wt\cH},D'_{\ov{\wt\cH}}$ on the one hand, and with $D''_{\wt\cH},D''_{\ov{\wt\cH}}$ on the other hand, \ie satisfies, for local sections $u,v$ of $H\subset\Gamma(\bS,\cHS)$,
\begin{align*}
d'h_\bS(u,\ov v)&=h_\bS\big((D'_E+\hbm\theta'_E)u,\ov v\big)+ h_\bS\big(u,(D'_{\ov E}-\hbm\theta'_{\ov E})\ov v\big),\\
d''h_\bS(u,\ov v)&=h_\bS\big((D''_E+\hb\theta''_E)u,\ov v\big)+ h_\bS\big(u,(D''_{\ov E}-\hb\theta''_{\ov E})\ov v\big).
\end{align*}
From this and from flatness properties of the connections $\wt D',\dots$, which is a consequence of the existence of a $\cR_\cX$ or a $\cR_{\ov\cX}$-structure, we get all relations needed for the harmonicity of $h$.
\end{proof}

\begin{remarque}\label{rem:baseorth}
Keep notation as in the proof above. Let $\varepsilong$ be a basis of $\cH$ which is orthonormal for $h_\bS$. Then $\varepsilong$ is contained in $H$ (and therefore is an orthonormal basis for $h$). Indeed, it defines bases $\varepsilong^\vee$ and $\ov{\varepsilong}$ of $\cH^\vee$ and $\ov\cH$ respectively. The orthonormality property exactly means that these two bases coincide near $\bS$, hence define a basis of $\wt\cH$. Consequently, $\ov{\varepsilong}$ is contained in $\ov H$ and therefore $\varepsilong$ is contained in $H$.
\end{remarque}

\subsection{Hodge theory for smooth twistor structures}\label{subsec:Hodgetw}
Let $X$ be a compact K\"ahler manifold with K\"ahler form $\omega$, and let $(\cT,\cS)$ be a smooth polarized twistor structure on $X$ of weight~$w$, with $\cT=(\cH',\cH',C)$ and $\cS=(S',S'')=((-1)^w\id,\id)$. Denote by $f:X\to{\mathrm{pt}}$ the constant map.

\begin{H-Stheoreme}\label{th:HodgeSimpson}
The direct image $\big(\oplus_jf_\dag^j\cT,\cL_\omega\big)$ is a graded polarized Lefschetz twistor structure of weight~$w$ and type $\varepsilon=1$.
\end{H-Stheoreme}

We will be more explicit on the polarization later on. We will restrict to the case $w=0$. The general case follows easily by changing $C$ to $\thb^wC$ and $S'$ to $(-1)^wS'$.

\begin{proof}
Let us first recall the results of \cite[\T2]{Simpson92} concerning this point. Let $(H,D_V,h)$ be a harmonic bundle on $X$ as in \T\ref{subsub:a}, with associated operators $D'_E$, $D''_E$, $\theta'_E$ and $\theta''_E$. Put $\index{$dcurlcurl$@$\ccD_0$, $\ccD_\infty$, $\ccD_{\hb}$}\ccD_\infty=D'_E+\theta''_E$ and $\ccD_0=D''_E+\theta'_E$, so that $D_V=\ccD_\infty+\ccD_0$ (notice that $\ccD_\infty,\ccD_0$ are not of type $(1,0)$ or $(0,1)$). The main observation is that the K\"ahler identities
\[
\Delta_{D_V}=2\Delta_{\ccD_\infty}=2\Delta_{\ccD_0},
\]
are satisfied for the Laplacian, and that the Lefschetz operator $L=\omega\wedge{}$ commutes with these Laplacians. For any $\hb_o\in\CC$, let $\Delta_{\hb_o}$ be the Laplacian of $\ccD_{\hb_o}\defin\ccD_0+\hb_o\ccD_\infty$. Then, the K\"ahler identities proved in \loccit for $\ccD_0$ and $\ccD_\infty$ (which are denoted there respectively by $D''$ and $D'_{\mathbf{K}}$) imply that
\begin{equation}\label{eq:Deltao}
\Delta_{\hb_o}=(1+\module{\hb_o}^2)\Delta_{\ccD_0}.
\end{equation}

It follows that the spaces $\Harm^k_{\hb_o}(H)$ of $\Delta_{\hb_o}$-harmonic sections are independent of $\hb_o$ and that the harmonic sections are closed with respect to any $\ccD_{\hb_o}$. Moreover, $\Harm^{\cbbullet}_{\hb_o}(H)$ is equal to the cohomology of the complex $\Gamma\big(X,(\cE_X^{\cbbullet}\otimes_{\cC_X^\infty}H,\ccD_{\hb_o})\big)$.

The Lefschetz operator induces a $\sld$-structure on the space of harmonic sections.

Consider the complex
\[
\Big(\CC[\hb]\ootimes_\CC\cE_X^{\cbbullet}\ootimes_{\cC_X^\infty}H; \ccD_0+\hb\ccD_\infty\Big).
\]
The restriction of this complex to each $\hb_o\in\Omega_o$ gives the previous complex. We may ``rescale'' this complex by the isomorphism
\begin{equation}\label{eq:rescale}
\begin{split}
\CC[\hb]\ootimes_{\CC}\cE_X^{p,q}\ootimes_{\cC_X^\infty}H&\To{\index{$iota$@$\iota$ (rescaling)}\iota}\hb^{-p}\CC[\hb]\ootimes_{\CC}\cE_X^{p,q}\ootimes_{\cC_X^\infty}H\subset\CC[\hb,\hbm]\ootimes_{\CC}\cE_X^{p,q}\ootimes_{\cC_X^\infty}H\\
\eta^{p,q}\otimes m&\mto\hb^{-p}\eta^{p,q}\otimes m
\end{split}
\end{equation}
so that the differential $\ccD_0+\hb\ccD_\infty$ is changed into $D_{\cH}=D'_\cH+D''_\cH$, where $D'_\cH=D'_E+\hbm\theta'_E$ and $D''_\cH=D''_E+\hb\theta''_E$ are the restriction of $D'_{\ov{\wt\cH}}$ and $D''_{\ov{\wt\cH}}$ respectively to the chart $\Omega_0$. We have
\[
\cO_\cX\ootimes_{\cO_X[\hb]}\iota\Big(\CC[\hb]\ootimes_\CC\cE_X^{\cbbullet}\ootimes_{\cC_X^\infty}H; \ccD_0+\hb\ccD_\infty\Big)= \big(\cE^{n+\bbullet}_\cX\otimes_{\cC^{\infty,\an}_\cX}\cH,D_\cH\big).
\]
It follows that, for $\hb_o\neq0$, the restricted complex $\big(\cE^{n+\bbullet}_\cX\otimes_{\cC^{\infty,\an}_\cX}\cH,D_\cH\big)_{|\hb=\hb_o}$ is nothing but the de~Rham complex of the flat holomorphic bundle $(V_{\hb_o},D'_E+\hbm_o\theta'_E)$, with $V_{\hb_o}=\ker(D''_E+\hb_o\theta''_E)$ (\cf \T\ref{subsub:a}), the cohomology of which is the local system $\ccV_{\hb_o}=\ker \nabla_{\hb_o}$; and for $\hb_o=0$, the restricted complex is the Dolbeault complex $\big(\Omega^{\cbbullet}_X\otimes_{\cO_X}E,\theta'_E\big)$.

\medskip
We now come back to the proof of the theorem. Recall that, as $D''_\cH$ defines a complex structure on $\cH$, we have
\[
(\Omega^{n+\bbullet}_\cX\otimes_{\cO_\cX}\cH',D'_{\cH'})\isom(\cE^{n+\bbullet}_\cX\otimes_{\cC^{\infty,\an}_\cX}\cH,D_\cH).
\]
Therefore, after \eqref{eq:lefrigimdir}, we may compute $f_\dag\cT$ as follows:
\[
f^j_\dag\cT=\Big(\cH^{-j}f_*\big(\cE^{n+\bbullet}_\cX\otimes_{\cC^{\infty,\an}_\cX} \cH,D_\cH\big), \cH^{j}f_*\big(\cE^{n+\bbullet}_\cX\otimes_{\cC^{\infty,\an}_\cX}\cH, D_\cH), C^j\Big)
\]
and $C^j=f^j_\dag h_\bS$ is the natural sesquilinear pairing
\[
f_\dag^jh_\bS:H^{-j}\big(X;(\cE_{\cX|\bS}^{n+\bbullet}\otimes\cHS,D_\cH)\big)\ootimes_{\cO_\bS} \ov{H^j\big(X;(\cE_{\cX|\bS}^{n+\bbullet}\otimes\cHS,D_\cH)\big)} \to\cO_\bS
\]
induced by
\[
\big[\eta^{n-j}\otimes m\big]\otimes\ov{\big[\eta^{n+j}\otimes \mu\big]}\mto \frac{\varepsilon(n+j)}{(2i\pi)^n}\int_Xh_\bS(m,\ov\mu)\,\eta^{n-j}\wedge\ov{\eta^{n+j}}.
\]

\subsubsection*{Strictness of $f_\dag\cH'$}
The cohomology $\cH^jf_\dag\cH'$ is a coherent $\cO_{\Omega_0}$-module. Denote by $h^j$ its generic rank. Notice that, for any $\hb_o\in \Omega_0$, the dimension of the space $H^j\big(X,(\cE^{n+\bbullet}_\cX\otimes_{\cC^{\infty,\an}_\cX}\cH,D_\cH)_{|\hb=\hb_o}\big)$ is equal to that of the space of $\Delta_{\hb_o}$-harmonic $(n+j)$-forms, hence is \emph{independent} of $\hb_o$. It is therefore equal to $h^j$ for any $\hb_o$. Consider now the exact sequence
\begin{multline*}
\cdots\to\cH^jf_\dag\cH'\To{\hb-\hb_o}\cH^jf_\dag\cH' \to H^j\big(X,(\cE^{n+\bbullet}_\cX\otimes\cH)_{|\hb=\hb_o}\big)\\
\To{k_j}\cH^{j+1}f_\dag\cH'\to\cdots
\end{multline*}
If $k_i=0$ for some $i$, then $\cH^if_\dag\cH'$ is locally free at $\hb_o$ and $k_{i-1}=0$, so $\cH^jf_\dag\cH'$ is locally free at $\hb_o$ for any $j\leq i$. As $\cH^jf_\dag\cH'=0$ for $j\gg0$, this shows that $\cH^jf_\dag\cH'$ is locally free for any $j$.

\subsubsection*{Twistor condition}
We want to prove that, for any $j$, the sesquilinear pairing $f_\dag^jC$ is nondegenerate and defines a gluing of weight~$j$. According to the strictness property above, it is enough to show the nondegeneracy after restricting to fibres $\hb=\hb_o$ for $\hb_o\in\bS$. Remark also that $\ov{D_\cH}_{|\hb=\hb_o}=D_{\ov H,\hb_o}$. The Hermitian metric $h$ induces therefore a nondegenerate pairing of flat bundles compatible with the differential:
\[
h:(H,D_{H,\hb_o})\otimes(\ov H,D_{\ov H,-\hb_o})\to (\cC_X^\infty,d).
\]
Poincar\'e duality applied to the de~Rham complex of these flat bundles gives the nondegeneracy of $(f_\dag^jC)_{\hb=\hb_o}$.

\begin{remarque*}
The nondegeneracy can also be obtained as a consequence of the positivity proved below, without referring to Poincar\'e duality.
\end{remarque*}

Consider now the inclusion
\[
\iota:\pi^{-1}\cE_X^{n-j}\ootimes_{\cC^\infty_X}H\hto\cE_\cX^{n-j}\ootimes_{\cC^{\infty,\an}_\cX}\cH
\]
sending $\eta^{p,n-j-p}\otimes m$ to $\hb^{-p}\eta^{p,n-j-p}\otimes m$. We then have $\iota\big(\Harm_X^{n-j}(H)\big)\subset\ker D_\cH$. Moreover, we have seen that, for any $\hb_o$, the projection of $\iota\big(\Harm_X^{n-j}(H)\big)$
\begin{itemize}
\item
in the de~Rham space $H^{n-j}(X;(\cE_X^\cbbullet\otimes H,D_{V_{\hb_o}}))$ if $\hb_0\neq0$, and
\item
in the Dolbeault space $H^{n-j}(X;(\cE_X^\cbbullet\otimes H,D''_E+\theta'_E))$ if $\hb_o=0$,
\end{itemize}
is an isomorphism of $\CC$-vector spaces. It follows that $\iota\big(\Harm_X^{n-j}(H)\big)$ is a lattice in $H^{-j}(X,(\cE^{n+\bbullet}_\cX\otimes\cH,D_\cH))$, \ie that we have
\[
H^{-j}\big(X,(\cE^{n+\bbullet}_\cX\otimes\cH,D_\cH)\big)=\cO_{\Omega_0}\ootimes_\CC\iota\big(\Harm_X^{n-j}(H)\big).
\]
To get the twistor condition, it is now enough to show that $f_\dag^jC$ induces a pairing
\[
\iota\big(\Harm_X^{n-j}(H)\big)\ootimes_\CC\ov{\iota\big(\Harm_X^{n+j}(H)\big)}\to\hb^j\CC.
\]
This follows from the fact that, for sections $\eta_1^{p,q}\otimes m$ of $\cE_X^{p,q}\otimes H$ and $\eta_2^{n-q,n-p}\otimes\mu$ of $\cE_X^{n-q,n-p}\otimes H$ with $p+q=n-j$,
\[
h(m,\ov\mu)\iota(\eta_1^{p,q})\wedge\ov{\iota(\eta_2^{n-q,n-p})}=(-1)^{n-q}\hb^jh(m,\ov\mu)\eta_1^{p,q}\wedge\ov{\eta_2^{n-q,n-p}}.
\]

\subsubsection*{The Lefschetz morphism}
The condition on the Lefschetz morphism, defined on \T\ref{subsec:lef}, comes from the same property for $\omega\wedge{}$ on harmonic sections $\Harm_{X}^{n-j}(H)$, as $\iota(\omega)=\hbm\omega$.

\subsubsection*{Polarization}
We will follow the notation introduced in \T\ref{subsec:grlefschetz}. Put $\cT_j=f_\dag^{-j}\cT$. This is a twistor structure of weight~$-j$. We have $\cT_j=(\cH^j,\cH^{-j},C_j)$ with $C_j=C^{-j}$, where $\cH^j$ stands for $\cH^{j}f_*\big(\cE^{n+\bbullet}_\cX\otimes_{\cC^{\infty,\an}_\cX}\cH,D_\cH\big)$. Moreover we have $C^*_{-j}=C_j$ (\cf Lemma \ref{lem:imdiradjointcomp}\eqref{lem:imdiradjointcomp2}). Hence we may put $\cS_j=(\id,\id):\cT_j\to\cT^*_{-j}$. It clearly satisfies $\cS_j^*=\cS_{-j}$. Moreover, $\cL_\omega$ is skewadjoint with respect to $\cS_j$ as in \eqref{eq:skew}, because by construction we have $\cL_\omega^*=-\cL_\omega$. Let us verify the positivity condition on the primitive part. We are in the situation of Example \ref{ex:lefpol}, with $L''_\omega=\hbm\omega\wedge=\iota(\omega)\wedge$ and $\varepsilon=1$.

Consider first a primitive section $\eta^{p,q}\otimes m_{p,q}$ of $\cE_X^{p,q}\otimes H$ with $p+q=n-j$. Then, by definition, $\eta^{p,q}$ is a primitive $(p,q)$-form. We have $\iota(\eta^{p,q}\otimes m_{p,q})=\hb^{-p}\eta^{p,q}\otimes m_{p,q}$. Taking notation of Example \ref{ex:lefpol}, we want to show that
\[
\thb^jC^{-j}\big(\iota(\omega^j\wedge\eta^{p,q}\otimes m_{p,q}),\ov{\iota(\eta^{p,q}\otimes m_{p,q})}\big)>0.
\]
This amounts to showing that
\[
(-1)^pi^j\frac{\varepsilon(n-j)}{(2i\pi)^n}\int_Xh(m_{p,q},\ov{m_{p,q}})\,\eta^{p,q}\wedge\ov{\eta^{p,q}}\wedge\omega^j>0.
\]
This classically follows from the primitivity of $\eta^{p,q}$, because, denoting by $\star$ the Hodge operator, we have $\varepsilon(n-j)i^{p-q}\ov{\eta^{p,q}}\wedge\omega^j=j!\,\star\ov{\eta^{p,q}}$ (see \eg \cite[\T8.C]{Demailly96}) and $h$ is positive definite.

By decomposing any primitive section of $\cE_X^{p,q}\otimes H$ with respect to an orthonormal basis of $H$, we get the positivity statement for it.

Given any primitive harmonic section in $\Harm_X^{n-j}(H)$, we apply the previous result to any of its $(p,q)$ component, with $p+q=n-j$, to get the positivity.
\end{proof}

\chapterspace{-3}
\chapter{Specializable $\cR_\cX$-modules}\label{chap:spe}
One of the main tools in the theory of polarized Hodge Modules \cite{MSaito86} is the notion of nearby cycles (or specialization) extended to $\cD$-modules. It involves the notions of Bernstein polynomial and Malgrange-Kashiwara filtration, denoted by $V$. The purpose of this chapter is to introduce a category of $\cR_\cX$-modules (or $\RTriples$) for which a good notion of specialization may be defined.

In \T\ref{sec:Vfil}, we recall with details the basic properties of the $V$-filtrations for $\cR_\cX$-modules. We follow \cite{Bibi87, M-S86}.

In \T\ref{sec:reviewspe}, we briefly review the construction of the Malgrange-Kashiwara filtration for coherent $\cD_X$-modules (see \eg \cite{M-S86}). We keep notation of \T\ref{subsec:Vfil}. Recall that this filtration was introduced by M\ptbl Kashiwara \cite{Kashiwara83} in order to generalize previous results by B\ptbl Malgrange \cite{Malgrange83} to arbitrary regular holonomic $\cD$-modules. The presentation we give here comes from various published sources (\eg \cite{Bibi87,M-S86,MSaito86}) and from an unpublished letter of B\ptbl Malgrange to P\ptbl Deligne dated january 1984.

\section{$V$-filtrations}\label{sec:Vfil}

\subsection{}\label{subsec:Vfil}
Let $X'$ be a complex manifold and let $X$ be an open set in $\CC\times X'$. We denote by $t$ the coordinate on $\CC$, that we also view as a function on $X$, and by $\partiall_t$ the corresponding vector field. We set $X_0=t^{-1}(0)\subset X$ (which is open in $X'$) and we denote by $\cI_{X_0}$ (\resp $\cI_{\cX_0}$) its ideal in $\cO_X$ (\resp in $\cO_\cX$).

Denote by $\index{$vr$@$V_{\bbullet}\cR_\cX$}V_{\bbullet}\cR_\cX$ the increasing filtration indexed by $\ZZ$ associated with $X_0$: for any $(x,\hb)\in\cX$,
\[
V_k\cR_{\cX,(x,\hb)}=\{P\in\cR_{\cX,(x,\hb)}\mid P\cdot\cI_{\cX_0,(x,\hb)}^j\subset \cI_{\cX_0,(x,\hb)}^{j-k}\;\forall j\in \ZZ\}
\]
where we put $\cI^\ell=\cO_\cX$ if $\ell\leq0$. In any local coordinate system $(x_2,\ldots,x_n)=x'$ of~$X'$, the germ $P\in\cR_\cX$ is in $V_k\cR_\cX$ iff
\begin{itemize}
\item
$P=\sum_{\bmj=(j_1,\bmj')}a_{\bmj}(t,x',\hb)(t\partiall_t)^{j_1}\partiall_{x'}^{\bmj'}$, if $k=0$:
\item
$P=t^{\module{k}}Q$ with $Q\in V_0\cR_\cX$, if $k\in -\NN$ (\ie $V_k\cR_\cX=t^{\module{k}}V_0\cR_\cX=V_0\cR_\cX\cdot t^{\module{k}}$);
\item
$P=\sum_{0\leq j\leq k}Q_j\partiall_{t}^{j}$ with $Q_j\in V_0\cR_\cX$, if $k\in\NN$ (\ie $V_k\cR_\cX=\sum_{j=0}^k\partiall_t^jV_0\cR_\cX=\sum_{j=0}^kV_0\cR_\cX\cdot\partiall_t^j$).
\end{itemize}

\medskip\noindent
Set $V_k\cO_\cX=V_k\cR_\cX\cap\cO_\cX$. This is nothing but the $\cI_{\cX_0}$-adic filtration on $\cO_\cX$. The following can be proved exactly as for $\cD_X$-modules (see \eg \cite{M-S86}):
\begin{itemize}
\item
$V_k\cR_\cX\cdot V_\ell\cR_\cX\subset V_{k+\ell}\cR_\cX$ with equality for $k,\ell\leq 0$ or $k,\ell\geq 0$.
\item
$V_k\cR_{\cX\moins\cX_0}=\cR_{\cX\moins\cX_0}$ for any $k\in\ZZ$.
\item
$\lefpar \cap_k{V_k\cR_\cX}\rigpar_{|\cX_0}=\{0\}$.
\end{itemize}

\begin{definition}
Let $\cM$ be a left $\cR_\cX$-module. A \emph{$V$-filtration} of $\cM$ is an increasing filtration $U_\bbullet\cM$ indexed by $\ZZ$, which is exhaustive, and such that, for any $k,\ell\in\ZZ$, we have $V_k\cR_\cX\cdot U_\ell\cM\subset U_{k+\ell}\cM$.
\end{definition}

\begin{Remarques}
\begin{enumerate}
\item
We will identify the sheaf of rings $\gr_0^V\cR_\cX\defin V_0\cR_\cX\big/V_{-1}\cR_\cX$, which is supported on $\cX_0$, with the ring $\cR_{\cX_0}[t\partiall_t]$, still denoting by $t\partiall_t$ the class of $t\partiall_t$ in $\gr_0^V\cR_\cX$. In particular, $\cR_{\cX_0}$ is a subring of $\gr_0^V\cR_\cX$. The class of $t\partiall_t$ commutes with any section of $\cR_{\cX_0}$.
\item
Given a holomorphic function $f:X'\to\CC$ on a complex manifold $X'$ and a $\cR_{\cX'}$-module $\cM$, we will usually denote by $i_f:X'\hto X=\CC\times X'$ the inclusion of the graph of $f$, by $t$ the coordinate on $\CC$, and we will consider $V$-filtrations on the $\cR_\cX$-module $i_{f,+}\cM$.
\end{enumerate}
\end{Remarques}

\Subsection{Coherence}\label{subsec:coh}
\subsubsection*{Coherence of the Rees sheaf of rings}
Introduce the Rees sheaf of rings $R_V\cR_\cX=\oplus_kV_k\cR_\cX\cdot q^k$, where $q$ is a new variable, and similarly $R_V\cO_\cX=\oplus_kV_k\cO_\cX\cdot q^k$, which is naturally a $\cO_\cX$-module. Let us recall some basic coherence properties of these sheaves on $\cX$.

Let $\cK$ be a compact polycylinder in $\cX$. Then $R_V\cO_\cX(\cK)=R_V(\cO_\cX(\cK))$ is Noetherian, being the Rees ring of the $\cI_{\cX_0}$-adic filtration on the Noetherian ring $\cO_\cX(\cK)$ (Theorem of Frisch). Similarly, as $\cO_{\cX,(x,\hb)}$ is flat on $\cO_\cX(\cK)$ for any $(x,\hb)\in\cK$, the ring $(R_V\cO_\cX)_{(x,\hb)}=R_V\cO_\cX(\cK)\otimes_{\cO_\cX(\cK)}\cO_{\cX_{(x,\hb)}}$ is flat on $R_V\cO_\cX(\cK)$.

Let us show that $R_V\cO_\cX$ is coherent on $\cX$. Let $\cU$ be any open set in $\cX$ and let $\varphi:(R_V\cO_\cX)^q_{|\cU}\to(R_V\cO_\cX)^p_{|\cU}$ be any morphism. Let $\cK$ be a polycylinder contained in $\cU$. Then, $\ker\varphi(\cK)$ is finitely generated over $R_V\cO_\cX(\cK)$ by noetherianity and, if $\cV$ is the interior of $\cK$, we have $\ker\varphi_{|\cV}=\ker\varphi(\cK)\otimes_{R_V\cO_\cX(\cK)}(R_V\cO_\cX)_{|\cV}$ by flatness. So $\ker\varphi_{|\cV}$ is finitely generated, proving the coherence of $R_V\cO_\cX$.

Before considering $R_V\cR_\cX$, consider the sheaf $\cO_\cX[\tau,\xi_2,\dots,\xi_n]$ equipped with the $V$-filtration for which $\tau$ has degree $1$, $\xi_2,\dots,\xi_n$ have degree $0$, and inducing the $V$-filtration on $\cO_\cX$. Firstly, forgetting $\tau$, we have $R_V(\cO_\cX[\xi_2,\dots,\xi_n])=(R_V\cO_\cX)[\xi_2,\dots,\xi_n]$. Secondly, $V_k(\cO_\cX[\tau,\xi_2,\dots,\xi_n])=\sum_{j\geq0}V_{k-j}(\cO_\cX[\xi_2,\dots,\xi_n])\tau^j$ for any $k\in\ZZ$, hence we have a surjective morphism
\begin{align*}
R_V\cO_\cX[\xi_2,\dots,\xi_n]\otimes_\CC\CC[\tau']&\To{} R_V(\cO_\cX[\tau,\xi_2,\dots,\xi_n])\\
V_\ell\cO_\cX[\xi_2,\dots,\xi_n] q^\ell\tau^{\prime j}&\mto V_\ell\cO_\cX[\xi_2,\dots,\xi_n] \tau^jq^{\ell+j}.
\end{align*}
If $\cK\subset\cX$ is any polycylinder, then $\big(R_V\cO_\cX(\cK)\big)[\tau',\xi_2,\dots,\xi_n]$ is Noetherian. Therefore, $R_V(\cO_\cX[\tau,\xi_2,\dots,\xi_n])(\cK)$ is Noetherian.

As $R_V\cR_\cX$ can be filtered (by the degree of the operators) in such a way that, locally on $\cX$, $\gr R_V\cR_\cX$ is isomorphic to $R_V(\cO_\cX[\tau,\xi_2,\dots,\xi_n])$, this implies that, if~$\cK$ is any sufficiently small polycylinder, then $R_V\cR_\cX(\cK)$ is Noetherian. Using the previous results and standard arguments, one concludes that $R_V\cR_\cX$ is coherent.

\subsubsection*{Good $V$-filtrations}
Let $(\cM,U_\bbullet \cM)$ be a $V$-filtered $\cR_\cX$-module. The filtration is \emph{good} if, for any compact set $\cK\subset \cX$, there exists $k_0\geq 0$ such that, in a neighbourhood of $\cK$, we have for all $k\geq k_0$
\[
U_{-k}\cM=t^{k-k_0}U_{-k_0}\cM\qqbox{and}U_k\cM=\sum_{0\leq j\leq k-k_0}\partiall_{t}^{j}U_{k_0}\cM,
\]
and any $U_\ell \cM$ is $V_0\cR_\cX$-coherent.

The filtration $U_{\bbullet}\cM$ is good if and only if the Rees module $\oplus_kU_k\cM\cdot q^k$ is coherent over $R_V\cR_\cX$. Equivalently, there should exist, locally on $\cX$, a presentation $\cR_\cX^b\to\cR_\cX^a\to\cM\to0$, inducing for each $k\in\ZZ$ a presentation \hbox{$U_k\cR_\cX^b\!\to\! U_k\cR_\cX^a\!\to\! U_k\cM\!\to\!0$}, where the filtration on the free modules $\cR_\cX^a,\cR_\cX^b$ are obtained by suitably shifting $V_\bbullet\cR_\cX$ on each summand. In particular, we get

\begin{lemme}\label{lem:tinj}
Locally on $\cX$, there exists $k_0$ such that, for any $k\leq k_0$, $t:U_{-k}\cM\to U_{-k-1}\cM$ is bijective.
\end{lemme}

\begin{proof}
Indeed, using a presentation of $\cM$ as above, it is enough to show the lemma for $\cR_\cX^a$ with a filtration as above, and we are reduced to consider each summand $\cR_\cX$ with a shifted standard $V$-filtration $U_\bbullet\cR_\cX$. There, we may choose $k_0$ such that $U_{k_0}\cR_\cX=V_0\cR_\cX$.
\end{proof}

In a similar way we get:

\begin{lemme}\label{lem:ttors}
Let $\cU$ be a coherent $V_0\cR_\cX$-module and let $\cT$ be its $t$-torsion subsheaf, \ie the subsheaf of local sections locally killed by some power of $t$. Then, locally on $\cX$, there exists $\ell$ such that $\cT\cap t^\ell\cU=0$.
\end{lemme}

\begin{proof}
Consider the $t$-adic filtration on $V_0\cR_\cX$, \ie the filtration $V_j\cR_\cX$ with $j\leq0$. Then the filtration $t^{-j}\cU$ is good with respect to it, and locally we have a surjective morphism $(V_0\cR_\cX)^n\to\cU$ which is strict with respect to the $V$-filtration. Its kernel~$\cK$ is coherent and comes equipped with the induced $V$-filtration, which is good. In particular, locally on $\cX$, there exists $j_0\leq0$ such that $V_{j+j_0}\cK=t^{-j}V_{j_0}\cK$ for any $j\leq0$. For any $j\leq0$ we therefore have locally an exact sequence
\[
(V_j\cR_\cX)^m\to(V_{j+j_0}\cR_\cX)^n\to t^{-(j+j_0)}\cU\to0.
\]
As $t:V_k\cR_\cX\to V_{k-1}\cR_\cX$ is bijective for $k\leq0$, we conclude that $t:t^{-j_0}\cU\to t^{-j_0+1}\cU$ is so, hence $\cT\cap t^{-j_0}\cU=0$.
\end{proof}

\begin{proposition}\label{prop:V-coh}
If $\cN$ is a coherent $\cR_\cX$-submodule of $\cM$ and $U_\bbullet \cM$ is a good filtration of $\cM$, then the $V$filtration $U_\bbullet \cN\defin \cN\cap U_\bbullet \cM$ is also good.
\end{proposition}

\begin{proof}
It is now standard (it follows from coherence properties of the Rees module $\oplus_kU_k\cM\cdot q^k$, see \eg \cite{Mebkhout87}).
\end{proof}

\begin{Remarques}\label{rem:leftrightspe}
\begin{enumerate}
\item\label{rem:leftrightspe1}
It is straightforward to develop the theory in the case of right $\cR_\cX$-modules. If $U_\bbullet(\cM)$ is a $V$-filtration of the left module $\cM$, then $U_\bbullet(\omega_\cX\otimes_{\cO_\cX}\cM)\defin \omega_\cX\otimes_{\cO_\cX}U_\bbullet(\cM)$ is the corresponding filtration of the corresponding right module. This correspondence is compatible with taking the graded object with respect to $U_\bbullet$. The operator $-\partiall_{t}t$ (acting on the left) corresponds to $t\partiall_t$ (acting on the right).

\item\label{rem:leftrightspe2}
Given an increasing filtration $U_\bbullet$ (lower indices), we define the associated decreasing filtration (upper indices) by $U^k=U_{-k-1}$. If $b(-(\partiall_tt+k\hb))\cdot \gr^U_k\cM=0$ for all $k\in\ZZ$, we have $b'(t\partiall_t-\ell\hb)\cdot \gr_U^\ell\cM=0$ for all $\ell\in\ZZ$, if we put $b'(s)=b(-s)$.
\end{enumerate}
\end{Remarques}

\subsection{$V$-filtration and direct images}\label{subsec:Vfildirim}

The purpose of this section is to establish the compatibility between taking a direct image and taking a graded part of a $V$-filtered $\cR_\cX$-module. We will give an analogue of Proposition 3.3.17 of \cite{MSaito86}.

\begin{definition}\label{def:monodromic}
Let $\cM$ be a left $\cR_\cX$-module equipped with an exhaustive increasing filtration $U_{\bbullet}\cM$ indexed by $\ZZ$ such that $V_k\cR_\cX\cdot U_\ell \cM\subset U_{k+\ell}\cM$ for any $k,\ell\in \ZZ$. We say that $(\cM,U_\bbullet\cM)$ is \emph{monodromic} if, locally on $\cX$, there exists a \emph{monic} polynomial $b(s)\in\CC[\hb][s]$ such that
\begin{enumerate}
\item
$b(-(\partiall_tt+k\hb))\cdot \gr^U_k\cM=0$ for all $k\in\ZZ$,
\item
$\gcd(b(s-k\hb),b(s-\ell\hb))\in\CC[\hb]\moins\{0\}$ for all $k\neq\ell$.
\end{enumerate}
For right $\cR_\cX$-modules, we use the convention of Remark~\ref{rem:leftrightspe}\eqref{rem:leftrightspe1}.
\end{definition}

\begin{theoreme}\label{th:imdirspe}
Let $f:X\to Y$ be holomorphic map between complex analytic manifolds and let $t\in\CC$ be a new variable. Put $F=f\times\id:X\times\CC\to Y\times\CC$. Let $\cM$ be a right $\cR_{\cX\times\CC}$-module equipped with a $V$-filtration $U_\bbullet\cM$ (relative to the function $t:X\times\CC\to\CC$). Then $U_\bbullet\cM$ defines canonically and functorially a $V$-filtration $U_\bbullet\cH^i(F_\dag\cM)$.

Assume that $F$ is proper on the support of $\cM$.
\begin{enumerate}
\item\label{th:imdirspe1}
If $\cM$ is good and $U_\bbullet\cM$ is a good $V$-filtration, then $U_\bbullet\cH^i(F_\dag\cM)$ is a good $V$-filtration.
\item\label{th:imdirspe2}
If moreover $(\cM,U\bbullet\cM)$ is monodromic and $f_\dag\gr^U\cM$ is strict, then one has a canonical and functorial isomorphism of $\cR_\cY$-modules ($k\in\ZZ$)
$$
\gr_k^U\big(\cH^iF_\dag\cM\big)=\cH^i\big(f_\dag\gr_k^U\cM\big),
$$
$\gr^U\big(\cH^iF_\dag\cM\big)$ is monodromic and strict.
\end{enumerate}
\end{theoreme}

\begin{remarque}
In the last assertion, we view $\gr_k^U\cM$ as a right $\cR_\cX$-module, and $f_\dag$ is defined as in \T\ref{subsec:imdir}. By functoriality, the action of $t\partiall_t$ descends to $\cH^i(f_\dag\gr_k^U\cM)$.
\end{remarque}

\begin{proof}
We will use the isomorphism $F_\dag=f_\dag$ for $\cM$ (see Remark \ref{rem:imdir}\eqref{rem:imdir2}), \ie we take the direct image viewing $\cM$ as a $\cR_{\cX\times\CC/\CC}$ equipped with a compatible action of $\partiall_t$. Put $\cN^{\cbbullet}=f_\dag\cM$. This complex is naturally filtered by $U_\bbullet\cN^{\cbbullet}\defin f_\dag U_\bbullet\cM$. Therefore, we define the filtration on its cohomology by
\[
U_\bbullet\cH^i(F_\dag\cM)=U_\bbullet\cH^i(f_\dag\cM)\defin\image \big[\cH^i(f_\dag U_\bbullet\cM)\to\cH^i(f_\dag\cM)\big].
\]
Notice that, for any $j$, $f_\dag U_j\cM$ is the direct image of $U_j\cM$ viewed as a $\cR_{\cX\times\CC/\CC}$-module, on which we put the natural action of $t\partiall_t$.

The relation with the Rees construction is given by the following lemma:

\begin{lemme}\label{lem:Hfiltre}
Let $(\cN^{\cbbullet},U_\bbullet\cN^{\cbbullet})$ be a $V$-filtered complex of $\cR_{\cY\times\CC}$-modules. Put $U_j\cH^i(\cN^{\cbbullet})\defin\image \big[\cH^i(U_j\cN^{\cbbullet})\to\cH^i(\cN^{\cbbullet})\big]$. Then we have
\[
\cH^i(R_U\cN^{\cbbullet})\big/q\text{\rm-torsion}=R_U\cH^i(\cN^{\cbbullet}).
\]
In particular, if $R_U\cN^{\cbbullet}$ has $\cR_{\cY\times\CC}$-coherent cohomology, then $U_\bbullet\cH^i(\cN^{\cbbullet})$ is a good $V$-filtration.
\end{lemme}

\begin{proof}
One has a surjective morphism of graded modules $\cH^i(R_U\cN^{\cbbullet})\to R_U\cH^i(\cN^{\cbbullet})$, by definition, and this morphism induces an isomorphism after tensoring with $\CC[q,q^{-1}]$.
\end{proof}

\begin{lemme}\label{lem:Vcoh}
If $\cM$ is good, then any coherent $V_0\cR_\cX$-submodule is good.
\end{lemme}

\begin{proof}
As a coherent $V_0\cR_\cX$-submodule of $\cM$ induces on any subquotient of $\cM$ a coherent $V_0\cR_\cX$-submodule, we may reduce to the case where $\cM$ has a good filtration. It is then enough to prove that any coherent $V_0\cR_\cX$-submodule $\cN$ of $\cM$ is contained in such a submodule having a good filtration. If $\cF$ is a $\cO_\cX$-coherent submodule of~$\cM$ which generates $\cM$, then $\cN$ is contained in $V_k\cR_\cX\cdot \cF$ for some $k$, hence the result.
\end{proof}

This lemma allows one to apply Grauert's coherence theorem to each $U_j$, in order to get that each $f_\dag U_j\cM$ has $V_0\cR_\cX$-coherent cohomology under the properness assumption. We conclude that for each $i,j$, $U_j\cH^if_\dag\cM$ is $V_0\cR_\cX$-coherent.

In order to end the proof of \eqref{th:imdirspe1}, we need to prove that each $U_\bbullet\cH^if_\dag\cM$ is a good $V$-filtration. We will compute directly the Rees module associated with this filtration, in order to get its coherence. Let us first consider the analogue of Lemma \ref{lem:Vcoh}.

\smallskip
Keep notation of \T\ref{subsec:coh}. The graded ring $R_V\cR_\cX$ is filtered by the degree in the derivatives $q\partiall_{x_j}$ and the degree-zero term of the filtration is $R_V\cO_\cX$, with $V_k\cO_\cX=\cO_\cX$ for $k\geq0$ and ${}=t^{-k}\cO_\cX$ for $k\leq0$.

Let $(\cM,U_\bbullet\cM)$ be a $V$-filtered right $\cR_\cX$-module and let $R_U\cM$ be the associated Rees module. We therefore have the notion of a good filtration on $R_U\cM$ (by coherent graded $R_V\cO_\cX$-submodules). If $R_U\cM$ has a good filtration (or equivalently if $R_U\cM$ is generated by a coherent graded $R_V\cO_\cX$-module), it is $R_V\cR_\cX$-coherent and has a left resolution by coherent ``induced'' graded $R_V\cR_\cX$-modules, of the form $G\otimes_{R_V\cO_\cX}R_V\cR_\cX$, where $G$ is graded $R_V\cO_\cX$-coherent. We may even assume (by killing the $q$-torsion) that each term $G\otimes_{R_V\cO_\cX}R_V\cR_\cX$ has no $q$-torsion, or in other words that it takes the form $R_U(L\otimes_{\cO_\cX}\cR_\cX)$, where $L$ is $\cO_\cX$-coherent, having support contained in $\supp\cM$, and equipped with a good $V$-filtration (\ie a good $\cI_{\cX_0}$-adic filtration) and $U_\bbullet(L\otimes_{\cO_\cX}\cR_\cX)$ is defined in the usual way.

We say that $R_U\cM$ is \emph{good} if, in the neighbourhood of any compact set $\cK\subset\cX$, $R_U\cM$ is a finite successive extension of graded $R_V\cR_\cX$-modules having a good filtration.

\begin{lemme}
Assume that $\cM$ is a good $\cR_\cX$-module and let $U_\bbullet\cM$ be a good $V$-filtration of $\cM$. Then $R_U\cM$ is a good graded $R_V\cR_\cX$-module.
\end{lemme}

\begin{proof}
Fix a compact set $\cK\subset\cX$. First, it is enough to prove the lemma when $\cM$ has a good filtration in some neighbourhood of $\cK$, because a good $V$-filtration $U_\bbullet\cM$ induces naturally on any subquotient $\cM'$ of $\cM$ a good $V$-filtration, so that $R_U\cM'$ is a subquotient of $R_U\cM$.

Therefore, assume that $\cM$ is generated by a coherent $\cO_\cX$-module $\cF$, \ie $\cM=\cR_\cX\cdot\cF$. Consider the $V$-filtration $U'_\bbullet\cM$ generated by $\cF$, \ie $U'_\bbullet\cM=V_\bbullet\cR_\cX\cdot\cF$. Then, clearly, $R_V\cO_\cX\cdot\cF=\oplus_k V_k\cO_\cX\cdot\cF q^k$ is a coherent graded $R_V\cO_\cX$-module which generates $R_{U'}\cM$.

If the filtration $U''_\bbullet\cM$ is obtained from $U'_\bbullet\cM$ by a shift by $-\ell\in\ZZ$, \ie if $R_{U''}\cM=q^\ell R_{U'}\cM\subset\cM[q,q^{-1}]$, then $R_{U''}\cM$ is generated by the $R_V\cO_\cX$-coherent submodule $q^\ell R_V\cO_\cX\cdot\cF$.

On the other hand, let $U''_\bbullet\cM$ be a good $V$-filtration such that $R_{U''}\cM$ has a good filtration. Then any good $V$-filtration $U_\bbullet\cM$ such that $U_k\cM\subset U''_k\cM$ for any $k$ satisfies the same property, because $R_U\cM$ is thus a coherent graded submodule of $R_{U''}\cM$, so a good filtration on the latter induces a good filtration on the former.

As any good $V$-filtration $U_\bbullet\cM$ is contained, in some neighbourhood of $\cK$, in the good $V$-filtration $U'_\bbullet\cM$ suitably shifted, we get the lemma.
\end{proof}
To end the proof of Part \eqref{th:imdirspe1}, it is therefore enough to prove it for induced modules $\cM=L\otimes_{\cO_\cX}\cR_\cX$, with $L$ coherent over $\cO_\cX$ and $F_{\supp L}$ proper. We will indicate it when $f:X=Y\times Z\to Y$ is the projection. We then have
\begin{align*}
U_j(L\otimes_{\cO_{\cX\times\CC}}\cR_{\cX\times\CC})&=U_j\Big[(L\otimes_{f^{-1}\cO_{\cY\times\CC}}f^{-1}\cR_{\cY\times\CC})\ootimes_{f^{-1}\cO_{\cY\times\CC}}\cR_{\cX\times\CC/\cY\times\CC}\Big]\\
&=U_j(L\otimes_{f^{-1}\cO_{\cY\times\CC}}f^{-1}\cR_{\cY\times\CC})\ootimes_{f^{-1}\cO_{\cY\times\CC}}\cR_{\cX\times\CC/\cY\times\CC},
\end{align*}
because the $V$-filtration on $\cR_{\cX\times\CC/\cY\times\CC}$ is nothing but the $t$-adic filtration. Now, we have
\begin{align*}
f_\dag U_j(L\otimes_{\cO_{\cX\times\CC}}\cR_{\cX\times\CC})&= \bR f_*U_j(L\otimes_{f^{-1}\cO_{\cY\times\CC}}f^{-1}\cR_{\cY\times\CC})\\
&= U_j(\bR f_* L\otimes_{\cO_{\cY\times\CC}}\cR_{\cY\times\CC}),
\end{align*}
if we filter the complex $\bR f_* L$ by subcomplexes $\bR f_*U_j(L)$ and we filter the tensor product as usual. By Grauert's theorem applied to coherent $R_V\cO_{\cX\times\CC}$-sheaves, $\bR f_* R_U L$ is $R_V\cO_{\cY\times\CC}$-coherent, hence $f_\dag R_U(L\otimes_{\cO_\cX}\cR_\cX)$ is $R_V\cR_{\cY\times\CC}$-coherent. After Lemma \ref{lem:Hfiltre}, we get \ref{th:imdirspe}\eqref{th:imdirspe1}.\qed

\medskip
In order to get Part \eqref{th:imdirspe2} of the theorem, we will first prove:

\begin{proposition}\label{prop:saito}
Let $(\cN^{\cbbullet},U_\bbullet\cN^{\cbbullet})$ be a $V$-filtered complex of $\cR_{\cY\times\CC}$-modules. Assume that
\begin{enumerate}
\item\label{prop:saito1}
the complex $\gr^U\cN^{\cbbullet}$ is strict and monodromic,
\item\label{prop:saito2}
there exists $j_0$ such that for all $j\leq j_0$ and all $i$, the left multiplication by $t$ induces an isomorphism $t:U_j\cN^i\isom U_{j-1}\cN^i$,
\item\label{prop:saito3}
There exists $i_0\in\ZZ$ such that, for all $i\geq i_0$ and any $j$, one has \hbox{$\cH^i(U_j\cN^{\cbbullet})=0$}.
\end{enumerate}
Then for any $i,j$ the morphism $\cH^i(U_j\cN^{\cbbullet})\to\cH^i(\cN^{\cbbullet})$ is injective. Moreover, the filtration $U_\bbullet\cH^i(\cN^{\cbbullet})$ defined by
\[
U_j\cH^i(\cN^{\cbbullet})=\image\lefcro \cH^i(U_j\cN^{\cbbullet})\to\cH^i(\cN^{\cbbullet})\rigcro
\]
satisfies $\gr^{U}\cH^i(\cN^{\cbbullet})=\cH^i(\gr^{U}\cN^{\cbbullet})$.
\end{proposition}

\begin{proof}
It will have three steps.
\subsubsection*{First step} This step proves a formal analogue of the conclusion of the proposition. Put
\[
\wh{U_j\cN^{\cbbullet}}=\liml{\ell}U_j\cN^{\cbbullet}/U_\ell\cN^{\cbbullet}\qqbox{and}\wh\cN^{\cbbullet}=\limr{j}\wh{U_j\cN^{\cbbullet}}.
\]
Under the assumption of Proposition \ref{prop:saito}, we will prove the following:
{\def\theenumi{\alph{enumi}}
\begin{enumerate}
\item\label{itema}
For all $k\leq j$, $\wh{U_k\cN^{\cbbullet}}\to\wh{U_j\cN^{\cbbullet}}$ is injective (hence, for all $j$, $\wh{U_j\cN^{\cbbullet}}\to\wh\cN^{\cbbullet}$ is injective) and $\wh{U_j\cN^{\cbbullet}}/\wh{U_{j-1}\cN^{\cbbullet}}=U_j\cN^{\cbbullet}/U_{j-1}\cN^{\cbbullet}$.
\item\label{itemb}
For any $k\leq j$, $\cH^i(U_j\cN^{\cbbullet}/U_k\cN^{\cbbullet})$ is strict.
\item\label{itemc}
$\cH^i(\wh{U_j\cN^{\cbbullet}})=\varprojlim_{\ell}\,\cH^i(U_j\cN^{\cbbullet}/U_\ell\cN^{\cbbullet})$.
\item\label{itemd}
$\cH^i(\wh{U_j\cN^{\cbbullet}})\to\cH^i(\wh\cN^{\cbbullet})$ is injective.
\item\label{iteme}
$\cH^i(\wh\cN^{\cbbullet})=\varinjlim_{j}\,\cH^i(\wh{U_j\cN^{\cbbullet}})$.
\end{enumerate}}

Define $U_j\cH^i(\wh\cN^{\cbbullet})=\image\lefcro \cH^i(\wh{U_j\cN^{\cbbullet}})\to\cH^i(\wh\cN^{\cbbullet})\rigcro$. Then the statements \eqref{itema} and \eqref{itemd} imply that
\[
\gr_{j}^{U}\cH^i(\wh\cN^{\cbbullet})=\cH^i(\wh{U_j\cN^{\cbbullet}}/\wh{U_{j-1}\cN^{\cbbullet}})=\cH^i(\gr_{j}^{U}\cN^{\cbbullet}).
\]

\medskip
For $\ell<k<j$ consider the exact sequence of complexes
\[
0\to U_k\cN^{\cbbullet}/U_\ell\cN^{\cbbullet}\to U_j\cN^{\cbbullet}/U_\ell\cN^{\cbbullet}\to U_j\cN^{\cbbullet}/U_k\cN^{\cbbullet}\to 0.
\]
As the projective system $(U_j\cN^{\cbbullet}/U_\ell\cN^{\cbbullet})_\ell$ trivially satisfies the Mittag-Leffler condition (ML), the sequence remains exact after passing to the projective limit, so we get an exact sequence of complexes
\[
0\to\wh{U_k\cN^{\cbbullet}}\to\wh{U_j\cN^{\cbbullet}}\to U_j\cN^{\cbbullet}/U_k\cN^{\cbbullet}\to 0,
\]
hence \eqref{itema}.

Let us show by induction on $n\geq 1$ that, for all $i$ and $j$,
\begin{enumerate}
\item[(b)${}_n$]
$\cH^i(U_{j}\cN^{\cbbullet}/U_{j-n}\cN^{\cbbullet})$ is strict (hence \eqref{itemb});
\end{enumerate}
Indeed, (b)${}_1$ follows from Assumption \ref{prop:saito}\eqref{prop:saito1}. Remark also that, by induction on $n\geq1$, \ref{prop:saito}\eqref{prop:saito1} implies that, for any $n,\ell,i$, $\cH^i(U_\ell/U_{\ell-n})$ is killed by $\prod_{k=\ell-n+1}^\ell b(\partiall_tt+k\hb)$.

For $n\geq2$, consider the exact sequence
\begin{multline*}
\cdots\to\cH^i(U_{j-1}/U_{j-n})\to\cH^i(U_{j}/U_{j-n})\to\cH^i(U_{j}/U_{j-1})\\
\To{\psi}\cH^{i+1}(U_{j-1}/U_{j-n})\to\cdots
\end{multline*}

Any local section of $\im\psi$ is then killed by $b(\partiall_tt+j\hb)$ and $\prod_{k=j-n+1}^{j-1}b(\partiall_tt+k\hb)$, hence by a nonzero holomorphic function of $\hb$. By strictness (b)${}_{n-1}$ applied to $\cH^{i+1}(U_{j-1}/U_{j-n})$, this implies that $\psi=0$, so the previous sequence of $\cH^i$ is exact and $\cH^i(U_{j}/U_{j-n})$ is also strict, hence (b)${}_n$.

By the same argument, we get an exact sequence, for all $\ell<k<j$,
\begin{equation}\label{eq:HiVex}
0\to\cH^i(U_k\cN^{\cbbullet}/U_\ell\cN^{\cbbullet})\to\cH^i(U_j\cN^{\cbbullet}/U_\ell\cN^{\cbbullet})\to\cH^i(U_j\cN^{\cbbullet}/U_k\cN^{\cbbullet})\to 0.
\end{equation}
Consequently, the projective system $(\cH^i(U_j\cN^{\cbbullet}/U_\ell\cN^{\cbbullet}))_\ell$ satisfies (ML), so we get \eqref{itemc} (see \eg \cite[Prop\ptbl 1.12.4]{K-S90}). Moreover, taking the limit on $\ell$ in the previous exact sequence gives, according to (ML), an exact sequence
\[
0\to\cH^i(\wh{U_k\cN^{\cbbullet}})\to\cH^i(\wh{U_j\cN^{\cbbullet}})\to\cH^i(U_j\cN^{\cbbullet}/U_k\cN^{\cbbullet})\to 0,
\]
hence \eqref{itemd}. Now, \eqref{iteme} is clear.

\subsubsection*{Second step}
For any $i,j$, denote by $\cT_j^i\subset\cH^i(U_j\cN^{\cbbullet})$ the $t$-torsion subsheaf of $\cH^i(U_j\cN^{\cbbullet})$. We will now prove that it is enough to show that there exists $j_0$ such that, for each $i$ and each $j\leq j_0$,
\begin{equation}\label{eq:HiVinj}
\cT_j^i=0.
\end{equation}

Assume that \eqref{eq:HiVinj} is proved (step $3$). Let $j\leq j_0$ and let $\ell\geq j$. Then, by definition of a $V$-filtration, $t^{\ell-j}$ acts by $0$ on $U_\ell\cN^{\cbbullet}/U_j\cN^{\cbbullet}$, so that the image of $\cH^{i-1}(U_\ell\cN^{\cbbullet}/U_j\cN^{\cbbullet})$ in $\cH^i(U_j\cN^{\cbbullet})$ is contained in $\cT_j^i$, and thus is zero. We therefore have an exact sequence for any $i$:
\[
0\to\cH^i(U_j\cN^{\cbbullet})\to\cH^i(U_\ell\cN^{\cbbullet})\to\cH^i(U_\ell\cN^{\cbbullet}/U_j\cN^{\cbbullet})\to0.
\]
Using \eqref{eq:HiVex}, we get for any $\ell$ the exact sequence
\[
0\to\cH^i(U_{\ell-1}\cN^{\cbbullet})\to\cH^i(U_\ell\cN^{\cbbullet})\to\cH^i(\gr^U_\ell\cN^{\cbbullet})\to0.
\]
This implies that $\cH^i(U_\ell\cN^{\cbbullet})\to\cH^i(\cN^{\cbbullet})$ is injective. Put $U_\ell\cH^i(\cN^{\cbbullet})=\image\lefcro \cH^i(U_\ell\cN^{\cbbullet})\to\cH^i(\cN^{\cbbullet})\rigcro$. We thus have, for any $i,\ell\in\ZZ$,
\[
\gr_{\ell}^{U}\cH^i(\cN^{\cbbullet})=\cH^i(\gr_{\ell}^{U}\cN^{\cbbullet}).
\]

\subsubsection*{Third step: proof of \eqref{eq:HiVinj}}
Remark first that, according to \ref{prop:saito}\eqref{prop:saito2}, the multiplication by $t$ induces an isomorphism $t:\wh{U_j\cN^{\cbbullet}}\to\wh{U_{j-1}\cN^{\cbbullet}}$ for $j\leq j_0$, and that \eqref{itemd} in Step one implies that, for all $i$ and all $j\leq j_0$, the multiplication by $t$ on $\cH^i(\wh{U_j\cN^{\cbbullet}})$ is injective.

The proof of \eqref{eq:HiVinj} is done by decreasing induction on $i$. It clearly hods for $i\geq i_0$ (given by \ref{prop:saito}\eqref{prop:saito3}). Assume that, for any $j\leq j_0$, we have $\cT_j^{i+1}=0$. We have (after \ref{prop:saito}\eqref{prop:saito2}) an exact sequence of complexes, for any $\ell\geq 0$,
$$
0\to U_j\cN^{\cbbullet}\To{t^\ell}U_j\cN^{\cbbullet}\to U_j\cN^{\cbbullet}\big/ U_{j-\ell}\cN^{\cbbullet}\to 0.
$$
As $\cT_j^{i+1}=0$, we have, for any $\ell\geq0$ an exact sequence
\[
\cH^i(U_j\cN^{\cbbullet})\To{t^\ell} \cH^i(U_j\cN^{\cbbullet})\to\cH^i(U_j\cN^{\cbbullet}/U_{j-\ell}\cN^{\cbbullet})\to0,
\]
hence, according to Step one,
\[
\cH^i(\wh{U_j\cN^{\cbbullet}})/\cH^i(\wh{U_{j-\ell}\cN^{\cbbullet}})=\cH^i(U_j\cN^{\cbbullet}/U_{j-\ell}\cN^{\cbbullet})=\cH^i(U_j\cN^{\cbbullet})/t^\ell\cH^i(U_j\cN^{\cbbullet}).
\]
According to Lemma \ref{lem:ttors}, for $\ell$ big enough (locally on $\cX$), the map $\cT_j^i\to\cH^i(U_j\cN^{\cbbullet})/t^\ell\cH^i(U_j\cN^{\cbbullet})$ is injective. It follows that $\cT_j^i\to\cH^i(\wh{U_j\cN^{\cbbullet}})$ is injective too. But we know that $t$ is injective on $\cH^i(\wh{U_j\cN^{\cbbullet}})$ for $j\leq j_0$, hence $\cT_j^i=0$, thus concluding Step $3$.
\end{proof}

We apply the proposition to $\cN^{\cbbullet}=f_\dag\cM$ equipped with $U_\bbullet\cN^{\cbbullet}=f_\dag U_\bbullet\cM$ to get \ref{th:imdirspe}\eqref{th:imdirspe2}. That Assumption \eqref{prop:saito1} in the proposition is satisfied follows from the assumptions in \ref{th:imdirspe}\eqref{th:imdirspe2}. Assumption \eqref{prop:saito2} is a consequence of the fact that $U_\bbullet\cM$ is a good $V$-filtration and Lemma \ref{lem:tinj}. Last, Assumption \eqref{prop:saito3} is satisfied because $f$ has finite cohomological dimension.
\end{proof}

\subsection{Regularity}\label{subsec:regularity}
We keep notation of \T\ref{subsec:Vfil}. We can identify the sheaf $\cR_{\cX/\CC}$ of differential operators relative to the function $t$ (constructed from the sheaf $\cD_{X/\CC}$ by the Rees procedure) as the subsheaf of $V_0\cR_\cX$ of operators commuting with $t$.

We say that the $V$-filtered $\cR_\cX$-module $(\cM,U_\bbullet\cM)$ is \emph{regular along $\cX_0$} if, for all $k\in\ZZ$, $U_k\cM$ is $\cR_{\cX/\CC}$-coherent near $\cX_0$. If such a condition is satisfied for some good filtration $U_\bbullet\cM$, it is satisfied for any. In an exact sequence, the extreme terms are regular along $\cX_0$ if and only if the middle term is so.

By an argument analogous to that of Lemma \ref{lem:Vcoh}, and applying Grauert's theorem, one proves that, in the situation of Theorem \ref{th:imdirspe}, if $\cM$ is good and regular along $\cX\times\{0\}$, then $F_+\cM$ is regular along $\cY\times\{0\}$.

\section{Review on specializable $\cD_X$-modules}\label{sec:reviewspe}
We keep notation of \T\ref{subsec:Vfil}. A coherent left $\cD_X$-module $M$ is said to be \emph{specializable} along $X_0$ if any local section $m$ of $M$ has a Bernstein polynomial $\index{$bm$@$b_m$}b_m(s)\in\CC[s]\moins\{0\}$ such that $b_m(-\partial_tt)m\in V_{-1}(\cD_X)\cdot m$ (the filtration $V_\bbullet(\cD_X)$ is defined as in \T\ref{subsec:Vfil}; we usually assume that $b_m$ is minimal for this property).

An equivalent definition is that there exists, locally on $X$, a good $V$-filtration $U_\bbullet M$ and a polynomial $\index{$bu$@$b_U$}b_U(s)\in\CC[s]\moins\{0\}$ such that, for any $k\in\ZZ$, we have
\begin{equation}\tag{$*$}\label{eq:star}
b_U(-(\partial_tt+k))\cdot\gr_k^UM=0.
\end{equation}
Indeed, in one direction, take the $V$-filtration generated by a finite number of local generators of $M$; in the other direction, use that two good filtrations are locally comparable.

If we decompose $b_U(s)$ as a product $b_1(s)b_2(s)$ then, putting
\begin{equation}\tag{$**$}\label{eq:starstar}
U'_k=U_{k-1}+b_1(-(\partial_tt+k))U_k,
\end{equation}
we get a new good $V$-filtration $U'_\bbullet$ with polynomial $b_{U'}=b_1(s-1)b_2(s)$. Therefore, an equivalent definition of specializability is that there exists, locally on $X$, a good $V$-filtration $U_\bbullet M$ and a polynomial $b_U(s)\in\CC[s]\moins\{0\}$ satisfying
\begin{equation}\label{eq:roots}
\text{the roots of $b_U$ do not differ by a nonzero integer},
\end{equation}
such that \eqref{eq:star} is satisfied.

For such a good filtration $U_\bbullet M$ and any $\alpha\in\CC$, put
\[
\psi_{t,\alpha}^UM=\ccup_n\ker[(\partial_tt+\alpha)^n:\gr^UM\to\gr^UM].
\]
We then have $\psi_{t,\alpha}^UM=0$ unless $\alpha\in b_U^{-1}(0)+\ZZ$, and
$\gr^UM=\oplus_\alpha\psi_{t,\alpha}^UM$.

For any $k\in\ZZ$, there are $\cD_{X_0}$-linear morphisms
\[
t:\gr_k^{U}M\to\gr_{k-1}^{U}M \qqbox{and} -\partial_t:\gr_k^{U}M\to\gr_{k+1}^{U}M.
\]
These morphisms are compatible with the decomposition with respect to the generalized eigenvalues of $-\partial_tt$ and induce morphisms
\[
t:\psi^U_{t,\alpha}M\to\psi^U_{t,\alpha-1}M\qqbox{and} -\partial_t:\psi^U_{t,\alpha}M\to\psi^U_{t,\alpha+1}M
\]
for any $\alpha\in\CC$. The first one is an isomorphism if $\alpha\neq0$ and the second one if $\alpha\neq-1$, as $\partial_tt$ (\resp $t\partial_t$) is invertible on $\psi^U_{t,\alpha}M$ if $\alpha\neq0$ (\resp $\alpha\neq-1$). We denote by $\can:\psi^U_{t,-1}M\to\psi^U_{t,0}M$ the morphism induced by $-\partial_t$ and by $\var:\psi^U_{t,0}M\to\psi^U_{t,-1}M$ the morphism induced by $t$.

\medskip
If $U_\bbullet M$ is any good $V$-filtration of $M$ defined on some open set of $X$, with Bernstein polynomial $b_U$, then any other good $V$-filtration $U'_\bbullet M$, defined on this open set or on any subset of it, has a Bernstein polynomial $b_{U'}$, and this polynomial satisfies $b_{U'}^{-1}(0)\subset b_U^{-1}(0)+\ZZ$. If we assume that $U_\bbullet M$ satisfies \eqref{eq:roots}, then any other good $V$-filtration $U'_\bbullet M$ defined on the same (or on a smaller) domain, and satisfying $b_{U'}^{-1}(0)\subset b_U^{-1}(0)$, is equal to $U_\bbullet M$.

Consequently, if $M$ is specializable, given any section $\sigma$ of the projection $\CC\to\CC/\ZZ$, there exists a unique good $V$-filtration $U^{\sigma}_\bbullet M$, globally defined on $X$, such that any local Bernstein polynomial $b_{U^\sigma}$ satisfies $b_{U^\sigma}^{-1}(0)\subset \image\sigma$. Any morphism between specializable $\cD_X$-modules is strictly compatible with the filtration $U^\sigma$.

Let $\ell:\CC\to\RR$ be a $\RR$-linear form such that $\ell(\ZZ)\subset\ZZ$. It defines a relation $\leq_{\ell}$ on~$\CC$: $\alpha_1\leq_{\ell}\alpha_2$ iff $\ell(\alpha_1)\leq\ell(\alpha_2)$. One usually takes $\ell(\alpha)=\reel(\alpha)$, but we will need below (see Proposition \ref{prop:restrhbo}) to consider various such linear forms.

Let $m$ be a local section of $M$. If $b_m$ is the Bernstein polynomial of $m$, we define the $\ell$-order of $m$ as $\ord_\ell(m)=\max\{\ell(\alpha)\mid b_m(\alpha)=0\}$. Define the $V$-filtration by the $\ell$-order $V_\bbullet^{(\ell)}M$ by the following property: a local section $m$ is in $V_k^{(\ell)}M$ iff $\ord_\ell(m)\leq k$. If $M$ is specializable, this filtration is \emph{good}. It is the filtration associated to the section of $\CC\to\CC/\ZZ$ which has image in $\{s\mid\ell(s)\in[0,1[{}\}$.

It will be convenient, later on, to view this filtration as indexed by $\RR$ with a discrete set of jumps, corresponding to the zeros of the possible $b_m$. Let us recall this notion. Let $A_\RR$ be a finite set in $\RR$ and put $\Lambda_\RR=A_\RR+\ZZ$. A \emph{good $V$-filtration} of $M$ indexed by $\Lambda_\RR$ is by definition a family ${}^{(a)}U_\bbullet M$ ($a\in\Lambda_\RR$) of good $V$-filtrations indexed by $\ZZ$ which satisfy the following properties:
\begin{itemize}
\item
${}^{(a)}U_\bbullet M\subset {}^{(b)}U_\bbullet M$ if $a\leq b$,
\item
${}^{(a+1)}U_\bbullet M={}^{(a)}U_{\bbullet+1} M$.
\end{itemize}
For any $a\in\Lambda_\RR$, one then defines $U_aM\defin {}^{(a)}U_0M$.
If $<a$ denotes the largest element of $\Lambda_\RR$ which is strictly smaller than $a$, then one puts $\gr_a^UM=U_aM/U_{<a}M$.

\medskip
If $U_\bbullet M$ and $U'_\bbullet M$ are two good $V$-filtrations satisfying \eqref{eq:roots}, then there are isomorphisms $\psi^U_{t,\alpha}M\isom\psi^{U'}_{t,\alpha}M$ which are compatible with $t$ and $-\partial_t$. Indeed, by the uniqueness above, $U$ and $U'$ may be related by a finite sequence of transformations of type \eqref{eq:starstar} for which, at each step, $b_1$ and $b_2$ do not have any common root. It is thus enough to prove the assertion when $U$ and $U'$ are as in \eqref{eq:starstar}, and $\gcd(b_1,b_2)=1$.

In such a situation, we have an exact sequence
\[
0\to U_{k-1}/U'_{k-1}\to U'_k/U'_{k-1}\to U'_k/U_{k-1}\to 0.
\]
On the one hand, the natural morphism $U'_k/U_{k-1}\to U_k/U_{k-1}$ is injective with image equal to $\ker b_2(-\partial_tt+k)$, as $\gcd(b_1,b_2)=1$. On the other hand, $U_{k-1}/U_{k-2}\to U_{k-1}/U'_{k-1}$ is onto and induces an isomorphism $\ker b_1(-\partial_tt+k-1)\isom U_{k-1}/U'_{k-1}$. The assertion follows.

As a consequence, taking a section $\sigma$ as above, the modules $\psi_{t,\alpha}^{U^\sigma}M$ are globally defined, and are independent of $\sigma$ up to a canonical isomorphism. They are equipped with the action of a nilpotent operator, locally obtained as the action of $-\partial_tt$. We denote them simply by $\psi_{t,\alpha}M$. Notice however that, to define $\can$ and $\var$, one needs an equation $\{t=0\}$ for $X_0$ and a corresponding vector field $\partial_t$.

\medskip
Any coherent sub or quotient module of a specializable $\cD_X$-module is so. For a specializable $\cD_X$-module,
\begin{enumerate}
\item\label{item:1}
$\can:\psi_{t,-1}M\to\psi_{t,0}M$ is onto iff $M$ has no coherent quotient $\cD_X$-module supported on $X_0$,
\item\label{item:2}
$\var:\psi_{t,0}M\to\psi_{t,-1}M$ is injective iff $M$ has no coherent sub $\cD_X$-module supported on $X_0$,
\item\label{item:3}
$\psi_{t,0}M=\im\can\oplus\ker\var$ iff $M=M'\oplus M''$ with $M'$ satisfying \ref{item:1} and \ref{item:2} and $M''$ supported on $X_0$.
\end{enumerate}

\section{The category $\cS^2(X,t)$}\label{subsec:S2t}
\subsection{}
Keep notation of \T\ref{subsec:Vfil}. We will work with increasing filtrations. To get a decreasing filtration $U^{\cbbullet}$ from an increasing one $U_\bbullet$, put $U^\beta=U_{-\beta-1}$ (see Remark \ref{rem:leftrightspe}\eqref{rem:leftrightspe2}).

We will introduce the Malgrange-Kashiwara filtration in the setting of $\cR_\cX$-modules. When the set $A$ below is contained in $\RR$, the presentation can be simplified, as the Malgrange-Kashiwara filtration is then defined \emph{globally} with respect to $\hb$, and not only locally. For $A\subset\CC$ general, the definitions below are suggested by Corollary \ref{cor:biquardRX}.

The strictness assumption is important, as emphasized yet in Theorem \ref{th:imdirspe}: for Hodge modules, it means a good behaviour of the Hodge filtration under the operation of taking nearby or vanishing cycles. Moreover, it important to notice that, under a strictness assumption, the ``nearby cycles'' $\psi_{t,\alpha}\cM$ are defined \emph{globally} with respect to~$\hb$.

\begin{definition}\label{def:spe}
A coherent left $\cR_\cX$-module $\cM$ is said to be \emph{specializable} along $\cX_0$ if there exists, locally on $X$, a finite subset $A\subset\CC$ and for any local section $m$ of~$\cM$, there exists a polynomial $\index{$bm$@$b_m$}b_m(s)=\prod_{\alpha\in A}\prod_{\ell\in\ZZ} (s-(\alpha+\ell)\star\hb)^{\nu_{\alpha,\ell}}$ satisfying
\[
\begin{array}{rcll}
b_m(-\partiall_tt)\cdot m&\in& V_{-1}\cR_\cX\cdot m&\text{(for left modules)},\\
m\cdot b_m(t\partiall_t)&\in& m\cdot V_{-1}\cR_\cX&\text{(for right modules)}.
\end{array}
\]
\end{definition}

An equivalent definition is that there exists, locally on $\cX$, a good $V$-filtration $U_\bbullet\cM$ and a polynomial $\index{$bu$@$b_U$}b_U(s)$ of the same kind, such that, for any $k\in\ZZ$, we have
\begin{equation}\tag{$*$}\label{eq:spe*}
b_U(-(\partiall_tt+k\hb))\cdot\gr_k^U\cM=0,\qquad \text{\resp\quad }\gr_k^U\cM\cdot b_U(t\partiall_t-k\hb)=0.
\end{equation}
Indeed, in one direction, take the $V$-filtration generated by a finite number of local generators of $\cM$; in the other direction, use that two good filtrations are locally comparable.

Remark that any coherent $\cR_\cX$-submodule or quotient of $\cM$ is specializable if $\cM$ is so (take the induced filtration, which is good, by Proposition \ref{prop:V-coh}). In particular, the category of specializable $\cR_\cX$-modules is abelian.

If we decompose $b_U(s)$ as a product $b_1(s)b_2(s)$ then, putting
\begin{equation}\tag{$**$}\label{eq:spe**}
U'_k=U_{k-1}+b_1(-(\partiall_tt+k\hb))U_k,
\end{equation}
we get a new good $V$-filtration $U_\bbullet$ with polynomial $b_{U'}=b_1(s-1)b_2(s)$. Therefore, an equivalent definition of specializability is that there exists locally a good $V$-filtration~$U$ with polynomial $b_U(s)=\prod_{\alpha\in A}(s-\alpha\star\hb)^{\nu_\alpha}$ for some integers $\nu_\alpha$ and $A$ is contained in the image of a section $\sigma$ of the projection $\CC\to\CC/\ZZ$. In other words, $(\cM,U_\bbullet\cM)$ is monodromic (\cf Definition \ref{def:monodromic}) with a particular form for $b$.

The constructions made in \T\ref{sec:reviewspe} may be applied word for word here, \emph{provided that we avoid singular points with respect to $\Lambda\defin A+\ZZ$} (\cf \T\ref{subsec:defstar}). Indeed, we need such a nonsingularity assumption to get that $\gcd(b_1,b_2)$ is invertible when $b_1=0$ and $b_2=0$ do not have common components. In the neighbourhood of such singular points (in particular in the neighbourhood of $0$), we will need the constructions of the next subsections.

\medskip
The choice of the generating set $A$ may be changed. Put $\Lambda=A+\ZZ$ and, for any $\hb_o\in\Omega_0$ and $\alpha\in\Lambda$, set $\ell_{\hb_o}(\alpha)=\alpha'-\imhb_o\alpha''$ (recall that $\imhb_o\defin\im\hb_o$, so that $\ell_{\hb_o}(\alpha)=\reel(\alpha'+i\hb_o\alpha'')$; remark also that $\ell_{\hb_o}(\alpha+1)=\ell_{\hb_o}(\alpha)+1$, \cf\T\ref{subsec:defstar}); set also $\Lambda(\hb_o)=\{\alpha\in\Lambda\mid \ell_{\hb_o}(\alpha)\in[0,1[\}$. Then, the specializability of $\cM$ is equivalent to the local existence of a good $V$-filtration $U^{(\hb_o)}_\bbullet\cM$ with polynomial $b_{U^{(\hb_o)}}(s)$ having roots in $\Lambda(\hb_o)$.

\subsection{}\label{subsec:IR}
It will be convenient to work with filtrations indexed by $I_\RR+\ZZ$ for some finite set $I_\RR\subset\RR$, that we now define. Let $\cM$ be a coherent $\cR_\cX$-module and $I_\RR\subset\RR$ be a finite set contained in the image of a section of $\RR\to\RR/\ZZ$. A \emph{good filtration} of $\cM$ indexed by $I_\RR+\ZZ$ consists of a family ($a\in I_\RR$) of filtrations $\aU_{\bbullet}\cM$ indexed by $\ZZ$, which are good with respect to $V_\bbullet\cR_\cX$ and such that one has
\begin{equation}\label{eq:boite}
a+k\leq b+\ell\implique \aU_{k}\cM\subset {}^{(b)}\!U_{\ell}\cM.
\end{equation}
We denote $U_{a+k}\cM=\aU_{k}\cM$ and $\gr_{a}^{U}\cM=U_a\cM/U_{<a}\cM$. We may also view $U_\bbullet\cM$ as a filtration indexed by $\RR$ with jumps at $I_\RR+\ZZ$ at most.

Saying that $\cM$ is specializable is then equivalent to saying that, near any $(x_o,\hb_o)\in\cX$, it has a good filtration indexed by $\ell_{\hb_o}(\Lambda)=\ell_{\hb_o}(\Lambda(\hb_o))+\ZZ$ such that, for any $a\in\RR$,
\begin{equation}\label{eq:hospe}
\prod_{\substack{\alpha\in\Lambda\\ \ell_{\hb_o}(\alpha)=a}} (-\partiall_tt-\alpha\star\hb)^{\nu_\alpha}\cdot\gr_a^{U^{(\hb_o)}}\cM=0,
\end{equation}
where the integers $\nu_\alpha$ only depend on $\alpha\mod\ZZ$. Remark that the set of indices (hence the order of the filtration) depends on the point $\hb_0$. This is suggested by Corollary \ref{cor:biquard} below.

\begin{lemme}\label{lem:uniciteMK}
Assume that $\cM$ has, in the neighbourhood of any point $(x_o,\hb_o)\in\cX_0$, a good filtration $U^{(\hb_o)}_\bbullet\cM$ indexed by $\ell_{\hb_o}(\Lambda)$, satisfying {\rm \eqref{eq:hospe}} and such that $\gr_{a}^{U^{(\hb_o)}}\cM$ is a strict $\cR_{\cX_0}$-module for any $a\in \RR$. Then,
\begin{enumerate}
\item\label{lem:uniciteMK1}
for any coherent submodule $\cN\subset\cM$, the filtration $U^{(\hb_o)}_\bbullet\cN\defin U^{(\hb_o)}_\bbullet\cM\cap\cN$ is a good filtration satisfying the same properties;
\item\label{lem:uniciteMK2}
such a filtration is unique; it is therefore globally defined on some neighbourhood of $X\times (\RR+i\imhb_o)$; it will be denoted by $\index{$vhb$@$V^{(\hb_o)}_\bbullet\cM$}V^{(\hb_o)}_\bbullet\cM$ and be called \emph{the Malgrange-Kashiwara filtration} of $\cM$ along $X\times (\RR+i\imhb_o)$;

\item\label{lem:uniciteMK3}
this filtration is equal to the filtration by the $\ell_{\hb_o}$-order along $\cX_0$;

\item\label{lem:uniciteMK4}
any morphism $\varphi:\cM\to\cN$ between such $\cR_\cX$-modules is compatible with the $V^{(\hb_o)}$-filtration;

\item\label{lem:uniciteMK5}
the construction of the $V^{(\hb)}$-filtration is locally constant with respect to $\hb$;

\item\label{lem:uniciteMK6}
near any $\hb_o\in\Omega_0$ one has, for any $a\in\ell_{\hb_o}(\Lambda)$,
$$
\gr_a^{V^{(\hb_o)}}\cM=\ooplus_{\ell_{\hb_o}(\alpha)=a}\psi^{(\hb_o)}_{t,\alpha}\cM,$$
with
$$
\index{$psihb$@$\psi^{(\hb_o)}_{t,\alpha}\cM$, $\psi_{t,\alpha}\cM$}\psi^{(\hb_o)}_{t,\alpha}\cM\defin\textstyle\bigcup_n \ker\big[(\partiall_tt+\alpha\star\hb)^n:\gr_{\ell_{\hb_o}(\alpha)}^{V^{(\hb_o)}}\cM\to\gr_{\ell_{\hb_o}(\alpha)}^{V^{(\hb_o)}}\cM\big];
$$

\item\label{lem:uniciteMK7}
for any $\alpha\in\Lambda$, there exists a strict coherent $\cR_{\cX_0}$-module $\psi_{t,\alpha}\cM$ equipped with a nilpotent endomorphism $\rN$ such that, near any $\hb_o\in\Omega_0$,
\[
\big(\psi_{t,\alpha}\cM,\rN\big)\simeq
\big(\psi^{(\hb_o)}_{t,\alpha}\cM,-(\partiall_tt+\alpha\star\hb)\big).
\]
\end{enumerate}
\end{lemme}

\begin{Proof}
\eqref{lem:uniciteMK1} The filtration $U^{(\hb_o)}_\bbullet\cN\defin U^{(\hb_o)}_\bbullet \cM\cap \cN$ is good (by Proposition \ref{prop:V-coh}) and clearly satisfies \eqref{eq:hospe}. For each $a\in \ell_{\hb_o}(\Lambda)$, $\gr_{a}^{U^{(\hb_o)}}\cN$ is a submodule of $\gr_{a}^{U^{(\hb_o)}}\cM$, hence is also strict.

\medskip
\eqref{lem:uniciteMK2} Let $U^{(\hb_o)}$ and $U^{\prime(\hb_o)}$ be two such filtrations on $\cM$, that we may assume to be indexed by the same set $\ell_{\hb_o}(\Lambda)$. Locally, there exists $\ell\in \NN$ such that, for all $a\in \ell_{\hb_o}(\Lambda)$, one has
$$
U^{\prime(\hb_o)}_{a-\ell}\subset U^{(\hb_o)}_{a} \subset U^{\prime(\hb_o)}_{a+\ell}\subset U^{(\hb_o)}_{a+2\ell}.
$$
Let $m$ be a local section of $U^{(\hb_o)}_a\cM$ and assume that there exists $b\in{}]a, a+\ell]$ such that $m\in U^{\prime(\hb_o)}_b\moins U^{\prime(\hb_o)}_{<b}$. Then, there exist polynomials $B_U(s)$ and $B_{U'}(s)$, where $B_U$ has roots $\alpha\star\hb$ with $\ell_{\hb_o}(\alpha)\leq a$ and $B_{U'}$ has roots $\alpha\star\hb$ with $\ell_{\hb_o}(\alpha)=b$, such that $B_U(-\partiall_tt)\cdot m\in U^{(\hb_o)}_{<b-\ell}\subset U^{\prime(\hb_o)}_{<b}\cM$ and $B_{U'}(-\partiall_tt)\cdot m\in U^{\prime(\hb_o)}_{<b}\cM$. Hence, there exists $p(\hb)\in\CC[\hb]\moins\{0\}$ such that $p(\hb)m\in U^{\prime(\hb_o)}_{<b}$. As $\gr_{b}^{U^{\prime(\hb_o)}}\cM$ is strict, this implies that $m\in U^{\prime(\hb_o)}_{<b}$, a contradiction. This shows that $m\in U^{\prime(\hb_o)}_a$. Exchanging the roles of $U$ and $U'$ gives the uniqueness.

\medskip
\eqref{lem:uniciteMK3} Let us first define the filtration by the $\ell_{\hb_o}$-order for a filtered module satisfying \eqref{eq:hospe}. Let $(x_o,\hb_o)\in \cX_0$ and $m\in\cM_{(x_o,\hb_o)}$. According to the proof of \eqref{lem:uniciteMK1}, we will assume that, locally at $(x_o,\hb_o)$, the section $m$ generates $\cM$.

There exists a minimal polynomial $b_m(s)$ as in Definition \ref{def:spe} such that a relation $b_m(-\partiall_tt)\cdot m\in V_{-1}\cR_\cX\cdot m$ is satisfied. The \emph{$\ell_{\hb_o}$-order} of $m$ along $\cX_0$ is the largest $a\in\RR$ for which there exists $\alpha$ with $\ell_{\hb_o}(\alpha)=a$, such that $\alpha\star\hb$ is a root of $b_m$. This defines an increasing filtration of $\cM$, called \emph{the filtration by the order}. It is not necessarily good \emph{a~priori}. Denote it by $V^{\prime(\hb_o)}_\bbullet\cM$ and denote by $V^{(\hb_o)}_\bbullet\cM$ the filtration obtained in \eqref{lem:uniciteMK2}.

We clearly have $V^{(\hb_o)}_a\cM\subset V^{\prime(\hb_o)}_a\cM$ for each $a\in \RR$. Let us prove the reverse inclusion. Assume that $m$ has $\ell_{\hb_o}$-order $a$. If $m\in V^{(\hb_o)}_b\cM-V^{(\hb_o)}_{<b}\cM$ with $b\geq a$, there exists some integer $\ell\geq 0$ such that
\[
b_m(-\partiall_tt)\cdots b_m(-\partiall_tt+\ell\hb)\cdot m\in V^{(\hb_o)}_{<b}\cM.
\]
We also have
\[
\prod_{\ell_{\hb_o}(\alpha)=b}(-\partiall_tt-\alpha\star\hb)^{\nu_\alpha}\cdot m\in V^{(\hb_o)}_{<b}\cM,
\]
for some $\nu_\alpha\geq 1$. If $b>a$, we deduce from both relations that there exists $p(\hb)\in\CC[\hb]\moins\{0\}$ such that $p(\hb)\cdot m\in V^{(\hb_o)}_{<b}\cM$ and, by strictness of $\gr_b^{V^{(\hb_o)}}\cM$, we have $m\in V^{(\hb_o)}_{<b}\cM$, a contradiction.

\medskip
\eqref{lem:uniciteMK4}
It is clear that any morphism between two $\cR_\cX$-modules satisfying \eqref{eq:hospe} is compatible with the filtration by the $\ell_{\hb_o}$-order. The assertion follows then from \eqref{lem:uniciteMK3}.

\medskip
\eqref{lem:uniciteMK5}
We fix a compact set in $X$. The constants $\epsilon,\eta$ below will be relative to this compact set. We say that the filtration $V^{(\hb)}_\bbullet\cM$ is locally constant at $\hb_o$ with respect to $\hb$ if, for any $a\in\RR$ and any $\varepsilon>0$ sufficiently small, there exists $\eta=\eta(\hb_o,a,\varepsilon)$ such that, for any $\hb\in\Delta_{\hb_o}(\eta)$ (disc of radius $\eta$ centered at $\hb_o$), we have
\begin{equation}\label{eq:Vlocconst}
V^{(\hb)}_{a-\varepsilon}(\cM_\hb)=V^{(\hb_o)}_{a-\varepsilon}(\cM)_\hb,\quad V^{(\hb)}_{a+\varepsilon}(\cM_\hb)=V^{(\hb_o)}_{a+\varepsilon}(\cM)_\hb.
\end{equation}
By considering Fig\ptbl \ref{fig:lines}, one shows that this property is true for the filtration by the $\ell_{\hb}$-order, hence the result by \eqref{lem:uniciteMK3}.
\begin{figure}[htb]
\setlength{\unitlength}{1mm}
\begin{center}
\begin{picture}(60,60)(0,0)
\put(0,15){\vector(1,0){60}}
\put(20,0){\vector(0,1){60}}
\put(59,11.5){$\imhb$}
\put(20,50){\circle*{1}}\put(17,50){$1$}
\put(20,37){\circle*{1}}
\put(20,34){\circle*{1}}
\put(20,31){\circle*{1}}

\put(14,38){$\scriptstyle a+\varepsilon$}
\put(18,32){$\scriptstyle a$}
\put(14,29){$\scriptstyle a-\varepsilon$}

\put(0,10){\line(1,1){50}}
\put(0,30){\line(3,-2){45}}
\put(20,51.7){\line(3,-2){45}}\put(20,51.7){\line(-3,2){11}}
\put(0,40){\line(4,-1){65}}

\multiput(-1,37)(4,0){16}{\line(1,0){2}}
\multiput(-1,34)(4,0){16}{\line(1,0){2}}
\multiput(-1,31)(4,0){16}{\line(1,0){2}}

\multiput(24.1,15)(0,1){19}{\circle*{.5}}
\put(24.1,15){\circle*{1}}\put(24.5,16.5){$\imhb_o$}

\multiput(27,38)(0,1){7}{\circle*{.5}}
\multiput(21.5,31)(0,1){14}{\circle*{.5}}
\put(21.5,45){\vector(1,0){5.5}}
\put(27,45){\vector(-1,0){5.5}}
\put(23.5,42.5){$\eta$}

\put(45,52){$\imhb\mto\ell_{\hb}(\alpha_1)$}
\put(40,23){$\imhb\mto\ell_{\hb}(\alpha_2)$}
\put(40,5){$\imhb\mto\ell_{\hb}(\alpha_3)$}
\put(-11,54){$\imhb\mto\ell_{\hb}(\alpha_3+1)$}
\end{picture}
\caption{The set $\{\alpha\in \Lambda\mid\ell_{\hb_o}(\alpha)=a\}$ is equal to $\{\alpha_1,\alpha_2\}$.}\label{fig:lines}
\end{center}
\end{figure}
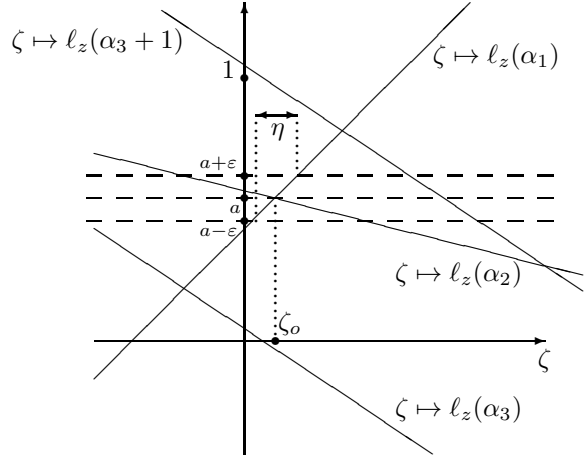

\medskip
\eqref{lem:uniciteMK6} Locally on $\cX$, there exists $\nu_\alpha$ such that $\psi_{t,\alpha}^{(\hb_o)}\cM=\ker (\partiall_tt+\alpha\star\hb)^{\nu_\alpha}$. The roots of the minimal polynomial of $-\partiall_tt$ on $\gr_a^{V^{(\hb_o)}}\cM$ are the $\alpha\star\hb$ for $\alpha$ such that $\ell_{\hb_o}(\alpha)=a$. By Lemma \ref{lem:imhb}, the various $\alpha\star\hb_o$ corresponding to distinct $\alpha$ are distinct, hence, if $\hb$ is sufficiently close to $\hb_o$, the roots $\alpha\star\hb$ remain distinct. Therefore one has a B\'ezout relation between the polynomials $(\partiall_tt+\alpha\star\hb)^{\nu_\alpha}$ for $\alpha\in\ell_{\hb_o}^{-1}(a)$, and the decomposition follows.

\medskip
\eqref{lem:uniciteMK7} Fix $\hb_o\in\Omega_0$, set $a=\ell_{\hb_o}(\alpha)$ and fix $\varepsilon,\eta>0$ as in \eqref{lem:uniciteMK5} so that, moreover, for any $\hb\in\Delta_{\hb_o}(\eta)$, we have $\ell_\hb(\alpha)\in{}]a-\varepsilon,a+\varepsilon]$. If $\epsilon$ is small enough, we have $\gr_a^{V^{(\hb_o)}}\cM_\hb=V^{(\hb_o)}_{a+\epsilon}\cM_\hb/V^{(\hb_o)}_{a-\epsilon}\cM_\hb$ for any $\hb\in\Delta_{\hb_o}(\eta)$. By~\eqref{lem:uniciteMK5}, we can compute $\psi_{t,\alpha}^{(\hb_o)}(\cM)_\hb$ using the filtration $V^{(\hb)}_\bbullet\cM_\hb$:
\begin{multline*}
\psi_{t,\alpha}^{(\hb_o)}(\cM)_\hb=\ker\big[-(\partiall_tt+\alpha\star\hb)^{\nu_\alpha}:V_{a+\varepsilon}^{(\hb)}(\cM_\hb)\big/V_{a-\varepsilon}^{(\hb)}(\cM_\hb)\\
\to V_{a+\varepsilon}^{(\hb)}(\cM_\hb)\big/V_{a-\varepsilon}^{(\hb)}(\cM_\hb)\big].
\end{multline*}
Fix $\hb\in\Delta_{\hb_o}(\eta)$. By construction, if $\beta\in\Lambda$ is such that $\ell_\hb(\beta)\in{}]a-\epsilon,a+\epsilon]$, then $\ell_{\hb_o}(\beta)=\ell_{\hb_o}(\alpha)$ (on Fig\ptbl \ref{fig:lines}, each such $\beta$ corresponds to a line cutting the intersection of the horizontal strip of width $2\epsilon$ and the vertical strip of width $\eta$). By Lemma \ref{lem:imhb}, if $\beta\neq\alpha$, we then have $\beta\star\hb_o\neq\alpha\star\hb_o$ and this remains true for any $\hb$ near $\hb_o$, so that $\hb\mto(\beta-\alpha)\star\hb$ is locally invertible. For $b\in{}]a-\epsilon,a+\epsilon]$ and $b\neq\ell_\hb(\alpha)$, the operator $-(\partiall_tt+\alpha\star\hb)$ acting on $\gr_b^{V^{(\hb)}}\cM_\hb=\oplus_{\beta\mid\ell_\hb(\beta)=b}\psi^{(\hb)}_{t,\beta}\cM_\hb$ (see \eqref{lem:uniciteMK6} above for the decomposition) may be written, on each summand, as the sum of the nilpotent operator $-(\partiall_tt+\beta\star\hb)$ and of the invertible operator $[(\beta-\alpha)\star\hb]\id$. It is therefore invertible. This implies that
\[
\psi_{t,\alpha}^{(\hb_o)}(\cM)_\hb\subset V_{\ell_{\hb}(\alpha)}^{(\hb)}(\cM_\hb)\big/V_{a-\varepsilon}^{(\hb)}(\cM_\hb)\subset V_{a+\varepsilon}^{(\hb)}(\cM_\hb)\big/V_{a-\varepsilon}^{(\hb)}(\cM_\hb),
\]
and that the natural projection (which is now defined) from $\psi_{t,\alpha}^{(\hb_o)}(\cM)_\hb$ to $\gr_{\ell_\hb(\alpha)}^{V^{(\hb)}}\cM_\hb$ is an isomorphism $\lambda_{\hb,\hb_o}$ onto $\psi_{t,\alpha}^{(\hb)}(\cM_\hb)$.

In order to show that the various $\psi_{t,\alpha}^{(\hb)}(\cM)$ glue as a sheaf $\psi_{t,\alpha}(\cM)$ on $\cX_0$, it is enough to show that the family of isomorphisms $\lambda_{\hb',\hb}$ is transitive, \ie $\lambda_{\hb'',\hb}=\lambda_{\hb'',\hb'}\circ\lambda_{\hb',\hb}$ whenever $\hb'\in\Delta_\hb$ and $\hb''\in\Delta_{\hb'}\cap\Delta_{\hb}$. This follows from \eqref{eq:Vlocconst} which holds for $\hb,\hb',\hb''$.
\end{Proof}

\begin{remarques}\label{rem:psi}
The strictness assumption in Lemma \ref{lem:uniciteMK} is used only near the points of $\Sing\Lambda(\cM)$. More precisely, the lemma holds without any strictness assumption if we extend the coefficients by tensoring over $\cO_{\Omega_0}$ by the ring $\cO_{\Omega_0}(*\Sing\Lambda(\cM))$ obtained by inverting polynomials vanishing on $\Sing\Lambda$. If we allow to divide by such polynomials, the proof gives the uniqueness of the Malgrange-Kashiwara filtration, as well as its coincidence with the filtration by the order.

\smallskip
Let $\cM,\cN$ be as in Lemma \ref{lem:uniciteMK}.
\begin{enumerate}
\item\label{rem:psi0}
If $\varphi:\cM\to\cN$ is any morphism, then $\gr_{\ell_{\hb_o}(\alpha)}^{V^{(\hb_o)}}\varphi$ sends $\psi_{t,\alpha}\cM_{\hb_o}$ in $\psi_{t,\alpha}\cN_{\hb_o}$. It is globally defined. Denote by $\psi_{t,\alpha}\varphi:\psi_{t,\alpha}\cM\to\psi_{t,\alpha}\cN$ the morphism induced locally by $\gr_{\ell_{\hb_o}(\alpha)}^{V^{(\hb_o)}}\varphi$.

\item\label{rem:psi1}
What prevents the filtration $V_{\cbbullet}^{(\hb)}$ to be defined independently of $\hb$ is that $V_a^{(\hb)}(\cM_{\hb})$ jumps for values of $a$ depending on $\imhb$, more precisely for $a=\alpha'-\imhb\alpha''$, $\alpha\in \Lambda$. If $\Lambda\subset\RR$ (\eg if the local monodromy of the perverse sheaf corresponding to $\Xi_{\DR}(\cM)$ around $0$ is quasi-unipotent or unitary), then the $V$-filtration is globally defined on $\cM$. Otherwise, the natural order on the set $\ell_{\hb}(\Lambda)$ may depend on $\hb$, as shown on Fig\ptbl \ref{fig:lines}. One may compare the various $V_a^{(\hb)}(\cM_{\hb})$ for $\hb$ in a neighbourhood of $\hb_o$ if no jump occurs at $a$ for these germs. This explains that, given any $a$, we compare the various $V_{a\pm\varepsilon}^{(\hb)}(\cM_\hb)$ (see Fig\ptbl \ref{fig:lines}).

If $\Lambda\subset\RR$, then the previous lemma is simpler, as $\ell_{\hb_o}(\alpha)=\alpha$ for any $\alpha\in\Lambda$. The order does not depend on $\hb_o$ and the $V$-filtration is globally defined. Moreover, $\psi_{t,\alpha}\cM=\gr_\alpha^V\cM$ for any $\alpha\in\Lambda$.

In Corollary~\ref{cor:biquard}, we will view $\ell_\hb(\alpha)$ as a growth order indexing a parabolic filtration. This growth order depends on $\imhb$. A similar problem occurs when defining the Stokes structure of a meromorphic connection of one variable, as for instance in \cite[p\ptbl 56]{Malgrange91}.
\item\label{rem:Ncanvar1}
The nilpotent endomorphism
\[
-(\partiall_tt+\alpha\star\hb):\psi_{t,\alpha}\cM\to\psi_{t,\alpha}\cM
\]
is denoted by $\index{$n$@$\rN$}\rN$. There exists a unique increasing filtration $\index{$mn$@$\rM(\rN)_\bbullet$}\rM(\rN)_\bbullet$ of $\psi_{t,\alpha}\cM$ by $\cR_{\cX_0}$-submodules, indexed by $\ZZ$, such that, for any $\ell\geq 0$, $\rN$ maps $\rM_k$ into $\rM_{k-2}$ for all $k$ and $\rN^\ell$ induces an isomorphism $\gr_\ell^{\rM}\isom\gr_{-\ell}^{\rM}$ for any $\ell\geq 0$. It is called the \emph{monodromy filtration} of $\rN$ (\cf \cite[\T1.6]{Deligne80}). Each $\gr^{\rM}\psi_{t,\alpha}\cM$ has a Lefschetz decomposition, with basic pieces the primitive parts ($\ell\geq 0$)
\[
\index{$pgrm$@$P\gr_\ell^{\rM}$}P\gr_\ell^{\rM}\psi_{t,\alpha}\cM\defin\ker\lefcro \rN^{\ell+1}:\gr_\ell^{\rM}\psi_{t,\alpha}\cM\to\gr_{-\ell-2}^{\rM}\psi_{t,\alpha}\cM\rigcro.
\]
Notice however that $\gr_\ell^{\rM}\psi_{t,\alpha}\cM$ need not be strict, even if $\psi_{t,\alpha}\cM$ is so.

\item\label{rem:psia}
By an argument similar to that of the lemma, one shows that $t:\psi_{t,\alpha}\to\psi_{t,\alpha-1}\cM$ is injective for $\alpha\neq0$ and $\partiall_t:\psi_{t,\alpha}\cM\to\psi_{t,\alpha+1}\cM$ is so for $\alpha\neq-1$. Notice that, by Lemma \ref{lem:tinj}, we also get that $t:V^{(\hb_o)}_{a}\cM\to V^{(\hb_o)}_{a-1}\cM$ is injective for $a<0$, because $V^{(\hb_o)}_\bbullet\cM$ is good.

Notice also that, for any $\alpha$ with $\reel\alpha\neq0$, the morphism $t:\psi_{t,\alpha}\cMS\to\psi_{t,\alpha-1}\cMS$ is an isomorphism. Indeed, we already know that $t$ is injective. Moreover, as $\reel\alpha\neq0$, we have $\alpha\star\hb/\hb\neq0$ for any $\hb\in\bS$, and therefore $t\partiall_t$ is invertible on $\psi_{t,\alpha-1}\cMS$, hence the surjectivity of $t$.

\item\label{rem:psib}
Recall that if $\hb_o\not\in i\RR$, then $\alpha\star\hb_o=0\implique\alpha=0$. For such a $\hb_o$ and for $\alpha\neq0$, $\partiall_tt$ is therefore invertible on $\psi_{t,\alpha}\cM$ in a neighbourhood of $\hb_o$. This implies that, in a neighbourhood of such a $\hb_o$, we may replace ``injective'' with ``isomorphism'' in \eqref{rem:psia} above.

\item\label{rem:Ncanvar2}
There is a diagram of morphisms
$$
\xymatrix{
\psi_{t,-1}\cM\ar@/^/[rr]^-{\index{$can$@$\can$}\can=-\partiall_t}&&\psi_{t,0}\cM\ar@/^/[ll]^-{\index{$var$@$\var$}\var=t}
}
$$
where $\var$ is induced by the action of $t$ and $\can$ by that of $\partiall_t$.

\item
For right $\cR_\cX$-modules, $\rN$ is defined as the endomorphism induced by the right action of $t\partiall_t-\alpha\star\hb$, $\can$ by that of $\partiall_t$ and $\var$ by that of $t$.
\end{enumerate}
\end{remarques}

\begin{corollaire}\label{cor:Bernstein}
Assume that $\cM$ is as in Lemma \ref{lem:uniciteMK} and let $\hb_o\in\Omega_0$. Let $K$ be a compact set of $X$ and let $\Omega$ be an open neighbourhood of $\hb_o$ such that $V_\bbullet^{(\hb_o)}\cM$ exists on $K\times \Omega$. Let $W\subset K$ be an open set. If $m\in\Gamma(W\times\Omega,\cM)$, there exists a finite set $A(m)\subset \CC$ and $\nu:A(m)\to\NN^*$ such that
\begin{enumerate}
\item
$\gamma\in A(m)\implique\psi_{t,\gamma}\cM_{|W_0\times\Omega}\neq0$,
\item
$\nu(\gamma)\leq{}$the order of nilpotency of $\rN$ on $\psi_{t,\gamma}\cM_{|W_0\times\Omega}$,
\end{enumerate}
so that, putting $b_m(s)=\prod_{\gamma\in A(m)}(s-\gamma\star\hb)^{\nu(\gamma)}$, one has $b_m(-\partiall_tt)\cdot m\in V_{-1}\cR_\cX\cdot m$.
\end{corollaire}

In other words, a Bernstein polynomial exists in a more global setting than in Definition~\ref{def:spe}.

\begin{proof}
One can assume that $\cM=\cR_\cX\cdot m$, according to Lemma \ref{lem:uniciteMK}\eqref{lem:uniciteMK1}. Use then that the good filtration $V_\bbullet\cR_\cX\cdot m$ is comparable to the Malgrange-Kashiwara filtration $V_\bbullet^{(\hb_o)}\cM$ on $K\times\Omega$.
\end{proof}

\pagebreak[2]
\begin{Definitions}[Strict specializability]\label{def:strictspe}
\begin{enumerate}
\item
A specializable $\cR_\cX$-module is said to be \emph{strictly specializable} along $\cX_0$ if one can find, locally near any point $(x_o,\hb_o)\in\cX_0$, a good filtration $V^{(\hb_o)}_\bbullet\cM$ satisfying \eqref{eq:hospe} and such that moreover
\begin{enumerate}
\item\label{def:strictspea}
for every $a\in \RR$, $\gr_{a}^{V^{(\hb_o)}}\cM$ is a strict $\cR_{\cX_0}$-module (hence $V^{(\hb_o)}_\bbullet\cM$ is the Malgrange-Kashiwara filtration of $\cM$ near $\hb_o$);
\item\label{def:strictspeb}
$t:\psi_{t,\alpha}\cM\to\psi_{t,\alpha-1}\cM$ is onto for $\ell_{\hb_o}(\alpha)<0$;
\item\label{def:strictspec}
$\partiall_t:\psi_{t,\alpha}\cM\to\psi_{t,\alpha+1}\cM$ is onto for $\ell_{\hb_o}(\alpha)\geq-1$ but $\alpha\neq-1$.
\end{enumerate}

\item\label{def:Ncanvar3}
A morphism $\varphi:\cM\to\cN$ between two strictly specializable $\cR_\cX$-modules is \emph{strictly specializable} if, for any $\alpha\in\Lambda$, the morphisms $\psi_{t,\alpha}\varphi$ are strict.

\item\label{def:Ncanvar4}
The category $\index{$scurl2xt$@$\cS^2(X,t)$}\cS^2(X,t)$ has strictly specializable $\cR_\cX$-modules as objects and strictly specializable morphisms as morphisms.

\item\label{def:Ncanvar5}
Let $f:X\to\CC$ be an analytic function and let $\cM$ be a $\cR_\cX$-module. Denote by $i_{f}:X\hto X\times \CC$ the graph inclusion. We say that $\cM$ is strictly specializable along $f=0$ if $i_{f,+}\cM$ is strictly specializable along $X\times \{0\}$. We then set $\index{$psif$@$\psi_{f,\alpha}\cM$}\psi_{f,\alpha}\cM=\psi_{t,\alpha}(i_{f,+}\cM)$. These are coherent $\cR_\cX$-modules. If $f=t$ is induced by a projection, we have, by an easy verification, $\psi_{t,\alpha}(i_{f,+}\cM)=i_{f,+}(\psi_{f,\alpha}\cM)$ for any $\alpha$.
\end{enumerate}
\end{Definitions}

\begin{Remarques}\label{rem:strictspe}
\begin{enumerate}
\item\label{rem:strictspe0}
As we have seen in Lemma \ref{lem:uniciteMK}\eqref{lem:uniciteMK6}, Condition \eqref{def:strictspea} implies that, for every $a\in \RR$, we have a local decomposition $\gr_{a}^{V^{(\hb_o)}}\cM=\oplus_{\alpha\in \Lambda,\; \ell_{\hb_o}(\alpha)=a}\psi_{t,\alpha}\cM$.

\item\label{rem:strictspe2}
Notice that, according to \ref{rem:psi}\eqref{rem:psia}, we may replace ``onto'' with ``an isomorphism'' in \ref{def:strictspe}\eqref{def:strictspeb} and \eqref{def:strictspec}. Therefore, there is, for any $k\in\ZZ$ and near $\hb_o$, a preferred isomorphism $\psi_{t,\alpha}\cM\isom\psi_{t,\alpha+k}\cM$, obtained by suitably composing \ref{def:strictspe}\eqref{def:strictspeb} and \ref{def:strictspe}\eqref{def:strictspec}; this isomorphism is \emph{not} globally defined with respect to $\hb$, unless $\alpha$ is \emph{real}.

Moreover, locally, $t:V^{(\hb_o)}_a\cM\to V^{(\hb_o)}_{a-1}\cM$ is an isomorphism for $a<0$: indeed, it is an isomorphism for $a\ll0$, as the filtration is $V$-good; apply then \ref{def:strictspe}\eqref{def:strictspeb}.

\item\label{rem:strictspe2b}
Conditions \ref{def:strictspe}\eqref{def:strictspeb} and \eqref{def:strictspec} are automatically satisfied if we restrict to $\hb\not\in \Sing(\Lambda)$, as we remarked in \ref{rem:psi}\eqref{rem:psib}. These are really conditions near the singular points with respect to $\Lambda$, as defined in \T\ref{subsec:defstar}.

\item\label{rem:strictspe1}
Assume that $\cM$ and $\cN$ are strictly specializable along $\cX_0$. If $\varphi:\cM\to\cN$ is any $\cR_\cX$-linear morphism, it induces a morphism $\psi_{t,\alpha}\varphi:\psi_{t,\alpha}\cM\to\psi_{t,\alpha}\cN$, by Remark \ref{rem:psi}\eqref{rem:psi0}. According to \ref{def:strictspe}\eqref{def:strictspeb} or \ref{def:strictspe}\eqref{def:strictspec} (and to \eqref{rem:strictspe2} above), in order to get the strict specializability of $\varphi$, it is enough to verify strictness of $\psi_{t,\alpha}\varphi$ for one representative $\alpha$ of each class in $\Lambda/\ZZ$, and for $\alpha=-1,0$. Assume now that $\varphi$ is strictly specializable. Then we have locally $\gr_a^{V^{(\hb_o)}}\varphi=\oplus_{\alpha\in\Lambda,\,\ell_{\hb_o}(\alpha)=a} \psi_{t,\alpha}\varphi$, which is therefore also strict.

\item\label{rem:strictspe3}
Let $U^{(\hb_o)}_\bbullet\cM$ be a good $V$-filtration indexed by $\ell_{\hb_o}(\Lambda)$ which satisfies \eqref{eq:hospe}, \ref{def:strictspe}\eqref{def:strictspeb} and \eqref{def:strictspec}, but satisfies \ref{def:strictspe}\eqref{def:strictspea} for $a\leq0$ only. Then one shows by induction that $\partiall_t:\gr_{a}^{U^{(\hb_o)}}\to\gr_{a+1}^{U^{(\hb_o)}}\cM$ is an isomorphism if $a>-1$ and therefore \ref{def:strictspe}\eqref{def:strictspea} is satisfied for any $a$.
\end{enumerate}
\end{Remarques}

\begin{lemme}\label{lem:morstrictspe}
Let $\varphi:\cM\to\cN$ be a strictly specializable morphism between strictly specializable modules. Then $\varphi$ is $V$-strict, $\ker\varphi$ and $\coker\varphi$ are strictly specializable and, for each $\alpha$, one has
\[
\psi_{t,\alpha}\ker\varphi=\ker\psi_{t,\alpha}\varphi\quad\text{and}\quad \psi_{t,\alpha}\coker\varphi=\coker\psi_{t,\alpha}\varphi.
\]
\end{lemme}

\begin{proof}
Fix $\hb_o$ and forget about the exponent $(\hb_o)$. As indicated in Remark \ref{rem:psi}, the result holds after inverting polynomials of the variable $\hb$ vanishing on $\Sing\Lambda(\cM)\cup\Sing\Lambda(\cN)$. The assumption on $\varphi$ is made to control the behaviour near the singular set.

Let us prove the $V$-strictness of $\varphi$. We have to show that, for any $a\in \RR$, we have $\im\varphi\cap V_a\cN=\varphi(V_a\cM))$. As both filtrations $\im\varphi\cap V_\bbullet\cN$ and $\varphi(V_\bbullet\cM))$ of $\im\varphi$ are good (Artin-Rees for the first one), there exists $k\in\NN$ such that, for any $a\in\RR$ we have $\im\varphi\cap V_{a-k}\cN\subset \varphi(V_a\cM))$. Therefore, if, for any $a$, the morphism $V_a/V_{a-k}(\varphi)$ is strict for the induced $V$-filtrations on $V_a/V_{a-k}(\cM)$ and $V_a/V_{a-k}(\cN)$, then $\varphi$ is $V$-strict.

Let us now show, by induction on the length of the induced $V$-filtration, that, for any $a',a$ with $a'<a$, the morphism $V_a/V_{a'}(\varphi)$ is $V$-strict and that $\coker V_a/V_{a'}(\varphi)$ is strict. This is by assumption if the length is one. Let $a''\in{}]a',a[$ be a jumping index and let $n$ be a local section of $\varphi(V_a\cM)\cap V_{a''}\cN+V_{a'}\cN$. There exists a polynomial $p(\hb)$ such that $p(\hb)\cdot n$ is a local section of $\varphi(V_{a''}\cM)+V_{a'}\cN$, as indicated above. Hence the class of $p(\hb)n$ in $\coker V_a/V_{a''}(\varphi)$ is zero. By induction, the previous module is strict, hence the class of $n$ itself vanishes, that is, $n\in \varphi(V_a''\cM)+V_{a'}\cN$. In other words, $V_a/V_{a'}(\varphi)$ is $V$-strict. As a consequence, we have $\gr_{a''}^V\coker V_a/V_{a'}(\varphi)=\coker \gr_{a''}^V\varphi$ for any $a''\in{}]a',a]$, which is strict by assumption, hence, by Lemma \ref{lem:Wstrict}\eqref{lem:Wstrict1}, $\coker V_a/V_{a'}(\varphi)$ is strict. This gives the $V$-strictness of $\varphi$.

Put on $\ker\varphi$ and $\coker\varphi$ the filtration naturally induced by $V$. This is a good filtration satisfying \eqref{eq:hospe}. We know that it satisfies \ref{def:strictspe}\eqref{def:strictspea} on $\ker\varphi$ (Lemma \ref{lem:uniciteMK}\eqref{lem:uniciteMK1}), so we call it $V_{\bbullet}\ker\varphi$. By the $V$-strictness of $\varphi$ that we have just  proved, we have an exact sequence
\[
0\to\gr_{a}^{V}\ker\varphi\to\gr_{a}^{V}\cM\To{\gr_{a}^{V}{\varphi}}\gr_{a}^{V}\cN\to\gr_{a}^{U}\coker\varphi\to 0.
\]
By assumption and Remark \ref{rem:strictspe}\eqref{rem:strictspe1}, $\coker\gr_{a}^{V}{\varphi}$ is strict, so $U_\bbullet\coker\varphi$ also satisfies \ref{def:strictspe}\eqref{def:strictspea}. By Lemma \ref{lem:uniciteMK}, it is equal to the $V$-filtration on $\coker\varphi$. Now, \ref{def:strictspe}\eqref{def:strictspeb} and \eqref{def:strictspec} are clear.
\end{proof}

\begin{proposition}\label{prop:canvar}
Let $\cM$ be a strictly specializable $\cR_\cX$-module.
{\def\theenumi{\alph{enumi}}
\begin{enumerate}
\item\label{prop:canvara}
If $\cM=\cM'\oplus \cM''$, then $\cM'$ and $\cM''$ are strictly specializable.
\item\label{prop:canvarb}
If $\cM$ is supported on $\cX_0$, then $V_{<0}\cM=0$ and $\cM=i_+V_0\cM$.
\item\label{prop:canvarc}
The following properties are equivalent:
{\def\theenumii{\arabic{enumii}}
\begin{enumerate}
\item\label{prop:canvarc1}
$\var:\psi_{t,0}\cM\to\psi_{t,-1}\cM$ is injective,
\item\label{prop:canvarc2}
$\cM$ has no proper sub-$\cR_\cX$-module supported on $\cX_0$,
\item\label{prop:canvarc3}
$\cM$ has no proper subobject in $\cS^2(X,t)$ supported on $\cX_0$.
\end{enumerate}}
\item\label{prop:canvard}
If $\can:\psi_{t,-1}\cM\to\psi_{t,0}\cM$ is onto, then $\cM$ has no proper quotient satisfying \ref{def:strictspe}\eqref{def:strictspea} and supported on $\cX_0$.
\item\label{prop:canvare}
The following properties are equivalent:
{\def\labelenumii{\textrm{(\theenumii${}'$)}}
\def\theenumii{\arabic{enumii}}
\begin{enumerate}
\item\label{prop:canvare1}
$\psi_{t,0}\cM=\im\can\oplus \ker\var$,
\item\label{prop:canvare2}
$\cM=\cM'\oplus \cM''$ with $\cM'$ satisfying \eqref{prop:canvarc} and \eqref{prop:canvard} and $\cM''$ supported on~$\cX_0$.
\end{enumerate}}
\end{enumerate}}
\end{proposition}

In \eqref{prop:canvarb}, one should put an exponent $(\hb_o)$; however, as a consequence of the proof, the lattice $\Lambda$ is then contained in $\ZZ$, and therefore the various filtrations $V_\bbullet^{(\hb)}\cM$ glue as a global $V$-filtration, so that the statement is not ambiguous.

\begin{proof}
We work locally near $\hb_o$ and forget the exponent $(\hb_o)$ in the notation. Because of Properties \ref{def:strictspe}\eqref{def:strictspeb} and \eqref{def:strictspec} and Remark \ref{rem:strictspe}\eqref{rem:strictspe0}, we may replace $\psi_{t,0}$ with $\gr_0^{V}$ and $\psi_{t,-1}$ with $\gr_{-1}^V$ in \ref{prop:canvar}\eqref{prop:canvarc}, \eqref{prop:canvard} and \eqref{prop:canvare}.

\smallskip
\eqref{prop:canvara} For $\cN=\cM'$ or $\cM''$, put $U_a\cN=V_a\cM\cap\cN$. Then $U_\bbullet\cM\defin U_\bbullet\cM'\oplus U_\bbullet\cM''$ is a good filtration of $\cM$ satisfying \eqref{eq:hospe}. As $\gr_a^{U}\cN$ is a submodule of $\gr_a^{V}\cM$, it is strict. From Lemma \ref{lem:uniciteMK} one concludes that $U_\bbullet\cM=V_\bbullet\cM$ and it follows that $\cM'$ and $\cM''$ are strictly specializable with $U_\bbullet\cM'=V_\bbullet\cM'$ and $U_\bbullet\cM''=V_\bbullet\cM''$.

\smallskip
\eqref{prop:canvarb} As $t$ is injective on $V_{<0}\cM$, one has $V_{<0}\cM=0$. Similarly, $\gr_a^{V}\cM=0$ for $a\not\in\NN$. As $t$ is injective on $\gr_{k}^{V}\cM$ for $k\geq 1$ (\cf Remark \ref{rem:psi}\eqref{rem:psia}), one has
\[
V_0\cM=\ker \lefcro t:\cM\to\cM\rigcro.
\]
As $\partiall_t:\gr_{k}^{V}\cM\to\gr_{k+1}^{V}\cM$ is an isomorphism for $k\geq 0$, one gets
\[
\cM=\ooplus_{k\geq 0}V_0\cM\partiall_{t}^{k}=i_+V_0\cM.
\]

\eqref{prop:canvarc} $\eqref{prop:canvarc1}\ssi\eqref{prop:canvarc2}$: It is enough to show that the mappings
$$\xymatrix@C=-.5cm@R=.5cm{
&\ker\lefcro t:V_0\cM\to V_{-1}\cM\rigcro\ar@{^{ (}->}[dl]\ar[dr]&\\
\ker\lefcro t:\cM\to\cM\rigcro&&\ker\lefcro t:\gr_{0}^{V}\cM\to\gr_{-1}^{V}\cM\rigcro
}
$$
are isomorphisms. It is clear for the right one, since $t:V_{<0}\cM\to V_{<-1}\cM$ is an isomorphism. For the left one this follows from the fact that $t$ is injective on $\gr_a^{V}\cM$ for $a\neq0$ according to Remark \ref{rem:psi}\eqref{rem:psia}.

\enlargethispage{3mm}%
\smallskip
$\eqref{prop:canvarc2}\ssi\eqref{prop:canvarc3}$: let us verify $\Leftarrow$ (the other implication is clear). Let $\cT$ denote the $t$-torsion submodule of $\cM$ and $\cT'$ the submodule generated by
\[
\cT_0\defin\ker \lefcro t:\cM\to\cM\rigcro.
\]
\begin{assertion*}
$\cT'$ is a subobject of $\cM$ in $\cS^2(X,t)$.
\end{assertion*}

This assertion gives the implication $\Leftarrow$ because, by assumption, $\cT'=0$, hence $t:\cM\to\cM$ is injective, so $\cT=0$.

\begin{proof}[Proof of the assertion]
Let us show first that $\cT'$ is $\cR_\cX$-coherent. As we remarked above, we have $\cT_0=\ker[t:\gr_0^V\cM\to\gr_{-1}^V\cM]$. Now, $\cT_0$ is the kernel of a linear morphism between $\cR_{\cX_0}$-coherent modules, hence is also $\cR_{\cX_0}$-coherent. It follows that $\cT'$ is $\cR_\cX$-coherent.

Let us now show that $\cT'$ is strictly specializable. Notice that $\cT_0$ is strict because it is equal to
\[
\ker\lefcro t:\gr_{0}^{V}\cM\to\gr_{-1}^{V}\cM\rigcro.
\]
Let $U_\bbullet\cT'$ be the filtration induced by $V_\bbullet\cM$ on $\cT'$. One shows as in \eqref{prop:canvarb} that \hbox{$U_{<0}\cT'\!=\!0$} and $\gr_a^{U}\cT'=0$ for $a\not\in\NN$. Let us show by induction on $k$ that
\[
U_k\cT'=\cT_0+\partiall_t\cT_0+\cdots+\partiall_{t}^{k}\cT_0.
\]
Denote by $U'_k\cT'$ the right hand term. The inclusion $\supset$ is clear. Let $(x_o,\hb_o)\in\cX_0$, $m\in U_k\cT'_{(x_o,\hb_o)}$ and let $\ell\geq k$ such that $m\in U'_\ell\cT'_{(x_o,\hb_o)}$. If $\ell>k$ one has $m\in\cT'_{(x_o,\hb_o)}\cap V_{\ell-1}\cM_{(x_o,\hb_o)}$ hence $t^\ell m\in V_{-1}\cM_{(x_o,\hb_o)}\cap \cT'_{(x_o,\hb_o)}=0$. Put
\[
m=m_0+\partiall_tm_1+\cdots+\partiall_{t}^\ell m_\ell,
\]
with $tm_j=0$ ($j=0,\ldots,\ell$). One then has $t^\ell\partiall_{t}^\ell m_\ell=0$ and, as
$$
t^\ell\partiall_{t}^\ell m_\ell=\prod_{j=0}^\ell(\partiall_tt-j\hb)\cdot m_\ell=(-1)^\ell\ell!\hb^\ell m_\ell
$$
and $\cT_0$ is strict, one concludes that $m_\ell=0$, hence $m\in U'_{\ell-1}\cT'_{(x_o,\hb_o)}$. This implies the other inclusion.

As $\gr_a^{U}\cT'$ is strict (because it is contained in $\gr_a^{V}\cM$), one deduces from Remark \ref{rem:psi}\eqref{rem:psia} that $\partiall_t:\gr_{k}^{U}\cT'\to\gr_{k+1}^{U}\cT'$ is injective for $k\geq 0$. The previous computation shows that it is onto, hence $\cT'$ is strictly specializable and $U_\bbullet\cT'$ is its Malgrange-Kashiwara filtration.

\enlargethispage{3mm}%
It is now enough to prove that the injective morphism $\gr_{0}^{U}\cT'\to\gr_{0}^{V}\cM$ is strict. But its cokernel is identified with the submodule $\im[t:\gr_{0}^{V}\cM\to\gr_{-1}^{V}\cM]$ of $\gr_{-1}^{V}\cM$, which is strict.
\end{proof}

\eqref{prop:canvard} If $\can$ is onto, then $\cM=\cR_\cX\cdot V_{<0}\cM$. If $\cM$ has a $t$-torsion quotient $\cT$ satisfying \ref{def:strictspe}\eqref{def:strictspea}, then $V_{<0}\cT=0$, so $V_{<0}\cM$ is contained in $\ker[\cM\to\cT]$ and so is $\cR_\cX\cdot V_{<0}\cM=\cM$, hence $\cT=0$.

\smallskip
\eqref{prop:canvare} $\eqref{prop:canvare1}\implique\eqref{prop:canvare2}$: Put
\[
U_0\cM'=V_{<0}\cM+\partiall_tV_{-1}\cM\qqbox{and}\cT_0=\ker \lefcro t:\cM\to\cM\rigcro.
\]
The assumption $(1')$ is equivalent to $V_0\cM=U_0\cM'\oplus \cT_0$. Define
\[
U_k\cM'=V_k\cR_\cX\cdot U_0\cM'\qqbox{and}U_k\cM''=V_k\cR_\cX\cdot \cT_0
\]
for $k\geq 0$. As $V_k\cM=V_{k-1}\cM+\partiall_tV_{k-1}\cM$ for $k\geq 1$, one has $V_k\cM=U_k\cM'+U_k\cM''$ for $k\geq 0$. Let us show by induction on $k\geq 0$ that this sum is direct. Fix $k\geq 1$ and let $m\in U_k\cM'\cap U_k\cM''$. Write
$$
m=m'_{k-1}+\partiall_t n'_{k-1}= m''_{k-1}+\partiall_t n''_{k-1}
$$
with $m'_{k-1},n'_{k-1}\in U_{k-1}\cM'$ and $m''_{k-1},n''_{k-1}\in U_{k-1}\cM''$. One has $\partiall_t [n'_{k-1}]= \partiall_t[n''_{k-1}]$ in $V_k\cM/V_{k-1}\cM$, hence, as
\[
\partiall_t:V_{k-1}\cM/V_{k-2}\cM\to V_k\cM/V_{k-1}\cM
\]
is bijective for $k\geq 1$, one gets $[n'_{k-1}]= [n''_{k-1}]$ in $V_{k-1}\cM/V_{k-2}\cM$ and by induction one deduces that both terms are zero. One concludes that \hbox{$m\!\in\! U_{k-1}\cM'\cap U_{k-1}\cM''\!=\!0$} by induction.

This implies that $\cM=\cM'\oplus \cM''$ with $\cM'\defin\cup_kU_k\cM'$ and $\cM''$ defined similarly. It follows from \eqref{prop:canvara} that both $\cM'$ and $\cM''$ are strictly specializable and the filtrations~$U_\bbullet$ above are their Malgrange-Kashiwara filtrations. In particular $\cM'$ satisfies \eqref{prop:canvarc} and~\eqref{prop:canvard}.

\smallskip
$\eqref{prop:canvare2}\implique \eqref{prop:canvare1}$: One has $V_{<0}\cM''=0$. Let us show that $\im\can=\gr_{0}^{V}\cM'$ and $\ker\var=\gr_{0}^{V}\cM''$. The inclusions $\im\can\subset\gr_{0}^{V}\cM'$ and $\ker\var\supset\gr_{0}^{V}\cM''$ are clear. Moreover $\gr_{0}^{V}\cM'\cap \ker\var=0$ as $\cM'$ satisfies \eqref{prop:canvarc}. Last, $\can:\gr_{-1}^{V}\cM'\to\gr_{0}^{V}\cM'$ is onto, as $\cM'$ satisfies \eqref{prop:canvard}. Hence $\gr_{0}^{V}\cM=\im\can\oplus \ker\var$.
\end{proof}

\begin{corollaire}[Kashiwara's equivalence]\label{cor:Kashiwaraequiv}
The functor $i_+$ induces an equivalence between the category of coherent strict $\cR_{\cX_0}$-modules (and strict morphisms) and the full subcategory $\cS_{\cX_0}^{2}(X,t)$ of $\cS^2(X,t)$ consisting of objects supported on $\cX_0$. An inverse functor is $\psi_{t,0}$.
\end{corollaire}

\begin{proof}
It follows from Proposition \ref{prop:canvar}\eqref{prop:canvarb}.
\end{proof}

\begin{proposition}[Strict specializability along $\{t^r=0\}$]\label{prop:strictspetr}
\hspace*{-2mm}
Assume that $\cM$~is~strictly spe\-cia\-liza\-ble along $\{t=0\}$. Put $f=t^r$ for $r\geq2$. Then $\cM$ is strictly spe\-cia\-liza\-ble along $\{f=0\}$ and, if we denote by $i:\{t=0\}\hto X$ the closed inclusion, we have $(\psi_{f,\alpha}\cM,\rN)=(i_+\psi_{t,r\alpha}\cM, \rN/r)$ for any $\alpha$ and an isomorphism
\[
\lefacc\xymatrix{\psi_{f,-1}\cM\ar@/^/[rr]^-{\can_f}&&\psi_{f,0}\cM\ar@/^/[ll]^-{\var_f}}\rigacc \simeq
i_+\lefacc\xymatrix@C.4cm{
\psi_{t,-r}\cM\ar@/^2.5pc/[rrrr]^-{\can_f\defin\can_t\circ(t^{r-1})^{-1}}&&\psi_{t,-1}\cM\ar@{-->}[ll]_-{t^{r-1}}^-{\sim}\ar@{-->}@/^/[rr]^-{\can_t}
&&\psi_{t,0}\cM
\ar@/^2.5pc/[llll]^-{\var_f\defin t^{r-1}\circ\var_t}\ar@{-->}@/^/[ll]^-{\var_t}
}\rigacc
\]
\end{proposition}

\begin{proof}
We fix $\hb_o$ and we forget about the exponent $(\hb_o)$, when working in a neighbourhood of $\hb_o$. We may write $i_{f,+}\cM=\oplus_{k\in\NN}\cM\otimes\partiall_u^k\delta$ as a $\cR_\cX\otimes_\CC\CC[u]\langle\partiall_u\rangle$-module, with
\begin{align*}
\partiall_u^k(m\otimes\delta)&=m\otimes\partiall_u^k\delta,\\
\partiall_t(m\otimes\delta)&=(\partiall_tm)\otimes\delta-(rt^{r-1}m)\otimes\partiall_u\delta,\\
u(m\otimes\delta)&=(t^rm)\otimes\delta,\\
\cO_\cX(m\otimes\delta)&=(\cO_\cX m)\otimes\delta,
\end{align*}
and with the usual commutation rules. For the sake of simplicity, we will write $\cR_{\cX\times\CC}$ instead of $\cR_\cX\otimes_\CC\CC[u]\langle\partiall_u\rangle$ and $V_0\cR_{\cX\times\CC}$ instead of $\cR_\cX\otimes_\CC\CC[u]\langle u\partiall_u\rangle$.

For $a\leq0$, put
\[
V_a i_{f,+}\cM\defin V_0\cR_{\cX\times\CC}\cdot\big(V_{ra}\cM\otimes\delta\big),
\]
and for $a>0$ define inductively
\[
V_a i_{f,+}\cM\defin V_{<a}i_{f,+}\cM+\partiall_uV_{a-1}i_{f,+}\cM.
\]
Assume that $a\leq0$. Using the relation
\[
(\partiall_uu+\alpha\star\hb)(m\otimes\delta)=\frac1r\big([(\partiall_tt+r\alpha\star\hb)m]\otimes\delta-\partiall_t(tm\otimes\delta)\big),
\]
one shows that, if
\[
\prod_{\ell_{\hb_o}(r\alpha)=ra}(\partiall_tt+(r\alpha)\star\hb)^{\nu_{r\alpha}}V_{ra}\cM\subset V_{<ra}\cM,
\]
then
\[
\prod_{\ell_{\hb_o}(\alpha)=a}(\partiall_uu+\alpha\star\hb)^{\nu_{r\alpha}}V_a i_{f,+}\cM\subset V_{<a}i_{f,+}\cM,
\]
thus \eqref{eq:hospe} for $a\leq0$.

If $m_1,\dots,m_\ell$ generate $V_{ra}\cM$ over $V_0\cR_\cX$, then $m_1\otimes\delta,\dots,m_\ell\otimes\delta$ generate $V_a i_{f,+}\cM$ over $V_0\cR_{\cX\times\CC}$, as follows from the relation
\[
(\partiall_ttm)\otimes\delta=(\partiall_tt-r\partiall_uu)(m\otimes\delta).
\]
It follows that $V_a i_{f,+}\cM$ is $V_0\cR_{\cX\times\CC}$-coherent for any $a\leq0$, hence for any $a$.

For any $a$ we clearly have
\[
V_{a-1}i_{f,+}\cM\subset uV_a i_{f,+}\cM \quad(\text{\resp}\quad V_{a+1}i_{f,+}\cM\subset V_{<a+1}i_{f,+}\cM+\partiall_uV_a i_{f,+}\cM)
\]
with equality if $a<0$ (\resp if $a\geq-1$), as an analogous property is true for $\cM$. Therefore, $V_{\bbullet}i_{f,+}\cM$ is a good $V$-filtration.

According to Remark \ref{rem:strictspe}\eqref{rem:strictspe3}, it is now enough to prove the second part of the proposition, hence we now assume that $a\leq0$.

We have $V_a i_{f,+}\cM=V_{<a}i_{f,+}\cM+ \sum_{k\geq0}\partiall_t^k(V_{ra}\cM\otimes\delta)$. One shows, by considering the degree in $\partiall_u$, that the natural map
\begin{align*}
\oplus_k\gr_{ra}^V\cM\partiall_t^k&\to\gr_a^Vi_{f,+}\cM\\
\oplus_k[m_k]\partiall_t^k&\mto
\Big[\sum_k\partiall_t^k(m_k\otimes\delta)\Big]
\end{align*}
is an isomorphism of $\cR_\cX$-modules. The desired assertions follow.
\end{proof}

\begin{proposition}[Restriction to $\hb=\hb_o$]\label{prop:restrhbo}
Let $\hb_o\in\Omega_0$ and let $\cM$ be a strictly specializable $\cR_\cX$-module with Malgrange-Kashiwara filtration $V^{(\hb_o)}_\bbullet\cM$ near $\hb_o$. Put $\index{$mhb$@$M_{\hb_o}$}M_{\hb_o}=\cM/(\hb-\hb_o)\cM$. \begin{enumerate}
\item
For any $a\in\RR$, we have, near $\hb_o$,
\[
V^{(\hb_o)}_a(\cM)\cap(\hb-\hb_o)\cM=(\hb-\hb_o)\cdot V^{(\hb_o)}_a(\cM).
\]
The filtration $U_\bbullet(M_{\hb_o})$ naturally induced by $V^{(\hb_o)}_\bbullet(\cM)$ on $\cM/(\hb-\hb_o)\cM$ is good with respect to $V_\bbullet\cD_X$ ($\hb_o\neq0$) or to $V_\bbullet\gr^F\cD_X$ ($\hb_o=0$) and, for any $a$,
\[
\gr_a^U(M_{\hb_o})= \gr_a^{V^{(\hb_o)}}\cM\Big/(\hb-\hb_o)\gr_a^{V^{(\hb_o)}}\cM.
\]
Moreover, $\gr_a^U(M_{\hb_o})$ is naturally decomposed as the direct sum $\oplus_{\alpha\mid\ell_{\hb_o}(\alpha)=a}\psi^{(\hb_o)}_{t,\alpha}M_{\hb_o}$, with
\[
\psi^{(\hb_o)}_{t,\alpha}M_{\hb_o}\defin\psi_{t,\alpha}\cM\big/(\hb-\hb_o)\psi_{t,\alpha}\cM.
\]
Last,
\par$\bbullet$
$t:\psi^{(\hb_o)}_{t,\alpha}M_{\hb_o}\to\psi^{(\hb_o)}_{t,\alpha-1}M_{\hb_o}$ is an isomorphism if $\ell_{\hb_o}(\alpha)<0$,
\par$\bbullet$
$\partiall_t:\psi^{(\hb_o)}_{t,\alpha}M_{\hb_o}\to\psi^{(\hb_o)}_{t,\alpha+1}M_{\hb_o}$ is an isomorphism if $\ell_{\hb_o}(\alpha)\geq-1$ but $\alpha\neq-1$.

We still denote by $\rN:\psi^{(\hb_o)}_{t,\alpha}M_{\hb_o}\to\psi^{(\hb_o)}_{t,\alpha}M_{\hb_o}$ the nilpotent endomorphism induced by~$\rN$, and similarly for $\can$ and $\var$.
\item
If $\hb_o\neq0$, then
\begin{enumerate}
\item
$M_{\hb_o}$ is specializable along $X_0$ as a $\cD_X$-module,
\item
$-(\partial_tt+\alpha\star\hb_o/\hb_o)=\rN/\hb_o$ is nilpotent on $\psi^{(\hb_o)}_{t,\alpha}M_{\hb_o}$,
\item
$\partial_t:\psi^{(\hb_o)}_{t,\alpha}M_{\hb_o}\to\psi^{(\hb_o)}_{t,\alpha+1}M_{\hb_o}$ is an isomorphism if $\ell_{\hb_o}(\alpha)\geq-1$ but $\alpha\neq-1$.
\end{enumerate}

\item
If $\hb_o\not\in\Sing\Lambda$, in particular if $\hb_o\not\in i\RR$, the induced $V$-filtration $U_\bbullet M_{\hb_o}$ is a Malgrange-Kashiwara filtration $V^{(L_{\hb_o})}_\bbullet M_{\hb_o}$ (for some $\RR$-linear form $L_{\hb_o}$ depending on $\hb_o$ only) and we have
$$\psi^{(\hb_o)}_{t,\alpha}M_{\hb_o}= \psi_{t,(\alpha\star\hb_o)/\hb_o}M_{\hb_o};$$
in particular, $\psi^{(\hb_o)}_{t,\alpha}M_{\hb_o}= \psi_{t,\alpha}M_{\hb_o}$ for $\alpha\in \RR$, \eg $\alpha=-1,0$; moreover, $\rN$ induces $\hb_o\rN_{\hb_o}$, $\can$ induces $\hb_o\can_{\hb_o}$ and $\var$ induces $\var_{\hb_o}$ (where $\rN_{\hb_o},\can_{\hb_o},\var_{\hb_o}$ are defined in \T\ref{sec:reviewspe} for $M_{\hb_o}$).
\end{enumerate}
\end{proposition}

\begin{proof}
We work locally near $\hb_o$ and forget $(\hb_o)$ in the notation. Let $m$ be a local section of $V_a(\cM)\cap(\hb-\hb_o)\cM$. Then $m=(\hb-\hb_o)n$ where $n$ is a local section of $V_b\cM$ for some $b$. If $b>a$, then $n$ induces a torsion element in $\gr_b^V\cM$, hence $n\in V_{<b}\cM$ by \ref{def:strictspe}\eqref{def:strictspea}. This gives the first assertion. The other assertions are clear (\cf \T\ref{sec:reviewspe}). Let us give more details when $\hb_o\not\in \Sing\Lambda$. If we fix $a$, the roots of the minimal polynomial of $-\partial_tt$ on $\gr_a^U(M_{\hb_o})$ are the $\gamma=(\alpha\star\hb_o)/\hb_o$ where $\alpha$ satisfies $\ell_{\hb_o}(\alpha)=a$. If $\hb_o\not\in \Sing\Lambda$, then $\gamma\mto a$ is a $\RR$-linear form on $\gamma$ (because $\alpha\mto(\alpha\star\hb_o)/\hb_o$ is a $\RR$-linear automorphism of $\CC$), that we denote by $L_{\hb_o}(\gamma)$. Notice that, as $\hb_o\not\in \Sing\Lambda$, we have $\gamma\in\ZZ$ only if $\alpha\in\ZZ$ and therefore $L_{\hb_o}(\ZZ)\subset\ZZ$. The filtration $U_\bbullet$ is then equal to the Malgrange-Kashiwara filtration associated with $L_{\hb_o}$ (see \T\ref{sec:reviewspe}).

If $\hb_o\in \Sing\Lambda\moins\{0\}$, the roots of the minimal polynomial of $-\partial_tt$ on $\gr_a^U(M_{\hb_o})$ may differ by a nonzero integer, hence the filtration $U_\bbullet$ is not useful to compute $\psi_t(M_{\hb_o})$.

On the other hand, if $\hb_o\in\RR^*$, we have $\ell_{\hb_o}=\reel$ and $\reel(\alpha\star\hb_o/\hb_o)=\reel(\alpha)$, so the roots of the minimal polynomial of $-\partial_tt$ on $\gr_a^U(M_{\hb_o})$ are the $\alpha$ for which $\reel(\alpha)=a$. In the case $\hb_o=1$ (corresponding to the functor $\Xi_{\DR}$), we also have a perfect correspondence with $\can$ and $\var$.
\end{proof}

\begin{remarque*}
With the only assumption of strict specializability, we cannot give general statements concerning the behaviour of Properties \eqref{prop:canvarb} to \eqref{prop:canvare} of Proposition \ref{prop:canvar} by restriction to $\hb=\hb_o\neq0$. We will come back on this in Proposition \ref{prop:restricanvar}.
\end{remarque*}

We may now reformulate Theorem \ref{th:imdirspe} for strictly specializable modules. Let us take notation used in this theorem.

\begin{theoreme}\label{th:imdirstrictspe}
Assume that $\cM$ is good and strictly specializable along $X\times\{0\}$, and that $F$ is proper on the support of $\cM$. Assume moreover that, for any $\alpha$, the complexes $f_\dag\psi_{t,\alpha}\cM$ are strict. Then the $\cR_{\cY\times\CC}$-modules $\cH^i(F_\dag\cM)$ are strictly specializable along $Y\times\{0\}$. Moreover, for any $\alpha$, we have a canonical and functorial isomorphism
\[
\psi_{t,\alpha}\cH^i(F_\dag\cM)=\cH^i(f_\dag\psi_{t,\alpha}\cM).
\]
\end{theoreme}

\begin{proof}
We may work locally near $\hb_o\in\Omega_0$ and we forget the exponent $(\hb_o)$. Notice that, because of Remark \ref{rem:strictspe}\eqref{rem:strictspe0}, the complex $f_\dag\gr_a^V\cM$ is strict for any $a\in\RR$. Denote by $U_\bbullet\cH^i(F_\dag\cM)$ the filtration induced by $\cH^i(f_\dag V_\bbullet\cM)$. We may apply to it the conclusions of Theorem \ref{th:imdirspe} (after extending it to the case of filtrations indexed by $I_\RR+\ZZ$). This filtration satisfies \eqref{eq:hospe} and, by assumption, \ref{def:strictspe}\eqref{def:strictspea}. It is therefore equal to the Malgrange-Kashiwara filtration (\cf Lemma \ref{lem:uniciteMK}). We hence have
\[
\gr_a^V\cH^i(F_\dag\cM)=\cH^i(f_\dag\gr_a^V\cM)= \ooplus_{\alpha'=a}\cH^i(f_\dag\psi_{t,\alpha}\cM).
\]
As the image of $\cH^i(f_\dag\psi_{t,\alpha}\cM)\to\gr_a^V\cH^i(F_\dag\cM)$ is clearly contained in $\psi_{t,\alpha}\cH^i(F_\dag\cM)$, it is therefore equal to it.

The canonical isomorphism $\cH^i(f_\dag\psi_{t,\alpha}\cM)\isom\psi_{t,\alpha}\cH^i(F_\dag\cM)$ that we have constructed is \emph{a~priori} only defined locally near $\hb_o$. However, it is locally independent of the choice of $\hb_o$: to see this, replace $\gr_a^V$ above with $V_{a+\varepsilon}^{(\hb_o)}/V_{a-\varepsilon}^{(\hb_o)}$ for $\varepsilon$ small enough, and argue as in the proof of Lemma \ref{lem:uniciteMK}\eqref{lem:uniciteMK7}. Therefore, it is globally defined with respect to $\hb$.

As $t:V_a\cM\to V_{a-1}\cM$ is an isomorphism for $a<0$, it also induces an isomorphism $\cH^i(f_\dag V_a\cM)\isom \cH^i(f_\dag V_{a-1}\cM)$. As $\cH^i(f_\dag V_a\cM)\to\cH^i(F_\dag\cM)$ is injective for any~$a$ and has image $V_a\cH^i(F_\dag\cM)$, this shows that $\cH^i(F_\dag\cM)$ satisfies \ref{def:strictspe}\eqref{def:strictspeb} and, by the same argument, \ref{def:strictspe}\eqref{def:strictspec}.
\end{proof}

\begin{remarque}
It is enough to verify the strictness condition of the theorem for those $\alpha$ such that $\reel(\alpha)\in[-1,0[$ and for $\alpha=0$: indeed, strictness is a local property with respect to $\hb$, and one may apply locally \ref{def:strictspe}\eqref{def:strictspeb} and \eqref{def:strictspec}.
\end{remarque}

\section{Localization and minimal extension across a hypersurface}\label{sec:minext}

\subsection{Localization of a strictly specializable $\cR_\cX$-module} \label{subsec:loc}
Consider the sheaf of rings $\cR_\cX[\tm]$. Notice that we have $\cR_\cX[\tm]=(V_0\cR_\cX)[\tm]$, as $\partiall_t=\tm (t\partiall_t)$. This ring has a $V$-filtration defined by $V_k\cR_\cX[\tm]=t^{-k}V_0\cR_\cX$. One can define the notion of a good $V$-filtration for a coherent $\cR_\cX[\tm]$-module $\index{$mcurltilde$@$\wt\cM$}\wt\cM$, as well as the notion of specializability. Then Lemma \ref{lem:uniciteMK} applies similarly, and shows the existence of a canonical $V_\bbullet^{(\hb_o)}$-filtration.

The situation simplifies here, as $t:V_a^{(\hb_o)}\wt\cM\to V_{a-1}^{(\hb_o)}\wt\cM$ is an isomorphism for any~$a$. It follows that $t:\psi_{t,\alpha}\wt\cM\to\psi_{t,\alpha-1}\wt\cM$ is an isomorphism for any $\alpha$, and we do not need to consider the $(\can,\var)$ diagram. Moreover, strict specializability reduces here to Condition \ref{def:strictspe}\eqref{def:strictspea}, as we are not interested in Condition \ref{def:strictspe}\eqref{def:strictspec}.

\begin{lemme}\label{lem:localstrictspe}
Let $\cM$ be a coherent $\cR_\cX$-module which is strictly specializable along $t=0$. Then $\wt\cM\defin\cR_\cX[\tm]\otimes_{\cR_\cX}\cM$ is a coherent $\cR_\cX[\tm]$-module which is strictly specializable along $t=0$. Moreover, the natural morphism $\cM\to\wt\cM$ induces, in the neighbourhood of any $\hb_o\in\Omega_0$, an isomorphism of $V_0\cR_\cX$-modules
\[
\forall a<0,\quad V_a^{(\hb_o)}\cM\isom V_a^{(\hb_o)}\wt\cM.
\]
\end{lemme}

\begin{proof}
We have $\wt\cM=\cO_\cX[\tm]\otimes_{\cO_\cX}\cM$ as a $\cO_\cX[\tm]$-module. Locally, the injective map $V_{0}^{(\hb_o)}\cM\hto\cM$ induces, by flatness of $\cO_\cX[\tm]$ over $\cO_\cX$, an injective map $\cO_\cX[\tm]\otimes_{\cO_\cX}V_{0}^{(\hb_o)}\cM\hto\wt\cM$, which is onto because $\cM=\sum_{k\geq0}\partiall_t^kV_{0}^{(\hb_o)}\cM$, so that
\[
\wt\cM=\sum_{k\geq0}t^{-k}\cO_\cX[\tm]\otimes t^k\partiall_t^kV_{0}^{(\hb_o)}\cM=\sum_{k\geq0}t^{-k}\cO_\cX[\tm]\otimes V_{0}^{(\hb_o)}\cM.
\]
As $\gr_0^{V^{(\hb_o)}}\cM$ is killed (locally) by some power of $t$, we also have $\wt\cM=\cO_\cX[\tm]\otimes_{\cO_\cX}V_{<0}^{(\hb_o)}\cM$. Put then $V_{a+k}^{(\hb_o)}\wt\cM=t^{-k}\otimes V_{a}^{(\hb_o)}\cM$ for any $a\in[-1,0[$ and any $k\in\ZZ$. This defines a filtration of $\wt\cM$, which has all the properties required for the Malgrange-Kashiwara filtration. As $t:V_a^{(\hb_o)}\cM\to V_{a-1}^{(\hb_o)}\cM$ is bijective for $a<0$, we get the required isomorphism.
\end{proof}

\begin{lemme}\label{lem:psiMMt}
Let $\cM$ be strictly specializable along $t=0$. Put \hbox{$\wt\cM=\cR_\cX[\tm]\!\smash{\ootimes_{\cR_\cX}}\!\cM$}. Then,
\begin{enumerate}
\item\label{lem:psiMMt1}
for any $\alpha\not\in\NN$, we have $\psi_{t,\alpha}\cM\subset\psi_{t,\alpha}\wt\cM$,
\item\label{lem:psiMMt2}
for any $\alpha$ with $\reel\alpha\in[-1,0[$, we have $\psi_{t,\alpha}\cMS=\psi_{t,\alpha}\wt\cMS$.
\end{enumerate}
\end{lemme}

\begin{Proof}
\begin{enumerate}
\item
We know, by Lemma \ref{lem:localstrictspe}, that the inclusion is an equality near $\hb_o$ if \hbox{$\ell_{\hb_o}(\alpha)\!<\!0$}. Fix now $\alpha\not\in\NN$. Locally near $\hb_o$, there exists $k\in\NN$ such that $\ell_{\hb_o}(\alpha-k)<0$. We have a commutative diagram
\[
\xymatrix{\psi_{t,\alpha}\cM\ar@{_{ (}->}[d]_-{t^k}\ar[r]&\psi_{t,\alpha}\wt\cM\ar[d]_-{\wr}^-{t^k}\\
\psi_{t,\alpha-k}\cM\ar@{=}[r]&\psi_{t,\alpha-k}\wt\cM
}
\]
where the left vertical arrow is injective, by Remark \ref{rem:psi}\eqref{rem:psia}, hence the result.

\item
Choose $k\in\NN^*$ big enough so that $\ell_{\hb}(\alpha-k)<0$ for any $\alpha$ with $\reel\alpha\in[-1,0[$ and any $\hb\in\bS$. Notice now that, for any $\ell\geq0$, $(\alpha-\ell)\star\hb\neq0$ for $\hb\in\bS$, as $\reel(\alpha-\ell)\neq0$ (\cf\T\ref{subsec:defstar}). Therefore, as $t\partiall_t+(\alpha-\ell)\star\hb$ is nilpotent on $\psi_{t,\alpha-\ell-1}\cMS$, $t\partiall_t$ is bijective on it and thus the map $t:\psi_{t,\alpha-\ell}\cMS\to\psi_{t,\alpha-\ell-1}\cMS$ is onto, so that the left vertical map in the diagram above, restricted to $\bS$, is onto, as was to be proved.\qedhere
\end{enumerate}
\end{Proof}

\begin{definition}[Nearby cycles]\label{def:nearby}
Let $\cM$ be strictly specializable along $t=0$ and put $\wt\cM=\cR_\cX[\tm]\otimes_{\cR_\cX}\cM$. For $\alpha$ such that $\reel\alpha\in[-1,0[$, put
\[
\index{$psimaj$@$\Psi_{t,\alpha}\cM$}\Psi_{t,\alpha}\cM\defin\psi_{t,\alpha}\wt\cM.
\]
\end{definition}

\begin{Remarques}\label{rem:nearby}
\begin{enumerate}
\item\label{rem:nearby1}
By the definition of strict specializability, we have, for any $\alpha\not\in\NN$, a \emph{local} (with respect to~$\hb$) isomorphism $\Psi_{t,\alpha}\cM\simeq\psi_{t,\alpha}\cM$, given by a suitable power of $t$ or of $\partiall_t$.
\item\label{rem:nearby2}
On the other hand, by Lemma \ref{lem:psiMMt}\eqref{lem:psiMMt1}, we have $\psi_{t,\alpha}\cM\subset\Psi_{t,\alpha}\cM$ and $\psi_{t,\alpha}\cMS=\Psi_{t,\alpha}\cMS$, for any $\alpha$ with $\reel\alpha\in[-1,0[$.
\item\label{rem:nearby3}
Last, if $\alpha$ is \emph{real} and in $[-1,0[$, we have $\psi_{t,\alpha}\cM=\Psi_{t,\alpha}\cM$, as $\ell_{\hb}(\alpha)=\alpha<0$ for any $\hb$.
\end{enumerate}
\end{Remarques}

\begin{remarque}[Strict specializability along $\{t^r=0\}$] \label{rem:strictspetr}
Proposition \ref{prop:strictspetr} (forgetting the assertion on $\can$ and $\var$) applies to strictly specializable $\cR_\cX[\tm]$-modules. Remark that, with the notation of \loccit, the action of $u$ is invertible on $i_{f,+}\wt\cM$. Therefore, we deduce that, if $\reel\alpha\in[-1,0[$, $(\Psi_{f,\alpha}\cM,\rN)=i_+(\Psi_{t,r\alpha-\lceil r\alpha\rceil}\cM,\rN/r)$, where $\lceil r\alpha\rceil\in \ZZ$ is such that $\reel(r\alpha-\lceil r\alpha\rceil)\in[-1,0[$.
\end{remarque}

\Subsection{Minimal extension across a hypersurface}\label{subsec:minext}
\begin{proposition}\label{prop:minext}
Let $\wt\cM$ be a strictly specializable $\cR_\cX[\tm]$-module. In the neighbourhood of any $\hb_o\in\Omega_0$, consider the $\cR_\cX$-submodule $\cM^{(\hb_o)}$ of $\wt\cM$ generated by $V_{<0}^{(\hb_o)}\wt\cM$. Then, the various $\cM^{(\hb_o)}$ glue as a coherent $\cR_\cX$-submodule $\cM$ of $\wt\cM$, which is the unique strictly specializable $\cR_\cX$-submodule of $\wt\cM$ satisfying
\begin{enumerate}
\item\label{prop:minext1}
$\cR_\cX[\tm]\otimes_{\cR_\cX}\cM=\wt\cM$,
\item\label{prop:minext2}
$\can$ is onto and $\var$ is injective.
\end{enumerate}
Moreover, the filtration defined by
\begin{equation}\tag*{(\protect\ref{prop:minext})$(*)$}
\label{prop:minext*}
V_a^{(\hb_o)}\cM=
\begin{cases}
V_a^{(\hb_o)}\wt\cM&\text{if }a<0,\\
\sum_{\substack{\ell\in\NN,a'<0\\ a'+\ell\leq a}} \partiall_t^{\ell}V_{a'}^{(\hb_o)}\wt\cM&\text{if }a\geq0.
\end{cases}
\end{equation}
is its Malgrange-Kashiwara filtration near $\hb_o$.
\end{proposition}

\begin{definition}[Minimal extension]\label{def:minext}
We call the $\cR_\cX$-submodule $\cM$ of $\wt\cM$ given by Proposition \ref{prop:minext} the \emph{minimal extension of $\wt\cM$ across $t=0$}.
\end{definition}

\begin{proof}[Proof of Proposition \ref{prop:minext}]
The question is local on $\cX$, so we work in the neighbourhood of some compact set in $X$ and on some disc $\Delta_{\hb_o}(\eta)$, on which $V_\bbullet^{(\hb_o)}\wt\cM$ exists. We denote by $\Lambda$ the set of indices of this filtration.

Let us first prove that $\cM$ does not depend on $\hb_o$. We have to prove that, if $\eta$ is small enough, then for any $\hb\in\Delta_{\hb_o}(\eta)$, the germ $\cM^{(\hb)}_\hb$ is equal to the germ $\cM^{(\hb_o)}_\hb$. The problem comes from the fact that $V_{<0}^{(\hb_o)}\wt\cM$ may not induce $V_{<0}^{(\hb)}\wt\cM$ at $\hb$. Fix $\epsilon>0$ such that
\begin{equation}\label{eq:minV1}
V_{<0}^{(\hb_o)}\wt\cM=V_{-\epsilon}^{(\hb_o)}\wt\cM.
\end{equation}
Then, after Lemma \ref{lem:uniciteMK}\eqref{lem:uniciteMK5} and \eqref{eq:Vlocconst}, we have, for $\eta>0$ small enough and any $\hb\in\Delta_{\hb_o}(\eta)$,
\begin{equation}\label{eq:minV2}
V_{<0}^{(\hb_o)}\wt\cM_\hb=V_{-\epsilon}^{(\hb_o)}\wt\cM_\hb=V_{-\epsilon}^{(\hb)}\wt\cM_\hb\subset V_{<0}^{(\hb)}\wt\cM_\hb.
\end{equation}
Therefore, $\cM^{(\hb_o)}_\hb\subset\cM^{(\hb)}_\hb$. In order to prove the reverse inclusion, it is enough to show the inclusion
\begin{equation}\label{eq:reverse}
V_{<0}^{(\hb)}\wt\cM_\hb\subset V_{-\epsilon}^{(\hb)}\wt\cM_\hb+\partiall_t V_{<-1}^{(\hb)}\wt\cM_\hb= V_{-\epsilon}^{(\hb)}\wt\cM_\hb+\partiall_tt V_{<0}^{(\hb)}\wt\cM_\hb,
\end{equation}
because if $\eta$ is small enough, we have $V_{<-1}^{(\hb)}\wt\cM_\hb\subset V_{<0}^{(\hb_o)}\wt\cM_\hb$; equivalently, it is enough to show that, for any $\alpha\in\Lambda$ and any $\hb\neq\hb_o$ near $\hb_o$ such that $\ell_{\hb_o}(\alpha)\in{}]-\epsilon,0[$, the operator $\partiall_tt$ is onto on $\psi_{t,\alpha}\wt\cM_\hb$. Recall (\cf Lemma \ref{lem:imhb}) that, if $\alpha$ is such that $\ell_{\hb_o}(\alpha)=0$, then $\alpha\star\hb_o=0$ if and only if $\alpha=0$. Therefore, for any $\epsilon>0$ small enough there exists $\eta>0$ such that
\[
(\ell_{\hb_o}(\alpha)\in{}]-\epsilon,0[,\ \hb\in\Delta_{\hb_o}(\eta)\text{ and } \alpha\in\Lambda\moins\{0\})\implique \alpha\star\hb\neq0.
\]
As $-(\partiall_tt+\alpha\star\hb)$ is nilpotent on $\psi_{t,\alpha}\wt\cM_\hb$, the previous choice of $\epsilon,\eta$ is convenient to get \eqref{eq:reverse}.

Clearly, each $\cM^{(\hb_o)}$ is $\cR_\cX$-coherent on the open set where it is defined. Let us now show that it is strictly specializable along $t=0$. Near $\hb_o$, $\cM=\cM^{(\hb_o)}$ comes equipped with a filtration $V_{\bbullet}^{(\hb_o)}\cM$ defined by \ref{prop:minext*}. This $V$-filtration is good, and $\cM$ is specializable. This filtration satisfies Properties \ref{def:strictspe}\eqref{def:strictspea}, \eqref{def:strictspeb} and \eqref{def:strictspec}: indeed, this follows from the strict specializability of $\wt\cM$ for \eqref{def:strictspea} with $a<0$ and for \eqref{def:strictspeb}; for \eqref{def:strictspea} with $a=0$, notice that $\gr_0^{V^{(\hb_o)}}\cM$ is identified with the image of $\partiall_t:\gr_{-1}^{V^{(\hb_o)}}\wt\cM\to\gr_0^{V^{(\hb_o)}}\wt\cM$ by construction, hence is strict, being contained in the strict module $\gr_0^{V^{(\hb_o)}}\wt\cM$; for \eqref{def:strictspec}, this follows from the definition of the $V$-filtration; apply then Remark \ref{rem:strictspe}\eqref{rem:strictspe3}. Therefore, $\cM$ is strictly specializable.

Similarly, $\partiall_t:\psi_{t,-1}\cM\to\psi_{t,0}\cM$ is onto, by construction, as $\psi_{t,0}\cM$ is identified with
\[
\image\big[\partiall_t:\psi_{t,-1}\wt\cM\to\psi_{t,0}\wt\cM\big].
\]
As $t:\psi_{t,0}\wt\cM\to\psi_{t,-1}\wt\cM$ is an isomorphism, we conclude that, for $\cM$, $\var$ is injective.

Let us end with the uniqueness statement. Let $\cN\subset\wt\cM$ satisfying \ref{prop:minext}\eqref{prop:minext1} and \eqref{prop:minext2}. Then, by Lemma \ref{lem:localstrictspe} and \ref{prop:minext}\eqref{prop:minext1}, $V_{<0}^{(\hb_o)}\cN=V_{<0}^{(\hb_o)}\wt\cN=V_{<0}^{(\hb_o)}\wt\cM$. As $\can$ is onto and as $\cN$ satisfies \ref{def:strictspe}\eqref{def:strictspec}, we have $\cN=\cR_\cX\cdot V_{<0}^{(\hb_o)}\cN$ near $\hb_o$, hence the desired uniqueness assertion.
\end{proof}

\section{Strictly S(upport)-decomposable $\cR_\cX$-modules} \label{sec:Sdecomp}
\begin{definition}\label{def:strictdec}
We say that a $\cR_\cX$-module $\cM$ is
\begin{itemize}
\item
\emph{strictly S-decomposable along $\cX_0$} if it is strictly specializable along $\cX_0$ and satisfies the equivalent conditions \ref{prop:canvar}\eqref{prop:canvare};
\item
\emph{strictly S-decomposable at $x_o\in X$} if for any analytic germ $f:(X,x_o)\to (\CC,0)$, $i_{f,+}\cM$ is strictly S-decomposable along $\cX\times \{0\}$ in some neighbourhood of $x_o$;
\item
\emph{strictly S-decomposable} if it is strictly S-decomposable at all points $x_o\in X$.
\end{itemize}
\end{definition}

\pagebreak[2]
\begin{Lemme}\label{lem:strictdec}
\begin{enumerate}
\item\label{lem:strictdec1}
If $\cM$ is strictly S-decomposable along $\{t=0\}$, then it is strictly S-decomposable along $\{t^r=0\}$ for any $r\geq1$.
\item\label{lem:strictdec2}
If $\cM=\cM_1\oplus\cM_2$, then $\cM$ is strictly S-decomposable of some kind if and only if $\cM_1$ and $\cM_2$ are so.
\item\label{lem:strictdec3}
Assume that $\cM$ is strictly S-decomposable and $Z$ is pure dimensional. Then
the following conditions are equivalent:
\begin{enumerate}
\item\label{lem:strictdec3a}
for any analytic germ $f:(X,x_o)\to(\CC,0)$ such that $f^{-1}(0)\cap Z$ has everywhere codimension one in $Z$, $i_{f,+}\cM$ satisfies both conditions \ref{prop:canvar}\eqref{prop:canvarc} and~\eqref{prop:canvard};
\item\label{lem:strictdec3b}
near any $x_o$, there is no coherent submodule of $\cM$ with support having codimension $\geq1$ in $Z$;
\item\label{lem:strictdec3c}
near any $x_o$, there is no nonzero morphism $\varphi:\cM\to\cN$, with $\cN$ strictly S-decomposable at $x_o$, such that $\im\varphi$ is supported in codimension $\geq1$ in $Z$.
\end{enumerate}
\end{enumerate}
\end{Lemme}

\begin{definition}\label{def:strictsupp}
Let $Z$ be a pure dimensional closed analytic subset of $X$ and let $\cM$ be strictly S-decomposable. We say that $\cM$ has \emph{strict support $Z$} if the equivalent conditions of \ref{lem:strictdec}\eqref{lem:strictdec3} are satisfied.
\end{definition}

\begin{proof}[Proof of Lemma \ref{lem:strictdec}]
The first point is a direct consequence of Proposition \ref{prop:strictspetr} and the second one is clear. For the third one, let us show for instance $\eqref{lem:strictdec3a}\iff\eqref{lem:strictdec3c}$. Let $\varphi:\cM\to\cN$, with $\cN$ strictly S-decomposable at $x_o$, such that $\im\varphi\subset f^{-1}(0)$. Then \ref{prop:canvar}\eqref{prop:canvard} implies that $\im\varphi=0$. Conversely, given $f$ such that $f^{-1}(0)$ has everywhere codimension one in $Z$, decompose $i_{f,+}\cM$ as in \ref{prop:canvar}\eqref{prop:canvare}. Then \eqref{lem:strictdec3c} implies that $\cM''=0$.
\end{proof}

We will now show that a strictly S-decomposable holonomic $\cR_\cX$-module may indeed be decomposed as the direct sum of holonomic $\cR_\cX$-modules having strict support. We first consider the local decomposition and, by uniqueness, we get the global one. It is important for that to be able to define \emph{a~priori} the strict components. They are obtained from the characteristic variety.

\begin{proposition}
Let $\cM$ be holonomic and strictly S-decomposable at $x_o$, and let $(Z_i,x_o)_{i\in I}$ be the minimal family of closed irreducible analytic germs $(Z_i,x_o)$ such that $\Char\cM\subset\cup_iT^*_{Z_i}X\times \Omega_0$ near $x_o$. There exists a unique decomposition $\cM_{x_o}=\oplus_{i\in I}\cM_{Z_i,x_o}$ of germs at $x_o$ such that $\cM_{Z_i,x_o}=0$ or has strict support $(Z_i,x_o)$.
\end{proposition}

\begin{proof}
We will argue by induction on $\dim\supp\cM$. First, we reduce to the case where the support $S$ of $\cM$ (see after Definition \ref{def:holo}) is irreducible. Let $S'$ be an irreducible component of $S$ at $x_o$ and let $S''$ be the union of all other ones. Let $f:(X,x_o)\to (\CC,0)$ be an analytic germ such that $S''\subset f^{-1}(0)$ and $(S',x_o)\not\subset f^{-1}(0)$. Then, according to \ref{prop:canvar}\eqref{prop:canvare}, near $x_o$, $\cM$ has a decomposition $\cM=\cM'\oplus \cM''$, with $\cM'$ supported on $S'$ and satisfying \ref{prop:canvar}\eqref{prop:canvarc} and \eqref{prop:canvard}, and $\cM''$ supported on $S''$.

Conversely, if we have any local decomposition $\cM=\oplus\cM_{S_i}$, with $(S_i,x_o)$ irreducible and $\cM_{S_i}$ (strictly S-decomposable after Lemma \ref{lem:strictdec}\eqref{lem:strictdec2}) having strict support~$S_i$, then $S_i\subset S'$ or $S_i\subset S''$ and $\cM'=\oplus_{S_i\not\subset S''}\cM_{S_i}$, $\cM''=\oplus_{S_i\subset S''}\cM_{S_i}$.

By induction on the number of irreducible components, we are reduced to the case where $(S,x_o)$ is irreducible. We may assume that $\dim S>0$.

Choose now a germ $f:(X,x_o)\to (\CC,0)$ which is nonconstant on $S$ and such that $f^{-1}(0)$ contains all components $Z_i$ except $S$. We have, as above, a unique decomposition $\cM=\cM'\oplus \cM''$ of germs at $x_o$, where $\cM'$ satisfies \ref{prop:canvar}\eqref{prop:canvarc} and \eqref{prop:canvard}, and $\cM''$ is supported on $f^{-1}(0)$, by Proposition \ref{prop:canvar}\eqref{prop:canvare}. Moreover, $\cM'$ and $\cM''$ are also strictly S-decomposable at $x_o$. We may apply the inductive assumption to~$\cM''$.

Let us show that $\cM'$ has strict support $S$ near $x_o$: if $\cM'_1$ is a coherent submodule of $\cM'$ supported on a strict analytic subset $Z\subset S$, then $Z$ is contained in the union of the components $Z_i$, hence $\cM'_1$ is supported in $f^{-1}(0)$, so is zero.
\end{proof}

By uniqueness of the local decomposition, we get:
\begin{corollaire}
Let $\cM$ be holonomic and strictly S-decomposable on $X$ and let $(Z_i)_{i\in I}$ be the minimal (locally finite) family of closed irreducible analytic subsets $Z_i$ such that $\Char\cM^\Cir\subset\cup_iT^*_{Z_i}X\times \CC^*$. There exists a unique decomposition $\cM=\oplus_i\cM_{Z_i}$ such that each $\cM_{Z_i}=0$ or has strict support $Z_i$.\qed
\end{corollaire}

A closed analytic irreducible subset $Z$ of $X$ such that $\cM_Z\neq0$ is called a \emph{strict component} of $\cM$.

\begin{corollaire}\label{cor:morphismestrictdecompo}
Let $\cM',\cM''$ be two holonomic $\cR_\cX$-module which are strictly S-decomposable and let $(Z_i)_{i\in I}$ be the family of their strict components.
Then any morphism $\varphi:\cM'_{Z_i}\to\cM''_{Z_j}$ vanishes identically if $Z_i\neq Z_j$.
\end{corollaire}

\begin{proof}
The image of $\varphi$ is supported on $Z_i\cap Z_j$, which is a proper closed analytic subset either of $Z_i$ or of $Z_j$, if $Z_i\neq Z_j$.
\end{proof}

The following will be useful:

\begin{corollaire}\label{cor:strictdecstrict}
Let $\cM$ be holonomic and strictly S-decomposable. Then $\cM$ is strict.
\end{corollaire}

\begin{proof}
The question is local, and we may assume that $\cM$ has strict support $Z$ with $Z$ closed irreducible analytic near $x_o$.

First, there exists an open dense set of $Z$ on which $\cM$ is strict. Indeed, by Kashiwara's equivalence on the smooth part of $Z$, we may reduce to the case where $Z=X$, and by restricting to a dense open set, we may assume that $\Char\cM$ is the zero section. Hence we are reduced to the case where $\cM$ is $\cO_\cX$-coherent. If $t$ is a local coordinate, notice that $\cM/t\cM=\psi_{t,-1}\cM$, as the filtration defined by $U_k\cM=t^{-k}\cM$ for $k\leq0$ and $U_k\cM=\cM$ for $k\geq0$ satisfies all properties of the Malgrange-Kashiwara filtration. Let $m$ be a local section of $\cM$ killed by $p(\hb)$. Then $m$ is zero in $\cM/t\cM$ by strict specializability. As $\cM$ is $\cO_\cX$-coherent, Nakayama's lemma implies that $m=0$.

Let now $m$ be a local section of $\cM$ near $x_o$ killed by some $p(\hb)$. Then $\cR_\cX\cdot m$ is supported by a strict analytic set of $Z$ near $x_o$ by the previous argument. As $\cM$ has strict support $Z$, we conclude that $m=0$.
\end{proof}

Let us end this paragraph with a result concerning sesquilinear pairings:

\begin{proposition}\label{prop:Cdecomposable}
Let $\cM',\cM''$ be two holonomic $\cR_\cX$-module which are strictly S-decomposable and let $(Z_i)_{i\in I}$ be the family of their strict components.
Then any sesquilinear pairing $C:\cM'_{Z_i|\bS}\otimes_{\cO_\bS}\ov{\cM''_{Z_j|\bS}}\to\Dbh{X}$ vanishes identically if \hbox{$Z_i\neq Z_j$}.
\end{proposition}

\begin{proof}
The assertion is local on $X\times\bS$, so we fix $\hb_o\in\bS$ and $x_o\in X$ and we work with germs at $(x_o,\hb_o)$. Assume for instance that $Z_i$ is not contained in $Z_j$ and consider an analytic function, that we may assume to be a local coordinate~$t$ by Kashiwara's equivalence, such that $t\equiv0$ on $Z_j$ and $t\not\equiv0$ on $Z_i$. Consider $C$ as a morphism $\cM'_{Z_i|\bS}\to\cHom_{\ov{\cR_{\cX|\bS}}}(\ov{\cM''_{Z_j|\bS}},\Dbh{X})$. Fix local $\cR_{\cX}$-generators $m''_1,\ldots,m''_\ell$ of $\cM''_{Z_j,(x_o,\hb_o)}$. By \ref{prop:canvar}\eqref{prop:canvarb}, there exists $q\geq 0$ such that $t^qm''_k=0$ for all $k=1,\ldots,\ell$. Let $m'\in\cM'_{Z_i,(x_o,\hb_o)}$ and let $p$ be the maximum of the orders of $C(m')(\ov{m''_k})$ on some neighbourhood of $(x_o,\hb_o)$. As $t^{p+1+q}/\ov t^q$ is $C^p$, we have, for any $k=1,\dots,\ell$,
\[
t^{p+1+q}C(m')(\ov{m''_k})=\frac{t^{p+1+q}}{\ov t^q}\cdot \ov t^qC(m')(\ov{m''_k})=0,
\]
hence $t^{p+1+q}C(m')\equiv0$. Applying this to generators of $\cM'_{Z_i,(x_o,\hb_o)}$ shows that all local sections of $C(\cM'_{Z_i,(x_o,\hb_o)})$ are killed by some power of $t$.

As $\cM'_{Z_i}$ has strict support $Z_i$, we know from Proposition~\ref{prop:canvar}\eqref{prop:canvard} that $V^{(\hb_o)}_{<0}\cM'_{Z_i,(x_o,\hb_o)}$ generates $\cM'_{Z_i,(x_o,\hb_o)}$ over $\cR_{\cX}$. It is therefore enough to show that $C(V^{(\hb_o)}_{<0}\cM'_{Z_i,(x_o,\hb_o)})=0$.

On the one hand, we have $C(V^{(\hb_o)}_{k}\cM'_{Z_i,(x_o,\hb_o)})=0$ for $k\ll0$: indeed, by Lemma~\ref{lem:tinj}, \hbox{$t:C(V^{(\hb_o)}_{k}\cM'_{Z_i,(x_o,\hb_o)})\to C(V^{(\hb_o)}_{k-1}\cM'_{Z_i, \hb_o})$} is an isomorphism for $k\ll0$, hence acts injectively on $C(V^{(\hb_o)}_{k}\cM'_{Z_i,(x_o,\hb_o)})$, therefore the conclusion follows, as $t$ is also nilpotent by the argument above.

Let now $k<0$ be such that $C(V^{(\hb_o)}_{k-1}\cM'_{Z_i,(x_o,\hb_o)})=0$, and let $m'$ be a section of $V^{(\hb_o)}_{k}\cM'_{Z_i,(x_o,\hb_o)}$; there exists $b(s)$ of the form $\prod_{\alpha\mid\ell_{\hb_o}(\alpha)\in[k,k-1[}(s-\alpha\star\hb)^{\nu_\alpha}$ such that $b(-\partiall_tt)m'\in V^{(\hb_o)}_{k-1}\cM'_{Z_i,(x_o,\hb_o)}$, hence $b(-\partiall_tt)C(m')=0$; on the other hand, we have seen that there exists $N$ such that $t^{N+1}C(m')=0$, hence, putting $B(s)=\prod_{\ell=0}^{N}(s-\ell\hb)$, it also satisfies $B(-\partiall_tt)C(m')=0$; notice now that $b(s)$ and $B(s)$ have no common root, so there exists $p(\hb)\in\CC[\hb]\moins\{0\}$ such that $p(\hb)C(m')=0$. According to \eqref{eq:sanstorsion}, we conclude that $C(m')=0$.
\end{proof}

\section{Specialization of a sesquilinear pairing}\label{sec:spesesqui}
\subsection{Sesquilinear pairing on nearby cycles}
We keep notation of \S\T\ref{subsec:sesqui} and \ref{subsec:Vfil}.
Let $\cM'$ and $\cM''$ be two objects of $\cS^2(X,t)$ and let
\[
C:\cMS'\ootimes_{\cO_\bS}\ov{\cMS''} \to\Dbh{X}
\]
be a sesquilinear pairing. In the following, we will assume that the $\cR_\cX$-modules are also good (in the applications they will be holonomic). The purpose of this paragraph is to define, for each $\alpha\in\CC$ a sesquilinear pairing
\[
\index{$psitc$@$\psi_{t,\alpha}C$}\psi_{t,\alpha}C:\psi_{t,\alpha}\cMS'\ootimes_{\cO_\bS}\ov{\psi_{t,\alpha}\cMS''}\to\Dbh{X_0},
\]
compatible with $\rN$, \ie such that, with obvious notation,
\begin{equation}\label{eq:psiCN}
\psi_{t,\alpha}C(\rN[m],\ov{[\mu]})=\thb^2\psi_{t,\alpha}C([m],\ov{\rN[\mu]}),
\end{equation}
where $\rN$ denotes the action of $-(\partiall_tt+\alpha\star\hb)$ on $\psi_{t,\alpha}$. Using the notation of Tate twist introduced in \eqref{eq:Tate} and the notion of morphism of $\RTriples(X)$ introduced in Definition \ref{def:RtriplesX}, we will put $\index{$ncurl$@$\cN$}\cN=(\rN',\rN'')$ with $\rN''=i\rN$ and $\rN'=-\rN''=-i\rN$, so that $\cN$ is a morphism of $\RTriples(X_0)$:
\begin{equation}\label{eq:NTate}
\cN:\psi_{t,\alpha}(\cM',\cM'',C)\to\psi_{t,\alpha}(\cM',\cM'',C)(-1)
\end{equation}
which satisfies $\cN^*=-\cN$.

\begin{Remarques}\label{rem:primitif}
\begin{enumerate}
\item\label{rem:primitif1}
Once such specializations are defined, we get, according to the compatibility with $\rN$, pairings
\begin{equation*}
\index{$psitcl$@$\psi_{t,\alpha,\ell}C$}\psi_{t,\alpha,\ell}C:\gr_{-\ell}^{\rM}\psi_{t,\alpha}\cMS' \ootimes_{\cO_\bS} \ov{\gr_\ell^{\rM}\psi_{t,\alpha}\cMS''}\to\Dbh{X_0}.
\end{equation*}
In other words, using the notion of graded Lefschetz $\RTriples$ introduced in Remark \ref{rem:grlefrtriples}, the graded object
\[
\big(\gr_\bbullet^{\rM}\psi_{t,\alpha}\cT,\gr_{-2}^{\rM}\cN\big)\defin \Big(\ooplus_\ell\big(\gr_{-\ell}^{\rM}\psi_{t,\alpha}\cM',\gr_{\ell}^{\rM}\psi_{t,\alpha}\cM'',\psi_{t,\alpha,\ell}C\big),\gr_{-2}^{\rM}\cN\Big)
\]
is a graded Lefschetz triple with $\varepsilon=-1$.

\item\label{rem:primitif3}
As $C$ is $\cR_{(X,\ov X),\bS}$-linear, it easily follows from the definition of $\psi_{t,\alpha}C$ that we know all $\psi_{t,\alpha}C$ as soon as we know them for $\reel(\alpha)\in[-1,0[$ and for $\alpha=0$, according to \ref{def:strictspe}\eqref{def:strictspeb} and \eqref{def:strictspec}.
\end{enumerate}
\end{Remarques}

In order to define the specialization of $C$, we will use the residue of a Mellin transform, that we consider now.

Let $\cM',\cM''$ be two objects of $\cS^2(X,t)$ and let $C:\cMS'\otimes_{\cO_\bS}\ov{\cMS''}\to\Dbh{X}$ be a sesquilinear pairing. Fix $(x_o,\hb_o)\in X_0\times\bS$. For local sections $m',m''$ of $\cM',\cM''$ defined in some neighbourhood of $(x_o,\hb_o)$ in $X\times\bS$ the distribution $C(m',\ov{m''})$ has some finite order $p$ on $\nb_{X\times\bS}(x_o,\hb_o)$. For $2\reel s>p$, the function $\mt^{2s}$ is $C^p$, so for any such $s$, $\mt^{2s}C(m',\ov{m''})$ is a section of $\Dbh{X}$ on $\nb_{X\times\bS}(x_o,\hb_o)$. Moreover, for any relative form $\psi$ of maximal degree with compact support on $\nb_{X\times\bS}(x_o,\hb_o)$, the function $s\mto \big\langle\mt^{2s}C(m',\ov{m''}),\psi\big\rangle$ is holomorphic on the half-plane $\{2\reel s>p\}$. We say that $\mt^{2s}C(m',\ov{m''})$ \emph{depends holomorphically on $s$} on $\nb_{X\times\bS}(x_o,\hb_o)\times\{2\reel s>p\}$.

Let $\chi(t)$ be a real $C^\infty$ function with compact support, which is ${}\equiv 1$ near $t=0$. In the following, we will consider differential forms $\psi=\varphi\wedge\chi(t)\,\itwopi dt\wedge d\ov t$, where $\varphi$ is a relative form of maximal degree on $X_0\times\bS$.

\begin{proposition}\label{prop:MellinC}
Let $\cM',\cM'',C$ be as above. Then, for any \hbox{$(x_o,\hb_o)\!\in\! X_0\times\bS$}, there exists an integer $L\geq0$ and a finite set of complex numbers~$\gamma$ satisfying $\psi_{t,\gamma}\cM'_{(x_o,\hb_o)}\neq0$ and $\psi_{t,\gamma}\cM''_{(x_o,-\hb_o)}\neq0$, such that, for any element $m'$ of $\cM'_{(x_o,\hb_o)}$ and $m''$ of $\cM''_{(x_o,-\hb_o)}$, the correspondence
\begin{equation}\tag*{(\protect\ref{prop:MellinC})$(*)$}
\label{eq:MellinC}
\varphi\mto\prod_\gamma\Gamma(s-\gamma\star\hb/\hb)^L\cdot\big\langle\mt^{2s}C(m',\ov{m''}),\varphi\wedge\chi(t)\,\itwopi dt\wedge d\ov t\big\rangle
\end{equation}
defines, for any $s\in\CC$, a section of $\Dbh{X_0}$ on $\nb_{X_0\times\bS}(x_o,\hb_o)$ which is holomorphic with respect to $s\in\CC$.
\end{proposition}

The proposition asserts that the distribution
\[
\varphi\mto \big\langle\mt^{2s}C(m',\ov{m''}),\varphi\wedge\chi(t)\,\itwopi dt\wedge d\ov t\big\rangle
\]
extends as a distribution on $\nb_{X_0\times\bS}(x_o,\hb_o)$ depending meromorphically on $s$, with poles along the sets $s=-k+\gamma\star\hb/\hb$ ($k\in\NN$), and with a bounded order. Notice that changing the function $\chi$ will modify the previous meromorphic distribution by a holomorphic one, as $\mt^{2s}$ is $C^\infty$ for any $s$ away from $t=0$. The proposition is a consequence of the following more precise lemma.

\begin{lemme}\label{lem:polesI}
Let $(x_o,\hb_o)\in X_0\times\bS$ and let $a_1,a_2\in\RR$. There exist $L\geq0$ and a~finite set of $\gamma$ satisfying
$$
\psi_{t,\gamma}\cM'_{(x_o,\hb_o)}\neq0,\quad \psi_{t,\gamma}\cM''_{(x_o,-\hb_o)}\neq0,\quad \ell_{\hb_o}(\gamma)\leq a_1,\quad \ell_{-\hb_o}(\gamma)\leq a_2,
$$
such that, for any sections $m'\in V_{a_1}^{(\hb_o)}\cM'_{(x_o,\hb_o)}$ and $m''\in V_{a_2}^{(-\hb_o)}\cM''_{(x_o,-\hb_o)}$, the correspondence
\begin{equation}\tag*{(\protect\ref{lem:polesI})$(*)$}\label{eq:polesI*}
\varphi\mto\prod_\gamma\Gamma(s-\gamma\star\hb/\hb)^L\cdot\big\langle\mt^{2s}C(m',\ov{m''}),\varphi\wedge\chi(t)\,\itwopi dt\wedge d\ov t\big\rangle
\end{equation}
defines, for any $s\in\CC$, a section of $\Dbh{X_0}$ on $\nb_{X_0\times\bS}(x_o,\hb_o)$ which is holomorphic with respect to $s\in\CC$.

Assume moreover that the class of $m'$ (\resp $m''$) in $\gr_{a_1}^{V^{(\hb_o)}}\cM'_{(x_o,\hb_o)}$ (\resp in $\gr_{a_2}^{V^{(-\hb_o)}}\cM''_{(x_o,-\hb_o)}$) is in $\psi_{t,\alpha_1}\cM'_{(x_o,\hb_o)}$ (\resp in $\psi_{t,\alpha_2}\cM''_{(x_o,-\hb_o)}$). Then the product of $\Gamma$ factors can be indexed by a set of $\gamma$ satisfying moreover
\begin{equation}\tag*{(\protect\ref{lem:polesI})$(**)$}\label{eq:polesI**}
2\reel(\gamma)<a_1+a_2\quad\text{or, if $\alpha_1=\alpha_2\defin\alpha$, }\gamma=\alpha.
\end{equation}
\end{lemme}

\begin{proof}
Let $b_{m'}(S)=\prod_{\gamma\in A(m')}(S-\gamma\star\hb)^{\nu(\gamma)}$ be the Bernstein polynomial of $m'$ (\cf Corollary~\ref{cor:Bernstein}), with $\nu(\gamma)$ bounded by the nilpotency index $L$ of $\rN$. It is enough to prove that $\prod_{\gamma\in A(m')}\Gamma(s-\gamma\star\hb/\hb)^{\nu(\gamma)}$ is a convenient product of $\Gamma$ factors. Indeed, arguing similarly for $m''$, one obtains that the product indexed by $A(m')\cap A(m'')$ is convenient. It is then easy to verify that Conditions \ref{eq:polesI**} on $\gamma$ are satisfied by any $\gamma\in A(m')\cap A(m'')$. Remark that $\ell_{\hb_o}(\gamma)+\ell_{-\hb_o}(\gamma)=2\reel(\gamma)$.

Notice first that, for any local section $Q$ of $V_0\cR_{\cX,(x_o,\hb_o)}$, and any $C^\infty$-form $\psi$ on $\nb_{X\times\bS}(x_o,\hb_o)$ with compact support, the form $(\mt^{2s}\psi)\cdot Q$ is $C^p$ with compact support if $2\reel s>p$. Applying this to the Bernstein relation $Q=b_{m'}(-\partiall_tt)- tP$ for $m'$, one gets
\begin{equation}\label{eq:extension}
\begin{aligned}
0&=\big\langle [b_{m'}(-\partiall_tt)-tP]C(m',\ov{m''}),\mt^{2s}\psi\big\rangle\\
&=\big\langle C(m',\ov{m''}),(\mt^{2s}\psi)\cdot[b_{m'}(-\partiall_tt)- tP]\big\rangle\\
&=b_{m'}(\hb s)\big\langle C(m',\ov{m''}),\mt^{2s}\psi\big\rangle + \big\langle C(m',\ov{m''}),\mt^{2s}t\eta\big\rangle
\end{aligned}
\end{equation}
for some $\eta$, which is a polynomial in $s$ with coefficients being $C^\infty$ with compact support contained in that of $\varphi$. As $\mt^{2s}t$ is $C^p$ for $2\reel s+1>p$, we can argue by induction to show that, for any $\psi$ and $k\in\NN$,
\begin{equation}\label{eq:extensionk}
s\mto b_{m'}(\hb(s-k+1))\cdots b_{m'}(\hb s)\big\langle\mt^{2s}C(m',\ov{m''}),\psi\big\rangle
\end{equation}
extends as a holomorphic function on $\{s\mid2\reel s>p-k\}$. Apply this result to $\psi=\varphi\wedge\chi(t)\,\itwopi dt\wedge d\ov t$ to get the lemma.
\end{proof}

\begin{remarque}\label{rem:extension}
The previous proof also applies if we only assume that $C$ is $\cR_{(X,\ov X),\bS}$-linear \emph{away from} $\{t=0\}$. Indeed, this implies that $[b_{m'}(-\partiall_tt)- tP]C(m',\ov{m''})$ is supported on $\{t=0\}$, and \eqref{eq:extension} only holds for $\reel s$ big enough, maybe $\gg p$. Then, \eqref{eq:extensionk} coincides with a holomorphic distribution defined on $\{s\mid2\reel s>p-k\}$ only for $\reel s\gg0$. But, by uniqueness of analytic extension, it coincides with it on $\reel s>p$.
\end{remarque}

A distribution on $X_0\times\bS/\bS$ which is continuous with respect to $\hb$ and holomorphic with respect to $\bS$ can be restricted as a distribution to sets of the form $s=\alpha\star\hb/\hb$. This restriction is continuous with respect to $\hb$. By a similar argument, the polar coefficients along $s=\alpha\star\hb/\hb $ of the meromorphic distribution \hbox{$\big\langle\mt^{2s}C(m',\ov{m''}),\cbbullet\wedge\chi(t)\,\itwopi dt\wedge d\ov t\big\rangle$} exist as \emph{semi-meromorphic} distributions on $\nb_{X_0\times\bS}(x_o,\hb_o)$ (\ie the exists a polynomial $p(\hb)$ such that, after multiplication by~$p(\hb)$, the distribution is continuous with respect to~$\hb$). The possible poles are the $\hb\in\bS$ such that there exists $\gamma$ as in \ref{eq:polesI**} with $(\gamma-n-\alpha)\star\hb=0$, $n\in\NN$ and $n\neq0$ if $\gamma=\alpha$.

\begin{lemme}\label{lem:coefshol}
Let $[m']$ be a local section of $\psi_{t,\alpha}\cM'$ near $(x_o,\hb_o)$ and $[m'']$ a local section of $\psi_{t,\alpha}\cM''$ near $(x_o,-\hb_o)$. Then, the polar coefficients of the distribution $\big\langle\mt^{2s}C(m',\ov{m''}),\cbbullet\wedge\chi(t)\,\itwopi dt\wedge d\ov t\big\rangle$ along $s=\alpha\star\hb/\hb$ do neither depend on the choice of the local liftings $m',m''$ of $[m'],[m'']$ nor on the choice of $\chi$, and take value in $\Dbh{X_0}$.
\end{lemme}

\begin{proof}
Indeed, any other local lifting of $m'$ can be written as $m'+\mu'$, where $\mu'$ is a local section of $V_{<\ell_{\hb_o}(\alpha)}^{(\hb_o)}\cM'$. By the previous lemma, $\big\langle\mt^{2s}C(\mu',\ov{m''}),\cbbullet\wedge\chi(t)\,\itwopi dt\wedge d\ov t\big\rangle$ is holomorphic along $s=\alpha\star\hb/\hb$. Notice also that a different choice of the function~$\chi$ does not modify the polar coefficients.

We want to show that the polar coefficients do not have poles in some neighbourhood of $\bS$. The possible poles of the polar coefficients, as we have seen above, are such that $(\gamma-n-\alpha)\star\hb=0$, with $n\in\NN$ and $n\neq0$ if $\gamma=\alpha$. Now, \ref{eq:polesI**} shows that the only possible $\gamma\neq\alpha$ are such that $\reel(\gamma)<\reel(\alpha)$, hence for any $\gamma,n$ that we have to consider, we have $\reel(\gamma-n-\alpha)<0$, hence $\neq0$. Now, there can be no $\hb\in\bS$ with $(\gamma-n-\alpha)\star\hb=0$ (by \T\ref{subsec:defstar}, we should have $\hb=\pm i$ and $\gamma-n-\alpha$ purely imaginary).
\end{proof}

According to this lemma, we get a sesquilinear pairing
\begin{equation}\label{eq:psiC}\index{$psitc$@$\psi_{t,\alpha}C$}
\begin{split}
\psi_{t,\alpha}\cMS'\ootimes_{\cO_\bS} \ov{\psi_{t,\alpha}\cMS''}&\To{\psi_{t,\alpha}C} \Dbh{X_0}\\
([m'],\ov{[m'']}) &\mto\res_{s=\alpha\star\hb/\hb} \big\langle\mt^{2s}C(m',\ov{m''}),\cbbullet\wedge\chi(t)\,\itwopi dt\wedge d\ov t\big\rangle,
\end{split}
\end{equation}
where $m',m''$ are local liftings of $[m'],[m'']$. The compatibility \eqref{eq:psiCN} of $\psi_{t,\alpha}C$ with~$\rN$ follows from $\ov t\partiall_{\ov t}\mt^{2s}=\thb^{-2}t\partiall_t\mt^{2s}$ (recall that $\alpha\star\hb/\hb$ is real).

\begin{definition}\label{def:psiC}
For $\reel\alpha\in[-1,0[$, the specialized sesquilinear pairing $\Psi_{t,\alpha}C$ is defined as $\psi_{t,\alpha}C$, according to Remark~\ref{rem:nearby}\eqref{rem:nearby2}.
\end{definition}

\begin{remarque}\label{rem:psitwist}
We have defined a functor $\psi_{t,\alpha}$, and similarly $\Psi_{t,\alpha}$ if $\reel\alpha\in[-1,0[$, from the category of strictly specializable objects of $\RTriples(X)$, \ie objects $\cT=(\cM',\cM'',C)$ such that $\cM'$ and $\cM''$ are in $\cS^2(X,t)$, to the category $\RTriples(X_0)$ by putting $\psi_{t,\alpha}\cT=(\psi_{t,\alpha}\cM',\psi_{t,\alpha}\cM'',\psi_{t,\alpha}C)$. This functor clearly commutes with any Tate twist by $k\in\hZZ$.
\end{remarque}

\begin{remarque}[Behaviour with respect to adjunction]\label{rem:psiadj}
As $\chi(t)\,\itwopi dt\wedge d\ov t$ is real, we have $\psi_{t,\alpha}(C^*)=(\psi_{t,\alpha}C)^*$. If $\cS:\cT\to\cT^*(-w)$ is a sesquilinear duality of weight~$w$ on $\cT$, then $\psi_{t,\alpha}\cS$ is a sesquilinear duality of weight~$w$ on $\psi_{t,\alpha}\cT$. As $\psi_{t,\alpha}S'$ and $\psi_{t,\alpha}S''$ commute with $\rN$, we have $\cN^*\circ\psi_{t,\alpha}\cS=-\psi_{t,\alpha}\cS\circ\cN$ (recall that $\cN=(-i\rN,i\rN)$). Then $\gr_\bbullet^{\rM}\psi_{t,\alpha}\cT$ is a graded triple in the sense given in \T\ref{subsec:dgtriples} and $\gr_\bbullet^{\rM}\psi_{t,\alpha}\cS$ is a (graded) sesquilinear duality of weight~$w$ on it. Last, we see that $\gr_{-2}^{\rM}\cN$ is skewadjoint with respect to $\gr_\bbullet^{\rM}\psi_{t,\alpha}\cS$; in other words, $\gr_\bbullet^{\rM}\psi_{t,\alpha}\cS$ is a Hermitian duality of $(\gr_\bbullet^{\rM}\psi_{t,\alpha}\cT,\gr_{-2}^{\rM}\cN)$. Consequently, \eqref{eq:grlefhermdual} defines a Hermitian duality on the primitive parts.

If for instance $\cM'=\cM''=\cM$ and $C^*=C$, so that $\cS=(\id,\id)$ is a Hermitian duality of weight~$0$, we have $\psi_{t,\alpha}(C)^*=\psi_{t,\alpha}(C)$, and we are in the situation of Example \ref{ex:lefpol}. The sesquilinear pairing on the primitive part
\[
\index{$ppsitcl$@$P\psi_{t,\alpha,_\ell}C$}P\psi_{t,\alpha,_\ell}C:P\gr_\ell^{\rM}\psi_{t,\alpha}\cMS\ootimes_{\cO_\bS} \ov{P\gr_{\ell}^{\rM}\psi_{t,\alpha}\cMS}\to\Dbh{X_0}
\]
is given by the formula
\begin{equation}\label{eq:psipolell}
P\psi_{t,\alpha,_\ell}C\defin\thb^{-\ell}\psi_{t,\alpha,\ell}C((i\rN)^\ell\cbbullet,\ov\cbbullet).
\end{equation}
\end{remarque}

\begin{remarque}\label{rem:restrsmooth}
Assume that $\cT$ is a smooth twistor. In particular, $C$ takes values in $C^\infty$ functions. Then $\psi_{t,-1}\cT$ is equal to the restriction of $\cT$ along $\{t=0\}$ as defined in Definition \ref{def:restrsmooth}.
\end{remarque}

\begin{remarque}\label{rem:psiramif}
Let $f=t^{r}$ for $r\geq1$. Let $i_{f}$ be the graph inclusion of $f$ and let $i:\{t=0\}\hto X$ be the closed inclusion. Similarly to Proposition \ref{prop:strictspetr}, one shows that $\psi_{f,\alpha}C=i_{+}\psi_{t,r\alpha}C$.
\end{remarque}

\subsection{Vanishing cycles and sesquilinear pairings}\label{subsec:vanishsesqui}
If $\cM'$ or $\cM''$ are supported on $X_0$, we have $\psi_{t,\alpha}C=0$ for any $\alpha$. We should therefore also define the ``vanishing cycle analogue'' $\phi_{t,0}C$ in order to recover an interesting sesquilinear form on $\psi_{t,0}\cM'\otimes_{\cO_\bS}\ov{\psi_{t,0}\cM''}$ in any case. We continue to assume in the following that $\cM',\cM''$ are strictly specializable along $\{t=0\}$.

\subsubsection*{The function $\Ichih$}
In the following, we always assume that $\hb$ varies in $\bS$. Let $\chih(\theta)$ be a $C^\infty$ function of the complex variable $\theta$ such that $\chih$ has compact support on $\CC$ and $\chih\equiv1$ near $\theta=0$. For $s$ such that $\reel s>0$, the function
\[
\Ichih(t,s,\hb)\defin\int e^{\ov t \hb/\ov\theta-t/\theta\hb}\module{\theta}^{2(s-1)}\chih(\theta)\itwopi d\theta\wedge d\ov\theta
\]
is continuous with respect to $t$ and holomorphic with respect to $s$ (notice that the exponent $\ov t \hb/\ov\theta-t/\theta\hb$ is purely imaginary, as $\hb\in\bS$). It also varies smoothly with respect to $\hb$. For any $p\in\NN$, the function $\Ichih$, when restricted to the domain $2\reel s>p$, is $C^p$ in $t$ and holomorphic with respect to $s$.

Define $I_{\chih,k,\ell}$ by replacing $\module{\theta}^{2(s-1)}$ with $\theta^k\ov\theta^\ell\module{\theta}^{2(s-1)}$ in the integral defining $\Ichih$; in particular, we have $\Ichih=I_{\chih,0,0}$ and $I_{\chih,k,k}(t,s,\hb)=\Ichih(t,s+k,\hb)$ for any $k\in\ZZ$.

\begin{remarque}
We can also use the coordinate $\tau=1/\theta$ to write $\Ichih(t,s,\hb)$ as
\[
\Ichih(t,s,\hb)=\int e^{\ov{t\tau}\hb-t\tau/\hb}\module{\tau}^{-2(s+1)}\chih(\tau)\itwopi d\tau\wedge d\ov\tau
\]
where now $\chih$ is $C^\infty$, is $\equiv1$ near $\tau=\infty$ and $\equiv0$ near $\tau=0$. It is the Fourier transform of $\module{\tau}^{-2(s+1)}\chih(\tau)$ up to a scaling factor $\hb$: put $\tau=\xi+i\eta$ and $t/2\hb=y+ix$; then $\Ichih(t,s,\hb)=\frac{1}{\pi}\int e^{-i(x\xi+y\eta)}\module{\tau}^{-2(s+1)}\chih(\tau)\,d\xi\wedge d\eta$.

If we denote by $\cF$ the Fourier transform with kernel $e^{\ov{t\tau}\hb-t\tau/\hb}\itwopi d\tau\wedge d\ov\tau$, then the inverse Fourier transform $\cF^{-1}$ has kernel $e^{-\ov{t\tau}\hb+t\tau/\hb}\itwopi dt\wedge d\ov t$.
\end{remarque}

For $\reel s$ large enough, using Stokes formula, we obtain
\begin{align*}
tI_{\chih,k-1,\ell}(t,s,\hb)&=-\hb (s+k) I_{\chih,k,\ell}(t,s,\hb)-\hb I_{\partial\chih/\partial\theta,k+1,\ell}(t,s,\hb)\\
\ov tI_{\chih,k,\ell-1}(t,s,\hb)&=-\ov\hb (s+\ell) I_{\chih,k,\ell}(t,s,\hb)-\ov\hb I_{\partial\chih/\partial\ov\theta,k,\ell+1}(t,s,\hb),
\end{align*}
with $I_{\partial\chih/\partial\theta,k+1,\ell},I_{\partial\chih/\partial\ov\theta,k,\ell+1}\in C^\infty(\CC\times\CC\times\bS)$, holomorphic with respect to $s\in\CC$. In particular we get
\[
\mt^2\Ichih(t,s-1,\hb)=-s^2\Ichih(t,s,\hb)+\cdots,
\]
where ``$\cdots$'' is $C^\infty$ in $(t,s,\hb)$ and holomorphic with respect to $s\in\CC$. This equality holds on $\reel s>1$. This allows one to extend $\Ichih$ as a $C^\infty$ function on $\{t\neq0\}\times\CC\times\bS$, holomorphic with respect to $s$.

For $\reel s>1$, we have
\[
\partiall_t\Ichih(t,s,\hb)=-I_{\chih,-1,0}(t,s,\hb)\quad\text{and}\quad
\ov{\partiall_t}\Ichih(t,s,\hb)=-I_{\chih,0,-1}(t,s,\hb),
\]
hence
\[
t\partiall_t\Ichih=\hb s\Ichih+\hb I_{\partial\chih/\partial\theta,1,0} \quad\text{and}\quad
\ov{t\partiall_t}\Ichih=\ov\hb s\Ichih+\ov\hb I_{\partial\chih/\partial\ov\theta,0,1}.
\]
By analytic extension, these equalities hold on $\{t\neq0\}\times\CC\times\bS$.

\subsubsection*{Definition of $\phi_{t,0}C$}

Let $m',m''$ be local sections (near $(x_o,\hb_o)$ and $(x_o,-\hb_o)$ with $\hb_o\in\bS$) of $V_0\cM',V_0\cM''$ lifting local sections $[m'],[m'']$ of $\psi_{t,0}\cM',\psi_{t,0}\cM''$. Using the previous properties of $\Ichih$, one shows as in Lemma \ref{lem:polesI} that, for any test form $\varphi$ on $\cX_0$ and any compactly supported $C^\infty$ function $\chi(t)$ such that $\chi\equiv1$ near $t=0$, the function
\[
s\mto\big\langle \Ichih(t,s,\hb)C(m',\ov{m''}),\varphi\wedge\chi(t)\itwopi dt\wedge d\ov t\big\rangle
\]
is holomorphic for $\reel s$ big enough and extends as a meromorphic function of $s$ with poles at most on $s=0$ and on sets $s=\gamma\star\hb/\hb$ with $\reel \gamma<0$.

We put
\begin{equation}
\big\langle\phi_{t,0}C([m'],\ov{[m'']}),\varphi\big\rangle\defin \res_{s=0}\big\langle C(m',\ov{m''}),\varphi\wedge\Ichih(t,s,\hb)\chi(t)\itwopi dt\wedge d\ov t\big\rangle.\index{$phif$@$\phi_{f,0}\cT$}
\end{equation}
This residue does not depend on the choice of $\chih$ and $\chi$, nor on the choice of the representatives $m',m''$ in $V_0\cM',V_0\cM''$ (\cf Lemma \ref{lem:coefshol}), and defines a section of $\Dbh{X_0}$. As we can take $\chi$ and $\chih$ real, one obtains that $\phi_{t,0}(C^*)=(\phi_{t,0}C)^*$. Arguing as for $\psi_{t,\alpha}C$, one gets the analogue of \eqref{eq:psiCN}. We define then
\begin{equation}\label{eq:defphiT}
\phi_{t,0}\cT=\phi_{t,0}(\cM',\cM'',C)\defin(\psi_{t,0}\cM',\psi_{t,0}\cM'',\phi_{t,0}C),
\end{equation}
and we have a morphism $\cN:\phi_{t,0}\cT\to\phi_{t,0}\cT(-1)$ in $\RTriples (X_0)$.

\begin{remarque}
Let us explain the definition of $\phi_{t,0}C$. Consider the one-variable distribution with compact support $\langle\chi C(m',\ov{m''}),\varphi\rangle$. Its Fourier transform is a distribution of the variable $\tau$, that we localize near $\tau=\infty$ and to which we apply the functor $\psi_{\theta,-1}$, to obtain $\phi_{t,0}C$. This procedure is similar to a microlocalization with respect to the variable $t$.
\end{remarque}

\subsubsection*{The morphisms $\cCan$ and $\cVar$}
We define
\[
\cCan=(\var,i\can)\quad\text{and}\quad \cVar=(-i\can,\var).
\]
Once they are known to be morphisms in $\RTriples(X_0)$, they clearly satisfy $\cVar\circ\cCan=\cN_{\psi_{t,-1}}$, $\cCan\circ\cVar=\cN_{\phi_{t,0}}$.

\begin{lemme}\label{lem:cCancVar}
The morphisms $\cCan$ and $\cVar$ are morphisms in $\RTriples(X_0)$:
\[
\psi_{t,-1}\cT\To{\cCan\index{$can$@$\cCan$}}\phi_{t,0}\cT(-1/2),\quad \phi_{t,0}\cT(1/2)\To{\cVar\index{$var$@$\cVar$}}\psi_{t,-1}\cT.
\]
\end{lemme}

\begin{remarque}[Behaviour with respect to adjunction]\label{rem:canvaradj}
Let $\cS:\cT\to\cT^*(-w)$ be a sesquilinear duality of weight~$w$. Using the canonical isomorphism \eqref{eq:adjcanonical} $(\phi_{t,0}\cT(1/2))^*\isom \phi_{t,0}\cT^*(-1/2)$ given by $(\id_{\psi_{t,0}\cM''},-\id_{\psi_{t,0}\cM'})$, we get a commutative diagram
\[
\xymatrix@C=2cm{
\psi_{t,-1}\cT\ar[r]^-{\psi_{t,-1}\cS}\ar[d]_-{\cCan}&\psi_{t,-1}\cT^*(-w) \ar[d]^-{\cVar^*}\\
\phi_{t,0}\cT(-1/2)\ar[r]^-{\psi_{t,0}\cS}&\phi_{t,0}\cT^*(-w-1/2)
}
\]
and an adjoint anticommutative diagram.
\end{remarque}

\begin{proof}[Proof of Lemma \ref{lem:cCancVar}]
Let us show that $\cCan$ is a morphism in $\RTriples(X_0)$, the proof for $\cVar$ being similar. Let $[m'_0]$ (\resp $[m''_{-1}]$) be a local section of $\psi_{t,0}\cM'$ (\resp $\psi_{t,-1}\cM''$). We have to show that
\begin{multline}\label{eq:psiphiC}
\res_{s=0}\hb\big\langle C(m'_0,\ov{m''_{-1}}),\varphi\wedge\ov\partiall_t\Ichih(t,s,\hb)\chi(t)\itwopi dt\wedge d\ov t\big\rangle\\
=\res_{s=-1}\big\langle C(m'_0,\ov{m''_{-1}}),\varphi\wedge t\mt^{2s}\chi(t)\itwopi dt\wedge d\ov t\big\rangle.
\end{multline}
We can replace $\ov\partiall_t\Ichih(t,s,\hb)$ with $-I_{\chih,1,0}(t,s-1,\hb)$, so that the left-hand term in \eqref{eq:psiphiC} is
\begin{equation}\label{eq:Ichih}
\res_{s=-1}\big\langle C(m'_0,\ov{m''_{-1}}),\varphi\wedge(-\hb I_{\chih,1,0}(t,s,\hb))\chi(t)\itwopi dt\wedge d\ov t\big\rangle.
\end{equation}

Denote by $T$ the one-variable distribution $\langle\chi C(m'_0,\ov{m''_{-1}}),\varphi\rangle$ obtained by integration in the $X_0$ direction. It has compact support by definition of $\chi$. Therefore, its Fourier transform $\cF T$ is a $C^\infty$-function of $\tau,\hb$, which has slow growth, as well as all its derivatives, when $\tau\to\infty$. The function in \eqref{eq:Ichih} is then written as
\begin{equation}\label{eq:FTtau}
-\hb\int\cF T(\tau,\hb)\cdot\tau^{-1}\module{\tau}^{-2(s+1)}\chih(\tau)\itwopi d\tau\wedge d\ov\tau.
\end{equation}
On the other hand, the function in the RHS of \eqref{eq:psiphiC} is
\begin{equation}\label{eq:FTt}
\begin{aligned}
\langle T,t\mt^{2s}\chi(t)\itwopi dt\wedge d\ov t\rangle&=\big\langle \cF T,\cF^{-1}(t\mt^{2s}\chi(t)\itwopi dt\wedge d\ov t)\big\rangle\\
&=\int\cF T(\tau,\hb)\cdot\cF^{-1}(t\mt^{2s}\chi)\itwopi d\tau\wedge d\ov\tau
\end{aligned}
\end{equation}
(in order to get this expression, we replace $\chi$ with $\chi^2$ in \eqref{eq:psiphiC}, which does not change the residue, as previously remarked).

\subsubsection*{The function $\wh I_{\chi,1,0}(\tau,s,\hb)$}
Let us state some properties of the function $\wh I_{\chi,1,0}(\tau,s,\hb)\defin \cF^{-1}(t\mt^{2s}\chi)$.
\begin{enumerate}
\item\label{item:whIchi1}
Denote by $\wh I_{\chi,k,\ell}(\tau,s,\hb)$ ($k,\ell\in\ZZ$) the function obtained by integrating $\mt^{2s}t^k\ov t^\ell$. Then, for any $s\in\CC$ with $\reel(s+1+(k+\ell)/2)>0$ and any $\hb\in\bS$, the function $(\tau,s,\hb)\mto \wh I_{\chi,k,\ell}(\tau,s,\hb)$ is $C^\infty$, depends holomorphically on $s$, and satisfies $\lim_{\tau\to\infty}\wh I_{\chi,k,\ell}(\tau,s,\hb)=0$ locally uniformly with respect to $s,\hb$ (apply the classical Riemann-Lebesgue lemma saying that the Fourier transform of a function in $L^1$ is continuous and tends to $0$ at infinity).

\item\label{item:whIchi2}
We have
\begin{equation}\label{eq:whIchi2}
\begin{aligned}
\tau\wh I_{\chi,k,\ell}&=-\hb(s+k)\wh I_{\chi,k-1,\ell}-\hb \wh I_{\partial\chi/\partial t,k,\ell}\quad&\partiall_\tau\wh I_{\chi,k,\ell}&=\wh I_{\chi,k+1,\ell}\\
\ov\tau\wh I_{\chi,k,\ell}&=-\ov\hb(s+\ell)\wh I_{\chi,k,\ell-1}-\ov\hb \wh I_{\partial\chi/\partial\ov t,k,\ell}&\ov\partiall_\tau\wh I_{\chi,k,\ell}&=\wh I_{\chi,k,\ell+1},
\end{aligned}
\end{equation}
where the equalities hold on the common domain of definition (with respect to $s$) of the functions involved. Notice that the functions $\wh I_{\partial\tau,k,\ell}$ and $\wh I_{\partial\chi/\partial\ov t,k,\ell}$ are $C^\infty$ on $\PP^1\times\CC\times\bS$, depend holomorphically on $s$, and are infinitely flat at $\tau=\infty$ (because $t^k\ov t^\ell\mt^{2s}\partial_{t,\ov t}\chi$ is $C^\infty$ in $t$ with compact support, and holomorphic with respect to $s$, so that its Fourier transform is in the Schwartz class, holomorphically with respect to~$s$).

It follows that, for $\reel(s+1)+(k+\ell)/2>0$, we have
\begin{equation}\label{eq:tdtwhIchi}
\begin{aligned}
\tau\partiall_\tau\wh I_{\chi,k,\ell}&=-\hb(s+k+1)\wh I_{\chi,k,\ell}-\hb \wh I_{\partial\chi/\partial t,k+1,\ell},\\
\ov{\tau\partiall_\tau}\wh I_{\chi,k,\ell}&=-\ov\hb(s+\ell+1)\wh I_{\chi,k,\ell}-\ov\hb \wh I_{\partial\chi/\partial\ov t,k,\ell+1}.
\end{aligned}
\end{equation}

\item\label{item:whIchi3}
Consider the variable $\theta=\tau^{-1}$ with corresponding derivation $\partial_\theta=-\tau^2\partial_\tau$, and write $\wh I_{\chi,k,\ell}(\theta,s,\hb)$ the function $\wh I_{\chi,k,\ell}$ in this variable. Then, for any $p\geq0$, any $s\in\CC$ with \hbox{$\reel (s+1+(k+\ell)/2)>p$} and any $\hb\in\bS$, all derivatives up to order $p$ of $\wh I_{\chi,k,\ell}(\theta,s,\hb)$ with respect to $\theta$ tend to $0$ when $\theta\to0$, locally uniformly with respect to $s,\hb$ (use \eqref{eq:tdtwhIchi} and \eqref{eq:whIchi2}); in particular, $\wh I_{\chi,k,\ell}(\tau,s,\hb)$ extends as a function of class $C^p$ on $\PP^1\times\{\reel (s+1+(k+\ell)/2)>p\}\times\bS$, holomorphic with respect to $s$.
\end{enumerate}

\medskip
The function $\wh I_{\chi,1,0}(\tau,s,\hb)$ is $C^\infty$ in $\tau$ and holomorphic in $s$ on $\{s\mid\reel s>-3/2\}$. Using the function $\wh\chi(\tau)$ as above, we conclude that the integral
\begin{equation}\label{eq:unmoinschi}
\int\cF T(\tau,\hb)\cdot\cF^{-1}(t\mt^{2s}\chi)(1-\chih(\tau))\itwopi d\tau\wedge d\ov\tau
\end{equation}
is holomorphic with respect to $s$ for $\reel s>-3/2$. It can thus be neglected when computing the residue at $s=-1$. The question reduces therefore to the comparison of $\wh I_{\chi,1,0}(\tau,s,\hb)$ and $\tau^{-1}\module{\tau}^{-2(s+1)}$ when $\tau\to\infty$.

Put $\wh J_{\chi,1,0}(\tau,s,\hb)=\tau\module{\tau}^{2(s+1)}\wh I_{\chi,1,0}(\tau,s,\hb)$. Then, by \eqref{eq:tdtwhIchi}, we have
\[
\tau\frac{\partial\wh J_{\chi,1,0}}{\partial \tau}=-\wh J_{\partial\chi/\partial t,1,0},\quad \ov{\tau}\frac{\partial\wh J_{\chi,1,0}}{\partial \ov{\tau}}=-\wh J_{\partial\chi/\partial\ov t,0,1},
\]
and both functions $\wh J_{\partial\chi/\partial t,1,0}$ and $\wh J_{\partial\chi/\partial\ov t,0,1}$ extend as $C^\infty$ functions, infinitely flat at $\tau=\infty$ and holomorphic with respect to $s\in\CC$. Put
\[
\whKchi(\tau,s,\hb)=-\int_0^1\big[\wh J_{\partial\chi/\partial t,1,0}(\lambda \tau,s,\hb)+\wh J_{\partial\chi/\partial\ov t,0,1}(\lambda \tau,s,\hb)\big]d\lambda.
\]
Then $\whKchi$ is of the same kind.

\begin{lemme}\label{lem:Bessel}
For any $s$ in the strip $\reel (s+1)\in{}]-1,-1/4[$, the function $\tau\mto\wh J_{\chi,1,0}(\tau,s,\hb)$ satisfies
\[
\lim_{\tau\to\infty}\wh J_{\chi,1,0}(\tau,s,\hb)=-\hb\,\frac{\Gamma(s+2)}{\Gamma(-s)}.
\]
\end{lemme}

\begin{proof}
We can assume that $\chi$ is a $C^\infty$ function of $\mt^2$, that we still write $\chi(\mt^2)$. For simplicity, we assume that $\chi\equiv1$ for $\mt\leq1$. Then the limit of $\wh J_{\chi,1,0}$ is also equal to the limit of the integral
\[
J(\tau,s,\hb)=\int_{\mt\leq1} e^{-\ov {t\tau} \hb+ t\tau/\hb} \, t\tau\module{t\tau}^{2(s+1)}\,\itwopi \frac{dt}{t}\wedge \frac{d\ov t}{\ov t}.
\]
By a simple change of variables, we have
\[
J(\tau,s,\hb)=\hb\int_{\module{u}\leq\module{\tau}} e^{2i\im u}u\module{u}^{2s}\itwopi du\wedge d\ov u.
\]
Using the Bessel function $J_{\pm1}(x)=\frac{1}{2\pi}\int_0^{2\pi}e^{-ix\sin y}e^{\pm iy}dy$, we can write
\begin{align*}
J(\tau,s,\hb)&=2\hb\int_{\rho\leq\module{\tau}} J_{-1}(2\rho)\rho^{2(s+1)} d\rho\\ &=-2^{-2(s+1)}\hb\int_{\rho\leq2\module{\tau}} J_1(\rho)\rho^{2(s+1)}d\rho,\quad\text{as } J_1=-J_{-1}.
\end{align*}
For $\reel (s+1)\in{}]-1,-1/4[$, the limit when $\module{\tau}$ of the previous integral is equal to \hbox{$2^{2(s+1)}\Gamma(s+2)/\Gamma(-s)$} (\cf \cite[\T13.24, p\ptbl 391]{Watson22}).
\end{proof}

From Lemma \ref{lem:Bessel}, we can write, on the strip $\reel (s+1)\in{}]-1/2,-1/4[$,
\begin{equation}\label{eq:IGamma}
\wh I_{\chi,1,0}(\tau,s,\hb)=-\hb\tau^{-1}\module{\tau}^{-2(s+1)} \,\frac{\Gamma(s+2)}{\Gamma(-s)}+K_{\chi}(\tau,s,\hb)
\end{equation}
where $K_{\chi}(\tau,s,\hb)=-\hb\tau^{-1}\module{\tau}^{-2(s+1)}\whKchi$ is $C^\infty$ on $\CC\times\CC\times\bS$, infinitely flat at $\tau=\infty$ and holomorphic with respect to~$s$. For any $p\geq0$, apply $(\partiall_\tau\partiall_{\ov\tau})^p$ to the previous equality restricted to $\tau\neq0$ (where both sides are $C^\infty$ in $\tau$ and holomorphic with respect to $s$; preferably, multiply both sides by $\wh\chi(\tau)$), to get, for $s$ in the same strip,
\[
\wh I_{\chi,1,0}(\tau,s+p,\hb)=-\hb\tau^{-1}\module{\tau}^{-2(s+p+1)} \,\frac{\Gamma(s+p+2)}{\Gamma(-s-p)}+(\partiall_\tau\partiall_{\ov\tau})^pK_{\chi}(\tau,s,\hb)
\]
where the last term remains infinitely flat at $\tau=\infty$. It follows that \eqref{eq:IGamma} remains true on any strip $\reel (s+1)\in{}]p-1/2,p-1/4[$ with $p\geq0$ and a function $K^{(p)}_\chi$ instead of~$K_\chi$.

Choose $p$ such that the two the meromorphic functions considered in \eqref{eq:psiphiC} are holomorphic on the strip $\reel (s+1)\in{}]p-1/2,p-1/4[$. The difference between $\Gamma(s+2)/\Gamma(-s)$ times the function in the LHS and the function in the RHS coincides, on this strip, with a holomorphic function on the half-plane $\{s\mid\reel s>-3/2\}$ (taking into account \eqref{eq:unmoinschi} and $K^{(p)}_\chi$). It is then equal to it on this whole half-plane, hence has residue $0$ at $s=-1$.
\end{proof}

Let us emphasize two cases:

\begin{enumerate}
\item
Denote by $i$ the inclusion $\{t=0\}\hto X$. Remark that, if $\cM'$ and $\cM''$ are in $\cS^2(X,t)$ and are supported on $\{t=0\}$, so that $\cM'=i_+\cM'_0$ and $\cM''=i_+\cM''_0$, then any sesquilinear pairing $C$ on $\cMS'\otimes_{\cO_\bS}\ov{\cMS''}$ is equal to $i_{++}C_0$ for some sesquilinear pairing $C_0$ on $\cM'_{0|\bS}\otimes_{\cO_\bS}\ov{\cM''_{0|\bS}}$. Indeed, by $\cR_{(X,\ov X),\bS}$-linearity, $C$ is determined by its restriction to $\cM'_{0|\bS}\otimes_{\cO_\bS}\ov{\cM''_{0|\bS}}$; conclude by using that $tC(m'_0,\ov{m''_0})=C(tm'_0,\ov{m''_0})=0$.

We have $\cM'_0=\psi_{t,0}\cM'$ and $\cM''_0=\psi_{t,0}\cM''$. Moreover:

\begin{lemme}\label{lem:phiCC0}
The pairing $\phi_{t,0}C$ is equal to $C_0$.
\end{lemme}

\begin{proof}
By definition, as $\chi(0)=1$, we have for $\reel s\gg0$,
\[
\big\langle C(m'_0,\ov{m''_0}),\varphi\wedge\Ichih(t,s,\hb)\chi(t)\itwopi dt\wedge d\ov t\big\rangle=\big\langle C_0(m'_0,\ov{m''_0}),\varphi\big\rangle\cdot \Ichih(0,s,\hb).
\]
As $\chih(0)=1$, we have $\res_{s=0}\Ichih(0,s,\hb)=1$.
\end{proof}

\item
Assume that $\can$ is onto.

\begin{lemme}\label{lem:phiC}
Let $m'_0,m''_0$ be local sections of $V_0\cM',V_0\cM''$ lifting local sections $[m'_0],[m''_0]$ of $\psi_{t,0}\cM',\psi_{t,0}\cM''$. Then
\[
\big\langle\phi_{t,0}C([m'_0],\ov{[m''_0]}),\cbbullet\big\rangle=\res_{s=0}\frac{-1}{s}\big\langle \mt^{2s}C(m'_0,\ov{m''_0}),\cbbullet\wedge\chi(t)\itwopi dt\wedge d\ov t\big\rangle.
\]
\end{lemme}

\begin{proof}
There exists a local section $[m''_{-1}]$ of $\psi_{t,-1}\cM''$ such that $[m''_0]=i\can[m''_{-1}]$. We have
\begin{align*}
\big\langle \mt^{2(s+1)}C(m'_0,&\ov{m''_0}),\cbbullet\wedge\chi(t)\itwopi dt\wedge d\ov t\big\rangle
= \big\langle \mt^{2s}C(tm'_0,\ov{tm''_0}),\cbbullet\wedge\chi(t)\itwopi dt\wedge d\ov t\big\rangle\\
&= \big\langle \mt^{2s}C(tm'_0,\ov{i\rN m''_{-1}}),\cbbullet\wedge\chi(t)\itwopi dt\wedge d\ov t\big\rangle\\
&= \big\langle (-i\ov{\partiall_tt}\mt^{2s})C(tm'_0,\ov{m''_{-1}}),\cbbullet\wedge\chi(t)\itwopi dt\wedge d\ov t\big\rangle+J(s,\hb)\\
&= -i\ov\hb(s+1)\big\langle \mt^{2s}C(tm'_0,\ov{m''_{-1}}),\cbbullet\wedge\chi(t)\itwopi dt\wedge d\ov t\big\rangle+J(s,\hb)\\
&= -\thb^{-1}(s+1)\big\langle \mt^{2s}C(tm'_0,\ov{m''_{-1}}),\cbbullet\wedge\chi(t)\itwopi dt\wedge d\ov t\big\rangle+J(s,\hb),
\end{align*}
where $J(s,\hb)$ is meromorphic with respect to $s$ and has no pole along $s=-1$. Therefore, as $\phi_{t,0}C([m'_0],\ov{[m''_0]})=\phi_{t,0}C([m'_0],\ov{i\can[m''_{-1}]})=\thb^{-1}\psi_{t,-1}C(\var[m'_0],\ov{[m''_{-1}]})$, we get
\begin{align*}
\big\langle\phi_{t,0}C([m'_0],\ov{[m''_0]}),\cbbullet\big\rangle&=\thb^{-1}\res_{s=-1}\big\langle \mt^{2s}C(tm'_0,\ov{m''_{-1}}),\cbbullet\wedge\chi(t)\itwopi dt\wedge d\ov t\big\rangle\\
&=\res_{s=-1}\frac{-1}{s+1}\big\langle \mt^{2(s+1)}C(m'_0,\ov{m''_0}),\cbbullet\wedge\chi(t)\itwopi dt\wedge d\ov t\big\rangle\\
&=\res_{s=0}\frac{-1}{s}\big\langle \mt^{2s}C(m'_0,\ov{m''_0}),\cbbullet\wedge\chi(t)\itwopi dt\wedge d\ov t\big\rangle.\qedhere
\end{align*}
\end{proof}
\end{enumerate}

\begin{corollaire}\label{cor:RTriplesdecomposables}
Let $\cT=(\cM',\cM'',C)$ be an object of $\RTriples(X)$. Assume that $\cM',\cM''$ are strictly specializable along $\{t=0\}$. The following properties are equivalent:
\begin{enumerate}
\item
$\phi_{t,0}\cT=\im\cCan\oplus\ker\cVar$ in $\RTriples(X)$,
\item
$\cT=\cT_1\oplus\cT_2$ in $\RTriples(X)$, with $\cT_2$ supported on $\{t=0\}$ and $\cT_1$ being such that its $\cCan$ is onto and its $\cVar$ is injective.
\end{enumerate}
\end{corollaire}

\begin{proof}
The part for $\cM',\cM''$ is Proposition \ref{prop:canvar}. That $\phi_{t,0}C$ decomposes is proved as in Proposition \ref{prop:Cdecomposable}.
\end{proof}

\subsection{Direct images and specialization of sesquilinear pairings}
We take the notation used in Theorem \ref{th:imdirspe}.

\begin{corollaire}\label{cor:imdirstrictspeRT}
Let $\cT=(\cM',\cM'',C)$ be an object of $\RTriples(X\times\CC)$. Assume that $\cM',\cM''$ satisfy the conditions in Theorem \ref{th:imdirstrictspe}. Then, for any~$\alpha$ with $\reel\alpha\in[-1,0[$, we have $\Psi_{t,\alpha}\cH^i(F_\dag\cT)=\cH^i(f_\dag\Psi_{t,\alpha}\cT)$. Moreover, we have $\phi_{t,0}\cH^i(F_\dag\cT)=\cH^i(f_\dag\phi_{t,0}\cT)$ and, with obvious notation, $\cCan_{\cH^i(F_\dag\cT)}=\cH^i(f_\dag\cCan)$, $\cVar_{\cH^i(F_\dag\cT)}=\cH^i(f_\dag\cVar)$.
\end{corollaire}

\begin{proof}
Apply Theorem \ref{th:imdirstrictspe} for $\cM',\cM''$ and $\cCan,\cVar$. It remains to controlling the behaviour of $\psi_{t,\alpha}C,\phi_{t,0}C$ under $F_\dag$. Now the result is a direct consequence of the definition of $\psi_{t,\alpha}C,\phi_{t,0}C$, as we can compute with local sections of $V_a\cM',V_a\cM''$, knowing that, for $\cM=\cM'$ or $\cM''$, the filtration $V_\bbullet\cH^i(F_\dag\cM)$ is $\cH^i(f_\dag V_\bbullet\cM)$.
\end{proof}

\section{Noncharacteristic inverse image}\label{sec:nonchar}

\subsection{Noncharacteristic and strictly noncharacteristic $\cR_\cX$-modules along a submanifold}
Let $\cM$ be a holonomic $\cR_\cX$-module with characteristic variety $\Char\cM$ contained in $\Lambda \times\Omega_0$, where $\Lambda\subset T^*X$ is Lagrangian. Let $Z\subset X$ be a submanifold of $X$ and denote by $i:Z\hto X$ the inclusion. We say that $\cM$ is \emph{noncharacteristic} along $Z$ if $T^*_{Z}X\cap\Lambda\subset T^*_XX$ for some choice of $\Lambda$ as above.

Locally on $Z$, we may choose a smooth map $\bmt=(t_1,\dots,t_p):X\to\CC^p$ such that $Z=\bmt^{-1}(0)$ and we may view $\bmt$ as a projection. We may therefore consider the sheaf $\cR_{\cX/\CC^p}$ of relative differential operators with respect to the projection $\bmt$. The following is classical and easy:

\begin{lemme}\label{lem:noncar}
If $\cM$ is noncharacteristic along $Z$, then $\cM$ is (locally on $Z$) $\cR_{\cX/\CC^p}$-coherent. If $Z$ has codimension one, then $\cM$ is regular along $Z$.
\end{lemme}

\begin{proof}
Indeed, if $\cM$ is noncharacteristic along $Z$, then any local good filtration $F_\bbullet\cM$ of $\cM$ as a $\cR_\cX$-module is such that $\gr^F\cM$ is $\gr^F\cR_{\cX/\CC^p}$-coherent.
\end{proof}

\begin{definition}
The $\cR_\cX$-module $\cM$ is \emph{strictly noncharacteristic} along $Z$ if it is noncharacteristic along $Z$ and the (ordinary) restriction $\cO_{\cZ}\otimes_{\cO_\cX}\cM_{|Z}$ is strict.
\end{definition}

\begin{remarque}\label{rem:iminvnoncar}
If $f:Z\to X$ is any morphism between smooth complex manifold, we may similarly define, for a holonomic $\cR_\cX$-module $\cM$, what ``(strictly) noncharacteristic with respect to $f$'' means. Decompose $f$ as en embedding followed by a projection. Then $\cM$ is always strictly non characteristic with respect to the projection. Hence, in practice, it is enough to check this property for embeddings.
\end{remarque}

\begin{lemme}\label{lem:strictnoncar}
Let $Z$ be a smooth hypersurface. Assume that $\cM$ is strictly noncharacteristic along $Z$. Then,
\begin{enumerate}
\item\label{lem:strictnoncar3}
$\cM$ is strictly specializable along $Z$,
\item\label{lem:strictnoncar4}
we have $\cO_{\cZ}\otimes_{\cO_\cX}\cM_{|Z}=\cO_{\cZ}\otimes^{\bL}_{\cO_\cX}\cM_{|Z}$.
\end{enumerate}
\end{lemme}

\begin{proof}
We may work locally on $X$. Assume that $Z$ is defined by a local equation $t=0$. As $\cM$ is noncharacteristic along $Z$, Lemma \ref{lem:noncar} shows that $\cM$ is specializable and regular along $t=0$, and we may choose as a good $V$-filtration the filtration given by $V_{-1}\cM=\cM$, $V_{-k-1}\cM=t^k\cM$ for $k\geq0$ and $\gr^V_a\cM=0$ for $a\not\in -\NN^*$ (here the filtration is independent of the choice of $\hb_o$ as the indices are real). The restriction $i^+\cM$ is equal to $\gr_{-1}^V\cM$ and is still holonomic with characteristic variety contained in $\Lambda_0\times\Omega_0$, where $\Lambda_0$ is the image of $\Lambda_{|Z}$ by the cotangent map $T^*i:T^*X_{|Z}\to T^*Z$. Moreover, the action of $t\partiall_t$ on $\gr_{-1}^V\cM$ vanishes, because for a local section $m$ of $\cM$, $\partiall_tm$ is also a local section of $\cM$.

With the assumption of strictness of $\cM/t\cM$, we conclude that $\cM$ is strictly specializable and regular along $Z$. In such a situation, the $V$-filtration defined above is the Malgrange-Kashiwara filtration and we have $\psi_{t,-1}\cM=\gr_{-1}^V\cM=\cM/t\cM$. Moreover, $t:\cM\to\cM$ is injective, because $\cM=V_{-1}\cM$ (\cf Remark \ref{rem:psi}\eqref{rem:psia}), \ie
$$
\cO_{\cZ}\otimes_{\cO_\cX}\cM_{|Z}=\cO_{\cZ}\otimes^{\bL}_{\cO_\cX}\cM_{|Z}.$$
Last, remark that $\can$ and $\var$ are both equal to $0$, as $\psi_{t,0}\cM=0$.
\end{proof}

\subsubsection*{An adjunction morphism}
Let $\pi:\wt X\to X$ be a proper analytic map between complex analytic manifolds, which is an isomorphism almost everywhere (say that $\pi$ is a \emph{proper modification of $X$}). Assume that $\pi=p\circ i$, where $i:\wt X\hto X\times\PP^1$ is a closed inclusion and $p:X\times\PP^1\to X$ is the projection. Let $\cM$ be an holonomic $\cR_\cX$-module which is \emph{strictly noncharacteristic} with respect to $\pi$. By Lemma \ref{lem:strictnoncar}, $\pi^+\cM$ is a holonomic $\cR_{\wt\cX}$-module.

\begin{lemme}\label{lem:adj}
Under these conditions, there is a natural adjunction morphism $\iota:\cM\to\pi_+^0\pi^+\cM$.
\end{lemme}

\begin{proof}
Put $n=\dim \wt X=\dim X$. The right $\cR_{\wt\cX}$-module associated with $\pi^+\cM$ is $\omega_{\wt\cX}\otimes_{\pi^{-1}\cO_{\cX}}\pi^{-1}\cM$. Using the contraction isomorphism \eqref{eq:contraction}, we identify the complex
\[
\Big(\omega_{\wt\cX}\ootimes_{\pi^{-1}\cO_\cX}\pi^{-1}\cM\Big)\ootimes_{\cR_{\wt\cX}}\Sp^\cbbullet_{\wt\cX\to\cX}=\Big(\omega_{\wt\cX}\ootimes_{\pi^{-1}\cO_\cX}\pi^{-1}\cM\Big)\ootimes_{\cO_{\wt\cX}}\wedge^\cbbullet\Theta_{\wt\cX}\ootimes_{\pi^{-1}\cO_\cX}\pi^{-1}\cR_\cX
\]
with the de~Rham complex
\[
\Omega_{\wt\cX}^{n+\cbbullet}\ootimes_{\pi^{-1}\cO_\cX}\pi^{-1}(\cM\otimes_{\cO_\cX}\cR_\cX)
\]
of the (inverse image of) the left $\cR_\cX$-module $\cM\otimes_{\cO_\cX}\cR_\cX$, the right $\cR_\cX$-structure being trivially induced by that on $\cR_\cX$. Using the isomorphism \eqref{eq:lefrig}, this complex is isomorphic to the de~Rham complex
\[
\Omega_{\wt\cX}^{n+\cbbullet}\ootimes_{\pi^{-1}\cO_\cX}\pi^{-1}(\cR_\cX\otimes_{\cO_\cX}\cM)\simeq\cE_{\wt\cX}^{n+\cbbullet}\ootimes_{\pi^{-1}\cO_\cX}\pi^{-1}(\cR_\cX\otimes_{\cO_\cX}\cM).
\]
We now have a morphism
\begin{multline*}
\cM^r=\omega_\cX\ootimes_{\cO_\cX}\cM\simeq\Omega^{n+\cbbullet}_\cX\ootimes_{\cO_\cX}(\cR_\cX\otimes_{\cO_\cX}\cM)\to (\pi_*\cE_{\wt\cX}^{n+\cbbullet})\ootimes_{\cO_\cX}(\cR_\cX\otimes_{\cO_\cX}\cM)\\
\to\pi_*\Big(\cE_{\wt\cX}^{n+\cbbullet}\ootimes_{\pi^{-1}\cO_\cX}\pi^{-1}(\cR_\cX\otimes_{\cO_\cX}\cM)\Big) \simeq \pi_+(\pi^+\cM)^r,
\end{multline*}
hence a morphism $\cM^r\to\pi^0_+(\pi^+\cM)^r$.
\end{proof}

\subsection{Noncharacteristic inverse image of a sesquilinear pairing} \label{subsec:noncarC}
Consider first the case of the inclusion of a smooth hypersurface $Z=\{t=0\}$. If $C$ is a sesquilinear pairing, then $\psi_{t,-1}C$ is defined by the formula of Definition \ref{def:psiC}. Notice that, by applying the same argument as in Lemma \ref{lem:polesI}, for any local sections $m'$ of $\cM'$ and $m''$ of $\cM''$, the function $\cI_{C(m',\ov{m''}),\varphi}^{(0)}(s)$ has at most simple poles at $s=-1,-2,\dots$ and no other poles.

If $\cT=(\cM',\cM'',C)$ is an object of $\RTriples(X)$ with $\cM',\cM''$ holonomic, and if $\cM',\cM''$ are strictly noncharacteristic along $X_0$, then one may define $i^+\cT$ as $\psi_{t,-1}\cT$.

The following result will be useful in the proof of Theorem \ref{th:imdirtwistor} (\cf \T\ref{sec:n1n0}\eqref{n0step3}).

\begin{proposition}\label{prop:sesquiunidef}
Assume that $\cM$ is strictly noncharacteristic along the smooth hypersurface $Z\subset X$. Let $j:X\moins Z\hto X$ denote the open inclusion. Given any sesquilinear pairing $C^o:j^*\cMS\otimes_{\cO_\bS}\ov{j^*\cMS}\to\Db_{X\moins Z\times\bS/\bS}$, there exists at most one sesquilinear pairing $C$ on $\cMS$ which extends~$C^o$.
\end{proposition}

\begin{proof}
The question is local on $X\times\bS$, so let $x_o\in Z$, $\hb_o\in\bS$. Assume that we have a coordinate system $(t,x')$ such that $Z=\{t=0\}$. We may consider relative differential forms of maximal degree, namely forms $\psi=a\cdot dx'\wedge d\ov x'$ where $a$ is a section of $\cC^\infty_{X_\RR\times\bS}$. Let $m',m''$ be sections of $\cM_{(x_o,\hb_o)},\cM_{(x_o,-\hb_o)}$. Assume that $C$ is a sesquilinear pairing on $\cM$. Then, for a relative differential form $\psi$ of maximal degree supported in $\nb(x_o,\hb_o)$, $\langle C(m',\ov{m''}),\psi\rangle$ is the section of $\Dbh{D}$, where $D=\{\mt<R\}$, defined by
\[
\eta(t,\hb)\itwopi dt\wedge d\ov t\mto\llangle C(m',\ov{m''}),\psi\wedge\eta(t,\hb)\itwopi dt\wedge d\ov t\rrangle.
\]
We can view $\langle C(m',\ov{m''}),\psi\rangle$ as a distribution on $D_\RR$ with values in the Banach space $C^0(\bS)$, using the formula above.

Assume that we have two extensions $C_1,C_2$ of $C^o$. For $m',m''$ as above, put $u=C_1(m',\ov{m''})-C_2(m',\ov{m''})$. It is then enough to prove that, for any relative form $\psi$, the distribution $\langle u,\psi\rangle$ on the disc is equal to $0$: indeed, for any test form $\varphi$ supported in $\nb(x_o,\hb_o)$, we can write $\varphi=\psi\wedge\chi$, for some relative form $\psi$ and some test form $\chi$ on the disc; then, $\langle u,\varphi\rangle=\big\langle\langle u,\psi\rangle,\chi\big\rangle=0$.

Denote by $p$ the order of the distribution $u$ on $\nb(x_o,\hb_o)$. We have
\begin{equation}\label{eq:enough}
\langle u,\psi\rangle=\sum_{0\leq a+b\leq p} \lambda_{a,b}(\psi)\partiall_t^a\partiall_{\ov t}^b\delta_0,
\end{equation}
for some $\lambda_{a,b}(\psi)\in\cC^0_{\bS,\hb_o}$, where $\delta_0\in\Dbh{D}$ denotes the Dirac distribution at $t=0$. Let us show that all $\lambda_{a,b}(\psi)$ vanish identically. This is true if $\psi$ vanishes at order $\geq p+1$ along $\{t=0\}$. We may therefore assume that $\psi$ does not depend on $t,\ov t$. Using the Bernstein equation for $m'$, we obtain, for a convenient choice of $N$,
\[
B(t\partiall_t)\cdot m'=\Big[\prod_{k=0}^N(t\partiall_t-k\hb)\Big]\cdot m'=t^{p+1}\sum_j(t\partiall_t)^j P_j(t,x',\partiall_x')\cdot m'.
\]
According to \eqref{eq:enough}, the coefficient of $\partiall_t^a\partiall_{\ov t}^b\delta_0$ in $B(t\partiall_t)\cdot\langle u,\psi\rangle$ is $\mu_a(\hb)\lambda_{a,b}(\psi)$, with $\mu_a(\hb)$ invertible on $\bS$. On the other hand, for any $j$, by \eqref{eq:enough} applied to $\psi\cdot P_j$, $t^{p+1}\llangle u,\psi\cdot P_j\rrangle=0$, hence $t^{p+1}(t\partiall_t)^j\llangle u,\psi\cdot P_j\rrangle=0$. Therefore, all $\lambda_{a,b}(\psi)$ vanish.
\end{proof}

\subsection{Specialization along two normally crossing divisors}
We will need the results below in the proof of Theorem \ref{th:imdirtwistor}, \T\ref{sec:n1n0}. The result is an adaptation of \cite[\S\T3.5.11--3.5.18]{MSaito86}.

Assume that we are in the following situation: let $Y=Y_1\cup Y_2$ be a normal crossing divisor in a smooth manifold $X$ and let $\cT=(\cM',\cM'',C)$ be an object of $\RTriples(X)$. Assume that $\cM',\cM''$ are holonomic with characteristic variety contained in $\Lambda\times\Omega_0$, for some Lagrangian variety $\Lambda\subset T^*X$. Assume also that $\cM',\cM''$ are strictly noncharacteristic along $Y_1$, $Y_2$ and $Z=Y_1\cap Y_2$ in a neighbourhood of~$Z$. We will work in local coordinates near a point of $Z$: we put $Y_1=\{x_1=0\}$, $Y_2=\{x_2=0\}$.

\begin{lemme}\label{lem:nonbichar}
Under these conditions, for $\cM=\cM'$ or $\cM''$ and near each point of $Z$,
\begin{itemize}
\item
$\cM$ is strictly specializable along $Y_1$and $Y_2$,
\item
$\psi_{x_j,\alpha}\cM$ vanishes for $\alpha\not\in -\NN^*$ and $j=1,2$,
\item
$\psi_{x_1,-1}\cM=\cM/x_1\cM$ is strictly specializable along $\{x_2=0\}$ and conversely,
\item
we have $\psi_{x_2,-1}\psi_{x_1,-1}\cM=\psi_{x_1,-1}\psi_{x_2,-1}\cM=\cM/x_1x_2\cM$,
\item
for any local sections $m'$ of $\cM'$ and $m''$ of $\cM''$, the two-variable Mellin transform $\langle C(m',\ov{m''}),\module{x_1}^{2s_1}\module{x_2}^{2s_2}\varphi(x)\rangle$ has only simple poles along lines $s_1=-1-k$, $s_2=-1-\ell$, $k,\ell\in\NN$ and no other poles,
\item
we have
$\psi_{x_2,-1}\psi_{x_1,-1}\cT=\psi_{x_1,-1}\psi_{x_2,-1}\cT$.\qed
\end{itemize}
\end{lemme}

We now will compute the specialization of $\cT$ along $x_1x_2=0$. Denote by $f$ the monomial $x_1x_2$ and let $i_f:X\to X\times\CC$ be the graph inclusion. Let $t$ be the coordinate on $\CC$, so that $i_f(X)=\{t-f=0\}$. Assume that $\cM'=\cM''=\cM$ and that $\cS=(\id,\id)$ is a Hermitian duality of weight~$0$ on $\cT$.

\begin{lemme}\label{lem:psi2dcn}
Under the same conditions, for $\cM=\cM'$ or $\cM''$, the $\cR_{\cX\times\CC}$-module $i_{f,+}\cM$ is strictly S-decomposable along $\{t=0\}$, we have $\psi_{f,\alpha}\cM=0$ for $\alpha\not\in\ZZ$ and there are local isomorphisms
\[
P\gr_\ell^{\rM}\psi_{f,-1}(\cM)\simeq
\begin{cases}
i_{Y_1,+}\psi_{x_1,-1}\cM\oplus i_{Y_2,+}\psi_{x_2,-1}\cM&\text{if }\ell=0,\\
i_{Z,+}\psi_{x_1,-1}\psi_{x_2,-1}\cM&\text{if }\ell=1.
\end{cases}
\]
Last, the sesquilinear pairing on the primitive part given by Formula \eqref{eq:psipolell} coincides with the corresponding specialization of $C$.
\end{lemme}

\begin{proof}
We will only insist on the computation of $\psi_{t,-1}C$, as the computation of $\psi_t(i_{f,+}\cM)$ is done in \cite{MSaito86}. We have $i_{f,+}\cM=\oplus_{k\in\NN}\cM\otimes\partiall_t^k\delta(t-x_1x_2)$ with the usual structure of a $\cR_{\cX\times\CC}$-module. By \loccit, we have $V_{-1}(i_{f,+}\cM)=\cR_\cX\cdot(\cM\otimes\delta)$ and, for $k\geq0$,
\[
V_{-1-k}(i_{f,+}\cM)=t^kV_{-1}(i_{f,+}\cM),\qquad
V_{-1+k}(i_{f,+}\cM)=\sum_{\ell\leq k}\partiall_t^\ell V_{-1}(i_{f,+}\cM).
\]
Moreover, $(t\partiall_t)^2$ vanishes on $\gr_{-1}^V(i_{f,+}\cM)$ and the monodromy filtration is given by
\begin{align*}
\rM_{-2}V_{-1}(i_{f,+}\cM)&=V_{-2}(i_{f,+}\cM),\\
\rM_{-1}V_{-1}(i_{f,+}\cM)&=t\partiall_t\cdot V_{-1}(i_{f,+}\cM) + V_{-2}(i_{f,+}\cM),\\
\rM_0V_{-1}(i_{f,+}\cM)&=\cR_\cX\cdot\big((x_1,x_2)\cM\otimes\delta\big),\\
\rM_1V_{-1}(i_{f,+}\cM)&=V_{-1}(i_{f,+}\cM).
\end{align*}
These formulas lead to the isomorphisms given in the lemma. It is also clear that $\can$ is onto, and one shows that $\var$ is injective, identifying therefore $\psi_{t,0}(i_{f,+}\cM)$ with $\gr_{-1}^{\rM}\psi_{t,-1}(i_{f,+}\cM)$. This gives the strict S-decomposability of $\cM$ along $\{f=0\}$.

Let us now compute Formula \eqref{eq:psipolell} for $\ell=0,1$. For that purpose, let $m',m''$ be local sections of $\cM$. They define local sections $m'\otimes\delta$ and $m''\otimes\delta$ of $V_{-1}(i_{f,+}\cM)$.

Assume first that $\ell=1$. For a local section $\varphi$ of $\cE_{X}^{n,n}$, we have to compute
\[
\thb^{-1}\res_{s=-1}\Big\langle i_{f,+}C\big((-it\partiall_t)m'\otimes\delta,\ov{m''\otimes\delta}\big), \mt^{2s}\chi(t)\varphi\wedge\itwopi dt\wedge d\ov t\Big\rangle.
\]
Then it is equal to (see also the computation of \eqref{eq:res} below)
\[
\res_{s=-1}(s+1)\big\langle C(m',\ov{m''}),\module{x_1x_2}^{2s}\chi(x_1x_2)\varphi\big\rangle,
\]
that we have to compare with $\res_{s_1=-1}\res_{s_2=-1}\!\big\langle C(m',\ov{m''}),\module{x_1}^{2s_1}\!\module{x_2}^{2s_2}\!\chi(x_1x_2)\varphi\big\rangle$. By Lemma \ref{lem:nonbichar}, both residues coincide.

Assume now that $\ell=0$ and take $m',m''\in(x_1,x_2)\cM$. If $m'=x_1m$ and $m''=x_2\mu$, then the function $\big\langle C(m,\ov{\mu}),\module{x_1x_2}^{2s+1}\chi(x_1x_2)\varphi\big\rangle$ has no pole at $s=-1$, after Lemma \ref{lem:nonbichar}. If for instance $m'=x_1m$ and $m''=x_1\mu$, then $\langle i_{Y_2,+}\psi_{x_2,-1}C(m,\ov\mu),\varphi\rangle$ and $\res_{s=-1}\big\langle C(m,\ov{\mu}),\module{x_1}^{2s+2}\module{x_2}^{2s}\chi(x_1x_2)\varphi\big\rangle$ coincide.
\end{proof}

\section{A local computation}\label{subsec:localcomput}
Let $\cT=(\cM',\cM'',C)$ be a \emph{smooth} twistor structure (\cf \T\ref{sec:smtwqc}) on a complex manifold $X$. The purpose of this paragraph is to compute the nearby cycles of $(\cM',\cM'',C)$ with respect to a function $f$ which takes the form $f(x_1,\dots,x_n)=x_1\cdots x_p$ for some local coordinates $x_1,\dots,x_n$ on $X$ and for some $p\geq1$. The goal is to show that, first, $\cT$ is strictly specializable along $f=0$, and to compute the primitive parts in terms of the restriction of $\cT$ to various coordinate planes, in the sense of Definition \ref{def:restrsmooth}. As, by definition, these restrictions are smooth twistor structures, this will imply that the primitive parts are so. The computation is local on $X$.

For $\ell\leq p$, denote by $\cJ_\ell$ the set of subsets $I\subset\{1,\dots,p\}$ having cardinal equal to $\ell$. For $I\in\cJ_\ell$, denote by $I^c$ its complement in $\{1,\dots,p\}$ and by $i_I$ the inclusion $\{x_{I}=0\}\hto X$.

\begin{proposition}\label{prop:localcomput}
Let $\cT$ be a smooth twistor structure of weight~$w$ on $X$. Then,
\begin{enumerate}
\item\label{prop:localcomput1}
the $\cR_\cX$-triple $\cT$ is regular and strictly specializable along $f=0$ (with a set $\Lambda$ of eigenvalues contained in $\ZZ$), and is strictly S-decomposable along $f=0$, with only one strict component;
\item\label{prop:localcomput2}
for $\alpha\in[-1,0[$, we have $\psi_{f,\alpha}\cT=0$ unless $\alpha=-1$ and, for any $\ell\geq0$, there is a functorial isomorphism
\[
\ooplus_{J\in\cJ_{\ell+1}}i_{J+}i_J^*\cT\isom P\gr^{\rM}_\ell\psi_{f,-1}(\cT)(\ell/2),
\]
where $\rM_\bbullet\psi_{f,\alpha}\cT$ denotes the monodromy filtration.
\end{enumerate}
\end{proposition}

\begin{proof}
Let us begin with the trivial smooth twistor $\cT=(\cO_\cX,\cO_\cX,C)$ where $C$ is trivially defined by $C(1,1)=1$. Except for the computation of $\psi_fC$, this is proved in \cite[\T3.6]{MSaito86}. We will recall some details in order to compute $\psi_fC$. We put $y_j=x_{p+j}$ for $j=1,\dots,n-p$. If $\delta$ denotes the $\cR_{\cX\times\CC}$-generator of $i_{f,+}\cO_\cX$, we have the following relations:
\[
t\delta=f(x)\delta,\quad x_i\partiall_{x_i}\delta=-(t\partiall_t+\hb)\delta,\quad \partiall_{y_j}\delta=0,\quad
t\Big(\prod_{i=1}^p\partiall_{x_i}\Big)\delta=(-t\partiall_t)^p\delta.
\]
This shows that $\delta\in V_{-1}(i_{f,+}\cO_\cX)$ and that $\gr_\alpha^V(i_{f,+}\cO_\cX)=0$ for $\alpha\not\in\ZZ$. Regularity along $\{t=0\}$ is also clear. It can be shown that \hbox{$V_{-1}(i_{f,+}\cO_\cX)=V_0(\cR_{\cX\times\CC})\cdot\delta$} and that any local section of $V_{-1}(i_{f,+}\cO_\cX)$ has a unique representative modulo $V_{-2}(i_{f,+}\cO_\cX)$ of the form
\[
\sum_{a\in\NN^p}\sum_{k=0}^{d_{a-\bf1}}g(x_{I^c(a)},y)x^{-a}(-t\partiall_t)^kP_{a-\bf1}(-t\partiall_t)(t\delta),
\]
with $P_{a-\bf1}(s)=\prod_{i=1}^p\prod_{1\leq\ell\leq a_i-1}(s-\ell \hb)$ (see \cite{Bibi96a}), $d_{a-\bf1}=\#\{i\mid a_i-1\geq0\}$, $I^c(a)=\{i\mid a_i=0\}$ and $g$ holomorphic in its variables. One may also show that sections of $P\gr_\ell^{\rM}\gr_{-1}^Vi_{f,+}\cO_\cX$ are uniquely represented by elements of the form
\[
\sum_{\substack{I\subset\{1,\dots,p\}\\ \module{I}=\ell+1}}\sum_{a_I\in(\NN^*)^I}g(x_{I^c},y)x_I^{-a_I}P_{a_I-\bf1_I}(-t\partiall_t)(t\delta)
\]
which can be rewritten as
\[
\sum_{\substack{I\subset\{1,\dots,p\}\\ \module{I}=\ell+1}}\sum_{b_I\in\NN^I}g(x_{I^c},y)\partial_{x_I}^{b_I}(x_I^{-\bf1_I}t\delta).
\]
It can be shown that, for $j\in I^c$,
\[
\partiall_{x_j}\big[g(x_{I^c},y)\partial_{x_I}^{b_I}(x_I^{-\bf1_I}t\delta)\big]=\partiall_{x_j}(g(x_{I^c},y))\partial_{x_I}^{b_I}(x_I^{-\bf1_I}t\delta)\mod \rM_{\ell-1}.
\]
One then gets an isomorphism
\begin{align*}
P\gr_\ell^{\rM}\gr_{-1}^Vi_{f,+}\cO_\cX&\isom\ooplus_{I\in\cJ_{\ell+1}}i_{I+}i_I^*i_{f,+}\cO_\cX\\
g(x_{I^c},y)\partial_{x_I}^{b_I}(x_I^{-\bf1_I}t\delta)&\mto g(x_{I^c},y)\partiall_{x_I}^{b_I}.
\end{align*}
One may compute similarly $P\gr_\ell^{\rM}\gr_0^Vi_{f,+}\cO_\cX$ (see \cite{MSaito86,Bibi87,Bibi96a}) and $\can,\var$, and prove in that way that $\cO_\cX$ is strictly specializable and strictly S-decomposable along $\{f=0\}$.

Let us now show that the previous isomorphism is an isomorphism of $\cR$-triples, once the left one is twisted by $(\ell/2)$, \ie the corresponding $C$ is multiplied by $\thb^{-\ell}$. Fix $J\in\cJ_{\ell+1}$ and, for any test form $\varphi$ with support contained in the fixed coordinate chart, put $\varphi=\varphi_J\wedge\prod_{j\in J}\itwopi dx_j\wedge d\ov x_j$. On the one hand, we have
\[
\big\langle i_{J+}i_J^*C(1,1),\varphi\big\rangle=\int_{\{x_j=0\mid j\in J\}}\varphi_J.
\]
On the other hand, according to Example \ref{ex:lefpol}, as $\rN''=i\rN=-it\partiall_t$, we have to compute
\begin{equation}\label{eq:res}
\thb^{-\ell}\res_{s=-1}\Big\langle i_{f,+}C\big((-it\partiall_t)^\ell x_J^{-\bf1_J}t\delta,\ov{x_J^{-\bf1_J}t\delta}\big),\mt^{2s}\chi(t)\varphi\wedge\itwopi dt\wedge d\ov t\Big\rangle.
\end{equation}
Remark first that
\begin{align*}
\eqref{eq:res}&=\thb^{-\ell}\res_{s=-1}\Big\langle i_{f,+}C\big(x_{J^c}\delta,\ov{x_{J^c}\delta}\big),(i\partiall_tt)^\ell(\mt^{2s}\chi(t))\varphi\wedge\itwopi dt\wedge d\ov t\Big\rangle\\
&=\res_{s=-1}\Big[(s+1)^\ell\Big\langle i_{f,+}C\big(x_{J^c}\delta,\ov{x_{J^c}\delta}\big),\mt^{2s}\chi(t)\varphi\wedge\itwopi dt\wedge d\ov t\Big\rangle\Big],
\end{align*}
as the term containing derivatives of $\chi$ will not create any residue. Putting $\tau=t-f$, we have
\[
\frac{\mt^{2s}\chi(t)\varphi\wedge\itwopi dt\wedge d\ov t} {\itwopi d\tau\wedge d\ov\tau}=\module{f}^{2s}\chi(f)\varphi,
\]
and, by definition of $i_{f,+}$,
\[
\eqref{eq:res}=\res_{s=-1}\Big[(s+1)^\ell\int_X\module{x_{J^c}}^2\module{f}^{2s}\chi(f)\varphi\Big].
\]
The integral has a pole of order $\ell+1$ at $s=-1$ produced by
\[
\int_X\module{x_J}^{2s}\varphi_J\wedge\prod_{j\in J}\itwopi dx_j\wedge d\ov x_j
\]
and the coefficient of the pole is $\int_{\{x_j=0\mid j\in J\}}\varphi_J$. This gives \eqref{prop:localcomput2} for $i_{f,+}\cO_\cX$.

\medskip
If now $\cT$ is any smooth twistor, remark that, for $\cM=\cM'$ or $\cM''$, we have $i_{f,+}\cM=\cM\otimes_{\cO_\cX}i_{f,+}\cO_\cX$ with its usual twisted structure of $\cR_\cX$-module, and that the action of $t$ and $\partiall_t$ comes from that on $i_{f,+}\cO_\cX$. As $\cM$ is assumed to be $\cO_\cX$-locally free, the filtration of $i_{f,+}\cM$ defined by $V_\alpha i_{f,+}\cM=\cM\otimes_{\cO_\cX}V_\alpha(i_{f,+}\cO_\cX)$ satisfies all properties of the Malgrange-Kashiwara filtration. It is then easy to deduce all assertions of the proposition for $\cM$ from the corresponding statement for $i_{f,+}\cO_\cX$.
\end{proof}

\chapter{Polarizable twistor $\cD$-modules}\label{chap:twistor}

\section{Definition of a twistor $\cD$-module} \label{subsec:defDtwist}
We will follow the inductive method of \cite[\T5.1]{MSaito86} to define the notion of (polarized) twistor $\cD$-module.

Let $X$ be a complex analytic manifold and let $w\in\ZZ$. We will define by induction on $d\in\NN$ the category $\MT_{\leq d}(X,w)$ of \emph{twistor $\cD_X$-modules of weight~$w$ on $X$, having a support of dimension $\leq d$}. This will be a subcategory of the category $\RTriples(X)$ introduced in \T\ref{subsec:Rtriples}. We will also define the full subcategory $\MTr_{\leq d}(X,w)$ of \emph{regular} twistor $\cD_X$-modules

\begin{definition}[Twistor $\cD$-modules]\label{def:twt}
The category $\index{$mt$@$\MT_{\leq d}(X,w)$, $\MTr_{\leq d}(X,w)$}\MT_{\leq d}(X,w)$ is the \emph{full subcategory} of $\RTriples(X)$ for which the objects are triples $(\cM',\cM'',C)$ satisfying:

\smallskip\noindent
(HSD)
$\cM',\cM''$ are \emph{holonomic}, \emph{strictly S-decomposable}, and have support of dimension $\leq d$.

\smallskip\noindent
$(\MT_{>0})$
For any open set $U\subset X$ and any holomorphic function $f:U\to\CC$, for any~$\alpha$ with $\reel(\alpha)\in[-1,0[$ and any integer $\ell\geq 0$, the triple
\[
\gr_\ell^{\rM}\Psi_{f,\alpha}(\cM',\cM'',C)\defin \lefpar \gr_{-\ell}^{\rM}\Psi_{f,\alpha}(\cM'),\gr_\ell^{\rM}\Psi_{f,\alpha}(\cM''),\gr_\ell^{\rM}\Psi_{f,\alpha}C\rigpar
\]
is an object of $\MT_{\leq d-1}(U,w+\ell)$.

\smallskip\noindent
$(\MT_0)$
For any zero-dimensional strict component $\{x_o\}$ of $\cM'$ or $\cM''$, we have
\[
(\cM'_{\{x_o\}},\cM''_{\{x_o\}},C_{\{x_o\}})=i_{\{x_o\}+}(\cH',\cH'',C_o)
\]
where $(\cH',\cH'',C_o)$ is a twistor structure of dimension $0$ and weight~$w$.
\end{definition}

Regular objects may be defined similarly:

\begin{definition}[Regular twistor $\cD$-modules]\label{def:regtwt}
The category $\MTr_{\leq d}(X,w)$ is the \emph{full subcategory} of $\RTriples(X)$ for which the objects are triples $(\cM',\cM'',C)$ satisfying:

\smallskip\noindent
(HSD)
$\cM',\cM''$ are \emph{holonomic}, \emph{strictly S-decomposable}, and have support of dimension $\leq d$.

\smallskip\noindent
(REG)
For any open set $U\subset X$ and any holomorphic function $f:U\to\CC$, the restrictions $\cM'_{\vert U},\cM''_{\vert U}$ are \emph{regular} along $\{f=0\}$.

\smallskip\noindent
$(\MT_{>0})$
For any open set $U\subset X$ and any holomorphic function $f:U\to\CC$, for any~$\alpha$ with $\reel(\alpha)\in[-1,0[$ and any integer $\ell\geq 0$, the triple
\[
\gr_\ell^{\rM}\Psi_{f,\alpha}(\cM',\cM'',C)\defin \lefpar \gr_{-\ell}^{\rM}\Psi_{f,\alpha}(\cM'),\gr_\ell^{\rM}\Psi_{f,\alpha}(\cM''),\gr_\ell^{\rM}\Psi_{f,\alpha}C\rigpar
\]
is an object of $\MTr_{\leq d-1}(U,w+\ell)$.

\smallskip\noindent
$(\MT_0)$
For any zero-dimensional strict component $\{x_o\}$ of $\cM'$ or $\cM''$, we have
\[
(\cM'_{\{x_o\}},\cM''_{\{x_o\}},C_{\{x_o\}})=i_{\{x_o\}+}(\cH',\cH'',C_o)
\]
where $(\cH',\cH'',C_o)$ is a twistor structure of dimension $0$ and weight~$w$.
\end{definition}

Clearly, $\MTr_{\leq d}(X,w)$ is a full subcategory of $\MT_{\leq d}(X,w)$. Let us justify all understatements made in the definition of the category $\MT(X,w)$. Notice that we have used the $\Psi$ functor of Definition~\ref{def:nearby}. Remark first:

\begin{proposition}\label{prop:twstrict}
If $(\cM',\cM'',C)$ is an object of $\MT_{\leq d}(X,w)$, then $\cM'$ and $\cM''$ are \emph{strict}, as well as $\gr_\ell^{\rM}\psi_{f,\alpha}\cM'$, $\gr_\ell^{\rM}\psi_{f,\alpha}\cM''$ for any analytic germ $f$, any $\alpha\in\CC$ and any $\ell\in\ZZ$. In particular, $\psi_{f,\alpha}\cM'$ and $\psi_{f,\alpha}\cM''$ are strict for any $\alpha\in\CC$.
\end{proposition}

\begin{proof}
Set $\cM=\cM'$ or $\cM''$. The strictness of $\cM$ follows from (HSD), after Corollary \ref{cor:strictdecstrict}. The strictness of $\gr_\ell^{\rM}\Psi_{f,\alpha}\cM$ for $\reel\alpha\in[-1,0[$ is by definition. To get the strictness of $\gr_\ell^{\rM}\psi_{f,\alpha}\cM$ for any $\alpha\not\in\NN$, remark that the property is local with respect to $\hb$. Use then the filtration $V_\bbullet^{(\hb_o)}$ and its graded pieces, and use the isomorphisms $t$ or $\partiall_t$ to increase or decrease $\reel\alpha$, depending whether $\ell_{\hb_o}(\alpha)<0$ or $\ell_{\hb_o}(\alpha)\geq-1$, if $\alpha\not\in \NN$.

Let us show the strictness of $\gr_\ell^{\rM}\psi_{f,\alpha}\cM$ for $\alpha=0$ (hence for any $\alpha\in\NN$). We can assume that $\cM$ has strict support. If $f\equiv0$ on the support of $\cM$, then the monodromy filtration is trivial and the strictness of $\psi_{f,0}\cM$ is a consequence of the strictness of $\cM$, by Kashiwara's equivalence \ref{cor:Kashiwaraequiv}. Otherwise, we know by \cite[Lemma 5.1.12]{MSaito86} that $\var:(\psi_{f,0}\cM,\rM_\bbullet)\to (\psi_{f,-1}\cM,\rM_{\bbullet-1})$ is injective and strict, \ie induces an injective morphism after grading. Therefore, each $\gr_\ell^{\rM}\psi_{f,0}\cM$ is strict.

The strictness of $\psi_{f,\alpha}\cM$ follows then from Lemma \ref{lem:Wstrict}.
\end{proof}

Notice also that we have locally finite strict S-decompositions $\cM'=\oplus_{Z'}\cM'_{Z'}$ and $\cM''=\oplus_{Z''}\cM''_{Z''}$ where $Z'$ belongs to the set of strict irreducible components of $\cM'$ and $Z''$ to that of $\cM''$. For any open set $U\subset X$, the irreducible components of all $Z'\cap U$ form the set of strict components of $\cM'_{\vert U}$, and similarly for $\cM''$. For $f:U\to\CC$, we have $\psi_{f,\alpha}\cM'_{Z'_U}=0$ for any $\alpha\not\in\NN$ if $f$ vanishes identically on the strict component $Z'_U$ of $\cM'_{\vert U}$, and has support of codimension one in $Z'_U$ otherwise. The support of $\psi_{f,\alpha}\cM'_{\vert U}$ has therefore dimension $\leq d-1$.

\smallskip
According to Proposition \ref{prop:Cdecomposable}, the component $C_{Z',Z''}$ of $C$ on $\cM'_{Z'|\bS}\otimes_{\cO_\bS}\ov{\cM''_{Z''|\bS}}$ vanishes unless $Z'=Z''$. We denote therefore by $C_Z$ the component of $C$ when $Z=Z'=Z''$ is a common strict component of $\cM'$ and $\cM''$. We thus have a S-decomposition
\begin{equation}\label{eq:MTdecompo}
(\cM',\cM'',C)=\oplus_Z(\cM'_Z,\cM''_Z,C_Z)
\end{equation}
indexed by the set of strict components of $\cM'$ or $\cM''$. We will see below (Corollary \ref{cor:MTdec}) that the set of strict components is the same for $\cM'$ and $\cM''$, and that each $(\cM'_Z,\cM''_Z,C_Z)$ is a twistor $\cD$-module of weight~$w$.

\smallskip
With such a notation, $(\MT_0)$ is concerned with the zero-dimensional strict components, which are not seen by $(\MT_{>0})$. Assume for instance that we work with right $\cR_\cX$-modules. Take local coordinates $x_1,\ldots,x_n$ at $x_o$. Then $(\MT_0)$ says that $\cM'_{\{x_o\}}=\cH'\otimes_\CC\CC[\partiall_{x_1},\ldots,\partiall_{x_n}]$, $\cM''_{\{x_o\}}=\cH''\otimes_\CC\CC[\partiall_{x_1},\ldots,\partiall_{x_n}]$ and $C_{\{x_o\}}$ is obtained by $\cR_{(X,\ov X),\bS}$-linearity from its restriction to $\cHS'\otimes_{\cO_\bS}\ov{\cHS''}$. There, it is equal to $C_o\cdot \delta_{x_o}$, where $\delta_{x_o}$ denotes the Dirac current at $x_o$ and $C_o:\cHS'\otimes_{\cO_\bS}\ov{\cHS''}\to\cO_\bS$ is given by \eqref{eq:twC0}.

\smallskip
It is easy to see now that the set of zero-dimensional strict components is the same for $\cM'$ and $\cM''$: if $\{x_o\}$ is not a strict component of $\cM''$ for instance, then $\cM''_{\{x_o\}}=0$ and thus $\cH''=0$. As $C_o$ is nondegenerate, this implies that $\cH'=0$, therefore $\cM'_{\{x_o\}}=0$ and $\{x_o\}$ is not a strict component of $\cM'$.

\smallskip
We will now give the basic properties of twistor $\cD$-modules.

\subsection{Locality}
For any open set $U\subset X$, there exists a natural restriction functor (and a regular analogue)
\[
\MT_{\leq d}(X,w)\To{\rho_U}\MT_{\leq d}(U,w).
\]
Moreover, if $(\cM',\cM'',C)$ is any object of $\RTriples(X)$ such that, for any open set~$U$ of a covering of $X$, $(\cM',\cM'',C)_{\vert U}$ is an object of $\MT_{\leq d}(U,w)$, then $(\cM',\cM'',C)$ is an object of the category $\MT_{\leq d}(X,w)$.

\Subsection{Stability by direct summand}
\begin{proposition}\label{prop:directsum}
If $(\cM',\cM'',C)=(\cM'_1\oplus \cM'_2,\cM''_1\oplus \cM''_2,C_1\oplus C_2)$ is an object of $\MT_{\leq d}(X,w)$, then each $(\cM'_i,\cM''_i,C_i)$ ($i=1,2$) also. Moreover, regularity is conserved.
\end{proposition}

\begin{proof}
The property of holonomicity restricts to direct summands, as well as the property of strict specializability (\ref{prop:canvar}\eqref{prop:canvara}) and, as Property \ref{prop:canvar}\eqref{prop:canvare} also restricts to direct summands, it follows that strict S-decomposability restricts too. It is easy to see that (REG) restricts to direct summands. Then argue by induction on $d$ for $(\MT_{>0})$. For $(\MT_0)$, use the fact that a direct summand of a trivial holomorphic vector bundle (on $\PP^1$) remains trivial.
\end{proof}

\begin{corollaire}\label{cor:MTdec}
If $(\cM',\cM'',C)$ is an object of the category $\MT_{\leq d}(X,w)$, then the strict components of $\cM'$ and $\cM''$ are the same and the S-decomposition \eqref{eq:MTdecompo} holds in $\MT_{\leq d}(X,w)$. Moreover, $\MT_{\leq d}(X,w)$ is the direct sum of the full subcategories $\MT_{\langle Z\rangle}(X,w)$ consisting of objects having strict support on the irreducible closed analytic subset $Z\subset X$ of dimension $\leq d$.
\end{corollaire}

\begin{proof}
Assume that there is a strict component $Z'$ of $\cM'$ which is not a strict component of $\cM''$. Then we have an object $(\cM'_{Z'},0,0)$ in $\MT_{\leq d}(X,w)$, according to the previous proposition. Argue now by induction on $\dim Z'$, the case $\dim Z'=0$ having being treated above. Let $f$ be the germ of any smooth function such that $f^{-1}(0)\cap Z'$ has codimension one in $Z'$. By induction, we have $\psi_{f,\alpha}\cM'_{Z'}=0$ for any $\alpha\not\in\NN$. By Kashiwara's equivalence on some open dense set of $Z'$, we may assume that $Z'=X$, that $\cM'_{Z'}$ is $\cO_\cX$-coherent, and we may choose for $f$ a local coordinate $t$, so that we conclude that $\cM'_{Z'}/t\cM'_{Z'}=0$. By Nakayama's lemma, we have $\cM'_{Z'}=0$ near $t=0$, hence $\cM'_{Z'}=0$ by definition of the strict support. This gives a contradiction.

The remaining statement is easy.
\end{proof}

\subsection{Kashiwara's equivalence}\label{subsubsec:KashequivMT}
Let $i$ denote the inclusion of $X$ as a closed analytic submanifold of the analytic manifold $X'$. Then the functor $i_+$ induces an equivalence between $\MT(X,w)$ and $\MT_X(X',w)$ (objects supported on $X$), which induces an equivalence between the regular subcategories.

\begin{remarque}\label{rem:ramif}
It follows from Remarks \ref{rem:strictspetr} and \ref{rem:psiramif} that, if $(\text{MT})_{>0}$ is satisfied for some holomorphic function $f$, it is satisfied for all $f^{r}$, $r\geq1$. Therefore, it is enough in practice to verify $(\text{MT})_{>0}$ for holomorphic functions which are not a power. A similar reduction holds for strict S-decomposability in (HSD), according to Lemma \ref{lem:strictdec}\eqref{lem:strictdec1}.
\end{remarque}

\Subsection{Generic structure of twistor $\cD_X$-modules}
\begin{proposition}\label{prop:genericMT}
Let $(\cM',\cM'',C)$ be an object of $\MT(X,w)$ having strict support on the irreducible closed analytic set $Z\subset X$. Then there exists an open dense set $Z'\subset Z$ and a smooth twistor structure $(\cH',\cH'',C)$ of weight~$w$ on $Z'$, such that $(\cM',\cM'',C)_{\vert Z'}=i_{Z'+}(\cH',\cH'',C)$.
\end{proposition}

\begin{proof}
Restrict first to a smooth open set of $Z$ and apply Kashiwara's equivalence to reduce to the case where $Z=X$. On some dense open set of $X$, the characteristic variety of $\cM'$ and $\cM''$ is contained in the zero section. By Proposition \ref{prop:CK}\eqref{prop:CK3}, $\cM'$ and $\cM''$ are $\cO_\cX$-locally free on some dense open set $X'$ of $X$, and by \ref{prop:CK}\eqref{prop:CK2}, putting $\cM=\cM'$ or $\cM''$, we have $\cM^\Cir=V_{-1}\cM^\Cir$ with respect to any local coordinate~$t$. Consequently, $\psi_{t,\alpha}\cM$ is supported on $\{\hb=0\}$ if $\alpha\not\in-\NN^*$, hence vanishes because of strictness (\cf Proposition \ref{prop:twstrict}). We thus have $\cM/t\cM=\gr_{-1}^V\cM$ and the monodromy filtration on it is trivial. Moreover, this is a $\cO_{\cX_0}$-locally free module and, after Remark \ref{rem:restrsmooth}, $\psi_{t,-1}C$ is the restriction of $C$ to $t=0$, so we may continue until we reach a twistor structure of dimension $0$. Therefore, $(\cM',\cM'',C)$ is a smooth twistor structure according to the definition given in \T\ref{sec:smtwqc}.
\end{proof}

\Subsection{Morphisms}
\begin{proposition}\label{prop:wsupw}
There is no nonzero morphism (in $\RTriples(X)$) from an object in the category $\MT(X,w)$ to an object in $\MT(X,w')$ if $w>w'$.
\end{proposition}

\begin{proof}
Let $\varphi:(\cM'_1,\cM''_1,C_1)\to (\cM'_2,\cM''_2,C_2)$ be such a morphism. According to Corollary \ref{cor:MTdec}, we may assume that both have the irreducible closed analytic set $Z$ as their strict support. As the result is clear for smooth twistor structures (there is no nontrivial morphism from $\cO_{\PP^1}(w)$ to $\cO_{\PP^1}(w')$ if $w>w'$), it follows from Proposition \ref{prop:genericMT} that the support of $\im\varphi$ is strictly smaller than $Z$. By definition of the strict support (\cf Definition \ref{def:strictsupp}), this implies that $\im\varphi=0$.
\end{proof}

\begin{proposition}\label{prop:abel}
The categories $\MT(X,w)$ and $\MTr(X,w)$ are abelian, all morphisms are strict and strictly specializable.
\end{proposition}

\begin{proof}
It is analogous to that of \cite[Prop\ptbl 5.1.14]{MSaito86}. Let us indicate it for $\MTr(X,w)$, the case of $\MT(X,w)$ being similar. Introduce the subcategory $\index{$mtw$@$\MTWr(X,w)$}\MTWr(X,w)$ of $\RTriples(X)$, the objects of which are triples with a finite filtration $W_\bbullet$ indexed by $\ZZ$ such that, for each $\ell$, $\gr_\ell^W$ is in $\MTr(X,w+\ell)$. The morphisms in $\MTWr(X,w)$ are the morphisms of $\RTriples(X)$ which respect the filtration $W$. Consider both properties:
\begin{enumerate}
\item[(a$_d$)] $\MTr_{\leq d}(X,w)$ abelian, all morphisms are strict and strictly specializable;
\item[(b$_d$)] $\MTWr_{\leq d}(X,w)$ abelian and morphisms are strict and strictly compatible with the filtration $W$.
\end{enumerate}

Remark first that $(\mathrm{a}_0)$ follows from Kashiwara's equivalence of \T\ref{subsubsec:KashequivMT} and the corresponding result in dimension $0$ (\cf\T\ref{subsubsec:smtw0}).

\medskip\noindent
$(\mathrm{a}_d)\implique(\mathrm{b}_d)$. Notice first that, by Proposition \ref{prop:twstrict} and Lemma \ref{lem:Wstrict}, the objects in $\MTWr_{\leq d}(X,w)$ are strict $\cR_\cX$-modules. According to Proposition \ref{prop:wsupw} and \cite[Lemme 5.1.15]{MSaito86}, $(\mathrm{a}_d)$ implies that the category $\MTWr_{\leq d}(X,w)$ is abelian and that morphisms are strictly compatible with $W$. Using Lemma \ref{lem:Wstrict} once more, we conclude that all morphisms are strict.

\medskip\noindent
$(\mathrm{b}_{d-1})\implique(\mathrm{a}_d)$ for $d\geq1$. The question is local. Let $\varphi=(\varphi',\varphi''):(\cM'_1,\cM''_1,C_1)\to (\cM'_2,\cM''_2,C_2)$ be a morphism of pure twistor $\cD$-modules of weight~$w$. According to Proposition \ref{prop:wsupw}, we may assume that all the $\cR_\cX$-modules involved have strict support $Z$ (closed irreducible analytic subset of $X$) of dimension $d$. We will first show that $\ker\varphi$ and $\coker\varphi$ are also strictly specializable, S-decomposable and have strict support $Z$.

Let $f$ be the germ of an analytic function not vanishing identically on $Z$, that we may assume to be a local coordinate $t$, using the graph embedding of $f$ and according to Kashiwara's equivalence of \T\ref{subsubsec:KashequivMT}. By $(\mathrm{b}_{d-1})$, $\psi_{t,\alpha}\varphi$ is strict if $\alpha\not\in\NN$, according to \ref{def:strictspe}\eqref{def:strictspeb} and \eqref{def:strictspec}. We will show below that
\begin{itemize}
\item
$\psi_{t,0}\varphi'$ and $\psi_{t,0}\varphi''$ are strict (hence so are $\psi_{t,k}\varphi'$ and $\psi_{t,k}\varphi''$ for any $k\in\NN$, according to \ref{def:strictspe}\eqref{def:strictspec}),
\item
$\can$ is onto for $\ker\varphi'$ and $\ker\varphi''$, and
\item
$\var$ is injective for $\coker\varphi'$ and $\coker\varphi''$.
\end{itemize}
The first assertion will be enough to show that $\varphi'$ and $\varphi''$ are strictly specializable, hence $\ker\varphi',\ldots,\coker\varphi''$ are also strictly specializable (Lemma \ref{lem:morstrictspe}). The two other assertions will insure that these modules satisfy Properties \ref{prop:canvar}\eqref{prop:canvarc} and \eqref{prop:canvard}, hence are strictly S-decomposable along $\{t=0\}$ and have neither sub nor quotient module supported on $Z\cap\{f=0\}$. Applying this for any such $f$ implies that $\ker\varphi',\dots,\coker\varphi''$ are strictly S-decomposable and have strict support $Z$. Now, $\ker\varphi',\ldots,\coker\varphi''$ are clearly holonomic and regular along $\{t=0\}$, hence they are also strict (\cf Corollary \ref{cor:strictdecstrict}). We now have obtained that $\varphi$ is strict and strictly specializable.

\medskip
Let us come back to the proof of the previous three assertions. As $\var$ is injective for $\cM'$ and $\cM''$, we identify $\psi_{t,0}\varphi'$ to the restriction of $\psi_{t,-1}\varphi'$ on $\im \rN\subset\psi_{t,-1}\cM'_2$, and similarly for $\varphi''$. By the inductive assumption, the morphism
\[
\cN:(\psi_{t,-1}\cM'_k,\psi_{t,-1}\cM''_k,\psi_{t,-1}C_k)\to (\psi_{t,-1}\cM'_k,\psi_{t,-1}\cM''_k,\psi_{t,-1}C_k)(-1)
\]
is strict, for $k=1,2$ and $\im \cN$ is an object of $\MTWr_{\leq d-1}(X_0,w)$. Using once more this inductive assumption, the restriction of $\psi_{t,-1}\varphi$ on $\im \cN$ is strict, hence the first point.

In order to show the other assertions, consider the following diagram of exact sequences (and the similar diagram for $\varphi'$):
\[
\xymatrix@C=.79cm{
0\to\psi_{t,-1}\ker\varphi''\ar[r]\ar[d]^-{\can}&
\psi_{t,-1}\cM''_1\ar[rr]^-{\psi_{t,-1}\varphi''}
\ar@{>>}[d]^-{\can}\ar@/_1.5pc/[dd]_(.27){\rN_1} |\hole&& \psi_{t,-1}\cM''_2\ar[r]\ar@{>>}[d]^-{\can}\ar@/_1.5pc/[dd]_(.27){\rN_2} |\hole&
\psi_{t,-1}\coker\varphi''\ar@{>>}[d]^-{\can}\to0\\
0\to\psi_{t,0}\ker\varphi''\ar[r]\ar@{^{ (}->}[d]^-{\var}&\psi_{t,0}\cM''_1\ar[rr]^-{\psi_{t,0}\varphi''}\ar@{^{ (}->}[d]^-{\var} && \psi_{t,0}\cM''_2\ar[r]\ar@{^{ (}->}[d]^-{\var}& \psi_{t,0}\coker\varphi''\ar[d]^-{\var}\to0
\\0\to\psi_{t,-1}\ker\varphi''\ar[r]&\psi_{t,-1}\cM''_1\ar[rr]^-{\psi_{t,-1}\varphi''}&& \psi_{t,-1}\cM''_2\ar[r]&\psi_{t,-1}\coker\varphi''\to0
}
\]
We have to prove that the left up $\can$ is onto and that the right down $\var$ is injective. This amounts to showing that $\im \rN_1\cap\ker\psi_{t,-1}\varphi''=\rN_1(\ker\psi_{t,-1}\varphi'')$ (because this is equivalent to $\im \can\cap\ker\psi_{t,0}\varphi''=\can(\ker\psi_{t,0}\varphi'')$) and $\im \rN_2\cap\im\psi_{t,-1}\varphi''=\rN_2(\im\psi_{t,-1}\varphi'')$. This follows from the lemma below applied to the germs of the various sheaves.

\begin{lemme}\label{lem:morozov}
Let $E_1,E_2$ be two $\ZZ$-modules, equipped with nilpotent endomorphisms $\rN_1,\rN_2$. Let $\lambda:E_1\to E_2$ be a morphism commuting with $\rN_1,\rN_2$, which is strictly compatible with the corresponding monodromy filtrations $\rM(\rN_1),\rM(\rN_2)$. Then
$$
\im \rN_1\cap\ker\lambda=\rN_1(\ker\lambda) \qqbox{and} \im \rN_2\cap \im\lambda=\rN_2(\im\lambda).
$$
\end{lemme}

\begin{proof}
By the characteristic property of the monodromy filtration and by the strict compatibility of $\lambda$, we have $\rM(\rN_1)\cap\ker\lambda=\rM(\rN_{1|\ker\lambda})$ and $\rM(\rN_2)\cap\im\lambda=\rM(\rN_{2|\im\lambda})$. Moreover, each $\gr_\ell^{\rM}\lambda$ decomposes with respect to the Lefschetz decomposition. It follows that the property of the lemma is true at the graded level.

Let us show the first equality, the second one being similar. By the previous remark we have
\[
\im\rN_1\cap\ker\lambda\cap \rM(\rN_1)_\ell=\rN_1\big(\ker\lambda\cap \rM(\rN_1)_{\ell+2}\big)+\im\rN_1\cap\ker\lambda\cap \rM(\rN_1)_{\ell-1},
\]
and we may argue by induction on $\ell$ to conclude.
\end{proof}

To end the proof of $(\mathrm{b}_{d-1})\implique(\mathrm{a}_d)$, it remains to be proved that $\ker\varphi$ and $\coker\varphi$ satisfy $(\MT_{>0})$. It follows from the abelianity of $\MTWr_{\leq d-1}(X,w)$ and from the strict specializability of $\varphi$ that $\Psi_{t,\alpha}\ker\varphi$ and $\Psi_{t,\alpha}\coker\varphi$ (with $\reel\alpha\in[-1,0[$) are in $\MTWr_{\leq d-1}(X,w)$ and, as we have seen in Lemma \ref{lem:morozov}, the weight filtration is the monodromy filtration. This gives $(\MT_{>0})$, concluding the proof of Proposition \ref{prop:abel}.
\end{proof}

\begin{corollaire}\label{cor:Lefschetzdecmorphism}
Given any morphism $\varphi:\cT_1\to\cT_2$ between objects of $\MT(X,w)$ and any germ $f$ of holomorphic function on $X$, then, for any $\alpha\not\in\NN$, the specialized morphism $\Psi_{f,\alpha}\varphi$ is strictly compatible with the monodromy filtration $\rM_\bbullet$ and, for each $\ell\in\ZZ$, $\gr_\ell^{\rM}\Psi_{f,\alpha}\varphi$ decomposes with respect to the Lefschetz decomposition, \ie
$$
\gr_\ell^{\rM}\Psi_{f,\alpha}\varphi=
\begin{cases}
\ooplus_{k\geq0}\cN^kP\gr_{\ell+2k}^{\rM}\Psi_{f,\alpha}\varphi&(\ell\geq0),\\
\ooplus_{k\geq0}\cN^{k-\ell}P\gr_{-\ell+2k}^{\rM}\Psi_{f,\alpha}\varphi&(\ell\leq0).
\end{cases}
$$
In particular we have
\[
\gr_\ell^{\rM}\Psi_{t,\alpha}\ker\varphi=\ker\gr_\ell^{\rM}\Psi_{t,\alpha}\varphi
\]
and similarly for $\coker$, where, on the left side, the filtration $\rM_\bbullet$ is that induced naturally by $\rM_\bbullet\Psi_{t,\alpha}\cT_1$ or, equivalently, the monodromy filtration of $\cN$ acting on $\Psi_{t,\alpha}\ker\varphi=\ker\Psi_{t,\alpha}\varphi$.\qed
\end{corollaire}

\begin{corollaire}\hspace*{-2mm}
If $\cT$ is in $\MT_{\leq d}(X,w)$, then the Lefschetz decomposition for $\gr_\ell^{\rM}\Psi_{t,\alpha}\cT$ (with $\reel(\alpha)\in[-1,0[$) holds in $\MT_{\leq d-1}(X,w+\ell)$.
\end{corollaire}

\begin{proof}
Indeed, $\cN:\gr_\ell^{\rM}\Psi_{t,\alpha}\cT\to\gr_{\ell-2}^{\rM}\Psi_{t,\alpha}\cT(-1)$ is a morphism in the category $\MT_{\leq d-1}(X,w+\ell)$, which is abelian, so the primitive part is an object of this category, and therefore each term of the Lefschetz decomposition is also an object of this category.
\end{proof}

\subsection{Graded Lefschetz twistor $\cD_X$-modules}\label{subsec:grleftw}
Given $\varepsilon=\pm1$, we may define the category $\index{$mtl$@$\MLT_{\leq d}(X,w;\varepsilon)$, $\MLTr_{\leq d}(X,w;\varepsilon)$}\MLT_{\leq d}(X,w;\varepsilon)$ of graded Lefschetz twistor $\cD_X$-modules as in \T\ref{subsec:grlefschetz}: the objects are pairs $(\cT,\cL)$, with $\cT=\oplus\cT_j$, and $\cT_j$ are objects of \hbox{$\MT_{\leq d}(X,w-\varepsilon j)$}; $\cL$ is a graded morphism $\cT_j\to\cT_{j-2}(\varepsilon)$ of degree $-2$, such that, for $j\geq0$, $\cL^j:\cT_j\to\cT_{-j}(\varepsilon j)$ is an isomorphism. Notice that, by Proposition \ref{prop:abel}, $P\cT_j$ is an object of $\MT_{\leq d}(X,w-\varepsilon j)$ and the Lefschetz decomposition of $\cT_j$ holds in $\MT_{\leq d}(X,w-\varepsilon j)$; moreover, the category $\MLT_{\leq d}(X,w;\varepsilon)$ is abelian, any morphism is graded with respect to the Lefschetz decomposition, and moreover is strict and strictly specializable, as follows from Proposition \ref{prop:abel}.

More generally, for any $k\geq0$ and $\varepsilong=(\varepsilon_1,\dots,\varepsilon_k)=(\pm1,\dots,\pm1)$, we may define the category $\MLT_{\leq d}(X,w;\varepsilong)$ of $k$-graded Lefschetz twistor $\cD_X$-modules: the objects are tuples $(\cT,\bcL)$, with $\bcL=(\cL_1,\dots,\cL_k)$, $\cT=\oplus_{\bmj\in\ZZ^k}\cT_{\bmj}$, each $\cT_{\bmj}$ is an object in $\MT_{\leq d}\big(X,w-\sum_i\varepsilon_i j_i\big)$, the morphisms $\cL_i$ should pairwise commute, be of $k$-degree $(0,\dots,-2,\dots,0)$ and for any $\bmj$ with $j_i\geq0$, $\cL_i^{j_i}$ should induce an isomorphism from $\cT_{\bmj}$ to the component where $j_i$ is replaced with $-j_i$; the primitive part $P\cT_{\bmj}$, for $j_1,\dots,j_k\geq0$, is the intersection of the $\ker\cL_i^{j_i+1}$ and we have a Lefschetz multi-decomposition, with respect to which any morphism is multi-graded. The category is abelian, and any morphism is strict and strictly specializable.

\begin{lemme}\label{lem:Ppsi}
Let $(\cT,\bcL)$ be an object of the category $\MLT_{\leq d}(X,w;\varepsilong)$. Then, for any $\alpha$ with $\reel(\alpha)\in[-1,0[$, the specialized object $\big(\gr_\bbullet^{\rM}\Psi_{t,\alpha}\cT,(\gr_\bbullet^{\rM}\Psi_{t,\alpha}\bcL,\cN)\big)$ is an object of $\MLT_{\leq d-1}\big(X,w;(\varepsilong,-1)\big)$ and $P_{\bcL}\gr_\ell^{\rM}\Psi_{t,\alpha}\cT_{\bmj}=\gr_\ell^{\rM}\Psi_{t,\alpha}P_{\bcL}\cT_{\bmj}$, where $P_{\bcL}$ denotes the multi-primitive part with respect to $\bcL$.
\end{lemme}

\begin{proof}
The lemma is a direct consequence of the strict compatibility of the $\Psi_{t,\alpha}\cL_i$ with the monodromy filtration $\rM(\cN)$, as follows from Proposition \ref{prop:abel}.
\end{proof}

\begin{lemme}
The category $\MLT_{\leq d}(X,w;\varepsilong)$ has an inductive definition \hbox{analogous} to that of $\MT_{\leq d}(X,w)$, where one replaces the condition $(\MT_{>0})$ with the condition $(\MLT_{>0})$, asking that $\big(\gr_\bbullet^{\rM}\Psi_{t,\alpha}(\cT,\bcL),\rN\big)$ is an object of $\MLT_{\leq d-1}\big(X,w;(\varepsilong,-1)\big)$, and the condition $(\MT_0)$ with the analogous property $(\MLT_0)$.
\end{lemme}

\begin{proof}
According to the previous lemma, it is enough to show that, if $(\cT,\bcL)$ satisfies the inductive conditions, then it is an object of $\MLT_{\leq d}(X,w;\varepsilong)$. This is done by induction on $d$, the case $d=0$ being easy. One shows first that each $\cT_{\bmj}$ is in $\MT_{\leq d}(X,w-\sum\varepsilon_ij_i)$ for any $\bmj$ and that $\Psi_{t,\alpha}\cL_i^{j_i}$ is an isomorphism from $\Psi_{t,\alpha}\cT_{\bmj}$ to $\Psi_{t,\alpha}\cT_{j_1,\dots,-j_i,\dots,j_k}$ for any $i=1,\dots,k$, any $\bmj$ with $j_i\geq0$, any local coordinate $t$ and any $\alpha\in\CC$. Considering the decomposition with respect to the support, one deduces that $\cL_i^{j_i}$ is an isomorphism from $\cT_{\bmj}$ to $\cT_{j_1,\dots,-j_i,\dots,j_k}$.
\end{proof}

\begin{remarque}[Regularity]
Similar results hold for the category $\MLTr_{\leq d}(X,w;\varepsilon)$ of graded Lefschetz regular twistor $\cD_X$-modules and its multi-graded analogues.
\end{remarque}

\subsection{Vanishing cycles}

Let $\cT$ be an object of $\MT_{\leq d}(X,w)$. By definition, for any locally defined analytic function $f$, the object $(\psi_{f,-1}\cT,\rM_\bbullet(\cN))$ is an object of $\MTW_{\leq d}(X,w)$.

\begin{corollaire}[Vanishing cycles, \cf {\cite[Lemme 5.1.12]{MSaito86}}] \label{cor:vanishingcycles}
For such a $\cT$, the object $(\phi_{f,0}\cT,\rM_\bbullet(\cN))$ is in $\MTW_{\leq d}(X,w)$ and $(\gr_\bbullet^{\rM}\phi_{f,0}\cT,\gr^{\rM}_{-2}\cN)$ is an object of $\MLT_{\leq d}(X,w;-1)$. Moreover, the morphisms $\cCan,\cVar$ are filtered morphisms
\begin{align*}
(\psi_{f,-1}\cT,\rM_{\bbullet}(\cN))&\To{\cCan}(\phi_{f,0}\cT(-1/2), \rM_{\bbullet-1}(\cN))\\
(\phi_{f,0}\cT(-1/2), \rM_{\bbullet-1}(\cN))&\To{\cVar} (\psi_{f,-1}\cT(-1), \rM_{\bbullet-2}(\cN)),
\end{align*}
hence are morphisms in $\MTW(X,w)$, and similarly for $\gr^{\rM}_{-1}\cCan$ and $\gr^{\rM}_{-1}\cVar$.
\end{corollaire}

\begin{proof}
We can assume that $\cT$ has strict support on an irreducible closed analytic subset $Z$ of $X$. If $f\equiv0$ on $Z$, then the result follows from Kashiwara's equivalence and Lemma \ref{lem:phiCC0}.

Assume now that $f\not\equiv0$ on $Z$. The object $\phi_{f,0}\cT(-1/2)$ is equipped with a filtration $W_\bbullet\phi_{f,0}\cT$ naturally induced by $\rM_\bbullet(\cN)\psi_{f,-1}\cT$. As such, according to Lemma \ref{lem:cCancVar}, it is identified with the image of $\cN:\big(\psi_{f,-1}\cT,\rM_\bbullet(\cN)\big)\to\big(\psi_{f,-1}\cT(-1),\rM_{\bbullet-2}(\cN)\big)$, hence is an object of $\MTW(X,w)$, because this category is abelian.

The result now follows from \cite[Lemme 5.1.12]{MSaito86}, which gives in particular that $W_\bbullet\phi_{f,0}\cT=\rM_{\bbullet-1}(\cN)\phi_{f,0}\cT$.
\end{proof}

\begin{remarque}[Regularity]
Starting with an object $\cT$ of $\MTr(X,w)$, we conclude that $(\gr_\bbullet^{\rM}\phi_{f,0}\cT,\cN)$ is an object of $\MLTr(X,w)$.
\end{remarque}

\Subsection{Behaviour with respect to the functors $\Xi_{\DR}$ and $\Xi_{\Dol}$}

We can now give a statement more precise than Proposition \ref{prop:restrhbo} concerning the restriction to $\hb=\hb_o$, and in particular the behaviour of the monodromy filtration and the property of S-decomposability.

\begin{proposition}[Restriction to $\hb=\hb_o$]\label{prop:restricanvar}
Let $(\cM',\cM'',C)$ be an object of $\MT(X,w)$. Put $\cM=\cM'$ or $\cM''$. Fix $\hb_o\in\Omega_0$ and put $M_{\hb_o}=\cM/(\hb-\hb_o)\cM$.
\begin{enumerate}
\item\label{prop:restricanvar3}
If $(\cM',\cM'',C)$ is in $\MTr(X,w)$ and $\hb_o\neq0$, then $M_{\hb_o}$ is a regular holonomic $\cD_X$-module.
\item\label{prop:restricanvar2}
Let $f:U\to\CC$ be a holomorphic function on some open set $U$.
\begin{enumerate}
\item\label{prop:restricanvar2a}
For any $\alpha\not\in\NN$, the restriction to $\hb=\hb_o$ of the monodromy filtration $\rM_\bbullet \psi_{f,\alpha}\cM$ of $\rN$ is the monodromy filtration of its restriction $\rN$ on $\psi_{f,\alpha}^{(\hb_o)} M_{\hb_o}$.
\item\label{prop:restricanvar2b}
Assume that $f$ is a projection $t$, that $\can$ is onto and $\var$ is injective. Then, $\can:\psi_{t,-1}^{(\hb_o)}M_{\hb_o}\to\psi_{t,0}^{(\hb_o)}M_{\hb_o}$ is onto and $\var:\psi_{t,0}^{(\hb_o)}M_{\hb_o}\to\psi_{t,-1}^{(\hb_o)}M_{\hb_o}$ is injective.
\end{enumerate}
\item\label{prop:restricanvar1}
Assume that $\hb_o\not\in\Sing\Lambda$. Then $M_{\hb_o}$ is a strictly S-decomposable holonomic $\cD_X$-module. If $\cM$ has strict support $Z$ (irreducible closed analytic subset of $X$), then so has the restriction to $\hb=\hb_o$.
\end{enumerate}
\end{proposition}

\begin{Proof}
\eqref{prop:restricanvar3}
Using the definition of regularity as in \cite[Def\ptbl (3.1.12)]{Mebkhout89}, one shows by induction on $\dim\supp M_{\hb_o}$ that $M_{\hb_o}$ is regular.

\smallskip\eqref{prop:restricanvar2a}
This is a consequence of the strictness of $\gr_k^{\rM}\psi_{t,\alpha}\cM$ proved in Proposition \ref{prop:twstrict}. Indeed, by strictness, the filtration $\rM_\bbullet\psi_{f,\alpha}^{(\hb_o)}M_{\hb_o}$ naturally induced by $\rM_\bbullet\psi_{f,\alpha}\cM$ satisfies $\gr_k^{\rM}\psi_{f,\alpha}^{(\hb_o)}M_{\hb_o}=\gr_k^{\rM}\psi_{f,\alpha}\cM/(\hb-\hb_o)\gr_k^{\rM}\psi_{f,\alpha}\cM$, and then satisfies the characteristic properties of the monodromy filtration of the restriction of $\rN$.

\smallskip\eqref{prop:restricanvar2b}
That $\can$ remains onto is clear. In order to show that $\var$ remains injective by restriction, we will use that $\psi_{t,-1}\cM\big/\im \rN$ is strict: indeed, $\rN:\psi_{t,-1}\cM\to\psi_{t,-1}\cM$ is part of a morphism in $\MTW(X,w)$, hence its cokernel is strict. This implies that
\[
\im \rN\big/(\hb-\hb_o)\im \rN\to\psi_{t,-1}\cM\big/(\hb-\hb_o)\psi_{t,-1}\cM
\]
is injective and therefore
\[
\im \rN\cap(\hb-\hb_o)\psi_{t,-1}\cM= (\hb-\hb_o)\im \rN.
\]

Let $m$ be a local section of $\psi_{t,0}\cM$ such that $tm\in(\hb-\hb_o)\psi_{t,-1}\cM$. As $\can$ is onto, there exists a local section $m'$ of $\psi_{t,-1}\cM$ such that $m=-\partiall_tm'$. Then $\rN m'\in(\hb-\hb_o)\psi_{t,-1}\cM$. By the strictness property above, we have $\rN m'=\rN(\hb-\hb_o)m''$ for some local section $m''$ of $\psi_{t,-1}\cM$, and hence $t[m-(\hb-\hb_o)(-\partiall_tm'')]=0$. As $\var$ is injective, we have $m\in(\hb-\hb_o)\psi_{t,0}\cM$, as was to be proved.

\smallskip
\eqref{prop:restricanvar1}
Assume that $\hb_o\not\in\Sing\Lambda$ and that $\cM$ has strict support $Z$. We will show that $M_{\hb_o}$ is strictly S-decomposable and has strict support $Z$ (the definition of these notions for $\cD_X$-modules are given in \cite{MSaito86}; they are also obtained by doing $\hb=1$ in the corresponding definitions for $\cR_\cX$-modules). Let $f:(X,x_o)\to(\CC,0)$ be an analytic germ which is nonconstant on $(Z,0)$. Using the graph embedding of $f$, we may assume that $f$ is a coordinate $t$. By Proposition \ref{prop:restrhbo}, we have $\psi_{t,-1}^{(\hb_o)}M_{\hb_o}=\psi_{t,-1}M_{\hb_o}$ and $\psi_{t,0}^{(\hb_o)}M_{\hb_o}=\psi_{t,0}M_{\hb_o}$, $\can$ restricts to $\hb_o\can_{\hb_o}$ and $\var$ restricts to $\var_{\hb_o}$. Therefore, the conclusion follows from \eqref{prop:restricanvar2b}.
\end{Proof}

\section{Polarization}

\begin{definition}[Polarization]\label{def:polarization}
A \emph{polarization} of an object $\cT$ of $\MT_{\leq d}(X,w)$ is a ses\-qui\-li\-near Hermitian duality $\cS:\cT\to\cT^*(-w)$ of weight~$w$ (\cf Definition \ref{def:hermw}) such that:

\smallskip\noindent
$(\MTP_{>0})$
for any open set $U\subset X$ and any holomorphic function $f:U\to\CC$, for any $\alpha$ with $\reel(\alpha)\in[-1,0[$ and any integer $\ell\geq 0$, the morphism $(P\gr_\bbullet^{\rM}\Psi_{f,\alpha}\cS)_\ell$ induces a polarization of $P_\ell\Psi_{f,\alpha}\cT$,

\smallskip\noindent
$(\MTP_0)$
for any zero-dimensional strict component $\{x_o\}$ of $\cM'$ or $\cM''$, we have $\cS=i_{\{x_o\}+}\cS_o$, where $\cS_o$ is a polarization of the zero-dimensional twistor structure $(\cH',\cH'',C_o)$.
\end{definition}

\begin{Remarques}\label{rem:pol}
\begin{enumerate}
\item\label{rem:pol1}
Notice that Condition $(\MTP_{>0})$ is meaningful because of Remark \ref{rem:psiadj}.

\item\label{rem:pol2}
Conditions $(\MTP_{>0})$ and $(\MTP_0)$ imply that the components $S'$ and $S''=(-1)^wS'$ of $\cS$ are isomorphisms $\cM''\isom\cM'$: indeed, one may assume that $\cT$ has only one strict component; by induction on the dimension, using a local coordinate, one obtains that $S'$ is an isomorphism on a dense open set of the support; by definition of the strict support, $S'$ is thus an isomorphism.
\end{enumerate}
\end{Remarques}

We will denote by $\index{$mtp$@$\MT_{\leq d}(X,w)^\rp$, $\MTr_{\leq d}(X,w)^\rp$}\MT_{\leq d}(X,w)^\rp$ the full subcategory of $\MT_{\leq d}(X,w)$ of polarizable objects, and similarly for $\MTr_{\leq d}(X,w)^\rp$. According to Proposition \ref{cor:morphismestrictdecompo}, we have a S-decomposition
\begin{equation}\label{eq:MTPdecompo}
(\cM',\cM'',C,\cS)=\oplus_Z(\cM'_Z,\cM''_Z,C_Z,\cS_Z).
\end{equation}
The following proposition is easy:
\pagebreak[2]
\begin{Proposition}
\begin{enumerate}
\item
In the situation of Proposition \ref{prop:directsum}, if a polarization $\cS$ is the direct sum of two morphisms $\cS_1$ and $\cS_2$, then each $\cS_i$ is a polarization of $(\cM'_i,\cM''_i,C_i)$.
\item
Corollary \ref{cor:MTdec} holds for $\MT_{\leq d}(X,w)^\rp$ or $\MTr_{\leq d}(X,w)^\rp$.
\item
Kashiwara's equivalence of \T\ref{subsubsec:KashequivMT} holds for $\MT(X,w)^\rp$ or $\MTr(X,w)^\rp$.\qed
\end{enumerate}
\end{Proposition}

\Subsection{Semi-simplicity}
\begin{proposition}\label{prop:semisimplicity}
If $\cT_1$ is a subobject (in the category $\MT(X,w)$) of a polarized object $(\cT,\cS)$, then $\cS$ induces a polarization $\cS_1$ of $\cT_1$ and $(\cT_1,\cS_1)$ is a direct summand of $(\cT,\cS)$ in $\MT(X,w)^\rp$. In particular, the category $\MT(X,w)^\rp$ is semisimple (all objects are semisimple and morphisms between simple objects are zero or isomorphisms).
\end{proposition}

\begin{proof}
By induction on the dimension of the support, the result being clear if the support has dimension $0$ (see Fact \ref{fact:pol}). We may also assume that $\cT$ has strict support a closed irreducible analytic subset $Z$ of $X$. Put $\cS=(S',S'')$ and $S=S'$ or $S''$. Consider then the exact sequences
\[
\xymatrix{
0&\cT^*_1(-w)\ar[l]&\cT^*(-w)\ar[l]&\cT^*_2(-w)\ar[l]&0\ar[l]\\
0\ar[r]&\cT_1\ar[r]\ar[u]^-{\cS_1}&\cT\ar[r]\ar[u]^-{\cS}_-\wr&\cT_2\ar[r]&0
}
\]
where $\cT_2$ is the cokernel, in the abelian category $\MT(X,w)$, of $\cT_1\hto\cT$. We want to show first that $\cS_1$ is an isomorphism. This is a local statement. Take a local coordinate $t$ such that $Z\not\subset\{t=0\}$ and apply $\Psi_{t,\alpha}$ to the previous diagram ($\reel(\alpha)\in[-1,0[$). According to Corollary \ref{cor:Lefschetzdecmorphism}, the exact sequences in the following diagram remain exact, if $\rM_\bbullet$ denotes the monodromy filtration of $\rN$:
\[
\xymatrix{
0&\gr_\ell^{\rM}\Psi_{t,\alpha}\cT^*_1(-w)\ar[l]&\gr_\ell^{\rM}\Psi_{t,\alpha}\cT^*(-w)\ar[l]& \gr_\ell^{\rM}\Psi_{t,\alpha}\cT^*_2(-w)\ar[l]&0\ar[l]\\
0\ar[r]&\gr_\ell^{\rM}\Psi_{t,\alpha}\cT_1\ar[r]\ar[u]^-{\gr_\ell^{\rM}\Psi_{t,\alpha}\cS_1}&\gr_\ell^{\rM}\Psi_{t,\alpha}\cT\ar[r]\ar[u]^-{\gr_\ell^{\rM}\Psi_{t,\alpha}\cS}_-\wr&\gr_\ell^{\rM}\Psi_{t,\alpha}\cT_2\ar[r]&0
}
\]
Using the inductive assumption, we conclude that each $\gr_\ell^{\rM}\Psi_{t,\alpha}\cS_1$ is an isomorphism, hence $\Psi_{t,\alpha}\cS_1$ too. Arguing now as in Remark \ref{rem:pol}\eqref{rem:pol2}, we conclude that $\cS_1$ is an isomorphism.

We now have a decomposition $\cM'=\cM'_2\oplus S(\cM''_1)$ and $\cM''=\cM''_1\oplus S^{-1}(\cM'_2)$ and we have by definition a decomposition $S=S_1\oplus S_2$, where $S_2$ is the isomorphism such that $S_2^{-1}$ is the restriction of $S^{-1}$ to $\cM'_2$.

It remains to be proved that we have a decomposition $C=C_1\oplus C_2$. By definition, we have $C(m'_2,\ov{n''_1})=0$ for local sections $m'_2,n''_1$ of $\cM'_{2|\bS}$ and $\cM''_{1|\bS}$ respectively. It is enough to show that $C(Sm''_1,\ov{S^{-1}n'_2})=0$ for local sections $m''_1,n'_2$ of $\cM''_{1|\bS}$ and $\cM'_{2|\bS}$ respectively. This is a direct consequence of the fact that $\cS$ is Hermitian.
\end{proof}

\begin{remarque}[Regularity]
The same result holds with regular objects.
\end{remarque}

\subsection{Polarized graded Lefschetz twistor $\cD$-modules}
Let $(\cT,\bcL)$ be an object of $\MLT(X,w;\varepsilong)$. A \emph{polarization} $\cS$ is a graded isomorphism $\cS:\cT\to\cT^*(-w)$ which is Hermitian, \ie satisfying $\cS_{\bmj}^*=(-1)^j\cS_{-\bmj}$, such that each $\cL_i$ is skew-adjoint with respect to $\cS$ (\ie $\cL_i^*\circ\cS_{\bmj}=-\cS_{\bmj-2{\bf1}_i}\circ\cL_i$ for any $i=1,\dots,k$ and any $\bmj$) and that, for each $\bmj$ with nonnegative components, the induced morphism
\[
\cS_{-\bmj}\circ\cL_1^{j_1}\cdots\cL_k^{j_k}:P_{\bcL}\cT_{\bmj}\to (P_{\bcL}\cT_{\bmj})^*(-w+\textstyle\sum\varepsilon_ij_i)
\]
is a polarization of the object $P_{\bcL}\cT_{\bmj}$ of $\MT(X,w-\sum\varepsilon_ij_i)$.

\begin{lemme}\label{lem:semisimplicitelef}
The categories $\MLT(X,w;\varepsilong)^\rp$ and $\MLTr(X,w;\varepsilong)^\rp$ have an inductive definition as in Definition \ref{def:polarization}.
\end{lemme}

\begin{proof}
This directly follows from the commutativity of $P_{\bcL}$ and $\gr_\ell^{\rM}\Psi_{t,\alpha}$ shown in Lemma \ref{lem:Ppsi}.
\end{proof}

We also have, using Remark \ref{rem:factpol} in dimension $0$:

\begin{lemme}\label{lem:conclprop}
The conclusion of Proposition \ref{prop:semisimplicity} holds for $\MLT(X,w;\varepsilong)^\rp$ and $\MLTr(X,w;\varepsilong)^\rp$.\qed
\end{lemme}

\begin{corollaire}\label{cor:psi0}
Let $(\cT,\cS,\bcL)$ be an object of $\MLT(X,w;\varepsilong)^\rp$ with strict support $Z$. Let $f:U\to\CC$ be a holomorphic function $\not\equiv0$ on $Z$. Then $\big(\gr_{\bbullet}^{M}\phi_{f,0}(\cT,\cS,\bcL),\cN\big)$ is an object of $\MLT(X,w+1;\varepsilong,-1)^\rp$. A similar result holds for regular objects.
\end{corollaire}

\begin{proof}
Apply the Lefschetz analogue of Corollary \ref{cor:vanishingcycles} and Lemma \ref{lem:conclprop}.
\end{proof}

\begin{proposition}\label{prop:declefcoh}
The conclusions of Propositions \ref{prop:declef} and \ref{prop:lefcoh} remain valid for graded Lefschetz (regular) twistor $\cD_X$-modules.
\end{proposition}

\begin{proof}
We will give the proof for nonregular objects, the regular case being similar. Let us begin with Proposition \ref{prop:declef}. First, we remark that $c(\cT_{j+1}),\ker v\subset\cT'_j$ are objects of $\MT(X,w-\varepsilon j)$, according to Proposition \ref{prop:abel}.

Let us show that $\im c$ and $\ker v$ are subobjects of $\cT'$ in $\MLT(X,w;\varepsilon)$. We may assume that $\cT,\cT'$ have strict support $Z$. Choose a local coordinate $t$ such that $\codim_Z(\{t=0\}\cap Z)=1$. We know that $c$ or $v$ and $\gr_\ell^{\rM}\Psi_{t,\alpha}$ commute (Proposition \ref{prop:abel}). It follows that, by induction, $\coker\cL^j:c(\cT_{j+1})\to c(\cT_{-j+1})$ is supported in $\{t=0\}$, hence is equal to $0$, as $\cT'$ has strict support $Z$. Argue similarly for $\ker v$.

By Lemma \ref{lem:conclprop}, $\im c$ and $\ker v$ decompose as direct sums of simple objects in $\MLT(X,w;\varepsilon)$, so their intersection is an object in the same category. By the same argument as above, using induction on the dimension, the intersection $\im c\cap \ker v$ vanishes. Similarly, the direct summand of $\im c\oplus\ker v$ in $\cT'$ is an object of $\MLT(X,w;\varepsilon)$ and also vanishes by induction. We therefore have a decomposition $\cT'=\im c\oplus\ker v$ in $\MLT(X,w;\varepsilon)$.

\medskip
Let us now consider Proposition \ref{prop:lefcoh}. So, let $\big((\cT_{j_1,j_2})_{\bmj\in\ZZ^2},\cL_1,\cL_2\big)$ be an object of $\MLT(X,w;\varepsilon_1,\varepsilon_2)^\rp$ with a polarization $\cS$. Let $d:\cT_{j_1,j_2}\to\cT_{j_1-1,j_2-1}(\varepsilon_1+\varepsilon_2)$ be a differential in $\RTriples(X)$, which commutes with $\cL_1$ and $\cL_2$ and is selfadjoint with respect to $\cS$. As both source and target of $d$ are in $\MT(X,w-\varepsilon_1j_1-\varepsilon_2j_2)$, $d$ is a morphism in this category, hence is strict and strictly specializable (Proposition \ref{prop:abel}) and we have, for any germ $f$ of holomorphic function any $\alpha$ with $\reel(\alpha)\in[-1,0[$ and any $\ell\geq0$,
\[
P\gr_\ell^{\rM}\Psi_{f,\alpha}(\ker d/\im d)=\ker(P\gr_\ell^{\rM}\Psi_{f,\alpha} d)/\im (P\gr_\ell^{\rM}\Psi_{f,\alpha} d)
\]
(\cf Corollary \ref{cor:Lefschetzdecmorphism}). By induction on the dimension of the support, we may assert that $\big(P\gr_\ell^{\rM}\Psi_{f,\alpha}(\ker d/\im d),P\gr_\ell^{\rM}\Psi_{f,\alpha}\bcL,P\gr_\ell^{\rM}\Psi_{f,\alpha}\cS\big)$ is an object of the category $\MLT(X\cap f^{-1}(0),w+\ell,\varepsilong)$ and we conclude with Lemma \ref{lem:semisimplicitelef}.
\end{proof}

\begin{corollaire}[Degeneration of a spectral sequence]\label{cor:suitespecdeg}
Let $\cT^{\cbbullet}$ be an object of the category $D^+(\RTriples(X))$ equipped with a Hermitian duality $\cS:\cT^{\cbbullet}\to(\cT^{\cbbullet})^*(-w)$ and with $\cL:\cT^{\cbbullet}\to\cT^{\cbbullet}[2](1)$ and $\cN:\cT^{\cbbullet}\to\cT^{\cbbullet}(-1)$ which commute and are skewadjoint with respect to $\cS$. Assume that $\cN$ is nilpotent and and that each term $E_1^{i,j-i}$ of the spectral sequence associated to the monodromy filtration of $\rM_\bbullet(\cN)$ is part of an object
\[
\ooplus_{i,j}\big(E_1^{i,j-i}=\cH^j(\gr_{-i}^{\rM}\cT^{\cbbullet}),\cH^j\gr_{-i}^{\rM}\cN,\cH^j\gr_{-i}^{\rM}\cL,\cH^j\gr_{-i}^{\rM}\cS\big)
\]
of $\MLT(X,w;-1,1)^\rp$. Then,
\begin{enumerate}
\item
the spectral sequence degenerates at $E_2$,
\item
the filtration $W_\bbullet\cH^j(\cT^{\cbbullet})$ naturally induced by $\rM_\bbullet\cT^{\cbbullet}$ is the monodromy filtration $\rM_\bbullet$ associated to $\cH^j\cN:\cH^j(\cT^{\cbbullet})\to\cH^j(\cT^{\cbbullet})$,
\item
the object
\[
\ooplus_{i,j}\big(\gr_{-i}^{\rM}\cH^j(\cT^{\cbbullet}),\gr_{-i}^{\rM}\cH^j\cN,\gr_{-i}^{\rM}\cH^j\cL,\gr_{-i}^{\rM}\cH^j\cS\big)
\]
is an object of $\MLT(X,w;-1,1)^\rp$.
\end{enumerate}
A similar result holds for regular objects.
\end{corollaire}

\begin{proof}
We know (Lemma \ref{lem:d1}) that the differential $d_1$ is selfadjoint with respect to $\cH^j\gr_{-i}^{\rM}\cS$. Moreover, $d_1:\cH^j(\gr_{-i}^{\rM}\cT^{\cbbullet})\to\cH^j(\gr_{-i}^{\rM}\cT^{\cbbullet})$ is a morphism between objects in $\MT(X,w+j-i)$. From the analogue of Proposition \ref{prop:lefcoh}, we deduce that $(E_2^{i,j-i})$ is part of an object of $\MLT(X,w;-1,1)^\rp$. Now, one shows inductively that $d_r=0$ for any $r\geq2$, by applying Proposition \ref{prop:wsupw}. This gives the result.
\end{proof}

\Subsection{A conjecture}\label{subsec:conj}
We now restrict our discussion to regular objects. For the nonregular case, we lack at the moment of results in dimension one.

\begin{theoreme}\label{th:ssimple}
Assume that $X$ is a complex projective manifold. The functor which associates to each object $(\cT,\cS)$ in $\MTr(X,w)^\rp$ the regular holonomic $\cD_X$-module $\Xi_{\DR}\cM''$ (restriction to $\hb_o=1$) takes values in the category of semisimple regular holonomic $\cD_X$-modules.
\end{theoreme}

\begin{proof}
We may assume that $\cT$ has strict support the irreducible closed analytic subset $Z\subset X$. We know that $\Xi_{\DR}\cM''$ is regular holonomic and has strict support $Z$, according to Proposition \ref{prop:restricanvar}. This means that there exists a dense Zariski open set $Z^o\subset Z$ and a local system $\cL$ of finite dimensional $\CC$-vector spaces on $Z^o$ such that the de~Rham complex of $\Xi_{\DR}\cM''$ is isomorphic to the intersection complex $\IC^{\cbbullet}\cL$ up to a shift. We want to show that the local system $\cL$ is semisimple. We will argue by induction on $\dim Z$, starting from $\dim Z=1$.

\smallskip
The case when $Z$ is a \emph{smooth} curve is a corollary of Theorem~\ref{th:ssimplecourbes} proved later in Chapter~\ref{chap:curves}, as a consequence of results of C\ptbl Simpson and O\ptbl Biquard.

\smallskip
Let us now consider the case when $Z$ is singular. Denote by $\nu:\wt Z\to Z$ the normalization. We may assume that $Z^o$ is an open set in $\wt Z$. It is therefore enough to construct an object $(\wt\cT,\wt\cS)$ in $\MTr(\wt Z,w)^\rp$ which coincides with $(\cT,\cS)$ on~$Z^o$. The singular points of $Z$ being isolated, this is then a local problem on $Z$, as we may glue local solutions to this problem with the solution $(\cT,\cS)_{|Z^o}$ on $Z^o$. The noncharacteristic inverse image $p^+(\cT,\cS)$ by the projection $p:\wt Z\times X\to X$ is an object of $\MTr(\wt Z\times X,w)^\rp$. Choose a family of local equations $t_1,\dots,t_n$ of the graph $G(\nu)\subset \wt Z\times X$ of $\nu:\wt Z\to X$. Then the object $P\gr_0^{\rM}\Psi_{t_1}\big(\cdots P\gr_0^{\rM}\Psi_{t_n}p^+(\cT,\cS)\big)$ gives a local solution to the problem.

\smallskip
Assume now that $\dim Z\geq2$. According to \cite[th\ptbl 1.1.3(ii)]{H-L85}, it is enough to show that the restriction of $\cL$ to a generic hyperplane section of $Z^o$ is semisimple, because for such a hyperplane, $\pi_1(Z^o\cap H)\to\pi_1(Z^o)$ is onto. Now, the (noncharacteristic) restriction of $(\cT,\cS)$ to a generic hyperplane $H$ still belongs to $\MTr(X\cap H,w)^\rp$, because it can be locally expressed as $P\gr_0^{\rM}\psi_{t,-1}(\cT,\cS)$ for a local equation $t$ of $H$. We therefore get by induction the semisimplicity of $\cL_{|Z^o\cap H}$, hence of $\cL$.
\end{proof}

\begin{conjecture}\label{conj:ssimple}
The functor above is an equivalence.
\end{conjecture}

This assertion should also hold when $X$ is compact and K\"ahler. Its proof would give a proof of the conjecture of M\ptbl Kashiwara recalled in the introduction, for semisimple perverse sheaves or regular holonomic $\cD$-modules.

\begin{remarque}\label{rem:ssimple}
When the conjecture holds, the functor sends a simple object in the first category in a simple object of the second one. This will be the case when $X$ is a compact Riemann surface, as a consequence of Theorem \ref{th:ssimplecourbes}.
\end{remarque}

\Subsection{Polarizable Hodge $\cD$-modules and polarizable twistor $\cD$-modules}\label{subsec:HodgeDM}
One may develop a theory of Hodge $\cD$-modules along the lines of this chapter. We will indicate the main steps.

In dimension $0$, polarized Hodge $\cD$-modules correspond to polarized complex Hodge structures as in \T\ref{subsec:chsmtw0}.

In general, replace the category of $\cR_\cX$-modules with the category of graded $R_F\cD_X$-modules, the morphisms being graded. Strict objects correspond to $\cD_X$-modules equipped with a good filtration (by the Rees construction). In order to define graded $\RTriples$, consider sesquilinear pairings $C$ taking values in $\Db_{X_\RR}\otimes_\CC\CC[\hb,\hbm]$.

In the definition of specializable graded $R_F\cD_X$-modules, one should insist on the fact that the $V$-filtration is graded. Therefore, if $\cM=R_FM$ for some well-filtered $\cD_X$-module $M$, and if all $\psi_{f,\alpha}\cM$ are strict, then $M$ is specializable along $\{f=0\}$ and $\psi_{f,\alpha}\cM=R_F\psi_{f,\alpha}M$, where the filtration $F$ on $\psi_{f,\alpha}M$ is naturally induced from that of $M$ as in \cite{MSaito86}.

At this point, the definition of a complex Hodge $\cD$-module is obtained by working in the category of graded $\RTriples$ when considering Definition \ref{def:regtwt}. It is very similar to the category considered by M\ptbl Saito.

The polarization is introduced as a graded isomorphism between both $R_F\cD_X$-modules entering in the definition of a complex Hodge $\cD$-module.

The graded analogue of Conjecture \ref{conj:ssimple} asserts that the category of complex Hodge $\cD$-modules is equivalent to that of admissible variations of polarized complex Hodge structures. It would follow from a direct comparison with the category of complex Hodge modules constructed by M\ptbl Saito, by the results of \cite{MSaito86,MSaito87}.

Another approach is indicated in Chapter~\ref{chap:int}.

\chapter{Polarizable regular twistor $\cD$-modules on~curves} \label{chap:curves}
In this chapter we will prove:

\begin{theoreme}\label{th:ssimplecourbes}\label{thssimplecourbes}
Conjecture \ref{conj:ssimple} is true when $X$ is a compact Riemann surface.
\end{theoreme}

This is nothing but a reformulation of some of the main results in \cite{Simpson90}. Nevertheless, we will give details on the reduction to this result, as this is not completely straightforward. Moreover, we will use the more precise description given in \cite{Biquard97}.

We will begin with the detailed computation of a basic example when $\dim X=1$. It corresponds to ``nilpotent orbits'' in dimension one. It was considered in detail in \cite{Simpson90} and \cite{Biquard97}.

\section{A basic example}\label{subsec:basic}

Let $\beta\in\CC$ and put $\beta'=\reel\beta$, $\beta''=\im\beta$. Let $V^o$ be a $\CC$-vector space of dimension $d$ equipped with a $\sld$-triple $(\rY,\rX,\rH)$ and with a positive definite Hermitian form such that $\rX^*=\rY$, $\rY^*=\rX$ and $\rH^*=\rH$. Fix an orthonormal basis $\bmv^o=(v_{1}^{o},\ldots,v_{d}^{o})$ of eigenvectors for $\rH$ and let $w_i\in\ZZ$ be the eigenvalue of $\rH$ corresponding to $v_{i}^{o}$. It will be convenient to assume that the basis $\bmv^o$ is obtained as follows: fix an orthonormal basis $v_1^o,\dots,v_k^o$ of $\ker\rX$ made with eigenvectors of $\rH$; for any $j=1,\dots,k$, consider the vectors
\begin{equation}\label{eq:vjl}
v_{j,\ell}^o=c_{j,\ell}\rY^\ell v_j^o,
\end{equation}
for $\ell=0,\dots,w_j$, where $c_{j,\ell}$ is some positive constant; the basis $(v_{j,\ell}^o)_{j,\ell}$ is orthogonal, and one can choose $c_{j,\ell}$ (with $c_{j,0}=1$) such that the basis is orthonormal.

All along this section \ref{subsec:basic}, we denote by $X$ the unit disc centered at $0$ with coordinate~$t$ and by $X^*$ the punctured disc $X\moins\{t=0\}$. Let $H=\cC^\infty_X\otimes_\CC V^o$ be the trivial $C^\infty$-bundle on $X$ and let $\bmv=(v_1,\ldots,v_d)$ be the basis such that $v_i=1\otimes v_{i}^{o}$. We still denote by $\rY,\rX,\rH$ the corresponding matrices in this basis and by $H$ the restriction of the bundle $H$ to the punctured disc $X^*$. Put on $H$ the logarithmic connection $D_V$ such that
\begin{align*}
D''_V\bmv&=0,\\
D'_V\bmv&=\bmv\cdot(\rY+\beta\id)\dlt.
\end{align*}
Put $\Lt=\module{\log t\ov t}$ as in \T\ref{subsec:defLt} and let $\varepsilong$ be the basis obtained from $\bmv$ by the change of basis of matrix
\[
P=\mt^{-\beta}\Lt^{-\rHsd}e^{\rX},
\]
that is,
\[
(\varepsilon_1,\ldots,\varepsilon_d)=(v_1,\ldots,v_d)\cdot P(t).
\]
Put on $H_{|X^*}$ the Hermitian metric $h$ such that $\varepsilong$ is an orthonormal basis. Put $\cH=\cC_{\cX^*}^{\infty,\an}\otimes_{\cC_{X^*}^\infty}H$. It has a basis $\bme(\hb)=\bmv\cdot R_\hb\in\Gamma(\cX^*,\cH)$ with
\[
R_\hb(t)=\Lt^{-\rHsd}\lefpar \mt^{i\beta''}e^{-\rX}\rigpar^{(\hb-1)}\Lt^{\rHsd}, \quad \text{for } t\in X^*,\ \hb\in \Omega_0.
\]
Put $\bme=\bme(0)$. One also has $\bme(\hb)=\bme\cdot Q_\hb$ with
\[
Q_\hb(t)=\Lt^{-\rHsd}\lefpar \mt^{i\beta''}e^{-\rX}\rigpar^\hb \Lt^{\rHsd}.
\]
The metric $h$ and the connection $D_V$ on $H$ allow to define operators $D'_E$, $D''_E$, $\theta'_E$ and $\theta''_E$ (see \T\ref{subsub:a} and \cite{Simpson92}).

\begin{proposition}\label{prop:basic}
For any $\beta\in\CC$ we have:
\begin{enumerate}
\item\label{prop:basic1}
the metric $h$ on $H$ is harmonic;
\item\label{prop:basic2}
the basis $\bme(\hb)$ is holomorphic with respect to the holomorphic structure on $\cH$ defined by $D''_E+\hb\theta''_E$;
\item\label{prop:basic3}
the action of $\hb D'_E+\theta'_E$ defines a left $\cR_{\cX}$-module structure on the free $\cO_{\cX}[1/t]$-module $\wt\cM\defin\cO_{\cX}[1/t]\cdot\bme(\hb)\subset j_*\cH$ (where $j:X^*\hto X$ denotes the inclusion) and $\wt\cM$ is strict and strictly specializable (\cf \T\ref{subsec:loc});
\item\label{prop:basic4}
the minimal extension $\cM$ of $\wt\cM$ across $t=0$ is strict holonomic;
\item\label{prop:basic5}
using notation of \T\ref{subsec:smtwqc} on $X^*$, the sesquilinear pairing $h_\bS:\cH'_{|\bS}\otimes_{\cO_\bS}\ov{\cH'_{|\bS}}\to\cC_{\cX^*}^{\infty,\an}$ extends to a sesquilinear pairing
\[
C:\cMS\ootimes_{\cO_\bS}\ov{\cMS} \to\Dbh{X}.
\]
\end{enumerate}
\end{proposition}

\begin{proof}
We will use the following identities:
\begin{equation}\label{eq:sl2}
\begin{split}
\Lt^{\pm \rHsd}\rY\Lt^{\mp \rHsd}&=\Lt^{\mp1}\rY\\
\Lt^{\pm \rHsd}\rX\Lt^{\mp \rHsd}&=\Lt^{\pm1}\rX\\
e^{\rY}\rH e^{-\rY}&=\rH+2\rY\\
e^{\rX}\rH e^{-\rX}&=\rH-2\rX\\
e^{\rX}\rY e^{-\rX}&= \rY+\rH-\rX
\end{split}
\end{equation}
and, from\eqref{eq:diffLt},
\begin{equation}\label{eq:at}
\Lt^{\rHsd}t\;\frac{\partial\lefpar \Lt^{-\rHsd}\rigpar}{\partial t} =\Lt^{\rHsd}\ov t\;\frac{\partial\lefpar \Lt^{-\rHsd}\rigpar}{\partial \ov t} =\frac{\rHsd}{\Lt}.
\end{equation}
Write in the basis $\varepsilong$
$$
D'_V\varepsilong=\varepsilong\cdot M'\dlt,\quad D''_V\varepsilong=\varepsilong\cdot M''\dltb.
$$
One has $M'=\beta\id+P^{-1}\rY P+P^{-1}t\partial_tP$ and $M''=P^{-1}\ov t\partial_{\ov t}P$. According to the previous identities, one gets
\begin{align}
M'&=\frac{\beta}{2} \id+\frac{\rY-\rHsd}{\Lt}\label{eq:M'}\\
M''&=-\frac{\beta}{2}\id+\frac{\rX+\rHsd}{\Lt}\label{eq:M''}\\
\theta'_E&=\frac12(M'+M^{\prime\prime*})\dlt=\lefpar i\frac{\beta''}{2}\id+\frac{\rY}{\Lt}\rigpar \dlt\label{eq:theta'}\\
\theta''_E&=\frac12(M^{\prime*}+M'')\dltb=\lefpar -i\frac{\beta''}{2}\id+\frac{\rX}{\Lt}\rigpar \dltb\label{eq:theta''}\\
D''_E\varepsilong&=(D''_V-\theta''_E)\varepsilong=\varepsilong\cdot\lefpar -\frac{\beta'}{2}\id+\frac{\rHsd}{\Lt}\rigpar \dltb.\notag
\end{align}
Now, the matrix of $D''_E+\hb\theta''_E$ in the basis $\bme(\hb)$ is zero (which gives the second point): indeed, we have $\bme(\hb)=\bmv\cdot R_\hb=\varepsilong\cdot P^{-1}R_\hb$ and
\begin{equation}\label{eq:PRh}
P^{-1}R_\hb=\mt^{\beta'+i\hb\beta''}e^{-\hb \rX}\Lt^{\rHsd}\defin A_\beta(t,\hb).\index{$abthb$@$A_\beta(t,\hb)$}
\end{equation}
The matrix of $D''_E+\hb\theta''_E$ in $\varepsilong$ is
\[
\lefpar -\frac{\beta'+i\hb\beta''}{2}\id+\frac{(\rHsd+\hb \rX)}{\Lt}\rigpar \dltb,
\]
hence in the basis $\bme(\hb)$ the coefficient of $\dltb$ is
\begin{multline}
-\frac{\beta'+i\hb\beta''}{2} \id+\Lt^{-\rHsd}e^{\hb \rX}\frac{\rHsd+\hb \rX}{\Lt}e^{-\hb \rX}\Lt^{\rHsd}
\\ +\lefpar \frac{\beta'+i\hb\beta''}{2}\rigpar\id-\frac{\rHsd}{\Lt}.
\end{multline}
But
\begin{align*}
e^{\hb \rX}(\rHsd+\hb \rX)e^{-\hb \rX}&=e^{\hb \ad \rX}(\rHsd+\hb \rX)\\
&=\hb \rX+e^{\hb \ad \rX}(\rHsd)=\hb \rX+\rHsd-\hb \rX=\rHsd,\notag
\end{align*}
hence the result.

On the other hand, the matrix of $\hb D'_E+\theta'_E$ in the basis $\bme(\hb)$ is computed in the same way: it is equal to
\[
\lefcro \Big( i\frac{\beta''}{2}(1+\hb^2)+\hb\beta'\Big) \id+\rY\rigcro\dlt=\lefcro \lefpar \beta\star\hb\rigpar \id+\rY\rigcro\dlt.
\]
Therefore, $\cO_{\cX}[1/t]\cdot\bme(\hb)$ is a left $\cR_{\cX}$-module, which gives \eqref{prop:basic3}. By definition, one has
\begin{equation}\label{eq:tpartiallt}
t\partiall_t\bme(\hb)=\bme(\hb)\cdot\lefcro \lefpar \beta\star\hb\rigpar \id+\rY\rigcro.
\end{equation}
Moreover, putting $\hb=0$ shows that the matrix of $\theta'_E$ in the basis $\bme=\bme(0)$ is holomorphic, hence $h$ is harmonic. We have obtained \ref{prop:basic}\eqref{prop:basic1} and \eqref{prop:basic2}.

\medskip
Consider the filtration $U_k\wt\cM=t^{-k}\cO_\cX\cdot\bme(\hb)$. This is a good filtration with respect to $V_\bbullet\cR_\cX[1/t]$ and, for any $k\in\ZZ$, we have $\gr_k^U\wt\cM\simeq\cO_{\Omega_0}^d$. Moreover, putting $\alpha=-\beta-1$, the operator $\partiall_tt+(\alpha+k)\star\hb$ is nilpotent on $\gr_k^U\wt\cM$. Near any $\hb_o\in\Omega_0$, the filtration $V_\bbullet^{(\hb_o)}\wt\cM\defin U_{\bbullet-\ell_{\hb_o}(\alpha)}\wt\cM$ satisfies all the properties of the Malgrange-Kashiwara filtration, hence is equal to it. This shows \ref{prop:basic}\eqref{prop:basic3}. The lattice $\Lambda$ considered after Definition \ref{def:spe} reduces here to $\alpha+\ZZ$.

Although the filtration $V_\bbullet^{(\hb_o)}\wt\cM$ is only locally defined with respect to $\hb_o$ when the imaginary part $\alpha''$ is not $0$, the module $\psi_{t,\alpha+k}\wt\cM$ (for $k\in\ZZ$) is globally identified with $U_k/U_{k-1}$.

\medskip
Consider the minimal extension $\cM$ of $\wt\cM$ across $t=0$. By the results of \T\ref{subsec:minext}, it is strictly specializable along $t=0$. It is strict, because $\wt\cM$ is so. Near $\hb_o\in\Omega_0$, $\cM$ is the $\cR_\cX$-submodule of $\wt\cM$ generated by $t^{-k_{\hb_o}}\bme(\hb)$, where $k_{\hb_o}\in\ZZ$ is chosen such that $k_{\hb_o}+\ell_{\hb_o}(\alpha)\in[-1,0[$.

$\cM$ is holonomic because its characteristic variety is contained in $(T^*_DD\cup T^*_0D)\times \Omega_0$ (by an extension argument, reduce to the case where $\rY=0$). Remark also that the support $\Sigma(\cM)$ of $\Xi_{\Dol}\cM$ is contained in the curve $t\tau=i\beta''/2$: near $\hb_o$, the classes of $t^{-k_{\hb_o}}e_j(\hb)$ ($j=1,\ldots,d$) generate $\Xi_{\Dol}\cM$ over $\cO_D[TD]$; these classes satisfy the equation
\[
\det\lefpar (t\tau-i\beta''/2)\id-\rY\rigpar \cdot t^{-k_{\hb_o}}e_j(0)=0\quad(j=1,\ldots,d).
\]

Let us end by proving \eqref{prop:basic5}. It will simplify notation and not be restrictive to assume that $\reel\beta\in{}]-1,0]$, that is, putting $\alpha=-\beta-1$, $\reel\alpha\in[-1,0[$. Fix $\hb_o\in\Omega_0$. In the local basis $t^{-k_{\hb_o}}\bme(\hb)$ of $\cM$ near $\hb_o$, the matrix $\bC^{(\hb_o)}=(C^{(\hb_o)}_{ij})$ with entries $C^{(\hb_o)}_{ij}=h_\bS(t^{-k_{\hb_o}}e_i(\hb),\ov{t^{-k_{\hb_o}}e_j(\hb)})$ defines a $\ccR_{X^*,\ov X^*,\bS}$ sesquilinear pairing.

Let $A(\hb)$ be the matrix defined by \eqref{eq:PRh}, so that $\bme(\hb)=\varepsilong\cdot A(\hb)$, where $\varepsilong$ is the $h$-orthonormal basis constructed previously. Put
\[
\wt\bC={}^t\!A(\hb)\cdot \ov{A(\hb)}.
\]
Formula \eqref{eq:PRh} shows that
\begin{align}\label{eq:Cehb}
\wt\bC&=\mt^{2(\beta\star\hb)/\hb}\cdot \Lt^{\rHsd}e^{-\hb \rY} e^{\rX/\hb} \Lt^{\rHsd}\defin\mt^{2(\beta\star\hb)/\hb}\cdot B.
\end{align}
The matrix $\bC^{(\hb_o)}$ can be written as $\mt^{-2k_{\hb_o}}\wt\bC$. When $t\to0$, each $\vert C^{(\hb_o)}_{ij}\vert$ behaves therefore like $\mt^{2\ell_{\hb_o}(\beta-k_{\hb_o})}\Lt^k$ for some $k\in\ZZ$. By definition, we have $\ell_{\hb_o}(\beta-k_{\hb_o})\in{}]-1,0]$, so that $C^{(\hb_o)}_{ij}$ is $L^1_\loc$ and defines a distribution depending analytically on $\hb$ near $\hb_o\in\bS$.

We hence \emph{define} $C$ on $\cMS\otimes_{\cO_\bS}\ov{\cMS}$ as the unique (if it exists) $\ccR_{X,\ov X,\bS}$-linear pairing such that $C(t^{-k_{\hb_o}}e_i(\hb),\ov{t^{-k_{\hb_o}}e_j(\hb)})$ is the $L^1_\loc$ extension of $C^{(\hb_o)}_{ij}$ as a distribution. Uniqueness is clear, as we are given $C$ on generators. It will also allows us to glue along $\bS$ the various local constructions. It remains to prove the existence.

If one chooses a basis $\bmv^o=(v_{j,\ell})_{j,\ell}$ as in \eqref{eq:vjl}, the matrix $\bC^{(\hb_o)}$ is block-diagonal, $\bC^{(\hb_o)}=\oplus_j\bC^{(\hb_o)}_{(j)}$. We can therefore easily reduce to the case where $\rY$ has only one Jordan block.

Under such an assumption, $\cM$ is $\cR_\cX$-generated by $t^{-k_{\hb_o}}e_1(\hb)$. As we have $t^{-k_{\hb_o}}e_{\ell+1}(\hb)=[t\partiall_t-(\beta-k_{\hb_o})\star\hb]^\ell t^{-k_{\hb_o}}e_1(\hb)$, we first have to verify that we indeed have, as distributions,
$$
[t\partiall_t-(\beta-k_{\hb_o})\star\hb]^{\ell'}\ov{[t\partiall_t-(\beta-k_{\hb_o})\star\hb]}^{\ell''}C^{(\hb_o)}_{11}=C^{(\hb_o)}_{1+\ell',1+\ell''}.
$$
We know that this holds on $X^*$ as $C^\infty$ functions. It then holds as distributions for the $L^1_\loc$ extensions (\cf Example~\ref{ex:distr}).

Notice now that we have a local presentation of $\cM$ (recall that $\reel\beta\in{}]-1,0]$):
\begin{align*}
\cR_\cX\To{{}\cdot[t\partiall_t-(\beta-k_{\hb_o})\star\hb]^d} \cR_\cX\to\cM\to0&\quad\text{if }\beta\not\in\ZZ,\\
\cR_\cX\To{{}\cdot\partiall_t(t\partiall_t)^{d-1}}\cR_\cX\to\cM\to0&\quad\text{if }\beta=0.
\end{align*}
Indeed, we have a surjective morphism of the cokernel to $\cM$, and it is enough to show that the cokernel has no $t$-torsion, which can be seen easily. Therefore, $C$ will be well defined if we show that $C^{(\hb_o)}_{11}$ satisfies $[t\partiall_t-(\beta-k_{\hb_o})\star\hb]^dC^{(\hb_o)}_{11}=0$ (when $\beta\not\in\ZZ$) or $\partiall_t(t\partiall_t)^{d-1}C^{(\hb_o)}_{11}=0$ (when $\beta=0$).

By construction, this holds on $X^*$, so that we can write as $C^\infty$ functions on $X^*$:
$$
C^{(\hb_o)}_{1,1}=C(t^{-k_{\hb_o}}e_1,\ov{t^{-k_{\hb_o}}e_1})=\mt^{2((\beta-k_{\hb_o})\star\hb)/\hb}\sum_{k=0}^{d-1}a_k\Lt^k/k!
$$
for some integers $a_k$, with $a_{d-1}=1$. By definition of the extension $C^{(\hb_o)}_{1,1}$, this also holds as $L^1_\loc$ functions on $X$. Apply now Example~\ref{ex:distr}.
\end{proof}

\begin{Remarques}
\begin{enumerate}
\item
In particular, we have $\bC^{(\hb_o)}_{i,j}=0$ for $i+j>d$ and, in the expression \eqref{eq:Cehb} for $\bC^{(\hb_o)}$, the coefficients of the negative powers of $\Lt$ vanish. This can also be seen using the relation
\[
\big(e^{\rX/\hb}e^{-\hb \rY}e^{\rX/\hb}\big)\cdot \rH\cdot \big(e^{\rX/\hb}e^{-\hb \rY} e^{\rX/\hb}\big)^{-1}=-\rH.
\]
\item
As we assume that $\reel\beta\in{}]-1,0]$, the entries of $\wt\bC$ take value in $\Dbh{X}$. Moreover, the sesquilinear pairing $\psi_{t,\alpha}C$ on $\Psi_{t,\alpha}\cM$ can be directly computed by using~$\wt\bC$.
\end{enumerate}
\end{Remarques}

We will end this paragraph by proving:

\begin{proposition}\label{prop:basicpolarisation}
Put $\alpha=-\beta-1$. Then, for any $\ell\geq0$, the object
\[
\big(P\gr_\ell^{\rM}\Psi_{t,\alpha}(\cM),P\gr_\ell^{\rM}\Psi_{t,\alpha}(\cM), P\psi_{t,\alpha,_\ell}C,\cS=\id)
\]
is a polarized twistor structure of weight~$\ell$ in the sense of \T\ref{subsubsec:smtw0}.
\end{proposition}

\begin{proof}
Keeping notation as above, we have $\psi_{t,\gamma}\cM=0$ if $\gamma\not\equiv\alpha\mod\ZZ$ and $\Psi_{t,\alpha}\cM\simeq\cO_{\Omega_0}^d$ has basis $[\bme(\hb)]$. The matrix of $i\rN=-i[\partiall_tt+\alpha\star\hb]$ in this basis is $-i\rY$. We compute $\psi_{t,\alpha}C$ on $\Psi_{t,\alpha}\cM$ with the help of the matrix $\wt\bC$.

For simplicity, we will assume that $\rY$ has only one Jordan block, of size $d$. Therefore, $P\gr_\ell^{\rM}\Psi_{t,\alpha}(\cM)=0$ if $\ell\neq d-1$ and has dimension $1$ if $\ell=d-1$.

It is a matter of verifying (see Formula \eqref{eq:psipolell}) that the expression
\begin{equation}\label{eq:basicpositif}
\thb^{-(d-1)}\cdot\res_{s=\alpha\star\hb/\hb}\int\mt^{2s}\wt C((i\rN)^{d-1}e_1,\ov{e_1})\chi(t) \, \itwopi dt\wedge d\ov t,
\end{equation}
(for $\chi\in C^\infty_c$, $\chi\equiv 1$ near $t=0$) considered as a function of $\hb$, is a positive constant.

Put $\rY^{d-1}e_1=y_de_d$. As $i\rN$ acts as $-i\rY$, the expression \eqref{eq:Cehb} for $\wt\bC$ implies that
\[
\wt C_{d,1}=\mt^{2(\beta\star\hb)/\hb}(-\hb)^{d-1}y_d/(d-1)!,
\]
hence $\wt C((i\rN)^{d-1}e_1,\ov{e_1})=(-i)^{d-1}y_d\wt C_{d,1} =\mt^{2(\beta\star\hb)/\hb}\thb^{d-1}y_d^2/(d-1)!$. We therefore have
\[
\eqref{eq:basicpositif}= \frac{y_d^2}{(d-1)!}\cdot \res_{s=\alpha\star\hb/\hb}\int\mt^{2(s-\alpha\star\hb/\hb)}\chi(t)\,\itwopi \frac{dt}{t}\wedge\frac{d\ov t}{\ov t}.
\]
Now, use that
\[
\res_{s=0}\int\mt^{2s}\chi(t)\,\itwopi \frac{dt}{t}\wedge\frac{d\ov t}{\ov t}=1. \qedhere
\]
\end{proof}

\section{Review of some results of C\ptbl Simpson and O\ptbl Biquard}
In this section, $X$ denotes a compact Riemann surface and $P$ a finite set of point of $X$. We also set $X^*=X\moins P$ and we denote by $\cI_P$ the ideal of $P$ (as a reduced set).

\subsection{}
By a \emph{meromorphic bundle on $X$ with poles on $P$} we mean a locally free $\cO_X(*P)$-module of finite rank. Let $\wt V$ be such a bundle. A \emph{meromorphic connection on $\wt V$} is a $\CC$-linear morphism $\nabla:\wt V\to\Omega_X^1\otimes_{\cO_X}\wt V$ satisfying the usual Leibniz rule. We denote by $\wt M$ the meromorphic bundle $\wt V$ with connection $\nabla$, viewed as a left $\cD_X$-module. There is an equivalence between the category of semisimple regular holonomic $\cD_X$-modules with singularities at $P$ and strict support $X$, and the category of semisimple meromorphic bundles with connection having regular singularities at $P$: in one direction, associate to the $\cD_X$-module $M$ the meromorphic bundle $\wt M=\cO_X(*P)\otimes_{\cO_X}M$; in the other direction, use $\nabla$ to put the structure of a regular holonomic module on $\wt V$ and associate to $\wt M$ the \emph{minimal extension} $M\subset \wt M$, \ie the biggest $\cD_X$-submodule of $\wt M$ having no quotient supported in a finite set of points.

Notice also that, by the Riemann-Hilbert correspondence, these categories are equivalent to the category of semisimple representations of $\pi_1(X^*)$.

If we are given a decreasing filtration $\wt M^{\cbbullet}$ of $\wt M$, indexed by a finite set $B_\RR\subset{}]-1,0]$, by $\cO_X$-locally free submodules $\wt M^b$ such that $\wt M/\wt M^b$ is supported on $P$ for any $b\in B_\RR$ and on which the connection has at most logarithmic poles, we say, following \cite{Simpson90}, that $(\wt M,\wt M^{\cbbullet},\nabla)$ is a \emph{filtered regular meromorphic connection}. The filtration may be extended to indices in $B_\RR+\ZZ$ by putting $\wt M^{b+k}=\cI_P^k\wt M^b$.

\subsection{}\label{subsec:canonique}
Assume that $\wt M$ has only regular singularities. Consider the \emph{canonical} filtration of $\wt M$ (\ie the Malgrange-Kashiwara filtration, in dimension one): it is indexed by $A_\RR+\ZZ$ for some finite set $A\in\CC$, putting $A_\RR=\{\reel\alpha\mid\alpha\in A\}$. We will use the increasing version of it: each $V_a\wt M$ is a locally free $\cO_X$-module which coincides with $\wt M$ on $X^*$, and on which the connection acts with only simple poles, such that the eigenvalues of the residue have real part in $[-(a+1),-a[$. We will also consider the decreasing version $V^{\cbbullet}\wt M$, by putting $V^b=V_{-(b+1)}$ and $\gr^b_V=\gr_{-(b+1)}^V$ (see Remark \ref{rem:leftrightspe}\eqref{rem:leftrightspe2}), so that the eigenvalues of the residue of $\nabla$ on $V^b\wt M$ have real part in $[b,b+1[$.

The \emph{degree} of a filtered regular meromorphic connection $(\wt M,\wt M^{\cbbullet},\nabla)$ is defined as
\[
\deg(\wt M,\wt M^{\cbbullet},\nabla)=\deg \wt M^0+\sum_{x\in P} \sum_{b\in[0,1[}b\dim\gr^b\wt M_x.
\]
By the residue formula we have:
\begin{lemme}\label{lem:residu}
If $\wt M^{\cbbullet}$ is the canonical filtration $V^{\cbbullet}\wt M$ of $\wt M$, then \hbox{$\deg(\wt M,\wt M^{\cbbullet})\!=\!0$}.\qed
\end{lemme}

Say that a filtered regular meromorphic connection $(\wt M,\wt M^{\cbbullet})$ is \emph{stable} if any nonzero sub meromorphic connection $(\wt N,\nabla)$, equipped with the induced filtration $\wt N^{\cbbullet}=\wt N\cap \wt M^{\cbbullet}$, satisfies
\[
\frac{\deg(\wt N,\wt N^{\cbbullet})}{\rg\wt N}<\frac{\deg(\wt M,\wt M^{\cbbullet})}{\rg\wt M}.
\]

Owing to the fact that the filtration induced on $\wt N$ by the canonical filtration of $\wt M$ is the canonical filtration of $\wt N$, and according to the previous lemma, we get:

\begin{lemme}\label{lem:irredstable}
If $\wt M^{\cbbullet}$ is the canonical filtration of $\wt M$, then $(\wt M,\wt M^{\cbbullet})$ is stable if and only if $\wt M$ is irreducible.\qed
\end{lemme}

\subsection{}
Let $V$ be a holomorphic bundle on $X^{*}$ and let $h$ be a hermitian metric on $H=\cC^\infty_{X^*}\otimes_{\cO_{X^*}} V$. Say that $h$ is \emph{moderate} if the subsheaf $\wt V$ of $j_{*}V$ consisting of sections of $j_{*}V$, the $h$-norm of which has moderate growth near $P$, is a meromorphic bundle on $X$. If we are given a meromorphic extension $\wt M$ of $V$, we also say that $h$ is moderate with respect to $\wt M$ if $\wt V=\wt M$. The \emph{parabolic filtration} $\wt V^{\cbbullet}$ of $\wt V$ is then the filtration by the order of growth: in a local coordinate $t$ near $x_o\in P$,
\[
\wt V^b_{x_o}=\{\sigma\in j_{*}V_{x_{o}} \mid\lim_{t\to0}\mt^{-b+\varepsilon}\norme{\sigma(t)}_{h}=0\text{ for } \varepsilon>0\text{ and }\varepsilon\ll1\}.
\]
A criterion for the coherence of the parabolic bundles is given in \cite[Prop\ptbl 3.1]{Simpson90}, after \cite{C-G75}.

Let now $(H,D_{V})$ be a flat bundle on $X^{*}$. Following \loccit, say that $h$ is \emph{tame} with respect to $(H,D_V)$ if the $h$-norm of flat sections of $V$ grows at most polynomially near~$P$. If $h$ is tame and \emph{harmonic}, C\ptbl Simpson has shown \cite[Th\ptbl 2]{Simpson90} that $h$ is moderate, that each term of the parabolic filtration is $\cO$-coherent and is logarithmic with respect to the connection $\nabla$. The object $(H,D_{V},h)$ is then called a \emph{tame harmonic bundle}. It therefore defines a filtered meromorphic bundle with connection $(\wt M,\wt M^{\cbbullet})$.

We say that the tame harmonic bundle $(H,D_{V},h)$ has \emph{Deligne type} if the parabolic filtration $\wt M^{\cbbullet}$ is the canonical filtration of $\wt M$.

The category of tame harmonic bundles (morphisms are compatible with $D_{V}$ and bounded with respect to the metrics) is semisimple, as well as the full subcategory of tame harmonic bundles of Deligne type (\loccit, Th\ptbl 5).

We will now use:

\begin{theoreme}[\cite{Simpson90}]\label{th:simpson90}
Let $(\wt M,\wt M^{\cbbullet})$ be a filtered regular meromorphic connection. Then $(\wt M,\wt M^{\cbbullet})$ is poly-stable, each summand having degree $0$, if and only if there exists a Hermitian metric $h$ on $\wt M_{|X^*}$ which is tame with respect to $\wt M$, with associated parabolic filtration $\wt M^{\cbbullet}$, and such that $(\wt M_{|X^*},\nabla,h)$ is harmonic.\qed
\end{theoreme}

By the previous lemmas, if we assume that $\wt M^{\cbbullet}$ is the canonical filtration, then we may replace in the previous theorem the word ``poly-stable'' by ``semisimple'', and forget about the condition on the degree, which is automatically satisfied.

\begin{corollaire}\label{cor:simpson90}
The functor $(H,D_{V},h)\mto(\wt M,\wt M^{\cbbullet})$ from the category of tame harmonic bundles to that of log-filtered meromorphic connections induces an equivalence between the subcategory of tame harmonic bundles of Deligne type to that of semisimple meromorphic bundles with a regular connection (or equivalently, semisimple regular holonomic $\cD_{X}$-modules having strict support $X$).\qed
\end{corollaire}

\subsection{}\label{subsec:debutsimple}
Let $\wt M$ be a simple meromorphic bundle on $X$ with poles on $P$ and a regular connection; it is isomorphic, locally near each point of $P$, to a direct sum, indexed by $\beta\in\CC$ with $\reel\beta\in{}]-1,0]$, of meromorphic bundles with connection as in \T\ref{subsec:basic}. Denote by $(V,\nabla)$ the restriction of $\wt M$ to $X^*$. Let $\wt M^{\cbbullet}$ denote the (decreasing) Malgrange-Kashiwara filtration of $\wt M$. Hence $(\wt M,\wt M^{\cbbullet})$ is a polystable regular filtered meromorphic bundle with connection of degree $0$ (Lemmas \ref{lem:residu} and \ref{lem:irredstable}), to which we may apply Theorem \ref{th:simpson90}.

Choose a model (also called standard) metric $h^{\std}$ on $V$ which is equal, near the singular points, to a corresponding direct sum of metrics as in Prop\ptbl \ref{prop:basic}. The Malgrange-Kashiwara filtration $\wt M^{\cbbullet}$ of $\wt M$ can be recovered from this metric by measuring the order of growth of the norm of local sections of $\wt M$ (this is easily seen on the simple basic example).

The result of C\ptbl Simpson is then more precise than Theorem \ref{th:simpson90}: it asserts (see \cite[Th\ptbl 6(2)]{Simpson90}) that there exists a harmonic metric $h$ which is \emph{comparable} with $h^{\std}$ near $P$.

It will be even more convenient to use the construction made by O\ptbl Biquard in \cite{Biquard97} which gives a more precise description of the relationship between $h$ and $h^{\std}$.

We keep notation of \S\T\ref{subsub:a} and \ref{sec:smtwqc}. Denote by $B'\subset\{\beta\in\CC\mid\reel\beta\in{}[0,1[\}$ the finite set of eigenvalues of the residue of $\nabla$ on $\wt M^0/\wt M^1$, and let $B$ be the set obtained from $B'$ by adding $-1$ to any $\beta\in B'$ such that $\reel\beta\neq0$: for $\beta\in B$, we have $\reel\beta\in{}]-1,0]$.

The following result is also valid for general tame harmonic bundles or poly-stable filtered regular meromorphic connection of degree $0$ that we will not consider here.

\begin{theoreme}[{\cite[\S\T9 and 11]{Biquard97}}]\label{th:biquard}
Let $(H,D_V,h)$ be a tame harmonic bundle on $X^*$ of Deligne type, or equivalently, let $\wt M$ be a simple meromorphic regular connection on $X$ with singularities at $P$. For each puncture in $P$, there exists, on a small disc $D$ centered at this puncture, an $h$-orthonormal basis $\varepsilong$ of $\cC_{D^*}^\infty\otimes V$ on $D^*$ such that the matrix of $D_V$ in this basis can be written as
\[
M'\frac{dt}{t}+M''\frac{d\ov t}{\ov t}=(M^{\prime\std}+P')\frac{dt}{t}+(M^{\prime\prime\std}+P'')\frac{d\ov t}{\ov t},
\]
where $M^{\prime\std},M^{\prime\prime\std}$ are direct sums indexed by $\beta\in B$ of matrices \eqref{eq:M'} and \eqref{eq:M''}, and $P',P''$ are of some H\"older type (see \loccit).\qed
\end{theoreme}

\subsection{Sketch of the proof of Theorem \ref{th:ssimplecourbes}} \label{subsec:sketch}
In \T\ref{sec:first}, we prove a reconstruction result, namely, starting from a tame harmonic bundle of Deligne type $(H,D_{V},h)$, we associate to it a polarized regular twistor $\cD_{X}$-module $(\cT,\cS)_{h}=(\cM,\cM,C,\id)$ of weight~$0$ which coincides with that given by Lemma \ref{lem:twharm} on $X^{*}$ and such that $\Xi_{\DR}\cM=M$ is the minimal extension of the meromorphic bundle with connection defined by $h$.

In \T\ref{sec:second}, we show that the correspondence $(\cT,\cS)\mto(H,D_{V},h)$ on $X^{*}$ of Lemma \ref{lem:twharm}, starting from a polarizable regular twistor $\cD_{X}$-module of weight~$0$, gives rise to a tame harmonic bundle of Deligne type. Moreover, we show that $(\cT,\cS)$ and $(\cT,\cS)_{h}$ constructed in \T\ref{sec:first} are isomorphic, at least locally on~$X$.

Both results are enough to conclude. Indeed, it is enough to prove Theorem \ref{th:ssimplecourbes} for objects having strict support $X$. The functor ``restriction to $\hb=1$'' sends polarized regular twistor $\cD$-modules to semisimple regular holonomic $\cD$-modules, according to \T\ref{sec:second}. It is essentially surjective, according to \T\ref{sec:first}. By the equivalence of Lemma \ref{lem:twharm} and Corollary \ref{cor:simpson90}, any morphism $\varphi:M_{1}\to M_{2}$ lifts to a morphism $\psi:(\cT_{1},\cS_{1})_{|X^{*}}\to (\cT_{2},\cS_{2})_{|X^{*}}$ in a unique way. Moreover, it induces a morphism between the meromorphic extensions $\wt\cM_1\to\wt\cM_2$, which is compatible with the parabolic filtration constructed in Corollary \ref{cor:biquard}, because $\varphi$ respects the canonical filtration, \ie the parabolic filtration of the harmonic metric.

By the construction of \T\ref{subsec:exemplestricspe}, $\psi$ extends then as a morphism $(\cT_{1},\cS_{1})_{h}\to (\cT_{2},\cS_{2})_{h}$. Therefore, by the local isomorphism proved in \T\ref{sec:second}, $\psi$ also extends (in a unique way) as a morphism $(\cT_{1},\cS_{1})\to (\cT_{2},\cS_{2})$, hence the full faithfulness of the functor.\qed

\section[Proof of Theorem $\ref{thssimplecourbes}$, first part]{Proof of Theorem \ref{th:ssimplecourbes}, first part}\label{sec:first}

In this section, we will prove that any semisimple regular holonomic $\cD_X$-module is of the form $\Xi_{\DR}\cM''$ where $\cM''$ is part of an object of $\MTr(X,0)^\rp$.

Let $M$ be a \emph{simple} regular holonomic $\cD_X$-module and let $\wt M$ be the associated simple meromorphic bundle with regular connection as in \T\ref{subsec:debutsimple}. We fix a harmonic metric as given by Theorem \ref{th:biquard}, and we will work locally near a point of $P$, with a local coordinate $t$ centered at this point.

\Subsection{Construction of the $\cO_\cX[1/t]$-module $\wt\cM$ and the filtration $V_\bbullet\wt\cM$}\label{subsec:constr}
In this paragraph, $X$ will denote an open disc with coordinate $t$.

Denote by $j:X^*\times\Omega_0\hto X\times\Omega_0$ the open inclusion; denote as above by $\cX$ the product $X\times\Omega_0$. For $\hb\in\Omega_0$, recall that $\imhb$ denotes its imaginary part; denote also by $\Delta_{\hb_o}(\eta)$ the closed disc centered at $\hb_o$ of radius $\eta>0$; recall that $\rH\big(\Delta_{\hb_o}(\eta)\big)$ denotes the corresponding Banach space of holomorphic functions. For short, we will denote by $\wt\cM_\hb$ the germ at $(0,\hb)\in\cX$ of the $\cO_\cX$-module $\wt\cM$.

The first step consists in defining a $\cO_\cX[1/t]$-module $\wt\cM$ and the parabolic filtration $V^{\cbbullet}\wt\cM$ on it. As we will see below, the parabolic filtration is only \emph{locally} defined on~$\cX$. What is possible to define globally are the graded pieces, or even the $\psi_t^\beta\wt\cM$.

\medskip
Let us first state a consequence of Theorem \ref{th:biquard}.

\begin{Corollaire}\label{cor:biquard}
\begin{enumerate}
\item\label{cor:biquard1}
Let $\wt\cM\subset j_*\cH'$ be the subsheaf of local sections of $j_*\cH'$, the $\pi^*h$-norm of which has moderate growth along $\{0\}\times\Omega_0$. Then $\wt\cM$ is $\cO_\cX[1/t]$-locally free and is strictly specializable along $t=0$ (\cf \T\ref{sec:minext}).
\item\label{cor:biquard2}
For any $\hb_0\in\Omega_0$, let $(V^b\wt\cM_{\hb_o})_{b\in\RR}$ be the ``parabolic filtration'' of the germ $\wt\cM_{\hb_o}$ associated with $\pi^*h$ near $\hb_o$, \ie $V^b\wt\cM_{\hb_o}$ is the set of germs $m\in\wt\cM_{\hb_o}$ which satisfy $\mt^{-b+\varepsilon}\norme{m}_{\pi^*h}$ bounded near $(0,\hb_o)$ for any $\varepsilon>0$. Then, this filtration coincides with the Malgrange-Kashiwara filtration $V_{(\hb_o)}^b(\wt\cM_{\hb_o})$.

\item\label{cor:biquard3}
The Malgrange-Kashiwara filtration moreover satisfies:
\begin{enumerate}
\item\label{cor:biquard3a}
each $V_{(\hb_o)}^b(\wt\cM_{\hb_o})$ is $\cO_{\cX,(0,\hb_o)}$-locally free of finite rank,
\item\label{cor:biquard3b}
the monodromy filtration $\rM_\bbullet(\rN)$ (\cf Remark \ref{rem:psi}\eqref{rem:Ncanvar1}) of the nilpotent endomorphism $\rN=-(t\partiall_t-\beta\star\hb):\psi_{t}^\beta(\wt\cM_{\hb_o})\to\psi_{t}^\beta(\wt\cM_{\hb_o})$ is such that each graded piece $\gr_\ell^{\rM(\rN)}\psi_{t}^\beta(\wt\cM_{\hb_o})$ is $\cO_{\Omega_0,\hb_o}$-locally free for any $\ell$ and any $\beta$.
\end{enumerate}
\end{enumerate}
\end{Corollaire}

For the proof, we will need Lemma \ref{lem:killing} below, analogous to \cite[p\ptbl 79]{Biquard97}. Let us recall some notation: we put $\ccD'_\hb=\hb D'_E+\theta'_E$ and $\ccD''_\hb=D''_E+\hb\theta''_E$. By Theorem~\ref{th:biquard}, we write $\ccD'_\hb$ (\resp $\ccD''_\hb$) as the sum of $\ccD^{\prime\std}_\hb$ (\resp $\ccD^{\prime\prime\std}_\hb$) plus a perturbation which is controlled. We denote with the same letter the action of these connections on endomorphisms of $H$ with values in differential forms.

\begin{lemme}[Local killing of the perturbation $P''$]\label{lem:killing}
In the situation of Theorem \ref{th:biquard}, for any $\hb_o\in\Omega_0$, there exists $\eta>0$ and a matrix $Q^{(\hb_o)}(t,\hb)$ of functions $X^*\to\rH(\Delta_{\hb_o}(\eta))$ such that
\begin{enumerate}
\item\label{lem:killing1}
the $\rH(\Delta_{\hb_o}(\eta))$-norms of $\Lt^\delta Q^{(\hb_o)}(t,\hb)$, $\Lt^\delta\ccD''_{\hb} Q^{(\hb_o)}(t,\hb)$ and $\Lt^\delta\ccD'_{\hb} Q^{(\hb_o)}(t,\hb)$ remain bounded when $t\to0$, for some $\delta>0$;
\item\label{lem:killing2}
in the basis $\varepsilong^{\prime(\hb_o)}(\hb)\defin\varepsilong\cdot(\id+Q^{(\hb_o)}(t,\hb))$, the matrix of $\ccD''_{\hb}$ is the standard one, namely $\big[M^{\prime\prime\std}+(\hb-1)(M^{\prime\prime\std}+M^{\prime\std*})/2\big]d\ov t/\ov t$.
\end{enumerate}
\end{lemme}

\begin{proof}[Sketch of proof\,\footnotemark]\footnotetext{I thank O\ptbl Biquard for explaining me his proof and the referee for noticing a missing point in a previous proof.}

It is a variant of \loccit The basis $\varepsilong$ of Theorem \ref{th:biquard} is decomposed in subfamilies $\varepsilong_\beta$ for $\beta\in B$. Given any matrix $P$, we denote by $P_{\beta_i,\beta_j}$ its $(\beta_i,\beta_j)$-block. We will first use the following property of $P=P'$ or $P''$ given in \loccit: there exists $\delta>0$ such that
\begin{enumerate}
\item[(i)]
for any $\beta\in B$, the matrix function $\Lt^{1+\delta}P_{\beta,\beta}(t)$ is bounded on $X^*$.
\item[(ii)]
if $\beta_i\neq\beta_j\in B$, the matrix function $\Lt^{2+\delta}P_{\beta_i,\beta_j}(t)$ is bounded on $X^*$.
\end{enumerate}
(The second condition corresponds to the case denoted by $\ell\neq0$ in \loccit).

We will also need a better estimate, the proof of which will be indicated in \T\ref{app:biquard1}: there exists $\epsilon>0$ such that
\begin{enumerate}
\item[(iii)]
if $\beta''_i\neq\beta''_j$, the matrix function $\mt^{-\epsilon} P_{\beta_i,\beta_j}(t)$ is bounded on $X^*$.
\end{enumerate}
Let us now fix $\hb_o\in\Omega_o$. It will be convenient to work with the basis $\wt\varepsilong(\hb)=\varepsilong\cdot e^{-\hb\rX}$. In this basis, the matrix of $\ccD''_{\hb}=D''_E+\hb\theta''_E$ can be written as $[\Std''(\hb)+b''(t,\hb)]d\ov t/\ov t$ where
\begin{align*}
\Std''(\hb)&=\oplus_\beta\Big(\big[-(\beta'+i\hb\beta'')/2\big]\id+\rH/2\Lt\Big),\\
b''(t,\hb)&=e^{\hb\rX}\big[P''+(\hb-1)(P''+P^{\prime*})/2\big]e^{-\hb\rX}.
\end{align*}
We look for $\id+u(t,\hb)$ such that
\begin{equation}\tag*{(\protect\ref{lem:killing})$(*)$}\label{eq:uStd}
\ov t\,\frac{\partial u}{\partial\ov t}=-[\Std''(\hb),u(t,\hb)]-b''(t,\hb)\big(\id+u(t,\hb)\big).
\end{equation}
Then $\id+Q^{(\hb_o)}(t,\hb)=e^{-\hb\rX}\big(\id+u(t,\hb)\big)e^{\hb\rX}$ will be a solution to \ref{lem:killing}\eqref{lem:killing2}.

We can now argue as in \cite[pp\ptbl78--79]{Biquard97}. As both operators $\ad\Std''(\hb)$ and $\ad \rH$ commute, we can decompose any $d\times d$ matrix $u$ as $u=\oplus_{\gamma,\ell}u^{(\gamma,\ell)}$, where $u^{(\gamma,\ell)}$ belongs to the $\gamma(\hb)$-eigenspace of $\ad\Std''(\hb)$ and the $\ell$-eigenspace of $\ad\rH$. Denote by $E_{\delta,\epsilon}\big(D^*_R,\rH(\Delta_{\hb_o}(\eta))\big)$ the Banach space of matrices $u$ of functions $D^*_R\to\rH(\Delta_{\hb_o}(\eta))$ on the punctured disc of radius $R$ such that the entries of $u^{\gamma,\ell}$ satisfy $\big\Vert\Lt^\delta u^{\gamma,\ell}(t)\big\Vert_{\rH(\Delta_{\hb_o}(\eta))}$ bounded on $D^*_R$ and, if $\gamma(\hb)$ is not constant, $\big\Vert\mt^{-\epsilon} u^{\gamma,\ell}\big\Vert_{\rH(\Delta_{\hb_o}(\eta))}$ bounded.

For a matrix $v^{(\gamma,\ell)}$ in the $(\gamma,\ell)$ eigenspace of $(\ad\Std''(\hb),\ad\rH)$ and for $n\in\ZZ$, put
\[
T_n^{(\gamma,\ell)}v^{(\gamma,\ell)}=t^{-n}\mt^{-\reel\gamma}\Lt^{\ell/2}\int_{D_R^*}\frac{w^n}{\ov w}\,\module{w}^{\reel\gamma}\rL(w)^{-\ell/2}\cdot v^{(\gamma,\ell)}\itwopi\frac{dw\wedge d\ov w}{t-w}.
\]

Notice that $\gamma=\gamma(\hb)$ is a function of $\hb$ which takes the form $(\beta'_1-\beta'_2)+i\hb(\beta''_1-\beta''_2)$, for $\beta_1,\beta_2$ in $B$. Then, either $\hb\to\gamma(\hb)$ is not constant, \ie $\beta''_1\neq\beta''_2$, or it is constant and belongs to $]-1,1[$. We define $n_{\gamma(\hb_o),\ell}\in\ZZ$ by the property
\[
n_{\gamma(\hb_o),\ell}+\reel\gamma(\hb_o)
\begin{cases}
\in[-1,0[&\text{if }\gamma(\hb)\not\equiv0,\\
=-1&\text{if }\gamma(\hb)\equiv0\text{ and }\ell\geq-1,\\
=0&\text{if }\gamma(\hb)\equiv0\text{ and }\ell\leq-2.
\end{cases}
\]
The radius $\eta>0$ is chosen such that, for any $\gamma$ with $\gamma(\cbbullet)$ not constant and $\gamma(\hb_o)\in\ZZ$, then $\reel\gamma(\hb)+\epsilon\in{}]\reel\gamma(\hb_o),\reel\gamma(\hb_o)+1[$. Consider the operator
\[
\wt\cT:u\mto\oplus_{\gamma,\ell}T_{n_{\gamma,\ell}}^{(\gamma,\ell)}\big[(b''u+b'')^{(\gamma,\ell)}\big].
\]
We then obtain as in \loccit that, if $R$ is small enough, the operator $\wt\cT$ sends $E_{\delta,\epsilon}\big(D^*_R,\rH(\Delta_{\hb_o}(\eta))\big)$ into itself, and is contracting. The fixed point $u$ is solution of \ref{eq:uStd}, where the derivative is taken in the distributional sense on $D^*_R$.

Notice that the only assumption of a logarithmic decay of $u^{(\gamma,\ell)}$ would cause trouble if $\gamma(\hb_o)\in\ZZ$ and $\gamma(\hb)$ non constant. This is ruled out by the stronger decay in $\mt^\epsilon$ (Property (iii)) in such a case.

The first two properties of \ref{lem:killing}\eqref{lem:killing1} directly follow from the construction, as $\ccD''_\hb(\id+u)=-(\id+u)b''$.

The matrix of $\ccD'_\hb$ in the basis $\wt\varepsilong(\hb)$ can be written as $[\Std'(\hb)+b'(t,\hb)]dt/t$, with $\Std'(\hb)=\oplus_\beta\big[(\hb\beta'+i\beta'')/2\id+(\rY+\hb\rH/2)/\Lt\big]$ and $b'$ satisfying (i)--(iii). Notice that, for $u\in E_{\delta,\epsilon}\big(D^*_R,\rH_{\hb_o}(\Delta(\eta))\big)$, we have $\big\Vert\Lt^{1+\delta}[\Std'(\hb),u]\big\Vert_{\rH_{\hb_o}(\Delta(\eta))}$ bounded (and a similar property for $[\Std''(\hb),u]$). In order to obtain the desired estimate for $\ccD'_\hb Q^{(\hb_o)}$, it is therefore enough to show that $\Lt\, d'u(t,\hb)$ remains bounded when $t\to0$. This is obtained by an argument of elliptic regularity which closely follows that of \cite{Biquard97}. Details are given in \T\ref{app:biquard2}.
\end{proof}

\begin{Proof}[Proof of Corollary \ref{cor:biquard}]
\subsubsection*{Identification of the parabolic filtration}
Fix $\hb_o\in\Omega_0$ and work in a neighbourhood of $\hb_o$ as in the previous lemma. We denote by $Q^{(\hb_o)}$ the corresponding matrix. For $\beta\in B$, let $q_{\beta,\imhb_o}\in\ZZ$ be such that $\ell_{\hb_o}(q_{\beta,\imhb_o}+\beta)\defin q_{\beta,\imhb_o}+\beta'-\imhb_o\beta''$ belongs to $[0,1[$.

The basis $\varepsilong$ of Theorem \ref{th:biquard} is decomposed into subfamilies $\varepsilong_\beta$ for $\beta\in B$, so that $M^{\prime\std},M^{\prime\prime\std}$ are block-diagonal matrices. Recall that we set $A_{\beta}(t,\hb)\index{$abthb$@$A_\beta(t,\hb)$}=e^{-\hb X}\mt^{\beta'+i\hb\beta''}\Lt^{\rHsd}$ and we put $A(t,\hb)=\oplus_{\beta\in B}A_\beta(t,\hb)$. We also define similarly the diagonal matrices $\wt A_{\beta}(t,\hb)=\mt^{\beta'+i\hb\beta''}\Lt^{\rHsd}$ and $\wt A(t,\hb)=\oplus_{\beta\in B}\wt A_\beta(t,\hb)$. Last, put
\begin{equation}\label{eq:bme}
\index{$ehb$@$\bme^{(\hb_o)}$}\bme^{(\hb_o)}=\varepsilong\cdot(\id+Q^{(\hb_o)}(t))A(t,\hb)= \varepsilong^{\prime\prime(\hb_o)}\cdot \wt A(t,\hb)
\end{equation}
(defining therefore also $\varepsilong^{\prime\prime(\hb_o)}$), where $Q^{(\hb_o)}(t)$ is given by Lemma \ref{lem:killing}. For any $j$, we have $e^{(\hb_o)}_j=\mt^{\beta'_j+i\hb\beta''_j}\Lt^{w_j/2}\varepsilon^{\prime\prime(\hb_o)}_j$ for some $\beta_j\in B$ and $w_j\in\ZZ$.

According to \ref{lem:killing}\eqref{lem:killing2}, the computation following \eqref{eq:PRh} shows that the basis $\bme^{\prime(\hb_o)}$, defined by \begin{equation}\label{eq:bmeprime}
\index{$ehbp$@$\bme^{\prime(\hb_o)}$}e^{\prime(\hb_o)}_j=t^{q_{\beta_j,\imhb_o}}e^{(\hb_o)}_j,
\end{equation}
is holomorphic. As the base changes $\varepsilong\mto\bme$ and $\bme\mto\varepsilong$ have moderate growth along $\{0\}\times\Omega_0$, we conclude that, near $\hb_o$, the sheaf $\wt\cM$ is nothing but the $\cO_\cX[1/t]$-locally free sheaf generated by $\bme^{(\hb_o)}$ in $j_*\cH'$. This gives the first part of \eqref{cor:biquard1}. The strict specializability follows from the proof of \eqref{cor:biquard2} and \eqref{cor:biquard3} below.

\medskip
Notice that the basis $\bme^{(\hb_o)}$ may also be decomposed into subfamilies $$\bme^{(\hb_o)}_\beta=\{e^{(\hb_o)}_j\mid\beta_j=\beta\}\quad (\beta\in B)$$
so that
\begin{equation}\label{eq:bmebeta}
\bme^{(\hb_o)}_\beta=\Big[\varepsilong_\beta\cdot(\id+Q^{(\hb_o)}_{\beta,\beta})+\sum_{\beta'\neq\beta}\varepsilong_{\beta'}\cdot Q^{(\hb_o)}_{\beta',\beta}\Big]A_{\beta}.
\end{equation}
For $b\in\RR$, denote by $U^b_{(\hb_o)}\wt\cM$ the locally free $\cO_\cG$-module generated by the sections $t^{n_j}e^{(\hb_o)}_j$ with $n_j\in\ZZ$ such that $\ell_{\hb_o}(n_j+\beta_j)\in[b,b+1[$. Notice that, for any $k\in\ZZ$, we have $U^{b+k}_{(\hb_o)}\wt\cM=t^kU^b_{(\hb_o)}\wt\cM$.

We will now show that $U^{\cbbullet}_{(\hb_o)}\wt\cM$ induces on $\wt\cM_{\hb_o}$ the parabolic filtration (defined in the statement of Corollary \ref{cor:biquard}).

Formula \eqref{eq:bme} shows that each element $t^ke_j^{(\hb_o)}$ of $\bme^{(\hb_o)}$ has parabolic order equal to $\ell_{\hb_o}(k+\beta_j)$ exactly, according to the logarithmic decay of $Q^{(\hb_o)}(t)$ given by Lemma \ref{lem:killing}.

Conversely, let $m=\sum_jm_jt^{n_j}e^{(\hb_o)}_j\in\wt\cM_{\hb_o}$ with $m_j$ holomorphic and $m_j\not\equiv0\implique m_j(0,\hb)\not\equiv0$. We may assume that $m_j\not\equiv0\implique m_j(0,\hb)\neq0$ for $\hb\neq\hb_o$ in a neighbourhood of $\hb_o$. Put $b=\min_{j\mid m_j\neq0}\ell_{\hb_o}(n_j+\beta_j)$. We will prove that $m$ has order~$b$. Clearly, $m$ has order $\geq b$. It is then enough to prove that, for any $\epsilon>0$, $\mt^{-b-\epsilon}\norme{m}$ is not bounded in any neighbourhood of $(0,\hb_o)$. There exists only one $k\in\ZZ$ and $\beta\in B$ such that, for any $\imhb\in[\imhb_o,\imhb_o+\eta]$ with $\eta$ small enough, $\ell_{\hb}(k+\beta)$ achieves the minimum of $\{\ell_{\hb}(n_j+\beta_j)\mid m_j\neq0\}$. Notice that $\ell_{\hb}(k+\beta)\leq b$ if $\imhb\in[\imhb_o,\imhb_o+\eta]$ and that, unless $\beta$ is real, the inequality is strict if $\imhb\neq\imhb_o$ (see Fig\ptbl \ref{fig:lines} on page\kern2pt\pageref{fig:lines}). Denote by $J_1$ the set of $j$ for which this minimum is achieved. Put $w=\max_{j\in J_1}w_j$ and denote by $J$ the set of $j\in J_1$ for which this maximum is achieved. Let $\bmm_J(0,\hb)$ be the vector having entries $m_j(0,\hb)$ for $j\in J$ and $0$ otherwise and put $\wt\bmm_J(0,\hb)=e^{-\hb X}\bmm_J(0,\hb)$. Then an easy computation shows that
\begin{equation}\label{eq:norme}
\norme{m}_{\pi^*h}^2\underset{\substack{t\to0\\ \imhb\in[\imhb_o,\imhb_o+\eta]}}{\sim} \mt^{2\ell_\hb(k+\beta)}\Lt^w\norme{\wt\bmm_J(0,\hb)}^2.
\end{equation}
This gives the assertion.

\subsubsection*{Computation of $\Theta'_\hb$}
We now compute the matrix $\Theta'_\hb$ of $\ccD'_\hb$ in the basis $\bme^{\prime(\hb_o)}$. Notice first that, by flatness, the matrix of $\ccD'_\hb$ in the basis $\bme^{\prime(\hb_o)}$ takes the form
\begin{equation}\label{eq:matrix}
\Theta'_{\hb}=\Big[\big(\oplus_\beta\big[(q_{\beta,\imhb_o}+\beta)\star\hb\id+\rY_\beta\big]\big)+P(t,\hb)\Big]\frac{dt}{t},
\end{equation}
with $P(t,\hb)$ holomorphic on $\{t\neq0\}\times\nb(\hb_o)$.

\begin{enonce*}{Assertion}
The matrix $P(t,\hb)$ satisfies:
\begin{itemize}
\item
if $\beta_j\neq\beta_k$,
\begin{itemize}
\item
if $\ell_{\hb_o}(q_{\beta_j,\imhb_o}+\beta_j)\leq \ell_{\hb_o}(q_{\beta_k,\imhb_o}+\beta_k)$, then $P_{\beta_j,\beta_k}/t$ is holomorphic near $(0,\hb_o)$;
\item
if $\ell_{\hb_o}(q_{\beta_j,\imhb_o}+\beta_j)> \ell_{\hb_o}(q_{\beta_k,\imhb_o}+\beta_k)$, then $P_{\beta_j,\beta_k}$ is holomorphic near $(0,\hb_o)$;
\end{itemize}
\item
if $\beta_j=\beta_k$,
\begin{itemize}
\item
if $w_j\geq w_k-2$, then $P_{\beta_j,\beta_k}/t$ is holomorphic near $(0,\hb_o)$;
\item
if $w_j\leq w_k-3$, then $P_{\beta_j,\beta_k}$ is holomorphic near $(0,\hb_o)$.
\end{itemize}
\end{itemize}
\end{enonce*}

\begin{proof}[Proof of the assertion]
We will prove that the matrix $P(t,\hb)$ satisfies the following property, for some $\delta>0$ small enough and $\hb\in\nb(\hb_o)$:
\[
\lim_{t\to0}\Lt^\delta\mt^{(q_{\beta_j,\imhb_o}+\beta'_j-\imhb\beta''_j)- (q_{\beta_k,\imhb_o}+\beta'_k-\imhb\beta''_k)}\Lt^{1+(w_j-w_k)/2}P_{\beta_j,\beta_k}(t,\hb)=0,
\]
from which the assertion follows. Using notation of Lemma \ref{lem:killing}, this, in~turn, is equivalent to the fact that the matrix of $\ccD'_\hb$ in the basis $\wt\varepsilong(\hb)(\id+u)$ takes the form $[\Std'(\hb)+c'(t,\hb)]dt/t$ with $\big\Vert\Lt^\delta c'(t,\hb)dt/t\big\Vert_{\rH_{\hb_o}(\Delta(\eta))}$ bounded. We know, after \cite{Biquard97}, that the matrix of $\ccD'_\hb$ in the basis $\wt\varepsilong(\hb)$ takes the form \hbox{$[\Std'(\hb)+b'(t,\hb)]dt/t$}, where $b'dt/t$ satisfies the previous estimate and Lemma \ref{lem:killing} implies that the functions  $t\mto\big\Vert\Lt^\delta u(t,\hb)\big\Vert_{\rH_{\hb_o}(\Delta(\eta))}$ and $t\mto\big\Vert\Lt^\delta\ccD'_\hb u(t,\hb)\big\Vert_{\rH_{\hb_o}(\Delta(\eta))}$ are bounded. On the other hand, we have $c'(t,\hb)dt/t=(\id+u)^{-1}\ccD'_\hb u+b'dt/t$. Both terms in the right-hand side satisfy the desired property.
\end{proof}

\subsubsection*{End of the proof}
The assertion above shows that $U^b_{(\hb_o)}\wt\cM$ is stable under $t\partiall_t$. Moreover, the matrix of $t\partiall_t$ on $U_{(\hb_o)}^0/U^1_{(\hb_o)}$ is block-lowertriangular, if we order the families $\bme^{\prime(\hb_o)}_\beta$ with respect to the $\ell_{\hb_o}$-order, each block corresponding to a value of $\ell_{\hb_o}(q_{\beta,\imhb_o}+\beta)$. Each diagonal block is itself block-diagonal with respect to the various $\beta$ with $\ell_{\hb_o}(q_{\beta,\imhb_o}+\beta)$ fixed. This shows that each graded piece $\gr^b_{U_{(\hb_o)}}\wt\cM$ is $\cO_{\Omega}$-locally free near $\hb_o$, hence $U_{(\hb_o)}^\cbbullet\wt\cM$ satisfies the properties of Lemma \ref{lem:uniciteMK}. Therefore, $\wt\cM$ is strictly specializable along $t=0$ and $U_{(\hb_o)}^\cbbullet\wt\cM=V_{(\hb_o)}^\cbbullet\wt\cM$.

Last, let $W_\bbullet\gr^b_{V_{(\hb_o)}}\wt\cM$ be the weight filtration of $\oplus_\beta \rY_\beta$. It also follows from the previous properties that it is decomposed with respect to $\oplus\psi_t^{q_{\beta,\imhb_o}+\beta}\wt\cM$ and each summand $W_\bbullet\psi_t^{q_{\beta,\imhb_o}+\beta}\wt\cM$ satisfies the properties characterizing the monodromy filtration of $\rN=t\partiall_t-(q_{\beta,\imhb_o}+\beta)\star\hb$. It is therefore equal to it and $\rN:\psi_t^{q_{\beta,\imhb_o}+\beta}\wt\cM\to\psi_t^{q_{\beta,\imhb_o}+\beta}\wt\cM$ is conjugate to $\rY_\beta$. This proves \ref{cor:biquard}\eqref{cor:biquard3}.
\end{Proof}

\begin{Remarques}\label{rem:although}
\begin{enumerate}
\item
\label{rem:although2}
Although the bases $\bme^{(\hb_o)}$ are defined only locally with respect to $\hb$, the classes of their elements in the various bundles $\gr_\ell^{\rM}\psi_t^\beta\wt\cM$ are globally defined on $\Omega_0$. Notice that these bundles are holomorphically trivial on $\Omega_0$. More precisely, for any $\ell\in\ZZ$ and any $\beta\in B$, there exists a basis $\bme^o_{\beta,\ell}$ of $\gr_\ell^{\rM}\psi_t^\beta\wt\cM$ as a $\cO_{\Omega_0}$-module, uniquely determined from the basis $\varepsilong$ of Theorem \ref{th:biquard}, such that $(\bme^o_{\beta,\ell})_{\beta\in B,\,\ell\in\ZZ}$ lifts locally to bases $\bme^{(\hb_o)}$.

In order to prove this statement, it is enough to characterize the inverse image in $V_{(\hb_o)}^{\ell_{\hb_o}(\beta)}\wt\cM$ of $\gr_\ell^{\rM}\psi_t^\beta\wt\cM$ near $\hb_o$.

First, a section $m$ of $V_{(\hb_o)}^{\ell_{\hb_o}(\beta)}\wt\cM$ near $(0,\hb_o)$ has a class in $\gr^{\ell_{\hb_o}(\beta)}_{V_{(\hb_o)}}\wt\cM$ contained in $\psi_t^\beta\wt\cM$ if and only if, for any $j$, denoting by $n_j\in\ZZ\cup\{+\infty\}$ the order of its coefficient $m_j$ on $e_j^{(\hb_o)}$ along $t=0$, we have
$\ell_{\hb_o}(n_j+\beta_j)>\ell_{\hb_o}(\beta)$ if $\beta_j\neq\beta$ and at least one $m_j(0,\hb)\not\equiv0$ with $\beta_j=\beta$.

When this is satisfied, the class of $m$ has $\rM$-order $\ell$ if and only if, for any $j$ with $\beta_j=\beta$, $m_j(0,\hb)\equiv0$ if $w_j>\ell$ and $m_j(0,\hb)\not\equiv0$ for at least one $j$ with $\beta_j=\beta$ and $w_j=\ell$.

Both assertions are proved as in the proof of Corollary \ref{cor:biquard}. This criterion may be translated in terms of order of growth: the class of $m$ is in $\psi_t^\beta\wt\cM$ iff for any $\hb'_o\neq\hb_o$ and sufficiently near $\hb_o$, the norm $\norme{m_{|\hb=\hb'_o}}_{h}$ has growth order $\ell_{\hb'_o}(\beta)$, and this class has $\rM$-order $\ell$ iff moreover this norm grows as $\mt^{\ell_{\hb'_o}(\beta)}\Lt^{\ell/2}$.

With this criterion in mind, one observes that two meromorphic bases given by Formula \eqref{eq:bme} for different choices of $Q^{(\hb_o)}$ give rise to the same classes in the various $\gr^{\rM}_\ell\psi_t^\beta\wt\cM$, so these classes can be glued when $\hb_o$ varies in $\Omega_0$. By definition, the basis $\bme^o$ is adapted to the Lefschetz decomposition of $\gr_{-2}^{\rM}\rY_\beta$, hence to that of $\gr_{-2}^{\rM}\rN$. In particular, the whole basis of $\gr^{\rM}\psi_t^\beta\wt\cM$ can be recovered from its primitive vectors.

\item\label{rem:although3}
We also have an asymptotic estimate of $\norme{m}_{\pi^*h}$ as in Formula \eqref{eq:norme} for $\imhb\in[\imhb_o-\eta,\imhb_o]$. It uses a maybe distinct pair $(\beta_-,k_-)$, subject however to the relation $\ell_{\hb_o}(k_-+\beta_-)=\ell_{\hb_o}(k+\beta)$.

By restriction to $\hb=\hb_o$, we also get a criterion for a nonzero local section $m^o$ of $\wt M_{\hb_o}$: Denote by $U^\cbbullet M_{\hb_o}$ the filtration induced by $V^\cbbullet_{(\hb_o)}\wt\cM$ (\cf Proposition \ref{prop:restrhbo}). Denote by $\rN$ the corresponding nilpotent endomorphism on $\gr^b_U\wt M_{\hb_o}$, which is the direct sum of the various $-(\partiall_tt+\beta\star\hb_o)$ for $\beta$ such that $\ell_{\hb_o}(\beta)=b$, and let $\rM_\bbullet\gr^b_U\wt M_{\hb_o}$ be its monodromy filtration, that we lift as $\rM_\bbullet U_b\wt M_{\hb_o}$. Restrict the basis $\bme^{(\hb_o)}$ to a basis $\bme^o$ of $\wt M_{\hb_o}$, and write $m^o=\sum_jm^o_j(t)t^{n_j}e_j^o$. We assume that \hbox{$m^o_j\!\not\equiv\!0\implique m_j^o(0)\!\neq\!0$}. Put $b=\min_{j\mid m_j\neq0}\ell_{\hb_o}(n_j+\beta_j)$. We denote by $J^o_1$ the set of $j$ such that the minimum is achieved, we put $w=\max_{j\in J^o_1}w_j$ and we denote by $J^o$ the set of $j\in J_1^o$ such that the maximum is achieved. We then have the restriction of Formula \eqref{eq:norme}:
\begin{equation}\tag*{(\protect\ref{eq:norme})$_{\hb_o}$}\label{eq:normehbo}
\norme{m^o}_h^2\underset{t\to0}{\sim} \mt^{2b}\Lt^w\norme{\wt\bmm^o_J(0)}^2.
\end{equation}
Moreover, $m^o$ is in $U^b\wt M_{\hb_o}$ and its class in $\gr_w^{\rM}\gr^b_UM_{\hb_o}$ is nonzero.

\item\label{rem:although4}
For any $\hb_o\in\Omega_0$, it is possible to find a matrix $Q(\hb)$ which is holomorphic with respect to $\hb$ and is invertible for any $\hb$ near $\hb_o$, so that, after changing the basis $\bme^{(\hb_o)}$ by the matrix $Q(\hb)$, the residue at $t=0$ of the matrix of $\hb D_E+\theta'_E$ is block-diagonal, each block corresponding to eigenvalues $\beta\star\hb$ taking the same value at $\hb_o$, and is lower-triangular, each diagonal sub-block corresponding to a given $\beta$.

Assume moreover that $\hb_o\not\in\Sing\Lambda$, in particular $\hb_o\neq0$. Then, there exists a neighbourhood of $\hb_o$ such that, for any $\hb$ in this neighbourhood and any $\beta_1,\beta_2\in B$, the differences $(\beta_1-\beta_2)\star\hb/\hb\in\ZZ\iff\beta_1=\beta_2$.

One may therefore apply the same arguments as in the theory of regular meromorphic connections of one variable with a parameter to find, near any $(0,\hb_o)$ with $\hb_o\not\in\Sing\Lambda$, a basis $\wt\bme^{(\hb_o)}=\bme^{(\hb_o)}\cdot \ccP(t,\hb)$ of $\wt\cM$, with $\ccP(0,\hb_o)$ invertible, in which the matrix of $\hb D_E+\theta'_E$ is equal to
\[
\ooplus_\beta[(\beta\star\hb)\id+\rY_\beta]\frac{dt}{t}.
\]

If now $\hb_o\in\Sing\Lambda\moins\{0\}$, we first apply a classical ``shearing transformation'', which is composed by successive rescalings by powers of $t$ and invertible matrices depending holomorphically on $\hb$ for $\hb$ near $\hb_o$, so that, for any two eigenvalues $\beta_1\star\hb,\beta_2\star\hb$ of the ``constant'' part (\ie depending on $\hb$ only) of the matrix of $\hb D_E+\theta'_E$, the difference is not a nonzero integer. Then it is possible to find a base change $\ccP(t,\hb)$ as above, such that the new matrix of $\hb D_E+\theta'_E$ does not depend on $t$, is holomorphic with respect to $\hb$ and is lower triangular with eigenvalues $\beta\star\hb$ for the new set of $\beta$'s.

Notice that in both cases, $\wt\cM$ is, locally with respect to $\hb$, an extension of rank one $\cO_{\cX}[1/t]$-modules with connection.
\end{enumerate}
\end{Remarques}

\Subsection{Construction of the $\cR_\cX$-module $\cM$ and the filtration $V_\bbullet\cM$}\label{subsec:exemplestricspe}

Let us return to the global setting on the Riemann surface $X$. The $\cR_X$-module $\cM$ is defined as the minimal extension of $\wt\cM$ across its singular set $P$ (\cf Def\ptbl\ref{def:minext}).

\begin{corollaire}\label{cor:biquardRX}
$\cM$ is a regular holonomic $\cR_\cX$-module which is strictly specializable along $t=0$ and has strict support equal to $X$.
\end{corollaire}

\begin{proof}
Let us show that $\cM$ is good. On $\cX\moins \cP$, $\cM=\wt\cM$ is $\cO_\cX$-coherent. Let $x\in P$. If $K$ is any compact set in $\{x\}\times\Omega_0$, then $\cM_{|K}$ is generated by the sheaf of local sections $m$ which satisfy $\norme{m}_{\pi^*h}\leq C\mt^{-N}$ for $N=N_K$ large enough, as this sheaf contains $V_{\hb_o}^{>-1}\wt\cM_{\hb_o}$ for any $\hb_o\in K$. This sheaf is $\cO_\cX$-coherent, as follows from Corollary \ref{cor:biquard}, hence $\cM$ is good.

By definition, on $\Delta_{\hb_o}(\eta)$, the $V$-filtration $V_{(\hb_o)}^{\cbbullet}\cM$, restricted to indices in $\ZZ$, is good; hence $\cM$ is regular along $\{x\}\times \Omega_0$ for any $x\in P$, as each $V_{(\hb_o)}^{\cbbullet}\cM$ is $\cO_\cX$-coherent.

Last, as $\cM$ is $\cO$-coherent on $\cX\moins \cP$, its characteristic variety is the zero section on $\cX\moins \cP$, hence is contained in $(T^*_XX\cup T^*_PX)\times\Omega_0$; in other words, $\cM$ is holonomic.

The strict specializability along $t=0$ follows from Proposition \ref{prop:minext} and the fact that $\cM$ is a minimal extension implies that it has strict support equal to $X$.
\end{proof}

We will end this paragraph by proving:

\begin{lemme}
The restriction to $\hb=1$ of the $\cR_\cX$-module $\cM$ is equal to $M$.
\end{lemme}

\begin{proof}
Consider first the restriction of $\wt\cM$ to $\hb=1$. We will show that is is equal to $\wt M$ (with connection). On the open set $X^*=X\moins P$, this was shown in the proof of Lemma \ref{lem:twharm}. We know by \ref{cor:biquard}\eqref{cor:biquard1} that $\wt\cM/(\hb-1)\wt\cM$ is $\cO_X(*P)$-locally free. It is therefore enough to prove that $\wt\cM/(\hb-1)\wt\cM$ and $\wt M$ define the same meromorphic extension, or equivalently that $\wt\cM/(\hb-1)\wt\cM\subset\wt M$. By the construction of the metric~$h$, $\wt M$ is the subsheaf of $j_*V$ (where $j:X^*\hto X$ is the open inclusion) of sections, the $h$-norm of which has moderate growth. Now, Formula \eqref{eq:bme} computed at $\hb_o=1$ shows that the $h$-norm of the basis $\bme^{(1)}_{|\hb=1}$ has moderate growth, implying therefore the required inclusion.

Notice that the same formula shows that the filtration $V^{\cbbullet}_{(\hb_o=1)}\wt\cM$ restricts to the canonical filtration $V^{\cbbullet} M$ of $\wt M$. Recall now that the minimal extension $M\subset\wt M$ is the sub-$\cD_X$-module of $\wt M$ generated by $V^{>-1}\wt M$. By construction, we therefore have $\cM/(\hb-1)\cM\subset M$. As $M$ is simple, both $\cD_X$-modules coincide.
\end{proof}

\Subsection{Construction of the sesquilinear pairing $C$}\label{subsec:constrC}

We have to define a sesquilinear pairing $C$ on $\cMS\otimes_{\cO_\bS}\ov{\cMS}$ with values in $\Dbh{X}$, which has to extend $h_{\bS}$ defined on $X^*$. We will construct $C$ with values in the sheaf $\Db^\an_{\cX|\bS}$ of distributions which are holomorphic with respect to $\hb$.

We first define $C$ locally with respect to $\bS$ and then show that it glues along $\bS$. So we fix $\hb_o\in\bS$ and a compact disc $\Deltag$ neighbourhood of $\hb_o$, that we assume to be small enough. We denote by $\Deltag^{\!\circ}$ its interior. In the following, we will denote by $\rH(\Deltag)$ the Banach space of continuous functions on $\Deltag$ which are holomorphic on its interior.

Consider the basis $\bme^{\prime(\hb_o)}$ of $V_{(\hb_o)}^{>-1}\cM$ defined near $\hb_o$ by a formula analogous to \eqref{eq:bme}: we have $e_j^{\prime(\hb_o)}=t^{n_j}\mt^{\beta'_j+i\hb\beta''_j}\Lt^{w_j/2}\varepsilon_j^{\prime\prime(\hb_o)}$, with $\ell_{\hb_o}(n_j+\beta_j)\in{}]-1,0]$. Consider similarly the basis $\bme^{\prime(-\hb_o)}$ of $V_{(-\hb_o)}^{>-1}\cM$ defined near $-\hb_o$, with $\ell_{-\hb_o}(\nu_k+\beta_k)\in{}]-1,0]$. Notice that, as $\module{\hb_o}=1$, we have $\sigma(\hb_o)=-\hb_o$ and $\im\sigma(\hb_o)=-\imhb_o$, using the notation of \T\ref{subsubsec:conj}.

The entries of the matrix of $\hdel$ in these bases may be written as
\begin{equation}\label{eq:estimate}
t^{n_j}\mt^{\beta_j'+i\hb\beta_j''}\cdot \ov t^{\nu_k}\mt^{\beta_k'+i\beta_k''/\hb}\cdot \Lt^{(w_j+w_k)/2}\cdot a_{j,k}(t),
\end{equation}
with $\norme{a_{j,k}(t)}_{\rH(\Deltag)}\to0$ if $\beta_j\neq\beta_k$ and $\norme{a_{j,k}(t)}_{\rH(\Deltag)}$ locally bounded, when $t\to0$. As we have $\ell_{\hb_o}(n_j+\beta_j)>-1$ and $\ell_{-\hb_o}(\nu_k+\beta_k)>-1$, we also have $\reel(n_j+\beta_j'+i\hb\beta_j'')>-1$ and $\reel(\nu_k+\beta_k'+i\beta_k''/\hb)>-1$ for any $\hb\in\Deltag$, if $\Deltag$ is small enough. Therefore, the entries of the matrix of $\hdel$ are in $L^1_{\loc}(X,\rH(\Deltag))$, hence, as $V_{(\hb_o)}^{>-1}\cM$ is $\cO_\cX$-locally free, $\hdel$ defines a sesquilinear pairing
\[
\hdel:V_{(\hb_o)}^{>-1}\cM_{|\Deltag}\otimes_{\cO_{\Deltag}}\ov{V_{(-\hb_o)}^{>-1}\cM_{|\sigma(\Deltag)}}\to L^1_{\loc}(X,\rH(\Deltag))
\]
by $\cO_\cX\otimes\cO_{\ov \cX}$-linearity.

As $L^1_{\loc}(X,\rH(\Deltag))$ is contained in the subspace $\ker\ov\partial_\hb$ of $L^1_{\loc}(X\times\Deltag^{\!\circ})$, this defines similarly $\Cdelc$ on $V_{(\hb_o)}^{>-1}\cM_{|\Deltag^{\!\circ}}\otimes_{\cO_{\Deltag^{\!\circ}}}\ov{V_{(-\hb_o)}^{>-1}\cM_{|\sigma(\Deltag^{\!\circ})}}$ by $\cO_\cX\otimes\cO_{\ov \cX}$-linearity.

\begin{Assertion}\label{ass:1}
$\hdel$ is $V_0\cR_{\cX|\Deltag}\otimes_{\cO_{\Deltag}}\ov{V_0\cR_{\cX|\sigma(\Deltag)}}$-linear and $\Cdelc$ is $V_0\cR_{\cX|\Deltag^{\!\circ}}\otimes_{\cO_{\Deltag^{\!\circ}}}\ov{V_0\cR_{\cX|\sigma(\Deltag^{\!\circ})}}$-linear.
\end{Assertion}

\begin{Assertion}\label{ass:2}
$\Cdelc$ can be extended by $\cR_{(X,\ov X),\Deltag^{\!\circ}}$-linearity as a well-defined sesquilinear pairing on $\cM_{|\Deltag^{\!\circ}}\otimes_{\cO_{\Deltag^{\!\circ}}}\ov{\cM_{|\sigma(\Deltag^{\!\circ})}}$.
\end{Assertion}

At this stage, we cannot prove the assertions, as we do not have any information on the derivatives of the functions $a_{j,k}$ introduced above. We need a more precise expression for $\hdel$, that we will derive in Lemma \ref{lem:B} below, using techniques analogous to that of \cite{B-M87}, that we recall in \T\ref{sec:Mellin}.

Fix $\hb_o\in \bS$, and let $\Deltag$ be as above. We will also assume that $\Deltag$ is small enough so that $\Deltag\cap\Sing\Lambda\subset\{\pm i\}$. Let $m$ be a section of $V_{(\hb_o)}^{>-1}\cM$ on $W\times\Deltag$ and $\mu$ a local section of $V_{(-\hb_o)}^{>-1}\cM$ on $W\times\sigma(\Deltag)$, for some open set $W\subset X$. If $b_m$ denotes the Bernstein polynomial of $m$, we consider the set $A(m)$ introduced in Corollary \ref{cor:Bernstein}, and take a \emph{minimal} subset $A'(m)\subset A(m)$ such that $A(m)\subset A'(m)-\NN$. Put $A'(m,\mu)=A'(m)\cap A'(\mu)$, so that $A'(m,\mu)-\NN$ contains all the $\gamma$ such that $\gamma\star\hb/\hb$ is a root of $b_m$ \emph{and} $b_\mu$. We also set $B'(m,\mu)=\{-\alpha-1\mid\alpha\in A'(m,\mu)\}$.

\begin{lemme}\label{lem:B}
There exist integers $\ell_0\geq0$ and $N\geq0$, such that, for any local sections $m,\mu$ as above, we have
\begin{equation}\tag*{(\protect\ref{lem:B})$(*)$}\label{eq:B*}
(\hb+1/\hb)^N\hdel(m,\ov\mu)=\sum_{\beta\in B'(m,\mu)}\sum_{\ell=0}^{\ell_0} f_{\beta,\ell}(t)\mt^{2(\beta\star\hb)/\hb}\frac{\Lt^\ell}{\ell!}
\end{equation}
on some punctured neighbourhood $X^*$ of $0\in X$, where the $f_{\beta,\ell}$ are $C^\infty$ functions $X\to\rH(\Deltag)$.
\end{lemme}

The proof of the lemma will be given below.

\begin{proof}[Proof of Assertion \ref{ass:1}]
Let us prove the $V_0\cR_{\cX|\Deltag}$-linearity of $\hdel$, the conjugate linearity being obtained similarly. Let $m,\mu$ be as in the lemma. The real part of the exponents in \ref{eq:B*} are $>-1$, as $\reel (\beta\star\hb_o)/\hb_o=\ell_{\hb_o}(\beta)+\ell_{-\hb_o}(\beta)$ for any $\hb_o\in\bS$. For a section $P$ of $V_0\cR_{\cX|\Deltag}$, we have to show that $u\defin P\cdot \hdel(m,\ov\mu)-\hdel(Pm,\ov\mu)=0$. But on the one hand, $u$ is supported on $\{t=0\}$ as $\hdel$ is known to be $V_0\cR_{\cX|\Deltag}$-linear away from $\{t=0\}$ by the correspondence of Lemma \ref{lem:twharm}. On the other hand, using \ref{eq:B*} and \eqref{eq:tdtu}, $(\hb+1/\hb)^Nu$ is $L^1_{\loc}$. Therefore, $(\hb+1/\hb)^Nu=0$ and, by \eqref{eq:sanstorsion}, $u=0$.

Notice that Formula \ref{eq:B*}, when restricted to $\Deltag^{\!\circ}$ and applied to the bases $\bme^{(\hb_o)},\bme^{(-\hb_o)}$, gives coefficients $f_{\beta,\ell}$ which are in $C^\infty(X^*\times\Deltag^{\!\circ})$ and holomorphic with respect to $\hb$. This still holds for any local sections $m,\mu$ on $X\times\Deltag^{\!\circ}$, by $\cO_\cX\otimes\cO_{\ov\cX}$-linearity. Then the argument for $\hdel$ can be applied to $\Cdelc$.
\end{proof}

\begin{proof}[Proof of Assertion \ref{ass:2}]
Any local section $m$ of $\cM$ on $X\times\Deltag^{\!\circ}$ can be written, by definition, as $\sum_{j\geq0}\partiall_t^jm_j$, where $m_j$ are local sections of $V_{(\hb_o)}^{>-1}\cM_{\hb_o}$ on $X\times\Deltag^{\!\circ}$. Therefore, in order to define $\Cdelc(m,\ov\mu)$, we write $m=\sum_{j\geq0}\partiall_t^jm_j$, $\mu=\sum_{k\geq0}\partiall_t^k\mu_k$, and put
\[
\Cdelc(m,\ov\mu)=\sum_{j,k\geq0}\partiall_t^j\partiall_{\ov t}^k \Cdelc(m_j,\ov\mu_k).
\]
This will be well defined if we prove that
\[
\sum_{j\geq0}\partiall_t^jm_j=0\implique \sum_{j\geq0}\partiall_t^j\Cdelc(m_j,\ov\mu_k)=0\quad\text{for any }\mu_k\in V_{(-\hb_o)}^{>-1}\cM_{-\hb_o},
\]
and a conjugate statement.

Set $u=\sum_{j\geq0}\partiall_t^j\Cdelc(m_j,\ov\mu_k)$ and $\wt u=(\hb+1/\hb)^N u$ for $N$ big enough. After Lemma \ref{lem:B}, there exists $N\geq0$ and a set $B'$ satisfying Properties \eqref{eq:propB1} and \eqref{eq:propB2} in Example \ref{ex:distr}, such that $\wt u$ can be written as
\[
\wt u=\sum_{\substack{\beta\in B'\\ \ell\in\NN}}\Big(\sum_{j\geq1}c_{\beta,\ell,j}\partiall_t^ju_{\beta,\ell}+ g_{\beta,\ell}(t)u_{\beta,\ell}\Big),
\]
where $c_{\beta,\ell,j}$ are constants and $g_{\beta,\ell}$ are $C^\infty$ on $X\times\Deltag^{\!\circ}$ and holomorphic with respect to $\hb$. The condition $\sum_{j\geq0}\partiall_t^jm_j=0$ implies that there exists $L\geq0$ such that $t^L\wt u=0$, as $\Cdelc$ is $V_0\cR_\cX$-linear. By taking $L$ large enough, using the freeness of the family $(u_{\beta,\ell})$ and \eqref{eq:tdtu}, one obtains that the coefficients $g_{\beta,\ell}$ and $c_{\beta,\ell,j}$ vanish identically, hence $\wt u\equiv0$. This implies that $u\equiv0$.
\end{proof}

Let us now show that the construction of $C$ glues along $\bS$. So, fix $\hb_o$ and $\Deltag$ as above, and denote by $C_{(\hb_o)}$ the sesquilinear pairing constructed above. Let $\hb\in\bS\cap\Deltag^{\!\circ}$. We will show that $C_{(\hb)}$ and $C_{(\hb_o)}$ coincide on $\cM_{\hb}\otimes\ov{\cM_{-\hb}}$.

By \eqref{eq:minV1} and \eqref{eq:minV2}, they coincide on $V_{(\hb)}^{-1-\epsilon}\cM_{\hb}\otimes \ov{V_{(-\hb)}^{-1-\epsilon}\cM_{-\hb}}$, for some $\epsilon>0$, as both are $L^1_{\loc}$ there and as they coincide with $h_{\bS}$ away from $t=0$. As both are $\cR_{(X,\ov X),\hb}$-linear, they also coincide on
\[
\big[V_{(\hb)}^{-1-\epsilon}\cM_{\hb}+\partiall_tt V_{(\hb)}^{>-1}\cM_{\hb}\big]\otimes\big[\ov{V_{(-\hb)}^{-1-\epsilon}\cM_{-\hb}+\partiall_tt V_{(-\hb)}^{>-1}\cM_{-\hb}}\big]
\]
which contains $V_{(\hb)}^{>-1}\cM_{\hb}\otimes\ov{V_{(-\hb)}^{>-1}\cM_{-\hb}}$ by \eqref{eq:minV2}. Last, by $\cR_{(X,\ov X),\hb}$-linearity, they coincide on $\cM_{\hb}\otimes\ov{\cM_{-\hb}}$.\qed

\begin{proof}[Proof of Lemma \ref{lem:B}]
We will use the method developed in \cite{B-M87}, using the existence of a ``good operator'' for $m$ or $\mu$, which is a consequence of the regularity property (REG). Fix local sections $m,\mu$ as above. By \eqref{eq:estimate}, $\hdel(m,\ov\mu)$ can be considered as a distribution on $X$ with values in $\rH(\Deltag)$ which, restricted to $X^*$, is $C^\infty$. Denote by~$p$ its order in some neighbourhood of $t=0$ on which we work. As $m,\mu$ are fixed, put, for any $N\geq 0$,
\[
\cI_{\chi,N}^{(k)}(s)=
\begin{cases}
(\hb+1/\hb)^N\langle \hdel(m,\ov\mu),\mt^{2s}t^k\chi(t)\,\itwopi dt\wedge d\ov t\rangle&\text{if }k\geq0,\\[8pt]
(\hb+1/\hb)^N\langle \hdel(m,\ov\mu),\mt^{2s}\ov t^{\module{k}}\chi(t)\,\itwopi dt\wedge d\ov t\rangle&\text{if }k\leq0,
\end{cases}
\]
for every function $\chi\in\cC_c^\infty(X,\rH(\Deltag))$. Then, for any such $\chi$, on the open set $2\reel s+\module{k}>p$, the function $s\mapsto \cI_{\chi,N}^{(k)}(s)$ takes values in $\rH(\Deltag)$ and is holomorphic. In the following, we fix $R\in{}]0,1[$ and we assume that $\chi\equiv0$ for $\mt>R$. We will be mainly interested to the case where $\chi\equiv1$ near $t=0$. Applying Theorem \ref{th:B-M} to the family $(\cI_{\chi,N}^{(k)}(s))_k$ will give the result. We will therefore show that the family $(\cI_{\chi,N}^{(k)}(s))_k$ satisfies the necessary assumptions for some $N$ big enough. We will assume that $k\geq0$, the case $k\leq0$ being obtained similarly, exchanging the roles of $m$ and $\mu$. Arguing exactly as in Lemma \ref{lem:polesI} and Remark \ref{rem:extension}, we find that $\cI_{\chi,N}^{(k)}(s)$ extends as a meromorphic function of $s\in\CC$ with values in $\rH(\Deltag)$, with at most poles along the sets $s=\gamma\star\hb/\hb$, with $\gamma\in A(m,\mu)-\NN$. Moreover, we may choose $N$ big enough so that all polar coefficients of $\cI_{\chi,N}^{(k)}(s)$ (for any $k$) take values in $\rH(\Deltag)$. We will fix such a $N$ and we forget it in the notation below. We put $\wthdel=(\hb+1/\hb)^N\hdel$.

By the regularity property (REG), there exists an integer $d$ and a relation
\[
(-\partiall_tt)^d\cdot m=\Big(\sum_{j=0}^{d-1}a_j(t)(-\partiall_tt)^j\Big)\cdot m
\]
for some sections $a_j$ of $\cO_{X\times\Deltag}$. It follows that, for any $\chi$ as above, we have
\[
(s+k)^d\cdot \cI_\chi^{(k)}(s)=\sum_{j=0}^{d-1}(s+k)^j\cI_{\chi_j}^{(k)}(s),
\]
where the $\chi_j$ only depend on the $a_\ell$ and $\chi$: as $\hdel$ is \emph{a~priori} $\cR_{(X,\ov X),\Deltag}$-linear away from \hbox{$t=0$} only, this equality only holds for $\reel s\gg0$; by uniqueness of analytic continuation, it holds for any $s$. Applying the same reasoning for $\mu$ we get
\[
s^d\cdot\cI_\chi^{(k)}(s)=\sum_{j=0}^{d-1}s^j\cI_{\psi_j}^{(k)}(s).
\]
For $n',n''\geq1$, we therefore have
\[
(s+k)^{dn'}s^{dn''}\cdot \cI_\chi^{(k)}(s)=\sum_{j'=0}^{(d-1)n'}\sum_{j''=0}^{(d-1)n''}(s+k)^{j'}s^{j''}\cI_{\chi_{n',n'',j',j''}}^{(k)}(s),
\]
where the functions $\chi_{n',n'',j',j''}$ only depend on $\chi, n',n'',j',j''$ and the $a_\ell$. As the current $\wthdel(m,\ov\mu)$ has order $p$, there exists a constant $C_{\chi,n', n'',j',j''}$ such that, if $2\reel s+\module{k}>p$, we have
\[
\norme{\cI_{\chi_{n',n'',j',j''}}^{(k)}(s)}_{\rH(\Deltag)}\leq C_{\chi,n', n'',j',j''}\cdot(1+\module{s}+\module{k})^pR^{2\reel s+\module{k}}.
\]
Let $n\geq0$. By summing the various inequalities that we get for $n'+n''\leq n+p$, we get the existence of a constant $C_{\chi, n}$ such that, if $2\reel s+\module{k}>p$, we have
\[
(1+\module{s}+\module{k})^n\norme{\cI_\chi^{(k)}(s)}_{\rH(\Deltag)}\leq C_{\chi,n}\cdot R^{2\reel s+\module{k}}.
\]
Let us now extend this for $2\reel s+\module{k}>p-q$ for any $q\in\NN$. For this purpose, consider now a Bernstein relation
\[
\prod_{\ell=0}^q b_m(-\partiall_tt+\ell\hb)\cdot m=t^q P(t,t\partiall_t)\cdot m
\]
and put $\delta=\deg b_m$. Fix $\chi$ as above. We therefore have
\[
\prod_{\ell=0}^qb_m\big(\hb(s+\ell)\big)\cdot \cI_\chi^{(k)}(s)=\sum_{j=0}^{q\delta-1}(s+k+q)^j\cI_{\chi_{q,j}}^{(k+q)}(s)
\]
for some $C^\infty$ functions depending on $\chi$, $q$ and $j$. By the previous reasoning, we get
\[
(1+\module{s}+\module{k}+q)^n\norme{\textstyle\prod_{\ell=0}^qb_m\big(\hb(s+\ell)\big)\cdot \cI_\chi^{(k)}(s)}_{\rH(\Deltag)}\leq C(n,q,\chi)\cdot R^{2\reel s+\module{k}+q}.
\]

Applying Theorem \ref{th:B-M}, we now find that there exists a finite family $(f_{\beta,\ell})_{\beta,\ell},g$ of $C^\infty$ functions $X\to\rH(\Deltag)$, $g$ being infinitely flat at $t=0$ and $\beta\in B(m,\mu)$, such that we have on $X^*$, for $\chi\equiv1$ near $t=0$,
\[
\chi(t)\wthdel(m,\ov\mu)=\Big(\sum_{\beta,\ell} f_{\beta,\ell}(t)\mt^{2(\beta\star\hb)/\hb}\frac{\Lt^\ell}{\ell!}+g(t)\Big).\qedhere
\]
\end{proof}

\begin{remarque}\label{rem:apresB}
We now give a more explicit statement when $m$ and $\mu$ locally lift sections of $\psi_{t,\alpha}\cM$. In such a case, arguing as for \ref{eq:polesI**}, we have for some $\ell\leq\ell_0$, putting $\beta=-\alpha-1$,
\begin{align*}
\chi(t)\wthdel(m,\ov\mu)&=\sum_{k=0}^{\ell} \wt c_{\beta,k}\mt^{2(\beta\star\hb)/\hb}\frac{\Lt^k}{k!}\\
&+t\sum_{k=0}^{\ell} \wt f_{\beta,k}(t)\mt^{2(\beta\star\hb)/\hb}\frac{\Lt^k}{k!}+\sum_{\reel\gamma>\reel\beta}\sum_{k=0}^{\ell_0} f_{\gamma,k}(t)\mt^{2(\gamma\star\hb)/\hb}\frac{\Lt^k}{k!}\,.
\end{align*}
The integer $\ell+1$ is smaller than or equal to the index of nilpotency of $\rN$ on $[m]$ or $[\mu]$. Moreover, each $\wt c_{\beta,k}$ is divisible by $(\hb+1/\hb)^N$, as the polar coefficients of the function $\cI_{\chi,0}^{(0)}(s)$ along $s=\alpha\star\hb$ take value in $\rH(\Deltag)$ for $\Deltag$ small enough (\cf Lemma \ref{lem:coefshol}).

If $m,\mu$ locally lift sections of $P\gr_\ell^{\rM}\psi_{t,\alpha}\cM$, then, applying $(t\partiall_t-\beta\star\hb)^\ell$ to both terms and noticing that, for $\hb\in\bS$, we have $\reel(\gamma\star\hb/\hb)=\reel\gamma$,
gives
\begin{equation}\tag*{(\protect\ref{rem:apresB})$(*)$}\label{eq:apresB*}
\chi(t)\wthdel((t\partiall_t-\beta\star\hb)^\ell m,\ov\mu)=\wt c_{\beta,w}\mt^{2(\beta\star\hb)/\hb}(1+\wt r(t))
\end{equation}
and $\wt r(t)$ tends to $0$ faster than some positive power of $\mt$.
\end{remarque}

\Subsection{Proof of the twistor properties for $(\cM,\cM,C,\id)$}
We will show that the properties of Definitions \ref{def:regtwt} and \ref{def:polarization} are satisfied for the object $(\cM,\cM,C,\id)$ when one takes specialization along a coordinate $t$. This is enough, according to Remark \ref{rem:ramif}.

Corollary \ref{cor:biquardRX} shows that $\cM$ satisfies Properties (HSD) and (REG) of the category $\MTr_{\leq 1}(X,0)$ of Def\ptbl\ref{def:regtwt}.

We will now show that, for any $\beta\in B$ and $\ell\in\NN$, the sesquilinear form $P\Psi_{t,\ell}^\beta C$ defines a polarized twistor structure on $(P\Psi_{t,\ell}^\beta\cM,P\Psi_{t,\ell}^\beta\cM)$ with polarization $\cS=(\id,\id)$.

By Definition \ref{def:nearby}, we can replace $\cM$ with $\wt\cM$ and use the basis $\bme^o$ (with its primitive vectors) introduced in Remark \ref{rem:although}\eqref{rem:although2} to compute $P\Psi_{t,\ell}^\beta C$. We will do the computation locally on small compact sets $\Deltag$ as in \T\ref{subsec:constrC}, so that we may lift the primitive vectors $e_i^o$ to $e_i^{(\hb_o)}$. If we only use the dominant term in Formula \eqref{eq:bme} when computing $P\Psi_{t,\ell}^\beta \Cdelc(e_i^{(\hb_o)},\ov{e_j^{(-\hb_o)}})$, we recover the computation made in Proposition \ref{prop:basicpolarisation}. Therefore, in order to conclude, we only need to show that the non dominant term in \eqref{eq:bme} does not contribute to $P\Psi_{t,\ell}^\beta \Cdelc(e_i^{(\hb_o)},\ov{e_j^{(-\hb_o)}})$. Here also, the estimate given by \eqref{eq:estimate}, Formula \eqref{eq:bme} and Lemma \ref{lem:killing} is not strong enough to eliminate the non dominant term in the computation of $P\psi_{t,\ell}^\beta \Cdelc$, and we use Lemma \ref{lem:B} and Remark \ref{rem:apresB}.

If $e_i^{(\hb_o)},e_j^{(-\hb_o)}$ locally lift primitive sections of weight~$\ell$, Formulas \eqref{eq:estimate} and \ref{eq:apresB*} for $\wtCdelc((t\partiall_t-\beta\star\hb)^\ell e_i^{(\hb_o)},\ov{e_j^{(-\hb_o)}})$ give the same dominant term $\wt c_{\beta,w}\mt^{2(\beta\star\hb)/\hb}$. As $\mt^{-\delta}\wt r(t)\to0$ for some $\delta>0$, the remaining term in \ref{eq:apresB*} does not contribute to the computation of the residue in Definition \ref{def:psiC}. Dividing now by $(\hb+1/\hb)^N$ gives the desired result.\qed

\section[Proof of Theorem $\ref{thssimplecourbes}$, second part]{Proof of Theorem \ref{th:ssimplecourbes}, second part}\label{sec:second}

\subsection{}\label{subsec:second}
Let $(\cT,\cS)$ be a polarized regular twistor $\cD_X$-module of weight~$0$. We will assume that it has strict support $X$ and that $\cM'=\cM''\defin\cM$, $\cS=(\id,\id)$. Its restriction to the complement $X^*=X\moins P$ of a finite set of points is therefore a smooth twistor $\cD_X$-module of weight~$0$ (Proposition \ref{prop:genericMT}) and corresponds to a harmonic bundle $(H,D_V,h)$ on $X^*$ in such a way that $C=h_\bS$ (Lemma \ref{lem:twharm}). Put $M=\Xi_{\DR}(\cM')$ and $\wt M=\cO_X(*P)\otimes_{\cO_X}M$. By definition we have $M_{|X^*}=V=\ker D''_V$. We will show that the Hermitian metric $h$ is tame with respect to $\wt M$ and that its parabolic filtration is the canonical filtration of $\wt M$.

\medskip
We now work locally near a point of $P$, restricting $X$ when necessary. We will mainly work with $\wt\cM$ instead of $\cM$. Remark first:
\begin{lemme}\label{lem:family1}
For any $\hb_o\in\Omega_0$ and any $a<0$, the $\cO_\cX$-module $V_a^{(\hb_o)}\wt\cM$ is locally free on the open set where it is defined.
\end{lemme}

\begin{proof}
Recall that, for $a<0$, we have $V_a^{(\hb_o)}\wt\cM=V_a^{(\hb_o)}\cM$ near any $\hb_o\in\Omega_0$ (\cf Lemma \ref{lem:localstrictspe}).

By the regularity assumption, $V_a^{(\hb_o)}\wt\cM$ is $\cO_\cX$-coherent. Moreover, on $t\neq0$, it is $\cO_\cX$-locally free, being there equal to $\cM$, which is there a smooth twistor $\cD$-module. As $t$ is injective on $V_a^{(\hb_o)}\cM$ when $a<0$ (\cf Remark \ref{rem:psi}\eqref{rem:psia}), it follows that $V_a^{(\hb_o)}\cM$ has no $\cO_\cX$-torsion. It is therefore enough to show that $V_a^{(\hb_o)}\cM/tV_a^{(\hb_o)}\cM$ is $\cO_{\Omega_0}$-locally free. By Proposition \ref{prop:twstrict}, each $\psi_{t,\alpha}\cM$ is $\cO_{\Omega_0}$-locally free. As $\gr_a^{V^{(\hb_o)}}\cM=\oplus_{\alpha\mid\ell_{\hb_o}(\alpha)=a}\psi_{t,\alpha}\cM$ near $\hb_o$, it follows that $\gr_a^{V^{(\hb_o)}}\cM$, hence $V_a^{(\hb_o)}\cM/tV_a^{(\hb_o)}\cM$, is $\cO_{\Omega_0}$-locally free.
\end{proof}

The proof will now consist in constructing bases $\bme$ and $\varepsilong$ as in \T\ref{sec:first}. However, it will go in the reverse direction. Indeed, we will first construct local bases $\bme^{(\hb_o)}$ of $V_{<0}^{(\hb_o)}\wt \cM$ near any $\hb_o\in\bS$. We then construct local bases $\varepsilon^{(\hb_o)}$ using Formula \eqref{eq:bme} as a definition. Then we glue together these local bases and we recover an orthonormal basis of the bundle $H$ with respect to the harmonic metric $h$, in such a way that a formula like \eqref{eq:bme} still holds.

\subsection{}
We will first recover precise formulas for the sesquilinear pairing $C$, as in Lemma \ref{lem:B}, starting from the definition of a regular twistor $\cD$-module. Fix $\hb_o\in\bS$ and let $\Deltag=\Delta_{\hb_o}(\eta)$ be a small closed disc centered at $\hb_o$ on which the $V_\bbullet^{(\hb_o)}\cM$ is defined. On $X^*$, $C=h_\bS$ takes values in $\cCh{\cX^*}$. If the radius $\eta$ is small enough, we may assume that $\hdel$ is a $C^\infty$ function $X^*\to\rH(\Deltag)$, that can be written as a power series $\sum_nc_n(t)(\hb-\hb_o)^n$, where each $c_n(t)$ is a $C^\infty$ function on $X^*$.

\begin{lemme}\label{lem:apresB}
Assume that $m$ lifts a section of $\psi_{t,\alpha_1}\cM_{|\Deltag}$ and $\mu$ a section of $\psi_{t,\alpha_2}\cM_{|\sigma(\Deltag)}$. Then,
\begin{enumerate}
\item\label{lem:apresB1}
if $\alpha_1=\alpha_2\defin\alpha$ and $\beta\defin-\alpha-1$, we have, for some $\ell\in\NN$ and any $t\in X^*$,
\begin{equation}\tag*{(\protect\ref{lem:apresB})$(*)$}\label{lem:apresB*}
\hdel(m,\ov\mu)=\mt^{2(\beta\star\hb)/\hb}\Lt^\ell(c_{m,\ov\mu}+r(t))
\end{equation}
with $c_{m,\ov\mu}\in\rH(\Deltag)$, $r(t)$ takes values in $\rH(\Deltag)$ and there exists $\delta>0$ such that $\lim_{t\to0}\Lt^\delta\norme{r(t)}_{\rH(\Deltag)}=0$;
\item\label{lem:apresB2}
if $\alpha_1\neq\alpha_2$, there exists $\delta>0$ such that, for any $\ell\in\ZZ$,
\begin{equation}\tag*{(\protect\ref{lem:apresB})$(**)$}\label{lem:apresB**}
\lim_{t\to0}\mt^{-\delta}\norme{\mt^{-(\beta'_1+i\beta''_1\hb)}\mt^{-(\beta'_2+i\beta''_2/\hb)}\Lt^{-\ell}\hdel(m,\ov\mu)}_{\rH(\Deltag)}=0.
\end{equation}
\end{enumerate}
\end{lemme}

\begin{proof}
Lemma \ref{lem:B} and Remark \ref{rem:apresB} apply in this situation, as the argument only uses Lemma \ref{lem:polesI} and the regularity property (REG). If $\alpha_1=\alpha_2=\alpha$, we get that, near $t=0$, one has, for some $\ell\geq0$ and some $N\geq0$,
\[
(\hb+1/\hb)^N\hdel(m,\ov\mu)=\mt^{\beta\star\hb/\hb}\Lt^\ell((\hb+1/\hb)^Nc(m,\ov\mu)+\wt r(t)),
\]
where $c(m,\ov\mu)$ is in $\rH(\Deltag)$ and there exists $\delta>0$ such that $\lim_{t\to0}\norme{\wt r(t)}_{\rH(\Deltag)}=0$. Therefore, $\wt r(t)=(\hb+1/\hb)^Nr(t)$ where $r$ takes values in $\rH(\Deltag)$, and (by the maximum principle) we also have $\lim_{t\to0}\norme{r(t)}_{\rH(\Deltag)}=0$.

If $\alpha_1\neq\alpha_2$, then the same argument as in Remark \ref{rem:apresB} shows that, near $t=0$ and if the radius $\eta>0$ of $\Deltag$ is small enough, we have
\[
(\hb+1/\hb)^N\hdel(m,\ov\mu)=\sum_{\gamma}\sum_{k=0}^{\ell_0}\wt f_{\gamma,k}(t)\mt^{2(\gamma\star\hb)/\hb}\Lt^k,
\]
where the first sum is taken for $\gamma\in\Lambda$ such that $2\reel\gamma>\ell_{\hb_0}(\beta_1)+\ell_{-\hb_o}(\beta_2)$, and $\wt f_{\gamma,k}:X\to\rH(\Deltag)$ are $C^\infty$. This shows that \ref{lem:apresB**} holds for $(\hb+1/\hb)^N\hdel(m,\ov\mu)$, and one concludes as above that it holds for $\hdel(m,\ov\mu)$.
\end{proof}

\subsection{Construction of local bases $\bme^{(\hb_o)}$ of $V_{<0}^{(\hb_o)}\wt\cM$}\label{subsec:constrlocbases}

By definition, for any $\alpha$ with $\reel\alpha\in[-1,0[$ and $\ell\in\NN$, the triple $(P\gr_\ell^{\rM}\Psi_{t,\alpha}\cM,P\gr_\ell^{\rM}\Psi_{t,\alpha}\cM,P\Psi_{t,\alpha,\ell}C)$ is a twistor structure of weight~$0$ with polarization $(\id,\id)$ (\cf Remark \ref{rem:psiadj}). Choose therefore a basis $\bme^o_{\alpha,\ell,\ell}$ of $P\gr_\ell^{\rM}\Psi_{t,\alpha}\cM= P\gr_\ell^{\rM}\Psi_{t,\alpha}\wt\cM$ (\cf Remark \ref{rem:nearby}\eqref{rem:nearby2}) which is orthonormal for $P\Psi_{t,\alpha,\ell}C$, when restricted to $\nb(\bS)$: this is possible according to the twistor condition (\cf \T\ref{sec:twst0}). Extend the basis $\bme^o_{\alpha,\ell,\ell}$ as a basis $\bme^o_{\alpha,\ell}=(\bme^o_{\alpha,\ell,w})_{w\in\ZZ}$ of $\gr^{\rM}\psi_{t,\alpha}\wt\cM$ for which the matrix $-\rY_\alpha$ of $\rN$ is as in the basic example of \T\ref{subsec:basic}.

Fix now $\hb_o\in\Omega_0$. Locally near $\hb_o$, lift the family $\bme^o_{\alpha,\ell,\ell}$, defined as above, as a family $\wt\bme^{(\hb_o)}_{\alpha,\ell,\ell}$ of local sections of $V^{(\hb_o)}_{\ell_{\hb_o}(\alpha)}\wt\cM_{\hb_o}$. Similarly, extend the family $\wt\bme^{(\hb_o)}_{\alpha,\ell,\ell}$ in a family $\wt\bme^{(\hb_o)}_{\alpha,\ell}=(\wt\bme^{(\hb_o)}_{\alpha,\ell,w})_{w\in\ZZ}$ so that, putting $\beta=-\alpha-1$, the lift of $(-\rN)^j\wt e^o$, where $e^o$ is any element of $\bme^o_{\alpha,\ell,\ell}$ and $j\leq\ell$, is $(t\partiall_t-\beta\star\hb)^j\wt e$.

Last, for any $q\in\ZZ$, put $\wt\bme^{(\hb_o)}_{\alpha+q,\ell}=t^{-q}\wt\bme^{(\hb_o)}_{\alpha,\ell}$. By choosing $q_\alpha$ so that $\ell_{\hb_o}(\alpha+q_\alpha)\in[-1,0[$, we get a basis of $V_{<0}^{(\hb_o)}\cM$ near $\hb_o$ (this $\cO_\cX$-module is known to be $\cO_\cX$-locally free, \cf Lemma \ref{lem:family1}), as it induces a basis of $V_{<0}^{(\hb_o)}\wt\cM/tV_{<0}^{(\hb_o)}\wt\cM$.

Let $\wt\bme^{\prime(\hb_o)}$ be any other local $\cO_\cX$-basis of $V_{<0}^{(\hb_o)}\cM$ inducing $\bme^o$ on $\oplus_\alpha\gr^{\rM}\psi_{t,\alpha}^{(\hb_o)}\wt\cM$ near $\hb_o$. Define the matrix $P$ by the equality $\wt\bme^{\prime(\hb_o)}=\wt\bme^{(\hb_o)}\cdot(\id+P(t,\hb))$, and recall that the matrix $A(t,\hb)$ was defined just before Formula \eqref{eq:bme}. Using the definition of the monodromy filtration and of $\rH$ and denoting by $\Deltag$ a sufficiently small closed disc centered at $\hb_o$, one easily gets:

\begin{lemme}\label{lem:Pthb}
There exists $\delta>0$ such that $\lim_{t\to0}\Lt^\delta\norme{A^{-1}PA}_{\rH(\Deltag)}=0$.\qed
\end{lemme}

Define then, for $t\neq0$, the basis $\wt\varepsilong^{(\hb_o)}$ by $\wt\bme^{(\hb_o)}=\wt\varepsilong^{(\hb_o)}\cdot A(t,\hb)$. Notice that, if $\wt\varepsilong^{\prime(\hb_o)}=\wt\varepsilong^{(\hb_o)}(\id+Q(t,\hb))$ is defined similarly from another basis $\wt\bme^{\prime(\hb_o)}$ then, according to Lemma \ref{lem:Pthb}, we have
\begin{equation}\label{eq:Qthb}
\lim_{t\to0}\Lt^\delta\norme{Q}_{\rH(\Deltag)}=0.
\end{equation}

\subsection{Orthonormality with respect to $C$}
Fix now $\hb_o\in\bS$. Recall then that $\sigma(\hb_o)=-\hb_o$. Using Lemma \ref{lem:apresB}, arguing first at the level of primitive vectors as in Remark \ref{rem:apresB}, one gets:

\begin{lemme}\label{lem:Rthb}
The matrix $\bC^{(\hb_o)}$ of $\hdel$ in the bases $\wt\varepsilong^{(\pm\hb_o)}$ takes the form $\id+R(t,\hb)$, and there exists $\delta>0$ such that $\lim_{t\to0}\Lt^\delta\norme{R(t,\hb)}_{\rH(\Deltag)}=0$.\qed
\end{lemme}

In other words, the pair of local bases $\wt\varepsilong^{(\pm\hb_o)}$ is asymptotically orthonormal for $\hdel$, with speed a negative power of $\Lt$. Notice that, if $\wt\varepsilong^{\prime(\pm\hb_o)}$ is defined similarly from another basis $\wt\bme^{\prime(\hb_o)}$ then, according to Lemma \ref{lem:Pthb}, $\wt\varepsilong^{\prime(\pm\hb_o)}$ has the same property (maybe with a different $\delta$).

\Subsection{Globalization of the asymptotically orthonormal local bases}\label{subsec:glob}
For any $r\in[0,1]$, let $\bS_r$ denote the circle of radius $r$ in $\Omega_0$. For any such $r$, cover $\bS_r$ by a finite number of open discs $\Delta_{\hb_o}$ ($\hb_o\in\bS_r$) on the closure of which the previous construction applies. One can assume that the intersection of any three distinct such open discs is empty. On the intersection $\Delta_{ij}$ of two open sets $\Delta_i$ and $\Delta_j$, the base change $\wt\varepsilong^{(i)}=\wt\varepsilong^{(j)}\cdot (\id+Q_{ij}(t,\hb))$ satisfies $\lim_{t\to0}\Lt^\delta\norme{Q_{ij}(t,\hb)}_{\rH(\Delta_j)}=0$ for some $\delta>0$, according to Lemma \ref{lem:Pthb} and \eqref{eq:Qthb}.

\begin{lemme}\label{lem:thGrauert}
There exists $\delta>0$ and $C^0$ matrices $R_i(t,\hb)$ on $D^*\times\Delta_i$, holomorphic with respect to $\hb$, such that $\lim_{t\to0}\Lt^\delta\norme{R_i(t,\hb)}_{\rH(\Delta_j)}=0$ and $\id+Q_{ij}(t,\hb)=(\id+R_{j}(t,\hb))(\id+R_{i}(t,\hb))^{-1}$ on $X^*\times\Delta_{ij}$.
\end{lemme}

\begin{proof}
We consider the family $(\id+Q_{ij}(t,\hb))_{ij}$ as a cocycle relative to the covering $(\Delta_i)_i$ with values in the Banach Lie subgroup of $\GL_d(C^0(X^*))$ of matrices $U$ such that $U-\id\in\Mat_d(C^0_\delta(X^*))$, where $d$ is the size of the matrices $Q$ (generic rank of $\cM$) and $C^0_\delta(X^*)$ is the Banach algebra of continuous functions $\varphi$ on $X^*$ such that $\norme{\Lt^\delta\varphi}_\infty<+\infty$, with the corresponding norm $\norme{\cbbullet}_{\infty,\delta}$. By changing $\delta$, one can moreover assume that $\norme{Q_{ij}}_{\infty,\delta}<1$ for all $i,j$. Clearly, this cocycle can be deformed continuously to the trivial cocycle. Hence, by standard results (\cf \cite{Cartan49}), the Banach bundle defined by this cocycle is topologically trivial. It follows then from a generalization of Grauert's theorem to Banach bundles, due to L.~Bungart (\cf \cite{Bungart68}, see also \cite{Leiterer90}) that this bundle is holomorphically trivial. The trivialization cocycle takes the form given in the lemma if one chooses a smaller $\delta$.
\end{proof}

According to Lemma \ref{lem:thGrauert}, the basis $\wt\varepsilong^{(r)}$ defined by $\wt\varepsilong^{(r)}=\wt\varepsilong^{(i)}\cdot(\id+R_{i}(t,\hb))$ is globally defined on some open neighbourhood $\nb(\bS_r)$.

Let $\Delta_0$ be the closed disc centered at $0$ and of radius $1$ in $\Omega_0$. Cover $\Delta_0$ by a finite number of $\nb(\bS_r)$ on which the previous construction applies. Assume that the intersection of three distinct open sets is empty. Apply the previous argument to get a basis $\wt\varepsilong$ globally defined on $\Delta_0$ such that, according to Lemma \ref{lem:Rthb}, the matrix $\bC$ of $h_\bS$ takes the form $\id+R(t,\hb)$ on $X^*\times\nb(\bS)$ with $\lim_{t\to0}\Lt^\delta\norme{R(t,\hb)}_{\infty}=0$.

\begin{lemme}\label{lem:rigidP1}
If $X$ is small enough, there exists a $d\times d$ matrix $S(t,\hb)$ such that
\begin{itemize}
\item
$S(t,\hb)$ is continuous on $X^*\times\nb(\Delta_0)$ and holomorphic with respect to $\hb$,
\item
$\lim_{t\to0}\norme{S(t,\hb)}_{\rH(\Delta_0)}=0$,
\item
$\id+R(t,\hb)=\big(\id+S^*(t,\hb)\big)\cdot\big(\id+S(t,\hb)\big).$
\end{itemize}
\end{lemme}

As usual, we denote by $S^*$ the adjoint matrix of $S$ (where conjugation is taken as in \T\ref{subsubsec:conj}).

Define the continuous basis $\varepsilong$ of $\wt\cM_{X^*\times\nb(\Delta_0)}$ by $\varepsilong=\wt\varepsilong'\cdot(\id+\ov S(t,\hb))^{-1}$. In this basis, the matrix of $C$ is equal to $\id$, after the lemma. This means that $\varepsilong$ is a continuous basis of $H$ and that it is orthonormal for the harmonic metric associated with $(\cT,\cS)$ on $X^*$. Moreover, near any $\hb_o\in\bS$, the base change between the local basis $\wt\bme^{(\hb_o)}$ and $\varepsilong$ take the form of Formula \eqref{eq:bme} (we only get here that the corresponding $Q^{(\hb_o)}$ tends to $0$ with $t$, and not the logarithmic speed of decay). Arguing as in Remark \ref{rem:although}\eqref{rem:although3}, we conclude in particular that, when $\hb_o=1$, the $V$-filtration of $\wt M$ is equal to the parabolic filtration defined the metric $h$ obtained from $C$, as was to be proved\footnote{I thank the referee for indicating that a previous proof was not complete.}.\qed

\begin{proof}[Sketch of proof of Lemma \ref{lem:rigidP1}]
We view $\id+R(t,\hb)$ as a family, parametrized by~$X$, of cocycles of $\PP^1$ relative to the covering $\big(\nb(\Delta_0),\ov{\nb(\Delta_0)}\big)$. At $t=0$, this cocycle is equal to the identity, hence the corresponding bundle is trivial. By the rigidity of the trivial bundle on $\PP^1$, this remains true for any $t$ small enough. More precisely, arguing for instance as in \cite[Lemme 4.5]{Malgrange83db}, there exist unique invertible matrices $\wt\Sigma'$ continuous on $X\times\nb(\Delta_0)$ and $\Sigma''$ continuous on $X\times\ov{\nb(\Delta_0)}$, both holomorphic with respect to $\hb$, such that $\Sigma''(\cbbullet,\infty)\equiv\id$ and $\id+R(t,\hb)=\Sigma''\cdot\wt\Sigma'$. Set $\Sigma_0(t)=\wt\Sigma'(t,0)$, which is invertible, and $\Sigma'=\Sigma_0^{-1}\wt\Sigma'$; hence $\Sigma'(t,0)\equiv\id$ and $\id+R(t,\hb)=\Sigma''\cdot\Sigma_0\cdot\Sigma'$. Moreover, these matrices are all equal to $\id$ at $t=0$, by uniqueness.

Recall now that, as we have reduced to $\cS=(\id,\id)$ and the weight equal to $0$, the sesquilinear pairing satisfies $C^*=C$, hence $R(t,\hb)$ satisfies $R^*=R$. We thus have $\id+R(t,\hb)=\Sigma^{\prime*}\cdot\Sigma_0^*\cdot\Sigma^{\prime\prime*}$.

As $\Sigma^{\prime*}(t,\infty)\equiv\id$, we have, by uniqueness of such a decomposition, $\Sigma^{\prime*}=\Sigma''$. Moreover, $\Sigma_0^*=\Sigma_0$ (in the usual sense, as no $\hb$ is involved).

As $\Sigma_{0\mid t=0}=\id$, we can write locally $\Sigma_0=T^*T$ and put $\id+S=T\Sigma'$.
\end{proof}

It remains to explain that $(\cT,\cS)$ is locally isomorphic to $(\cT,\cS)_h$ constructed in \T\ref{sec:first}, where $(H,h)$ is defined in \T\ref{subsec:second}, as asserted in the sketch of \T\ref{subsec:sketch}. What we have done above is to show that, starting from a polarized twistor $\cD$-module $(\cT,\cS)$ on a curve $X$, we can recover the properties that have been proved by Simpson and Biquard for $(H,h)$. As we have similar bases $\bme^{(\hb_o)}$ on which we can compute the connection and the $h$-norm of which has moderate growth at $t=0$, locally uniformly with respect to $\hb$, we conclude that the localization $\wt\cM$ of $\cM$ is equal to $\wt\cM$ defined by Cor\ptbl\ref{cor:biquard}\eqref{cor:biquard1}. Then, the corresponding minimal extensions $\cM$ are the same. Last, the gluings $C$ coincide away from the singular point $P$, hence they coincide locally near any $\hb_o\in\bS$ on $V_{<0}^{(\hb_o)}$, as they take values in $L^1_{\loc}$ there, and therefore they coincide on $\cM$ by $\cR_{(X,\ov X),\bS}$-linearity.\qed

\setcounter{section}{0}
\let\oldthesection\thesection
\renewcommand{\thesection}{\thechapter.\Alph{section}}
\section{Mellin transform and asymptotic expansions}
\label{sec:Mellin}

We will recall here, with few minor modifications, some results of \cite{B-M89}. Fix a finite set $B\subset\CC$ such that no two complex numbers in $B$ differ by a nonzero integer. Let $\Deltag$ be any compact set in $\Omega_0\moins\{0\}$, which is the closure of its interior. We keep the notation of \T\ref{num:RX} concerning $\rH(\Deltag)$. Notice that the set $K_B=\{(\beta\star\hb)/\hb\mid\beta\in B\text{ and } \hb\in\Deltag\}$ is compact.

\begin{definition}
For a $C^\infty$ function $f:\CC^*\to\rH(\Deltag)$ satisfying
\begin{enumerate}
\item[(i)]
$f\equiv0$ for $\mt\geq R$ (for some $R>0$),
\item[(ii)]
$f$ has moderate growth at $t=0$, \ie there exists $\sigma_0\in\RR$ such that
\[
\lim_{t\to0}\mt^{2\sigma_0}\norme{f(t)}_{\rH(\Deltag)}=0,
\]
\end{enumerate}
the \emph{Mellin transform with parameters $k',k''\in\NN$} is
\[
\cI^{(k',k'')}_f(s)=\int_\CC\mt^{2s}t^{k'}\ov t^{k''} f(t)\,\itwopi dt\wedge d\ov t.
\]
\end{definition}

Notice that $\cI^{(k',k'')}_f(s)$ is holomorphic in the half-plane $2\reel s+k'+k''>\sigma_0-2$. It is clearly enough, up to a translation of $s$ by an integer, to consider the functions $\cI^{(k,0)}_f(s)$ and $\cI^{(0,k)}_f(s)$ for $k\in\NN$. For the sake of simplicity, we will denote, for any $k\in\ZZ$,
\[
\cI_f^{(k)}(s)=\begin{cases}
\cI^{(k,0)}_f(s)&\text{if }k\geq0,\\
\cI^{(0,\module{k})}_f(s)&\text{if }k\leq0.
\end{cases}
\]

\begin{definition}\label{def:nils}
Let $B$ be as above and $R\in\RR_+$.
\begin{enumerate}
\item\label{nils1}
Recall (\cf \T\ref{subsec:defLt}) that we put $\Lt=\Lt=\module{\log\mt^2}$. A $C^\infty$ function $f:\CC^*\to\rH(\Deltag)$ has \emph{Nilsson type $B$} at $t=0$ if there exist $L\in\NN$ and, for any $\beta\in B$ and $\ell\in[0,L]\cap\NN$, $C^\infty$ functions $f_{\beta,\ell}:\CC\to\rH(\Deltag)$ such that $f$ can be written on $\CC^*$ as
\begin{equation*}
f(t)=\sum_{\beta\in B}\sum_{\ell=0}^{L} f_{\beta,\ell}(t)\mt^{2(\beta\star\hb)/\hb}\Lt^\ell.
\end{equation*}
\item\label{mellin2}
Let $\cI=\big(\cI^{(k)}(s)\big)_{k\in\ZZ}$ be a family of meromorphic functions of $s\in\CC$ with values in $\rH(\Deltag)$, \ie of the form $\varphi(s)/\psi(s)$ with $\varphi,\psi:\CC\to\rH(\Deltag)$ holomorphic and $\psi\not\equiv0$. We say that $\cI$ has type $(B,R)$ if there exist
\begin{itemize}
\item
a polynomial $b(s)$ equal to a finite product of terms $s+n+(\beta\star\hb)/\hb$ with $n\in\NN$ and $\beta\in B$,
\item
for any $N\in\NN$, any $\sigma\in\RR^+$, a constant $C(N,\sigma,R)$,
\end{itemize}
such that, for any $k\in\ZZ$, $\cI^{(k)}$ satisfies on the half plane $\reel s>-\sigma-1-\module{k}/2$
\begin{equation*}
(1+\module{k}+\module{s})^{N}\Big\Vert\Big(\textstyle\prod_{\nu\in[0,\sigma]\cap\NN}b(s+\module{k}+\nu)\Big)\cI^{(k)}(s)\Big\Vert_{\rH(\Deltag)}\leq C(N,\sigma,R)R^{2\reel s+\module{k}}.
\end{equation*}
\end{enumerate}
\end{definition}

\begin{remarques}
Let $\cI$ be a family having type $(B,R)$.
\begin{enumerate}
\item
Notice that, $K_B$ being compact, there exists $\sigma_0$ such that $\cI^{(k)}(s)$ is holomorphic in the half plane $\reel s>\sigma_0-\module{k}/2$ for any $k\in\ZZ$; in particular, given a half plane $\reel s>\sigma$, there exist only a finite number of $k\in\ZZ$ for which $\cI^{(k)}$ has a pole on this half plane; moreover, the possible poles of $\cI^{(k)}(s)$ are $s=-\module{k}-n-(\beta\star\hb)/\hb$ for some $n\in\NN$ and $\beta\in B$, and the order of the poles is bounded by some integer $L$.
\item
If $\cI$ has type $(B,R)$ with some polynomial $b$, it also has type $(B,R)$ with any polynomial of the same kind that $b$ divides. This can be used to show that the sum of two families having type $(B,R)$ still has type $(B,R)$.
\item
The polar coefficients of $\cI^{(k)}$ at its poles may, as functions of $\hb$, have poles for some purely imaginary values of $\hb$. However, there exists a polynomial $\lambda(\hb)$ with poles at $i\RR$ at most, such that all possible polar coefficients of $\lambda\cdot\cI^{(k)}$, for any $k\in\ZZ$, are in $\rH(\Deltag)$.
\end{enumerate}
\end{remarques}

\begin{theoreme}[\cite{B-M87,B-M89}]\label{th:B-M}
The Mellin transform $f\mapsto\big(\cI^{(k)}_f(s)\big)_{k\in\ZZ}$ gives a one-to-one correspondence between $C^\infty$ functions $f:\CC^*\to\rH(\Deltag)$ of Nilsson type $B$ at $t=0$ and having support in $\mt\leq R$, and families of meromorphic functions $\big(\cI^{(k)}(s)\big)_{k\in\ZZ}:\CC\to\rH(\Deltag)$ of type $(B,R)$ having polar coefficients in $\rH(\Deltag)$.
\end{theoreme}

\begin{proof}
We will only indicate the few modifications to be made to the proof given in \loccit Given a family $\big(\cI^{(k)}(s)\big)_{k\in\ZZ}$ of type $(B,R)$, each $\cI^{(k)}(s)$ being holomorphic in a half plane $\reel s>\sigma_0-\module{k}/2$, there exists a $C^\infty$ function $f:\CC^*\to\rH(\Deltag)$ which has moderate growth at $t=0$ and vanishes for $\mt\geq R$, such that $\cI^{(k)}(s)=\cI^{(k)}_f(s)$ for all $k\in\ZZ$ and $\reel s>\sigma_0$. As the polar coefficients are in $\rH(\Deltag)$, one may construct as in \loccit, using Borel's lemma, a function $g(t)=\sum_{\beta,\ell}g_{\beta,\ell}(t)\mt^{2(\beta\star\hb)/\hb}\Lt^\ell$, the functions $g_{\beta,\ell}:\CC\to\rH(\Deltag)$ being $C^\infty$, such that the family $\big(\cI^{(k)}_{f-g}(s)\big)_{k\in\ZZ}$ has type $(B,R)$ and all functions $\cI^{(k)}_{f-g}(s)$ are entire.

Fix $\hb_o\in\Deltag$ and choose an increasing sequence $(\sigma_i)_{i\in\NN}$ with $\lim_i\sigma_i=+\infty$, such that no line $\reel s=\sigma_i$ contains a complex number of the form $n-(\beta\star\hb_o)/\hb_o$, for $n\in\ZZ$ and $\beta\in B$. Thus, there is a neighbourhood $V(\hb_o)\subset\Deltag$ such that the distance between the lines $\reel s=\sigma_i$ and the set $\{n-(\beta\star\hb)/\hb\mid n\in\ZZ,\; \beta\in B,\; \hb\in V(\hb_o)\}$ is bounded from below by a positive constant. Let $H_{i,k}$ be the half plane $\reel s>-\sigma_i-1-\module{k}/2$ and let $D$ be the union of small discs centered at these points $n-(\beta\star\hb_o)/\hb_o$ (small enough so that they do not cut the lines $\reel s=\sigma_i$). On $H_{i,k}\moins D$, $\module{\prod_{\nu\in[0,\sigma_i]\cap\NN}b(s+\module{k}+\nu)}$ is bounded from below, so we have an estimation on this open set:
\begin{equation*}
(1+\module{k}+\module{s})^{N}\norme{\cI^{(k)}_{f-g}(s)}_{\rH(V(\hb_o))}\leq C'(N,\sigma_i,R)R^{2\reel s+\module{k}}
\end{equation*}
By the maximum principle, $\cI^{(k)}_{f-g}(s)$ being entire, this estimation holds on $H_{i,k}$. This implies that $\norme{f-g}_{\rH(V(\hb_o))}=O(\mt^{\sigma_i})$ for any $i$. By compactness of $\Deltag$, we have $\norme{f-g}_{\rH(\Deltag)}=O(\mt^{\sigma})$ for any $\sigma>0$, so $f-g$ is $C^\infty$ and infinitely flat at $t=0$. We may therefore change some $g_{\beta,\ell}$ to get the decomposition of $f$.
\end{proof}

\begin{remarque}\label{rem:polesI}
It is not difficult to relate the order of the poles of $\cI^{(k)}_f(s)$ with the integer $L$ in Definition \ref{def:nils}\eqref{nils1}. In particular, if $\cI^{(k)}_f(s)/\Gamma(s+\module{k}+1)$ has no pole for any $k$, then one can choose $L=0$ and $B=\{0\}$ in \ref{def:nils}\eqref{nils1}.
\end{remarque}

\section{Some results of O\ptbl Biquard}\label{app:biquard}
The results in this section are direct consequences of \cite{Biquard97}, although the lemma in \T\ref{app:biquard1} is not explicitly stated there. They were explained to me by O\ptbl Biquard, whom I thank.

\subsection{Better estimate of some perturbation terms}\label{app:biquard1}
We indicate here how to obtain \emph{a posteriori} better estimates for the perturbation terms $P',P''$ of Theorem~\ref{th:biquard}, that is, Property (iii) used in the proof of Lemma \ref{lem:killing}.

One starts with the holomorphic bundle $E$ on the disc $X$ equipped with the harmonic metric $h$ and the holomorphic Higgs field $\theta'_E$. The Higgs field takes the form $\theta^{\prime\std}_E+R(t)dt/t$, where $\theta^{\prime\std}_E$ is the Higgs field for the standard metric $h^{\std}$ as in the basic example of \T\ref{subsec:basic}. Recall that the $i\beta''/2$ are the eigenvalues of $\theta^{\prime\std}_E$. It is proved in \cite{Biquard97} that the perturbation $R(t)$ satisfies $R(0)=0$ (this argument was used in the proof of the assertion after \eqref{eq:matrix}). Moreover, the metric $h$ can be written as the product $h=h^{\std}(\id+v)$ where $v$ is a section of the H\"older space $C^{2+\vartheta}_\delta$ for any $\vartheta\in[0,1[$ (\cf \cite[p\ptbl77]{Biquard97}). Using an argument similar to that of \cite[Th\ptbl1 (Main estimate)]{Simpson90}, one gets:

\begin{enonce*}{\lemmname}[O\ptbl Biquard]
If $\beta''_i\neq\beta''_j$, then the component $v_{\beta_i,\beta_j}$ of $v$ and its logarithmic derivatives with respect to $t$ are $O(\mt^\eta)$ for some $\eta>0$.\qed
\end{enonce*}

\noindent
(An analogous statement also holds for general parabolic weights.)

To prove Property (iii), one argues then as follows. First, the operators $D''_E$ and $D''_{E^{\std}}$ coincide, as the holomorphic bundles $E$ and $E^{\std}$ coincide. Denote by $\varepsilong^{\std}$ an $h^{\std}$-orthonormal basis of $E$ and by $\varepsilong$ an $h$-orthonormal basis obtained from $\varepsilong^{\std}$ by the Gram-Schmidt process. Then the base change $\ccP$ from $\varepsilong^{\std}$ to $\varepsilong$ satisfies $\module{\ccP_{\beta_i,\beta_j}}=O(\mt^\eta)$ for some $\eta>0,$ if $\beta''_i\neq\beta''_j$, after the lemma. It follows that, in the basis $\varepsilong$, the matrix of $D''_E$ is obtained from the standard matrix by adding a perturbation term satisfying (i)--(iii) of the proof of Lemma \ref{lem:killing}. By adjunction, the same property holds for $D'_E$.

Consider now the Higgs field. As $R(t)$ is holomorphic and $R(0)=0$, we see, arguing as in the proof of the assertion after \eqref{eq:matrix} but in the reverse direction, that the matrix of $\theta'_E$ in the basis $\varepsilong^{\std}$ differs from that of $\theta^{\prime\std}_E$ by a perturbation term which is $O(\mt^\epsilon)$ for some $\epsilon>0$. Then, according to the lemma, in the basis $\varepsilong$ the matrix of $\theta'_E$ differs from the standard one \eqref{eq:theta'} by a perturbation term satisfying (i)--(iii). By adjunction, the same property holds for $\theta''_E$.\qed

\subsection{Elliptic regularity}\label{app:biquard2}
We give details on the argument of elliptic regularity used at the end of the proof of Lemma \ref{lem:killing}. Recall that $D_R$ denotes the disc of radius $R<1$ equipped with the Poincar\'e metric. We will now use the supplementary property that $P=P',P''$ satisfies, according to \cite{Biquard97}:
\begin{enumerate}
\item[(iv)]
$P$ belongs to the H\"older space $C^\vartheta_{1+\delta}$ for any $\vartheta\in[0,1[$, that is, $\Lt^{1+\delta}P$ belongs to the H\"older space $C^\vartheta(D^*_R)$, where the distance is taken with respect to the Poincar\'e metric.
\end{enumerate}
(In fact, this is also true for the first derivatives of $P$, and we can replace $1+\delta$ with $2+\delta$ for the components $P_{\beta_i,\beta_j}$ with $\beta_i\neq\beta_j$, but we will not use these properties.)

For each component $u^{(\gamma(\hb),\ell)}$ of the matrix $u$ obtained in the first part of Lemma \ref{lem:killing}, we have
\begin{equation}\label{eq:ugammaell}
d''u^{(\gamma(\hb),\ell)}+\gamma(\hb)u^{(\gamma(\hb),\ell)}\frac{d\ov t}{\ov t}+ \frac{\ell}{2}u^{(\gamma(\hb),\ell)}\frac{d\ov t}{\Lt\ov t}\in C^\vartheta_\delta\big(D_R,\rH(\Delta_{\hb_o}(\eta))\big),
\end{equation}
after (iv) and Formula \ref{eq:uStd}, taking into account that the norm of $d\ov t/\ov t$ is $\Lt$. Moreover, we know that $u^{(\gamma(\hb),\ell)}\in L^\infty_{\delta}\big(D_R,\rH(\Delta_{\hb_o}(\eta))\big)$ and that, if $\gamma(\hb)\not\equiv0$, that $u^{(\gamma(\hb),\ell)}\in L^\infty_{1+\delta}\big(D_R,\rH(\Delta_{\hb_o}(\eta))\big)$.

The idea in \loccit is to apply regularity properties of an elliptic differential operator locally in the upper half-plane with the hyperbolic metric and to use the homogeneity property of this operator with respect to the isometries of the hyperbolic plane to extend the corresponding inequalities to a fundamental domain covering the punctured disc $D^*_R$ (recall that we assume that $R<1$).

For $t\in D^*_R$, we put $t=e^{-\tau}$ with $\tau=s-i\theta$ and $s=\Lt$. The metric on the half-plane $\bbH=\{s>0\}$ is $(ds^2+d\theta^2)/s^2$. The covering $\bbH\to D^*_1$ is denoted by $q$. On a fixed hyperbolic ball $B_{a_o}(A_o)\subset\bbH$ of radius $a_o\leq1$ and center $(e^{A_o},0)$ with $A_o\in\RR$, we have, by a standard property of the Cauchy kernel for $(*)$ and by elliptic regularity of the operator $d''+\gamma(\hb)d\ov\tau$ for $(**)$ (see \eg \cite[Chap\ptbl XI]{Taylor81}), two inequalities for any $\vartheta\in{}]0,1[$,
$$
\norme{v}_{C^\vartheta(B_{a_o/2}(A_o))}\leq C_0(\vartheta)\big(\norme{(d''+\gamma(\hb)d\ov \tau)(v)}_{L^\infty(B_{a_o}(A_o))}+\norme{v}_{L^\infty(B_{a_o}(A_o))}\big)
\leqno(*)
$$
and
\begin{multline*}\tag*{$(**)$}
\norme{d'v}_{C^\vartheta(B_{a_o/4}(A_o))}\\
\leq C_1(\vartheta)\big(\norme{(d''+\gamma(\hb)d\ov \tau)(v)}_{C^\vartheta(B_{a_o/2}(A_o))}+\norme{v}_{C^\vartheta(B_{a_o/2}(A_o))}\big).
\end{multline*}
Here, we use that $s=\Lt$ and $s^{-1}$ are bounded on the fixed balls, so that $\norme{\cbbullet}_{C^\vartheta(B_{a_o}(A_o))}$ computed in the hyperbolic metric is  comparable with the same expression computed with the Euclidean metric, where the usual elliptic inequalities apply.

If $\varphi:\bbH\to\bbH$ is a hyperbolic holomorphic isometry, we have $\varphi^*(d''v+\gamma(\hb)vd\ov\tau)=(d''+\gamma(\hb)d\ov\tau)(v\circ\varphi)$. If we choose such an isometry sending the point $(e^{A_o},0)$ to $(e^{A},0)$, we obtain that the inequalities $(*)$ and $(**)$ also apply on the balls of radii $a_o/4,a_o/2,a_o$ centered at $(e^{A},0)$, with the \emph{same} constants $C_0(\vartheta)$ and $C_1(\vartheta)$.

Moreover, given $\delta>0$, there exist constants $c_1=c_1(\vartheta,a_o)$ and $c_2=c_2(\vartheta,a_o)$ such that, for any $A\in\RR$, we have
\[
c_1e^{\delta A}\norme{\cbbullet}_{C^\vartheta(B_{a_o}(A))}\leq \norme{\cbbullet}_{C^\vartheta_\delta(B_{a_o}(A))}\leq c_2e^{\delta A}\norme{\cbbullet}_{C^\vartheta(B_{a_o}(A))}
\]
(see \cite[p\ptbl54]{Biquard97}).

Choose now a sequence $A_n$ and $R',R''>0$ such that
\[
D^*_{R''}\subset q\Big(\tbigcup_n B_{a_o/4}(A_n)\Big)\subset D^*_{R'}\subset q\Big(\tbigcup_n B_{a_o}(A_n)\Big)\subset D^*_R.
\]
For a function $u$ on $D^*_R$ with values in $\rH(\Delta_{\hb_o}(\eta))$, we have
\begin{align*}
\norme{u}_{C^\vartheta_\delta(D^*_{R'})}&\leq c_2(\vartheta,a_o/2)\sup_n e^{\delta A_n}\norme{u\circ q}_{C^\vartheta(B_{a_o/2}(A_n)}\\
&\leq C_0(\vartheta)c_2(\vartheta,a_o/2)\sup_n e^{\delta A_n} \Big(\norme{(d''+\gamma(\hb)d\ov \tau)(u\circ q)}_{L^\infty(B_{a_o}(A_n))}\\[-4pt]
&\hspace*{7cm}+\norme{u\circ q}_{L^\infty(B_{a_o}(A_n))}\Big)\\
&\leq C'_0 \big(\norme{(d''+\gamma(\hb)d\ov \tau)(u)}_{L^\infty_\delta(D^*_R)}+\norme{u}_{L^\infty_\delta(D^*_R)}\big),
\end{align*}
with $C'_0=c_1(\vartheta,a_o)c_2(\vartheta,a_o/2)C_0(\vartheta)$. Similarly,
\[
\norme{d'u}_{C^\vartheta_\delta(D^*_{R''})}\leq C'_1 \big(\norme{(d''+\gamma(\hb)d\ov \tau)(u)}_{C^\vartheta_\delta(D^*_{R'})}+\norme{u}_{C^\vartheta_\delta(D^*_{R'})}\big).
\]

Coming back to $u^{(\gamma,\ell)}$, we have $\norme{u^{(\gamma,\ell)}}_{L^\infty_\delta(D^*_R)}<+\infty$ by construction and thus $\norme{(d''+\gamma(\hb)d\ov \tau)(u^{(\gamma,\ell)})}_{L^\infty_\delta(D^*_R)}<+\infty$ after the $L^\infty_\delta$ version of \eqref{eq:ugammaell} (recall that $\norme{d\ov t/\Lt\ov t}=1$). The first inequality above gives $\norme{u^{(\gamma,\ell)}}_{C^\vartheta_\delta(D^*_{R'})}<+\infty$ for any $\vartheta\in{}]0,1[$. Now, \eqref{eq:ugammaell} implies that $\norme{(d''+\gamma(\hb)d\ov \tau)(u^{(\gamma,\ell)})}_{C^\vartheta_\delta(D^*_{R'})}<+\infty$ and the second inequality above gives then $\norme{d'u^{(\gamma,\ell)}}_{C^\vartheta_\delta(D^*_{R''})}<+\infty$, so that, in particular, $\norme{d'u^{(\gamma,\ell)}}_{L^\infty_\delta(D^*_{R''})}<+\infty$, which is the desired statement.\qed

\let\thesection\oldthesection
\chapterspace{-2}
\chapter[The decomposition theorem]{The decomposition theorem\\ for polarizable regular twistor $\cD$-modules}\label{chap:decomp}

\section{Statement of the main results and proof of the Main Theorems}\label{sec:main}

\begin{theoreme}\label{th:imdirtwistor}\label{thimdirtwistor}
Let $f:X\to Y$ be a projective morphism between complex analytic manifolds and let $(\cT,\cS)$ be an object of $\MTr(X,w)^\rp$. Let $c$ be the first Chern class of a relatively ample line bundle on $X$ and let $\cL_c$ be the corresponding Lefschetz operator. Then $(\oplus_if_\dag^i\cT,\cL_c,\oplus_if_\dag^i\cS)$ is an object of $\MLTr(Y,w;1)^\rp$.
\end{theoreme}

\begin{remarque}
At the moment, Theorem \ref{th:imdirtwistor} is proved for \emph{regular} twistor $\cD$-modules only. Regularity is used in the proof of the case where $\dim X=1$ and $f$ is the constant map (case denoted by $\eqref{th:imdirtwistor}_{(1,0)}$ in \T\ref{sec:imdirtwistor10} below). The reduction to this case, done in \S\T\ref{sec:nmnp} and \ref{sec:n1n0}, does not use the regularity assumption. It seems reasonable to expect that the techniques of \cite{Bibi98} and \cite{B-B03} may be extended to obtain $\eqref{th:imdirtwistor}_{(1,0)}$, hence Theorem \ref{th:imdirtwistor}, in the nonregular case as well.
\end{remarque}

\begin{theoreme}\label{th:smoothtw}\label{thsmoothtw}
Let $X$ be a complex manifold and let $(\cT,\cS)$ be a smooth polarized twistor structure of weight~$w$ on $X$, in the sense of \T\ref{subsec:smtwqc}. Then $(\cT,\cS)$ is an object of $\MTr(X,w)^\rp$.
\end{theoreme}

\begin{proof}[Proof of Main Theorem \ref{thm:decomp} (see Introduction)]
Let $X$ be a smooth complex projective variety and let $\cF$ be an irreducible local system on $X$. Let $M=\cO_X\otimes_\CC\cF$ be the corresponding $\cD_X$-module. It follows from the theorem of K\ptbl Corlette and C\ptbl Simpson (\cf Lemma \ref{lem:twharm}) that $M$ underlies a smooth polarized twistor structure of weight~$0$ and, from Theorem \ref{th:smoothtw}, that $M=\Xi_{\DR}\cM''$ where $\cM''$ is the second term of an object $(\cT,\cS)$ of $\MTr(X,0)^\rp$.

Let $U$ be an open set of $X$ and let $f:U\to Y$ be a proper morphism to a complex analytic manifold $Y$. Then $f$ is projective. Any ample line bundle on $X$ will be relatively ample with respect to $f$. From Theorem \ref{th:imdirtwistor}, we conclude that $(\oplus_if_\dag^i\cT_{|U},\cL_c,\oplus_if_\dag^i\cS_{_U})$ is an object of $\MLTr(Y,0,1)^\rp$. In particular, each $f_\dag^i\cT_{|U}$ is strictly S-decomposable. By \cite{Deligne68}, we also know that $f_\dag\cT_{|U}\simeq\oplus_if_\dag^i\cT_{|U}$.

By restricting to $\hb=1$, we get \eqref{thm:decomp1} and \eqref{thm:decomp2} in Main Theorem \ref{thm:decomp}. We also get \eqref{thm:decomp3}, according to Proposition \ref{prop:restricanvar}.

Assume now that $U=X$ and $Y$ is projective. By the Lefschetz decomposition, we may apply Theorem \ref{th:ssimple} to conclude that $f_\dag^iM$ is semisimple as a regular holonomic $\cD_Y$-module.
\end{proof}

\begin{proof}[Proof of Main Theorem \ref{thm:vanishcycles} (see Introduction)]
The argument is similar. Starting with an irreducible local system $\cF$ on $X$, we may assume that it underlies a object of $\MTr(X,0)^\rp$. So does $\gr_\ell^{\rM}\psi_{f,\alpha}\cF$ by Proposition \ref{prop:restricanvar}, and we may apply Theorem \ref{th:ssimple} to get semisimplicity, hence Theorem \ref{thm:vanishcycles}.
\end{proof}

The remaining part of the chapter is devoted to the proof of Theorems \ref{th:imdirtwistor} and \ref{th:smoothtw}. We closely follow \cite[\T5.3]{MSaito86}. The proof of Theorem \ref{th:imdirtwistor} is by induction on the pair
$$(\dim\supp\cT,\dim\supp f(\cT)).$$
As $\eqref{th:imdirtwistor}_{(0,0)}$ is clear, it will be enough to prove
\begin{itemize}
\item
$\eqref{th:imdirtwistor}_{(1,0)}$ when $\supp\cT$ is smooth (Section \ref{sec:imdirtwistor10}),
\item
$\eqref{th:imdirtwistor}_{(n,m)} \implique \eqref{th:imdirtwistor}_{(n+1,m+1)}$ (Section \ref{sec:nmnp}),
\item
$\eqref{th:imdirtwistor}_{(\leq(n-1),0)}$ and ($\eqref{th:imdirtwistor}_{(1,0)}$ with $\supp\cT$ smooth) $\implique \eqref{th:imdirtwistor}_{(n,0)}$ for $n\geq1$ (Section \ref{sec:n1n0}).
\end{itemize}

\section{Proof of $\eqref{thimdirtwistor}_{(1,0)}$ when $\supp\cT$ is smooth} \label{sec:imdirtwistor10}
Let $X$ be a compact Riemann surface and let $f:X\to\mathrm{pt}$ be the constant map. We will argue as in \cite{Zucker79}. Let $(\cT,\cS)=(\cM,\cM,C,\id)$ be a polarized twistor left $\cD_X$-module of weight~$0$ on $X$ (\cf Remark \ref{rem:hermdual}). For the proof of $\eqref{th:imdirtwistor}_{(1,0)}$, it will be enough to assume that it has strict support $X$. Denote by $\wt\cM$ the localization of~$\cM$ with respect to the singular set $P$. As $(\cT,\cS)_{|X^*}$ is a smooth twistor structure, it corresponds to a flat bundle with harmonic metric $h$, and we have a metric $\pi^*h$ on $\wt\cM_{\cX^*}$ (\cf \T\ref{subsec:constr}).

Fix a complete metric on $X^*$ which is equivalent to the Poincar\'e metric near each puncture $x_o\in P$, with volume form $\vol$. Extend it as a metric on the bundle $\Theta_{\cX^*}$, and therefore on $\cE_{\cX^*}^1$ so that $\norme{dt/\hb}=\norme{d\ov t}=\mt\Lt$. On the other hand, put on $\cX^*$ the product metric (Poincar\'e on $X^*$, Euclidean on $\Omega_0$) to compute the volume form $\cvol$. Recall that, in a local coordinate $t$ near a puncture, we have
\begin{equation}\label{eq:L2Poincare}
\module{t}^{a}\Lt^{w/2}\in L^2(\vol)\iff a>0\text{ or } a=0\text{ and } w\leq0.
\end{equation}

Put $M_{\hb_o}=\cM/(\hb-\hb_o)\cM$.

\subsection{The holomorphic $L^2$ de~Rham complex}
\label{subsec:holomDR}
Denote by $\wt\cM_{(2)}\subset\wt\cM$ the submodule consisting of local sections which are locally $L^2$ with respect to $\pi^*h$. Formula \eqref{eq:norme}, together with \eqref{eq:L2Poincare}, shows that (by switching from $\beta$ to $\alpha$) a local section $m=\sum_jm_jt^{n_j}e_j^{(\hb_o)}$ of $\wt\cM$ is in $\wt\cM_{(2), \hb_o}$ if and only if, for any $j$, either $\ell_{\hb_o}(n_j+\beta_j)>0$ (and therefore $\ell_{\hb}(n_j+\beta_j)>0$ for any $\hb$ near $\hb_o$) or $n_j+\beta_j=-1$ (that is, $\ell_{\hb}(n_j+\beta_j)=-1$ for any $\hb$ near $\hb_o$) and $w_j\leq0$. According to Remark \ref{rem:although}\eqref{rem:although2}, this is equivalent to saying that $m$ is a local section of $V^{(\hb_o)}_{-1}\wt\cM$, the class of which in $\gr_{-1}^{V^{(\hb_o)}}\wt\cM$ is contained in $\rM_0\psi_{t,-1}\wt\cM$, where $\rM_\bbullet$ denotes the monodromy filtration of $\rN$. In particular, we have a natural inclusion $\wt\cM_{(2)}\subset\cM$, globally defined with respect to~$\hb$.

Similarly, the sheaf $(\Omega_{\cX}^1\otimes_{\cO_X}\wt\cM)_{(2)}$ of $L^2$ local sections of $\Omega_{\cX}^1\otimes_{\cO_X}\wt\cM$ consists, near $\hb_o$, of local sections $(dt/t\hb)\otimes m$, where $m$ is a local section of $V^{(\hb_o)}_{-1}\wt\cM$ the class of which in $\gr_{-1}^{V^{(\hb_o)}}\wt\cM$ is contained in $\rM_{-2}\psi_{t,-1}\wt\cM$. We also have a natural inclusion $(\Omega_{\cX}^1\otimes_{\cO_X}\wt\cM)_{(2)}\subset \Omega_{\cX}^1\otimes_{\cO_X}\cM$, as follows from the following lemma:

\begin{lemme}\label{lem:yetseen}
Let $m$ be a local section of $V^{(\hb_o)}_{-1}\wt\cM$, the class of which in $\gr_{-1}^{V^{(\hb_o)}}\wt\cM$ is contained in $\rM_{-2}\gr_{-1}^{V^{(\hb_o)}}\wt\cM$. Then $m/t$ is a local section of $\cM$.
\end{lemme}

\begin{proof}
We denote here by $\rM_\bbullet\gr_{-1}^{V^{(\hb_o)}}\wt\cM$ the direct sum of the monodromy filtrations of $\rN$ on each $\psi_{t,\alpha}\wt\cM$ with $\ell_{\hb_o}(\alpha)=-1$. We want to show that $m/t=\partiall_tn_1+n_2$, where $n_1$ is a local section of $V^{(\hb_o)}_{-1}\wt\cM$ and $n_2$ a local section of $V^{(\hb_o)}_{<0}\wt\cM$. As $t:V^{(\hb_o)}_{<0}\wt\cM\to V^{(\hb_o)}_{<-1}\wt\cM$ is bijective, this is equivalent to finding $n_1$ as above and a local section $n_3$ of $V^{(\hb_o)}_{<-1}\wt\cM$ such that $m=t\partiall_tn_1+n_3$. Equivalently, the class of $m$ in $\gr_{-1}^{V^{(\hb_o)}}\wt\cM$ has to be contained in $\image t\partiall_t$. For the component in $\psi_{t,-1}\wt\cM$, this is equivalent to saying that it is contained in $\image(\rN)$. Now, that $\rM_{-2}\psi_{t,-1}\wt\cM$ is contained in $\image(\rN)$ is a general property of the monodromy filtration: indeed, we have $\rM_{-2}=\rN(\rM_0)$, as $\rN:\gr_k^{\rM}\to\gr_{k-2}^{\rM}$ is onto for any $k\leq0$. For the component in any $\psi_{t,\alpha}\wt\cM$ with $\alpha\neq-1$, notice that $t\partiall_t:\psi_{t,\alpha}\wt\cM\to\psi_{t,\alpha}\wt\cM$ is an isomorphism, by Lemma \ref{lem:imhb}, and therefore sends bijectively $\rM_{-2}\psi_{t,\alpha}\wt\cM$ to itself.
\end{proof}

The description above clearly shows that the connection $\nabla$ on $\wt\cM$ sends $\wt\cM_{(2)}$ into $(\Omega_{\cX}^1\otimes_{\cO_X}\wt\cM)_{(2)}$ and defines the \emph{holomorphic $L^2$ de~Rham complex} $\DR(\wt\cM)_{(2)}$. There is a natural inclusion of complexes
\[
\index{$dr2$@$\DR(\wt\cM)_{(2)}$}\DR(\wt\cM)_{(2)}\hto\DR\cM.
\]

\begin{remarque}[Restriction to $\hb=\hb_o$]\label{rem:restrhboL2}
Note that $\wt M_{\hb_o}=\wt\cM/(\hb-\hb_o)\wt\cM$ is also the localization of $M_{\hb_o}$ along $P$, as $\cO_{\cX}[*P]$ is $\cO_{\cX}$-flat. A similar computation can be made for the subsheaf of local sections of $\wt M_{\hb_o}$ the norm of which is $L^2$ with respect to the metric $h$ and the volume form $\vol$. Denoting by $U_\bbullet M_{\hb_o}$ the filtration induced by $V_\bbullet^{(\hb_o)}\cM$ (\cf Proposition \ref{prop:restrhbo}), this subsheaf is equal, according to Formula \ref{eq:normehbo} in Remark \ref{rem:although}\eqref{rem:although3}, to $\rM_0U_{-1}M_{\hb_o}$, which is the inverse image in $U_{-1}M_{\hb_o}$ of $\rM_0\gr_{-1}^{U}M_{\hb_o}$. Similarly, $(\Omega_X^1\otimes_{\cO_X}\wt M_{\hb_o})_{(2)}= \Omega_X^1\langle\log P\rangle\otimes \rM_{-2}U_{-1}M_{\hb_o}$.

Notice however that the derivative of a section of $\rM_0U_{-1}M_{\hb_o}$ may not be $L^2$. If we denote by $\wt M_{\hb_o,(2)}$ the subsheaf of $L^2$ sections of $\wt M_{\hb_o}$ with $L^2$ derivative, we identify $\wt M_{\hb_o,(2)}$ with the subsheaf of local sections of $\rM_0U_{-1}M_{\hb_o}$ the class of which in $\rM_0\gr_{-1}^UM_{\hb_o}$ have $\alpha$-components in $\rM_{-2}\psi_{t,\alpha}^{(\hb_o)}M_{\hb_o}$ for any $\alpha\neq-1$ such that $\ell_{\hb_o}(\alpha)=-1$.

According to Proposition \ref{prop:restricanvar}\eqref{prop:restricanvar2a} and the previous remarks, we have a natural inclusion of complexes
\[
\DR(\wt\cM)_{(2)}/(\hb-\hb_o)\DR(\wt\cM)_{(2)}\hto\index{$dr2hb$@$\DR(\wt M_{\hb_o})_{(2)}$}\DR(\wt M_{\hb_o})_{(2)}.
\]
Moreover, for the same reason, the holomorphic $L^2$ de~Rham complex $\DR(\wt M_{\hb_o})_{(2)}$ (or the Dolbeault complex if \hbox{$\hb_o=0$}) is a subcomplex of $\DR M_{\hb_o}$: indeed, it is clear by the same Proposition that the term in degree $-1$ is the restriction to $\hb=\hb_o$ of a submodule of $\cM$; for the term in degree $0$, use Lemma \ref{lem:yetseen}.
\end{remarque}

\begin{Proposition}\label{prop:DRL2}
\begin{enumerate}
\item\label{prop:DRL21}
The natural inclusions of complexes $\DR(\wt\cM)_{(2)}\hto\DR\cM$ and, for any $\hb_o$, $\DR(\wt M_{\hb_o})_{(2)}\hto\DR M_{\hb_o}$, are quasi-isomorphisms.
\item\label{prop:DRL22}
The natural inclusion of complexes
\[
\DR(\wt\cM)_{(2)}/(\hb-\hb_o)\DR(\wt\cM)_{(2)}\hto\DR(\wt M_{\hb_o})_{(2)}
\]
is a quasi-isomorphism.
\end{enumerate}
\end{Proposition}

\begin{proof}
As the morphisms are globally defined with respect to~$\hb$, the question is local near any $\hb_o$. It is also enough to work locally near a singular point of $\cM$ in $X$, with a local coordinate~$t$. We will therefore use the $V_\bbullet^{(\hb_o)}$-filtration.

Let us begin with the first part of \eqref{prop:DRL21}. We then have
$$
\DR\cM=\{0\to V_{-1}^{(\hb_o)}\cM\To{-\nabla} \Omega^1_\cX\otimes V_0^{(\hb_o)}\cM\to 0\},
$$
because locally this complex is nothing but
\[
0\to V_{-1}^{(\hb_o)}\cM\To{-\partiall_t}V_0^{(\hb_o)}\cM\to 0,
\]
and for any $a>-1$, the morphism $\partiall_t:\gr_a^{V^{(\hb_o)}}\cM\to\gr_{a+1}^{V^{(\hb_o)}}\cM$ is an isomorphism (\cf Definition \ref{def:strictspe}\eqref{def:strictspec}).

As $\cM$ is contained in $\wt\cM$ (being its minimal extension across $t=0$), we have isomorphisms $t:V_{<0}^{(\hb_o)}\cM\isom V_{<-1}^{(\hb_o)}\cM$ and $t:V_0^{(\hb_o)}\cM\isom tV_0^{(\hb_o)}\cM\subset V_{-1}^{(\hb_o)}\cM$ and, by Formula \ref{prop:minext*}, we have $V_0^{(\hb_o)}\cM=\partiall_t V_{-1}^{(\hb_o)}\cM+V_{<0}^{(\hb_o)}\cM$. Therefore, the previous complex is quasi-isomorphic to the complex
\[
0\to V_{-1}^{(\hb_o)}\cM\To{-t\partiall_t}t\partiall_tV_{-1}^{(\hb_o)}\cM+ V_{<-1}^{(\hb_o)}\cM\to 0.
\]
We therefore have an exact sequence of complexes:
\[
0\to\DR(\wt\cM)_{(2)}\to\DR\cM\to C^{\cbbullet}\to0,
\]
where $C^{\cbbullet}$ is the complex
\begin{multline*}
0\to\Big(\ooplus_{\substack{\alpha\neq-1\\\ell_{\hb_o}(\alpha)=-1}}\psi_{\alpha}^{(\hb_o)}\cM\Big)\oplus\big(\psi_{-1}^{(\hb_o)}\cM/\rM_0\psi_{-1}^{(\hb_o)}\cM\big)\\
\To{-t\partiall_t}\Big(\ooplus_{\substack{\alpha\neq-1\\\ell_{\hb_o}(\alpha)=-1}}t\partiall_t\psi_{\alpha}^{(\hb_o)}\cM\Big)\oplus\big(\rN\psi_{-1}^{(\hb_o)}\cM/\rM_{-2}\psi_{-1}^{(\hb_o)}\cM\big)\to0
\end{multline*}
In order to prove the first statement of Proposition \ref{prop:DRL2}, we are left with showing that
{\def\theenumi{\alph{enumi}}
\begin{enumerate}
\item\label{prop:DRL2a}
$t\partiall_t:\psi_{\alpha}^{(\hb_o)}\cM\to t\partiall_t\psi_{\alpha}^{(\hb_o)}\cM$ is an isomorphism for any $\alpha\neq-1$ such that $\ell_{\hb_o}(\alpha)=\nobreak-1$,
\item\label{prop:DRL2b}
the morphism
\[
\rN:\big(\psi_{-1}^{(\hb_o)}\cM/\rM_0\psi_{-1}^{(\hb_o)}\cM\big)\to\big(\rN\psi_{-1}^{(\hb_o)}\cM/\rM_{-2}\psi_{-1}^{(\hb_o)}\cM\big)
\]
is an isomorphism.
\end{enumerate}}

For \eqref{prop:DRL2a}, use Lemma \ref{lem:imhb} to conclude that $t\partiall_t$ is invertible on the corresponding $\psi_\alpha$'s.

For \eqref{prop:DRL2b}, the surjectivity is clear. The injectivity follows from the injectivity of $\rN:\gr_k^{\rM}\to\gr_{k-2}^{\rM}$ for $k\geq1$ and the equality $\rM_{-2}=\rN(\rM_0)$ that we have yet seen in the proof of Lemma \ref{lem:yetseen}.

The second part of \ref{prop:DRL2}\eqref{prop:DRL21} is proved similarly: consider the complex $C^{\prime\cbbullet}$ analogous to $C^{\cbbullet}$, where one replaces $\psi_{\alpha}^{(\hb_o)}\cM$ with $\psi_{\alpha}^{(\hb_o)}\cM/\rM_{-2}\psi_{\alpha}^{(\hb_o)}\cM$ and $t\partiall_t\psi_{\alpha}^{(\hb_o)}\cM$ with $t\partiall_t\psi_{\alpha}^{(\hb_o)}\cM/\rM_{-2}\psi_{\alpha}^{(\hb_o)}\cM$; then, as we have seen in Remark \ref{rem:restrhboL2}, the analogue of \eqref{prop:DRL2a} remains true for $C^{\prime\cbbullet}$; by restriction to $\hb=\hb_o$ one gets the desired statement.

Now, \ref{prop:DRL2}\eqref{prop:DRL22} is obtained from \ref{prop:DRL2}\eqref{prop:DRL21} by restricting to $\hb=\hb_o$.
\end{proof}

\subsection{The $L^2$ complex}\label{subsec:L2hb}
Keep notation of \T\ref{subsec:Hodgetw}. We will define the $L^2$ complex with $\hb_o$ fixed, and then for $\hb$ varying near $\hb_o\in\Omega_0$. Recall that we put $\ccD_\hb=\ccD_0+\hb\ccD_\infty$ and $\ccD_{\hb_o}=\ccD_0+\hb_o\ccD_\infty$.

Fix $\hb_o\in\Omega_0$. Let us recall the standard definition of the $L^2$ complex $\cLd^{\cbbullet}(H,\ccD_{\hb_o})$ on $X$. Denote by $j:X^*\hto X$ the inclusion, and by $L^1_{\loc,X^*}$ the sheaf on $X^*$ of $L^1_{\loc}$ functions. Then $\cLd^j(H,\ccD_{\hb_o})$ is the subsheaf of $j_*\big[L^1_{\loc,X^*}\otimes_{\cC^\infty_{X^*}}\cE_{X^*}^j\otimes H\big]$ defined by the following properties: a section $v$ of the previous sheaf is a section of $\index{$l2hd$@$\cLd^\bbullet(H,\ccD_{\hb_o})$}\cLd^j(H,\ccD_{\hb_o})$ iff
\begin{itemize}
\item
the norm $\norme{v_{|X^*}}$ of $v_{|X^*}$ with respect to the metric $h$ on $H$ and the metric on $\cE_{X^*}^j$ induced by the Poincar\'e metric on $X^*$---which is a local section of $j_*L^1_{\loc,X^*}$---is a section of $L^2_{\loc,X}(\vol)$,
\item
$\ccD_{\hb_o} v_{|X^*}$, which is taken in the distributional sense, is a section of $L^1_{\loc,X^*}\otimes_{\cC^\infty_{X^*}}\cE_{X^*}^{j+1}\otimes H$, and $\norme{\ccD_{\hb_o} v_{|X^*}}$ is a local section of $L^2_{\loc, X}(\vol)$.
\end{itemize}
The $L^2$ condition may be read on the coefficients of $v_{|X^*}$ in any $h$-orthonormal basis of $H$ or, more generally, in any $L^2$-adapted basis in the sense of Zucker. According to \cite[Lemma 4.5]{Zucker79}, the restriction to $\hb=\hb_o$ of the basis $\varepsilong^{\prime\prime(\hb_o)}$ introduced in Formula \eqref{eq:bme} is $L^2$-adapted, and therefore, as $\bme^{(\hb_o)}$ is obtained from $\varepsilong^{\prime\prime(\hb_o)}$ by a rescaling (the matrix $\wt A$ in \eqref{eq:bme} being diagonal), it is also $L^2$-adapted as well as the basis $\bme^{\prime(\hb_o)}$ defined in \eqref{eq:bmeprime}.

\medskip
We now consider the case where $\hb$ belongs to some neighbourhood of $\hb_o$. Let $U$ be an open set of $X$, put $U^*=U\cap X^*$. In the following, we only consider functions on open sets of $\cX$, the restriction of which to $\cX^*$ is $L^1_{\loc}$ with respect to the volume form $\cvol$.

Say that such a function $\varphi$ on some open set $U\times\Omega\subset\cX$ is holomorphic with respect to $\hb$ if its restriction to $U^*\times\Omega$ satisfies $\ov\partial_\hb\varphi=0$.

Say that a local section $v$ of $j_*\big(L^1_{\loc,\cX^*}\otimes_{\cC^\infty_{\cX^*}}\pi^{-1}\cE_X\otimes\pi^{-1}H\big)$ is a section of $\index{$l2dhcurl$@$\cLd^{\cbbullet}(\cH,\ccD_\hb)$}\cLd^{\cbbullet}(\cH,\ccD_\hb)$ if
\begin{enumerate}
\item\label{cond:L21}
$v_{|\cX^*}$ is holomorphic with respect to $\hb$, \ie in any local basis of $H$ the coefficients $v_j$ of $v_{|\cX^*}$ satisfy $\ov\partial_\hb v_j=0$; it is therefore meaningful to consider $v(\cbbullet,\hb)_{|X^*}$;
\item\label{cond:L22}
the norm $\norme{v(\cbbullet,\hb)_{|X^*}}$ of $v(\cbbullet,\hb)_{|X^*}$ (with respect to the metric $h$ on $H$ and the metric on $\cE_{X^*}^j$ induced by the Poincar\'e metric on $X^*$) is $L^2$ locally uniformly in $\hb$, \ie on any open set $U\times\Omega$ on which $v$ is defined, and any compact set $K\subset\Omega$, we have $\sup_{\hb\in K}\norme{v(\cbbullet,\hb)}_{2,\vol}<+\infty$;
\item\label{cond:L23}
$\ccD_\hb v_{|\cX^*}$, which is taken in the distributional sense, takes values in $L^1_{\loc}$-sections and satisfies \eqref{cond:L22} (and clearly \eqref{cond:L21}, as $\ccD_\hb$ commutes with $\ov\partial_\hb$).
\end{enumerate}
The $L^2$ condition may be read on the coefficients in any $L^2$-adapted basis. Formulas \eqref{eq:bme} and \eqref{eq:bmeprime} show that, near any $(x,\hb_o)$ with $x\in P$ (the singular set of $\cM$ in~$X$), the bases $\bme^{(\hb_o)}$ and $\bme^{\prime(\hb_o)}$ are $L^2$-adapted.

In particular, for any $\hb_o\in\Omega_0$, there is a restriction morphism $$\cLd^{\cbbullet}(\cH,\ccD_\hb)_{|\hb=\hb_o}\to\cLd^{\cbbullet}(H,\ccD_{\hb_o}).$$

\begin{theoreme}\label{th:L2}
Let $(\cT,\cS)=(\cM,\cM,C,\id)$ be a polarized twistor left $\cD_X$-module of weight~$0$ on $X$. Then the natural inclusions of complexes
\[
\DR(\wt\cM)_{(2)}\to\cLd^{1+\cbbullet}(\cH,\ccD_\hb)
\]
and, for any $\hb_o\in\Omega_o$,
\[
\DR(\wt M_{\hb_o})_{(2)}\to\cLd^{1+\cbbullet}(H,\ccD_{\hb_o})
\]
are quasi-isomorphisms compatible with restriction to $\hb=\hb_o$.
\end{theoreme}

We will use the rescaling \eqref{eq:rescale} to define the inclusion, as explained in \T\ref{subsec:endthL2}. Let us begin with the analysis of the $L^2$ complexes $\cLd^{\cbbullet}(\cH,\ccD_\hb)$ and $\cLd^{\cbbullet}(H,\ccD_{\hb_o})$. We will use the basis $\bme^{\prime(\hb_o)}$ or its restriction to $\hb=\hb_o$ that we denote in the same way. The norm of any element $e^{\prime(\hb_o)}_{\beta,\ell,k}$ in the subfamily $\bme^{\prime(\hb_o)}_{\beta,\ell}$ ($\beta\in B$, $\ell\in\ZZ$) is $\module{t}^{\ell_{\hb}(q_{\beta,\imhb_o}+\beta)}\Lt^{\ell/2}(1+o(1))$, where $o(1)=\Lt^{-\delta}$ for some $\delta>0$. Recall that $q_{\beta,\imhb_o}$ is chosen so that $\ell_{\hb_o}(q_{\beta,\imhb_o}+\beta)\in[0,1[$ (\cf proof of Corollary \ref{cor:biquard}). The proof of Theorem \ref{th:L2} will occupy \S\T\ref{subsec:Hardy}--\ref{subsec:endthL2}. It will be given for $\hb$ variable. The compatibility with restriction to $\hb=\hb_o$ will be clear.

\begin{Remarques}
\begin{enumerate}
\item
We will not follow here the proof given by Zucker (a dichotomy Dolbeault/Poincar\'e lemma) for two reasons: the first one is that we do not know how to adapt the proof for the Poincar\'e lemma near points $\hb_o\in\Sing\Lambda$; the second one is that we want to put a parameter $\hb$ in the proof, and have a proof which is ``continuous with respect to $\hb$'', in particular near $\hb=0$. However, the proof of the Poincar\'e lemma given by Zucker does not ``degenerate'' to the proof of the Dolbeault lemma.
\item
Notice that $\DR(\wt M_{\hb_o})_{(2)}=\DR M_{\hb_o}$ has cohomology in degree $-1$ only if $\hb_o\not\in\Sing\Lambda$, as a consequence of Proposition \ref{prop:restricanvar}\eqref{prop:restricanvar1}, but this may not be true if $\hb_o\in\Sing\Lambda$.
\end{enumerate}
\end{Remarques}

We put $\ccD'_\hb=\hb D'_E+\theta'_E$ and $\ccD''_\hb=D''_E+\hb\theta''_E$, so that $\ccD_\hb=\ccD'_\hb+\ccD''_\hb$. Recall that the basis $\bme^{\prime(\hb_o)}$ is holomorphic with respect to $\ccD''_\hb$, and that $\ccD'_\hb$ acts by $\hb d'+\Theta'_{\hb}$, where the matrix $\Theta'_{\hb}$ is defined by \eqref{eq:matrix}. Taking the notation of \eqref{eq:matrix} we will decompose $\Theta'_{\hb}$~as
\[
\Theta'_{\hb}=\Theta'_{\hb,\diag}+\Theta'_{\hb,\nilp},\quad\text{with}\quad
\begin{cases}
\Theta'_{\hb,\diag}=\oplus_\beta(q_{\beta,\imhb_o}+\beta)\star\hb\id\,\dfrac{dt}{t}\\
\Theta'_{\hb,\nilp}=\big[\rY+P(t,\hb)\big]\dfrac{dt}{t},
\end{cases}
\]
with $\rY=(\oplus_\beta \rY_\beta)$. We will view $\rY$ and $P$ as operators and not matrices. We index the basis $\bme^{\prime(\hb_o)}$ by $\beta,\ell,k$, where $\ell$ denotes the weight with respect to $\rY$ and $k$ is used to distinguish the various elements having the same $\beta,\ell$. Then we have, for any $\beta,\ell,k$, $\rY(e^{\prime(\hb_o)}_{\beta,\ell,k})=e^{\prime(\hb_o)}_{\beta,\ell-2,k}$ and, according to the assertion after \eqref{eq:matrix},
\begin{equation}\label{eq:Ptilde}
\begin{aligned}
P(e^{\prime(\hb_o)}_{\beta,\ell,k})&=\sum_{\substack{\gamma\neq\beta\\ \ell_{\hb_o}(q_{\gamma,\imhb_o}+\gamma)\leq\ell_{\hb_o}(q_{\beta,\imhb_o}+\beta)}}\hspace*{-1cm}\sum_\lambda\sum_\kappa tp_{\beta,\gamma,\ell,\lambda,k,\kappa}(t,\hb)e^{\prime(\hb_o)}_{\gamma,\lambda,\kappa}\\
&+\sum_{\substack{\gamma\neq\beta\\ \ell_{\hb_o}(q_{\gamma,\imhb_o}+\gamma)>\ell_{\hb_o}(q_{\beta,\imhb_o}+\beta)}}\hspace*{-1cm}\sum_\lambda\sum_\kappa p_{\beta,\gamma,\ell,\lambda,k,\kappa}(t,\hb)e^{\prime(\hb_o)}_{\gamma,\lambda,\kappa}\\
&+
\sum_{\lambda\geq\ell-2}\sum_\kappa tp_{\beta,\ell,\lambda,k,\kappa}(t,\hb)e^{\prime(\hb_o)}_{\beta,\lambda,\kappa}+
\sum_{\lambda\leq\ell-3}\sum_\kappa p_{\beta,\ell,\lambda,k,\kappa}(t,\hb)e^{\prime(\hb_o)}_{\beta,\lambda,\kappa},
\end{aligned}
\end{equation}
where the functions $p_{\beta,\gamma,\ell,\lambda,k,\kappa}(t,\hb),p_{\beta,\ell,\lambda,k,\kappa}(t,\hb)$ are holomorphic. It follows that the $L^2$ norm of $\Theta'_{\hb,\nilp}$ is bounded. As a consequence, in the basis $\bme^{\prime(\hb_o)}$, the $L^2$ condition on derivatives under $\ccD_\hb$ can be replaced with an $L^2$ condition on derivatives with respect to the diagonal operator
\begin{equation}\label{eq:diagop}
\ccD_{\hb,\diag}\defin d''+\hb d'+\Theta'_{\hb,\diag}.
\end{equation}

\subsection{Hardy inequalities}\label{subsec:Hardy}
We will recall here the basic results that we will need concerning Hardy inequalities in $L^2$.

Let $I={}]A,B[{}\subset\RR$ be a nonempty open interval and $v,w$ two functions (weights) on $I$, which are measurable and almost everywhere positive and finite. Let $u$ be a $C^1$ function on $I$.

\begin{theoreme}[$L^2$ Hardy inequalities, \cf \eg {\cite[Th\ptbl1.14]{O-K90}}]
There is an inequality
\[
\big\Vert u\cdot w\big\Vert_2\leq C\big\Vert u'\cdot v \big\Vert_2,
\]
with
\[
C=
\begin{cases}
\dpl\sup_{x\in I}\int_x^Bw(t)^2\,dt\cdot\int_A^xv(t)^{-2}\,dt&\text{if $\dpl\lim_{x\to A_+}u(x)=0$,}\\[8pt]
\dpl\sup_{x\in I}\int_A^xw(t)^2\,dt\cdot\int_x^Bv(t)^{-2}\,dt&\text{if $\dpl\lim_{x\to B_-}u(x)=0$.}
\end{cases}
\]
\qed
\end{theoreme}

\begin{corollaire}\label{cor:Hardy}
Take $A=0$ and $B=R>0$, with $R<1$. Let $(b,k)\in\RR\times\ZZ$. Assume that $(b,k)\neq(0,1)$. Given $g(r)$ continuous on $[0,R]$, put
\[
f(r)=\begin{cases}
\dpl\int_0^rg(\rho)\,d\rho&\text{if $b<0$ or if $b=0$ and $k\geq2$,}\\[8pt]
\dpl\int_{\min(R,e^{-k/2b})}^rg(\rho)\,d\rho&\text{if $b>0$ or if $b=0$ and $k\leq0$.}
\end{cases}
\]
(In the second case, we set $e^{-k/2b}=+\infty$ if $b=0$ and $k\leq0$; this is the limit case of $b>0$ and $k\leq0$ when $b\to0_+$.)

Then there exists a constant $C(R,b,k)>0$, such that $\lim_{b\to0_+}C(R,b,k)<+\infty$ when $k\neq1$ and $\lim_{b\to0_-}C(R,b,k)<+\infty$ when $k\geq2$, such that the following inequality holds (we consider $L^2([0,R];dr/r)$ norms)
\begin{align*}
\big\Vert f(r)\cdot r^b\rL(r)^{k/2-1}\big\Vert_{2,dr/r}&
\leq C(R,b,k) \big\Vert g(r)\cdot r^b\rL(r)^{k/2-1}\cdot r\rL(r)\big\Vert_{2,dr/r}\\
&= C \big\Vert g(r)\cdot r^{b+1}\rL(r)^{k/2}\big\Vert_{2,dr/r}.
\end{align*}
\end{corollaire}

\begin{proof}
We will choose the following weight functions with respect to the measure $dr$: $w(r)=r^{b-1/2}\rL(r)^{k/2-1}$ and $v(r)=r^{b+1/2}\rL(r)^{k/2}$.
Assume first $b>0$.
\begin{enumerate}
\item
If $R\leq e^{-k/2b}$, \ie $k/2b\leq\rL(R)$, we have thus $\lim_{r\to R_-}f(r)=0$. We will show the finiteness of
\[
\sup_{r\in[0,R]}\int_0^r \rho^{2b-1}\rL(\rho)^{k-2}\,d\rho\cdot \int_r^R\rho^{-2b-1}\rL(\rho)^{-k}\,d\rho.
\]
After the change of variable $y=\rL(\rho)$ we have to estimate
\[
\int_x^{+\infty}e^{-2by}y^k\,\frac{dy}{y^2}\cdot\int_{\rL(R)}^xe^{2by}y^{-k}\,dy.
\]
The function $e^{-2by}y^k$ is decreasing on $]\rL(R),+\infty[$, hence the first integral is bounded by $e^{-2bx}x^{k-1}$, and the second one by $e^{2bx}x^{-k}(x-\rL(R))$. Here, we can therefore take $C(R,b,k)=1$.

\item
Assume now $k/2b>\rL(R)>0$, and $k\neq1$, hence $k\geq2$. We have
\[
f(r)=\int_{e^{-k/2b}}^rg(\rho)\,d\rho.
\]
We will have to consider the two cases $r\in{}]0,e^{-k/2b}[$ and $r\in{}]e^{-k/2b},R[$.
\begin{enumerate}
\item
If $r\in{}]0,e^{-k/2b}[$, we have to bound the same expression as above, in the same conditions, after replacing $\rL(R)$ with $k/2b$, and we can therefore choose $C=1$.
\item
If $r\in{}]e^{-k/2b},R[$, we want to estimate
\[
\int_{\rL(R)}^xe^{-2by}y^{k-2}\,dy\cdot \int_x^{k/2b}e^{2by}y^{-k}\,dy.
\]
The function $e^{-2by}y^k$ is increasing, and an estimation analogous to the previous one would give a constant $C$ having a bad behaviour when $b\to0_+$. We will thus use that $e^{-2by}$ is decreasing. As $k\geq2$, we have, for the second integral,
\[
\int_x^{k/2b}e^{2by}y^{-k}\,dy\leq \frac{e^k}{k-1}\,(x^{1-k}-(k/2b)^{1-k})\leq \frac{e^kx^{1-k}}{k-1}.
\]
For the first one we have
\[
\int_{\rL(R)}^xe^{-2by}y^{k-2}\,dy\leq e^{-2b\rL(R)}\frac{(x^{k-1}-\rL(R)^{k-1})}{k-1}=R^{2b}\frac{(x^{k-1}-\rL(R)^{k-1})}{k-1},
\]
so that the product is bounded by
\[
\frac{e^kR^{2b}}{(k-1)^2}\big(1-(\rL(R)/x)^{k-1}\big)\leq\frac{e^kR^{2b}}{(k-1)^2}\defin C(R,b,k),
\]
which has a good behaviour when $b\to0_+$.
\end{enumerate}
\item
If $k/2b>\rL(R)>0$, and $k=1$, the trick of part (b) above does not apply, because the function $\log$ is increasing, though the function $x^{1-k}$ was decreasing if $k\geq2$. The constant that we get does not have a good behaviour when $b\to0_+$, which could be expected, as we do not have a good Hardy inequality if $b=0$ and $k=1$.
\end{enumerate}

When $b=0$ and $k\neq1$, the proof above degenerates to give the corresponding inequalities.

Consider now the case where $b<0$.
\begin{enumerate}
\item
Assume that $e^{(k-2)/2\module{b}}\geq R$, \ie $k\geq2(1-\module{b}\rL(R))$, which is satisfied in particular whenever $k\geq2$. Then the function $e^{-2by}y^{k-2}$ is increasing on $]\rL(R),+\infty[$. An upper bound of
\[
\int_{\rL(R)}^xe^{-2by}y^{k-2}\,dy\cdot\int_x^{+\infty}e^{2by}y^{2-k}\,\frac{dy}{y^2}
\]
is given by
\[
(x-\rL(R))e^{-2bx}x^{k-2}\cdot e^{2bx}x^{-k+2}x^{-1}=(1-\rL(R)/x)\leq1.
\]
This situation degenerates when $b\to0_-$, to the analogous situation when $b=0$ and $k\geq2$.
\item
If $e^{(k-2)/2\module{b}}< R$, \ie $k<2(1-\module{b}\rL(R))$
we find a constant with a bad behaviour when $b\to0_-$.\qedhere
\end{enumerate}
\end{proof}

\subsection{Proof of Theorem \ref{th:L2}: vanishing of $\cH^2$} \label{subsec:vanishingH2}

As in \cite{Zucker79}, we will express $L^2$ functions of $t=re^{i\theta}$, holomorphic in $\hb$, as Fourier series $\sum_{n\in\ZZ}u_n(r,\hb)e^{in\theta}$, where each $u_n$ is holomorphic with respect to $\hb$ when $\hb$ is not fixed. The $L^2$ condition will be recalled below.

\begin{lemme}[Image of $\ccD''_\hb$]\label{lem:imdbar}
Let $\varphi=\sum_{\beta,\ell,k}\varphi_{\beta,\ell,k}(t,\hb)e^{\prime(\hb_o)}_{\beta,\ell,k}\,\frac{dt}{t}\wedge\frac{d\ov t}{\ov t}$ be a local section of $\cLd^{2}(\cH,\ccD_\hb)$, with $\varphi_{\beta,\ell,k}(t,\hb)=\sum_{n\in\ZZ}\varphi_{\beta,\ell,k,n}(r,\hb)e^{in\theta}$. Then, if all coefficients $\varphi_{\beta,\ell,k,0}(r,\hb)$ vanish identically for any $\beta,\ell$ such that $\ell_{\hb_o}(q_{\beta,\imhb_o}+\beta)=0$ and $\ell\leq-1$, there exists a local section $\psi$ of $\cLd^{(1,0)}(\cH,\ccD_\hb)$ such that $\varphi=\ccD_\hb\psi=\ccD''_\hb\psi$.

The same assertion holds after restriction to $\hb=\hb_o$.
\end{lemme}

Let us now be more precise on the $L^2$ condition: $f(t,\hb)e^{\prime(\hb_o)}_{\beta,\ell,k}$ is a local section of $\cLd^{0}(\cH)$ iff $\module{f(t,\hb)}\module{t}^{\ell_{\hb}(q_{\beta,\imhb_o}+\beta)}\Lt^{\ell/2}\Lt^{-1}$ is in $L^2(d\theta\,dr/r)$, locally uniformly with respect to $\hb$. A similar statement holds for sections in $\cLd^{1}$ or $\cLd^{2}$ expressed with logarithmic forms $dt/t$ and $d\ov t/\ov t$, using that the norm of each of these form is $\Lt$.

\begin{proof}[Proof of Lemma \ref{lem:imdbar}]
We have to solve (up to sign) $\ov t\partial_{\ov t} \psi_{\beta,\ell,k}=\varphi_{\beta,\ell,k}$, if $\psi_{\beta,\ell,k}$ is the coefficient of $e^{\prime(\hb_o)}_{\beta,\ell,k}\,\sfrac{dt}{t}$ in $\psi$. We will argue as in \cite[Propositions 6.4 and 11.5]{Zucker79}. Put $u=\psi_{\beta,\ell,k}$ and $v=\varphi_{\beta,\ell,k}$, and $u=\sum_nu_n(r,\hb)e^{in\theta}$, $v=\sum_nv_n(r,\hb)e^{in\theta}$. We are then reduced to solving for any $n\in\ZZ$
\[
\frac{\partial}{\partial r}(r^{-n}u_n(r,\hb))=2r^{-n-1}v_n(r,\hb)
\]
with an $L^2$ condition. If we put $f_n(r,\hb)=r^{-n}u_n(r,\hb)$ and $g_n(r,\hb)=2r^{-n-1}v_n(r,\hb)$, this condition reads
\[
\norme{f_n(r,\hb)r^{n+\ell_\hb(q_{\beta,\imhb_o}+\beta)}\rL(r)^{\ell/2}}_{2,dr/r}\leq C \norme{g_n(r,\hb)r^{n+1+\ell_\hb(q_{\beta,\imhb_o}+\beta)}\rL(r)^{\ell/2+1}}_{2,dr/r},
\]
uniformly for $\hb$ near $\hb_o$. If $n+\ell_{\hb_o}(q_{\beta,\imhb_o}+\beta)\neq0$, this remains true for $\hb$ near $\hb_o$ and we apply Corollary \ref{cor:Hardy}. If $n+\ell_{\hb_o}(q_{\beta,\imhb_o}+\beta)=0$, we thus have $\ell_{\hb_o}(q_{\beta,\imhb_o}+\beta)=0$ and $n=0$, as $\ell_{\hb_o}(q_{\beta,\imhb_o}+\beta)\in[0,1[$. If moreover $\ell\geq0$, we also apply Corollary \ref{cor:Hardy}, as the constant can be chosen uniformly with respect to $\hb$ for $\hb$ near $\hb_o$ (notice that, if $\ell_{\hb}(q_{\beta,\imhb_o}+\beta)\not \equiv0$, the sign of this function changes at $\hb=\hb_o$, hence the condition on $\ell$).

If we consider the case where $\hb=\hb_o$ is fixed, then the only condition is $\ell\neq-1$ when $\ell_{\hb_o}(q_{\beta,\imhb_o}+\beta)=0$.
\end{proof}

We will now show that the $L^2$ complexes $\cLd^{\cbbullet}(\cH,\ccD_\hb)$ and $\cLd^{\cbbullet}(H,\ccD_{\hb_o})$ have no $\cH^2$. By Lemma \ref{lem:imdbar} we are reduced to showing that, for any $\beta$ with \hbox{$\ell_{\hb_o}(q_{\beta,\imhb_o}+\beta)=0$}, any section $f(r,\hb)e'_{\beta,\ell,k}\,\frac{dt}{t}\wedge\frac{d\ov t}{\ov t}$ of $\cLd^{2}(\cH)$ with $\ell\leq-1$ belongs to the image of~$\ccD_\hb$. We will distinguish two cases.
\begin{enumerate}
\item
If $(q_{\beta,\imhb_o}+\beta)\star\hb_0\neq0$, we remark that
\[
d''\Big(f(r,\hb)\,\frac{dt}{t}\Big)=-\ov t\partial_{\ov t}(f)\,\frac{dt}{t}\wedge\frac{d\ov t}{\ov t}= -\tfrac12r\partial_r(f)\,\frac{dt}{t}\wedge\frac{d\ov t}{\ov t}
\]
and a similar computation for $d'(f(r,\hb)\sfrac{d\ov t}{\ov t})$, so that
\[
(\hb d'+d'')\Big(f(r,\hb)\hb\frac{dt}{t}+f(r,\hb)\frac{d\ov t}{\ov t}\Big)=0.
\]
We will prove the assertion by increasing induction on $\ell$. It is true for $\ell\ll0$, as the section is the equal to $0$. For arbitrary $\ell\leq-1$, it is enough to show that $\Theta'_{\hb,\nilp}\big(f(r,\hb)e'_{\beta,\ell,k}\sfrac{d\ov t}{\ov t}\big)$ belongs to $\im\ccD_\hb$. This is true by induction for the part $f(r,\hb)\rY(e'_{\beta,\ell,k})\,\frac{dt}{t}\wedge\frac{d\ov t}{\ov t}$. Consider now the image under $P(t,\hb)\sfrac{dt}{t}$ and its coefficient on $e'_{\gamma,\lambda,\kappa}\,\frac{dt}{t}\wedge\frac{d\ov t}{\ov t}$. By \eqref{eq:Ptilde},
\begin{itemize}
\item
if $\gamma\neq\beta$, the coefficient is $p(t,\hb)f(r,\hb)$ if $\ell_{\hb_o}(q_{\gamma,\imhb_o}+\gamma)>\ell_{\hb_o}(q_{\beta,\imhb_o}+\beta)=0$, and $tp(t,\hb)f(r,\hb)$ if $\ell_{\hb_o}(q_{\gamma,\imhb_o}+\gamma)\leq\ell_{\hb_o}(q_{\beta,\imhb_o}+\beta)=0$, for some holomorphic function $p(t,\hb)$. In the first case, we have $\ell_{\hb_o}(q_{\gamma,\imhb_o}+\gamma)\in{}]0,1[$, so we apply Lemma \ref{lem:imdbar}. In the second case, we have $\ell_{\hb_o}(q_{\gamma,\imhb_o}+\gamma)=0$ and the coefficient in the Fourier series of $tp(t,\hb)f(r,\hb)$ corresponding to $n=0$ is zero, so we can apply the same lemma.
\item
If $\gamma=\beta$, the coefficient is $p(t,\hb)f(r,\hb)$ if $\lambda\leq\ell-3$ and $tp(t,\hb)f(r,\hb)$ if $\lambda\geq\ell-2$. In the second case, the coefficient in the Fourier series of $tp(t,\hb)f(r,\hb)$ corresponding to $n=0$ is zero, so we can apply Lemma \ref{lem:imdbar}. In the first case, we apply the same argument to $\big(p(t,\hb)-p(0,\hb)\big)f(r,\hb)$, and the inductive assumption for $p(0,\hb)f(r,\hb)$.
\end{itemize}

\item
If $(q_{\beta,\imhb_o}+\beta)\star\hb_0=0$, knowing that $\ell_{\hb_o}(q_{\beta,\imhb_o}+\beta)=0$ it follows from Lemma \ref{lem:imhb} that we have in fact $q_{\beta,\imhb_o}+\beta=0$, and therefore $q_{\beta,\imhb_o}=0$ and $\beta=0$ by definition of $q_{\beta,\imhb_o}$, hence $(q_{\beta,\imhb_o}+\beta)\star\hb\equiv0$. As we also have $\ell_{\hb}(q_{\beta,\imhb_o}+\beta)\equiv0$, we can reduce to $\ell=-1$ (see the proof of Lemma \ref{lem:imdbar}). We write
\begin{multline*}
f(r,\hb)e'_{0,-1,k}\,\frac{dt}{t}\wedge\frac{d\ov t}{\ov t}=\ccD_\hb\Big(f(r,\hb)e'_{0,1,k}\hb\frac{dt}{t}+f(r,\hb)e'_{0,1,k}\frac{d\ov t}{\ov t}\Big)\\
-f(r,\hb)P(t,\hb)e'_{0,1,k}\,\frac{dt}{t}\wedge\frac{d\ov t}{\ov t},
\end{multline*}
as $\rY e'_{0,1,k}=e'_{0,-1,k}$. It is therefore enough to show that \hbox{$f(r,\hb)P(t,\hb)e'_{0,1,k}\,\frac{dt}{t}\wedge\frac{d\ov t}{\ov t}$} belongs to the image of $\ccD_\hb$, which is proven as above using Lemma \ref{lem:imdbar} and \eqref{eq:Ptilde}.\qed
\end{enumerate}

\subsection{Proof of Theorem \ref{th:L2}: computation of $\cH^1$} \label{subsec:computationH1}

By the result of \T\ref{subsec:vanishingH2}, the $L^2$ complex is quasi-isomorphic to the complex
\begin{equation}\label{eq:complexered}
0\to\cLd^0(\cH,\ccD_\hb)\To{\ccD_\hb}\ker\ccD_\hb^1\to0.
\end{equation}

\begin{lemme}\label{lem:supprdtb}
Any local section $\psi\,dt/t+\varphi\,d\ov t/\ov t$ of $\ker\ccD^1_\hb\subset\cLd^1(\cH,\ccD_\hb)$ can be written as the sum of a term in $\im\ccD_\hb$ and a term in $\cLd^{(1,0)}(\cH)\cap\ker\ccD_\hb^1$.
\end{lemme}

\begin{proof}
Write
\[
\varphi=\sum_{\beta,\ell,k,n}\varphi_{\beta,\ell,k,n}(r,\hb)e^{in\theta}\cdot e^{\prime(\hb_o)}_{\beta,\ell,k},
\]
and put $\varphi_{(0,0)}$ (\resp $\varphi_{\neq(0,0)}$) the sum of terms for which $\ell_{\hb_o}(q_{\beta,\imhb_o}+\beta)=0$ and $n=0$ (\resp the sum of the other terms).
We claim that there exists a section $\eta_{\neq(0,0)}$ of $\cLd^0(\cH,\ccD_\hb)$ such that $\ccD''_\hb\eta_{\neq(0,0)}=\varphi_{\neq(0,0)}d\ov t/\ov t$. First, the existence of $\eta_{\neq(0,0)}$ in $\cLd^0(\cH,\ccD''_\hb)$, \ie without taking care of $\ccD'_\hb\eta_{\neq(0,0)}$, is obtained as in Lemma \ref{lem:imdbar}. We wish to show that $\ccD'_\hb\eta_{\neq(0,0)}$ belongs to $\cLd^{(1,0)}(\cH)$ or, as we have seen, that $\ccD'_{\hb,\diag}\eta_{\neq(0,0)}$ belongs to $\cLd^{(1,0)}(\cH)$. It is therefore enough to show that the ${\neq}(0,0)$-part of $\ccD'_\hb\eta_{\neq(0,0)}$ belongs to $\cLd^{(1,0)}(\cH)$. We have
\begin{align*}
\ccD''_\hb\big(\ccD'_\hb\eta_{\neq(0,0)}\big)&=-\ccD'_\hb\big(\ccD''_\hb\eta_{\neq(0,0)}\big)=-\ccD'_\hb\Big(\varphi_{\neq(0,0)}\frac{d\ov t}{\ov t}\Big)\\
&=\ccD'_\hb\Big((\varphi_{(0,0)}-\varphi)\frac{d\ov t}{\ov t}\Big)=\ccD''_\hb\Big(\psi\frac{dt}{t}\Big)+\ccD'_\hb\Big(\varphi_{(0,0)}\frac{d\ov t}{\ov t}\Big)\\
&=\ccD''_\hb\Big(\psi_{\neq(0,0)}\frac{dt}{t}\Big)+\Big(\Theta'_{\hb,\nilp}\varphi_{(0,0)}\frac{d\ov t}{\ov t}\Big)+(0,0)\text{-terms}.
\end{align*}
By considering the ${\neq}(0,0)$-part we get
\[
\ccD''_\hb\big(\ccD'_\hb\eta_{\neq(0,0)}\big)_{\neq(0,0)}=\ccD''_\hb\Big(\psi_{\neq(0,0)}\frac{dt}{t}\Big)+\Big(\Theta'_{\hb,\nilp}\varphi_{(0,0)}\frac{d\ov t}{\ov t}\Big)_{\neq(0,0)}.
\]
According to Lemma \ref{lem:imdbar}, the second term of the right-hand side is a section of $\ccD''_\hb\big(\cLd^{(1,0)}(\cH,\ccD''_\hb)\big)$. We conclude that there exists $\nu\in\cLd^{(1,0)}(\cH,\ccD''_\hb)$ such that
\[
\big(\ccD'_\hb\eta_{\neq(0,0)}\big)_{\neq(0,0)}=\nu+\omega
\]
with $\omega=\sum_{\beta,\ell,k}\omega_{\beta,\ell,k}e^{\prime(\hb_o)}_{\beta,\ell,k}\,dt$ and each $\omega_{\beta,\ell,k}$ is a distribution such that $\partial_{\ov t}\omega_{\beta,\ell,k}=0$, \ie $\ccD''_\hb\omega=0$; hence each $\omega_{\beta,\ell,k}$ is a holomorphic function. Notice now that each $\omega_{\beta,\ell,k}e^{\prime(\hb_o)}_{\beta,\ell,k}\,dt$ belongs to $\cLd^{(1,0)}(\cH,\ccD''_\hb)$. Indeed, the $L^2$ condition is that
\[
\module{t\omega_{\beta,\ell,k}}\big\Vert e^{\prime(\hb_o)}_{\beta,\ell,k}\big\Vert=r^{\ell_\hb(q_{\beta,\imhb_o}+\beta)+1}\rL(r)^{\ell/2}\module{\omega_{\beta,\ell,k}}\in L^2_{\loc}(d\theta\,dr/r),
\]
which is clearly satisfied as
\[
\int_0^Rr^{2\ell_\hb(q_{\beta,\imhb_o}+\beta)+2}\rL(r)^\ell\,\frac{dr}{r}<+\infty
\]
if $\hb$ is sufficiently close to $\hb_o$ so that $\ell_\hb(q_{\beta,\imhb_o}+\beta)>-1$ (recall that $\ell_{\hb_o}(q_{\beta,\imhb_o}+\beta)\in[0,1[$). We can now conclude that $\ccD'_\hb\eta_{\neq(0,0)}$ is a section of $\cLd^{(1,0)}(\cH)$ and that
\[
\psi\,\frac{dt}{t}+\varphi\,\frac{d\ov t}{\ov t}=\ccD_\hb\eta_{\neq(0,0)}+\wt \psi\,\frac{dt}{t}+\varphi_{(0,0)}\,\frac{d\ov t}{\ov t};
\]
in other words, we are reduced to the case where $\varphi=\varphi_{(0,0)}$.

Let $\gamma$ be such that $\ell_{\hb_o}(q_{\gamma,\imhb_o}+\gamma)=0$. We claim that, for any $\lambda\in\ZZ$ and $\kappa$, the coefficient $(P\cdot \varphi)_{\gamma,\lambda,\kappa,0}$ is a linear combination of terms $\varphi_{\gamma,\ell,k,0}$ with $\ell\geq\lambda+3$, coefficients being holomorphic functions of $\hb$ in a neighbourhood of $\hb_o$.

Indeed, by \eqref{eq:Ptilde} and by the assumption $\varphi=\varphi_{(0,0)}$, we have
\[
\big(P\cdot \varphi_{\beta,\ell,k,0}e^{\prime(\hb_o)}_{\beta,\ell,k}\big)_{\gamma,\lambda,\kappa}=
\begin{cases}
tp(t,\hb)\varphi_{\beta,\ell,k,0}&\text{if }\beta\neq\gamma\text{ or } \beta=\gamma\text{ and }\ell\leq\lambda+2\\
p(t,\hb)\varphi_{\beta,\ell,k,0}&\text{if }\beta=\gamma\text{ and } \ell\geq\lambda+3,
\end{cases}
\]
for some holomorphic function $p(t,\hb)$ (depending on $\beta,\gamma,\ell,\lambda,k,\kappa$). Therefore, a nonzero coefficient not depending on $\theta$ in the Fourier expansion appears only in the second case, by taking $p(0,\hb)\varphi_{\beta,\ell,k,0}$.

Consider now the component on $e^{\prime(\hb_o)}_{\gamma,\lambda,\kappa}\,\frac{dt}{t}\wedge\frac{d\ov t}{\ov t}$ of the relation $\ccD''_\hb(\psi\,dt/t)+\ccD'_\hb(\varphi\,d\ov t/\ov t)=0$, when $\ell_{\hb_o}(q_{\gamma,\imhb_o}+\gamma)=0$. We get
\begin{multline}\label{eq:0-1k0}
\frac12r\partial_r\psi_{\gamma,\lambda,\kappa,0}(r,\hb)\\
=(\frac{\hb}{2}r\partial_r+\gamma\star\hb)(\varphi_{\gamma,\lambda,\kappa,0}(r,\hb))+\varphi_{\gamma,\lambda+2,\kappa,0}(r,\hb)+(P\cdot \varphi)_{\gamma,\lambda,\kappa,0}.
\end{multline}
We know, by Lemma \ref{lem:imhb}, that either $\gamma\star\hb_o\neq0$, so $\gamma\star\hb$ is invertible for $\hb$ near $\hb_o$, or $\gamma=0$ and hence $\gamma\star\hb\equiv0$. We will show by decreasing induction on $\lambda$ that
\[
\varphi_{\gamma,\lambda,\kappa,0}(r,\hb)e^{\prime(\hb_o)}_{\gamma,\lambda,\kappa}\,\frac{d\ov t}{\ov t}= \ccD''_\hb\big(\eta_{\gamma,\lambda,\kappa,0}(r,\hb)e^{\prime(\hb_o)}_{\gamma,\lambda,\kappa}\big)=\frac12r\partial_r\big(\eta_{\gamma,\lambda,\kappa,0}(r,\hb)\big)e^{\prime(\hb_o)}_{\gamma,\lambda,\kappa}\,\frac{d\ov t}{\ov t}
\]
for some section $\eta_{\gamma,\lambda,\kappa,0}(r,\hb)e^{\prime(\hb_o)}_{\gamma,\lambda,\kappa}$ of $\cLd^0(\cH,\ccD_\hb)$. This will be enough to conclude the proof of the lemma.

Assume first that $\gamma\star\hb_o\neq0$. Then \eqref{eq:0-1k0} allows us to write
\[
(\gamma\star\hb)\varphi_{\gamma,\lambda,\kappa,0}=\frac12r\partial_r\Big[\psi_{\gamma,\lambda,\kappa,0}-\hb\varphi_{\gamma,\lambda,\kappa,0}-\sum_{m\geq2}c_m(\hb)\eta_{\gamma,\lambda+m,\kappa,0}\Big],
\]
with $c_m(\hb)$ holomorphic near $\hb_o$ and $c_2\equiv1$. Denote by $(\gamma\star\hb)\eta_{\gamma,\lambda,\kappa,0}$ the term between brackets. The $L^2$ condition to be satisfied is
\begin{equation}\label{eq:etagamma}
\big\vert\eta_{\gamma,\lambda,\kappa,0}(r,\hb)\big\vert r^{\ell_\hb(q_{\gamma,\imhb_o}+\gamma)}\rL(r)^{\lambda/2-1}\in L^2_{\loc}(dr/r)
\end{equation}
and we have by assumption
\begin{align*}
\big\vert\psi_{\gamma,\lambda,\kappa,0}(r,\hb)\big\vert r^{\ell_\hb(q_{\gamma,\imhb_o}+\gamma)}\rL(r)^{\lambda/2}&\in L^2_{\loc}(dr/r)\\
\big\vert\varphi_{\gamma,\lambda,\kappa,0}(r,\hb)\big\vert r^{\ell_\hb(q_{\gamma,\imhb_o}+\gamma)}\rL(r)^{\lambda/2}&\in L^2_{\loc}(dr/r)\\
\big\vert\eta_{\gamma,\lambda+m,\kappa,0}(r,\hb)\big\vert r^{\ell_\hb(q_{\gamma,\imhb_o}+\gamma)}\rL(r)^{(\lambda+m)/2-1}&\in L^2_{\loc}(dr/r)\quad\text{for any }m\geq2.
\end{align*}
Clearly, these conditions imply \eqref{eq:etagamma}. Moreover, we have
\begin{align*}
\ccD'_{\hb,\diag}\big(\eta_{\gamma,\lambda,\kappa,0}e^{\prime(\hb_o)}_{\gamma,\lambda,\kappa}\big)&=\Big(\frac{\hb}2r\partial_r+\gamma\star\hb\Big)(\eta_{\gamma,\lambda,\kappa,0})e^{\prime(\hb_o)}_{\gamma,\lambda,\kappa}\,\frac{dt}{t}\\
&=\Big[\psi_{\gamma,\lambda,\kappa,0}-\sum_{m\geq2}c_m(\hb)\eta_{\gamma,\lambda+m,\kappa,0}\Big]e^{\prime(\hb_o)}_{\gamma,\lambda,\kappa}\,\frac{dt}{t},
\end{align*}
and this satisfies the $L^2$ condition as $m\geq2$.

Consider now the case where $\gamma=0$. Then \eqref{eq:0-1k0} allows us to write, with a small change of notation,
\[
\varphi_{0,\lambda,\kappa,0}=\frac12r\partial_r\Big[\psi_{0,\lambda-2,\kappa,0}-\hb\varphi_{0,\lambda-2,\kappa,0}-\sum_{m\geq1}c_m(\hb)\eta_{0,\lambda+m,\kappa,0}\Big],
\]
Denote by $\eta_{0,\lambda,\kappa,0}$ the term between brackets. We have by assumption
\begin{align*}
\big\vert\psi_{0,\lambda-2,\kappa,0}(r,\hb)\big\vert \rL(r)^{(\lambda-2)/2}&\in L^2_{\loc}(dr/r)\\
\big\vert\varphi_{0,\lambda-2,\kappa,0}(r,\hb)\big\vert \rL(r)^{(\lambda-2)/2}&\in L^2_{\loc}(dr/r)\\
\big\vert\eta_{0,\lambda+m,\kappa,0}(r,\hb)\big\vert \rL(r)^{(\lambda+m)/2-1}&\in L^2_{\loc}(dr/r)\quad\text{for any }m\geq1.
\end{align*}
Clearly, these conditions imply the $L^2$ condition for $\eta_{0,\lambda,\kappa,0}e^{\prime(\hb_o)}_{0,\lambda,\kappa}$, namely
\[
\big\vert\eta_{0,\lambda,\kappa,0}(r,\hb)\big\vert \rL(r)^{\lambda/2-1}\in L^2_{\loc}(dr/r).
\]
Moreover, we have
\[
\ccD'_{\hb,\diag}\big(\eta_{0,\lambda,\kappa,0}e^{\prime(\hb_o)}_{0,\lambda,\kappa}\big)=\frac12r\partial_r\big(\eta_{0,\lambda,\kappa,0}\big)e^{\prime(\hb_o)}_{0,\lambda,\kappa}\,\frac{dt}{t}=\varphi_{0,\lambda,\kappa,0}e^{\prime(\hb_o)}_{0,\lambda,\kappa}\,\frac{dt}{t},
\]
which is $L^2$ by assumption.
\end{proof}

\subsection{End of the proof of Theorem \ref{th:L2}}\label{subsec:endthL2}
First, it easily follows from Lemma \ref{lem:supprdtb} that the natural inclusion of the complex
\[
0\to\ker\ccD_\hb^{\prime\prime0}\To{\ccD^{\prime0}_\hb}\cLd^{(1,0)}(\cH)\cap\ker\ccD^1_\hb\to0
\]
in the complex \eqref{eq:complexered} is a quasi-isomorphism. Notice that $\cLd^{(1,0)}(\cH)\cap\ker\ccD^1_\hb=\cLd^{(1,0)}(\cH)\cap\ker\ccD^{\prime\prime1}_\hb$. A section of each of the sheaves in this complex is therefore holomorphic away from $t=0$, by the usual Dolbeault-Grothendieck lemma. The $L^2$ condition implies that the coefficients in the bases $\bme^{\prime(\hb_o)}$ or $\bme^{\prime(\hb_o)}\,dt/t$ have moderate growth at \hbox{$t=0$}, hence are meromorphic. The complex $\DR(\wt\cM)_{(2)}$ is therefore isomorphic to the previous complex shifted by $1$, the morphism being given by the identity on the term of degree $-1$ and by the multiplication by $\hb$ (the rescaling) on the term of degree~$0$.\qed

\Subsection{End of the proof of Theorem \ref{thimdirtwistor} on a Riemann surface}

\subsubsection*{Strictness of $f_\dag\cT$}
Let $(\cT,\cS)=(\cM,\cM, C,\id)$ be a polarized twistor regular $\cD_X$-module of weight~$0$, as in Theorem \ref{th:imdirtwistor}.

\begin{corollaire}\label{cor:strict}
The complex $f_\dag\cM$ is strict.
\end{corollaire}

\begin{proof}
Argue as in the proof of Hodge-Simpson Theorem~\ref{th:HodgeSimpson}. For any $\hb_o\in\Omega_0$, $f_\dag M_{\hb_o}$ has finite dimensional cohomology. It follows from Proposition \ref{prop:DRL2}\eqref{prop:DRL21} that each $\DR(\wt M_{\hb_o})_{(2)}$ has finite dimensional hypercohomology on $X$ and, by Theorem \ref{th:L2}, so does $\cLd^{1+\bbullet}(H,\ccD_{\hb_o})$. One may therefore apply Hodge Theory to the complete manifold $X^*$ and the Laplace operator $\Delta_{\hb_o}$. As this operator is essentially independent of $\hb_o$ by \eqref{eq:Deltao}, $\hb_o$-harmonic sections are $\hb$-harmonic and $\ccD_{\hb}$-closed for any $\hb$. Moreover, they have finite dimension.

As $f_\dag\cM$ has $\cO_{\Omega_0}$-coherent cohomology, we may now conclude with the same argument as in the smooth case (see the part of the proof of the Hodge-Simpson Theorem \ref{th:HodgeSimpson} concerning strictness).
\end{proof}

\begin{Remarques}
\begin{enumerate}
\item
Let $M$ be an \emph{irreducible} holonomic $\cD_X$-module with regular singularities. We have $\DR M_{|X^*}=L[1]$, where $L$ is an \emph{irreducible} local system on $X^*$, and $\DR M=j_*L[1]$ is the intersection complex attached to the shifted local system $\DR M_{X^*}=L[1]$, if $j:X^*=X\moins P\hto X$ denotes the inclusion. In particular, $\DR M$ has cohomology in degree $-1$ only. Consequently, if $M\neq\cO_X$, \ie if $L\neq\CC_{X^*}$, we have
\[
\bH^{-1}(X,\DR M)=H^0(X,j_*L)=0,
\]
and, by Poincar\'e duality and using that the dual $M^\vee$ is also irreducible, we have $\bH^1(X,\DR M)=0$. Therefore, $f_\dag M=f_\dag^0M$.

\item
Let now $\cT=(\cM,\cM,C,\id)$ be a polarized twistor regular $\cD_X$-module of weight~$0$. If we assume that $\cT$ is simple and has strict support $X$, and if $\cT$ is not equal to the twistor $\cD$-module associated to $\cO_\cX$, then $f_\dag\cT=f_\dag^0\cT$: indeed, by Remark \ref{rem:ssimple}, we know that the restriction $M$ of $\cM$ at $\hb_o=1$ is simple and not equal to $\cO_X$; then, by the remark above, the restriction $f_\dag M$ at $\hb_o=1$ of $f_\dag\cM$ has cohomology in degree $0$ only; as $f_\dag\cM$ is strict, it must have cohomology in degree $0$ only. Notice also that all sections of $f_\dag^0\cM$ are primitive with respect to the Lefschetz morphism associated to any $C^\infty$ metric on $X$.
\end{enumerate}
\end{Remarques}

\subsubsection*{The twistor condition}
We can assume that $\cT$ is simple and has strict support~$X$, and we also assume that $\cT$ is not equal to the twistor $\cD$-module associated to $\cO_\cX$ (otherwise the result is clear). We want to prove that the twistor condition is satisfied for $f_\dag\cT=f_\dag^0\cT$. First, the harmonic sections $\Harm^1(X,H)$ with respect to any $\Delta_{\hb_o}$ form a lattice in $\bH^0\big(X,\DR(\wt\cM)_{(2)}\big)$: this is proved as in \T\ref{subsec:Hodgetw}, replacing the $C^\infty$ de~Rham or Dolbeault complex with the $L^2$ complex, and using Theorem \ref{th:L2}. Then the polarization property is proved as in \T\ref{subsec:Hodgetw}, with the simplification that all harmonic sections are primitive.\qed

\section{Proof of $\eqref{thimdirtwistor}_{(n,m)} \implique \eqref{thimdirtwistor}_{(n+1,m+1)}$}\label{sec:nmnp}
Let $f:X\to Y$ be a projective morphism between complex manifolds and let $(\cT,\cS)=(\cM,\cM,C,\id)$ be an object of $\MTr(X,0)^p$ (it is easy to reduce to the case of weight~$0$ by a Tate twist). Assume that it has strict support a closed irreducible analytic set $Z\subset X$. Put $n+1=\dim Z$ and $m+1=\dim f(Z)$. Fix a relatively ample line bundle on $X$ and let $c$ be its Chern class.

We know by Corollary \ref{cor:holo} that the $f^i_\dag\cM$ are holonomic. Let $t$ be a holomorphic function on an open set $V\subset Y$ and put $g=t\circ f:f^{-1}(V)\to\CC$. We assume that $\{t=0\}\cap f(Z)$ has everywhere codimension one in $f(Z)$. We will now show that the $f^i_\dag\cM$ are strictly S-decomposable along $t=0$, by proving that the other conditions for a twistor object are also satisfied.

Consider first $\oplus_if^i_\dag\Psi_{g,\alpha}\cT$ for $\reel(\alpha)\in[-1,0[$ (\resp $\oplus_if^i_\dag\phi_{g,0}\cT$) with its nilpotent endomorphism $f_\dag^i\cN$ and monodromy filtration $\rM_\bbullet(f_\dag^i\cN)$, and the nilpotent Lefschetz endomorphism $\cL_c$.

\begin{claim}\label{claim}
For any $\alpha$ with $\reel\alpha\in[-1,0[$, the object
\[
\ooplus_{i,\ell}\big[\gr_\ell^{\rM}f_\dag^i\Psi_{g,\alpha}(\cT), \gr_\ell^{\rM}f_\dag^i\Psi_{g,\alpha}(\cS), f_\dag^i\cN,\cL_c\big]
\]
is an object of $\MLTr(V,w;-1,1)^\rp$. Similarly,
\[
\ooplus_{i,\ell}\big[\gr_\ell^{\rM}f_\dag^i\phi_{g,0}(\cT), \gr_\ell^{\rM}f_\dag^i\phi_{g,0}(\cS), f_\dag^i\cN,\cL_c\big]
\]
is an object of $\MLTr(V,w;-1,1)^\rp$.
\end{claim}

\begin{proof}[Sketch of proof]
By the inductive assumption for $\alpha\neq0$ and using Corollary \ref{cor:psi0} if $\alpha=0$, we know that
\[
\ooplus_{i,\ell}\big[f_\dag^i\gr_\ell^{\rM}\Psi_{g,\alpha}(\cT), f_\dag^i\gr_\ell^{\rM}\Psi_{g,\alpha}(\cS), f_\dag^i\cN,\cL_c\big]
\]
(\resp\ ...) is an object of $\MLTr(V,w;-1,1)^\rp$. Then we may apply Corollary \ref{cor:suitespecdeg}.
\end{proof}

As a consequence, we get the strictness of $f_\dag^i\Psi_{g,\alpha}\cM$ for any $\alpha$ with $\reel\alpha\in[-1,0[$, hence that of $f_\dag^i\psi_{g,\alpha}\cM$ for any $\alpha\not\in\NN$, as this is a local property with respect to~$\hb$. Similarly, we get the strictness of $f_\dag^i\psi_{g,0}\cM$, hence that of $f_\dag^i\psi_{g,\alpha}\cM$ for any \hbox{$\alpha\in\NN$}.

\enlargethispage{5mm}%
Applying Theorem \ref{th:imdirstrictspe}, we conclude that the $f^i_\dag\cM$ are regular and strictly specializable along $t=0$. By Corollary \ref{cor:imdirstrictspeRT}, we have $f_\dag^i\Psi_{g,\alpha}\cT=\Psi_{t,\alpha}f_\dag^i\cT$ for any $\alpha$ with $\reel\alpha\in[-1,0[$ and $f_\dag^i\phi_{g,0}\cT=\phi_{t,0}f_\dag^i\cT$.

Now, Condition $(\MLT_{>0})$ along $t=0$ is satisfied for $\oplus_{i}f^{i}_{\dag}\cT$, because of the claim. According to Remark \ref{rem:canvaradj}, strict S-decomposability along $t=0$ follows then from Proposition \ref{prop:declefcoh}, \ie the analogue of Proposition \ref{prop:declef}: indeed, as $f_\dag^i\cN$ commutes with $\cL_c$, it is enough to prove the S-decomposability of the primitive (relative to $\cL_c$) modules $Pf_\dag^i\cM$; apply then Proposition \ref{prop:declefcoh} to the objects $\gr^{\rM}\psi_{g,-1}Pf_\dag^i(\cT,\cS)=\gr^{\rM}Pf_\dag^i\psi_{g,-1}(\cT,\cS)$ and $\gr^{\rM}\phi_{g,0}Pf_\dag^i(\cT,\cS)=\gr^{\rM}Pf_\dag^i\phi_{g,0}(\cT,\cS)$, which are polarized graded Lefschetz objects according to the claim, to get that $\gr^{\rM}\phi_{g,0}Pf_\dag^i\cT=\im\gr\cCan\oplus\ker\gr\cVar$; use then that $\can:(\psi_{g,-1}Pf_\dag^i\cM,\rM_\bbullet)\to(\phi_{g,0}Pf_\dag^i\cM,\rM_{\bbullet-1})$ and $\var:(\phi_{g,0}Pf_\dag^i\cM,\rM_\bbullet)\to(\psi_{g,-1}Pf_\dag^i\cM,\rM_{\bbullet-1})$ are strictly compatible with the monodromy filtrations (\cf \cite[Lemma 5.1.12]{MSaito86}) to get that $\phi_{g,0}Pf_\dag^i\cM=\im\can\oplus\ker\var$, hence the S-decomposability of $f_\dag^i\cM$ and then, as in Proposition \ref{prop:Cdecomposable}, that of $f_\dag^i\cT$.\qed

\section[Proof of $\eqref{thimdirtwistor}_{(\leq(n-1),0)}$ \& ($\eqref{thimdirtwistor}_{(1,0)}$ with $\supp\cT$ smooth) $\implique \eqref{thimdirtwistor}_{(n,0)}$]{Proof of $\eqref{thimdirtwistor}_{(\leq(n-1),0)}$ and ($\eqref{thimdirtwistor}_{(1,0)}$ with $\supp\cT$ smooth) $\implique \eqref{thimdirtwistor}_{(n,0)}$ for $n\geq1$}\label{sec:n1n0}

We will argue as in \cite[\T5.3.8]{MSaito86} by using a Lefschetz pencil. Let $(\cT,\cS)$ be a polarized regular twistor $\cD$-module of weight~$w$ on a smooth complex projective variety and let $c$ be the Chern class of an ample line bundle on $X$. Assume that $\cT$ has strict support $Z$, which is an irreducible closed $n$-dimensional algebraic subset of $X$ ($n\geq1$). It is not restrictive to assume that $c$ is very ample, so that, by Kashiwara's equivalence, we may further assume that $X=\PP^N$ and $c=c_1(\cO_{\PP^N}(1))$. Choose a generic pencil of hyperplanes in $\PP^N$ and denote by $\wt X\subset X\times\PP^1$ the blowing up of $\PP^N$ along the axis $A$ of the pencil. We have the following diagram:
\begin{equation}
\begin{array}{c}
\xymatrix@C=1.5cm{
A\times\PP^1=\hspace*{-1.8cm}&\hbox{\raisebox{3pt}{$\wt A$}}\ar[d]\ar@<-2pt>@{^{ (}->}[r]&\hbox{\raisebox{4pt}{$\wt X$}}\ar@<-2pt>@{^{ (}->}[r]^-{\wt\imath}\ar[d]_-{\pi}\ar@/^2.7pc/[drr] ^(.8)f&X\times\PP^1\ar[dl]^-p\ar[dr]&\\
&A\ar@<-2pt>@{^{ (}->}[r]&X&&\PP^1
}
\end{array}
\end{equation}
Notice also that the restriction of $\pi$ to any $f^{-1}(t)$ is an isomorphism onto the corresponding hyperplane in $X$. Put $c'=c_1(\cO_{\PP^1}(1))$. Using Remark \ref{rem:hermdual}, we will assume that $\cT$ has weight~$w=0$ and that $(\cT,\cS)=((\cM,\cM,C),(\id,\id))$.

The proof will take five steps:
\begin{enumerate}
\item\label{n0step1}
We show that $\pi^+(\cT,\cS)$ satisfies (HSD), (REG), $(\MT_{>0})$ and $(\MTP_{>0})$ along $f^{-1}(t)$, for any $t\in\PP^1$.

Choose the pencil generic enough so that the axis $A$ of the pencil is noncharacteristic with respect to $\cM$. Then the inclusion $\wt\imath:\wt X\hto X\times\PP^1$ is noncharacteristic for $p^+\cM$ (\cf \T\ref{sec:nonchar}): this is clear away from $\wt A$; if the characteristic variety $\Lambda$ of $\cM$ is contained in a union of sets $T^*_{Z_i}X\times\Omega_0$, with $Z_i\subset Z$ closed, algebraic and irreducible, then, if $A$ is noncharacteristic with respect to each $Z_i$, so is $\wt A$ ---and therefore $\wt X$ near any point of $\wt A$--- with respect to each $Z_i\times\PP^1$. The characteristic variety of $\pi^+\cM$ is contained in the union of sets $T^*_{\wt Z_i}\wt X\times\Omega_0$, where $Z_i$ is the blow-up of $Z_i$ along $A\cap Z_i$.

Moreover, for any $t\in\PP^1$, the inclusion $A\times\{t\}\hto\wt X$ is noncharacteristic for $\pi^+\cM=\wt\imath^+p^+\cM$: by the choice of $A$, for any $Z_i$ as above, the intersection of $T^*_{A\times\{t\}}(X\times\PP^1)$ with $T^*_{Z_i\times\PP^1}(X\times\PP^1)$ is contained in the zero-section of $T^*(X\times\PP^1)$. As we have $T^*_{A\times\{t\}}(X\times\PP^1)=(T^*\wt\imath)^{-1} T^*_{A\times\{t\}}\wt X$, it follows that $T^*_{A\times\{t\}}\wt X\cap T^*_{\wt Z_i}\wt X\subset T^*_{\wt X}\wt X$.

This implies that, for any $t\in\PP^1$, the inclusion $f^{-1}(t)\hto\wt X$ is noncharacteristic for $\pi^+\cM$ near any point $(x_o,t)\in A\times\{t\}$.

Therefore, near each point $x_o$ of the axis of the pencil, $\psi_{f-t} \pi^+(\cT,\cS)$ is identified with $\psi_{g}(\cT,\cS)$, where $g=0$ is a local equation of the hyperplane $f=t$ near~$x_o$: argue as in the beginning of the proof of Lemma \ref{lem:strictnoncar} to show that $\pi^+\cM$ is specializable along $f=t$ and that there exists a good $V$-filtration for which $\gr_{-1}^V\wt\cM=i_{f^{-1}(t)}^+\wt\cM$; this module is equal to $i_{g^{-1}(0)}^+\cM=\psi_{g}\cM$ as $\cM$ is assumed to be \emph{strictly} noncharacteristic with respect to $g=0$; it follows that $\pi^+\cM$ is so with respect to $f=t$; a similar argument is used to identify the sesquilinear pairings; the identification of the sesquilinear dualities causes no problem, as they all are equal to $(\id,\id)$.

Using the identification above near the axis, and the properties assumed for $(\cT,\cS)$ on and away from the axis, we get all properties for $\pi^+(\cT,\cS)$ along any fibre $f^{-1}(t)$ (regularity follows from Lemma \ref{lem:noncar}). This concludes \eqref{n0step1}.

\item\label{n0step2}
As $A$ cuts $Z$ in codimension $2$, the support of $\pi^+\cM$ is the blow-up $\wt Z$ of $Z$ and the fibres of $f_{|\wt Z}$ all have dimension $n-1$. According to Step~\eqref{n0step1} and to Assumption $\eqref{thimdirtwistor}_{(n-1,0)}$, we can argue as in \T\ref{sec:nmnp} to obtain that $\big({\oplus}_if^i_+\pi^+(\cT,\cS),\cL_c\big)$ is an object of the category $\MLTr(\PP^1,w;1)^\rp$ with $w=0$. Denote by $a_\bbullet$ the constant map on the space $\cbbullet$. Then, by assumption (\ie by the result of \T\ref{sec:imdirtwistor10}), $\big({\oplus}_j\oplus_ia_{\PP^1+}^jf^i_+\pi^+(\cT,\cS),\cL_c,\cL_{c'}\big)$ is a polarized bigraded Lefschetz twistor structure of weight~$w=0$. It follows from Lemma \ref{lem:n1n0} that
$$
\big({\oplus}_k(\oplus_{i+j=k}a_{\PP^1+}^jf^i_+\pi^+(\cT,\cS)),\cL_c+\cL_{c'}\big)
$$
is a polarized graded Lefschetz twistor structure of weight~$w=0$. By using the same arguments as in \cite{Deligne68}, one shows that the Leray spectral sequence
\[
\oplus_j\oplus_ia_{\PP^1+}^jf^i_+\pi^+\cT\implique \oplus_ka_{\wt X,+}^k\pi^+\cT
\]
degenerates at $E_2$. Therefore, the Leray filtration $\Ler^\cbbullet a_{\wt X,+}^k\pi^+\cT$ attached to this spectral sequence satisfies in particular the following properties:
\begin{enumerate}
\item\label{num:Lera}
$\gr^j_{\Ler}a_{\wt X,+}^k\pi^+\cT=a_{\PP^1+}^jf^{k-j}_+\pi^+\cT=0$ for $j\neq-1,0,1$;
\item\label{num:Lerc}
$\cL_c^{k+1}:\Ler^1a_{\wt X,+}^{-k}\pi^+\cT\to\Ler^1a_{\wt X,+}^{k+2}\pi^+\cT(k+1)$ is an isomorphism for any $k\geq0$ (because $\Ler^1=\gr^1_{\Ler}$);
\item\label{num:Lerd}
$\ker\cL_{c'}:a_{\wt X,+}^{k}\pi^+\cT\to a_{\wt X,+}^{k+2}\pi^+\cT(1)$ is contained in $\Ler^0a_{\wt X,+}^{k}\pi^+\cT$ for any~$k$ (because $\cL_{c'}:\gr_{\Ler}^{-1}a_{\wt X,+}^{k}\pi^+\cT\to \gr_{\Ler}^{1}a_{\wt X,+}^{k+2}\pi^+\cT(1)$ is an isomorphism).
\end{enumerate}

We conclude from \eqref{num:Lera} that the object $\big({\oplus}_ka_{\wt X,+}^k\pi^+\cT,\cL_c+\cL_{c'}\big)$ is an extension of graded Lefschetz twistor structures of weight~$w=0$ and, by Remark \ref{rem:extleftriples}, is itself such an object (this argument is similar to that used in \cite[Th\ptbl5.2]{G-N90}).

\item\label{n0step3}
We now prove that $\pi_+\pi^+\cT$ decomposes as a direct sum, one summand being~$\cT$, and moreover that $(\cT,\cS)$ is a direct summand of $(\pi^0_+\pi^+\cT,\pi^0_+\pi^+\cS)$. We follow the proof given in \cite[\T5.3.9]{MSaito86}.

Everything has to be done along $A$ only, as $\pi$ is an isomorphism outside of $A$. Let $g$ be a local equation of a hyperplane containing $A$. Then $\pi^+\cM$ is strictly noncharacteristic along both components of $g\circ\pi=0$ and their intersections, so we may apply Lemma \ref{lem:psi2dcn}. Arguing as in Claim \ref{claim} (this is permissible due to the inductive assumption $\eqref{thimdirtwistor}_{(\leq(n-1),0)}$, as the fibres of $\pi:\wt Z\to Z$ have dimension $\leq n-1$), we conclude that $(\oplus_i\pi^i_+\pi^+(\cT,\cS),\cL_{c'})$ satisfies (HSD), (REG), ($\MLT_{>0}$) and $(\MLTP_{>0})$ (see Lemma \ref{lem:semisimplicitelef}) along $g=0$. We may cover $A$ by finitely many open sets where we can apply the previous argument. After \cite{Deligne68}, the complex $\pi_+\pi^+\cT$ decomposes as $\oplus_i\pi^i_+\pi^+\cT[-i]$ and clearly $\pi^i_+\pi^+\cT$ is supported on $A$ if $i\neq0$. Notice that, as $\cL_{c'}^2=0$, $\pi^0_+\pi^+(\cT,\cS)=P\pi^0_+\pi^+(\cT,\cS)$ satisfies (HSD), (REG), ($\MT_{>0}$) and $(\MTP_{>0})$ along $g=0$. We will identify $(\cT,\cS)$ with a direct summand of it.

Put $\pi^0_+\pi^+\cT=(\cM_0,\cM_0,C_0)$. It decomposes as $(\cM_1,\cM_1,C_1)\oplus(\cM_2,\cM_2,C_2)$ with $\cM_2$ supported on $A$ and $\cM_1$ has no submodule nor quotient supported on $A$ (use the S-decomposability along any $g=0$ as above). After lemma \ref{lem:adj}, there is an adjunction morphism $\cM\to\cM_0$. This morphism is an isomorphism away from $A$, and is injective, as $\cM$ has no submodule supported on $A$. Its image is thus contained in $\cM_1$, and is equal to $\cM_1$, as $\cM_1$ has no quotient supported on $A$. Therefore, $\cM=\cM_1$. That $C=C_1$ follows from Proposition \ref{prop:sesquiunidef}, applied to any germ of hyperplane containing~$A$. It remains to consider the polarization: notice that $\pi^+\cS=(\id,\id)$, and that $\pi^0_+\pi^+\cS=(\id,\id)$, as $(\pi^0_+\pi^+\cT)^*=\pi^0_+(\pi^+\cT)^*=\pi^0_+\pi^+\cT$; hence the identification of the polarizations.

\item\label{n0step4}
As $\cL_{c'}$ vanishes on $\cT$, we conclude from Step \eqref{n0step3} that $(a_{X,+}\cT,\cL_c)$ is a direct summand of $(a_{\wt X,+}\pi^+\cT,\cL_c+\cL_{c'})$. From Step \eqref{n0step2} and \cite{Deligne68} we have a (noncanonical) decomposition $a_{\wt X,+}\pi^+\cT\simeq\oplus_ka_{\wt X,+}^k\pi^+\cT[-k]$. Therefore, this decomposition can be chosen to induce a decomposition $a_{X,+}\cT\simeq \oplus_ka_{X,+}^k\cT[-k]$. In particular, $(\oplus_ka_{X,+}^k\cT,\cL_c)$ is a graded Lefschetz twistor structure of weight~$w=0$, being a direct summand of the graded Lefschetz twistor structure $(\oplus_ka_{\wt X,+}^k\pi^+\cT,\cL_c+\cL_{c'})$.

\item\label{n0step5}
It remains to show the polarization property. In order to do so, we will use the Fact \ref{fact:pol} in its graded Lefschetz form given by Remark \ref{rem:factpol}. Denote by $P'a_{\PP^1,+}^0(\oplus_kf_+^k\pi^+\cT)$ the $\cL_c'$-primitive part of $a_{\PP^1,+}^0(\oplus_kf_+^k\pi^+\cT)$, that is, the kernel of $\cL_{c'}$ acting on the previous space. Then $P'_0\defin\big(P'a_{\PP^1,+}^0(\oplus_kf_+^k\pi^+\cT),\cL_c\big)$ remains a (simply) graded Lefschetz twistor structure of weight~$w=0$, polarized by the family of sesquilinear dualities $a_{\PP^1,+}^0f_+^k\pi^+\cS$.

According to Remark \ref{rem:factpol}, we get the desired property if we show that
\begin{enumerate}
\item\label{num:Pprimea}
$(\oplus_ka_{X,+}^k\cT,\cL_c)$ is a sub graded Lefschetz twistor structure of $P'_0$,
\item\label{num:Pprimeb}
the polarization of $P'_0$ induces the family $(a_{X,+}^k\cS)_k$.
\end{enumerate}

By definition, for $k\geq0$, $Pa_{X,+}^{-k}\cT$ is the kernel of $\cL_c^{k+1}$ acting on $a_{X,+}^{-k}\cT$. It follows from \eqref{num:Lerc} that $Pa_{X,+}^{-k}\cT\cap\Ler^1a_{\wt X,+}^{-k}\pi^+\cT=\{0\}$. On the other hand, $a_{X,+}^{-k}\cT$ is contained in $\Ler^0a_{\wt X,+}^{-k}\pi^+\cT$ as $\cL_{c'}$ vanishes on $a_{X,+}^{-k}\cT$ and according to \eqref{num:Lerd}. Therefore, $Pa_{X,+}^{-k}\cT$ is contained in $\gr^0_{\Ler}a_{\wt X,+}^{-k}\pi^+\cT$, and more precisely in the biprimitive part $\ker\cL_{c'}\cap\ker\cL_c^{k+1}\subset a_{\PP^1,+}^0f_+^{-k}\pi^+\cT$. This gives \eqref{num:Pprimea}.

As we assume that $\cS=(\id,\id)$ and $\cT=\cT^*$, the sesquilinear duality $a_{X,+}^k\cS$ is nothing but the identification $(a_{X,+}^k\cT)^*=a_{X,+}^{-k}\cT$ deduced from \eqref{eq:fdag*}. Similarly, the sesquilinear duality on $P'_0$ is induced from the identification $(a_{\PP^1,+}^0f_+^k\pi^+\cT)^*=a_{\PP^1,+}^0f_+^{-k}\pi^+\cT$. That the former identification is induced by the latter is a consequence of Lemma \ref{lem:imdiradjointcomp}\eqref{lem:imdiradjointcomp2} and \eqref{lem:imdiradjointcomp3}.\qed
\end{enumerate}

\section[Proof of Theorem $\ref{thsmoothtw}$]{Proof of Theorem \ref{th:smoothtw}}
We will prove it by induction on $\dim X$. The result is clear when $\dim X=0$. It is easy if $\dim X=1$: by Remark \ref{rem:ramif}, it is enough to verify Properties (HSD), (REG), $(\MT_{>0})$ and $(\MTP_{>0})$ along any coordinate; as nearby cycles reduce then to ordinary restriction, the result is clear.

Let $\dim X=n\geq2$ and let $(\cT,\cS)$ be a smooth twistor structure of weight~$w$ on $X$. It is enough to consider the case where $w=0$ and $\cS=(\id,\id)$. By the computation of Proposition \ref{prop:localcomput} and by induction on $\dim X$, we know that Properties (HSD), (REG), $(\MT_{>0})$ and $(\MTP_{>0})$ are satisfied along any function like $x_1\cdots x_p$ and, by Remark \ref{rem:ramif}, along any monomial $(x_1\cdots x_p)^r$.

\medskip
Let $t:U\to\CC$ be any nonconstant holomorphic function on a connected open set $U$ of $X$. Let $\pi:\wt U\to U$ be a resolution of singularities of $t$: there exist local coordinates near each point of $\wt U$ so that $t\circ\pi$ is a monomial when expressed in these coordinates. It is a projective morphism. Choose a relatively ample line bundle on~$\wt U$ and denote by $c$ its Chern class. Assume that Properties (HSD), (REG), $(\MT_{>0})$ and $(\MTP_{>0})$ are satisfied for the inverse image $\pi^+\cT$ along $t\circ\pi=0$. Then, by the argument of \T\ref{sec:nmnp}, they are satisfied for $\oplus_i\pi_+^i\pi^+\cT$ along $t=0$. In particular, $\pi_+^0\pi^+\cT$ is strictly S-decomposable along $t=0$. We denote by $\cT'=(\cM',\cM',C')$ its component not supported by $t=0$. Remark that, for $\reel(\alpha)\in[-1,0[$, we have $\Psi_{t,\alpha}\cT'=\Psi_{t,\alpha}\pi_+^0\pi^+\cT$. Remark also that $\cL_c$ acts by $0$ on $\cT'$, as it acts by $0$ on $\cT$ and $\cT=\cT'$ away from $t=0$. It follows that $(\cT',\cS'=(\id,\id))$ satisfies (HSD), (REG), $(\MT_{>0})$ and $(\MTP_{>0})$ along $t=0$.

The natural adjunction morphism $\cM\to\pi_+^0\pi^+\cM$ of Lemma \ref{lem:adj} is injective, as $\cM$ is $\cO_{\cX}$-locally free, therefore it is an isomorphism onto $\cM'$. Moreover, $C'=C+C'_1$, where $C'_1$ takes values in distributions supported on $t=0$. Therefore, for $\reel(\alpha)\in[-1,0[$, we have $\psi_{t,\alpha}C'=\psi_{t,\alpha}C$, and $(\cT,\cS)$ satisfies (HSD), (REG), $(\MT_{>0})$ and $(\MTP_{>0})$ along $t=0$.

\medskip
To end the proof, we now have to consider the case where the function $t$ is any monomial. By a multi-cyclic covering, we may reduce this monomial to a power $(x_1\cdots x_p)^r$ and we may apply the first part of the proof to the reduced monomial. We are therefore reduced to proving that, if $\pi$ is the covering
\begin{align*}
X=\CC^n&\To{\pi}Y=\CC^n\\
(x_1,x_2,\dots,x_n)&\mto(x_1^k=y_1,x_2,\dots,x_k)
\end{align*}
then, if $(\cT,\cS)$ is a smooth polarized twistor structure on $Y$, it is a direct summand of $\pi_+^0\pi^+(\cT,\cS)$. Indeed, we may assume (by induction on the number of cyclic coverings needed), that $\pi^+(\cT,\cS)$ satisfies (HSD), (REG), $(\MT_{>0})$ and $(\MTP_{>0})$ along $t\circ\pi=0$. We will conclude as above that $\pi_+^0\pi^+(\cT,\cS)$ does so along $t=0$, and therefore so does $(\cT,\cS)$.

Put $\pi_*\cO_\cX=\cO_\cY\oplus\wt\cO$, where $\wt\cO$ denotes the sheaf of functions having trace zero along $\pi$. Similarly, put
\[
\frac{dx_1}{\hb}\otimes\pi_*\cO_\cX=\Big(\frac{dx_1^k}{\hb}\otimes\cO_\cY\Big)\oplus \wt\Omega^1,\qquad \wt\Omega^1=\Big(\frac{dx_1^k}{\hb}\otimes\wt\cO\Big)\oplus \Big(\frac{dx_1}{\hb}\otimes\cO_{\cX'}[x_1]_{\leq k-2}\Big),
\]
where $\wt\Omega^1$ is the sheaf of relative $1$-forms having trace $0$ along $\pi$ and $\cO_{\cX'}[x_1]_{\leq k-2}$ is the sheaf of polynomials of degree $\leq k-2$ in $x_1$ with coefficients depending on $x_2,\dots,x_n$ only. Notice that the relative differential $d:\cO_\cY\oplus\wt\cO\to\Big(\dfrac{dx_1^k}{\hb}\otimes\cO_\cY\Big)\oplus \wt\Omega^1$ is diagonal with respect to the direct sum decomposition.

We will compute $\pi_+$ by using the diagram
\[
\xymatrix@C=.5cm@R=0cm{
X\ar@{^{ (}->}[r]^-{i}&\CC\times X\ar[r]^-{\varpi}&Y\\
(x_1,\dots,x_n)\ar@{|->}[r]&(x_1^k,x_1,\dots,x_n)\ar@{|->}[r]& (x_1^k,\dots,x_n).
}
\]
Then, $\pi^+\cM=\cO_\cX\otimes_{\pi^{-1}\cO_\cY}\pi^{-1}\cM$ with its left $\cR_\cX$-structure (given by $\partiall_{x_1}(1\otimes m)=kx_1^{k-1}\otimes\partiall_{y_1}m$, \cf\T\ref{subsec:iminv}), and $\pi_+\pi^+\cM$ is the complex
\[
\pi_*\cO_\cX\otimes_{\cO_\cY}\cM[\tau]\To{\nabla_{\varpi}}\frac{dx_1}{\hb}\otimes\pi_*\cO_\cX\otimes_{\cO_\cY}\cM[\tau],
\]
with
\[
\nabla_{\varpi}(f\otimes m\tau^j)=df\otimes m\tau^j+\frac{dx_1^k}{\hb}\otimes \big[f\otimes (\partiall_{y_1}m\tau^j-m\tau^{j+1})\big].
\]
As $\nabla_{\varpi}$ is compatible with the direct sum decomposition corresponding to the trace, we have a decomposition of $\cR_\cY$-modules
\[
\pi^0_+\pi^+\cM=\cM\oplus\wt{\pi^0_+\pi^+\cM}.
\]
We will now show that this decomposition is orthogonal with respect to $C'\defin\pi^0_+\pi^+C$. It is enough to show that, for compactly supported $(n,n)$ forms $\varphi(y_1,x')=\chi(y_1,x')dy_1\wedge d\ov y_1\wedge\prod_{j=2}^ndx_j\wedge d\ov x_j$,
\begin{equation}\label{eq:piplus1}
\Big\langle\varpi_+(i_+C')^0\Big(\frac{dx_1^k}{\hb}\otimes(1\otimes m),\ov{\frac{dx_1^k}{\hb}\otimes(f\otimes \mu)}\Big),\varphi(y_1,x')\Big\rangle=0
\end{equation}
if $f\in\wt\cO$ and, for $0\leq \ell\leq k-2$,
\begin{equation}\label{eq:piplus2}
\Big\langle\varpi_+(i_+C')^0\Big(\frac{dx_1^k}{\hb}\otimes(1\otimes m),\ov{\frac{dx_1}{\hb}\otimes(x_1^\ell g(x')\otimes \mu)}\Big),\varphi(y_1,x')\Big\rangle=0.
\end{equation}
The left-hand term of \eqref{eq:piplus1} is, up to constants,
\begin{multline*}
C(m,\ov\mu)\int\chi(y_1,x')\ov f(x_1,x')\frac{dx_1^k\wedge d\ov x_1^k\wedge dy_1\wedge d\ov y_1\wedge\prod_{j=2}^ndx_j\wedge d\ov x_j}{d(y_1-x_1^k)\wedge d\ov{(y_1-x^k)}}\\
=C(m,\ov\mu)\int\chi(y_1,x')\ov f(x_1,x')dx_1^k\wedge d\ov x_1^k\wedge\prod_{j=2}^ndx_j\wedge d\ov x_j=0
\end{multline*}
as $\tr_\pi\ov f=0$. The argument for \eqref{eq:piplus2} is similar.\qed

\chapter{Integrability}\label{chap:int}
This chapter is concerned with the notion of integrability of a twistor $\cD$-module, a notion which is directly inspired from \cite{Hertling01}, where it is called a \emph{CV-structure}\footnote{``CV'' is for Cecotti-Vafa.}.

We define the notion of integrability of a $\cR_\cX$-module. We analyze the behaviour of such a notion with respect to various functors, like direct image by a proper morphism, inverse image, specialization. This notion is then extended to the category $\RTriples$, \ie we define the notion of integrability of a sesquilinear pairing between integrable $\cR_\cX$-modules.
We also analyze its behaviour with respect to the previous functors extended to the category $\RTriples$. Last, we extend Theorems \ref{th:imdirtwistor} and \ref{th:smoothtw} to the corresponding categories of integrable objects.

It may seem \emph{a priori} that this notion is useless when the underlying manifold is projective or affine: a variation of smooth polarizable twistor structures on a compact K\"ahler manifold (\ie a flat holomorphic vector bundle with a harmonic metric) is integrable if and only if it underlies a variation of polarized Hodge structures. I~also conjecture that the same result holds for a flat holomorphic vector bundle on a punctured compact Riemann surface with a \emph{tame} harmonic metric. In any case, a consequence of this integrability property is that, in the tame or regular case, the eigenvalues of local monodromies have absolute value equal to one.

Nevertheless, this integrability property seems to be the right generalization of the notion of a variation of polarized Hodge structure when irregular singularities occur.

\section{Integrable $\cR_\cX$-modules and integrable triples}

\subsection{Integrable $\cR_\cX$-modules}\label{subsec:intRmod}
Let $\cM$ be a $\cO_\cX$-module equipped with a flat relative connection $\nabla_{\cX/\Omega_0}$ as above. We say that $\cM$ is \emph{integrable} if $\nabla_{\cX/\Omega_0}$ comes from a (absolute) flat meromorphic connection $\nabla$ having Poincar\'e rank one along $\hb=0$, \ie such that $\hb\nabla$ has coefficients in the sheaf of \emph{logarithmic} $1$-forms $\Omega^1_\cX\langle\log\{\hb=0\}\rangle$, \ie has holomorphic coefficients when expressed in any local basis $dx_1,\dots,dx_n,d\hb/\hb$.

In the following, we denote by $\index{$dartialhbb$@$\ovdhb$}\ovdhb$ the operator $\hb^2\dhb$, using here the notion of geometric conjugation of \T\ref{subsubsec:conj}. This should not be confused with the corresponding operator using the usual conjugation on $\hb$. The latter will not be used in this chapter. We consider the sheaf of rings $\cR_\cX\langle\ovdhb\rangle$ generated by $\cR_\cX$ and $\ovdhb$, with the following commutation relations:
\[
[\ovdhb,\partiall_{x_i}]=\hb\partiall_{x_i},\quad[\ovdhb,f(\hb,x)]=\hb^2\frac{\partial f}{\partial\hb}(\hb,x).
\]
Let $\cM$ be a $\cR_\cX$-module. We say that $\cM$ in \emph{integrable} if the $\cR_\cX$-structure extends to a $\cR_\cX\langle\ovdhb\rangle$-structure; in other words, if $\cM$ is equipped with a $\cO_X$-linear operator $\ovdhb:\cM\to\cM$, which satisfies the previous commutation relations with the action of $\cR_\cX$. The integrable $\cR_\cX$-modules are the $\cO_\cX$-modules equipped with an absolute flat connection having Poincar\'e rank one along $\hb=0$.

Notice that, if $\cM$ admits a $\ovdhb$-action, it admits a family of such actions parametrized by $\CC$: for $\lambda\in\CC$ and $m\in\cM$, put $\ovdhb\cbbullet_\lambda m=(\ovdhb-\lambda\hb)m$.

\begin{Exemples}
\begin{enumerate}
\item
Let $\cM$ be a locally free $\cO_\cX$-module equipped with a flat meromorphic connection $\nabla$ having pole along the divisor $\hb=0$ at most, and having Poincar\'e rank one there (\ie $\hb\nabla$ has logarithmic poles along $\hb=0$). Then $\cM$ is a coherent holonomic $\cR_\cX$-module (with characteristic variety equal to the zero section in the relative cotangent bundle $(T^*X)\times\Omega_0$). Moreover, it is integrable by definition.

Examples of such objects are constructed in \cite{D-S02a} by partial Fourier transform of regular holonomic modules on $X\times\Afu$, equipped with a lattice (\ie a $\cO$-coherent submodule) when a noncharacteristic assumption is satisfied.
\item
Let $M$ be a coherent $\cD_X$-module equipped with a good filtration $F_\bbullet M$. Consider the Rees module $R_FM\defin\oplus_kF_k\hb^k$ on the Rees ring $R_F\cD_X$. After tensoring over $\cO_X[\hb]$ by $\cO_\cX$, on gets a coherent $\cR_\cX$-module $\cM$, with $\Char M\times\Omega_0$ as characteristic variety. As $F_\bbullet M$ is increasing, there is a natural action of $\hb\partial_\hb$ on $R_FM$, hence an action of $\ovdhb$. Therefore, $\cM$ is an integrable coherent $\cR_\cX$-module.
\end{enumerate}
\end{Exemples}

\subsection{Integrability of a sesquilinear pairing}\label{subsec:intsesquilin}
A \emph{sesquilinear pairing} between two $\cR_\cX$-modules $\cM',\cM''$ is a $\cR_{(X\ov X),\bS}$-linear morphism
\[
C:\cMS'\otimes_{\cO_\bS}\ov{\cMS''}\to\Dbh{X}.
\]

Let $\cM$ be an integrable $\cR_\cX$-module. On $\cMS$, we can consider the action of $\hbdhb$, defined as the action of $(1/\hb)\cdot\ovdhb$.

Define the action of $\hb\partial_\hb$ on $\cC^{\infty}_{X\times\bS}$ using polar coordinates, namely, if $\hb=\module{\hb}e^{i\arghb}$, $\hb\partial_\hb\varphi(x,\arghb)=-\tfrac i2\partial\varphi/\partial{\arghb}$. We therefore have a natural action of $\hb\partial_\hb$ on the sheaf of distributions $\Db_{X_\RR\times\bS}$.

Let $\cM',\cM''$ be integrable $\cR_\cX$-modules and let $C:\cMS'\otimes_{\cO_\bS}\ov{\cMS''}\to\Dbh{X}$ be a sesquilinear pairing. We say that the sesquilinear pairing is \emph{integrable} if the following equation is satisfied in $\Db_{X_\RR\times\bS}$ (recall that $\Dbh{X}$ is naturally included in $\Db_{X_\RR\times\bS}$):
\begin{equation}\label{eq:wflat}
\hbdhb C(m',\ov{m''})=C(\hbdhb m',\ov{m''})-C(m',\ov{\hbdhb m''})
\end{equation}
for any local sections $m',m''$ of $\cMS',\cMS''$. Although the right-hand term is in $\Dbh{X}$, the left-hand term has \emph{a priori} a meaning in $\Db_{X_\RR\times\bS}$ only; the integrability condition implies that it belongs to $\Dbh{X}$.

\subsection{The category $\RdTriples(X)$}\index{$rtriplesint$@$\RdTriples$}
We say that an object $\cT=(\cM',\cM'',C)$ of $\RTriples(X)$ is integrable if $\cM',\cM''$ are integrable, \ie equipped with a $\ovdhb$-action, and $C$ is integrable, \ie compatible with it, \ie satisfying \eqref{eq:wflat}. There is a family, parametrized by $\CC$, of $\ovdhb$-actions on $\cT$: for any $\lambda\in\CC$, consider the action by $\ovdhb-\lambda\hb$ on $\cM'$ and the action of $\ovdhb-\ov\lambda\hb$ on $\cM''$. We say that such actions are \emph{equivalent}.

Let $\cT_1,\cT_2$ be two integrable triples, each one equipped with an equivalence class of $\ovdhb$-actions. We say that a morphism $\varphi:\cT_1\to\cT_2$ is integrable if it commutes with some representatives of the $\ovdhb$-actions. Then, for any representative of the $\ovdhb$-action on $\cT_1$, there is a unique representative of the $\ovdhb$-action on $\cT_2$ such that $\varphi$ commutes with both.

Given two objects $\cT_1,\cT_2$ in $\RTriples(X)$, denote by $\Hom(\cT_1,\cT_2)$ the set of morphisms in $\RTriples(X)$ between these two objects. If $\cT_1,\cT_2$ are objects of $\RdTriples(X)$, denote by $\Homint(\cT_1,\cT_2)\subset \Hom(\cT_1,\cT_2)$ the set of integrable morphisms between them. The category $\RdTriples(X)$ is abelian.

The adjunction functor is an equivalence in $\RdTriples(X)$. There is a notion of sesquilinear duality, which has to be a morphism in $\RdTriples(X)$ between $\cT$ and~$\cT^*$.

Notice that the Tate twist is compatible with integrability: if \eqref{eq:wflat} is satisfied by~$C$, it is satisfied by $\thb^{-2k}C$ for $k\in\hZZ$ if we change the choice of the $\ovdhb$-action on $\cM'$ and $\cM''$, and replace it with the action of $\ovdhb-k\hb$ and $\ovdhb+k\hb$ respectively.

\begin{exemple}\label{ex:sesqui}
Let $\cT$ be integrable. Then its adjoint $\cT^*$ is also integrable. Let $w\in\ZZ$ and let $\cS$ be a sesquilinear duality of weight~$w$ on $\cT$, \ie a morphism $\cT\to\cT^*(-w)$ in $\RTriples(X)$. We say that $\cS$ is integrable if $\cS$ is a morphism in $\Homint(\cT,\cT^*(-w))$.

Let $\cT$ be integrable and equipped with a sesquilinear duality $\cS$ of weight~$w$. There is a family, parametrized by $\RR$, of $\ovdhb$-actions for which $(\cT,\cS)$ satisfies the same properties: for any $\lambda\in\RR$, consider the action by $\ovdhb-\lambda\hb$ on $\cM'$ and the action of $\ovdhb-\lambda\hb$ on $\cM''$
\end{exemple}

\subsection{Integrability and direct images}

Let $f:X\to Y$ be a holomorphic map. The direct image functor for $\cR_\cX$-modules is defined in \ref{subsec:imdir}, mimicking the corresponding functor for $\cD_X$-modules. The direct image for objects in $\RTriples(X)$ is defined in \T\ref{subsec:imdirsesqui} of \loccit\ Integrability is well-behaved with respect to the direct image functor of (right) $\cR_\cX$-modules or triples:

\begin{proposition}\label{prop:imdirint}
Let $\cM$ be a right $\cR_\cX$-module which is (right) integrable. Then each right $\cR_\cY$-module $R^jf_\dag\cM$ is right integrable. If $(\cM',\cM'',C)$ is an object of $\RdTriples(X)$, then $\cH^jf_\dag(\cM',\cM'',C)$ is an object of $\RdTriples(Y)$ for any $j\in\ZZ$.
\end{proposition}

Notice that $\omega_\cX$ is an integrable right $\cR_\cX$-module, and that the usual right$\to$left transformations for $\cR_\cX$-modules also transforms right integrability into left integrability.

\begin{proof}
Remark first that $\cR_\cX$ is an integrable left $\cR_\cX$-module. The left action of $\ovdhb$ is locally defined by
\[
\ovdhb\big(\textstyle\sum_\alpha a_\alpha(x,\hb)\partiall_x^\alpha\big)=\sum_\alpha\big(\hb^2\partial a/\partial\hb +\hb\module{\alpha}a\big)\partiall_x^\alpha.
\]
On the other hand, the sheaf $\Theta_\cX$ (vector fields relative to the projection $\cX\to\Omega_0$ which vanish along $\hb=0$) is equipped with a left action of $\ovdhb$: simply put, in some local coordinate system $(x_1,\dots,x_n)$ on $X$, $\ovdhb(\partiall_{x_i})=\hb\partiall_{x_i}$. Similarly, the exterior product is equipped with such an action, such that
\[
\ovdhb(\partiall_{x_{i_1}}\wedge\cdots\wedge\partiall_{x_{i_k}})=k\hb\partiall_{x_{i_1}}\wedge\cdots\wedge\partiall_{x_{i_k}}.
\]
It follows that each term of the Spencer complex $(\Sp_\cX^{\cbbullet}(\cO_\cX),\delta)$ is a left $\cR_\cX\langle\ovdhb\rangle$-module. One checks that the differential $\delta$ commutes with the $\cR_\cX\langle\ovdhb\rangle$-action.

If $\cM$ is a right $\cR_\cX\langle\ovdhb\rangle$-module and $\cN$ is a left $\cR_\cX\langle\ovdhb\rangle$-module, then $\cM\otimes_{\cO_\cX}\cN$ is a right $\cR_\cX\langle\ovdhb\rangle$-module. Similarly, $\cM\otimes_{\cR_\cX}\cN$ remains equipped with a right action of $\ovdhb$ defined by
\[
(m\otimes n)\cdot\ovdhb=m\ovdhb\otimes n-m\otimes\ovdhb n.
\]

Let $f:X\to Y$ be a holomorphic map. The relative Spencer complex $\Sp_{\cX\to\cY}^{\cbbullet}(\cO_\cX)$ is a complex of left $\cR_\cX\langle\ovdhb\rangle$-modules and right $f^{-1}\cR_\cY$-modules.

If $\cM$ is a right $\cR_\cX\langle\ovdhb\rangle$-module, then $\cM\otimes_{\cR_\cX}\Sp_{\cX\to\cY}^{\cbbullet}(\cO_\cX)$ remains equipped with a right $f^{-1}\cR_\cY$-module structure and a right action of $\ovdhb$. It is in fact a complex of right $\cR_\cY\langle\ovdhb\rangle$-modules.

These properties remain true after taking a Godement resolution. Therefore, the action of $\ovdhb$ is compatible with the construction of direct images given in \T\ref{subsec:imdir}, hence the first part of the proposition.

The integrability of the various $f_\dag^jC:\cH^{-j}f_\dag\cM'\otimes_{\cO_\bS}\ov{\cH^jf_\dag\cM''}\to \Dbh{Y}$ is then easy to get.
\end{proof}

\begin{remarque}[Integrability of the Lefschetz morphism]\label{rem:intlef}
Let us consider the situation of \T\ref{subsec:lef}. We have a Lefschetz morphism $\cL_c:f_\dag^j\cT\to f_\dag^{j+2}\cT(1)$. By the previous proposition we know that, if $\cT$ is integrable, $f_\dag^j\cT$ and $f_\dag^{j+2}\cT$ are so. We claim that $\cL_c\in\Homint(f_\dag^j\cT,f_\dag^{j+2}\cT(1))$: in the case of a projection, for instance, $\cL_c:(-L_\omega,L_\omega)$, where $L_\omega$ is $\hbm\omega\wedge$ and $\omega$ is a closed real $(1,1)$ form on $X$ with class~$c$; use that $\hbm(\ovdhb-\hb)=\ovdhb\hbm$.
\end{remarque}

\section{Integrable smooth twistor structures}

\subsection{Preliminary remark}
Assume that $X=\mathrm{pt}$ and that $(\cM',\cM'',C)$ defines a twistor structure of weight~$0$, that is, $\cM',\cM''$ are locally free $\cO_{\Omega_0}$-modules of finite rank and $C$ takes values in $\cO_\bS$ (\cf \T\ref{subsubsec:smtw0}). Saying that $\cM',\cM''$ are \emph{integrable} means that they are equipped with a connection having a pole of order $\leq2$ at $0$, and no other pole, or equivalently, that they are equipped with a $\ovdhb$-action. Then, saying that $C$ is integrable means that \eqref{eq:wflat} is satisfied when $C$ is viewed as taking values in $\cC^\infty_\bS$, \emph{via} the restriction $\cO_\bS\to\cC^\infty_\bS$.

The matrix of $C$ in local bases of $\cM',\cM''$ which are horizontal with respect to $\ovdhb$ is therefore constant when restricted to $\bS$. As it is assumed to be holomorphic in some neighbourhood of~$\bS$, it is constant, and $C$ satisfies \eqref{eq:wflat} in $\cO_\bS$. In other words, if we view $C$ as a gluing between the dual bundle $\cM^{\prime\vee}$ and the conjugate bundle $\ov{\cM''}$ on some neighbourhood of $\bS$, $C$ is integrable if and only if the previous isomorphism is compatible with the connections. Conversely, such a property clearly implies integrability of~$C$.

We say that $(\cM',\cM'',C)$ is an \emph{integrable twistor structure of weight~$0$} if it is a twistor structure of weight~$0$, if $\cM',\cM''$ are integrable, and $C$ is integrable.

We say that $(\cM',\cM'',C)$ is an \emph{integrable twistor structure of weight~$w$} if it is obtained by a Tate twist $(-w/2)$ from one with weight~$0$ (\cf \T\ref{subsec:somedef} for the definition of the Tate twist in this context).

\begin{exemple*}
Let us show that a complex Hodge structure defines an integrable twistor structure. We take notation of \T\ref{subsec:chsmtw0}, and we assume for simplicity that $w=0$. We define the $\ovdhb$-action on $\CC[\hb,\hbm]\otimes_\CC \ov H$ as the one induced by the natural one on $\CC[\hb,\hbm]$. Let us show for instance that $\cH''$ is stable under this action. For $m_q\in F^{\prime\prime q}$, we have $\ovdhb m_q\hb^{-q}=-q m_q\hb^{-q+1}$ and we have to show that $m_q\in F^{\prime\prime q-1}$, which follows from the fact that $F^{\prime\prime\cbbullet}$ is decreasing. The other compatibilities with the $\ovdhb$-action are verified similarly.
\end{exemple*}

\subsection{Characterization of integrable twistor structures}
We assume that $X=\textup{pt}$. Recall (\cf \T\ref{subsubsec:smtw0}) that a twistor structure $(\cH',\cH'',C)$ of weight~$0$ defines a vector bundle $\wt\cH$ on $\PP^1$ which is isomorphic to the trivial bundle, obtained by gluing $\cH^{\prime\vee}$ with $\ov{\cH''}$ using $C$ in some neighbourhood of $\bS$. There is an equivalence between the category of twistor structures of weight~$0$ and the category of finite dimensional $\CC$-vector spaces; one functor is
\[
(\cH',\cH'',C)\mto\wt\cH\mto H\defin\ov{\Gamma(\PP^1,\wt\cH)},
\]
and the quasi-inverse functor is
\[
H\mto \wt\cH\defin\ov H\otimes_\CC\cO_{\PP^1}\mto (\wt\cH^\vee_{|\Omega_0}, \ov{\wt\cH_{|\Omega_\infty}},C),
\]
where $C$ comes from the identity $\id:\wt\cH_{|\bS}\to\wt\cH_{|\bS}$.

\begin{lemme}
The twistor structure $(\cH',\cH'',C)$ of weight~$0$ is integrable if and only if the corresponding bundle $\wt\cH$ is equipped with a meromorphic connection $\wt\nabla$ having a pole of Poincar\'e rank at most one at $0$ and at infinity, and no other pole.
\end{lemme}

\begin{proof}
Indeed, if $(\cH',\cH'',C)$ is integrable, the bundles $\cH'$ and $\cH''$ are equipped with a meromorphic connection having a pole of Poincar\'e rank at most one at $0$ and no other pole. Therefore so has $\cH^{\prime\vee}$. Similarly, $\ov{\cH''}$ has a connection with a pole of Poincar\'e rank at most one at infinity. Integrability means that, \emph{via} the gluing, both connections coincide on some neighbourhood of $\bS$. The converse is also clear.
\end{proof}

Let $(\cH',\cH'',C)$ be an integrable twistor structure of weight~$0$. By the correspondence above, we have $\cH'=H\otimes_\CC\cO_{\Omega_0}$. Integrability means that there exist endomorphisms $U_0,Q,U_\infty$ of $H$ such that, for any element $m$ of $H$, we have
\[
\ovdhb m=(U_0-\hb Q-\hb^2U_\infty)m\in\cH'.
\]
In $\cH''=H^\vee\otimes_\CC\cO_{\Omega_0}$ we have, for any $\mu\in H^\vee$,
\[
\ovdhb\mu=({}^tU_\infty-\hb\, {}^tQ-\hb^2\,{}^tU_0)\mu\in\cH''.
\]
The category of integrable twistor structures of weight~$0$ is therefore equivalent to the category of tuples $(H,U_0,Q,U_\infty)$ and the morphisms are
the homomorphisms of vector spaces which are compatible with $(U_0,Q,U_\infty)$.

Assume that $(\cH',\cH'',C)$ is equipped with a Hermitian duality $\cS$. We will suppose that $\cH'=\cH''$ and $\cS=(\id,\id)$. This defines a Hermitian pairing $h:H\otimes\ov H\to \CC$. The compatibility of $\cS$ with the $\ovdhb$-action means that $Q$ is self-adjoint with respect to $h$ and $U_\infty$ is the $h$-adjoint of $U_0$.

If $\cS$ is a polarization, \ie if $h$ is positive definite, the eigenvalues of $Q$ are real, and $Q$ is semisimple. We decompose $H$ as
\[
H=\ooplus_{\alpha\in[0,1[}\ooplus_{p\in\ZZ}H_{\alpha+p}
\]
with respect to the eigenvalues $\alpha+p$ of $Q$. If we put $H^{p,-p}=\oplus_{\alpha\in[0,1[}H_{\alpha+p}$, we get a polarized complex Hodge structure of weight~$0$ on $H$.

\begin{remarque}
According to Example \ref{ex:sesqui}, if we change $Q$ in $Q+\lambda\id$ with $\lambda\in\RR$, we get an equivalent $\ovdhb$-action on $((\cM',\cM'',C),\cS)$.
\end{remarque}

\Subsection{Characterization of integrable smooth polarizable twistor structures}
\label{subsec:smoothinttw}
Let $(H,D_V,h)$ be a harmonic flat bundle, with a positive Hermitian metric $h$ and a flat connection $D_V=D_E+\theta_E$, where $\theta_E$ is the Higgs field. It corresponds to a smooth polarized twistor structure $(\cH',\cH',C)$ of weight~$0$ with polarization $\cS=(\id,\id)$ by the following rule: consider the $\cC^{\infty,\an}_\cX$-module $\cH=\cC^{\infty,\an}_\cX\otimes_{\pi^{-1}\cC^\infty_X}\pi^{-1}H$, equipped with the $d''$ operator
\begin{equation}\label{eq:dprimeprimeH}
D''_{\cH}=D''_E+\hb\theta''_E.
\end{equation}
This defines a holomorphic subbundle $\cH'=\ker D''_{\cH}$. Moreover, it has the natural structure of $\cR_\cX$-module, using the flat connection
\begin{equation}\label{eq:dprimeH}
D'_\cH=D'_E+\hbm\theta'_E.
\end{equation}

The integrability property means that the connection on $\cH'$ comes from an integrable absolute connection, that we denote with the same letter, which has a pole of Poincar\'{e} rank at most one along $X\times\{0\}$. The connection thus takes the form
\begin{align*}
D'_{\cH}&=D'_E+d'_\hb+\hbm\theta'_E+\Big(\frac{U_0}{\hb^2}-\frac{Q}{\hb}-U_\infty\Big)d\hb\\
D''_{\cH}&=D''_E+d''_\hb+\hb\theta''_E,
\end{align*}
where $U_0,U_\infty,Q$ are endomorphisms of the $C^\infty$ bundle $H$ and $d_\hb$ means the differential with respect to $\hb$ only. The compatibility with the polarization means that $U_\infty$ is the $h$-adjoint of $U_0$ and $Q$ is self-adjoint. Knowing that the relative connection $D_E+\hbm\theta'_E+\hb\theta''_E$ is integrable, the integrability condition is equivalent to the following supplementary conditions:
\begin{equation}\label{eq:int-tw}
\left\{\begin{aligned}{}
[\theta'_E,U_0]&=0,\\
D''_E(U_0)&=0,\\
D'_E(U_0)-[\theta'_E,Q]+\theta'_E&=0,\\
D'_E(Q)+[\theta'_E,U_\infty]&=0,
\end{aligned}\right.
\end{equation}
as the other conditions are obtained by adjunction. In particular, $U_0$ is an endomorphism of the holomorphic bundle $E$, which commutes with the holomorphic Higgs field $\theta'_E$.

\begin{corollaire}
Let $(H,D_V,h)$ be a harmonic flat bundle. Then it is integrable if and only if there exist endomorphisms $U_0,Q$ of $H$, $Q$ being self-adjoint with respect to $h$, satisfying Equations \eqref{eq:int-tw}, where $U_\infty$ denotes the $h$-adjoint of $U_0$.\qed
\end{corollaire}

\begin{Remarques}
\begin{enumerate}
\item
Equations \eqref{eq:int-tw} are the equations defining a CV-structure in \cite{Hertling01}, if one forgets the real structure, \ie if one forgets Equations (2.50-52) and (2.59) in \loccit
\item
For an integrable smooth twistor structure, the various local systems \hbox{$\ker(D_E+\hbm_o\theta'_E+\hb_o\theta''_E)\subset H$}, for $\hb_o\in\CC^*$, are all isomorphic to $\cL\defin\ker D_V$.
\item
It is a consequence of the equations for a flat harmonic bundle that the Higgs field $\theta_E$ satisfies
\[
D_V(\theta'_E-\theta''_E)=0,
\]
and therefore defines a class in $H^1(X,\endom(\cL))$. Notice now that Equations \eqref{eq:int-tw} imply in particular that, putting $A=-(U_0-Q-U_\infty)$, we have
\[
\theta'_E-\theta''_E=D_V(A).
\]
\ie the class of $\theta'_E-\theta''_E$ in $H^1(X,\endom(\cL))$ is zero. Moreover, $Q$ is the self-adjoint part of $A$ and $-U_0+U_\infty$ is its skew-adjoint part.

\item
For instance, if the polarized smooth twistor structure is associated to a variation of polarized complex Hodge structures of weight~$0$, we have $U_0=0=U_\infty$, and $Q$ is the endomorphism equal to $p\id$ on $H^{p,-p}$.
\end{enumerate}
\end{Remarques}

\begin{corollaire}[Rigidity on a compact K\"ahler manifold]\label{prop:rigidite}
Let $(H,D_V,h)$ be an integrable flat harmonic bundle on a compact K\"ahler manifold $X$. Then the corresponding $U_0$ is constant, and $Q$ defines a grading, so that $(H,D_V,h)$ corresponds to a variation of polarized complex Hodge structures of weight~$0$.
\end{corollaire}

\begin{proof}
We know that $U_0$ is an endomorphism of the holomorphic Higgs bundle $(E,\theta'_E)$. By the equivalence of \cite[Cor\ptbl1.3]{Simpson92}, it corresponds to an endomorphism of the flat bundle $\ker D_V$. This bundle is semi-simple, hence can be written as $\oplus_j(V_j,D_{V_j})^{p_j}$, with $p_j\in\NN$, where each $(V_j,D_{V_j})$ is simple and $(V_j,D_{V_j})\not\simeq(V_k,D_{V_k})$ for $j\neq k$. Then, any morphism $(V_j,D_{V_j})\to(V_k,D_{V_k})$ is zero for $j\neq k$ and equal to $\mathrm{cst}\cdot\id$ for $j=k$. By the correspondence quoted above, the same property holds for $U_0$ on the stable summands of the polystable Higgs bundle $(E,\theta'_E)$. In particular, $U_0$ is constant, and so is $U_\infty$.

Equations \eqref{eq:int-tw} reduce then to
\[
D_E(Q)=0,\quad\text{and}\quad [\theta'_E,Q]=\theta'_E.
\]
The eigenvalues of $Q$ are thus constant and the eigenspace decomposition of $Q$ is stable by $D_E$. Let $H_{\alpha+p}$ be the eigenspace corresponding to the eigenvalue $\alpha+p$ of $Q$, $\alpha\in[0,1[$, $p\in\ZZ$. Then $\theta'_E(H_{\alpha+p})\subset H_{\alpha+p-1}\otimes\Omega_X^1$. If we put $H^{p,-p}=\oplus_{\alpha\in[0,1[}H_{\alpha+p}$, we get a variation of polarized complex Hodge structures of weight~$0$.
\end{proof}

\begin{conjecture}\label{conj:tameint}
Let $X$ be a compact Riemann surface, let $P\subset X$ be a finite set of points, and let $(V,\nabla_V)$ be a semisimple holomorphic flat bundle on $X\moins P$. Denote by $(H,D_V,h)$ the tame harmonic flat bundle associated with it as in \cite{Simpson90,Biquard97}. Then, if $(H,D_V,h)$ is integrable, the endomorphism $U_0$ is compatible with the parabolic filtration defined by $h$ near each puncture.
\end{conjecture}

With the same argument as in Proposition \ref{prop:rigidite} we get:

\begin{corollaire}[Rigidity on a punctured Riemann surface]
If Conjecture \ref{conj:tameint} is true, the corresponding integrable tame harmonic flat bundle corresponds to a variation of polarized complex Hodge structures of weight~$0$ on $X\moins P$.\qed
\end{corollaire}

\section{Integrability and specialization}
Let $X'$ be a complex manifold, let $X$ be an open set in $\CC\times X'$, and let $t$ be the coordinate on $\CC$. Put $X_0=t^{-1}(0)\cap X$. We use definitions of \T\ref{subsec:S2t}.

\Subsection{Specialization of integrable $\cR_\cX$-modules}
\begin{proposition}\label{prop:speint}
Let $\cM$ be a $\cR_\cX$-module which is strictly specializable along $\cX_0$. Assume that $\cM$ is integrable. Then, for any $a\in\RR$ and any $\hb_o\in\Omega_0$, we have $\ovdhb V^{(\hb_o)}_a\cM\subset V^{(\hb_o)}_a\cM$ and, for any $\alpha\in\CC$ such that $\ell_{\hb_o}(\alpha)=a$, we have $\ovdhb\psi_{t,\alpha}\cM\subset\ovdhb\psi_{t,\alpha}\cM$, where $\ovdhb$ is viewed as acting on $\gr_a^{V^{(\hb_o)}}\cM$; in other words, each $\psi_{t,\alpha}\cM$ is an integrable $\cR_{\cX_0}$-module.
\end{proposition}

\begin{proof}
We will need the following lemma:
\begin{lemme}\label{lem:critereVa}
A local section $m$ of $\cM$ near $(x,\hb_o)$ is in $V^{(\hb_o)}_a\cM$ iff is satisfies a relation
\[
B_a(-\partiall_tt)m=n
\]
where $n$ is a local section of $V^{(\hb_o)}_a\cM$ and $B_a(s)=\prod_\gamma(s-\gamma\star\hb)^{\nu_\gamma}$, the product being taken on a finite set of $\gamma$ such that $\ell_{\hb_o}(\gamma)\leq a$.
\end{lemme}

\begin{proof}
The ``only if'' part is clear.
Assume that $m$ is a local section $V^{(\hb_o)}_b\cM$ for some $b>a$ satisfying such a relation with the polynomial $B_a(s)$. Then the class of~$m$ in $\gr_b^{V^{(\hb_o)}}\cM$ is killed by $B_a(-\partiall_tt)$ and $B_b(-\partiall_tt)$, where $B_b(s)=\prod_\beta(s-\beta\star\hb)^{\nu_\beta}$, the product being taken on a finite set $\beta$ such that $\ell_{\hb_o}(\beta)=b$. Therefore, the class of $m$ is killed by a nonzero polynomial in $\hb$, and by strictness, the class of $m$ is zero in $\gr_b^{V^{(\hb_o)}}\cM$.
\end{proof}

Let $m$ be a local section of $V^{(\hb_o)}_a\cM$, and let $b_m(s)$ be the minimal polynomial such that $b_m(-\partiall_tt)m=tPm$ where $P$ is a section of $V_0\cR_\cX$. We know that $b_m$ is a product of terms $s-\gamma\star\hb$ with $\ell_{\hb_o}(\gamma)\leq a$.

The following lemma is easy to prove:
\begin{lemme}
Let $k\in\ZZ$ and let $P$ be a local section of $V_k\cR_\cX$. Then $[\ovdhb,P]$ is a local section of $V_k\cR_\cX$ (and does not depend on $\ovdhb$).\qed
\end{lemme}

We then have
\begin{align*}
b_m(-\partiall_tt)\ovdhb m&= \ovdhb b_m(-\partiall_tt)m+Qm\quad Q\in V_0\cR_\cX\\
&=\ovdhb tPm+Qm\\
&=tP\ovdhb m+Rm\quad R\in V_0\cR_\cX
\end{align*}
Therefore, there exists $k\geq0$ such that, if we put $B_k(s)=\prod_{\ell=0}^kb_m(s-\ell\hb)$, we have $B_k(-\partiall_tt)\ovdhb m\in V^{(\hb_o)}_a\cM$. Apply then Lemma \ref{lem:critereVa} to get that $\ovdhb m$ is a local section of $V^{(\hb_o)}_a\cM$. This gives the first part of Proposition \ref{prop:speint}.

Denote by $\ovdhb$ the induced operator on $\gr_a^{V^{(\hb_o)}}\cM$. We now want to show that, for any $\alpha\in\CC$ with $\ell_{\hb_o}(\alpha)=a$, $\cup_n\ker\big[(\partiall_tt+\alpha\star\hb)^n:\gr_a^{V^{(\hb_o)}}\cM\to\gr_a^{V^{(\hb_o)}}\cM\big]$ is stable by $\ovdhb$. The point is that $\ovdhb$ does not commute with $\partiall_tt+\alpha\star\hb$, but $[\ovdhb,(\partiall_tt+\alpha\star\hb)^n]$ is a polynomial in $\partiall_tt+\alpha\star\hb$ with polynomial coefficients in $\hb$, and therefore commutes with $\partiall_tt+\alpha\star\hb$. Let $m$ be a local section of $\gr_a^{V^{(\hb_o)}}\cM$ killed by $(\partiall_tt+\alpha\star\hb)^n$. Then
$$
(\partiall_tt+\alpha\star\hb)^n\ovdhb m=-[\ovdhb,(\partiall_tt+\alpha\star\hb)^n]m=\sum_jp_j(\hb)(\partiall_tt+\alpha\star\hb)^jm
$$
and certainly $(\partiall_tt+\alpha\star\hb)^{2n}\ovdhb m=0$.
\end{proof}

\begin{corollaire}\label{cor:extmin}
Let $\wt\cM$ be a strictly specializable $\cR_\cX[t^{-1}]$-module (as defined in \T\ref{sec:minext}) which is integrable. Then the minimal extension $\cM$ of $\wt\cM$ across $\cX_0$ is integrable.
\end{corollaire}

\begin{proof}
By definition, we have $V^{(\hb_o)}_{<0}\wt\cM=V^{(\hb_o)}_{<0}\cM$, therefore this is stable by the $\ovdhb$-action, according to the proposition. One shows similarly that all $V^{(\hb_o)}_a\cM$, defined in \loccit, are stable under the $\ovdhb$-action.
\end{proof}

\begin{remarque}[S-decomposability]
Assume that $\cM$ is strictly specializable along $\cX_0$ and integrable. Then the morphism $\var$ of Remark \ref{rem:psi}\eqref{rem:Ncanvar2} commutes with the $\ovdhb$-action, but the morphism $\can$ does not. However, $\im \can$ is stable by the $\ovdhb$-action, because $\ovdhb\partiall_t=\partiall_t(\ovdhb+\hb)$. Similarly, if $\cM$ is strictly decomposable along $\cX_0$, its strict components are integrable, as can be seen from the proof of Proposition \ref{prop:canvar}\eqref{prop:canvare}. As a consequence, if $\cM$ is integrable and strictly S-decomposable, its strict components are integrable.
\end{remarque}

\begin{remarque}[Local unitarity]
When working with twistor $\cD$-modules, we are led to consider the graded modules $\gr^{\rM}_\ell\psi_{t,\alpha}\cM$ with respect to the monodromy filtration $\rM_\bbullet(\rN)$ of the nilpotent endomorphism $\rN=-(\partiall_tt+\alpha\star\hb)$. \emph{A priori}, $\ovdhb$ is not compatible with the monodromy filtration, therefore we would need a new assumption to insure that this compatibility is satisfied. However, we will see below that when $\dim X=1$ and if all $\gr^{\rM}_\ell\psi_{t,\alpha}\cM$ are strict, this compatibility is automatically satisfied, as a consequence of the fact that the complex numbers $\alpha$ to be considered in the various Bernstein polynomials are real. We will see in \T\ref{subsec:inttwist} that this property extends to integrable twistor $\cD$-modules. We call it \emph{local unitarity}.

When a strictly specializable $\cR_\cX$-module is locally unitary, the various $V^{(\hb_o)}$-filtrations glue together when $\hb_o$ varies in $\Omega_0$ and we forget the exponent $\hb_o$. Moreover, if $\alpha=a$ is real, we then have $\gr_\alpha^V\cM=\psi_{t,\alpha}\cM$. Last, we have $\ell_{\hb_o}(\alpha)=\alpha$ and $\alpha\star\hb=\hb\alpha$.
\end{remarque}

\begin{lemme}\label{lem:unit}
Assume that $X$ is a disc with coordinate $t$. If $\cM$ is an integrable strictly specializable $\cR_\cX$-module such that each $\gr^{\rM}_\ell\psi_{t,\alpha}\cM$ is strict, it is locally unitary.
\end{lemme}

\begin{proof}
Fix $\alpha\in\CC$. As each $\gr^{\rM}_\ell\psi_{t,\alpha}\cM$ is $\cO_{\Omega_0}$-free, there exists a basis $\bme$ of $\psi_{t,\alpha}\cM$ for which the matrix $\rY$ of $\rN$ has the Jordan normal form, in particular is constant and nilpotent. Denote by $A(\hb)$ the matrix of $\ovdhb$ in this basis. Then
\[
\ovdhb\bme=\bme\cdot A(\hb),\quad -\partiall_tt\bme=\bme\cdot\big[(\alpha\star\hb)\id+\rY\big].
\]
Therefore,
\begin{align*}
-\ovdhb\partiall_tt\bme&=\bme\cdot [(\alpha\star\hb)A(\hb)+A(\hb)\rY +\hb^2\partial(\alpha\star\hb)/\partial\hb\id],\\
-\partiall_tt(\ovdhb+\hb)\bme&=\bme\cdot[(\alpha\star\hb)\id+\rY][A(\hb)+\hb\id].
\end{align*}
As the operators $\ovdhb$ and $\partiall_tt$ satisfy the commutation relation $\ovdhb\partiall_tt=\partiall_tt(\ovdhb+\hb)$, we must have
\[
\hb[\hb\partial(\alpha\star\hb)/\partial\hb-\alpha\star\hb]\id=[\rY,A(\hb)]+\hb\rY,
\]
thus, taking the trace, we get that $\alpha$ must be such that, for any $\hb\in\CC$, $\hb\partial(\alpha\star\hb)/\partial\hb=\alpha\star\hb$. But $\hb\partial(\alpha\star\hb)/\partial\hb=\alpha\star\hb+i\alpha''(\hb^2-1)/2$. Therefore, if $\psi_{t,\alpha}\cM\neq0$, $\alpha$ must be such that $\alpha''=0$, \ie $\alpha$ must be real.
\end{proof}

Let us now go back to $\dim X\geq1$.

\begin{lemme}\label{lem:locunitint}
If $\cM$ is strictly specializable along $t=0$ and locally unitary, then, if $\cM$ is integrable, so is each $\gr_\ell^{\rM}\psi_{t,\alpha}\cM$.
\end{lemme}

\begin{proof}
We now have $\ovdhb\rN=\rN(\ovdhb+\hb)$, hence the kernel filtration $\ker\rN^k$ and the image filtration $\im\rN^k$ of $\rN$ are stable by $\ovdhb$. As the monodromy filtration $\rM_\bbullet(\rN)$ is obtained by convolution of these two filtrations (\cf \cite[Remark (2.3)]{S-Z85}) it is also stable by $\ovdhb$.
\end{proof}

\subsection{Specialization of sesquilinear pairings}
The definition of specialization of a sesquilinear pairing involves the residue of a distribution depending meromorphically on a complex variable $s$ along a set having equation $s=\alpha\star\hb/\hb$, for a fixed complex number $\alpha$ and for $\hb$ varying in $\bS$. In general, the compatibility of taking the residue along such a set and the action of $\ovdhb$ is not clear. However, as soon as we assume local unitarity, \ie $\alpha\in\RR$, then $\alpha\star\hb/\hb=\alpha$ does not depend on $\hb$ and the compatibility is clearly satisfied. We therefore obtain:

\begin{proposition}\label{prop:locunitint}
Let $\cT$ be an object of $\RdTriples(X)$. Assume that the components $\cM',\cM''$ are strictly specializable and locally unitary along $t=0$. Then $\psi_{t,\alpha}\cT$ is integrable for any $\alpha\in\RR$. Moreover, the morphism $\cN:\psi_{t,\alpha}\cT\to\psi_{t,\alpha}\cT(-1)$ is integrable.
\end{proposition}

\begin{proof}
It remains to explain the integrability of $\cN$ defined by \eqref{eq:NTate}. We have $\cN=(\rN',\rN'')$ with $\rN''=i\hb(\partial_tt+\alpha)=-\rN'$. Then we argue as in Remark \ref{rem:intlef}, using that $\hb(\ovdhb+\hb)=\ovdhb\hb$.
\end{proof}

\section{Integrable polarizable regular twistor $\cD$-modules}

\subsection{A preliminary lemma on twistor $\cD$-modules}\label{subsec:prelimlem}
Let $(\cM',\cM'',C)$ be an object in $\MT_{\leq d}(X,w)$. Put $\cM=\cM'$ or $\cM''$. Let $f$ be holomorphic functions on some open set $U$ of $X$. Then $\cM$ is strictly specializable along $f=0$ and, for any $\alpha\in\CC$, $\psi_{f,\alpha}\cM$ is equipped with a nilpotent endomorphism $\rN$. Denote by $\rM_\bbullet(\rN)$ the corresponding monodromy filtration. Then each $\rM_\ell\psi_{f,\alpha}\cM$ is strict and, by definition of $\MT$, each $\gr_\ell^{\rM}\psi_{f,\alpha}\cM$ is also strict.

Let $g$ be another holomorphic function. By definition, each $\gr_\ell^{\rM}\psi_{f,\alpha}\cM$ is strictly specializable along $g=0$. By induction on $\ell$, this implies that each $\rM_\ell\psi_{f,\alpha}\cM$ is so, and, for any $\beta\in\CC$, we have exact sequences
\[
0\to\psi_{g,\beta}\rM_{\ell-1}\psi_{f,\alpha}\cM \to \psi_{g,\beta}\rM_\ell\psi_{f,\alpha}\cM\to \psi_{g,\beta}\gr^{\rM}_\ell\psi_{f,\alpha}\cM\to0.
\]

Denote by $\rM_\bbullet(\psi_{g,\beta}\rN)$ the monodromy filtration of the nilpotent endomorphism on $\psi_{g,\beta}\psi_{f,\alpha}\cM$. Then, according to the previous exact sequence and to the uniqueness of the monodromy filtration, we have
\[
\rM_\bbullet(\psi_{g,\beta}\rN)=\psi_{g,\beta}\rM_\bbullet\psi_{f,\alpha}\cM.
\]
In particular, each $\gr^{\rM}_\ell\psi_{g,\beta}\psi_{f,\alpha}\cM$ is strict, being equal to $\psi_{g,\beta}\gr^{\rM}_\ell\psi_{f,\alpha}\cM$.

Let now $f_1,\dots,f_p$ be holomorphic functions and let $\alpha_1,\dots,\alpha_p$ be complex numbers. Under the same assumption on $\cM$ we obtain similarly:

\begin{lemme}\label{lem:strictspespe}
For any $\ell\in\ZZ$ and for any $j=1,\dots,p$, $\psi_{f_j,\alpha_j}\cdots\psi_{f_1,\alpha_1}\gr^{\rM}_\ell\psi_{f,\alpha}\cM$ is strict and strictly specializable, and we have
$$
\gr^{\rM}_\ell\psi_{f_j,\alpha_j}\cdots\psi_{f_1,\alpha_1}\psi_{f,\alpha}\cM=\psi_{f_j,\alpha_j}\cdots\psi_{f_1,\alpha_1}\gr^{\rM}_\ell\psi_{f,\alpha}\cM. \eqno\qed
$$
\end{lemme}

\Subsection{Integrable twistor $\cD$-modules}\label{subsec:inttwist}

\begin{proposition}\label{prop:twlocunit}
Let $(\cM',\cM'',C)$ be an object in $\MT_{\leq d}(X,w)$. Assume that it is integrable, \ie is also an object of $\RdTriples(X)$. Then $\cM'$ and $\cM''$ are locally unitary.
\end{proposition}

\begin{proof}
Let $f$ be a holomorphic function defined in some open set $U\subset X$. Assume that there exists $\alpha\in\CC\moins\RR$ such that $\psi_{f,\alpha}\cM\neq0$ for $\cM=\cM'$ or $\cM=\cM''$. Let $\ell\in\ZZ$. By assumption, $\gr_\ell^{\rM}\psi_{f,\alpha}\cM$ is strictly S-decomposable. For any strict component $Z$ of its support, let $(\gr_\ell^{\rM}\psi_{f,\alpha}\cM)_Z$ be the corresponding direct summand. It is enough to show that each $(\gr_\ell^{\rM}\psi_{f,\alpha}\cM)_Z$ is zero, and also that its restriction to dense open set of $Z$ is zero. We can assume that the characteristic variety of $(\gr_\ell^{\rM}\psi_{f,\alpha}\cM)_Z$ is equal to $T^*_ZX\times\Omega_0$ near a general point $x_o$ of $Z$. Therefore, near such a point, by Kashiwara's equivalence Cor\ptbl\ref{cor:Kashiwaraequiv} and Prop\ptbl\ref{prop:CK}, $(\gr_\ell^{\rM}\psi_{f,\alpha}\cM)_Z$ is the direct image by the inclusion $Z\hto X$ of a locally free $\cO_\cZ$-module.

Let $f_1,\cdots,f_p$ be holomorphic functions near $x_o$ inducing a local coordinate system on $Z$. By induction on $p$, using Lemma \ref{lem:strictspespe} and Proposition \ref{prop:speint}, one shows that the $\cR_\cX$-module $\cN\defin \psi_{f_p,\alpha_p}\cdots\psi_{f_1,\alpha_1}\psi_{f,\alpha}\cM$ (with $\cM=\cM'$ or $\cM''$) is integrable. By Lemma \ref{lem:strictspespe}, for any $\ell\in\ZZ$, $\gr_\ell^{\rM}\cN$ is strict and is supported on $x_o$. By Kashiwara's equivalence Cor\ptbl\ref{cor:Kashiwaraequiv}, we can apply the same argument as in Lemma \ref{lem:unit} to conclude that $\gr_\ell^{\rM}\cN=0$ for any $\ell$, as we assume $\alpha\not\in\RR$. Therefore, applying once more Lemma \ref{lem:strictspespe}, we obtain that $\psi_{f_p,\alpha_p}\cdots\psi_{f_1,\alpha_1}(\gr_\ell^{\rM}\psi_{f,\alpha}\cM)_Z=0$. Near $x_o$, $\psi_{f_j}$ is nothing but the usual restriction to $f_j=0$, therefore the restriction of $(\gr_\ell^{\rM}\psi_{f,\alpha}\cM)_Z$ at $x_o$ is zero. But $(\gr_\ell^{\rM}\psi_{f,\alpha}\cM)_Z$ is (the direct image of) a locally free $\cO_\cZ$-module, hence, by Nakayama, $(\gr_\ell^{\rM}\psi_{f,\alpha}\cM)_Z=0$ near $x_o$, a contradiction.
\end{proof}

From Lemma \ref{lem:locunitint} and Proposition \ref{prop:locunitint} we get:

\begin{corollaire}
Let $(\cM',\cM'',C)$ be an object of $\MT_{\leq d}(X,w)$ and let $f$ be a holomorphic function on some open set $U$ of $X$. Then, for any $\alpha\in[-1,0[$ and any $\ell\in\ZZ$, the object $\gr_\ell^{\rM}\psi_{f,\alpha}(\cM',\cM'',C)$ of $\MT_{\leq d}(U,w+\ell)$ is integrable.\qed
\end{corollaire}

Notice that, according to Proposition \ref{prop:twlocunit}, we do not have to consider $\psi_{f,\alpha}$ for $\alpha\in\CC\moins\RR$, and that the two functors $\psi$ and $\Psi$ (\cf Definition \ref{def:nearby}) coincide.

We define the category of integrable twistor $\cD$-modules ${\MTint}_{\leq d}(X,w)$ as the subcategory of $\MT_{\leq d}(X,w)$ having integrable objects and integrable morphisms. By the previous corollary, it is stable by taking $\gr_\ell^{\rM}\Psi_{f,\alpha}$. It shares many properties of $\MT_{\leq d}(X,w)$ (\cf \T\ref{subsec:defDtwist}): it is abelian, it is local, it satisfies Kashiwara's equivalence, it is stable by direct summand in $\RdTriples(X)$. However, it is \emph{a priori} not stable by direct summand in $\RTriples(X)$ or in $\MT_{\leq d}(X,w)$.

The subcategory $\MTintr(X,w)$ of regular objects is defined similarly. Last, the category $\MLTintr(X,w)$ of graded Lefschetz objects is defined as in \T\ref{subsec:grleftw}.

\subsection{Integrable polarizable regular twistor $\cD$-modules}
Let $\cT$ be an integrable twistor $\cD$-module of weight~$w$ as defined above. We say that a polarization of $\cT$ is integrable if it is an integrable morphism $\cT\to\cT^*(-w)$.

It is now clear that the two main theorems of Chapter \ref{chap:decomp} have the following integrable counterpart:

\begin{theoreme}\label{th:imdirtwistorint}
Let $f:X\to Y$ be a projective morphism between complex analytic manifolds and let $(\cT,\cS)$ be an object of $\MTintr(X,w)^\rp$. Let $c$ be the first Chern class of a relatively ample line bundle on $X$ and let $\cL_c$ be the corresponding Lefschetz operator. Then $(\oplus_if_\dag^i\cT,\cL_c,\oplus_if_\dag^i\cS)$ is an object of $\MLTintr(Y,w;1)^\rp$.\qed
\end{theoreme}

\begin{theoreme}\label{th:smoothtwint}
Let $X$ be a complex manifold and let $(\cT,\cS)$ be an integrable smooth polarized twistor structure of weight~$w$ on $X$, in the sense of \T\ref{subsec:smoothinttw}. Then $(\cT,\cS)$ is an object of $\MTintr(X,w)^\rp$.\qed
\end{theoreme}

\chapter[Partial Fourier-Laplace~transform]{Monodromy at infinity and partial~Fourier-Laplace~transform}\label{chap:8}
\def\theenumi{\roman{enumi}}

In this chapter, we analyze the behaviour of polarized regular twistor $\cD$-modules under a partial (one-dimensional) Fourier-Laplace transform. We generalize to such objects the main result of \cite{Bibi96a}, comparing, for a given function $f$, the nearby cycles at $f=\infty$ and the nearby or vanishing cycles for the partial Fourier-Laplace transform in the $f$-direction (Theorem \ref{th:Fstrict}).

\section{Exponential twist}
\subsection{Exponential twist of an object of $\RTriples$}\label{subsec:twist}
Let $t:X\to\CC$ be a holomorphic function on the complex manifold $X$. If $\cM$ is a left $\cR_\cX$-module, \ie a $\cO_\cX$-module with a flat relative meromorphic connection $\nabla_{\cX/\Omega_0}$, the twisted $\cR_\cX$-module $\FcM=\cM\otimes \ccE^{-t/\hb}\index{$mcurlF$@$\FcM$}\index{$eth$@$\ccE^{-t/\hb}$}$ is defined as the $\cO_\cX$-module $\cM$ equipped with the twisted connection $e^{t/\hb}\circ\nabla_{\cX/\Omega_0}\circ e^{-t/\hb}$. If $\cM$ is integrable (\cf Chapter~\ref{chap:int}), then so is $\FcM$: just twist the absolute connection $\nabla$. Notice that, if $\nabla$ has Poincar\'e rank one, so has the twisted connection.

Let $C:\cMS'\otimes_{\cO_{\cX|\bS}}\ov{\cMS''}\to\Dbh{X}$ be a sesquilinear pairing. Then $\Fou C\defin \exp(\hb\ov t-t/\hb)C\index{$cf$@$\Fou C$}$ is a sesquilinear pairing $\FcM'_{|\bS}\otimes_{\cO_{\cX|\bS}}\ov{\FcM''_{|\bS}}\to\Dbh{X}$, \ie is $\cR_{(X,\ov X),\bS}$-linear.

If $\cT=(\cM',\cM'',C)$ is an object of $\RdTriples(X)$, then so is $\FcT\defin(\FcM',\FcM'',\Fou C)\index{$th$@$\FcT$}$. Exponential twist is compatible with Tate twist and adjunction (as $\hb\ov t-t/\hb=\hb\ov t+\ov\hb t$ is ``real'').

If $\varphi:\cT_1\to\cT_2$ is a morphism, then $\varphi$ induces a morphism $\varphi:\FcT_1\to\FcT_2$.
In particular, if $\cS$ is a sesquilinear duality of weight~$w$ on $\cT$, then $\cS$ induces a sesquilinear duality of the same weight on $\FcT$.

\subsection{Exponential twist of flat and Higgs bundles}\label{subsec:expflathiggs}
We will now give an explicit description of the exponential twist in the case of smooth triples, using the language of Higgs bundles. Let $H$ be a $C^\infty$-bundle on $X$ equipped with a flat connection $D_V=D'_V+d''$ and a Hermitian metric $h$. Denote by $V=\ker d''$ the corresponding holomorphic bundle, equipped with the holomorphic connection $\nabla_V$. Using the function $t$ we twist the connection $D_V$ and define
\begin{align*}
\Fou D_V\index{$dVF$@$\Fou D_V$}&=e^t\circ D_V\circ e^{-t},\quad \text{\ie } \Fou D'_V=D'_V-dt,\ \Fou D''_V=d'',\\
\Fou h\index{$hF$@$\Fou h$}&=e^{2\reel t}h.
\end{align*}

Using definitions in \cite{Simpson90,Simpson92}, we have:

\begin{lemme}
If the triple $(H,D_V,h)$ is harmonic on $X$, then so is the triple $(H,\Fou D_V,\Fou h)$.
\end{lemme}

The Higgs field is given by the formulas
\[
\index{$thEF$@$\Fou\theta_E$}\Fou\theta'_E=\theta'_E-dt,\quad\Fou\theta''_E=\theta''_E-d\ov t,
\]
and the metric connection $\Fou D_E=\Fou D'_E+\Fou D''_E$ by
\[
\index{$dEF$@$\Fou D_E$}\Fou D_E=e^{-\ov t}\circ D_E\circ e^{\ov t},\quad\text{\ie } \Fou D'_E=D'_E,\ \Fou D''_E= D''_E+ d\ov t.
\]

\begin{lemme}\label{lem:FouFou}
If $(\cM,\cM,C)$ denotes the smooth polarized twistor structure of weight~$0$ corresponding to the harmonic bundle $(H,D_V,h)$ then, using notation of \T\ref{subsec:twist}, the triple $(\FcM,\FcM,\Fou C)$ is the smooth polarized twistor structure of weight~$0$ corresponding to the harmonic bundle $(H,\Fou D_V,\Fou h)$, \emph{via} the correspondence of \T\ref{subsec:smtwqc}.
\end{lemme}

\begin{proof}
Consider the $\cC^{\infty,\an}_\cX$-module $\cH=\cC^{\infty,\an}_\cX\otimes_{\pi^{-1}\cC^\infty_X}\pi^{-1}H$, equipped with the $d''$ operator
\begin{equation}\label{eq:FouDprimeprimeH}
\Fou D''_{\cH}=\Fou D''_E+\hb\Fou\theta''_E.
\end{equation}
We get a holomorphic subbundle $\FcH'=\ker \Fou D''_{\cH}\subset\cH$ equipped with a $\cR_\cX$-action, \ie a relative connection $\Fou\nabla_{\cX/\Omega_0}$, obtained from the connection $\Fou D'_{\cH}=\Fou D'_E+\hbm\Fou\theta'_E$. We have by definition
\begin{equation}\label{eq:DprimeH}
\begin{split}
\Fou D'_\cH&=\exp(t/\hb)\circ D'_\cH\circ\exp(-t/\hb),\\
\Fou D''_\cH&=\exp\big((\hb-1)\ov t\big)\circ D''_\cH\circ \exp\big((1-\hb)\ov t\big).
\end{split}
\end{equation}
We have an isomorphism
\[
(\cH',e^{t/\hb}\circ\nabla_{\cX/\Omega_0}\circ e^{-t/\hb})\To{\cdot\exp\big((\hb-1)\ov t\big)} (\FcH',\Fou\nabla_{\cX/\Omega_0})
\]
and, \emph{via} this isomorphism, $\Fou C=\Fou h_{\FcH'_{|\bS}\otimes\ov{\FcH'_{|\bS}}}$ corresponds to $e^{\hb\ov t+\ov\hb t}\cdot h_{\cH'_{|\bS}\otimes\ov{\cH'_{|\bS}}}=e^{\hb\ov t+\ov\hb t}C$.
\end{proof}

\section{Partial Fourier-Laplace transform of $\cR_\cX$-modules} \label{sec:Fourier}

\subsection{The setting}\label{subsec:Fouriersettings}
We consider the product $\Afu\times\Afuh$ of two affine lines with coordinates $(t,\tau)$, and the compactification $\PP^1\times\wh\PP^1$, covered by four affine charts, with respective coordinates $(t,\tau)$, $(t',\tau)$, $(t,\tau')$, $(t',\tau')$, where we put $t'=1/t$ and $\tau'=1/\tau$. We denote by $\infty$ the divisor $\{t=\infty\}$ in $\PP^1$, defined by the equation $t'=0$, as well as its inverse image in $\PP^1\times\wh\PP^1$, and similarly we consider the divisor $\wh\infty\subset\wh\PP^1$. We will use the picture described below.
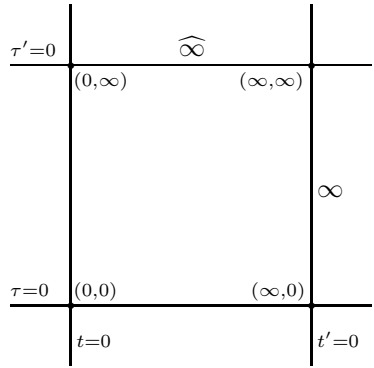
\begin{figure}[htb]
\setlength{\unitlength}{.8cm}
\begin{center}
\begin{picture}(6,6)(0,0)
\put(.93,.93){$\bbullet$}
\put(4.93,.93){$\bbullet$}
\put(.93,4.93){$\bbullet$}
\put(4.93,4.93){$\bbullet$}
\thicklines
\put(0,1){\line(1,0){6}}
\put(1,0){\line(0,1){6}}
\thinlines
\put(0,5){\line(1,0){6}}
\put(5,0){\line(0,1){6}}

\put(0,1.15){$\scst\tau=0$}
\put(0,5.15){$\scst\tau'=0$}
\put(1.1,.3){$\scst t=0$}
\put(5.1,.3){$\scst t'=0$}

\put(1.05,1.15){$\scst(0,0)$}
\put(4,1.15){$\scst(\infty,0)$}
\put(1.05,4.65){$\scst(0,\infty)$}
\put(3.8,4.65){$\scst(\infty,\infty)$}

\put(5.1,2.8){$\infty$}
\put(2.8,5.15){$\wh\infty$}
\end{picture}
\end{center}
\end{figure}

Let $Y$ be a complex manifold. We put $X=Y\times\PP^1$, $\wh X=Y\times\wh\PP^1$ and $Z=Y\times\PP^1\times\wh\PP^1$. The manifolds $X$ and $Z$ are equipped with a divisor (still denoted by) $\infty$, and $\wh X$ and~$Z$ are equipped with $\wh\infty$. We have projections
\begin{equation}\label{eq:projs}
\begin{array}{c}
\xymatrix{
&Z\ar[dl]_-p\ar[dr]^-{\wh p}&\\
X\ar[rd]_-q&&\wh X\ar[ld]^-{\wh q}\\
&Y&
}
\end{array}
\end{equation}

Let $\cM$ be a left $\cR_\cX$-module. We denote by $\wt\cM$ the localized module \hbox{$\cR_\cX[\infty]\otimes_{\cR_\cX}\cM$}. Then $p^+\wt\cM$ is a left $\cR_\cZ[*\infty]$-module. We consider its localization
$$
p^+\wt\cM[*\wh\infty]=\cR_\cZ[*(\infty\cup\wh\infty)]\otimes_{\cR_\cZ[*\infty]}p^+\wt\cM.
$$
We denote by $p^+\wt\cM[*\wh\infty]\otimes\ccE^{-t\tau/\hb}\index{$etth$@$\ccE^{-t\tau/\hb}$}$ the $\cO_\cZ[*(\infty\cup\wh\infty)]$-module $p^+\wt\cM[*\wh\infty]$ equipped with the twisted action of $\cR_\cZ$ described by the exponential factor: the $\cR_\cY$-action is unchanged, and, for any local section $m$ of $\cM$,
\begin{itemize}
\item
in the chart $(t,\tau)$,
\begin{equation}\label{eq:ttau}
\begin{split}
\partiall_t(m\otimes\ccE^{-t\tau/\hb})&=[(\partiall_t-\tau)m]\otimes\ccE^{-t\tau/\hb},\\
\partiall_\tau(m\otimes\ccE^{-t\tau/\hb})&=-tm\otimes\ccE^{-t\tau/\hb},
\end{split}
\end{equation}

\item
in the chart $(t',\tau)$,
\begin{equation}\label{eq:tprimetau}
\begin{split}
\partiall_{t'}(m\otimes\ccE^{-t\tau/\hb})&=[(\partiall_{t'}+\tau/t^{\prime2})m]\otimes\ccE^{-t\tau/\hb},\\
\partiall_\tau(m\otimes\ccE^{-t\tau/\hb})&=-m/t'\otimes\ccE^{-t\tau/\hb},
\end{split}
\end{equation}

\item
in the chart $(t,\tau')$,
\begin{equation}\label{eq:ttauprime}
\begin{split}
\partiall_t(m\otimes\ccE^{-t\tau/\hb})&=[(\partiall_t-1/\tau')m]\otimes\ccE^{-t\tau/\hb},\\
\partiall_{\tau'}(m\otimes\ccE^{-t\tau/\hb})&=tm/\tau^{\prime2}\otimes\ccE^{-t\tau/\hb},
\end{split}
\end{equation}

\item
in the chart $(t',\tau')$,
\begin{equation}\label{eq:tprimetauprime}
\begin{split}
\partiall_{t'}(m\otimes\ccE^{-t\tau/\hb})&=[(\partiall_{t'}+1/\tau't^{\prime2})m]\otimes\ccE^{-t\tau/\hb},\\
\partiall_{\tau'}(m\otimes\ccE^{-t\tau/\hb})&=m/t'\tau^{\prime2}\otimes\ccE^{-t\tau/\hb}.
\end{split}
\end{equation}
\end{itemize}

\begin{definition}\label{def:FcM}
The partial Fourier-Laplace transform $\wh\cM\index{$mcurlwh$@$\wh\cM$}$ of $\cM$ is the complex of $\cR_{\wh\cX}[*\wh\infty]$-modules
$$\wh p_+(p^+\wt\cM[*\wh\infty]\otimes\ccE^{-t\tau/\hb}).$$
\end{definition}

\subsection{Coherence properties}
We will give a criterion for the $\cR_{\wh\cX}[*\wh\infty]$-coherence of $\wh\cM$ when $\cM$ is $\cR_\cX$-coherent. As $\wh p$ is proper, it is enough to give a coherence criterion for $\cFcM\index{$mcurlFc$@$\cFcM$}\defin p^+\wt\cM[*\wh\infty]\otimes\ccE^{-t\tau/\hb}$.

\begin{proposition}\label{prop:Fouriercoh}
Let $\cM$ be a coherent $\cR_\cX$-module. Then $\cFcM$ is $\cR_\cZ[*\wh\infty]$-coherent. If moreover $\cM$ is good, then so is $\cFcM$, and therefore $\wh \cM=\wh p_+\cFcM$ is $\cR_{\wh\cX}[*\wh\infty]$-coherent.
\end{proposition}

\begin{proof}
The coherence is a local question near $t'=0$ (otherwise it is clearly satisfied) and it is enough to show that $\cFcM$ is locally finitely generated over $\cR_\cZ[*\wh\infty]$. Choose local generators $m_j$ of $\cM$ as a $\cR_\cX$-module. It is a matter of proving that, for any $k\in\NN$, $(\partiall_{t'}^km_j)\otimes\ccE^{-t\tau/\hb}$ and $t^{\prime-k}m_j\otimes\ccE^{-t\tau/\hb}$ belong to \hbox{$\cR_\cZ[*\wh\infty]\cdot (m_j\otimes\ccE^{-t\tau/\hb})$}.

Let us first compute in the chart $(t',\tau)$. We will use Formula \eqref{eq:tprimetau}. Up to a sign, the second term above is $\partiall_\tau^k(m_j\otimes\ccE^{-t\tau/\hb})$. The first one can be written as $(\partiall_{t'}^km_j)\otimes\ccE^{-t\tau/\hb}=(\partiall_{t'}-\tau\partiall_\tau^2)^k(m_j\otimes\ccE^{-t\tau/\hb})$. The computation in the chart $(t',\tau')$ is similar, using \eqref{eq:tprimetauprime}, as $\tau'$ acts in an invertible way.

The functor $\cM\mto \cFcM$ is exact and, for the property of being good, it is enough to show that if $\cL$ is a $\cO_\cX$-coherent submodule generating $\cM$ on a compact set $\cK\subset\cX$, then $p^*\cL[*\wh\infty]\otimes\ccE^{-t\tau/\hb}$ generates $\cFcM$ on $p^{-1}(\cK)$; this follows from the previous computation.
\end{proof}

\begin{remarque}\label{rem:alg}
When $\cM$ is good, we can compute the Fourier-Laplace transform in an algebraic way with respect to $t$ and $\tau$: we view $q_*\wt\cM$ as a coherent module over $q_*\cR_\cX[*\infty]=\cR_\cY[t]\langle\partiall_t\rangle$. Then $\wh q_*\wh\cM$ is the complex
\[
q_*\wt\cM[\tau]\To{\partiall_t-\tau}q_*\wt\cM[\tau],
\]
where the right-hand term is in degree $0$. In particular, the cohomology modules of this complex are $\cR_\cY[\tau]\langle\partiall_\tau\rangle$-coherent. (\Cf for instance \cite[Appendix A]{D-S02a} for an argument). Moreover, this complex has cohomology in degree $0$ only, and the cohomology is identified with $q_*\wt\cM$ as a $\cR_\cY$-module; the action of $\tau$ is induced by that of $\partiall_t$, and that of $\partiall_\tau$ by that of $-t$.
\end{remarque}

\begin{remarque}[Integrability of the Fourier-Laplace transform]
Let $\cM$ be a coherent $\cR_\cX$-module. Assume that $\cM$ is integrable (\cf \T\ref{subsec:intRmod}). Then $p^+\wt\cM\otimes\ccE^{-t\tau/\hb}$ is integrable as a $\cR_\cZ$-module. If moreover $\cM$ is good, then, using part of Proposition \ref{prop:imdirint}, we obtain the integrability of $\wh\cM$ as a $\cR_{\wh\cX}$-module.
\end{remarque}

\subsection{Fourier-Laplace transform of a sesquilinear pairing}\label{subsec:FC}
We will now forget the $\wh\infty$ divisor on $Z$ or $\wh X$, and still denote by $Z$ or $\wh X$ the sets $X\times\Afuh$ and $Y\times\Afuh$.

Assume that $\cM',\cM''$ are good $\cR_\cX$-modules. Let $C:\cMS'\otimes_{\cO_\bS}\ov{\cMS''}\to\Dbh{X}$ be a sesquilinear pairing. We will define a sesquilinear pairing between the corresponding Fourier-Laplace transforms:
\[
\index{$cwh$@$\wh C$}\wh C:\wh{\cMS'}\otimes_{\cO_\bS}\ov{\wh{\cMS''}}\to\Dbh{\wh X}.
\]

Firstly, define the sesquilinear pairing $p^+C:p^+\cMS'\otimes_{\cO_\bS}\ov{p^+\cMS''}\to\Dbh{Z}$ in the following way: local sections $m',m''$ of $p^+\cMS',p^+\cMS''$ can be written as $m'=\sum_i\phi_i\otimes m'_i$, $m''=\sum_j\psi_j\otimes m''_j$ with $\phi_i,\psi_j$ holomorphic functions on $\cZ$ and $m'_i,m''_j$ local sections of $\cMS',\cMS''$; put then
\begin{equation}\label{eq:Cmm}
\langle p^+C(m',\ov{m''}),\varphi\rangle\defin\sum_{i,j}\Big\langle C(m'_i,\ov{m''_j}),\int_p\phi_i\ov{\psi_j}\varphi\Big\rangle,
\end{equation}
for any $C^\infty$ (relative to $\bS$) form $\varphi$ on $Z\times\bS$ of maximal degree with compact support contained the open set of $\cZ$ where $m',m''$ are defined. That the previous expression does not depend on the decomposition of $m',m''$ and defines a sesquilinear pairing is easily verified: it is enough to show that, if $\sum_i\phi_i\otimes m'_i=0$, then the right-hand term in \eqref{eq:Cmm} vanishes; but, by flatness of $\cO_\cZ$ over $p^{-1}\cO_\cX$, the vector $\phi=(\phi_i)_i$ can be written as $\sum_ka_k\eta_k$, where each $\eta_k=(\eta_{k,i})_i$ is a vector of $\cO_{\cX}$-relations between the $m'_i$ in $\cM'$ and $a_k$ are local sections of $\cO_\cZ$; use then the $\cO_\cX$-linearity of $C$.

Secondly, extend $C$ as a sesquilinear pairing $\wt C$ on $\wt{\cMS'}\otimes_{\cO_\bS}\wt{\cMS''}$ with values in the sheaf of tempered distributions, that is, with poles along $\infty$. Define similarly $\wt{p^+C}$ (which is nothing but $p^+\wt C$). Such a distribution can be evaluated on forms $\varphi$ which are infinitely flat along $\infty$.

Remark that, for $\hb\in\bS$, we have $\big\vert e^{\hb\ov{t\tau}-t\tau/\hb}\big\vert=1$. The following lemma is standard (it is proved in the same way as one proves that the Fourier transform of a $C^\infty$ function with compact support is in the Schwartz class):

\begin{lemme}
Let $\varphi$ be a $C^\infty$ relative form of maximal degree on $Z\times\bS$ with compact support. Then $\int_pe^{\hb\ov{t\tau}-t\tau/\hb}\varphi$ is $C^\infty$ with compact support on $X\times\bS$ and is infinitely flat along $\infty$.\qed
\end{lemme}

For local sections $m',m''$ of $p^+\wt{\cMS'},p^+\wt{\cMS''}$ written as above and for $\varphi$ as in the lemma, it is meaningful to put
\[
\langle \cFC(m',\ov{m''}),\varphi\rangle\defin\sum_{i,j}\Big\langle \wt C(m'_i,\ov{m''_j}),\int_pe^{\hb\ov{t\tau}-t\tau/\hb}\phi_i\ov{\psi_j} \varphi\Big\rangle.
\]
This defines a sesquilinear pairing $\index{$cFc$@$\cFC$}\cFC:\cFcMS'\otimes_{\cO_\bS}\ov{\cFcMS''}\to\Dbh{Z}$. We can now define $\wh C=\wh p_\dag^0\cFC$.

\begin{remarque}[Behaviour with respect to adjunction]\label{rem:Fadjunction}
The formula above clearly implies that $(\cFC)^*={}^{\cF}\!(C^*)$. We hence have $(\wh C)^*=\wh{C^*}$.
\end{remarque}

It is possible to define $\wh C$ at the algebraic level considered in Remark \ref{rem:alg}. Notice first that $C$ defines a sesquilinear pairing $q_*\wt C$ on
$$
q_*\wt\cM'_{|\bS}\otimes_{\cO_\bS}\ov{q_*\wt\cM''_{|\bS}}
$$
which takes values in the sheaf on $Y_\RR\times\bS$ of distributions on $Y_\RR\times\Afu\times\bS$ which are \emph{tempered} with respect to the $t$-variable. Recall that $q_*\wt\cM',q_*\wt\cM''$ are $\cR_\cY[t]\langle\partiall_t\rangle$-modules and that their Fourier-Laplace transforms are the same objects viewed as $\cR_\cY[\tau]\langle\partiall_\tau\rangle$-modules \emph{via} the correspondence
\begin{equation}\label{eq:corrFourier}
\tau\longleftrightarrow \partiall_t,\quad \partiall_\tau\longleftrightarrow-t.
\end{equation}

Denote by $F$ the usual Fourier transform with kernel $\exp(\hb\ov {t\tau}-t\tau/\hb)\cdot\itwopi dt\wedge d\ov t$ for $\hb\in\bS$, sending $t$-tempered distributions on $Y_\RR\times\Afu\times\bS$ which are continuous with respect to $\hb$ to $\tau$-tempered distributions on $Y_\RR\times\Afuh\times\bS$ which are continuous with respect to $\hb$.

We define then $\wh q_*\wh C$ on
$$
\wh{q_*\wt\cM'_{|\bS}}\otimes_{\cO_\bS}\ov{\wh{q_*\wt\cM''_{|\bS}}}
$$
as the composition $F\circ q_*\wt C$. That $\wh q_*\wh C$ is $\cR_\cY[\tau]\langle\partiall_\tau\rangle\otimes\ov{\cR_\cY[\tau]\langle\partiall_\tau\rangle}$-linear follows from the fact that \eqref{eq:corrFourier} and its conjugate are the transformations that $F$ does.

\begin{lemme}
The analytization of $\wh q_*\wh C$ is equal to $\wh C$ defined as $\wh p_\dag^0\cFC$.\qed
\end{lemme}

\begin{remarque}[Integrability of the Fourier-Laplace transform of a sesquilinear pairing]
Let $\cM',\cM''$ be good $\cR_\cX$-modules which are integrable. Assume that the sesquilinear pairing $C$ is integrable (\cf \T\ref{subsec:intsesquilin}). Then $\cFC$ is integrable and, by Proposition \ref{prop:imdirint}, the pairing $\wh C$ is also integrable.
\end{remarque}

\section{Partial Fourier-Laplace transform and specialization}
As we are only interested in $\tau\neq\infty$, we continue to forget the divisor $\wh\infty$ and still denote by $Z$ or $\cZ$ the complement of this divisor.

\emph{A priori}, Proposition \ref{prop:Fouriercoh} does not restrict well to $\tau\!=\!\tau_o$ (of course, the problem is at $t'\!=\!0$). Indeed, we do not have a relation like \hbox{$\partiall_\tau(m\otimes\ccE^{-t\tau/\hb})=-m/t'\otimes\ccE^{-t\tau/\hb}$} to recover the polar part of $\wt\cM\otimes\ccE^{-t\tau_o/\hb}$ from the action of $\cR_\cX$. For instance, taking $\tau_o=0$, even with nice assumptions, $\wt\cM$ is not known to be $\cR_\cX$-coherent. We will introduce below an assumption which implies the $\cR_\cX$-coherence of $\wt\cM\otimes\ccE^{-t\tau_o/\hb}$ when $\tau_o\neq0$. For the coherence at $\tau=0$, we will need to consider the specialization at $\tau=0$ of $\cFcM$, and hence to prove first its strict specializability along $\tau=0$; for that, we will also need the same assumption. Let us introduce some notation.

Denote by $i_\infty$ the inclusion $Y\times\{\infty\}\hto X$. We will consider the functors $\psi_{\tau,\alpha}$ and $\psi_{t',\alpha}$ introduced in Lemma \ref{lem:uniciteMK}, as well as the functors $\Psi_{\tau,\alpha}$ and $\Psi_{t',\alpha}$ of Definition \ref{def:nearby}. We will denote by $\rN_\tau,\rN_{t'}$ the natural nilpotent endomorphisms on the corresponding nearby cycles modules. We denote by $\rM_\bbullet(\rN)$ the monodromy filtration of the nilpotent endomorphism $\rN$ and by $\gr\rN:\gr_\bbullet^{\rM}\to\gr_{\bbullet-2}^{\rM}$ the morphism induced by $\rN$. For $\ell\geq0$, $P\gr_\ell^{\rM}$ denotes the primitive part $\ker(\gr\rN)^{\ell+1}_{|\gr_\ell^{\rM}}$ of $\gr_\ell^{\rM}$ and $P\rM_\ell$ the inverse image of $P\gr_\ell^{\rM}$ by the natural projection $\rM_\ell\to\gr_\ell^{\rM}$. Recall that, in an abelian category, the primitive part $P\gr_0^{\rM}$ is equal to $\ker\rN/(\ker\rN\cap\im\rN)$. We will also denote by $\wt\cM_{\min}$ the minimal extension of $\wt\cM$ (\cf \T\ref{subsec:minext}).

\begin{proposition}\label{prop:Fstrictspe}
Assume that $\cM$ is strictly specializable and regular along $t'=0$ (\cf Definition~\ref{def:strictspe} and \T\ref{subsec:regularity}). Then,
\begin{enumerate}
\item\label{prop:Fstrictspe1}
for any $\tau_o\neq0$, the $\cR_\cX$-module $\wt\cM\otimes\ccE^{-t\tau_o/\hb}$ is $\cR_\cX$-coherent; it is also strictly specializable (but not regular in general) along $t'=0$, with a constant $V$-filtration, so that all \hbox{$\psi_{t',\alpha}(\wt\cM\otimes\ccE^{-t\tau_o/\hb})$} are identically $0$.
\end{enumerate}
Assume moreover that $\cM$ is strict. Then,
\begin{enumerate}\refstepcounter{enumi}
\item\label{prop:Fstrictspe2}
the $\cR_\cZ$-module $\cFcM\defin p^+\wt\cM\otimes\ccE^{-t\tau/\hb}$ is strictly specializable and regular along $\tau=\tau_o$ for any $\tau_o\in\Afuh$; it is equal to the \emph{minimal extension} of its localization along $\tau=0$;
\item\label{prop:Fstrictspe3}
if $\tau_o\neq0$, the $V$-filtration of $\cFcM$ along $\tau-\tau_o=0$ is given by
\[
V_k\cFcM=
\begin{cases}
\cFcM&\text{if }k\geq-1,\\
(\tau-\tau_o)^{-k+1}\cFcM&\text{if }k\leq-1;
\end{cases}
\]
we have
\[
\psi_{\tau-\tau_o,\alpha}\cFcM=
\begin{cases}
0&\text{if }\alpha\not\in-\NN-1,\\
\wt\cM\otimes\ccE^{-t\tau_o/\hb}&\text{if }\alpha\in-\NN-1.
\end{cases}
\]
\item\label{prop:Fstrictspe4}
If $\tau_o=0$, we have:
\begin{enumerate}
\item\label{prop:Fstrictspe4a}
for any $\alpha\neq-1$ with $\reel\alpha\in[-1,0[$, a functorial isomorphism on some neighbourhood of $\index{$dbf$@$\DD$}\DD\defin\{\module{\hb}\leq1\}$,
\[
\big(\Psi_{\tau,\alpha}\cFcM_{|\DD},\rN_\tau\big)\isom i_{\infty,+}\big(\psi_{t',\alpha}\wt\cM(-D_\alpha)_{|\DD},\rN_{t'}\big),
\]
where $D_\alpha$ is the divisor $1\cdot i$ if $\alpha'=-1$ and $\alpha''>0$, the divisor $1\cdot (-i)$ if $\alpha'=-1$ and $\alpha''<0$, and the empty divisor otherwise;
\item\label{prop:Fstrictspe4b}
for $\alpha=0$, a functorial isomorphism
\[
\big(\psi_{\tau,0}\cFcM,\rN_\tau\big)\isom i_{\infty,+}\big(\psi_{t',-1} \wt\cM,\rN_{t'}\big),
\]
\item\label{prop:Fstrictspe4c}
for $\alpha=-1$, two functorial exact sequences
\begin{gather*}
0\to i_{\infty,+}\ker\rN_{t'}\to\ker\rN_\tau\to\wt\cM_{\min}\to0\\
0\to\wt\cM_{\min}\to\coker\rN_\tau\to i_{\infty,+}\coker\rN_{t'}\to0,
\end{gather*}
inducing isomorphisms
\begin{align*}
i_{\infty,+}\ker\rN_{t'}&\isom\ker\rN_\tau\cap\im\rN_\tau\subset\ker\rN_\tau\\
\wt\cM_{\min}&\isom \ker\rN_\tau/(\ker\rN_\tau\cap\im\rN_\tau)\subset\coker\rN_\tau,
\end{align*}
such that the natural morphism $\ker\rN_\tau\to\coker\rN_\tau$ induces the identity on $\wt\cM_{\min}$.
\end{enumerate}
\end{enumerate}
\end{proposition}

\subsubsection*{Proof of \ref{prop:Fstrictspe}\eqref{prop:Fstrictspe1}}
Let us first prove the $\cR_\cX$-coherence of $\wt\cM\otimes\ccE^{-t\tau_o/\hb}$ when $\tau_o\neq0$. As this $\cR_\cX$-module is $\cR_\cX[*\infty]$-coherent by construction, it is enough to prove that it is locally finitely generated over $\cR_\cX$, and the only problem is at $t'=0$. We also work locally near $\hb_o\in\Omega_0$ and forget the exponent $(\hb_o)$ in the $V$-filtration along $t'=0$. Then, $\wt\cM=\cO_\cX[1/t']\otimes_{\cO_\cX}V_{<0}\cM$, equipped with its natural $\cR_\cX$-structure. By the regularity assumption, $V_{<0}\cM$ is $\cR_{\cX/\Afu}$-coherent, so we can choose finitely many $\cR_{\cX/\Afu}$-generators $m_i$ of $V_{<0}\cM$.

The regularity assumption implies that, for any $i$,
\[
t'\partiall_{t'} m_i\in \sum_j\cR_{\cX/\Afu}\cdot m_j.
\]
In $\wt\cM\otimes\ccE^{-t\tau_o/\hb}$, using \eqref{eq:tprimetau}, this is written as
\begin{equation}\label{eq:Vcoh}
(t'\partiall_{t'}-\tau_0/t')(m_i\otimes\ccE^{-t\tau_o/\hb})\in \sum_j\cR_{\cX/\Afu}\cdot (m_j\otimes\ccE^{-t\tau_o/\hb}),
\end{equation}
and therefore
\[
(\tau_o/t')(m_i\otimes\ccE^{-t\tau_o/\hb})\in \sum_jV_0\cR_\cX\cdot (m_j\otimes\ccE^{-t\tau_o/\hb}).
\]
It follows that $\wt\cM\otimes\ccE^{-t\tau_o/\hb}$ is $V_0\cR_\cX$-coherent, generated by the $m_i\otimes\ccE^{-t\tau_o/\hb}$. It is then obviously $\cR_\cX$-coherent. The previous relation also implies that $\tau_o(m_i\otimes\ccE^{-t\tau_o/\hb})\in t'\wt\cM\otimes\ccE^{-t\tau_o/\hb}$. Therefore, the constant $V$-filtration, defined by $V_a(\wt\cM\otimes\ccE^{-t\tau_o/\hb})=\wt\cM\otimes\ccE^{-t\tau_o/\hb}$ for any $a$, is good and has a Bernstein polynomial equal to $1$.

\subsubsection*{Proof of \ref{prop:Fstrictspe}\eqref{prop:Fstrictspe2} for $\tau_o\neq0$ and \ref{prop:Fstrictspe}\eqref{prop:Fstrictspe3}}
The analogue of Formula \eqref{eq:Vcoh} now reads
\[
(t'\partiall_{t'}+\tau\partiall_\tau)(m_i\otimes\ccE^{-t\tau/\hb})\in \sum_j\cR_{\cX/\Afu}\cdot(m_j\otimes\ccE^{-t\tau/\hb}).
\]
Therefore, the $\cR_{\cZ/\Afuh}$-module generated by the $m_j\otimes\ccE^{-t\tau/\hb}$ is $V_0\cR_\cZ$-coherent, where $V$ denotes the filtration relative to $\tau-\tau_o$. It is even $\cR_\cZ$-coherent if $\tau_o\neq0$, as $\tau$ is a unit near $\tau_o$, and this easily gives \ref{prop:Fstrictspe}\eqref{prop:Fstrictspe3}, therefore also \ref{prop:Fstrictspe}\eqref{prop:Fstrictspe2} when $\tau_o\neq0$.

\subsubsection*{Proof of \ref{prop:Fstrictspe}\eqref{prop:Fstrictspe2} for $\tau_o=0$}
Let us now consider the case where $\tau_o=0$. Then the previous argument gives the regularity of $\cFcM$ along $\tau=0$. We will now show the strict specializability along $\tau=0$. We will work near $\hb_o\in\Omega_0$ and forget the exponent $(\hb_o)$ in the $V$-filtrations relative to $\tau=0$ and to $t'=0$.

Out of $t'=0$ the result is easy: near $t=t_o$, formula \eqref{eq:ttau}, together with the strictness of $\cM$, implies that $\cFcM$ is strictly noncharacteristic along $\tau=0$, hence $\cFcM=V_{-1}\cFcM$ and $\psi_{\tau,-1}\cFcM=\cM$ (\cf \T\ref{sec:nonchar}).

We will now focus on $t'=0$. Denote by $V_\bbullet\cM$ the $V$-filtration of $\cM$ relative to $t'$ and put, for any $a\in[-1,0[$,
$$
V_{a+k}\wt\cM=\begin{cases}
V_a\cM&\text{if } k\leq0,\\
t^{\prime-k}V_a\cM&\text{if }k\geq0.
\end{cases}
$$
Each $V_a\wt\cM$ is a $V_0\cR_\cX$-coherent module and, by regularity, is also $\cR_{\cX/\Afu}$-coherent. We will now construct the $V$-filtration of $\cFcM$ along $\tau=0$. For any $a\in\RR$, put
\[
U_a\cFcM=\sum_{p\geq0}\partiall_{t'}^p \big[(p^*V_a\wt\cM)\otimes\ccE^{-t\tau/\hb}\big],
\]
\ie $U_a$ is the $\cR_{\cZ/\Afuh}$-module generated by $(p^*V_a\wt\cM)\otimes\ccE^{-t\tau/\hb}$ in $\cFcM$. Notice that, when we restrict to $t'\neq0$, we have for any $a\in\RR$,
\[
U_{a|t'\neq0}=\cFcM_{|t'\neq0}.
\]

\bgroup
\def\theenumi{\eqref{prop:Fstrictspe2}(\arabic{enumi})}
\def\labelenumi{\theenumi}
\begin{enumerate}
\item\label{prop:Fstrictspe21}
Clearly, $U_\bbullet$ is an increasing filtration of $\cFcM$ and each $U_a$ is $\cR_{\cZ/\Afuh}$-coherent for every $a\in\RR$.

\item\label{prop:Fstrictspe22}
$U_a$ is stable by $\tau\partiall_\tau$: indeed, for any local section $m$ of $V_a\wt\cM$, we have by \eqref{eq:tprimetau}:
\begin{align*}
(\tau\partiall_\tau)\partiall_{t'}^p(m\otimes\ccE^{-t\tau/\hb})&=\partiall_{t'}^p (\tau\partiall_\tau)(m\otimes\ccE^{-t\tau/\hb})\\
&=\partiall_{t'}^p\big[t'\partiall_{t'}(m\otimes\ccE^{-t\tau/\hb})-(t'\partiall_{t'}m)\otimes\ccE^{-t\tau/\hb}\big]\\
&=\partiall_{t'}^{p+1}(t'm\otimes\ccE^{-t\tau/\hb})-\partiall_{t'}^p\big[(\partiall_{t'}t'm)\otimes\ccE^{-t\tau/\hb}\big].
\end{align*}
The first term in the RHS is in $U_{a-1}$ and the second one is in $U_a$, as $V_a\wt\cM$ is stable by $\partiall_{t'}t'$.

\item\label{prop:Fstrictspe23}
For any $a\in\RR$, we have $U_{a+1}=U_a+\partiall_\tau U_a$: indeed, for $m$ as above, we have
\[
\partiall_\tau\cdot\partiall_{t'}^p(m\otimes\ccE^{-t\tau/\hb}) =-\partiall_{t'}^p\Big(\frac1{t'}m\otimes\ccE^{-t\tau/\hb}\Big)\in U_{a+1},
\]
hence $\partiall_\tau U_a\subset U_{a+1}$; applying this equality the in the other way gives the desired equality. This also shows that $\partiall_\tau:\gr_a^U\cFcM\to\gr_{a+1}^U\cFcM$ is an isomorphism for any $a\in\RR$.

\item\label{prop:Fstrictspe24}
For any $a\in\RR$, we have $\tau U_a\subset U_{a-1}$: indeed, one has, for $m$ as above
\begin{align*}
\tau(m\otimes\ccE^{-t\tau/\hb})&=t^{\prime2}\partiall_{t'}(m\otimes\ccE^{-t\tau/\hb})-(t^{\prime2}\partiall_{t'}m)\otimes\ccE^{-t\tau/\hb}\\
&=\partiall_{t'}(t^{\prime2}m\otimes\ccE^{-t\tau/\hb})-(\partiall_{t'}t^{\prime2}m)\otimes\ccE^{-t\tau/\hb};
\end{align*}
the first term of the RHS clearly belongs to $U_{a-2}$ and the second one to $U_{a-1}$.

\item\label{prop:Fstrictspe25}
Denote by $b_a(s)$ the minimal polynomial of $-\partiall_{t'}t'$ on $\gr_a^V\wt\cM$. Then, for $m$ as above, we have
\[
-(\partiall_{t'}t'+\tau\partiall_\tau)(m\otimes\ccE^{-t\tau/\hb})=-(\partiall_{t'}t'm)\otimes\ccE^{-t\tau/\hb}
\]
after \eqref{eq:tprimetau}. Therefore, we have $b_a(-[\partiall_{t'}t'+\tau\partiall_\tau])(m\otimes\ccE^{-t\tau/\hb})\in U_{<a}$. Using that $\partiall_{t'}t'(m\otimes\ccE^{-t\tau/\hb})=\partiall_{t'}(t'm\otimes\ccE^{-t\tau/\hb})\in U_{a-1}$ by definition, we deduce that $b_a(-\tau\partiall_\tau)(m\otimes\ccE^{-t\tau/\hb})\in U_{<a}$. Therefore, $b_a(-\tau\partiall_\tau)U_a\subset U_{<a}$.

\item\label{prop:Fstrictspe26}
We will now identify $U_a/U_{<a}$ with $\gr_a^V\wt\cM[\eta]\defin\CC[\eta]\otimes_\CC\gr_a^V\wt\cM$, where $\eta$ is a new variable. Notice first that both objects are supported on $\{t'=0\}$. Consider the map
\begin{align*}
V_a\wt\cM[\eta]&\to U_a\\
\sum_pm_p\eta^p&\mto\sum_p\partiall_{t'}^p(m_p\otimes\ccE^{-t\tau/\hb}).
\end{align*}
Its composition with the natural projection $U_a\to U_a/U_{<a}$ induces a surjective mapping $\gr_a^V\wt\cM[\eta]\to U_a/U_{<a}$. In order to show that it is injective, it is enough to show that, if $\sum_p\partiall_{t'}^p(m_p\otimes\ccE^{-t\tau/\hb})$ belongs to $U_{<a}$, then each $m_p$ belongs to $V_{<a}\wt\cM$. For that purpose, it is enough to work with an algebraic version of $U_a$, where ``$p^*$'' means ``$\otimes_\CC\CC[\tau]$''. Notice that, if a local section $\sum_{\ell=0}^r\tau^\ell(n_\ell\otimes\ccE^{-t\tau/\hb})$ of $\wt\cM[\tau]\otimes\ccE^{-t\tau/\hb}$ belongs to $U_a$, then the dominant coefficient $n_r$ is a local section of $V_{a+2r}\wt\cM$ (by using that $\partiall_{t'}(n\otimes\ccE^{-t\tau/\hb})=(\partiall_{t'}n)\otimes\ccE^{-t\tau/\hb}-\tau((n/t^{\prime2})\otimes\ccE^{-t\tau/\hb})$). Remark then that, using \eqref{eq:tprimetau}, $\sum_{p=0}^q\partiall_{t'}^p(m_p\otimes\ccE^{-t\tau/\hb})$ is a polynomial of degree $q$ in $\tau$ with dominant coefficient $\pm(\tau^q/t^{\prime2q})(m_q\otimes\ccE^{-t\tau/\hb})$. If the sum belongs to $U_{<a}$, this implies that $m_q/t^{\prime2q}\in V_{<a+2q}\wt\cM$, \ie $m_q\in V_{<a}\wt\cM$. Therefore, by induction on~$q$, all coefficients $m_p$ are local sections of $V_{<a}\wt\cM$, as was to be shown.

Let us describe the $\cR_\cX[\tau\partiall_\tau]$-module structure on $\gr_a^V\wt\cM[\eta]$ coming from the identification with $U_a/U_{<a}$. First, the $\cR_\cY$-module structure is the natural one on $\gr_a^V\wt\cM$, naturally extended to $\gr_a^V\wt\cM[\eta]$. Then one checks that
\begin{gather}\label{eq:UV1}
\partiall_{t'}\sum_pm_p\eta^p=\eta\sum_pm_p\eta^p,\quad
t'\sum_pm_p\eta^p=-\partiall_\eta\sum_pm_p\eta^p,\\
\label{eq:UV2}
\tau\partiall_\tau\sum_pm_p\eta^p=\sum_p(\partiall_t't')(m_p)\eta^p.
\end{gather}
In particular, any element is killed by some power of $t'$.

\item\label{prop:Fstrictspe27}
Consider the filtration $V_\bbullet\cFcM$ defined for $a\in[-1,0[$ and $k\in\ZZ$ by
\[
V_{a+k}\cFcM=
\begin{cases}
U_{a+1+k}&\text{if }k\geq0,\\
\tau^{-k}U_{a+1}&\text{if }k\leq0.
\end{cases}
\]
This is a $V$-filtration relative to $\tau$ on $\cFcM$, by \ref{prop:Fstrictspe21}, \ref{prop:Fstrictspe22}, \ref{prop:Fstrictspe23} and \ref{prop:Fstrictspe24}. It is good, by the equality in \ref{prop:Fstrictspe23} and because $\tau V_a=V_{a-1}$ for $a<0$ by definition.

On $\gr_a^V\cFcM$, for $a>-1$, there is a minimal polynomial of the right form for $-\partiall_\tau\tau$, by \ref{prop:Fstrictspe25} (here is the need for a shift by $1$ between $U$ and $V$), and strictness follows from \ref{prop:Fstrictspe26} and the strictness of $\gr_a^V\wt\cM$, which is by assumption.

It therefore remains to analyze $\gr_a^V\cFcM$ for $a\leq-1$.

\item\label{prop:Fstrictspe28}
We will analyze $\gr_{-1}^V\cFcM=U_{0}/\tau U_{<1}$ through the following two diagrams of exact sequences, where the non labelled maps are the natural ones:
\begin{equation}\label{eq:suiteexU}
\begin{array}{c}
\xymatrix@=5mm{
&0\ar[d]\\
&(U_{<0}\cap\tau\cFcM)/\tau U_{<1}\ar[d]\\
0\ar[r]& U_{<0}/\tau U_{<1}\ar[r]\ar[d]& U_{0}/\tau U_{<1}\ar[r]& U_{0}/U_{<0} \ar[r]& 0\\
&U_{<0}/(U_{<0}\cap\tau\cFcM)\ar[d]\\
&0
}
\end{array}
\end{equation}
and
\begin{equation}\label{eq:suiteexUc}
\begin{array}{c}
\xymatrix@=5mm{
&&&&0\ar[d]\\
&&&&U_{<0}/(U_{<0}\cap\tau U_1)\ar[d]\\
0\ar[r]&U_0/U_{<0}\ar[rr]^-{\tau\partiall_\tau}&&U_{0}/\tau U_{<1}\ar[r]& U_0/\tau U_1\ar[r]\ar[d]&0\\
&&&& U_0/(\tau U_1+U_{<0})\ar[d]\\
&&&&0
}
\end{array}
\end{equation}
Notice that, in \eqref{eq:suiteexUc}, $\tau\partiall_\tau$ is injective because it is the composition
\begin{equation}\label{eq:compo}
U_0/U_{<0}\To{\partiall_\tau} U_1/U_{<1}\To{\tau}U_{0}/\tau U_{<1},
\end{equation}
$\partiall_\tau$ is an isomorphism (\cf \ref{prop:Fstrictspe23}) and $\tau$ is injective, as it acts injectively on $\cFcM$. Recall that $(\gr_0^U\cFcM, \tau\partiall_\tau)$ is identified, by \ref{prop:Fstrictspe26}, with $i_{\infty,+}(\gr_0^V\wt\cM,\partiall_{t'}t')$. Notice also that $\tau\partiall_\tau$ vanishes on $U_{<0}/\tau U_{<1}$ (\resp on $U_0/\tau U_1$), as as $\partiall_\tau U_{<0}\subset U_{<1}$ (\resp $\partiall_\tau U_{0}\subset U_{1}$). It remains therefore to prove the strictness of $U_{<0}/\tau U_{<1}$ to get the desired properties for $\gr_{-1}^V\cFcM$. We denote by $\rN_{t'}$ the action of $-t'\partiall_{t'}$ on $\gr_{-1}^V\wt\cM$ (by strictness, $\ker\rN_{t'}$ is equal to the kernel of $-t'\partiall_{t'}$ acting on $\psi_{t',-1}\wt\cM\subset\gr_{-1}^V\wt\cM$). The strictness of $\gr_{-1}^V\cFcM$ follows then from the strictness of $i_{\infty,+}\psi_{t',-1}\wt\cM$, that of $\wt\cM_{\min}$ and the first two lines of the lemma below, applied to the diagram \eqref{eq:suiteexU}.

\begin{lemme}\label{lem:UU}
We have functorial isomorphisms of $\cR_\cX$-modules:
\begin{gather*}
U_{<0}/(U_{<0}\cap\tau U_1)= U_{<0}/(U_{<0}\cap\tau \cFcM)\isom\wt\cM_{\min}\\
i_{\infty,+}\ker\rN_{t'}\isom (U_{<0}\cap\tau\cFcM)/\tau U_{<1}\\
i_{\infty,+}\coker\rN_{t'}\isom U_0/(\tau U_1+U_{<0}).
\end{gather*}
\end{lemme}

\begin{proof}
For $m_0,\dots,m_p\in\wt\cM$, we can write
\begin{multline}\label{eq:mn}
m_0\otimes\ccE^{-t\tau/\hb}+\partiall_{t'}(m_1\otimes\ccE^{-t\tau/\hb})+\cdots+\partiall_{t'}^p(m_p\otimes\ccE^{-t\tau/\hb})\\
=n_0\otimes\ccE^{-t\tau/\hb}-\tau\Big[(n_1/t^{\prime2})\otimes\ccE^{-t\tau/\hb}+\cdots+\partiall_{t'}^{p-1}\big((n_p/t^{\prime2})\otimes\ccE^{-t\tau/\hb}\big)\Big]
\end{multline}
with
\begin{equation}\label{eq:changemn}
\begin{aligned}
n_p&=m_p&m_p&=n_p\\
n_{p-1}&=m_{p-1}+\partiall_{t'}m_p&m_{p-1}&=n_{p-1}-\partiall_{t'}n_p\\
&\hspace*{2mm}\vdots&&\hspace*{2mm}\vdots\\
n_1&=m_1+\partiall_{t'}m_2+\cdots+\partiall_{t'}^{p-1}m_p&m_1&=n_1-\partiall_{t'}n_2\\
n_0&=m_0+\partiall_{t'}m_1+\cdots+\partiall_{t'}^pm_p&m_0&=n_0-\partiall_{t'}n_1
\end{aligned}
\end{equation}

Sending an element to its constant term in its $\tau$ expansion gives an injective morphism $U_{<0}/(U_{<0}\cap \tau\cFcM)\to\wt\cM$. Formulas \eqref{eq:mn} and \eqref{eq:changemn} show that the image of this morphism is the $\cR_\cX$-submodule of $\wt\cM$ generated by $V_{<0}\wt\cM$: this is by definition the minimal extension of $\wt\cM$ across $t'=0$.

Let us show that
\begin{equation}\label{eq:UtauU}
U_{<0}\cap \tau U_1=U_{<0}\cap \tau\cFcM.
\end{equation}
Consider a local section of $U_{<0}\cap \tau\cFcM$, written as in \eqref{eq:mn}; it satisfies thus $m_0,\dots,m_p\in V_{<0}\wt\cM$ and $n_0=0$; then $\partiall_{t'}n_1=-m_0\in V_{<0}\wt\cM$. This implies that $n_1$ is a local section of $V_{-1}\wt\cM$: indeed, the condition on $n_1$ is equivalent to $t'\partiall_{t'}n_1\in V_{<-1}\wt\cM$; use then that, by strictness of $\gr_a^V\wt\cM$, $t'\partiall_{t'}$ acts injectively on $\gr_a^V\wt\cM$ if $a\neq-1$. Therefore, $(n_1/t^{\prime2})\otimes\ccE^{-t\tau/\hb}\in U_1$. We can now assume that $n_1=0$ and thus $\partiall_{t'}n_2\in V_{<0}\wt\cM$... hence \eqref{eq:UtauU}, and the first line of the lemma. Notice moreover that the class of each $n_j$ in $\gr_{-1}^V\wt\cM$ is in $\ker\rN_{t'}$.

Let $\eta$ be a new variable. We define a morphism
\[
\ker \rN_{t'}[\eta]\to U_{<0}/\tau U_{<1}
\]
by the rule
\begin{equation}\label{eq:sumnj}
\sum_{j\geq1}[n_j]\eta^{j-1}\mto -\tau\Big[(n_1/t^{\prime2})\otimes\ccE^{-t\tau/\hb}+\cdots+\partiall_{t'}^{p-1}\big((n_p/t^{\prime2})\otimes\ccE^{-t\tau/\hb}\big)\Big],
\end{equation}
by taking some lifting $n_j$ of each $[n_j]\in\ker \rN_{t'}\subset\gr_{-1}^V\wt\cM$ in $V_{-1}\wt\cM$.
\begin{itemize}
\item
This morphism is well defined: using \eqref{eq:mn}, write
\[
-\tau\partiall_{t'}^{j-1}\big((n_j/t^{\prime2})\otimes\ccE^{-t\tau/\hb}\big)= \partiall_{t'}^j\big(n_j\otimes\ccE^{-t\tau/\hb}\big)-\partiall_{t'}^{j-1}\big((\partiall_{t'}n_j)\otimes\ccE^{-t\tau/\hb}\big);
\]
that $[n_j]$ belongs to $\ker \rN_{t'}$ is equivalent to $t'\partiall_{t'}n_j\in V_{<-1}\wt\cM$; therefore, both $n_j$ and $\partiall_{t'}n_j$ belong to $V_{<0}\wt\cM$; moreover, if $n_j\in V_{<-1}\wt\cM$, so that $n_j/t^{\prime2}\in V_{<1}\wt\cM$, the image is in $\tau U_{<1}$.
\item
This morphism is injective: as we have seen in \ref{prop:Fstrictspe26}, the term between brackets in \eqref{eq:sumnj} belongs to $U_{<1}$ if and only if each $n_j/t^{\prime2}$ belongs to $V_{<1}\wt\cM$, \ie each $n_j$ is in $V_{<-1}\wt\cM$.
\item
The image of this morphism is equal to $(U_{<0}\cap \tau\cFcM)/\tau U_{<1}$: this was shown in the proof of \eqref{eq:UtauU}.
\end{itemize}
As in \ref{prop:Fstrictspe26}, we can identify $\ker \rN_{t'}[\eta]$ with $i_{\infty,+}\ker \rN_{t'}$ and the morphism is seen to be $\cR_\cX$-linear.

Let us now consider the third line of the lemma. We identify $U_0/(\tau U_1+U_{<0})$ with the cokernel of $\tau:\gr_1^U\to \gr_0^U$ or, equivalently, to that of $\tau\partiall_\tau:\gr_0^U\to\gr_0^U$. By~\ref{prop:Fstrictspe26}, it is identified with $i_{\infty,+}\coker \partiall_{t'}t'$ acting on $i_{\infty,+}\gr_0^V\wt\cM$. Use now the isomorphism $t':\gr_0^V\wt\cM\to\gr_{-1}^V\wt\cM$ to conclude.
\end{proof}

\item\label{prop:Fstrictspe29}
We will now prove that all the $\gr_a^V\cFcM$ for $a\leq-1$ are strict and have a Bernstein polynomial. In \ref{prop:Fstrictspe28} we have proved this for $a=-1$.

Choose $a<-1$. It follows from the definition of $V_\bbullet\cFcM$ that
\begin{equation}\label{eq:tauonto}
\tau:\gr_{a+1}^V\cFcM\to \gr_a^V\cFcM
\end{equation}
is onto. Therefore, by decreasing induction on $a$, we have a Bernstein relation on each $\gr_a^V\cFcM$. It remains to prove the strictness of such a module. This is also done by decreasing induction on $a$, as it is known to be true for any $a\in[-1,0[$. It is enough to show that \eqref{eq:tauonto} is also injective for any $a<-1$, and it is also enough to show that
\[
\partiall_\tau\tau:\gr_{a+1}^V\cFcM\to \gr_{a+1}^V\cFcM.
\]
is injective. If a section $m$ satisfies $\partiall_\tau\tau m=0$, it also satisfies $\prod(\alpha\star\hb)^{\nu_\alpha}m=0$, where the product is taken on a set of $\alpha\in\CC$ with $\ell_{\hb_o}(\alpha)=a+1<0$ and $\nu_\alpha\in\NN$. Such a set does not contain $0$ and the function $\hb\mto\prod(\alpha\star\hb)^{\nu_\alpha}$ is not identically $0$. By induction, $\gr_{a+1}^V\cFcM$ is strict. Therefore, $m=0$, hence the injectivity.

\item\label{prop:Fstrictspe210}
By construction, the filtration $V_\bbullet\cFcM$ satisfies moreover that
\begin{itemize}
\item[$\bbullet$]
$\tau:\gr_{a}^V\cFcM\to \gr_{a-1}^V\cFcM$ is onto for any $a<0$,
\item[$\bbullet$]
$\partiall_\tau:\gr_{a}^V\cFcM\to \gr_{a+1}^V\cFcM$ is onto for any $a\geq-1$.
\end{itemize}
This implies that both conditions \ref{def:strictspe}\eqref{def:strictspeb} and \eqref{def:strictspec} are satisfied, and that moreover the morphism $\can_\tau$ introduced in Remark \ref{rem:psi}\eqref{rem:Ncanvar2} is \emph{onto}. Notice also that the morphism $\var_\tau$ is injective: indeed, this means that $\tau:\gr_{0}^V\cFcM\to \gr_{-1}^V\cFcM$ is injective, or equivalently that $\tau:U_1/U_{<1}\to U_0/\tau U_{<1}$ is injective, which has been seen after \eqref{eq:compo}.

In other words, we have shown that $\cFcM$ is strictly specializable along $\tau=0$ in the sense of Definition \ref{def:strictspe} and that is is equal to the \emph{minimal extension} of its localization along $\tau=0$, as defined in \T\ref{subsec:minext}.
\end{enumerate}
\egroup

\subsubsection*{Proof of \ref{prop:Fstrictspe}\eqref{prop:Fstrictspe4}}
Now that $\cFcM$ is known to be strictly specializable along $\tau=0$, the $\cR_\cX$-modules $\psi_{\tau,\alpha}\cFcM$ are defined. We can compare them with $i_{\infty,+}\psi_{t',\alpha}\wt\cM$.

\bgroup
\def\theenumi{\eqref{prop:Fstrictspe4}(\arabic{enumi})}
\def\labelenumi{\theenumi}
\begin{enumerate}
\item\label{prop:Fstrictspe41}
For any $\hb_o\in\Omega_0$, we have a natural morphism, defined locally near $\hb_o$ (putting $a=\ell_{\hb_o}(\alpha)$)
\begin{equation}\label{eq:idefines}
\psi_{\tau,\alpha}\cFcM\hto\gr_a^V\cFcM\to\gr_{a+1}^U\cFcM\isom i_{\infty,+}\gr_{a+1}^V\wt\cM\xrightarrow[\ts\sim]{\ts~i_{\infty,+}t'~}i_{\infty,+}\gr_a^V\wt\cM,
\end{equation}
which takes values in $i_{\infty,+}\psi_{t',\alpha}\wt\cM$. One verifies that the various morphisms glue together in a well defined morphism $\psi_{\tau,\alpha}\cFcM\to i_{\infty,+}\psi_{t',\alpha}\wt\cM$.

\begin{lemme}\label{lem:morphisme}
Near any $\hb_o\in\DD$, the natural morphism $\psi_{\tau,\alpha}\cFcM\to\gr_{a+1}^U\cFcM$ ($a=\ell_{\hb_o}(\alpha)$) is injective for any $\alpha\in\CC\moins(-\NN^*)$ and, if $a\geq-1$, $\psi_{\tau,\alpha}\cFcM\to i_{\infty,+}\psi_{t',\alpha}\wt\cM$ is an isomorphism near $\hb_o$.
\end{lemme}

\begin{proof}
If $a>-1$, this has been proved in \ref{prop:Fstrictspe26}. Assume that $a=-1$ (and $\alpha\not\in-\NN^*$). If we decompose the horizontal sequence \eqref{eq:suiteexU} with respect to the eigenvalues of $-\tau\partiall_\tau$, we get that, for any $\alpha\neq-1$ with $\ell_{\hb_o}(\alpha)=-1$, the natural morphism
\[
\psi_{\tau,\alpha}\cFcM\to U_{0}/U_{<0}
\]
is injective and, according to \ref{prop:Fstrictspe26}, we have an isomorphism
\[
\psi_{\tau,\alpha}\cFcM\isom i_{\infty,+}\psi_{t',\alpha+1}\wt\cM.
\]

Assume now that $a<-1$. Let $k\geq0$ be such that $b=a+k\in[-1,0[$. We prove the result by induction on $k$, knowing that it is true for $k=0$. By induction, we have a commutative diagram
\[
\xymatrix{
\psi_{\tau,\alpha+1}\cFcM\ar@{^{ (}->}[r]\ar[d]^-\wr_-\tau& \gr_{a+2}^U\cFcM\\
\psi_{\tau,\alpha}\cFcM\ar[r]&\gr_{a+1}^U\cFcM\ar[u]_-{\partiall_\tau}^-\wr
}
\]
showing that the lower horizontal arrow is injective if and only if $\partiall_\tau\tau$ is injective on $\psi_{\tau,\alpha+1}\cFcM$, which follows from strictness if $(\alpha+1)\star\hb\not\equiv0$, that is, if $\alpha\neq-1$.
\end{proof}

\item\label{prop:Fstrictspe41b}Proof of \ref{prop:Fstrictspe}\eqref{prop:Fstrictspe4b}.
When $\alpha=0$, the proof follows from Lemma \ref{lem:morphisme}.

\item\label{prop:Fstrictspe41c}
Assume now that $\alpha\neq-1$ satisfies $\reel\alpha\in[-1,0[$. We wish to show that \eqref{eq:idefines} induces an isomorphism
\begin{equation}\label{eq:isoPsineqa}
\psi_{\tau,\alpha}\cFcM_{|\DD}\isom i_{\infty,+}\psi_{t',\alpha}\wt\cM(-D_\alpha)_{|\DD}.
\end{equation}
This is a local question with respect to $\hb\in\DD$.

Clearly, the image of $\psi_{\tau,\alpha}\cFcM\to \gr_{a+1}^U\cFcM$ is contained in $\ker[(\partiall_\tau\tau+\alpha\star\hb)^N:\gr_{a+1}^U\cFcM\to \gr_{a+1}^U\cFcM]$, for $N\gg0$ and is equal to this submodule if $a\geq-1$.

If $a<-1$ and if $k\geq1$ is such that $a+k\in[-1,0[$, the image is identified with
\[
\im(\tau^k\partiall_\tau^k):\ker(\partiall_\tau\tau+\alpha\star\hb)^N\to\ker(\partiall_\tau\tau+\alpha\star\hb)^N,
\]
and it is identified with the image of the multiplication by $\prod_{j=1}^{k}(\alpha+j)\star\hb$ on this module. For $j=1,\dots,k$, the number $\beta=\alpha+j$ satisfies $\reel\beta\geq0$, $\beta\neq0$ and $\ell_{\hb_o}(\beta)<0$. Then $\beta\star\hb=0$ has a solution in $\DD$ iff $\reel\beta=0$, and this solution is $\pm i$. This occurs iff $\reel\alpha=-1$ and $j=1$. In conclusion, the image of $\psi_{\tau,\alpha}\cFcM_{|\DD}$ in $i_{\infty,+}\psi_{\tau,\alpha}\wt\cM_{|\DD}$, is equal to the image of the multiplication by $(\alpha+1)\star\hb$ on $i_{\infty,+}\psi_{\tau,\alpha}\wt\cM_{|\DD}$. As we assume that $\ell_{\hb_o}(\alpha)<-1$, the divisor of $\hb\mto(\alpha+1)\star\hb$ coincides, near $\hb_o$, with the divisor $D_\alpha$, hence \eqref{eq:isoPsineqa}.

\item\label{prop:Fstrictspe42}
We now show that there is no difference between $\psi_{\tau,\alpha}\cFcM$ and $\Psi_{\tau,\alpha}\cFcM$ on some neighbourhood of $\DD$.

\begin{lemme}\label{lem:Psipsi}
Assume that $\alpha\neq-1$ and $\alpha'\defin\reel\alpha\in[-1,0[$. Then the natural inclusion (\cf Lemma \ref{lem:psiMMt}\eqref{lem:psiMMt1}) $\psi_{\tau,\alpha}\cFcM_\DD\hto\Psi_{\tau,\alpha}\cFcM_\DD$ is an isomorphism.
\end{lemme}

\begin{proof}
The question is local near points $\hb\in\DD$ such that $\ell_{\hb}(\alpha)\geq0$, otherwise the result follows from Lemma \ref{lem:localstrictspe}. Fix $\hb_o$ such that $\ell_{\hb_o}(\alpha)\geq0$ and let $k\geq1$ be such that $\ell_{\hb_o}(\alpha-k)\in[-1,0[$. We have a commutative diagram
\[
\xymatrix{
\psi_{\tau,\alpha}\cFcM\ar@{^{ (}->}[r]\ar[d]_-{\tau^k}&\Psi_{\tau,\alpha}\cFcM\ar[d]_-{\wr}^-{\tau^k}\\
\psi_{\tau,\alpha-k}\cFcM\ar[r]^-\sim&\Psi_{\tau,\alpha-k}\cFcM
}
\]
and, as $a\defin\ell_{\hb_o}(\alpha)$ and $a-k$ are $\geq-1$ and $\alpha\neq-1$, $\psi_{\tau,\alpha}\cFcM$ (\resp $\psi_{\tau,\alpha-k}\cFcM$) is contained in $\gr_{a+1}^U\cFcM$ (\resp in $\gr_{a+1-k}^U\cFcM$), using the local filtration $U$ near $\hb_o$. It follows (\cf \ref{prop:Fstrictspe23}) that $\partiall_\tau^k:\psi_{\tau,\alpha-k}\cFcM\to\psi_{\tau,\alpha}\cFcM$ is an isomorphism. Therefore, the image of $\psi_{\tau,\alpha}\cFcM$ in $\Psi_{\tau,\alpha}\cFcM$ is identified with the image of $\partiall_\tau^k\tau^k$ acting on $\Psi_{\tau,\alpha}\cFcM$. Using the nilpotent endomorphism $\rN_\tau=-(\partiall_\tau\tau+\alpha\star\hb)$, we write $\partiall_\tau^k\tau^k$ as \hbox{$(-1)^k(\rN_\tau+\alpha\star\hb)\cdots(\rN_\tau+(\alpha-k+1)\star\hb)$}. The proof of the lemma will be complete if we show that none of the $(\alpha-j)\star\hb_o$ ($j=0,\dots,k-1$) vanishes (assuming that $\hb_o\in\DD$).

Notice that $\beta\defin\alpha-j$ satisfies $\beta'<0$ and $\beta'-\imhb_o\beta''\geq0$. Assume that $\beta\star\hb_o=0$. By the previous conditions, we must have $\beta''\neq0$ and $\hb_o\neq0$, and the only possibility for $\hb_o$ is then $\hb_o=i\imhb_o$ and $\imhb_o=\frac{\beta'-\sqrt{\beta^{\prime2}+\beta^{\prime\prime2}}}{\beta''}$. Now, the condition $\beta'<0$ implies $\module{\imhb_o}>1$, so $\hb_o\not\in\DD$.
\end{proof}

\item Proof of \ref{prop:Fstrictspe}\eqref{prop:Fstrictspe4a}.
It follows from \eqref{eq:isoPsineqa} and Lemma \ref{lem:Psipsi} that we have a functorial isomorphism
\begin{equation}\label{eq:isoPsineq}
\Psi_{\tau,\alpha}\cFcM_{|\DD}\to i_{\infty,+}\psi_{\tau,\alpha}\wt\cM(-D_\alpha)_{|\DD}
\end{equation}
when $\alpha\neq-1$ satisfies $\reel\alpha\in[-1,0[$. This ends the proof of \ref{prop:Fstrictspe}\eqref{prop:Fstrictspe4} when $\alpha\neq-1$.

\item\label{prop:Fstrictspe43} Proof of \ref{prop:Fstrictspe}\eqref{prop:Fstrictspe4c}.
Let us now consider the case when $\alpha=-1$. The two exact sequences that we consider are the vertical exact sequences in \eqref{eq:suiteexU} and \eqref{eq:suiteexUc}, according to Lemma \ref{lem:UU}.

For the second assertion, notice first that, as the image of $\im\rN_\tau\cap\ker\rN_\tau$ in $\wt\cM_{\min}$ is supported on \hbox{$\{t'=0\}$}, it is zero by the definition of the minimal extension, hence we have an inclusion $\im\rN_\tau\cap\ker\rN_\tau\subset i_{\infty,+}\ker\rN_{t'}$. To prove $i_{\infty,+}\ker\rN_{t'}\subset\im\rN_\tau$, remark that the image of \eqref{eq:sumnj} is in $\tau (U_1/U_{<1})$, hence in $\tau\psi_{\tau,0}\cFcM$, that is, in $\im\var_\tau$, hence in $\im\rN_\tau$.

The last assertion is nothing but the identification $U_{<0}\cap\tau\cFcM=U_{<0}\cap\tau U_1$ of Lemma \ref{lem:UU}.
\qed
\end{enumerate}
\egroup

\section[Partial Fourier-Laplace transform of twistor $\cD$-modules]
{Partial Fourier-Laplace transform of regular twistor $\cD$-modules}
The main result of this chapter is:

\begin{theoreme}\label{th:Fstrict}
Let $(\cT,\cS)=(\cM',\cM'',C,\cS)$ be an object of $\MTr(X,w)^\rp$. Then, along $\tau=0$, $\wh\cM'$ and $\wh\cM''$ are strictly specializable, regular and S-decomposable (\cf Definition \ref{def:strictdec}). Moreover, $\Psi_{\tau,\alpha}(\wh\cT,\wh\cS)$, with $\reel\alpha\in[-1,0[$, and $\phi_{\tau,0}(\wh\cT,\wh\cS)$ induce, by grading with respect to the monodromy filtration $\rM_\bbullet(\rN_\tau)$, an object of $\MLTr(\wh X,w;-1)^\rp$.
\end{theoreme}

\noindent
(\Cf Chapter~\ref{chap:twistor} for the definition of the categories $\MTr$ and $\MLTr$.) In particular, all conditions of Definition \ref{def:regtwt} are satisfied along the hypersurface $\tau=0$.

This theorem is a generalization of \cite[Th\ptbl5.3]{Bibi96a}, without the $\QQ$-structure however. In fact, we give a precise comparison with nearby cycles of $(\cT,\cS)$ at $t=\infty$ as in \cite[Th\ptbl4.3]{Bibi96a}.

In order to prove Theorem \ref{th:Fstrict}, we need to extend the results of Proposition \ref{prop:Fstrictspe} to objects with sesquilinear pairings.

\subsection{``Positive'' functions of $\hb$}
Recall that we denote by $\DD$ the disc $\module{\hb}\leq1$ and by $\bS$ its boundary. Let $\lambda(\hb)$ be a meromorphic function defined in some neighbourhood of $\bS$. If the neighbourhood is sufficiently small, it has zeros and poles at most on $\bS$. We say that $\lambda$ is ``real'' if it satisfies $\ov\lambda=\lambda$, where (\cf \T\ref{subsubsec:conj}) $\ov\lambda(\hb)$ is defined as $c(\lambda(-1/c(\hb)))$ and $c$ is the usual complex conjugation. For instance, if $\alpha\in\CC$, the function $\hb\mto\alpha\star\hb/\hb$ is ``real''. If $\lambda(\hb)$ is ``real'' and if $\psi$ is a meromorphic function on $\CC$ which is real (in the usual sense, \ie $\psi c=c\psi$), then $\psi\circ\lambda$ is ``real''. In particular, for any $\alpha\in\CC^*$, the function $\hb\mto\Gamma(\alpha\star\hb/\hb)$ is ``real''.

\begin{lemme}\label{lem:positive}
Let $\lambda(\hb)$ be a ``real'' invertible holomorphic function in some neighbourhood of $\bS$. Then there exists an invertible holomorphic function $\mu(\hb)$ in some neighbourhood of $\DD$ such that $\lambda=\pm\mu\ov\mu$ in some neighbourhood of $\bS$. Moreover, such a function $\mu$ is unique up to multiplication by a complex number having modulus equal to $1$.
\end{lemme}

\begin{definition}\label{def:positive}
Let $\lambda$ be as in the lemma. We say that $\lambda$ is ``positive'' if $\lambda=\mu\ov\mu$, with $\mu$ invertible on $\DD$, and ``negative'' if $\lambda=-\mu\ov\mu$.
\end{definition}

\begin{remarque}\label{rem:positive}
Assume that $\lambda$ is a nonzero ``real'' \emph{meromorphic} function in some neighbourhood of $\bS$. Then $\lambda$ can be written as $\prod_i[(\hb-\hb_i)\ov{(\hb-\hb_i)}]^{m_i}\cdot h$ with $\hb_i\in\bS$, $h$~holomorphic invertible near $\bS$ and $\ov h=h$: indeed, one shows that, if $\hb_o\in\bS$, then $\ov{\hb-\hb_o}=(\hb+\hb_o)\cdot(-1/\hb_o\hb)$; therefore, if $\hb_o\in\bS$ is a pole or a zero of $\lambda$ with order $m_o\in\ZZ$, then $-\hb_o$ has the same order, hence the product decomposition of $\lambda$.

It follows from Lemma \ref{lem:positive} that $\lambda=\pm g\ov g$, with $g=\mu\prod_i(\hb-\hb_i)^m_i$, $\hb_i\in\bS$, $m_i\in\ZZ$ and $\mu$ holomorphic invertible on $\DD$. This decomposition is not unique, as one may change some $\hb_i$ with $-\hb_i$. The sign is also non uniquely determined, as we have, for any $\hb_o\in\bS$,
\[
-1=\Big(\frac{\hb-\hb_o}{\hb+\hb_o}\Big)\cdot \ov{\Big(\frac{\hb-\hb_o}{\hb+\hb_o}\Big)}.
\]
Nevertheless, the decomposition and the sign are uniquely defined (up to a multiplicative constant) if we fix a choice of a ``square root'' of the divisor of $\lambda$ so that no two points in the support of this divisor are opposed, and if we impose that the divisor of~$g$ is contained in this ``square root''. The sign does not depend on the choice of such a ``square root''. We say that $\lambda$ is ``positive'' if the sign is $+$, and ``negative'' if the sign is $-$.
\end{remarque}

\begin{proof}[Proof of Lemma \ref{lem:positive}]
One can write $\lambda=\nu\cdot\ov\mu$ with $\mu$ holomorphic invertible near $\DD$ and $\nu$ meromorphic in some neighbourhood of $\DD$ and having poles or zeros at $0$ at most. The function $c(\hb)=\nu/\mu=\ov\nu/\ov\mu$ defines a meromorphic function on $\PP^1$ with divisor supported by $\{0,\infty\}$. Thus, $c(\hb)=c\cdot\hb^k$ with $c\in\CC$ and $k\in\ZZ$, so $\lambda=c\hb^k\mu\ov\mu$. Moreover, the equality $\ov \lambda=\lambda$ implies that $c\in\RR$ and $k=0$. Changing notation for $\mu$ gives $h=\pm\mu\ov\mu$, with $\mu$ invertible on $\DD$.

For uniqueness, assume that $\mu\ov\mu=\pm1$ with $\mu$ holomorphic invertible in some neighbourhood of $\DD$. Then $\pm1/\ov\mu$ is also holomorphic in some neighbourhood of $\module{\hb}\geq1$, so $\mu$ extends as a holomorphic function on $\PP^1$ and thus is constant. This implies that $\mu\ov\mu=1$.
\end{proof}

\begin{lemme}\label{lem:Gamma}
Let $\alpha\in\CC$ be such that $\reel\alpha\in[0,1[$ and $\alpha\neq0$. Then the meromorphic function
\[
\lambda(\hb)=\frac{\Gamma(\alpha\star\hb/\hb)}{\Gamma(1-\alpha\star\hb/\hb)}
\]
is ``real'' and ``positive'' (it is holomorphic invertible near $\bS$ if $\reel\alpha\neq0$).
\end{lemme}

\begin{proof}
That this function is ``real'' has yet been remarked. The only possible pole/zero of $\lambda$ on $\bS$ is $\pm i$, which occurs if there exists $k\in\ZZ$ such that $\reel\alpha+k=0$. It is a simple pole (\resp a simple zero) if $k\geq0$ (\resp $k\leq-1$). As we assume $\reel\alpha\in[0,1[$, the only possibility is when $\reel\alpha=0$, with $k=0$ (hence a pole).

Write $\lambda(\hb)$ as $\Gamma(\alpha\star\hb/\hb)^2\cdot(1/\pi)\sin\pi(\alpha\star\hb/\hb)$. It is then equivalent to showing that $(1/\pi)\sin\pi(\alpha\star\hb/\hb)$ is ``positive'' for $\alpha$ as above.

Write $\alpha=\alpha'+i\alpha''$. The result is clear if $\alpha''=0$, as we then have $\alpha\star\hb/\hb=\alpha'\in{}]0,1[$. We thus assume now that $\alpha''\neq0$.

For any $\beta\in\CC$ with $\beta''\neq0$, we put $b=\dfrac{\beta'+\sqrt{\beta^{\prime2}+\beta^{\prime\prime2}}}{\beta''}$ and we can write
\[
\frac{\beta\star\hb}{\hb}= \frac{\beta''b}{2}\Big(1+\frac{i\hb}{b}\Big)\ov{\Big(1+\frac{i\hb}{b}\Big)}.
\]
If $\alpha$ is as above, we have $n-\alpha',n+\alpha'>0$ for any $n\geq1$ and we put for $n\geq0$
\[
b_n=-\frac{n-\alpha'+\sqrt{(n-\alpha')^2+\alpha^{\prime\prime2}}}{\alpha''},\quad c_n=\frac{n+\alpha'+\sqrt{(n+\alpha')^2+\alpha^{\prime\prime2}}}{\alpha''}.
\]
For $n\geq1$, we have $\module{b_n},\module{c_n}>1$ and
\begin{align*}
\frac{(n-\alpha)\star\hb}{\hb}&=\frac{n-\alpha'+\sqrt{(n-\alpha')^2+\alpha^{\prime\prime2}}}{2}\Big(1+\frac{i\hb}{b_n}\Big)\ov{\Big(1+\frac{i\hb}{b_n}\Big)},
\\
\frac{(n+\alpha)\star\hb}{\hb}&=\frac{n+\alpha'+\sqrt{(n+\alpha')^2+\alpha^{\prime\prime2}}}{2}\Big(1+\frac{i\hb}{c_n}\Big)\ov{\Big(1+\frac{i\hb}{c_n}\Big)}.
\end{align*}
The number
\[
c(\alpha)=\prod_{n\geq1}\frac{(n-\alpha'+\sqrt{(n-\alpha')^2+\alpha^{\prime\prime2}})(n+\alpha'+\sqrt{(n+\alpha')^2+\alpha^{\prime\prime2}})}{4n^2}
\]
is (finite and) positive. On the other hand, as $\dfrac{1}{b_n}+\dfrac{1}{c_n}=-\dfrac{\alpha'\alpha''}{n^2}+O(1/n^3)$, the infinite product
\[
\prod_{n\geq1}\Big(1+\frac{i\hb}{b_n}\Big)\Big(1+\frac{i\hb}{c_n}\Big)
\]
defines an invertible holomorphic function in some neighbourhood of $\DD$. Put
\[
g(\hb)=\Big(\frac{c(\alpha)(\alpha'+\sqrt{\alpha^{\prime2}+\alpha^{\prime\prime2}})}{2}\Big)^{1/2}\cdot \Big(1+\frac{i\hb}{c_0}\Big) \prod_{n\geq1}\Big(1+\frac{i\hb}{b_n}\Big)\Big(1+\frac{i\hb}{c_n}\Big).
\]
Then we have $(1/\pi)\sin\pi(\alpha\star\hb/\hb)=g(\hb)\ov g(\hb)$.
\end{proof}

If $\mu$ is a meromorphic function on some neighbourhood of $\DD$, we denote by $D_\mu$ its divisor on $\DD$. If $\cM$ is a $\cR_\cX$-module, we put $\cM(D_\mu)=\cO_\cX(D_\mu)\otimes_{\cO_\cX}\cM$ with its natural $\cR_\cX$-structure.

\begin{lemme}\label{lem:Dmu}
Let $(\cT,\cS)=(\cM',\cM'',C,\cS)$ be an object of $\MT(X,w)^{\rp}$. Then, for each $\mu$ as above, $(\cM'(D_\mu),\cM''(D_\mu),\mu\ov\mu C,\cS)$ is an object of $\MT(X,w)^{\rp}$ isomorphic to $(\cT,\cS )$.
\end{lemme}

\begin{remarque}
We only assume here that $\cM',\cM''$ are defined in some neighbourhood of $\DD$, and not necessarily on $\Omega_0$, as in Chapter~\ref{chap:twistor}. This does not change the category
$\MT(X,w)^{\rp}$.
\end{remarque}

\begin{proof}
The isomorphism is given by ${\cdot} \mu:\cM'(D_\mu){\to}\cM'$ and ${\cdot}(1/\mu):\cM''{\to}\cM''(D_\mu)$.
\end{proof}

\Subsection{Exponential twist and specialization of a sesquilinear pairing} \label{subsec:expspesesqui}
We now come back to our original situation of \T\ref{subsec:Fouriersettings}. Let $\cT=(\cM',\cM'',C)$ be an object of $\RTriples(X)$. We have defined the object $\cFcT=(\cFcM',\cFcM'',\cFC)$ of $\RTriples(Z)$. If we assume that $\cM',\cM''$ are strict and strictly specializable along $t'=0$, then $\cFcM',\cFcM''$ are strictly specializable along \hbox{$\tau=0$}. Then, for $\reel\alpha\in[-1,0[$, $\Psi_{\tau,\alpha}\cFcT$ is defined as in \T\ref{sec:spesesqui}. Recall (\cf \eqref{eq:NTate}) that we denote by $\cN_\tau:\Psi_{\tau,\alpha}\cFcT \to \Psi_{\tau,\alpha}\cFcT(-1)$ the morphism $(-i\rN_\tau,i\rN_\tau)$. If $\alpha=-1$ (more generally if $\alpha$ is real) we have $\Psi_{\tau,\alpha}\cFcT=\psi_{\tau,\alpha}\cFcT$. We also consider, as in \T\ref{subsec:vanishsesqui}, the vanishing cycle object $\phi_{\tau,0}\cFcT$.

The purpose of this section is to extend Proposition \ref{prop:Fstrictspe}\eqref{prop:Fstrictspe4} to objects of $\RTriples$. It will be convenient to assume, in the following, that $\cM'=\cM'_{\min}$ and $\cM''=\cM''_{\min}$; with such an assumption, we will not have to define a sesquilinear pairing on the minimal extensions used in Proposition \ref{prop:Fstrictspe}\eqref{prop:Fstrictspe4}, as we can use $C$.

\begin{proposition}\label{prop:Fspetw}
For $\cT$ as above, we have isomorphisms in $\RTriples(X)$:
\begin{align*}
\big(\Psi_{\tau,\alpha}\cFcT,\cN_\tau\big)
&\isom i_{\infty,+}\big(\Psi_{t',\alpha}\cT,\cN_{t'}\big),\quad \forall\alpha\neq-1 \text{ with }\reel\alpha\in[-1,0[,\\
(\phi_{\tau,0}\cFcT,\cN_{\tau})&\isom i_{\infty,+}\big(\psi_{t',-1} \cT,\cN_{t'}\big),
\end{align*}
and an exact sequence
\[
0\to i_{\infty,+}\ker \cN_{t'}\to \ker \cN_\tau\to\cT\to0
\]
inducing an isomorphism $P\gr_0^{\rM}\psi_{\tau,-1}\cFcT\isom\cT$.
\end{proposition}

\begin{corollaire}\label{cor:MTpsiF}
Assume that $\cT$ is an object of $\MTr(X,w)$ (\resp $(\cT,\cS)$ is an object of $\MTr(X,w)^\rp$). Then, for any $\alpha\in\CC$ with $\reel\alpha\in[-1,0[$, $(\Psi_{\tau,\alpha}\cFcT,\cN_\tau)$ induces by gradation an object of $\MLTr(X,w;-1)$ (\resp an object of $\MLTr(X,w;-1)^\rp$).
\end{corollaire}

\begin{proof}[Proof of Corollary \ref{cor:MTpsiF}]
Suppose that Proposition \ref{prop:Fspetw} is proved. Assume first that $\cT$ is an object of $\MTr(X,w)$. Then, by definition, $i_{\infty,+}\big(\gr_\bbullet^{\rM}\Psi_{t',\alpha}\cT,\gr_{-2}^{\rM} \cN_{t'}\big)$ is an object of $\MLTr(X,w;-1)$ for any $\alpha$ with $\reel\alpha\in[-1,0[$; therefore, so is $\big(\gr_\bbullet^{\rM}\Psi_{\tau,\alpha}\cFcT,\gr_{-2}^{\rM}\cN_\tau\big)$ for any such $\alpha\neq-1$. When $\alpha=-1$, as $\cFcM',\cFcM''$ are equal to their minimal extension along $\tau=0$ (\cf Proposition \ref{prop:Fstrictspe}) the morphism
\[
\cCan:\big(\psi_{\tau,-1}\cFcT,\rM_\bbullet(\cN_\tau)\big)\to\big(\phi_{\tau,0}\cFcT(-1/2),\rM_{\bbullet-1}(\cN_\tau)\big),
\]
(\cf \T\ref{subsec:vanishsesqui}) is onto. It is strictly compatible with the monodromy filtrations (\cf \cite[Lemme 5.1.12]{MSaito86}), and induces an isomorphism $P\gr_\ell^{\rM}\psi_{\tau,-1}\cFcT\isom P\gr_{\ell-1}^{\rM}\phi_{\tau,0}\cFcT(-1/2)$ for any $\ell\geq1$, hence an isomorphism
\[
P\gr_\ell^{\rM}\psi_{\tau,-1}\cFcT\isom i_{\infty,+}P\gr_{\ell-1}^{\rM}\psi_{t',-1}\cT(-1/2).
\]
By assumption on $\cT$, the right-hand term is an object of $\MTr(X,w+\ell)$, hence so is the left-hand term. Moreover, $P\gr_0^{\rM}\psi_{\tau,-1}\cFcT\simeq\cT$ is in $\MTr(X,w)$. This gives the claim when $\alpha=-1$.

In the polarized case, we can reduce to the case where $w=0$, $\cM'=\cM''$, $\cS=(\id,\id)$ and $C^*=C$. Then these properties are satisfied by the objects above, and the polarizability on the $\tau$-side follows from the polarizability on the $t'$-side.
\end{proof}

The proof of the proposition will involve the computation of a Mellin transform with kernel given by a function $I_\chi(t,s,\hb)$. We first analyze this Mellin transform.

\subsubsection*{The function $\Ichih(t,s,\hb)$}
Let $\chih\in C_c^\infty(\Afuh,\RR)$ be such that $\chih(\tau)\equiv1$ near $\tau=0$. For any $\hb\in\bS$, $t\in\Afu$ and $s\in\CC$ such that $\reel (s+1)>0$, put
\begin{equation}\label{eq:Ichiab}
\Ichih(t,s,\hb)=\int_{\Afuh}e^{\hb\ov{t\tau}-t\tau/\hb}\module{\tau}^{2s}\chih(\tau)\itwopi\,d\tau\wedge d\ov\tau.
\end{equation}
We also write $\Ichih(t',s,\hb)$ when working in the coordinate $t'$ on $\PP^1$. We will use the following coarse properties (they are similar to the properties described for the function $\wh I_\chi$ of \T\ref{subsec:vanishsesqui}).
\begin{enumerate}
\item\label{item:Ichi1}
Denote by $I_{\chih,k,\ell}(t,s,\hb)$ ($k,\ell\in\ZZ$) the function obtained by integrating $\module{\tau}^{2s}\tau^k\ov\tau^\ell$. Then, for any $s\in\CC$ with $\reel(s+1+(k+\ell)/2)>0$ and any $\hb\in\bS$, the function $(t,s,\hb)\to I_{\chih,k,\ell}(t,s,\hb)$ is $C^\infty$, depends holomorphically on $s$, and satisfies $\lim_{t\to\infty}I_{\chih,k,\ell}(t,s,\hb)=0$ locally uniformly with respect to $s,\hb$.
\item\label{item:Ichi2}
We have
\begin{align*}
tI_{\chih,k,\ell}&=\hb(s+k)I_{\chih,k-1,\ell}+\hb I_{\partial\chih/\partial\tau,k,\ell}&\partiall_tI_{\chih,k,\ell}&=-I_{\chih,k+1,\ell}\\
\ov tI_{\chih,k,\ell}&=\ov\hb(s+\ell)I_{\chih,k,\ell-1}+\ov\hb I_{\partial\chih/\partial\ov\tau,k,\ell}&\ov\partiall_tI_{\chih,k,\ell}&=-I_{\chih,k,\ell+1},
\end{align*}
where the equalities hold on the common domain of definition (with respect to $s$) of the functions involved. Notice that the functions $I_{\partial\chih/\partial\tau,k,\ell}$ and $I_{\partial\chih/\partial\ov\tau,k,\ell}$ are $C^\infty$ on $\PP^1\times\CC\times\bS$, depend holomorphically on $s$, and are infinitely flat at $t=\infty$.

It follows that, for $\reel(s+1)>0$, we have
\begin{equation}\label{eq:tdtIchih}
\begin{aligned}
t\partiall_t\Ichih&=-\hb(s+1)\Ichih+\hb I_{\partial\chih/\partial\tau,1,0},\\
\ov{t\partiall_t}\Ichih&=-\ov\hb(s+1)\Ichih+\ov\hb I_{\partial\chih/\partial\ov\tau,0,1}.
\end{aligned}
\end{equation}

\item\label{item:Ichi3}
Moreover, for any $p\geq0$, any $s\in\CC$ with \hbox{$\reel (s+1+(k+\ell)/2)>p$} and any $\hb\in\bS$, all derivatives up to order $p$ of $I_{\chih,k,\ell}(t',s,\hb)$ with respect to $t'$ tend to $0$ when $t'\to0$, locally uniformly with respect to $s,\hb$; therefore, $I_{\chih,k,\ell}(t,s,\hb)$ extends as a function of class $C^p$ on $\PP^1\times\{\reel (s+1+(k+\ell)/2)>p\}\times\bS$, holomorphic with respect to $s$.
\end{enumerate}

\subsubsection*{Mellin transform with kernel $\Ichih(t,s,\hb)$}
We will work near $\hb_o\in\bS$. For any local sections $\mu',\mu''$ of $\cM',\cM''$ and any $C^\infty$ relative form $\varphi$ of maximal degree on $X\times\bS$ with compact support contained in the open set where $\mu',\mu''$ are defined, the function
\[
(s,\hb)\mto\big\langle C(\mu',\ov{\mu''}),\varphi \Ichih(t,s,\hb)\big\rangle
\]
is holomorphic with respect to $s$ for $\reel s\gg0$ (according to \eqref{item:Ichi1}), continuous with respect to $\hb$. One shows as in Lemma \ref{lem:polesI}, using \eqref{item:Ichi3}, that it extends as a meromorphic function on the whole complex plane, with poles on sets $s=\alpha\star\hb/\hb$.

This result can easily be extended to local sections $\mu',\mu''$ of $\wt{\cM'},\wt{\cM''}$: indeed, this has to be verified only near $t=\infty$; there exists $p\geq0$ such that, in the neighbourhood of the support of $\varphi$, $t^{\prime p}\mu',t^{\prime p}\mu''$ are local sections of $\cM',\cM''$; apply then the previous argument to the kernel $\module{t}^{2p}\Ichih(t,s,\hb)$. In the following, we will write $\big\langle C(\mu',\ov{\mu''}),\varphi \Ichih(t,s,\hb)\big\rangle$ instead of $\big\langle C(t^{\prime p}\mu',t^{\prime p}\ov{\mu''}),\varphi\module{t}^{2p} \Ichih(t,s,\hb)\big\rangle$ near $t=\infty$.

\begin{lemme}\label{lem:outinfty}
Assume that $\varphi$ is compactly supported on $(X\moins\infty)\times\bS$. Then, for $\mu',\mu''$ as above, we have
\[
\res_{s=-1}\big\langle C(\mu',\ov{\mu''}),\varphi \Ichih(t,s,\hb)\big\rangle=\big\langle C(\mu',\ov{\mu''}),\varphi\big\rangle.
\]
\end{lemme}

\begin{proof}
The function $(s+1)\Ichih(t,s,\hb)$ can be extended to the domain $\reel(s+1)>-1/2$ as $C^\infty$ function of $(t,s,\hb)$, holomorphic with respect to~$s$: use \eqref{item:Ichi2} with $k=1$, $\ell=0$ to write $(s+1)\Ichih(t,s,\hb)=(t/\hb)I_{\chih,1,0}-I_{\partial\chih/\partial\tau,1,0}$. It is then enough to show that this $C^\infty$ function, when restricted to $s=-1$, is identically equal to $1$. It amounts to proving that, for any $t,\hb$, $\lim_{\substack{s\to-1\\\reel s>-1}}\big[(s+1)\Ichih(t,s,\hb)\big]=1$.
For $\reel s>-1$ we have $\Ichih(t,s,\hb)=J(t,s,\hb)+H_{\chih}(t,s,\hb)$, with
\[
J(t,s,\hb)=\int_{\module{\tau}\leq1}e^{-2i\im t\tau/\hb}\module{\tau}^{2s}\itwopi\,d\tau\wedge d\ov\tau,
\]
and $H_{\chih}$ extends as a $C^\infty$ function on $\Afu\times\CC\times\bS$, holomorphic with respect to~$s$. It is therefore enough to work with $J(t,s,\hb)$ instead of $\Ichih$. We now have
\begin{align*}
J(t,s,\hb)&=\module{t}^{-2(s+1)}\int_{\module{u}<\module{t}}e^{-2i\im u}\module{u}^{2s}\itwopi\,du\wedge d\ov u\\
&=\frac1\pi\,\module{t}^{-2(s+1)}\int_0^{2\pi}\hspace*{-2mm} \int_0^{\module{t}}e^{-2i\rho\sin\theta}\rho^{2s+1}\,d\rho d\theta.
\end{align*}
Now, integrating by part, we get
\[
\int_0^{\module{t}}e^{-2i\rho\sin\theta}\rho^{2s+1}\,d\rho=\frac{\module{t}^{2s+2}e^{-2i\module{t}\sin\theta}}{2s+2}+\frac{2i\sin\theta}{2s+2}\int_0^{\module{t}}e^{-2i\rho\sin\theta}\rho^{2s+2}\,d\rho,
\]
and the second integral is holomorphic near $s=-1$. Therefore,
\[
(s+1)J(t,s)=
\frac{\module{t}^{-2(s+1)}}{2\pi}
\int_0^{2\pi}\Big[\module{t}^{2s+2}e^{-2i\module{t}\sin\theta}
+2i\sin\theta\int_0^{\module{t}}e^{-2i\rho\sin\theta}\rho^{2s+2}\,d\rho\Big]d\theta.
\]
Taking $s\to-1$ gives
\[
\lim_{\substack{s\to-1\\\reel s>-1}}[(s+1)J(t,s)]=
\frac1{2\pi}\int_0^{2\pi}\Big[e^{-2i\module{t}\sin\theta}
+2i\sin\theta\int_0^{\module{t}}e^{-2i\rho\sin\theta}\,d\rho\Big] d\theta.
\]
Now,
\[
2i\sin\theta\int_0^{\module{t}}e^{-2i\rho\sin\theta}\,d\rho
=-\int_0^{\module{t}}\frac{d}{d\rho}\big(e^{-2i\rho\sin\theta}\big)d\rho
=1-e^{-2i\module{t}\sin\theta},
\]
hence $\lim_{\substack{s\to-1\\\reel s>-1}}[(s+1)J(t,s)]=1$.
\end{proof}

\begin{remarque}\label{rem:Jchih}
To simplify notation, we now put
\[
\Jchih(t,s,\hb)=\frac1{\Gamma(s+1)}\,\Ichih(t,s,\hb).
\]
Using \eqref{item:Ichi2} as in the previous lemma, one obtains that there exists a $C^\infty$ function on $\Afu\times\CC\times\bS$, holomorphic with respect to $s$, which coincides with $\Jchih$ when $\reel(s+1)>0$. This implies that, when the support of $\varphi$ does not cut $\infty$, the meromorphic function $s\mto\big\langle C(\mu',\ov{\mu''}),\varphi \Jchih(t,s,\hb)\big\rangle$ is entire.
\end{remarque}

We now work near $\infty$ with the coordinate $t'$. Assume that $\mu'$ is a local section of $V_{a_1+1}^{(\hb_o)}\wt\cM'$ and that $\mu''$ is a local section of $V_{a_2+1}^{(-\hb_o)}\wt\cM''$. Assume moreover that the class of $\mu'$ in $\gr_{a_1+1}^{(\hb_o)}\wt\cM'$ is in $\psi_{t',\alpha_1+1}\wt\cM'$, and that the class of $\mu''$ in $\gr_{a_2+1}^{(-\hb_o)}\wt\cM''$ is in $\psi_{t',\alpha_2+1}\wt\cM''$. Then on proves as in Lemma \ref{lem:polesI} that $\big\langle C(\mu',\ov{\mu''}),\varphi \Jchih(t',s,\hb)\big\rangle$ has poles on sets $s=\gamma\star\hb/\hb$ with $\gamma$ such that $2\reel\gamma<a_1+a_2$ or $\gamma=\alpha_1=\alpha_2$.

Let us then consider the case where $\alpha_1=\alpha_2\defin\alpha$. Then, if $\psi$ has compact support and vanishes along $t'=0$, the previous result shows that $\big\langle C(\mu',\ov{\mu''}),\psi \Jchih(t',s,\hb)\big\rangle$ has no pole along $s=\alpha\star\hb/\hb$. It follows that $\res_{s=\alpha\star\hb/\hb}\big\langle C(\mu',\ov{\mu''}),\varphi \Jchih(t',s,\hb)\big\rangle$ only depends on the restriction of $\varphi$ to $t'=0$; in other words, it is the direct image of a distribution on $t'=0$ by the inclusion $i_\infty$. We will identify this distribution with $\psi_{t',\alpha+1}C$. We will put
\[
i^*_\infty\varphi=\frac{\varphi_{|\infty}}{\itwopi dt'\wedge d\ov{t'}}.
\]

\begin{lemme}\label{lem:Jchi}
For any $\alpha\in\CC$ and $\mu',\mu''$ lifting local sections $[\mu'],[\mu'']$ of $\psi_{t',\alpha+1}\wt\cM'$, $\psi_{t',\alpha+1}\wt\cM''$, we have, when the support of $\varphi$ is contained in the open set where $\mu',\mu''$ are defined,
\[
\res_{s=\alpha\star\hb/\hb}\big\langle C(\mu',\ov{\mu''}),\varphi \Jchih(t',s,\hb)\big\rangle= \frac{1}{\Gamma(-\alpha\star\hb/\hb)}\, \big\langle \psi_{t',\alpha+1}C([\mu'],\ov{[\mu'']}),i^*_\infty\varphi\big\rangle.
\]
\end{lemme}

\begin{proof}
Let $\chi(t')$ be a $C^\infty$ function which has compact support and is $\equiv1$ near $t'=0$. As $\varphi-i^*_\infty\varphi\wedge\chi(t')\itwopi dt'\wedge d\ov{t'}$ vanishes along $t'=0$, the left-hand term in the lemma is equal to
\begin{equation}\label{eq:CJchih}
\res_{s=\alpha\star\hb/\hb}\big\langle C(\mu',\ov{\mu''}), \Jchih(t',s,\hb)i^*_\infty\varphi\wedge\chi(t')\itwopi dt'\wedge d\ov{t'}\big\rangle.
\end{equation}
On the other hand, by definition of $\psi_{t',\alpha+1}C$ (\cf \eqref{eq:psiC}), the right-hand term is equal to
\begin{equation}\label{eq:Cts}
\res_{s=\alpha\star\hb/\hb}\frac{1}{\Gamma(-s)}\big\langle C(\mu',\ov{\mu''}), \module{t'}^{2(s+1)}i^*_\infty\varphi\wedge\chi(t')\itwopi dt'\wedge d\ov{t'}\big\rangle.
\end{equation}

Put $\wtJchih(t,s,\hb)=\mt^{2(s+1)}\Jchih(t,s,\hb)$. Then, by \eqref{eq:tdtIchih} expressed in the coordinate~$t'$, we have
\[
t'\frac{\partial\wtJchih}{\partial t'}=-\wt J_{\partial\chih/\partial\tau,1,0},\quad \ov{t'}\frac{\partial\wtJchih}{\partial \ov{t'}}=-\wt J_{\partial\chih/\partial\ov\tau,0,1},
\]
and both functions $\wt J_{\partial\chih/\partial\tau,1,0}$ and $\wt J_{\partial\chih/\partial\ov\tau,0,1}$ extend as $C^\infty$ functions, infinitely flat at $t'=0$ and holomorphic with respect to $s\in\CC$. Put
\[
\wtKchih(t',s,\hb)=-\int_0^1\big[\wt J_{\partial\chih/\partial\tau,1,0}(\lambda t',s,\hb)+\wt J_{\partial\chih/\partial\ov\tau,0,1}(\lambda t',s,\hb)\big]d\lambda.
\]
Then $\wtKchih$ is of the same kind. Notice now that, for any $s\in\CC$ with $\reel(s+1)\in{}]0,1/4[$ and any $\hb\in\bS$, we have
\begin{equation}\label{eq:Bessel}
\lim_{t\to\infty}(\module{t}^{2(s+1)}\Ichih(t,s,\hb))=\frac{\Gamma(s+1)}{\Gamma(-s)}.
\end{equation}
[Let us sketch the proof of this statement. We assume for instance that $\chih\equiv1$ when $\module{\tau}\leq1$. We can replace $\Ichih(t,s,\hb)$ with
\[
\int_{\module{\tau}\leq1}e^{\hb\ov{t\tau}-t\tau/\hb}\module{\tau}^{2s}\itwopi\,d\tau\wedge d\ov\tau
\]
without changing the limit, and we are reduced to computing
\[
\frac1\pi\int_0^{2\pi}\int_0^\infty e^{-2i\rho\sin\theta}\rho^{2s+1}d\rho\,d\theta.
\]
Using the Bessel function $J_0(r)=\frac1{2\pi}\int_0^{2\pi}e^{-ir\sin\theta}d\theta$, this integral is written as
\[
2\int_0^\infty\rho^{2s+1}J_0(2\rho)d\rho=\frac1{2^{2s+1}}\int_0^\infty r^{2s+1}J_0(r)dr,
\]
and it is known (\cf \cite[\T13.24, p\ptbl 391]{Watson22}) that, on the strip $\reel(s+1)\in{}]0,1/4[$, the last integral is equal to $2^{2s+1}\Gamma(s+1)/\Gamma(-s)$.]

\enlargethispage{7mm}%
On this strip, we can therefore write $\wtJchih(t',s,\hb)=\wtKchih(t',s,\hb)+1/\Gamma(-s)$, by Taylor's formula. Hence, on the strip $\reel(s+1)\in{}]0,1/4[$ and, by applying $(\partiall_{t'}\ov\partiall_{t'})^p$ as for \eqref{eq:IGamma}, on any strip $\reel(s+1)\in{}]p,p+1/4[$ for $p\geq0$, we have
\[
\Jchih(t',s,\hb)=\frac{\module{t'}^{2(s+1)}}{\Gamma(-s)}+\Kchih^{(p)}(t',s,\hb),
\]
where $\Kchih^{(p)}(t',s,\hb)=(\partiall_{t'}\ov\partiall_{t'})^p[\module{t'}^{2(s+1)}\wtKchih(t',s,\hb)]$ satisfies the same properties as $\wtKchih$ does when $t'\to0$. This implies that the function
\[
s\mto \big\langle C(\mu',\ov{\mu''}), \Kchih^{(p)}(t',s,\hb)i^*_\infty\varphi\wedge\chi(t')\itwopi dt'\wedge d\ov{t'}\big\rangle
\]
is entire for any $\hb\in\bS$. Hence, there exists an entire function of $s$ such that the difference of the meromorphic functions considered in \eqref{eq:CJchih} and \eqref{eq:Cts}, when restricted to the strip $\reel(s+1)\in{}]p,p+1/4[$ (with $p$ large enough so that they are holomorphic on the strip), coincides with this entire function. This difference is therefore identically equal to this entire function of $s$, and \eqref{eq:CJchih} and \eqref{eq:Cts} coincide. This proves the lemma.
\end{proof}

\subsubsection*{Proof of Proposition \ref{prop:Fspetw}}
We will work near $\hb_o\in\bS$. By definition (\cf \eqref{eq:psiC}), given any local sections $[m'],[m'']$ of $\psi_{\tau,\alpha}\cFcM',\psi_{\tau,\alpha}\cFcM''$ and local liftings $m',m''$ in $V_{a'}\cFcM',V_{a''}\cFcM''$ with $a'=\ell_{\hb_o}(\alpha)$ and $a''=\ell_{-\hb_o}(\alpha)$, we have, for any $C^\infty$ relative form $\varphi$ of maximal degree on $X\times\bS$,
\begin{equation}\label{eq:psiFC}
\big\langle\psi_{\tau,\alpha}\cFC([m'],\ov{[m'']},\varphi\big\rangle=\res_{s=\alpha\star\hb/\hb}\big\langle\cFC(m',\ov{m''}),\varphi\module{\tau}^{2s}\chih(\tau)\itwopi\,d\tau\wedge d\ov\tau\big\rangle,
\end{equation}
where $\chih\equiv1$ near $\tau=0$. In particular, for sections $m',m''$ of the form $\mu'\otimes\ccE^{-t\tau/\hb}$, $\mu''\otimes\ccE^{-t\tau/\hb}$ with $\mu',\mu''$ local sections of $\wt\cM$, the definition of $\cFC$ implies that the right-hand term above can be written as
\begin{equation}\label{eq:psiFCmu}
\res_{s=\alpha\star\hb/\hb}\big\langle C(\mu',\ov{\mu''}),\varphi \Ichih(t,s,\hb)\big\rangle.
\end{equation}
[Here, we mean that both functions
$$
\big\langle C(\mu',\ov{\mu''}),\varphi \Ichih(t,s,\hb)\big\rangle\quad\text{and}\quad \big\langle\cFC(\mu'\otimes\ccE^{-t\tau/\hb},\ov{\mu''\otimes\ccE^{-t\tau/\hb}}),\varphi\module{\tau}^{2s}\chih(\tau)\itwopi\,d\tau\wedge d\ov\tau\big\rangle,
$$
\emph{a priori} defined for $\reel s\gg0$, are extended as meromorphic functions of $s$ on the whole complex plane.] Moreover, by $\cR_\cX$-linearity, it is enough to prove Proposition \ref{prop:Fspetw} for such sections.

\begin{proof}[Proof of Proposition \ref{prop:Fspetw} away from $\infty$]
This is the easy part of the proof. We only have to consider $\alpha=-1$ and, for $\varphi$ compactly supported on $(X\moins\infty)\times\bS$, we are reduced to proving that
\[
\res_{s=-1}\big\langle C(\mu',\ov{\mu''}),\varphi \Ichih(t,s,\hb)\big\rangle=\big\langle C(\mu',\ov{\mu''}),\varphi\big\rangle,
\]
for local sections $\mu',\mu''$ of $\wt\cM',\wt\cM''$. This is Lemma \ref{lem:outinfty}
\end{proof}

\begin{proof}[Proof of Proposition \ref{prop:Fspetw} near $\infty$ for $\alpha\neq-1,0$]
The question is local on $\DD$. We can compute \eqref{eq:psiFC} by using liftings of $m',m''$ in $\gr_{a'+1}^U\cFcM',\gr_{a''+1}^U\cFcM''$, according to \eqref{eq:idefines}. By $\cR$-linearity, we only consider sections $m'=t^{\prime-1}\mu'\otimes\ccE^{-t\tau/\hb}$, \hbox{$m''=t^{\prime-1}\mu''\otimes\ccE^{-t\tau/\hb}$}, where $\mu'$ is a local section of $V_{a'}\wt\cM'$ and $\mu''$ of $V_{a''}\wt\cM''$. According to \eqref{eq:psiFCmu}, we have
\[
\big\langle\psi_{\tau,\alpha}\cFC([m'],\ov{[m'']},\varphi\big\rangle= \res_{s=\alpha\star\hb/\hb}\big\langle C(t^{\prime-1}\mu',\ov{t^{\prime-1}\mu''}),\varphi \Ichih(t,s,\hb)\big\rangle,
\]
and, from Lemma \ref{lem:Jchi}, this is
\begin{multline*}
\frac{\Gamma(1+\alpha\star\hb/\hb)}{\Gamma(-\alpha\star\hb/\hb)}\big\langle\psi_{t',\alpha+1}C([t^{\prime-1}\mu'],\ov{[t^{\prime-1}\mu''}]),i_\infty^*\varphi\big\rangle\\
= \frac{\Gamma(1+\alpha\star\hb/\hb)}{\Gamma(-\alpha\star\hb/\hb)}\big\langle\psi_{t',\alpha}C([\mu'],\ov{[\mu''}]),i_\infty^*\varphi\big\rangle.
\end{multline*}
By Lemma \ref{lem:Gamma} and its proof, we have $\sfrac{\Gamma(1+\alpha\star\hb/\hb)}{\Gamma(-\alpha\star\hb/\hb)}=\mu\ov\mu$, with \hbox{$D_\mu=-D_\alpha$} (recall that $D_\alpha$ was defined in Proposition \ref{prop:Fstrictspe}\eqref{prop:Fstrictspe4a}), as we assume $\reel\alpha\in[-1,0[$. We then apply Lemma \ref{lem:Dmu}.
\end{proof}

\begin{proof}[Proof of Proposition \ref{prop:Fspetw} near $\infty$ for $\alpha=0$]
By the same reduction as above, we consider local sections $m'_0,m''_0$ of $V_0\cFcM',V_0\cFcM''$ of the form $m'_0=\mu'_1\otimes\ccE^{-t\tau/\hb}$, $m''_0=\mu''_1\otimes\ccE^{-t\tau/\hb}$, where $\mu'_1,\mu''_1$ are local sections of $V_1\wt\cM',V_1\wt\cM''$. By Lemma \ref{lem:phiC}, we have
\begin{align*}
\big\langle\phi_{\tau,0}\cFC([m'_0],\ov{[m''_0]}),\varphi\big\rangle
&=\res_{s=0}\frac{-1}{s}\big\langle \cFC(m'_0,\ov{m''_0}),\varphi\wedge\module{\tau}^{2s}\chih(\tau)\itwopi d\tau\wedge d\ov \tau\big\rangle\\
&=\res_{s=0}\frac{-1}{s}\big\langle C(\mu'_1,\ov{\mu''_1}),\varphi \Ichih(t,s,\hb)\big\rangle.
\end{align*}
By Lemma \ref{lem:Jchi}, using that $(-\Gamma(s+1)/s\Gamma(-s))_{|s=0}=1$, this expression is
\[
\res_{s=0}\big\langle C(\mu'_1,\ov{\mu''_1}),\module{t'}^{2(s+1)}\varphi\big\rangle =\res_{s=-1}\big\langle C(\mu'_{-1},\ov{\mu''_{-1}}),\module{t'}^{2s}\varphi\big\rangle
\]
if we put $\mu_{-1}=t^{\prime2}\mu_1$. This is $\big\langle \psi_{t',-1}C([\mu'_{-1}],\ov{[\mu''_{-1}]}),i_\infty^*\varphi\big\rangle$.
\end{proof}

\begin{proof}[Proof of Proposition \ref{prop:Fspetw} near $\infty$ for $\alpha=-1$]
Let us first explain how $\psi_{\tau,-1}\cFC$ is defined and how it induces a sesquilinear pairing on $P\gr_0^{\rM}\psi_{\tau,-1}\cFcM',P\gr_0^{\rM}\psi_{\tau,-1}\cFcM''$.

In order to compute $\psi_{\tau,-1}\cFC$, we lift local sections $[m'],[m'']$ of $\psi_{\tau,-1}\cFcM',\psi_{\tau,-1}\cFcM''$ in $U_0\cFcM',U_0\cFcM''$ and compute \eqref{eq:psiFC} for $\alpha=-1$. We know, by Lemma \ref{lem:polesI}, that this is well defined.

To compute the induced form on $P\gr_0^{\rM}$, we use \eqref{eq:suiteexU} and \eqref{eq:suiteexUc} and, arguing as above, we have to consider sections $m',m''$ of $U_{<0}\cFcM',U_{<0}\cFcM''$. We are then reduced to proving that, for local sections $\mu',\mu''$ of $V_{<0}\wt\cM',V_{<0}\wt\cM''$, we have
\[
\res_{s=-1} \frac{\Gamma(s+1)}{\Gamma(-s)}\big\langle C(\mu',\ov{\mu''}),\module{t'}^{2(s+1)}\varphi\big\rangle= \big\langle C(\mu',\ov{\mu''}),\varphi\big\rangle.
\]
By Lemma \ref{lem:polesI}, the meromorphic function $s\mto\big\langle C(\mu',\ov{\mu''}),\module{t'}^{2(s+1)}\varphi\big\rangle$ has poles along sets $s+1=\gamma\star\hb/\hb$ with $\reel\gamma<0$. For such a $\gamma$ and for $\hb\in\bS$, we cannot have $\gamma\star\hb/\hb=0$. Therefore, $s\mto\big\langle C(\mu',\ov{\mu''}),\module{t'}^{2(s+1)}\varphi\big\rangle$ is holomorphic near $s=-1$ and its value at $s=-1$ is $\big\langle C(\mu',\ov{\mu''}),\varphi\big\rangle$. The assertion follows.
\end{proof}

\subsection{Proof of Theorem \ref{th:Fstrict}}
We first reduce to weight~$0$, and assume that $w=0$. It is then possible to assume that $(\cT,\cS)=(\cM,\cM,C,\id)$. We may also assume that $\cM$ has strict support. Then, in particular, we have $\cM=\wt\cM_{\min}$, as defined above.

According to Corollary \ref{cor:MTpsiF} (and to Proposition \ref{prop:Fspetw} for $\phi_{\tau,0}$), we can apply the arguments given in \T\ref{sec:nmnp} to the direct image by $q$.\qed

Notice that we also get:

\begin{corollaire}\label{cor:Ftw}
Let $(\cT,\cS)=(\cM',\cM'',C,\cS)$ be an object of $\MTr(X,w)^\rp$. Then, we have isomorphisms in $\RTriples(X)$:
\tagdroite
\begin{align*}
\big(\Psi_{\tau,\alpha}\wh \cT,\cN_\tau\big)
&\isom \big(\Psi_{t',\alpha}\cT,\cN_{t'}\big),\quad \forall\alpha\neq-1 \text{ with }\reel\alpha\in[-1,0[,\\
(\phi_{\tau,0}\wh \cT,\cN_{\tau})&\isom \big(\psi_{t',-1} \cT,\cN_{t'}\big).\tag*{\qed}
\end{align*}
\taggauche
\end{corollaire}

\begin{remarque}
Let us indicate some shortcut to obtain the S-decomposability of $\wh\cM$ when $Y$ is reduced to a point. By Proposition \ref{prop:Fstrictspe}, we have exact sequences
\begin{gather*}
0\to\ker\rN_\tau\to\psi_{\tau,-1}\cFcM\To{\can_\tau}i_{\infty,+}\psi_{t',-1}\cM\to0,\\
0\to i_{\infty,+}\psi_{t',-1}\cM\To{\var_\tau}\psi_{\tau,-1}\cFcM\to\coker\rN_\tau\to0,
\end{gather*}
and
\begin{gather*}
0\to i_{\infty,+}\ker\rN_{t'}\to\ker\rN_\tau\to\cM\to0\\
0\to\cM\to\coker\rN_\tau\to i_{\infty,+}\coker\rN_{t'}\to0.
\end{gather*}
It follows that $\cH^1q_+\ker\can_\tau=\cH^1q_+\cM$ and $\cH^{-1}q_+\coker\var_\tau=\cH^{-1}q_+\cM$. By the first part of the proof, we then have exact sequences
\begin{gather*}
\psi_{\tau,-1}\wh\cM\To{\can_\tau}\psi_{\tau,0}\wh\cM=\psi_{t',-1}\cM\to\cH^1q_+\cM\to0\\
0\to\cH^{-1}q_+\cM\to\psi_{t',-1}\cM=\psi_{\tau,0}\wh\cM\To{\var_\tau}\psi_{\tau,-1}\wh\cM.
\end{gather*}
Therefore, if $q_+\cM$ has cohomology in degree $0$ only, $\wh\cM$ is a minimal extension along $\tau=0$. Such a situation occurs if $Y$ is reduced to a point, so that $X=\PP^1$: indeed, as $(\cT,\cS)$ is an object of $\MTr(\PP^1,0)^\rp$, we can assume that $\cT$ is simple (\cf Proposition \ref{prop:semisimplicity}); denote by $M$ the restriction of $\cM$ to $\hb=1$, \ie $M=\cM/(\hb-1)\cM$; by Theorem \ref{th:ssimplecourbes}, $M$ is an irreducible regular holonomic $\cD_{\PP^1}$-module;
\begin{itemize}
\item
if $M$ is not isomorphic to $\cO_{\PP^1}$, then $q_+M$ has cohomology in degree $0$ only [use duality to reduce to the vanishing of $\cH^{-1}q_+M$, which is nothing but the space of global sections of the local system attached to $M$ away from its singular points]; By Theorem \ref{th:imdirtwistor}, each cohomology $\cH^jq_+\cM$ is strict and its fibre at $\hb=1$ is $\cH^jq_+M$; therefore, $\cH^jq_+\cM=0$ if $j\neq0$;
\item
otherwise, $M$ is isomorphic to $\cO_{\PP^1}$ with its usual $\cD_{\PP^1}$ structure, and $\wh\cM$ is $\cO_{\cP^1}$, so $\wh\cM$ is supported on $\tau=0$ and $\psi_{\tau,-1}\wh\cM=0$;
\end{itemize}
in conclusion, the S-decomposability of $\wh\cM$ along $\tau=0$ is true in both cases.
\end{remarque}

\backmatter
\nocite{CIRM81,Mebkhout87,CIMPA90-2,K-M-S-S98}
\bibliographystyle{smfplain}
\bibliography{strings,sabbah-twistor_041220}
\renewcommand{\indexname}{Notation}
\printindex
\end{document}